\documentclass[final,3p,times]{elsarticle}
\usepackage{graphicx}
\usepackage{color,soul} 
\usepackage{array}
\usepackage{amsmath}
\usepackage{amsfonts}
\usepackage{amssymb}
\usepackage{psfrag}
\usepackage{epstopdf}
\usepackage{subcaption}

\usepackage{caption}
\usepackage{booktabs}
\usepackage[hidelinks]{hyperref}
\usepackage{latexsym}
\usepackage{kantlipsum}
\usepackage{balance}
\usepackage{lscape}
\usepackage{rotating}
\usepackage{multicol}
\usepackage{multirow,bigdelim}
\usepackage{color, colortbl}
\usepackage{threeparttable}
\usepackage{float}
\usepackage[ruled]{algorithm2e}
\usepackage{footnote}
\usepackage{enumitem}
\usepackage{comment}
\setlist{  
  listparindent=\parindent,
  parsep=0pt,
}

\usepackage{sfmath}

\usepackage{afterpage}

\usepackage{tikz}
\usetikzlibrary{patterns}
\usetikzlibrary{calc}

\definecolor{jetblue}{RGB}{35,120,255}
\definecolor{jetorange}{RGB}{255,255,45}
\definecolor{jetdarkred}{RGB}{195,35,125}

\newcommand{\parder}[2]{\frac{\partial #1}{\partial #2}}

\newcommand{\dparder}[2]{\dfrac{\partial #1}{\partial #2}}

\newcommand{\eval}[2][\right]{\relax\ifx#1\right\relax \left.\fi#2#1\rvert}
\newcommand{\abs}[1]{\left\lvert#1\right\rvert}

\newcommand{\bn}{\boldsymbol{n}}

\newcommand{\bx}{\boldsymbol{x}}

\newcommand{\bF}{\boldsymbol{F}}

\newcommand{\bM}{\boldsymbol{M}}
\newcommand{\bN}{\boldsymbol{N}}

\newcommand{\bR}{\boldsymbol{R}}

\makeatletter
\newcommand{\Rmnum}[1]{\expandafter\@slowromancap\romannumeral #1@}
\makeatother

\newcommand{\bK}{\boldsymbol{K}}

\newcommand{\bj}{\boldsymbol{j}}

\newcommand*\diff{\mathop{}\!\mathrm{d}}

\allowdisplaybreaks

\biboptions{sort&compress}
\begin{document}

\begin{frontmatter}
\title{Cahn--Hilliard phase-field modeling of tumor growth via locally adaptive isogeometric analysis with THB-splines}
\author[label1]{Dhiraj S. Bombarde\corref{cor1}} %
\ead{dhiraj.bombarde@unipv.it} %
\author[label2]{Carlotta Giannelli} %
\ead{carlotta.giannelli@unifi.it} %
\author[label1]{Alessandro Reali} %
\ead{alessandro.reali@unipv.it} %
\author[label3]{Guillermo Lorenzo} %
\ead{guillermo.lorenzo@udc.es} %
\cortext[cor1]{Corresponding author}
\address[label1]{Department of Civil Engineering and Architecture, University of Pavia, via Ferrata 3, 27100 Pavia, Italy}
\address[label2]{Dipartimento di Matematica e Informatica ``U. Dini'', Università degli Studi di Firenze, Viale Morgagni 67/A, 50134 Firenze, Italy}
\address[label3]{Group of Numerical Methods in Engineering, Department of Mathematics, School of Civil Engineering, and CITEEC, University of A Coruña, Campus de Elviña s/n, 15008 A Coruña, Spain}

\begin{abstract}
Predicting tumor dynamics in biological systems under physiologically relevant conditions using mathematical and computational models remains a challenging problem.
Continuum models based on phase-field (diffuse-interface) formulations have proven to be an effective modeling strategy to govern tumor dynamics and interactions of multiple species. 
Within this framework, the present work investigates a tumor growth model based on the Cahn-Hilliard (CH) equation.
The formulation involves a fourth-order differential operator that imposes higher continuity requirement on approximation spaces for a well-defined primal variational formulation. 
To address this challenge, we use isogeometric analysis (IGA), which inherently satisfies this requirement through spline-based basis functions and eliminates the need for mixed or auxiliary-variable approaches commonly used in standard finite element discretizations. 
Additionally, a locally adaptive IGA scheme with truncated hierarchical B-splines (THB-splines) is used to reduce computational cost while maintaining accuracy. %
The model is first evaluated on standard benchmark cases and then applied to an organ-scale, patient-specific geometric model of the breast reconstructed from magnetic resonance imaging (MRI) data. 
Our results show that the model reproduces known tumor morphologies, ranging from a spheroidal pattern to fingered growth. 
A series of numerical experiments further shows the diversity of tumor dynamics produced by different model parameter choices. 
Our findings demonstrate the predictive potential of the CH-based phase-field tumor growth model integrated with a locally adaptive IGA framework.
\end{abstract}
\end{frontmatter}

\section{Introduction}

Cancer places significant strain on healthcare systems, affecting both resources and patient well-being. 
As of 2022, cancer accounted for nearly 20 million new cases and 9.7 million deaths globally, contributing to about 16.8\% of all deaths and 22.8\% of deaths from noncommunicable diseases.
This burden is projected to reach 35 million new cases annually by 2050, a 77\% rise from 2022 \cite{Bray2024}.
Despite these concerning statistics, the outlook for cancer has improved in recent years. 
Increased awareness, widespread screening programs, and advances in early detection and treatment have led to significant improvements in patient outcomes over the past decades. 
As a result, survival rates for all cancers combined have improved substantially, rising from 49\% in the mid-1970s to 69\% in 2014--2020 in the United States \cite{Siegel2025}  and to approximately 40–80\% (depending on patient's age and sex) in 2000--2007 across Europe \cite{ECIS2025}.
With improved survival rates and advances in cancer care, the clinical focus is increasingly shifting towards personalized disease management, which requires a better understanding and anticipation of the long-term evolution of the disease after diagnosis and treatment.
Such efforts aim to improve therapeutic decision-making while enhancing quality of life without compromising survival outcomes \cite{Yankeelov2015, Karolak2018, Lorenzo2024}.

Mathematical and computational models representing the fundamental mechanisms of cancer provide an interpretable in-silico framework for investigating tumor responses under diverse conditions \cite{Anderson2018, Lorenzo2024, Desoyer2025}. 
Despite their potential, accurately predicting tumor dynamics in biological systems under physiologically relevant conditions remains a challenging problem. 
Mathematical oncology encompasses a wide spectrum of modeling approaches and associated mathematical frameworks, each designed to capture different aspects of tumor biology depending on the research objective and the level of biological detail required. %
These typically include ordinary differential equations (ODEs) for spatially averaged dynamics, partial differential equations (PDEs) for tissue-level spatial processes, discrete or agent-based models for individual cell behavior and heterogeneity, and hybrid frameworks for coupled multiscale phenomena, etc \cite{Bull2022, Desoyer2025}. 
Among these approaches, continuum models constitute a hallmark strategy to describe tissue-scale tumor dynamics, including  spatial interactions between cell types and with respect to tumor-driving substances (e.g., nutrients, growth factors, treatment drugs)  as well as evolving interfaces between tumor and healthy tissue or between intratumoral regions. 
In particular, phase-field (diffuse-interface) formulations provide an attractive continuum modeling framework to capture the spatiotemporal evolution of several interacting species subjected to multiple physical processes into a unified mathematical description within a physiologically relevant setting \cite{Wise2008,Garcke2016,Lima2016,Lorenzo2017,Garcke2021,Ebenbeck2021,Fritz2023}.

Phase-field models avoid the explicit tracking of evolving interfaces, naturally enforce complex boundary conditions across the interface, and allow realistic representation of diffuse interfaces, particularly in cases where phase boundaries are not sharply defined.
Many phase-field models are further grounded in a continuum thermodynamic framework based on balance laws, constitutive relations, and energy principles, enabling the systematic inclusion of multiple interacting fields and physical processes.
Although originally developed in the context of material science to describe phase separation in alloys and binary mixtures \cite{Cahn1958, Allen1979}, these advantages have motivated the use of phase-field models across a wide range of engineering problems, including solidification, solid-state phase transformations, grain growth and coarsening, microstructure evolution in thin films, surface-stress-induced pattern formation, fracture mechanics, crystal growth in the presence of strain, dislocation dynamic, electromigration, and multiphase flow (see~\cite{Chen2002, Tourret2022, Diehl2022, Zhao2023} and references therein). 
Tumor growth modeling is no exception, where phase-field models enables the simultaneous representation of complex tissue interfaces and interacting processes such as cell proliferation, necrosis, nutrient transport, and mechanical effects such as stress-induced mechanical inhibition of tumor growth and viscoelasticity \cite{Cristini2009, Hawkins-Daarud2013, ODEN2010, Lowengrub2010, Wise2008,Garcke2016,Lima2016,Lorenzo2017,Garcke2021,Fritz2023}.

Early phase-field tumor models were formulated as diffuse-interface regularizations of classical continuum free boundary problems \cite{Cristini2003}, introducing an interfacial thickness parameter to smooth the discontinuous interface over a finite width \cite{Frieboes2007,Wise2008, Bures2021}. 
This process is mathematically demanding, as proving convergence of diffuse-interface models to the corresponding free boundary problem in the limit of vanishing interfacial thickness is nontrivial. 
In contrast, several tumor models are derived directly from free-energy functionals within the framework of classical thermomechanics using Coleman–Noll-type arguments. %
In these models, constitutive relations depend on the variational derivatives of the free-energy and are coupled with the fundamental balance laws of mass, linear and angular momentum, and energy, making them thermomechanically and thermodynamically consistent \cite{HawkinsDaarud2012, Gomez2017}. 
A notable subset includes models based on the classical Allen-Cahn (AC) and Cahn-Hilliard (CH) equations, which have received considerable attention in recent years for modeling tumor dynamics \cite{Cristini2009, Hawkins-Daarud2013, ODEN2010, Lowengrub2010, Wise2008,Garcke2016,Lima2016,Lorenzo2017,Garcke2021,Fritz2023}. The present work focuses on the latter.
CH-based models are particularly attractive as they provide a diffuse-interface framework that maintains physically realistic interfaces over long times, prevents unphysical mass changes, and allows consistent coupling with nutrient diffusion, angiogenesis, or extracellular matrix evolution. In the sharp-interface limit, they align with moving boundary formulations, offering a rigorous basis for free-boundary comparisons \cite{Riva2025,Melchionna2018}.

Despite these advantages, numerical treatment to solve CH-based tumor growth models remains a challenge due to the fourth-order spatial differential operators inherent to this phase-field equation, which pose a higher continuity requirement (at least $C^1$) on approximation spaces for a well-defined primal variational formulation.
Traditional numerical approaches for higher-order operators, such as finite difference and spectral methods, are typically restricted to simple geometries. 
However, tumors grow into complex shapes, necessitating more geometrically flexible discretization framework. 
Consequently, the finite element analysis (FEA) has become the most widely used approach for such problems. 
Unfortunately, only a limited number of two-dimensional (2D) elements provide $C^1$ inter-element continuity, with no practical counterparts available in three-dimensions (3D).
As a result, FEA solutions of the CH-based models are relatively uncommon and are most often formulated using mixed methods %
at the cost of extra degrees of freedom (DOFs) in addition to primary variables \cite{Gomez2008}.

In the present work, we focus on a recently proposed CH-based tumor growth model \cite{Ebenbeck2021, Garcke2024} and propose a computational approach to solve it directly in its primal form, intead of using the originally proposed FEA formulation. 
To address the higher-continuity requirement of the underlying numerical framework, isogeometric analysis (IGA) is employed \cite{Hughes2005,Cottrell2009_book}. 
IGA offers a robust computational approach for solving PDEs within a unified framework through integration of computer-aided design (CAD) and FEA principles.
The foundational idea lies in the use of spline basis functions for both the geometric representation of the computational domain and the discretization of the governing PDEs within an isoparametric framework.
This %
eliminates the geometric approximation errors inherent in classical FEA, and allows a more coherent interaction between design and analysis.
Moreover, IGA significantly minimizes the effort needed in pre-processing to create an analysis-ready model, naturally achieves higher inter-element continuity, and reduces dependence on CAD during successive refinements \cite{Cottrell2009_book}.

Beyond the continuity requirement, accurately resolving the steep phase-field gradients confined to narrow interfacial regions requires sufficiently fine meshes, making the computational treatment of CH models a nontrivial task. 
Hence, to maintain computational tractability while balancing accuracy, local refinement with mesh adaptivity is incorporated within IGA.
Of note, these techniques constitute an area of active research; see \cite{Buffa2022} and references therein.
Among the available locally refinable spline spaces, the present work uses the local refinement capabilities of truncated hierarchical B-splines (THB-splines) to attain the accuracy of global refinement while drastically reducing both DOFs and computational efforts within an efficient suitably-graded adaptive meshing strategy \cite{giannelli2012, giannelli2016, bracco2018}.
The use of THB-splines is motivated by their key mathematical properties such as linear independence, non-negativity, variation-diminishing behavior, and partition of unity. 
In addition, their reduced support improves the sparsity of system matrices compared to standard hierarchical B-splines (HB-splines) and, in some cases, enhances matrix conditioning \cite{bracco2018,Buffa2022}.

Finally, this paper assesses the proposed strategy through a series of benchmark problems in 2D and 3D domains, demonstrating its capability to make large-scale simulations of tumor growth with complex morphologies computationally feasible. 
The model efficiently reproduces established morphologies, such as spheroidal growth and finger-like patterns, while also shows how parameter variations and sensitivity to initial conditions influence the resulting tumor shapes.
This is further extended to a patient-specific breast geometry reconstructed from magnetic resonance imaging (MRI) data, demonstrating the potential of the phase-field model and the computational approach in a clinically relevant setting.
To the best of the authors' knowledge, this represents the first attempt to investigate breast cancer growth using a higher-order CH-based model in a realistic, organ-scale 3D setting.

The rest of the paper is organized as follows: In Section~\ref{model_fromulation_numerical_methods}, the CH-based phase-field tumor growth model coupled with nutrient dynamics is presented, which includes the governing equations and their corresponding weak formulation.
In Section~\ref{adaptive_IGA_with_THB_splines}, the basic concepts of adaptive IGA with THB-splines are briefly recalled, including adaptive refinement and coarsening strategies that ensure suitably graded hierarchical meshes.
This is followed by the spatial discretization leading to the semi-discrete residual system and the description of a stable time-integration scheme for the nonlinear system.
In Section~\ref{Numerical_studies}, the model is assessed using 2D and 3D benchmark problems to study tumor morphologies and parameter effects, then it is extended to a patient-specific breast geometry from MRI to evaluate performance in a realistic anatomical setting. Section~\ref{Conclusion} provides a summary of the results, discusses the potential of the model, and outlines potential directions for future research. %
Finally, \ref{Appendix-flowchart} presents a pseudocode of the implementation framework, while \ref{Complementary_Numerical_Studies} provides supplementary numerical results.

\section{Model formulation} \label{model_fromulation_numerical_methods} %
This section outlines the CH-based phase-field framework for spatiotemporal tumor growth dynamics through proliferation (cell growth and division), apoptosis (programmed cell death), chemotaxis (tumor cell migration toward nutrient-rich regions), and coupled nutrient dynamics.
The governing equations are presented in strong form with relevant boundary and generic initial conditions, followed by their weak formulation, which provides the foundation for the numerical discretization in subsequent sections.

\subsection{CH-based tumor growth model}

Two field variables characterize the model: the tumor phase-field $\phi$ that follows a CH-type equation and the nutrient concentration $\sigma$ governed by a reaction-diffusion equation.
The model considers a two-phase mixture of healthy and tumor species, having mass densities $\rho_i$ and $\tilde{\rho}_i$ defined with respect to the total and pure component volume, respectively. 
The mass density of the material body is given as $\rho = \sum_{i=1}^2 \rho_i$, whereas mass and volume fractions of the components are defined as $c_i = \rho_i/\rho$ and $\phi_i = \rho_i/\tilde{\rho}_i$, respectively.
Given that $\sum_{i=1}^{2} \phi_i = 1 $ and $ \sum_{i=1}^{2} c_i = 1$, the order parameter (phase-field, $\phi$) is defined as $\phi = \phi_2 - \phi_1$ with $\phi \in [-1,1]$. 
By construction, $\phi = -1$ and $\phi= 1$ correspond to healthy and tumor tissue, respectively. 
Assuming constant densities and no excess volume due to mixing, the mass balance equations can be expressed solely in terms of the $\phi$ as
\begin{equation*}
    \parder{\phi}{t} = \nabla \cdot \left[ M(\phi) \nabla \mu \right] + S(\phi,\sigma),   
\end{equation*}
where $M(\phi) = M \left[\frac{1}{2} (1+\phi)^2\right]$ is the non-negative phase-dependent (degenerate) cell mobility, $M$ is the mobility parameter, $\mu$ is the chemical potential, and $S(\phi,\sigma)$ is the source term accounting for mass transfer due to nutrient-dependent proliferation and apoptosis \cite{Ebenbeck2021}.
A general free-energy density of the form
\begin{equation*}
    e(\phi, \nabla \phi, \sigma) = \mathcal{F}(\phi, \nabla \phi) + \mathcal{N}(\phi,\sigma),
\end{equation*}
is considered, where $\mathcal{F}(\phi, \nabla \phi): H^1(\Omega) \rightarrow{\mathbb{R}},\; \Omega \subset \mathbb{R}^d $ ($\Omega$ = open domain, $d$ = spatial dimension) is the Ginzburg--Landau free energy and $\mathcal{N}(\phi,\sigma)$ is the energy contribution from nutrient-dependent interactions with the tumor.
These terms can be defined as 
\begin{align}
    \mathcal{F}(\phi, \nabla \phi) &= Ef(\phi) + \frac{\lambda}{2} \abs{\nabla\phi}^2, \label{Ginzburg-Landau}\\
    \mathcal{N}(\phi,\sigma) &= \frac{\mathcal{X}_\sigma}{2}\abs{\sigma}^2 + \mathcal{X}_\phi \sigma (1-\phi) \label{nutrient_energy},
\end{align}
where $f(\phi) = \frac{1}{4} (1-\phi^2)^2$ is a double-well potential function. 
The constants are scaled as follows: $E = \frac{\beta_1}{\epsilon} $ and $\lambda = \beta_2\epsilon$, where $\epsilon>0$ controls the thickness of the diffuse interface separating the healthy and tumor phases and $\beta_1$, $\beta_2$ are positive non-dimensional constants.
In Eq.~\eqref{nutrient_energy}, the first term captures nutrient-induced energy effects, while the second term models chemotaxis. 
The parameters $\mathcal{X}_\sigma$ and $\mathcal{X}_\phi$ represent diffusion and chemotaxis parameter, respectively.

The chemical potential %
is obtained as the variational derivative of the free-energy functional $\int_\Omega e(\phi, \nabla \phi, \sigma) \diff \Omega$ with respect to $\phi$, resulting in  %
\begin{equation*}
     \mu =   \parder{e}{\phi} - \nabla \cdot \left( \parder{e}{\nabla\phi}\right)  = E f'(\phi)   - \lambda \Delta \phi - \mathcal{X}_\phi \sigma.
\end{equation*}
Finally, the source term $S(\phi,\sigma)$ is assumed to follow linear kinetics and is defined as 
\begin{equation*}
 S(\phi,\sigma) = \kappa (\mathcal{P} \sigma - \mathcal{A})h(\phi),    
\end{equation*}
where $h(\phi) = \frac{1}{2} (\phi + 1)$ varies linearly from 0 to 1 as $\phi$ transitions from -1 to 1, and $\kappa$ is the scale parameter for the source. 
Here, $\mathcal{P}$ and $\mathcal{A}$ are non-negative constants associated with proliferation and apoptosis, respectively. %

The nutrient dynamics is governed by the following balance law: %
\begin{equation*}
    \parder{\sigma}{t} =  -\nabla \cdot \bj^* -S^*(\phi,\sigma),    %
\end{equation*}
where $\bj^*$ is the diffusive flux representing nutrient transport, and $S^*(\phi,\sigma)$ is a sink term accounting for nutrient consumption.
Assuming a constant nutrient diffusion coefficient $\bar{D}$, $\bj^*$ is then evaluated as
\begin{equation*}
    \bj^* = -\bar{D} \nabla  \left(\parder{e}{\sigma}\right)  = -D \left( \nabla \sigma - \mathcal{X} \nabla \phi \right), 
\end{equation*}
where $D = \bar{D}\mathcal{X}_\sigma$ and $\mathcal{X} = \frac{\mathcal{X}_\phi}{\mathcal{X}_\sigma}$ are the effective diffusion and active transport coefficients, allowing later to be suppressed ($\mathcal{X}\xrightarrow{} 0$) while retaining chemotaxis.
The sink $S^*(\phi,\sigma)$ is assumed to follow linear kinetics, with $S^*(\phi,\sigma) = \mathcal{C} \sigma h(\phi)$, where $\mathcal{C}$ is the nutrient consumption rate.
Nutrient dynamics are further assumed quasi-static $\left(\parder{\sigma}{t} =0\right)$, which is justified by the fact that diffusion occurs on a much shorter time scale compared to cell proliferation.
Lastly, homogeneous Neumann boundary conditions are prescribed for $\phi$ and $\mu$, ensuring that the energy equality is satisfied.
For $\sigma$, the value is fixed on the boundary by a Dirichlet condition thereby representing nutrient supply from well-perfused neighboring healthy tissue. %
The resulting initial/boundary value problem is then summarized as follows: Let $\Omega \subset \mathbb{R}^d$ be an open domain of spatial dimension $d$ with a smooth boundary $\Gamma$ and closure $\bar{\Omega} = \Omega \cup \Gamma$.
Given the initial condition $\phi(\bx,0) = \phi_0(\bx) \in \mathbb{R},$ $\bx \in \bar{\Omega}$ and time $t \in \left[ 0, T \right]$, find $\phi(\bx,t)$ and $\sigma(\bx,t)$ $\in \mathbb{R} $ $: (\bx,t) \in \bar{\Omega} \times \left[0, T \right]$ such that  
\begin{align}
    \parder{\phi}{t} &= \nabla\cdot\left[ M(\phi) \nabla \left( E f'(\phi) - \lambda \Delta\phi - \mathcal{X}_\phi \sigma\right)\right] + \frac{\kappa}{2}(\mathcal{P} \sigma - \mathcal{A})(\phi + 1) && \text{in } \Omega \times \left[ 0, T \right], \label{CH_Eq}\\
    0 &= D \nabla\cdot\left( \nabla\sigma-\mathcal{X} \nabla\phi\right) - \frac{\mathcal{C}}{2} \sigma(\phi + 1) && \text{in } \Omega,    \label{Nutri_Eq}\\
    \nabla\phi\cdot\bn &= \nabla \left( E f'(\phi) - \lambda \Delta\phi - \mathcal{X}_\phi \sigma\right)\cdot\bn   = 0 && \text{on } \Gamma \times \left[ 0, T \right],  \label{FluxBC} \\
    \sigma &= \sigma_B && \text{on } \Gamma \times \left[ 0, T \right], \label{Dir_Nutri}\\
    \phi(\bx,0) &= \phi_0(\bx) && \text{in } \bar{\Omega} \label{IC_Tumor_Model},
\end{align}
where  $M(\phi) = M \left[\frac{1}{2} (1+\phi)^2\right]$ and $f'(\phi) = \left(  \phi^3 - \phi\right)$. 
The volume-averaged mixture velocity, as considered in \cite{Ebenbeck2021}, is excluded from the present formulation to maintain simplicity, though its inclusion is a natural extension of the present framework.

\subsection{Variational formulation of the coupled system}
We obtain the weak form of the coupled system (Eqs.~\eqref{CH_Eq}--\eqref{IC_Tumor_Model}) using a Bubnov--Galerkin method, where the discrete trial and test spaces coincide for each field variable.
Let $\mathbb{V} := H^2(\Omega)$ be the Hilbert space of functions with square-integrable first and second derivatives, 
and $\mathbb{V}_\phi \subset \mathbb{V}$ be the admissible space for $\phi$ and its variation $\delta\phi$ such that 
\begin{equation*}
    \phi(\cdot,t), \delta\phi(\cdot, t) \in \mathbb{V}_\phi \subset \mathbb{V}, \quad \forall t \in \left[0, T\right].   
\end{equation*}
Following a staggered approach, the variational form of Eq.~\eqref{CH_Eq} can be written as: For a given nutrient field $\sigma$ at $t$, find $\phi \in \mathbb{V}_\phi$ such that $\forall \: \delta \phi \in \mathbb{V}_\phi$
\begin{equation}\label{CH_Weak_Final}
\begin{split} 
\int_{\Omega} & \delta\phi  \parder{\phi}{t} \diff \Omega  + 
\int_{\Omega} \nabla\delta\phi \cdot \left[M(\phi) E f''(\phi)\nabla\phi \right] \diff \Omega -
\int_{\Gamma} \left(\nabla\delta\phi \cdot \bn\right) \lambda  M(\phi) \Delta\phi \diff \Gamma +
\int_{\Omega} \left(\Delta \delta\phi\right) \lambda M(\phi) \Delta\phi \diff \Omega 
\\ +& \int_{\Omega} \nabla \delta\phi \cdot \left[\lambda M'(\phi) \nabla\phi\Delta\phi\right] \diff \Omega - 
\int_{\Omega} \nabla \delta\phi \cdot \left[M(\phi) \mathcal{X}_\phi \nabla\sigma  \right] \diff \Omega -
\int_{\Omega}  \delta\phi \left[ \frac{\kappa}{2}(\mathcal{P} \sigma - \mathcal{A})(\phi + 1) \right] \diff \Omega = 0. %
\end{split}
\end{equation}
The essential boundary condition $\nabla \phi \cdot \bn = 0$ (Eq.~\eqref{FluxBC}) is imposed weakly using Nitsche’s method \cite{Zhao2017, Bracco2023} by adding the following boundary terms to the left hand side of Eq.~\eqref{CH_Weak_Final}:
\begin{equation}\label{Nitsche_Terms}
    - \int_{\Gamma}   \left(\Delta\delta\phi\right) M(\phi) \lambda \left( \nabla \phi \cdot \bn \right) \diff \Gamma 
    + \int_{\Gamma} \left(\nabla\delta\phi \cdot \bn\right) \frac{\varepsilon_n}{h_e} \left( \nabla\phi \cdot \bn\right) \diff \Gamma,
\end{equation}
where $\varepsilon_n > 0$ is a Nitsche penalty parameter, and $h_e$ is a characteristic length corresponding to the element size.

Next, the quasi-static second-order nutrient PDE (Eq.~\eqref{Nutri_Eq}) is formulated in a variational setting by considering the Hilbert space $\bar{\mathbb{V}} := H^1(\Omega)$. %
If $\mathbb{V}_\sigma \subset \bar{\mathbb{V}}$ denotes the admissible space for $\sigma$ and its variation $\delta\sigma$ such that $\sigma, \delta\sigma \in \mathbb{V}_\sigma \subset \bar{\mathbb{V}}$, then the  variational statement for Eq.~\eqref{Nutri_Eq} can be written as: For a given phase-field $\phi$, find $\sigma \in \mathbb{V}_\sigma$ such that $\forall \: \delta \sigma \in \mathbb{V}_\sigma$
\begin{equation}\label{Nutri_weak_form}
    -  \int_{\Omega} D \nabla\delta\sigma \cdot \nabla \sigma \diff \Omega 
    +  \int_{\Omega} D \mathcal{X} \nabla\delta\sigma \cdot \nabla \phi \diff \Omega
    -  \int_{\Omega} \delta\sigma  \left[ \frac{\mathcal{C}}{2} \sigma(\phi + 1) \right] \diff \Omega = 0. %
\end{equation}
Eqs.~\eqref{CH_Weak_Final}--\eqref{Nutri_weak_form} are then discretized in space using a Galerkin isogeometric formulation \cite{Cottrell2009_book} and in time using the generalized-$\alpha$ integration scheme \cite{jansen2000generalized,chung1993}.

\section{Adaptive IGA with THB-splines} \label{adaptive_IGA_with_THB_splines}
This section briefly recalls the basic concepts of adaptive IGA with THB-splines, including adaptive refinement and coarsening strategies that ensure suitably graded hierarchical meshes.

\subsection{B-splines}
A $p$\textsuperscript{th} degree univariate B-spline basis function $\hat{\mathtt{B}}_{i,p}(\xi)$, defined over a knot vector $\boldsymbol\Xi := \{\xi_j\}_{j = 0}^{n+p+1}$ with $\xi_j \in \mathbb{R}$ and $\xi_j \leq \xi_{j+1}$ for $j = 0,\: \ldots,\: n+p$, is defined recursively by the Cox--de Boor formula \cite{Piegl1997}
\begin{align*} \label{DE_boor_formula}
\text{for } p = 0,\quad \hat{\mathtt{B}}_{i,0}(\xi) &= \begin{cases}
    1 & \text{if } \xi_i \leq \xi < \xi_{i+1}\\
    0 & \text{Otherwise }
  \end{cases}, \\
\text{for } p \geq 1,\quad \hat{\mathtt{B}}_{i,p}(\xi) &= \frac{\xi - \xi_i}{\xi_{i+p} - \xi_i}\hat{\mathtt{B}}_{i,p-1}(\xi) + \frac{\xi_{i+p+1} - \xi}{\xi_{i+p+1} - \xi_{i+1}}\hat{\mathtt{B}}_{i+1,p-1}(\xi),
\end{align*}
where $\xi_i \in \boldsymbol\Xi$ and $n$ is the total number of basis functions.
A multivariate B-spline basis is obtained by tensor product of univariate functions. For degree $\mathbf{p}=(p_1,\ldots,p_d)$ defined over knot vectors $\boldsymbol{\Xi}^k$ with $k = 1,\ldots,d$, the multivariate B-spline basis function can be defined as
\begin{equation*}
\hat{\mathtt{B}}_{\textbf{i},\textbf{p}}(\boldsymbol{\xi}) = \prod_{k=1}^d \hat{\mathtt{B}}_{i_k,p_k}(\xi_k),
\end{equation*}
where $\textbf{i}=(i_1,\ldots,i_d)$ is a multi-index and $\boldsymbol{\xi}=(\xi_1,\ldots,\xi_d)$ denotes the vector of parametric coordinates.

\subsection{THB-splines}

\begin{figure}[h!]
\centering
\begin{subfigure}[t]{0.25\linewidth}
\centering
\begin{tikzpicture}[scale=.75]

\draw[fill = jetblue] (0,0) rectangle (1,5);
\draw[fill = jetblue] (4,0) rectangle (5,5);
\draw[fill = jetblue] (1,0) rectangle (4,1);
\draw[fill = jetblue] (1,4) rectangle (4,5);

\draw [pattern=north east lines,
    pattern color=darkgray] (0,0) rectangle (5,5);
\draw (0,0) grid (5,5);

\end{tikzpicture}
\subcaption{$\hat{\Omega}^0$, $\hat{G}^0$, and $\hat{\cal G}^0$
}
\end{subfigure}\hspace*{-.15cm}
\begin{subfigure}[t]{0.25\linewidth}
\centering
\begin{tikzpicture}[scale=.75]

\draw[fill = jetorange] (1,1) rectangle (4,2);
\draw[fill = jetorange] (0.99,2.99) rectangle (4,4);
\draw[fill = jetorange] (0.99,1.99) rectangle (2,3);
\draw[fill = jetorange] (2.99,1.99) rectangle (4,3);

\draw [pattern=north east lines,
    pattern color=darkgray] (1,1) rectangle (4,4);
\draw (0,0) grid (5,5);
\foreach \a in {1,3,5,7,9}
    \draw (0,\a/2) -- (5,\a/2); %
\foreach \a in {1,3,5,7,9}
    \draw (\a/2,0) -- (\a/2,5); %

\end{tikzpicture}
\subcaption{$\hat{\Omega}^1$, $\hat{G}^1$, and $\hat{\cal G}^1$ 
}
\end{subfigure}\hspace*{-.15cm}
\begin{subfigure}[t]{0.25\linewidth}
\centering
\begin{tikzpicture}[scale=.75]

\draw[fill = jetdarkred] (1.99,1.99) rectangle (3,3);

\draw [pattern=north east lines,
    pattern color=darkgray] (2,2) rectangle (3,3);
\draw (0,0) grid (5,5);
\foreach \a in {1,3,5,7,9}
    \draw (0,\a/2) -- (5,\a/2); %
\foreach \a in {1,3,5,7,9}
    \draw (\a/2,0) -- (\a/2,5); %
\foreach \a in {1,3,5,7,9,11,13,15,17,19}
    \draw (0,\a/4) -- (5,\a/4); %
\foreach \a in {1,3,5,7,9,11,13,15,17,19}
    \draw (\a/4,0) -- (\a/4,5); %

\end{tikzpicture}
\subcaption{$\hat{\Omega}^2$, $\hat{G}^2$, and $\hat{\cal G}^2$
}
\end{subfigure}\hspace*{-.15cm}
\begin{subfigure}[t]{.25\linewidth}
\centering
\begin{tikzpicture}[scale=.75]

\draw[fill = jetblue] (0,0) rectangle (1,5);
\draw[fill = jetblue] (4,0) rectangle (5,5);
\draw[fill = jetblue] (1,0) rectangle (4,1);
\draw[fill = jetblue] (1,4) rectangle (4,5);

\draw[fill = jetorange] (0.99,0.99) rectangle (4,2);
\draw[fill = jetorange] (0.99,2.99) rectangle (4,4);
\draw[fill = jetorange] (0.99,1.99) rectangle (2,3);
\draw[fill = jetorange] (2.99,1.99) rectangle (4,3);

\draw[fill = jetdarkred] (1.99,1.99) rectangle (3,3);

\draw (0,0) grid (5,5);
\foreach \a in {3,5,7}
    \draw (1,\a/2) -- (4,\a/2); %
\foreach \a in {3,5,7}
    \draw (\a/2,1) -- (\a/2,4); %
\foreach \a in {9,11}
    \draw (2,\a/4) -- (3,\a/4); %
\foreach \a in {9,11}
    \draw (\a/4,2) -- (\a/4,3); %

\end{tikzpicture}
\subcaption{$\hat{\cal Q}$
}
\end{subfigure}
\caption{Example of domains (hatched regions), grids, and sets of active elements ($\hat{\cal G}^0$- blue, $\hat{\cal G}^1$- yellow, and $\hat{\cal G}^2$- dark red) at levels 0 (a), 1 (b), and 2 (c) for $d=2$. The resulting hierarchical mesh is shown in (d), with the elements at levels 0, 1, and 2 shown in blue, yellow, and dark red, respectively.}
\label{fig:domains}
\end{figure}
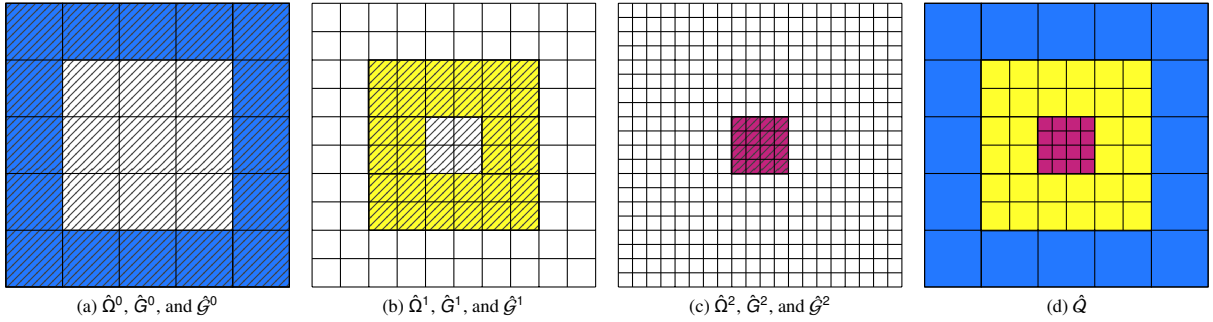

Consider a nested sequence of parametric domains $\hat{\Omega}^0\supseteq\ldots\supseteq\hat{\Omega}^{N-1}$, each contained in a closed hyper-rectangle $\mathtt{R} \in \mathbb{R}^d$, together with a nested sequence of $\mathbf{p}\textsuperscript{th}$ degree tensor-product spline spaces $\hat V^0 \subset \hat V^1 \subset \ldots \subset \hat V^{N-1}$. 
For each level $ \ell = 0,\ldots,N-1$, let $\hat{\mathcal{B}}^\ell$ and $\hat{G}^\ell$ denote the collection of associated B-spline basis of degree $\mathbf{p}$ and rectilinear grid, respectively. 
Open knot vectors are assumed in every parametric direction at level 0, with single knot multiplicities.
To ensure the nestedness of the spline spaces, dyadic mesh refinement is adopted between consecutive hierarchical levels\footnote{Note that the hierarchical spline framework can also accommodate more general refinement strategies \cite{giannelli2014}.}.
Furthermore, any newly inserted knot is assigned multiplicity one. Starting from $\hat{\Omega}^0=\Omega$ and assuming that each domain $\hat{\Omega}^\ell$, for $\ell=1,\ldots,N-1$ is defined as a union of elements from the previous level $\ell-1$, the hierarchical mesh $\hat{\mathcal{Q}}$ is defined as the collection of active elements across all levels, namely
\[
\hat{\mathcal{Q}} :=
\left\{ \hat{Q} \in \hat{\mathcal{G}}^\ell \;:\; \ell = 0,\ldots,N-1 \right\},
\]
where the set of active elements at level $\ell$ is given by
$
\hat{\mathcal{G}}^\ell :=
\left\{ 
\hat{Q} \in \hat{G}^\ell \;:\;
\hat{Q} \subset \hat{\Omega}^\ell 
\ \text{and} \ 
\hat{Q} \not\subset \hat{\Omega}^{\ell+1}
\right\}$. An example of hierarchical configuration on three refinement
levels in a bivariate case is shown in Fig.~\ref{fig:domains}. The plots depict the domains $\{\hat{\Omega}^\ell\}_{\ell=0}^2$, the rectilinear grids $\{\hat{{G}}^\ell\}_{\ell=0}^2$, active elements $\{\hat{{\mathcal{G}}}^\ell\}_{\ell=0}^2$, and the hierarchical mesh $\hat{\mathcal{Q}} = \{\hat{Q}\in {\cal G}^\ell: \ell=0,1,2\}$.

Let $\text{supp }\hat{\mathtt{B}}$ denote the intersection of the support of the B-spline basis $\hat{\mathtt{B}}$ with the level 0 domain $\hat{\Omega}^0$. Then, the HB-spline basis $\hat{\mathcal{H}}$ on the hierarchical mesh $\mathcal{\hat{Q}}$ is defined as:
\begin{equation*}
\hat{\mathcal{H}}(\hat{\mathcal{Q}}):=\left\lbrace \hat{\mathtt{B}}\in\hat{\mathcal{B}}^\ell: \text{supp }\hat{\mathtt{B}}\subseteq \hat{\Omega}^\ell \wedge \text{supp }\hat{\mathtt{B}} \not \subseteq \hat{\Omega}^{\ell+1}, \: \ell=0,\ldots,N-1, \: \hat{\Omega}^{N} = \emptyset \right\rbrace.
\label{eq:hierarchicalBasis}
\end{equation*}
Let $\hat{s} \in \hat V^\ell \subset \hat V^{\ell+1}$ and write $\hat{s}$ in terms of the refined B-spline basis at level $\ell+1$ as
\[
\hat{s}
=
\sum_{\hat{\mathtt{B}} \in \hat{\mathcal{B}}^{\ell+1}}
c_{\hat{\mathtt{B}}}^{\ell+1}(\hat{s})\, \hat{\mathtt{B}},
\]
where $c_{\hat{\mathtt{B}}}^{\ell+1}(\hat{s})$ denotes the coefficient of $\hat{s}$ associated with the basis function $\hat{\mathtt{B}}$. The truncation of $\hat{s}$ with respect to $\hat{\mathcal{B}}^{\ell+1}$ is defined as
\[
\text{trunc}^{\ell+1} \hat{s}
:=
\sum_{\substack{\hat{\mathtt{B}} \in \hat{\mathcal{B}}^{\ell+1} \\
\text{supp }\hat{\mathtt{B}} \not\subseteq \hat{\Omega}^{\ell+1}}}
c_{\hat{\mathtt{B}}}^{\ell+1}(\hat{s}) \, \hat{\mathtt{B}}.
\]
The THB-spline basis $\hat{\mathcal{T}}$ on the mesh $\hat{\mathcal{Q}}$ is then defined as
\begin{equation*}
\hat{\mathcal{T}}(\hat{\mathcal{Q}}):=
\left\lbrace
\textnormal{Trunc}^{\ell+1}\hat{\mathtt{B}}:\hat{\mathtt{B}}\in \hat{\mathcal{B}}^\ell\cap\hat{\mathcal{H}}(\hat{\mathcal{Q}}), \ell=0,\ldots,N-2
\right\rbrace
\bigcup
\left\{
\hat{\mathtt{B}}\in \hat{\mathcal{B}}^{N-1}\cap \hat{\mathcal{H}}(\hat{\mathcal{Q}})
\right\},
\end{equation*}
where $\textnormal{Trunc}^{\ell+1}\hat{\mathtt{B}}:=\textnormal{trunc}^{N-1}
(\textnormal{trunc}^{N-2}
(\ldots(\textnormal{trunc}^{\ell+1}(\hat{\mathtt{B}} ) )\ldots ) )$,
for any $\hat{\mathtt{B}}\in\hat{\mathcal{B}}^\ell \cap \hat{\mathcal{H}}(\hat{\mathcal{Q}})$.
Note that the two bases introduced above span the same space, $\text{span}\,\hat{\mathcal{H}} = \text{span}\,\hat{\mathcal{T}}$. Further details are provided in \citep{giannelli2012}, where the properties of non-negativity, linear independence, and partition of unity of the THB-spline basis were also established.
The successive truncation $\textnormal{Trunc}^{\ell+1}\hat{\mathtt{B}}$ of any HB-spline $\hat{\mathtt{B}}$ leads to a THB-spline whose support is either equal or reduced when compared to the one of the mother function $\hat{\mathtt{B}}$. 
The reduced support of THB-splines leads to improved sparsity patterns in the corresponding system matrices compared to HB-splines, and, in some cases, to enhanced conditioning~\cite{giannelli2012, giannelli2016, bracco2018}. The truncation mechanism underlying THB-splines also influences the design of refinement algorithms, enabling the construction of structured hierarchical meshes with a more localized graded refinement and, consequently, a reduced number of DOFs compared to HB-splines \cite{bracco2018,Buffa2022}.

\subsection{Error estimator}\label{Error_estimator_sec}

Before detailing the admissible THB-spline refinement and coarsening scheme, we first introduce the error estimator used to identify the regions where the mesh is refined or coarsened. %
For the present work, we use a gradient-based error estimator as follows:
\begin{equation}
    \varepsilon_{\phi} = \frac{ 1 }{\text{Vol}(\Omega_e)} \int_{\Omega_e} \| \nabla\phi \|^2 \diff \Omega, \label{error_estimator}
\end{equation}
where $\| \nabla\phi \|^2$ is the squared $L^2$-norm of the gradient of the phase field over each element $\Omega_e$, and $\text{Vol}(\Omega_e)$ is the corresponding elemental volume.
We introduce two tolerances, $\alpha$ and $\beta$, that define the marking strategy. 
Element marking is performed using the greatest appearing error utilization (GARU) strategy. In brief, let $\varepsilon^i_{\phi}$, $i = 1, \dots, n_{\mathrm{el}}$, denote the local error on element $i$. 
Then, an element is marked for refinement if $\varepsilon^i_{\phi} > \alpha \, \varepsilon_{\max}$ and for coarsening if $\varepsilon^i_{\phi} < \beta \, \varepsilon_{\max}$, where $\varepsilon_{\max} = \max \{\varepsilon^i_{\phi}\}$ and the thresholds $\alpha, \beta \in (0,1)$ with $\beta < \alpha$.

\subsection{Admissible THB-spline refinement and coarsening}\label{sec:admissible_meshes}
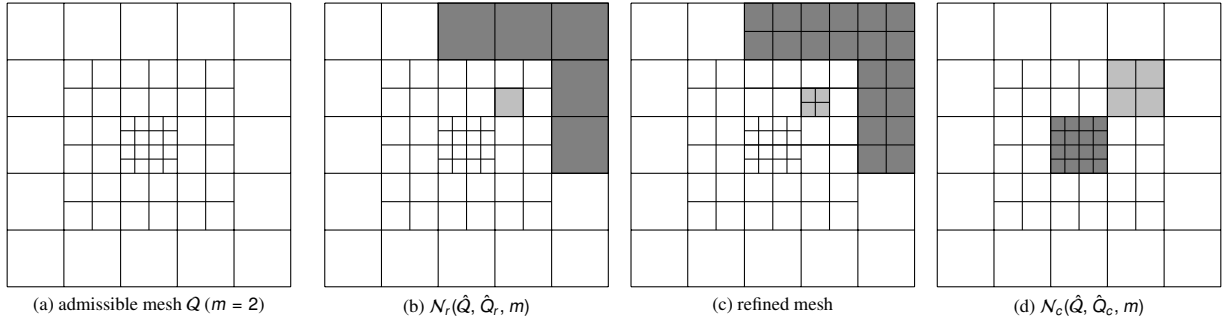
\begin{figure}[t!]
\centering
\begin{subfigure}[t]{.25\linewidth} 
\centering
\begin{tikzpicture}[scale=.75]
\draw (0,0) grid (5,5);
\foreach \a in {3,5,7}
    \draw (1,\a/2) -- (4,\a/2); %
\foreach \a in {3,5,7}
    \draw (\a/2,1) -- (\a/2,4); %
\foreach \a in {9,11}
    \draw (2,\a/4) -- (3,\a/4); %
\foreach \a in {9,11}
    \draw (\a/4,2) -- (\a/4,3); %
\end{tikzpicture}
\subcaption{admissible mesh ${\cal Q}$ ($m=2$)} \label{fig:neigh-a}
\end{subfigure}
\begin{subfigure}[t]{.25\linewidth} 
\centering
\begin{tikzpicture}[scale=.75]
\draw [fill=lightgray] (3,3) rectangle (3.5,3.5);
\draw [fill=gray] (2,4) rectangle (5,5);
\draw [fill=gray] (4,2) rectangle (5,4);
\draw (0,0) grid (5,5);
\foreach \a in {3,5,7}
    \draw (1,\a/2) -- (4,\a/2); %
\foreach \a in {3,5,7}
    \draw (\a/2,1) -- (\a/2,4); %
\foreach \a in {9,11}
    \draw (2,\a/4) -- (3,\a/4); %
\foreach \a in {9,11}
    \draw (\a/4,2) -- (\a/4,3); %
\end{tikzpicture}
\subcaption{ %
${\cal N}_r(\hat{\cal Q},\hat{Q}_r,m)$
} \label{fig:neigh-b}
\end{subfigure}\hspace{-0.15cm}
\begin{subfigure}[t]{.25\linewidth} 
\centering
\begin{tikzpicture}[scale=.75]
\draw [fill=lightgray] (3,3) rectangle (3.5,3.5);
\draw [fill=gray] (2,4) rectangle (5,5);
\draw [fill=gray] (4,2) rectangle (5,4);
\draw (0,0) grid (5,5);
\foreach \a in {3,5,7}
    \draw (1,\a/2) -- (4,\a/2); %
\foreach \a in {3,5,7}
    \draw (\a/2,1) -- (\a/2,4); %
\foreach \a in {5,7,9}
    \draw (2,\a/2) -- (5,\a/2); %
\foreach \a in {5,7}
    \draw (\a/2,4) -- (\a/2,5); %
\foreach \a in {9}
    \draw (\a/2,2) -- (\a/2,5); %
\foreach \a in {9,11}
    \draw (2,\a/4) -- (3,\a/4); %
\foreach \a in {9,11}
    \draw (\a/4,2) -- (\a/4,3); %
\foreach \a in {13}
    \draw (3,\a/4) -- (3.5,\a/4); %
\foreach \a in {13}
    \draw (\a/4,3) -- (\a/4,3.5); %
\end{tikzpicture}
\subcaption{refined mesh} \label{fig:neigh-c}
\end{subfigure}\hspace{-0.15cm}
\begin{subfigure}[t]{.25\linewidth} 
\centering
\begin{tikzpicture}[scale=.75]
\draw [fill=lightgray] (3,3) rectangle (4,4);
\draw [fill=gray] (2,2) rectangle (3,3);
\draw (0,0) grid (5,5);
\foreach \a in {3,5,7}
    \draw (1,\a/2) -- (4,\a/2); %
\foreach \a in {3,5,7}
    \draw (\a/2,1) -- (\a/2,4); %
\foreach \a in {9,11}
    \draw (2,\a/4) -- (3,\a/4); %
\foreach \a in {9,11}
    \draw (\a/4,2) -- (\a/4,3); %
\end{tikzpicture}
\subcaption{${\cal N}_c(\hat{\cal Q},\hat{Q}_c,m)$
} \label{fig:neigh-d}
\end{subfigure}
\caption{Example of refinement and coarsening neighborhoods for the hierarchical mesh shown in (a), with $\mathbf{p}=(2,2)$ and $m=2$. For the marked element $\hat{Q}_r$ (light gray) in (b), the refinement neighborhood ${\cal N}_r(\hat{\cal Q},\hat{Q}_r,m)$ is highlighted in gray, and the resulting refined mesh is displayed in (c); for the marked element $\hat{Q}_c$ (light gray) in (d), the coarsening neighborhood ${\cal N}_c(\hat{\cal Q},\hat{Q}_c,m)$ is highlighted in gray. All internal knots have multiplicity one.}
\label{fig:neigh}
\end{figure}

Given a hierarchical mesh $\hat{\mathcal{Q}}$ and a set of active elements marked for refinement or coarsening according to the criterion defined in Section~\ref{Error_estimator_sec}, the finer or coarser elements to be activated or deactivated are determined based on the prescribed refinement and coarsening rules.
An essential aspect for adaptive isogeometric methods is the construction of mesh configurations for which the number of basis functions that are nonzero on any given element remains uniformly bounded. 
A mesh $\hat{\mathcal{Q}}$ is said to be admissible of class $m$ if, for every element $\hat{Q} \in \hat{\mathcal{Q}}$, the basis functions in $\hat{\mathcal{T}}(\hat{\mathcal{Q}})$ that are nonzero on $\hat{Q}$ originate from at most $m$ consecutive refinement levels.
Refinement and coarsening of admissible meshes produce hierarchical meshes in which the number of THB-spline basis overlapping any given element is independent of the total number of refinement levels. An example of an admissible hierarchical mesh of class $m=2$, degree $\mathbf{p}=(2,2)$, and levels $\ell=0,1,2$ is shown in Fig.~\ref{fig:neigh-a}. 
Note that level-0 THB-spline basis functions overlap all elements in the outer ring near the boundary but vanish on the inner elements belonging to level 2.
If HB-splines were used on the same mesh instead, the active level-0 B-splines would take nonzero values on the inner elements at level 2. 
Although admissibility could still be enforced, this would be at the cost of introducing additional elements and a different mesh configuration. In this case, an intermediate ring consisting of 4 (rather than 2) level-1 elements would be required to obtain an admissible mesh of class $2$ for the HB-spline basis. 
For a comprehensive comparison of admissible mesh configurations for HB- and THB-splines, see \cite{bracco2018,Buffa2022, giannelli2025}, which highlights the advantage of THB-spline basis functions in generating less refined meshes and, consequently, a lower number of DOFs.

For THB-splines, refinement and coarsening algorithms that preserve the admissible nature of hierarchical meshes were developed in \cite{buffa2016,buffa2017b} and \cite{carraturo2019}, respectively. 
For any element $\hat{Q}_r\in \hat{\cal G}^\ell$ or $\hat{Q}_c\in \hat{\cal G}^\ell$ marked for refinement or coarsening, these algorithms rely on the concept of \emph{refinement and coarsening neighborhoods} on the hierarchical mesh $\hat{\mathcal{Q}}$ with respect to the admissibility class $m$, denoted as ${\cal N}_r(\hat{\mathcal{Q}},\hat{Q}_r,m)$ and ${\cal N}_c(\hat{\mathcal{Q}},\hat{Q}_c,m)$, and defined as follows
\begin{align*}
{\cal N}_r(\hat{\cal Q},\hat{Q}_r,m) &:=
\left\{
\hat{Q}'\in\hat{\cal G}^{\ell-m+1}: \exists\, \hat{Q}''\in \omega(\hat{Q},\ell-m+2), \hat{Q}'' \subset\hat{Q}'
\right\},\\
{\cal N}_c(\hat{\cal Q},\hat{Q}_c, m) &:= \left\{\hat{Q}'\in \hat{\cal G}^{\ell+m}: \exists\, \hat{Q}''\in\hat{\cal G}^{\ell+1} \text{ and } \hat{Q}''\subset \hat{Q}, \text{with } \hat{Q}'\in \omega(\hat{Q}'',\ell+1)\right\}
\end{align*}
where
\[
\omega(\hat{Q},k) :=
\left\{
\hat{Q}'\in \hat{G}^k : \exists\, \hat{\mathtt{B}} \in \hat{B}^k, 
\text{supp}\,\hat{\mathtt{B}} \,\cap\,\hat{Q}' \ne \emptyset
\text{ and } \text{supp}\,\hat{\mathtt{B}}\,\cap\, \hat{Q} \ne \emptyset
\right\}
\]
is the multilevel support extension consisting of the set of (active and not active) level-$k$ elements contained in the supports of level-$k$ B-splines that are nonzero on $\hat{Q}\in \hat{G}^\ell$, with $0\le k \le \ell$ and $\hat{\cal G}^\ell = \emptyset$ for $\ell<0$. 
Note that, while $\hat{Q}_r \in \hat{\mathcal{G}}^\ell$ is an active level-$\ell$ element to be replaced by its children upon refinement, $\hat{Q}_c \in \hat{\mathcal{G}}^\ell$ is a currently inactive element that may be reactivated through coarsening.

Defining the refinement patch for each marked level-$\ell$ element  $\hat{Q}_r$ through its refinement neighborhood ${\cal N}_r(\hat{\cal Q},\hat{Q}_r,m)$ produces admissible meshes that account for the truncation of any THB-splines introduced at level $\ell-m+1$ with respect to $\ell-m+2$. The idea is the following. When $\hat{Q}_r$  is refined into its children at level $\ell+1$, any THB-spline of level $\ell-m+1$ that is nonzero on $\hat{Q}_r$ must vanish on the children of $\hat{Q}_r$ activated after refinement. Consequently, we recursively mark for refinement also any element $\hat{Q}'\in \hat{\cal G}^{\ell-m+1}$ included in ${\cal N}_r({\cal Q},\hat{Q}_r,m)$. In this way elements of levels $\ell-m+1$ in the refinement neighborhoods of marked elements are recursively refined, 
forcing any THB-spline of level $\ell-m+1$  to be fully truncated with respect to level $\ell -m +2$ and then vanish on the new elements of $\ell+1$ \cite{buffa2016}. An example of admissible refinement in the bivariate case is shown in Fig.~\ref{fig:neigh-b}--\ref{fig:neigh-c}.

When coarsening is required instead of refinement, a similar strategy is adopted to determine ${\cal N}_c({\hat Q},\hat{Q}_c,m)$. In this case, the neighborhood consists of the active level-$(\ell+m)$ elements lying within the supports of active THB-splines that are nonzero on the children of $\hat{Q}_c$. A coarse element is reactivated only if all its children have been marked for coarsening. To ensure admissibility of the hierarchical mesh, a level-$\ell$ element is reactivated if and only if none of the level-$(\ell+m)$ elements in its coarsening neighborhood remains active. As in the admissible refinement procedure, mesh admissibility is enforced through a (recursive) coarsening process applied to elements in the refinement neighborhoods of the marked elements~\cite{carraturo2019}. Since coarsening entails a loss of information, the coarsening algorithm is inherently conservative: (1) an element is reactivated only if all its children are marked for coarsening, and (2) to preserve the admissibility condition, an element is reactivated only if its coarsening neighborhood is empty. Fig.~\ref{fig:neigh-d} illustrates a nonempty coarsening neighborhood, in which the presence of active level-2 elements prevents the coarsening of the marked level-0 element. In addition to \cite{carraturo2019}, the effectiveness of admissible refinement and coarsening with THB-spline discretization was recently investigated in \cite{Bracco2023,carraturo2024,ventavinuela2026}.

\subsection{Spatial discretization}\label{Spacial_discretization}

Let $\Omega \subset \mathbb{R}^d$ be the physical domain.
The isogeometric mapping $\mathbf{F} : \hat{\Omega}^0 \rightarrow \Omega$ maps the parametric domain $\hat{\Omega}^0$, where the mesh and basis functions are defined, onto the physical space.
For $\bx \in \Omega$, we write
\[
\bx = \mathbf{F}(\boldsymbol{\xi}) = \sum_i \mathbf{P}_i \, \hat{\tau}_i(\boldsymbol{\xi}),
\]
where $\mathbf{P}_i$ denotes the control point coordinates and $\hat{\tau}_i \in \hat{\mathcal{T}}^0$ are the basis functions of the initial tensor-product space. The corresponding physical basis functions are defined as $\tau_i := \hat{\tau}_i \circ \mathbf{F}^{-1}$.

To define the spline approximation over $\Omega$, we first define the domains and active elements at the different levels $\ell=0,\ldots,N-1$, together with the hierarchical mesh in the physical space ($\mathcal{Q}$) as
\[
\Omega^\ell = \mathbf{F}(\hat{\Omega}^\ell), 
\qquad 
\mathcal{G}^\ell = \{ Q \in \mathcal{Q} : \hat{Q} \in \hat{\mathcal{G}}^\ell \}, \qquad 
\mathcal{Q} = \{ Q = \mathbf{F}(\hat{Q}) : \hat{Q} \in \hat{\mathcal{Q}} \}.
\]

For a given hierarchical mesh, the basis functions $\tau_i(\bx)$ are organized into the system matrices following the Galerkin method.
The field variables $\phi$ and $\sigma$ introduced in Section~\ref{model_fromulation_numerical_methods} are approximated as follows: %
\begin{align}
    \phi(\bx, t) & = \bN(\bx)^T \tilde{\boldsymbol{\phi}}, \quad \delta\phi(\bx, t)  = \bN(\bx)^T  \delta\tilde{\boldsymbol{\phi}}, \label{phi_approx}\\
    \sigma(\bx) & =  \bN(\bx)^T  \tilde{\boldsymbol{\sigma}}, \quad \delta\sigma(\bx)  = \bN(\bx)^T  \delta\tilde{\boldsymbol{\sigma}}, \label{sigma_approx} 
\end{align}
where $\bN(\bx)$ is the vector of the basis functions $\tau_i(\bx)$, and $\tilde{\boldsymbol{\phi}}$ and $\tilde{\boldsymbol{\sigma}}$ are the corresponding vectors of the control variables.
The requirement of $H^2(\Omega)$ conforming spaces for $\phi$ is satisfied by the use of THB-spline basis with at least $C^1$ inter-element continuity, enabling a direct discretization of the fourth-order CH operator without introducing auxiliary variables.
Although the nutrient equation (Eq.~\eqref{Nutri_Eq}) only requires $H^1(\Omega)$-conforming spaces, an identical spline approximation space is used for $\sigma$. %
This ensures consistency, simplicity in the implementation of hierarchical refinement strategies, and numerical stability of the coupled multiphysics formulation in the staggered solution scheme.

Substitution of Eqs.~\eqref{phi_approx} and \eqref{sigma_approx} into Eqs.~\eqref{CH_Weak_Final} and \eqref{Nutri_weak_form}, together with the Nitsche's penalty terms, results in the following semi-discrete residual forms: %
\begin{align}
    \bR_\phi(\tilde{\boldsymbol{\phi}}, \dot{\tilde{\boldsymbol{\phi}}}, \tilde{\boldsymbol{\sigma}}) & = \bM \dot{\tilde{\boldsymbol{\phi}}} + \bM_n \tilde{\boldsymbol{\phi}} 
    + \bF_1(\tilde{\boldsymbol{\phi}}) + \bF_2(\tilde{\boldsymbol{\phi}}) + \bF_3(\tilde{\boldsymbol{\phi}}) - \bF_4(\tilde{\boldsymbol{\phi}}) - \bF_n(\tilde{\boldsymbol{\phi}})  
    - \bF_N(\tilde{\boldsymbol{\sigma}}) - \bF_S(\tilde{\boldsymbol{\phi}},\tilde{\boldsymbol{\sigma}}) = \bf{0}, \label{R_phi}\\
    \bR_\sigma(\tilde{\boldsymbol{\sigma}}, \tilde{\boldsymbol{\phi}}) &= \left(D \bK_{\nabla} + \bM^* \right) \tilde{\boldsymbol{\sigma}} - D\mathcal{X} \bK_\nabla \tilde{\boldsymbol{\phi}} = \bf{0}, \label{R_sigma}
\end{align}
where $\dot{\tilde{\boldsymbol{\phi}}} = \dparder{\tilde{\boldsymbol{\phi}}}{t}$ and the integrals correspond to the domain contributions are 
\begin{gather*}
    \bM = \int_{\Omega} \bN \bN^T \diff \Omega, \quad 
    \bF_1(\tilde{\boldsymbol{\phi}}) = \int_{\Omega} (\nabla \bN) M(\phi) Ef''(\phi) \nabla\bN^T \tilde{\boldsymbol{\phi}} \diff \Omega \quad
    \bF_2(\tilde{\boldsymbol{\phi}}) = \int_{\Omega} \Delta \bN \lambda M(\phi) \Delta\bN^T \tilde{\boldsymbol{\phi}} \diff \Omega, \\
    \bF_3(\tilde{\boldsymbol{\phi}}) = \int_{\Omega} \nabla \bN \lambda M'(\phi) \nabla \bN^T \tilde{\boldsymbol{\phi}} \Delta \bN^T \tilde{\boldsymbol{\phi}} \diff \Omega, \quad
    \bF_4(\tilde{\boldsymbol{\phi}}) = \int_{\Gamma} \nabla\bN \bn \lambda M(\phi) \Delta \bN^T \tilde{\boldsymbol{\phi}} \diff \Gamma, \quad
    \bF_N(\tilde{\boldsymbol{\sigma}}) = \int_{\Omega} \nabla\bN M(\phi) \mathcal{X}_\phi \nabla\bN^T \tilde{\boldsymbol{\sigma}} \diff \Omega, \\
    \bF_S(\tilde{\boldsymbol{\phi}},\tilde{\boldsymbol{\sigma}}) = \int_{\Omega} \bN S(\phi,\sigma) \diff \Omega, \quad
    \bK_{\nabla} = \int_{\Omega} \nabla\bN \nabla\bN^T \diff \Omega, \quad \bM^* =  \int_{\Omega} \bN \tilde{S}^*(\phi) \bN^T \diff \Omega.
\end{gather*}
The Nitsche's boundary integrals, enforcing $\nabla\phi\cdot\bn = \bf{0}$ weakly, are 
\begin{gather*}
    \bM_n = \int_{\Gamma} \nabla\bN \bn \frac{\varepsilon_n}{h_e} \left(\nabla \bN \bn\right)^T \diff \Gamma, \quad \bF_n(\tilde{\boldsymbol{\phi}}) = \int_{\Gamma} \Delta \bN \lambda M(\phi) (\nabla\bN \bn)^T  \tilde{\boldsymbol{\phi}} \diff \Gamma.
\end{gather*}
The associated functions in above expressions can be evaluated as 
\begin{gather*}
    M(\phi) = M \left[\dfrac{1}{2} (1+\bN^T \tilde{\boldsymbol{\phi}})^2\right] , \quad 
    M'(\phi) = M (1+\bN^T \tilde{\boldsymbol{\phi}}), \quad 
    f''(\phi) =  \left[ 3 (\bN^T \tilde{\boldsymbol{\phi}})^2 - 1 \right], \\
    S(\phi,\sigma) = \frac{\kappa}{2}(\mathcal{P} \bN^T\tilde{\boldsymbol{\sigma}} - \mathcal{A})(\bN^T\tilde{\boldsymbol{\phi}} + 1) ,\quad
    \tilde{S}^*(\phi) =  \frac{\mathcal{C}}{2} (\bN^T\tilde{\boldsymbol{\phi}} + 1).
\end{gather*}
After defining the semi-discrete residuals, Eq.~\eqref{R_phi} is advanced in time using the generalized-$\alpha$ method with the nonlinear system solved using Newton--Raphson (NR) iterations.
As Eq.~\eqref{R_sigma} is linear and quasi-static, it is solved directly at each time step.

\subsection{Time-stepping scheme}

The generalized-$\alpha$ method, which has proven to be effective for the CH equation \cite{Gomez2008, Bracco2023}, is employed to solve Eq.~\eqref{R_phi}. 
The framework can be formulated as follows:
Given $\tilde{\boldsymbol{\phi}}_{n}$, $\dot{\tilde{\boldsymbol{\phi}}}_{n}$, $\tilde{\boldsymbol{\sigma}}_{n}$, and timestep $\Delta t$, find $\tilde{\boldsymbol{\phi}}_{n+1}$, $\dot{\tilde{\boldsymbol{\phi}}}_{n+1}$, $\tilde{\boldsymbol{\phi}}_{n+\alpha_f}$, $\dot{\tilde{\boldsymbol{\phi}}}_{n+\alpha_m}$ such that
\begin{equation} \label{R_phi_time_march}
    \bR_\phi(\tilde{\boldsymbol{\phi}}_{n+\alpha_f}, \dot{\tilde{\boldsymbol{\phi}}}_{n+\alpha_m}, \tilde{\boldsymbol{\sigma}}_n) = \bf{0},
\end{equation}
where the vectors $\tilde{\boldsymbol{\phi}}_{n+\alpha_f}$, $\dot{\tilde{\boldsymbol{\phi}}}_{n+\alpha_m}$, and $\tilde{\boldsymbol{\phi}}_{n+1}$ are defined as
\begin{align}
    \tilde{\boldsymbol{\phi}}_{n+\alpha_f} & = \tilde{\boldsymbol{\phi}}_{n} + \alpha_f \left( \tilde{\boldsymbol{\phi}}_{n+1} - \tilde{\boldsymbol{\phi}}_{n} \right), \label{phi_nalpha}\\
    \dot{\tilde{\boldsymbol{\phi}}}_{n+\alpha_m} & = \dot{\tilde{\boldsymbol{\phi}}}_n + \alpha_m \left( \dot{\tilde{\boldsymbol{\phi}}}_{n+1} - \dot{\tilde{\boldsymbol{\phi}}}_n \right),  \label{phidot_nalpha}\\
    \tilde{\boldsymbol{\phi}}_{n+1} & = \tilde{\boldsymbol{\phi}}_{n} + \Delta t \dot{\tilde{\boldsymbol{\phi}}}_n + \gamma \Delta t \left( \dot{\tilde{\boldsymbol{\phi}}}_{n+1} - \dot{\tilde{\boldsymbol{\phi}}}_n \right). \label{phi_nplus1}
\end{align}
The parameters $\alpha_f$, $\alpha_m$, and $\gamma$ that define the method are chosen in accordance with \cite{Gomez2008, Bracco2023} to ensure accuracy and stability; they are given as
\begin{equation*}
    \alpha_f = \frac{1}{1+\rho_\infty}, \: \alpha_m = \frac{1}{2}\left( \frac{3 - \rho_\infty}{1+\rho_\infty} \right), \text{ and } \gamma = \frac{1}{2} +\alpha_m - \alpha_f,
\end{equation*}
where the spectral radius of the amplification matrix is $\rho_\infty = 0.5$.

To advance in time, Eq.~\eqref{R_phi_time_march} is solved at each time step using an iterative NR scheme due to its nonlinearity. The iteration begins with initial estimates $\tilde{\boldsymbol{\phi}}_{n+1}^0$ and $\dot{\tilde{\boldsymbol{\phi}}}_{n+1}^0$, forming the predictor stage as follows:
\begin{align}
    \tilde{\boldsymbol{\phi}}_{n+1}^0 &= \tilde{\boldsymbol{\phi}}_{n}, \label{predictor1}\\
    \dot{\tilde{\boldsymbol{\phi}}}_{n+1}^0 &= \frac{\gamma-1}{\gamma} \dot{\tilde{\boldsymbol{\phi}}}_{n}, \label{predictor2}
\end{align}
which are iteratively updated in the multicorrector stage using the following relations until the residual falls below a specified tolerance:
\begin{align}
     \dot{\tilde{\boldsymbol{\phi}}}_{n+1}^{(i+1)} &= \dot{\tilde{\boldsymbol{\phi}}}_{n+1}^{(i)} + \Delta\dot{\tilde{\boldsymbol{\phi}}}_{n+1}^{(i+1)}, \label{phidot_nplus1_i}\\
     \tilde{\boldsymbol{\phi}}_{n+1}^{(i+1)} &= \tilde{\boldsymbol{\phi}}_{n+1}^{(i)} + \gamma \Delta t \Delta\dot{\tilde{\boldsymbol{\phi}}}_{n+1}^{(i+1)}, \label{phi_nplus1_i}
\end{align}
where iterations are indexed with $i = 0, 1, 2, \dots, i_{\text{max}}$ and $\Delta\dot{\tilde{\boldsymbol{\phi}}}_{n+1}^{(i+1)}$ is evaluated by solving the following linear system:
\begin{equation}
    \bK_{\text{NR}}^{(i+1)} \Delta\dot{\tilde{\boldsymbol{\phi}}}_{n+1}^{(i+1)} = - \bR_\phi^{(i+1)}. \label{knablaphi_R}
\end{equation}
The tangent matrix $\bK_{\text{NR}}^{(i+1)}$ is computed using 
\begin{equation}
    \bK_{\text{NR}}^{(i+1)} = \alpha_m \parder{\bR_\phi(\tilde{\boldsymbol{\phi}}_{n+\alpha_f}^{(i+1)}, \dot{\tilde{\boldsymbol{\phi}}}_{n+\alpha_m}^{(i+1)}, \tilde{\boldsymbol{\sigma}}_n)}{\dot{\tilde{\boldsymbol{\phi}}}_{n+\alpha_m}^{(i+1)}} 
    + \alpha_f \gamma \Delta t \parder{\bR_\phi(\tilde{\boldsymbol{\phi}}_{n+\alpha_f}^{(i+1)}, \dot{\tilde{\boldsymbol{\phi}}}_{n+\alpha_m}^{(i+1)}, \tilde{\boldsymbol{\sigma}}_n)}{\tilde{\boldsymbol{\phi}}_{n+\alpha_f}^{(i+1)}}. \label{tangentMatK}
\end{equation}

Once convergence for $\tilde{\boldsymbol{\phi}}_{n+1}$ is achieved, it is used in the linear system $\bR_\sigma(\tilde{\boldsymbol{\sigma}}_{n+1}, \tilde{\boldsymbol{\phi}}_{n+1}) = \bf{0} $ (Eq.~\eqref{R_sigma}) to compute the nutrient concentration at the current time step. 

\subsection{Adaptive refinement and coarsening}
Given an initial tensor-product mesh $\mathcal{Q}_0$, the corresponding THB-spline basis space $\mathcal{T}_0$, an admissibility class $m$, a maximum number of hierarchical refinement levels $N$, and a tolerance $\alpha$ and $\beta$, the initial condition $ \phi(\bx,0) = \phi_0(\bx)  \text{ in } \Omega$ is projected onto $\mathcal{T}_0$ to obtain the initial solution vector $\tilde{\boldsymbol{\phi}}_0$.
To ensure that the initial mesh adequately resolves $\phi_0(\bx)$, we recursively apply the refinement and coarsening procedure driven by the error estimator described in Section~\ref{Error_estimator_sec}.
Elements are marked to attain prescribed admissibility according to the strategy detailed in Section~\ref{sec:admissible_meshes}, and the mesh is iteratively refined until the maximum level is reached. 
This initialization concentrates DOFs in regions of large error while simultaneously maintaining an admissible hierarchical mesh.
The refinement and coarsening passes are performed sequentially to avoid undesired oscillations in the solution.
Once the initial mesh and solution vector are established, the CH-based tumor growth model is advanced in time. 
To reduce computational overhead, refinement and coarsening are performed after prescribed intervals $n_{\mathrm{r}}$ and $n_{\mathrm{c}}$, respectively, rather than at every time step.
At each such interval, the mesh is first refined to capture evolving solution features, and subsequently coarsened to remove superfluous resolution. 
The solution from the previous time step is then projected onto the updated space, which serves as the initial condition for the subsequent time integration interval.
While coarsening, the solution is transferred from one hierarchical mesh to another using an $L^2$ projection \cite{ventavinuela2026}.
The pseudocode summarizing the implementation is presented in the Algorithm~\ref{seudocode_timemarchNR} in ~\ref{Appendix-flowchart}. %

\begin{figure}[!t]
\begin{subfigure}{0.185\columnwidth} \centering
\includegraphics[trim={0cm 0cm 0cm 0cm},clip,width=1\columnwidth]{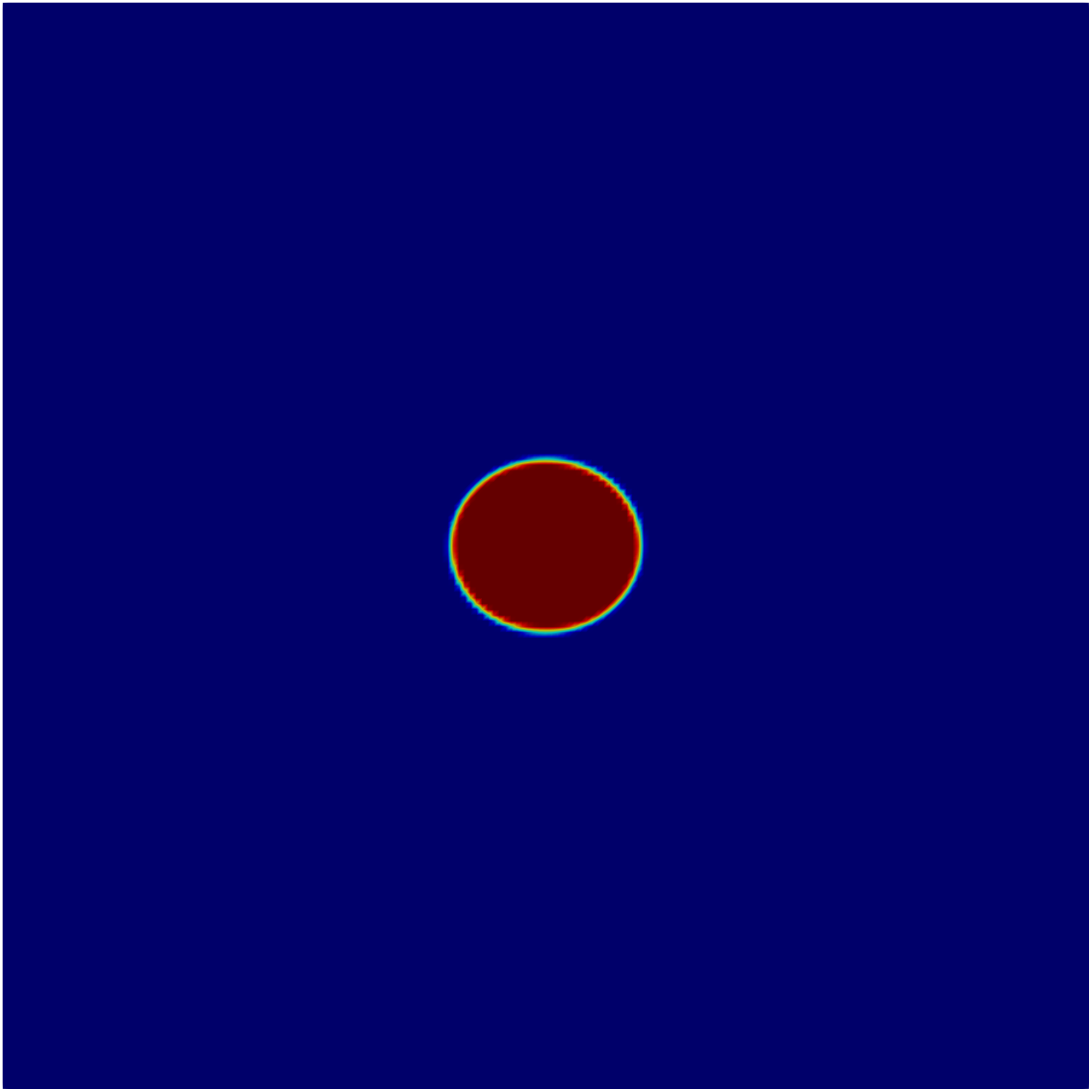}  \end{subfigure}
\begin{subfigure}{0.185\columnwidth} \centering
\includegraphics[trim={0cm 0cm 0cm 0cm},clip,width=1\columnwidth]{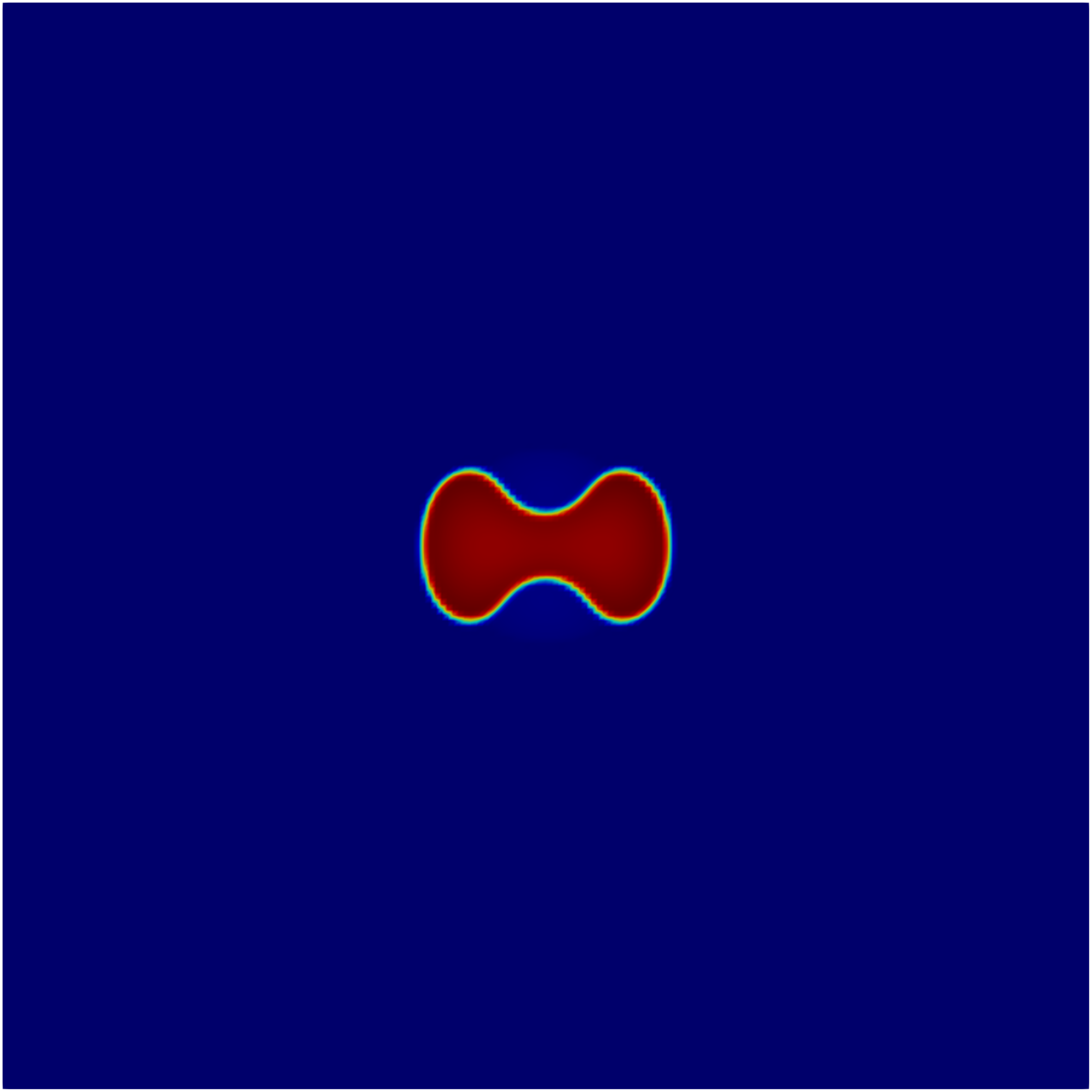}   \end{subfigure}
\begin{subfigure}{0.185\columnwidth} \centering
\includegraphics[trim={0cm 0cm 0cm 0cm},clip,width=1\columnwidth]{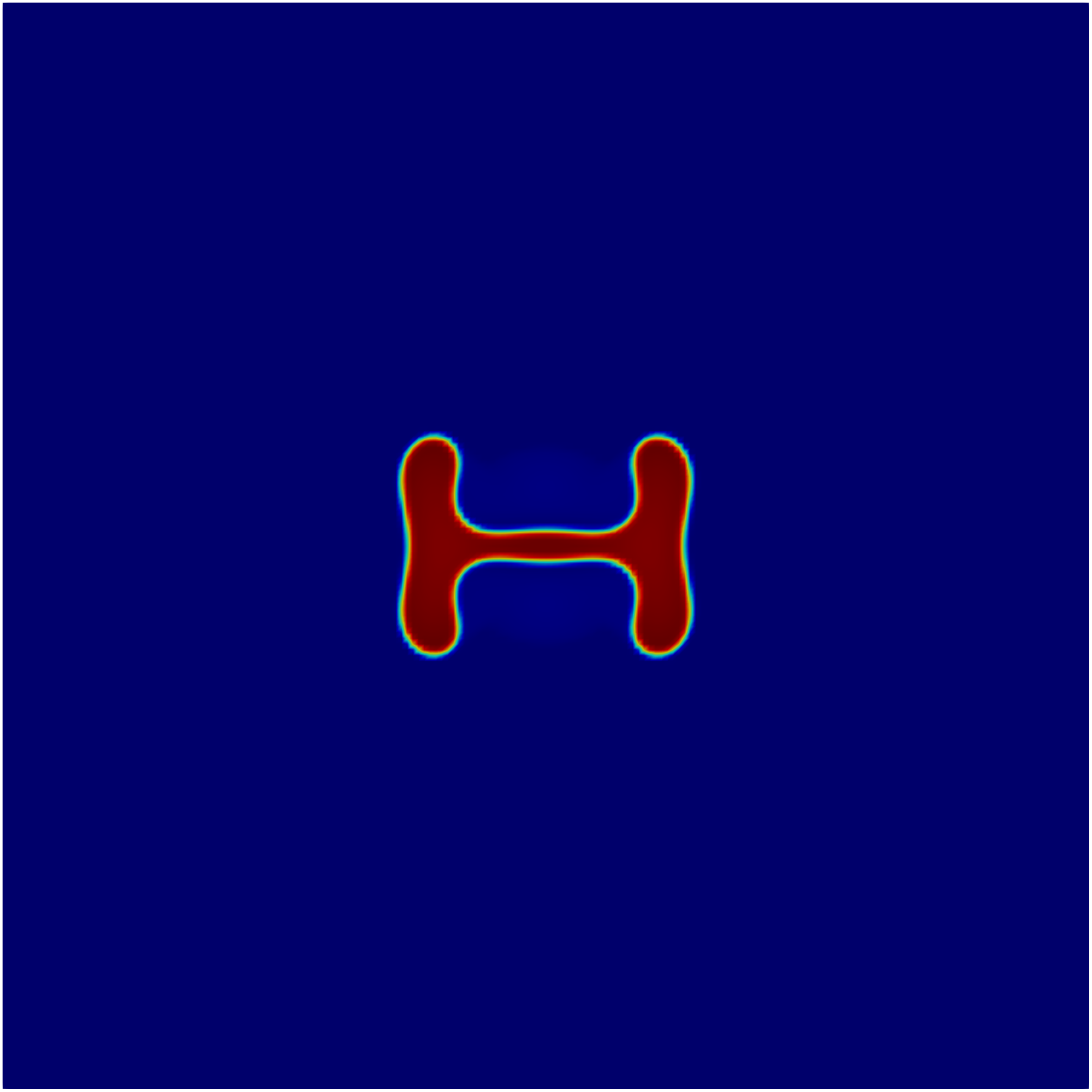} \end{subfigure}
\begin{subfigure}{0.185\columnwidth} \centering
\includegraphics[trim={0cm 0cm 0cm 0cm},clip,width=1\columnwidth]{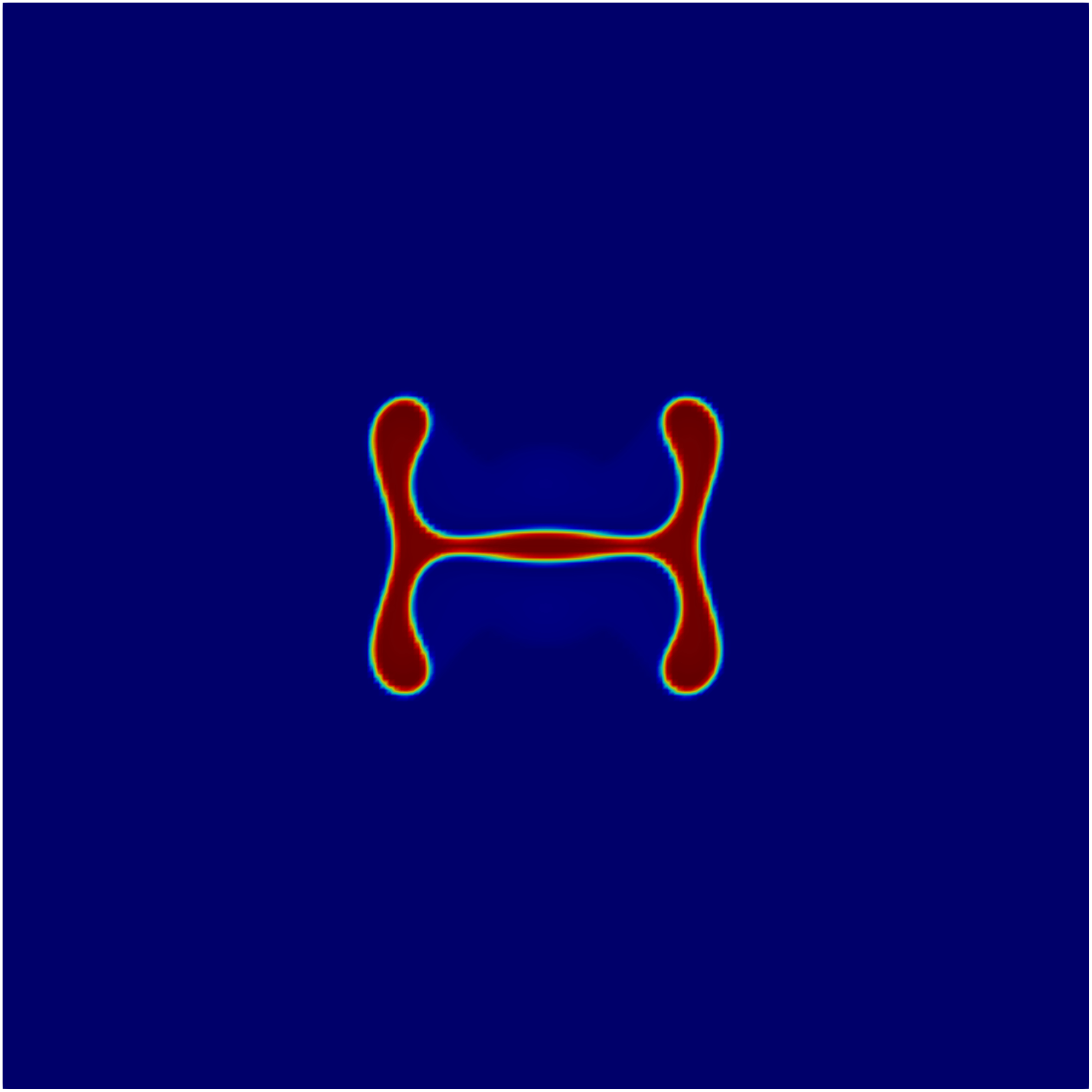} \end{subfigure}
\begin{subfigure}{0.185\columnwidth} \centering
\includegraphics[trim={0cm 0cm 0cm 0cm},clip,width=1\columnwidth]{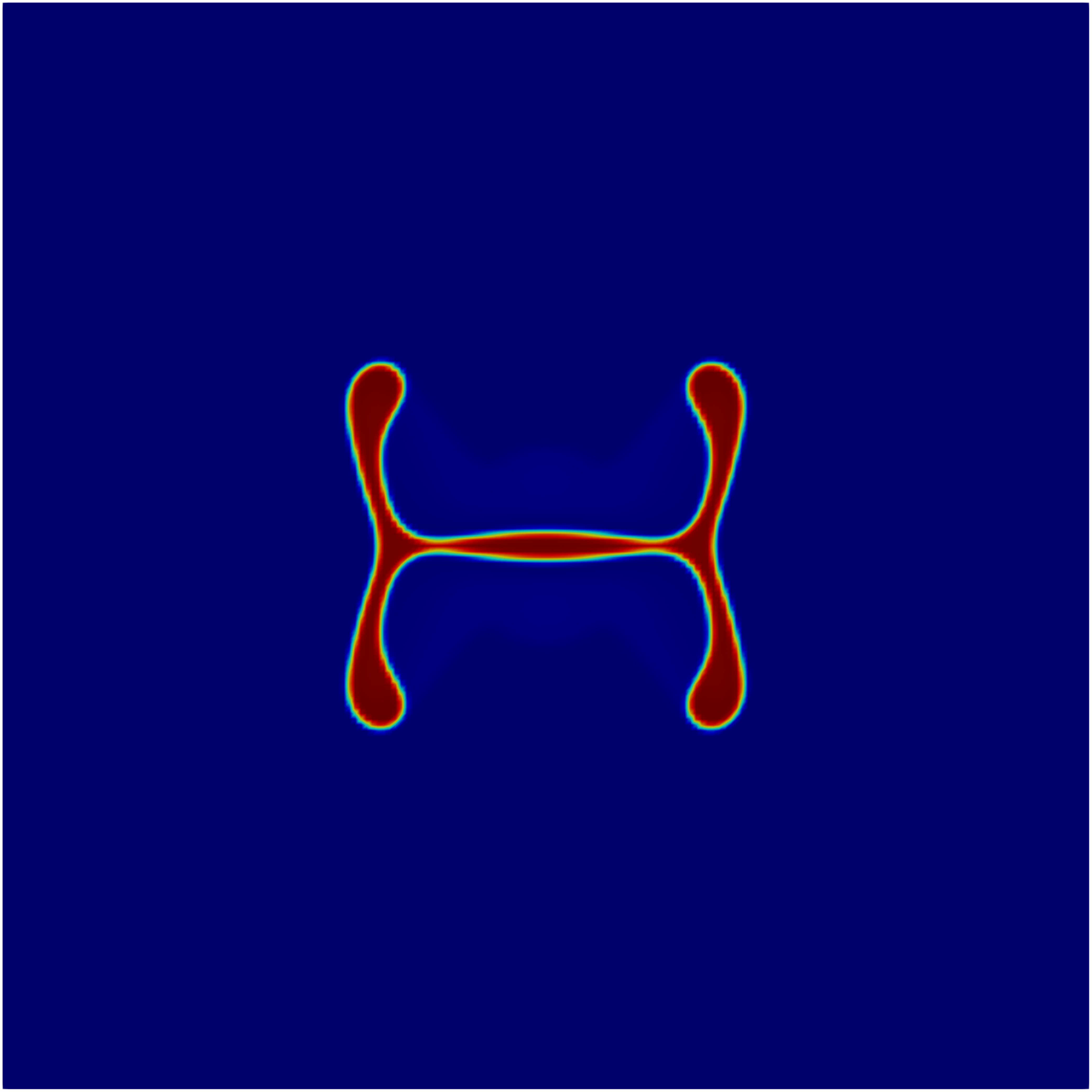}  \end{subfigure}
\begin{subfigure}{0.05\columnwidth} \centering
\includegraphics[trim={0cm 0cm 0cm 0cm},clip,width=1.2\columnwidth]{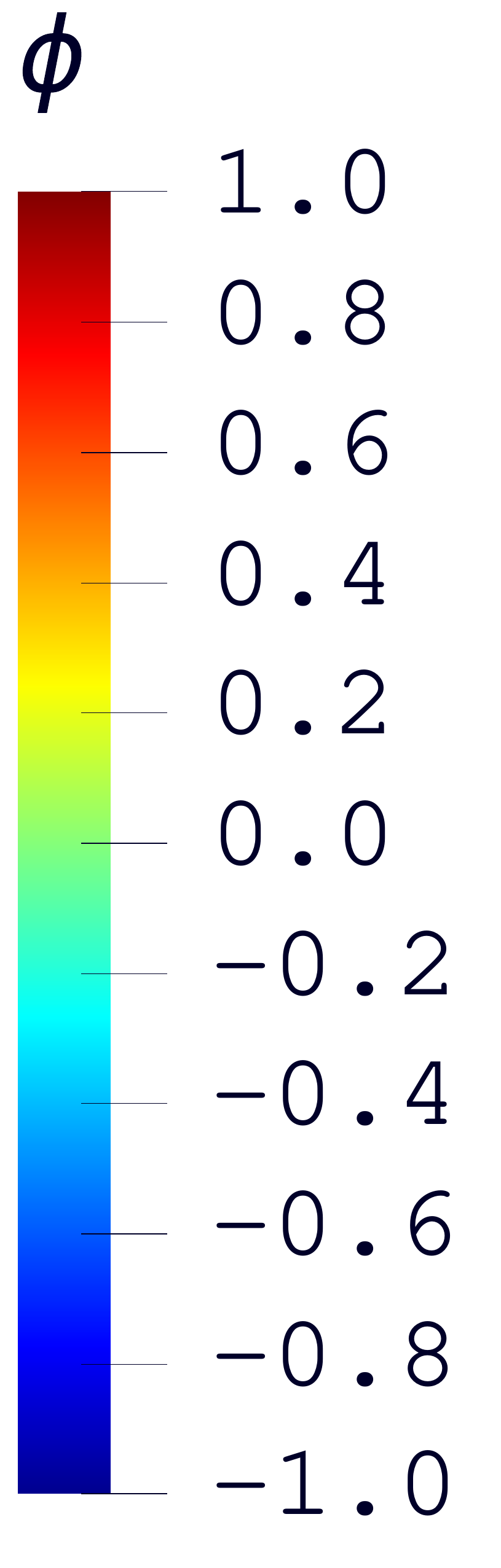} \end{subfigure}
\\
\begin{subfigure}{0.185\columnwidth} \centering
\includegraphics[trim={0cm 0cm 0cm 0cm},clip,width=1\columnwidth]{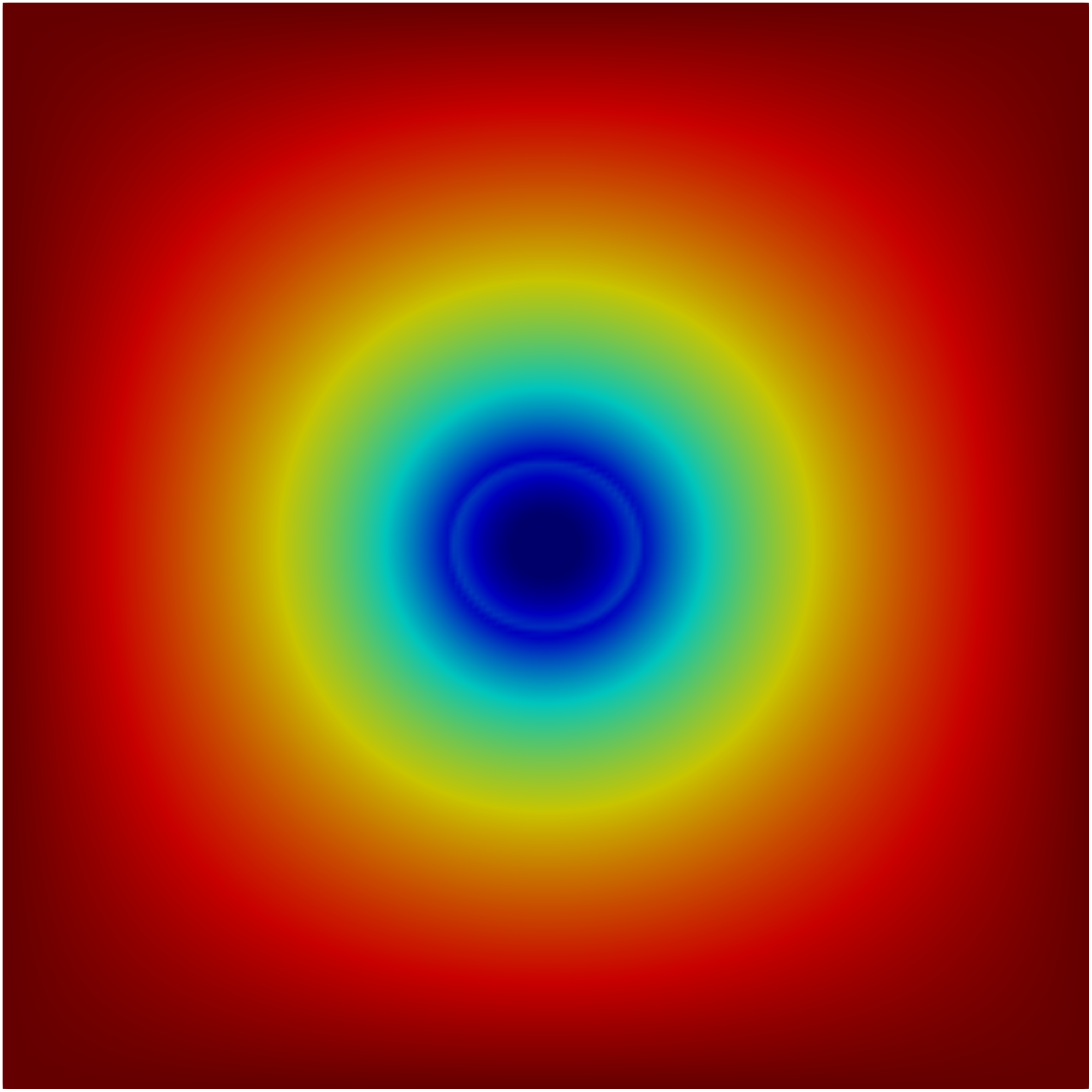} \caption{$t = 0$} \label{fig1_initial_condition}\end{subfigure}
\begin{subfigure}{0.185\columnwidth} \centering
\includegraphics[trim={0cm 0cm 0cm 0cm},clip,width=1\columnwidth]{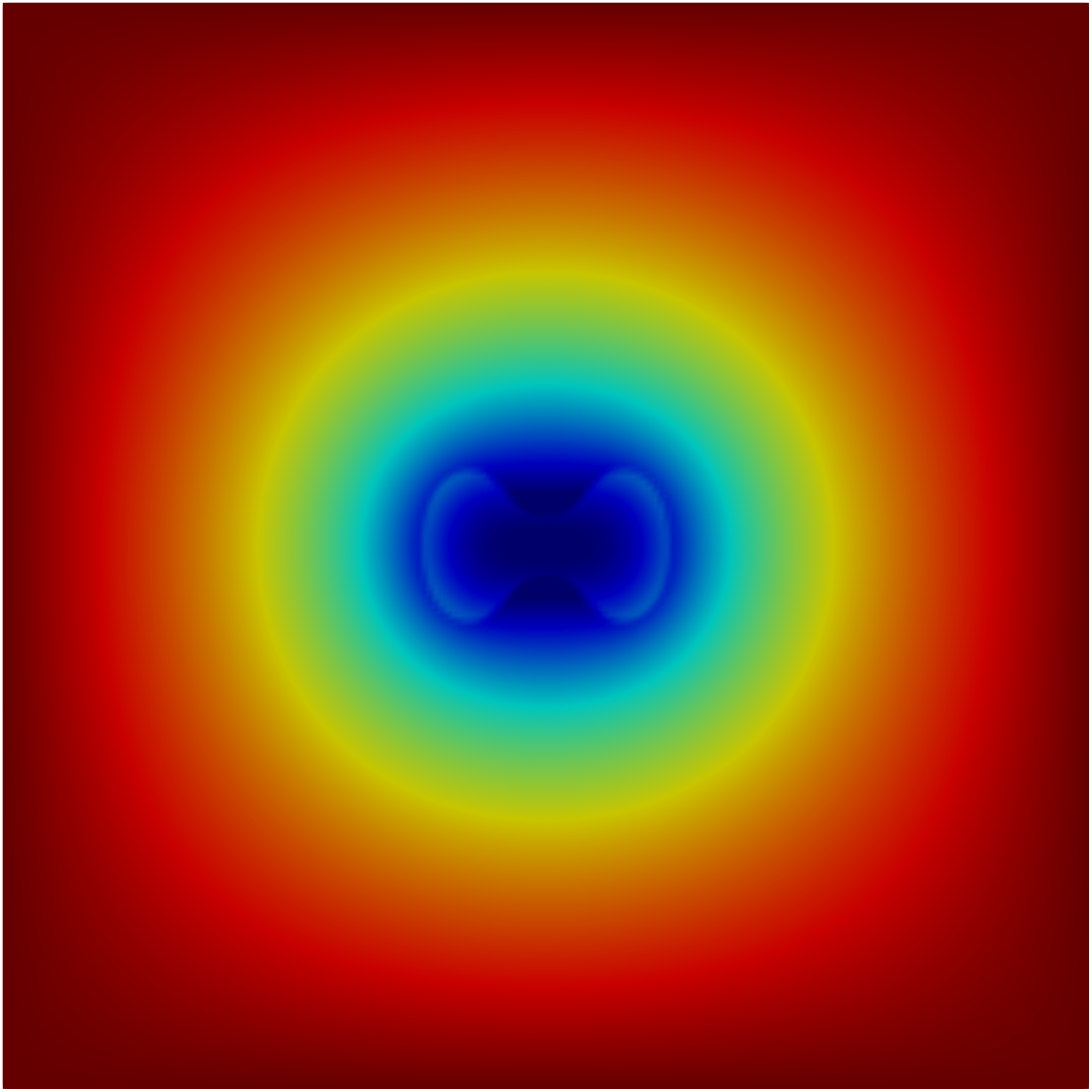} \caption{$t = 1$}  \end{subfigure}
\begin{subfigure}{0.185\columnwidth} \centering
\includegraphics[trim={0cm 0cm 0cm 0cm},clip,width=1\columnwidth]{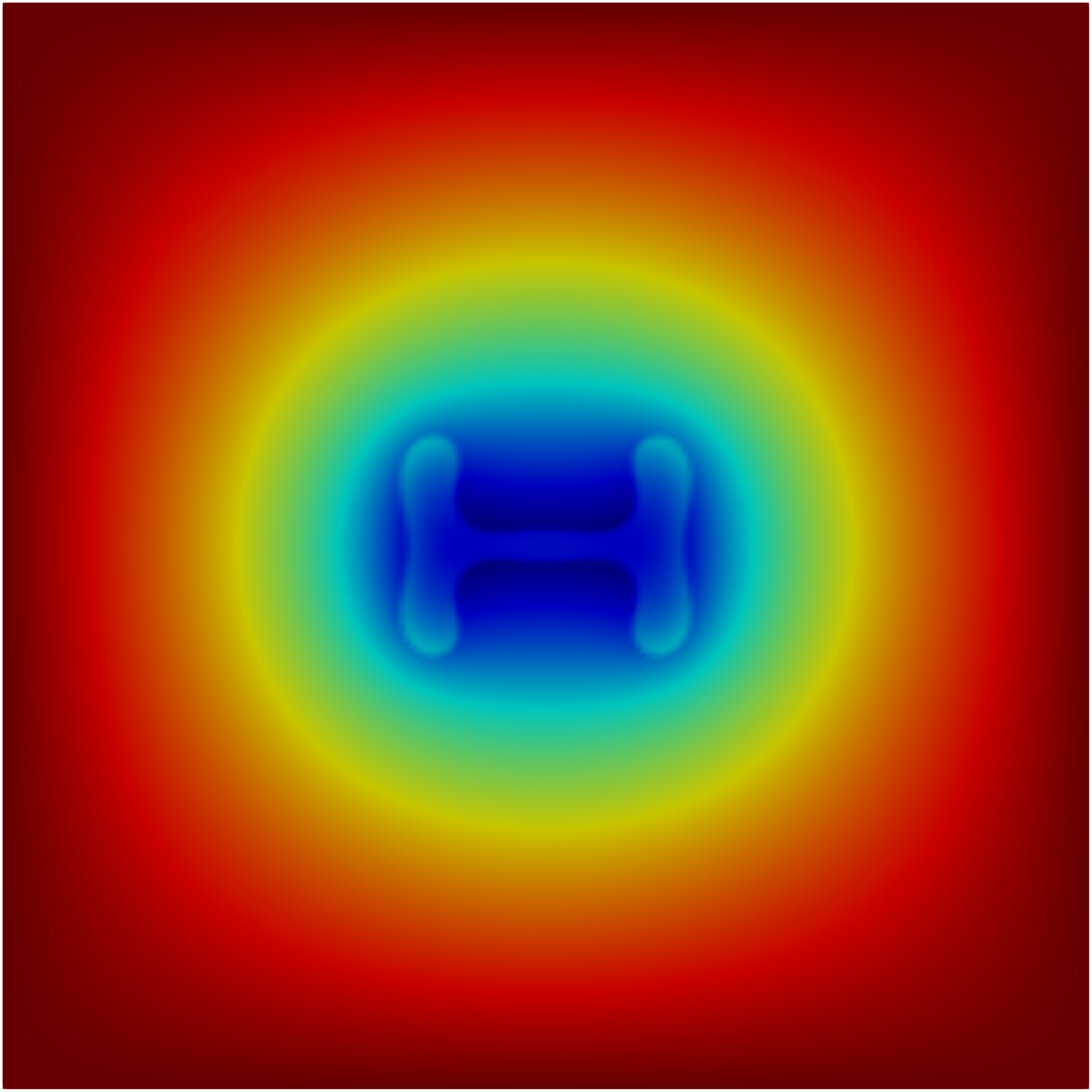} \caption{$t = 1.5$}\end{subfigure}
\begin{subfigure}{0.185\columnwidth} \centering
\includegraphics[trim={0cm 0cm 0cm 0cm},clip,width=1\columnwidth]{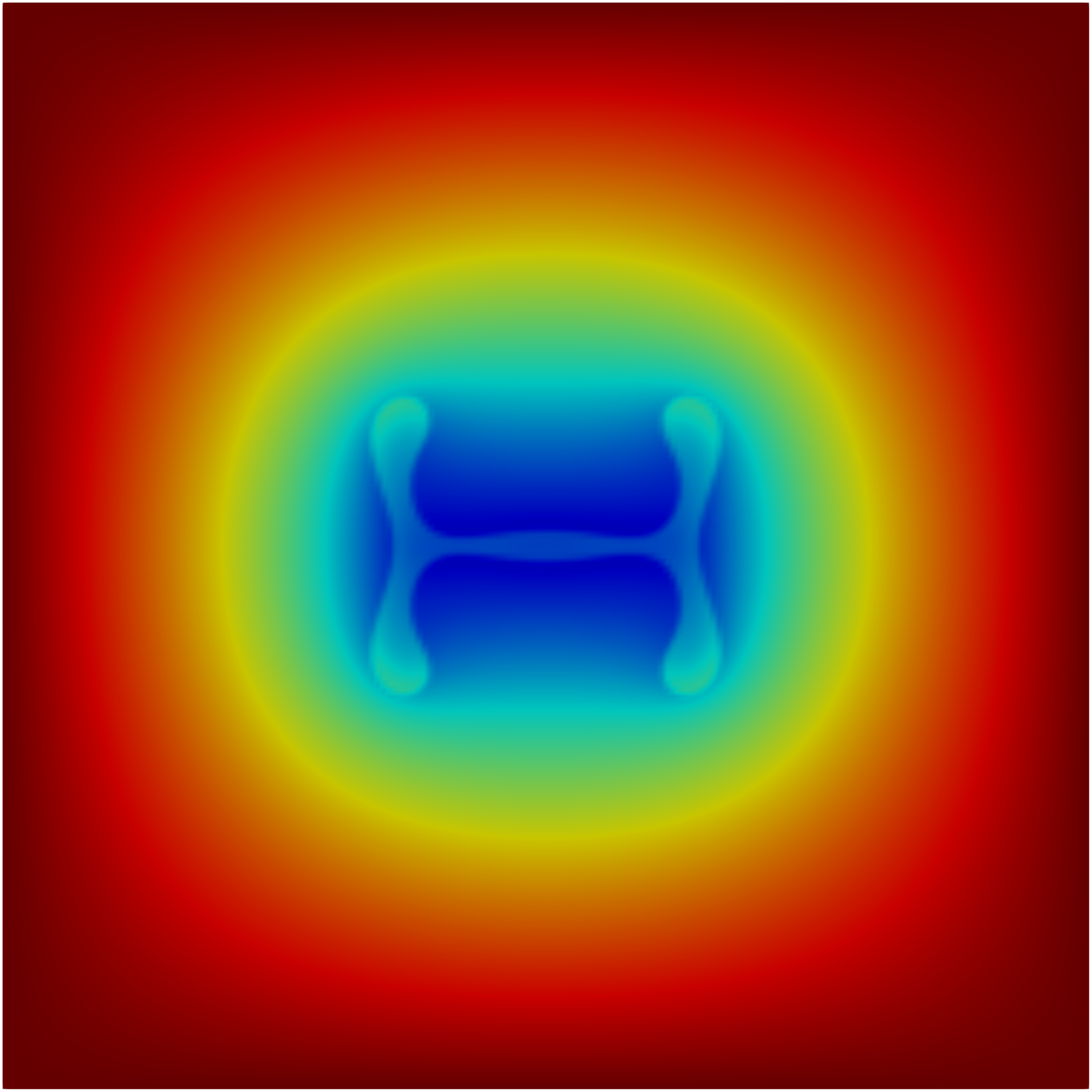} \caption{$t = 2$} \end{subfigure}
\begin{subfigure}{0.185\columnwidth} \centering
\includegraphics[trim={0cm 0cm 0cm 0cm},clip,width=1\columnwidth]{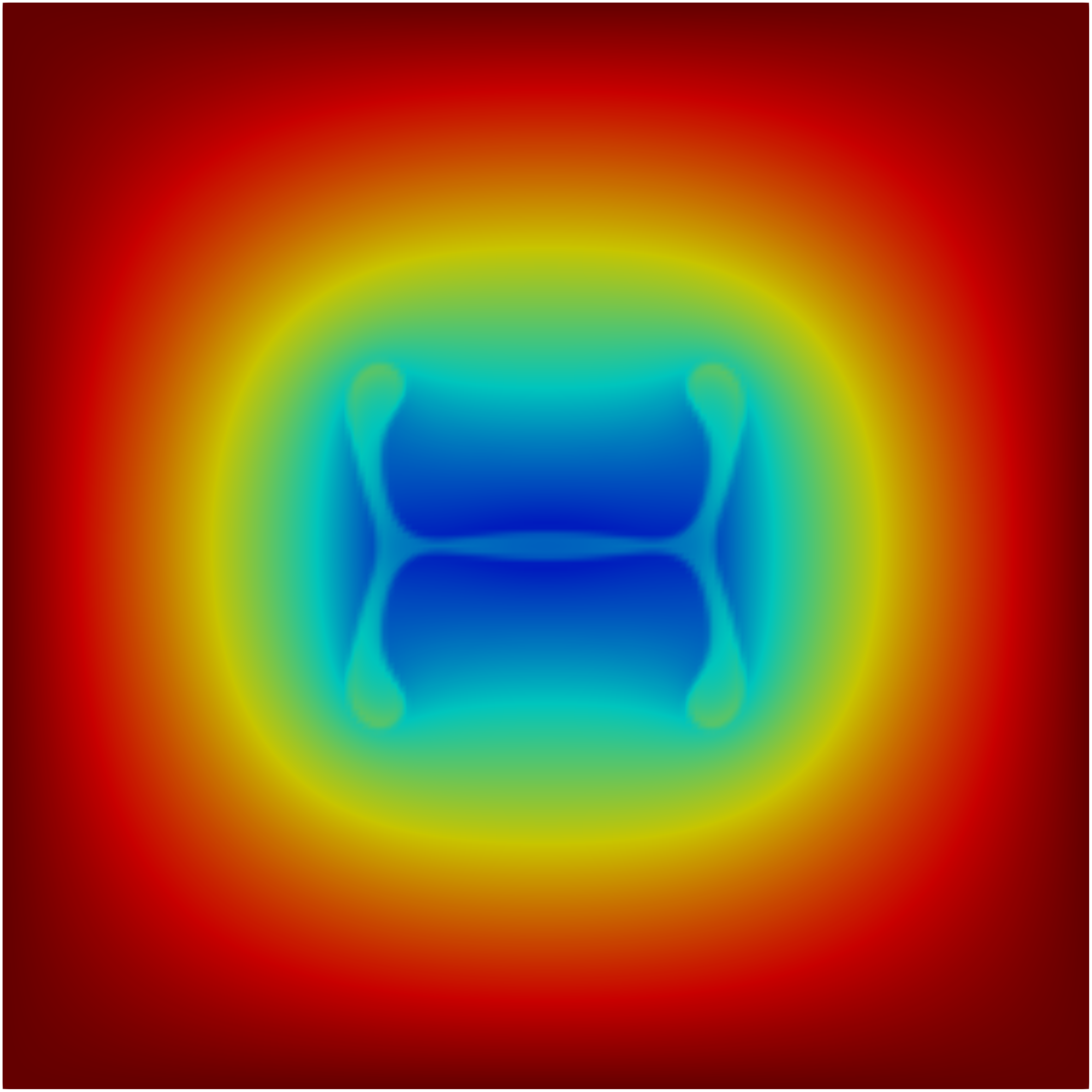} \caption{$t = 2.5$}\end{subfigure}
\begin{subfigure}{0.05\columnwidth} \centering
\includegraphics[trim={0cm -6cm 0cm 0cm},clip,width=1.2\columnwidth]{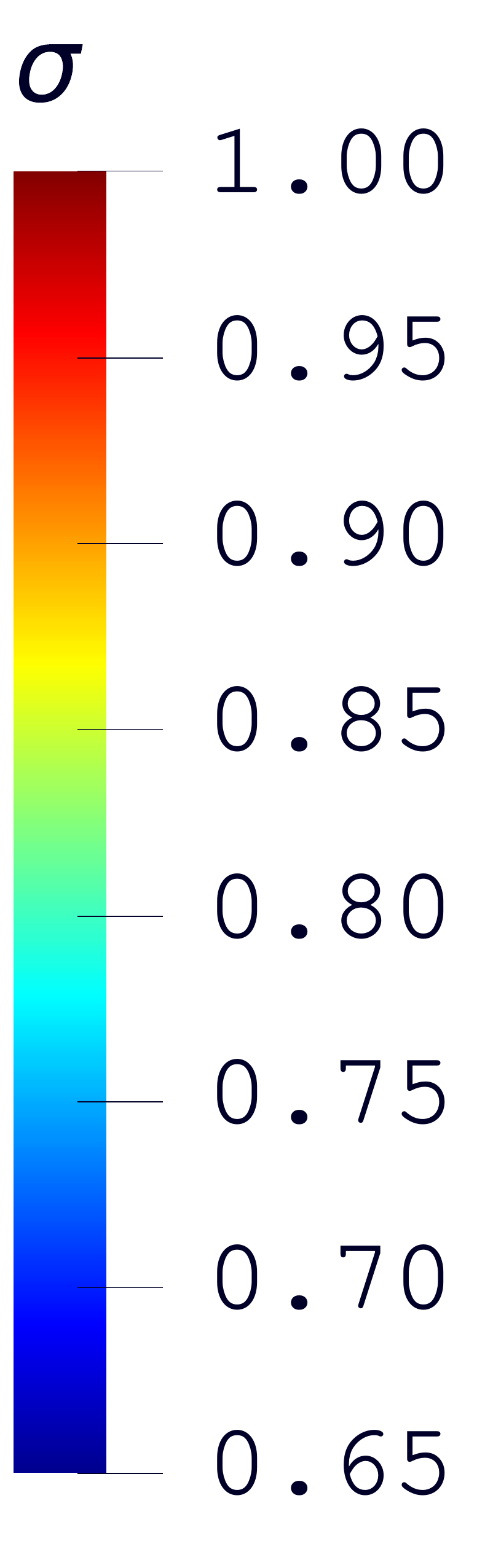} \end{subfigure}
\caption{Evolution of an elliptical tumor (top row) and corresponding nutrient concentration (bottom row) in a square domain $[-3,3]^2$ at time 0, 1, 1.5, 2 and 2.5. This figure shows our overkill solution computed with uniform tensor product B-splines with 1024 elements in each direction and $p=4$.}
\label{figure_1_final}
\end{figure}
\begin{figure}[!t]
\begin{subfigure}{1\columnwidth} \centering
\includegraphics[trim={0cm 0cm 0cm 0cm},clip,width=0.58\columnwidth]{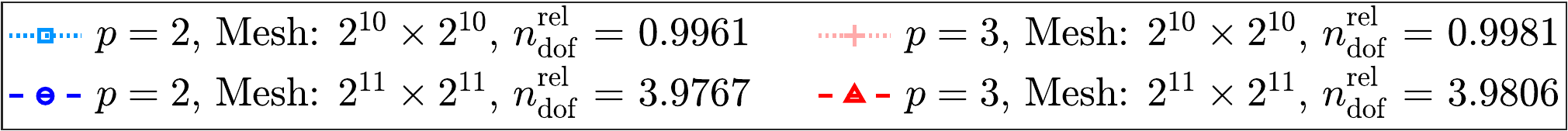}
\end{subfigure}
\begin{subfigure}{0.49\columnwidth} \centering
\includegraphics[trim={0cm 0cm 0cm 0cm},clip,width=1\columnwidth]{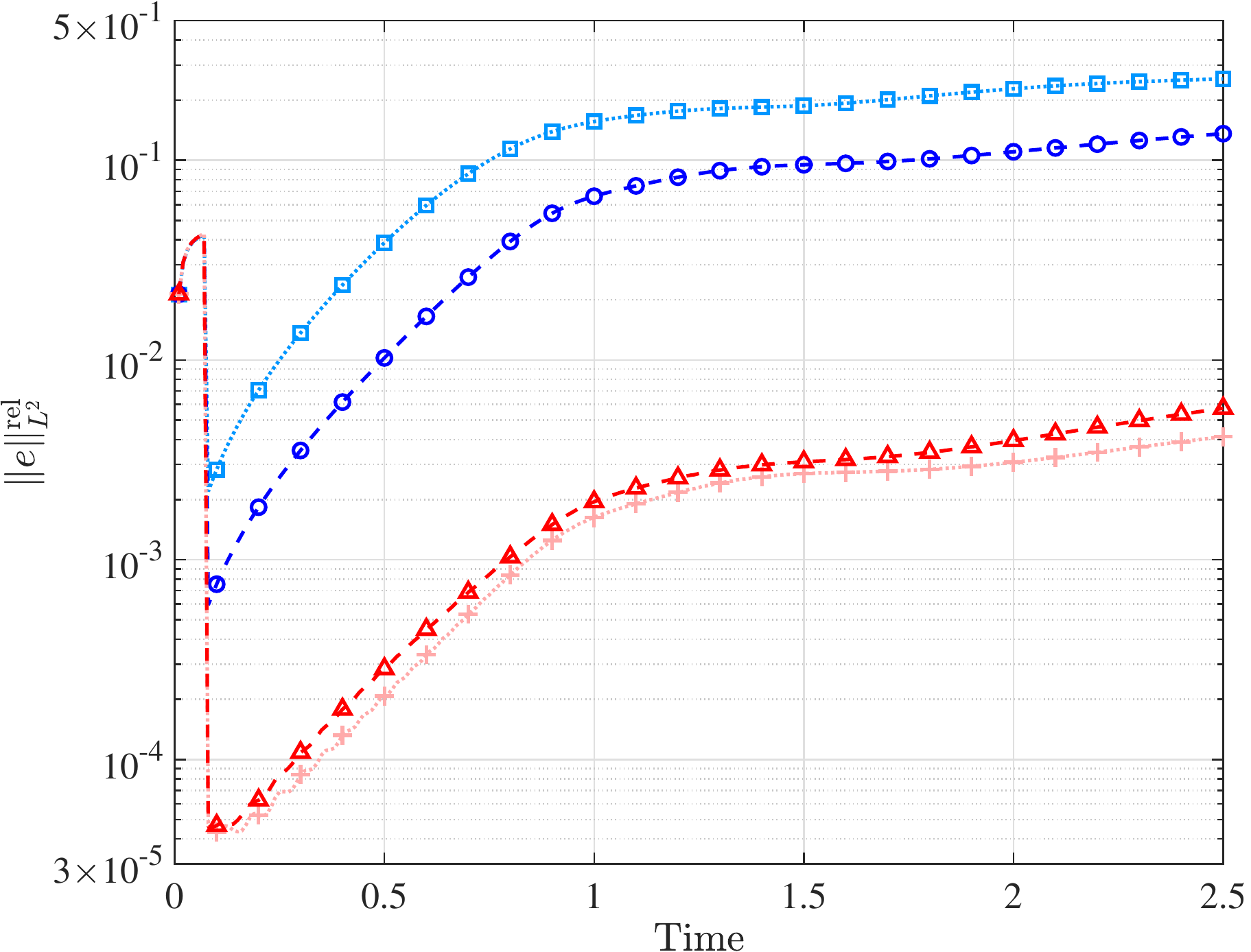}
\caption{$\| e \|_{L^2}^{\mathrm{rel}}$}
\label{}
\end{subfigure}
\begin{subfigure}{0.49\columnwidth} \centering
\includegraphics[trim={0cm 0cm 0cm 0cm},clip,width=1\columnwidth]{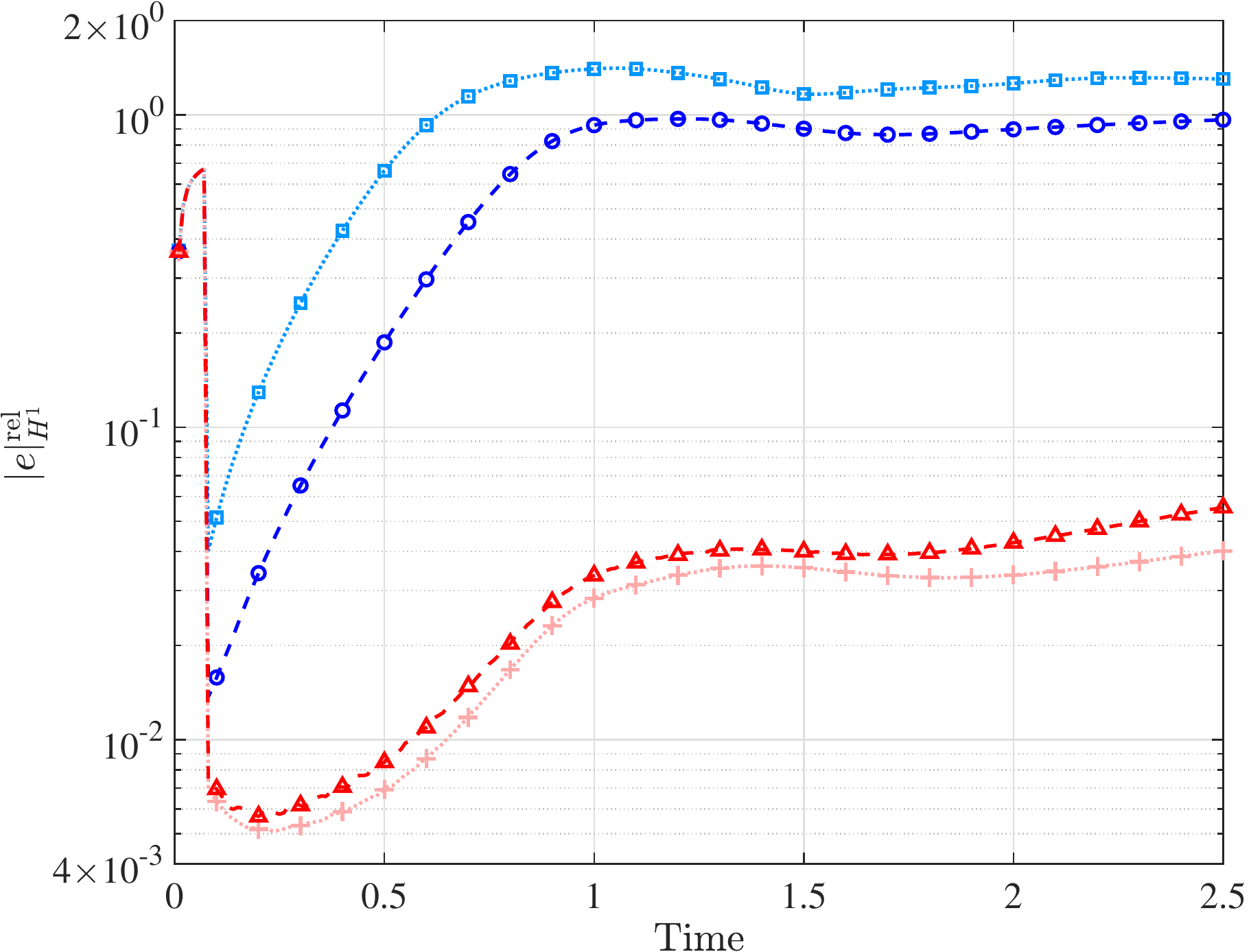}
\caption{$| e |_{H^1}^{\mathrm{rel}}$}
\label{}
\end{subfigure}
\caption{Relative errors $\| e \|_{L^2}^{\mathrm{rel}}$ and $| e |_{H^1}^{\mathrm{rel}}$ over the time period $\left[0, 2.5\right]$ for different tensor product mesh configurations with either $2^{10}$ or $2^{11}$ elements per side and polynomial degree either $p=2$ or $p=3$. The reference solution corresponds to the uniform tensor product B-spline mesh with $p = 4$ and $1024$ elements in each direction shown in Fig.~\ref{figure_1_final}.}
\label{figure_2_final}
\end{figure}
\begin{figure}[!t]\centering
\begin{subfigure}{1\columnwidth} \centering
\includegraphics[trim={0cm -0.35cm 0cm 0cm},clip,width=0.88\columnwidth]{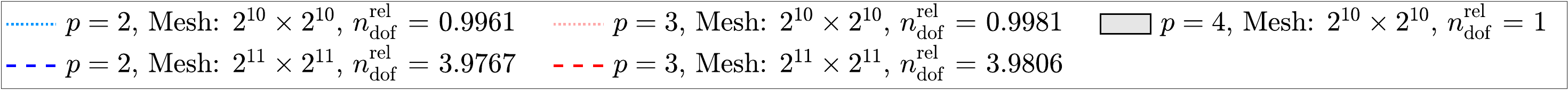}
\end{subfigure}
\begin{subfigure}{0.245\columnwidth} \centering
\includegraphics[trim={0cm 0cm 0cm 0cm},clip,width=1\columnwidth]{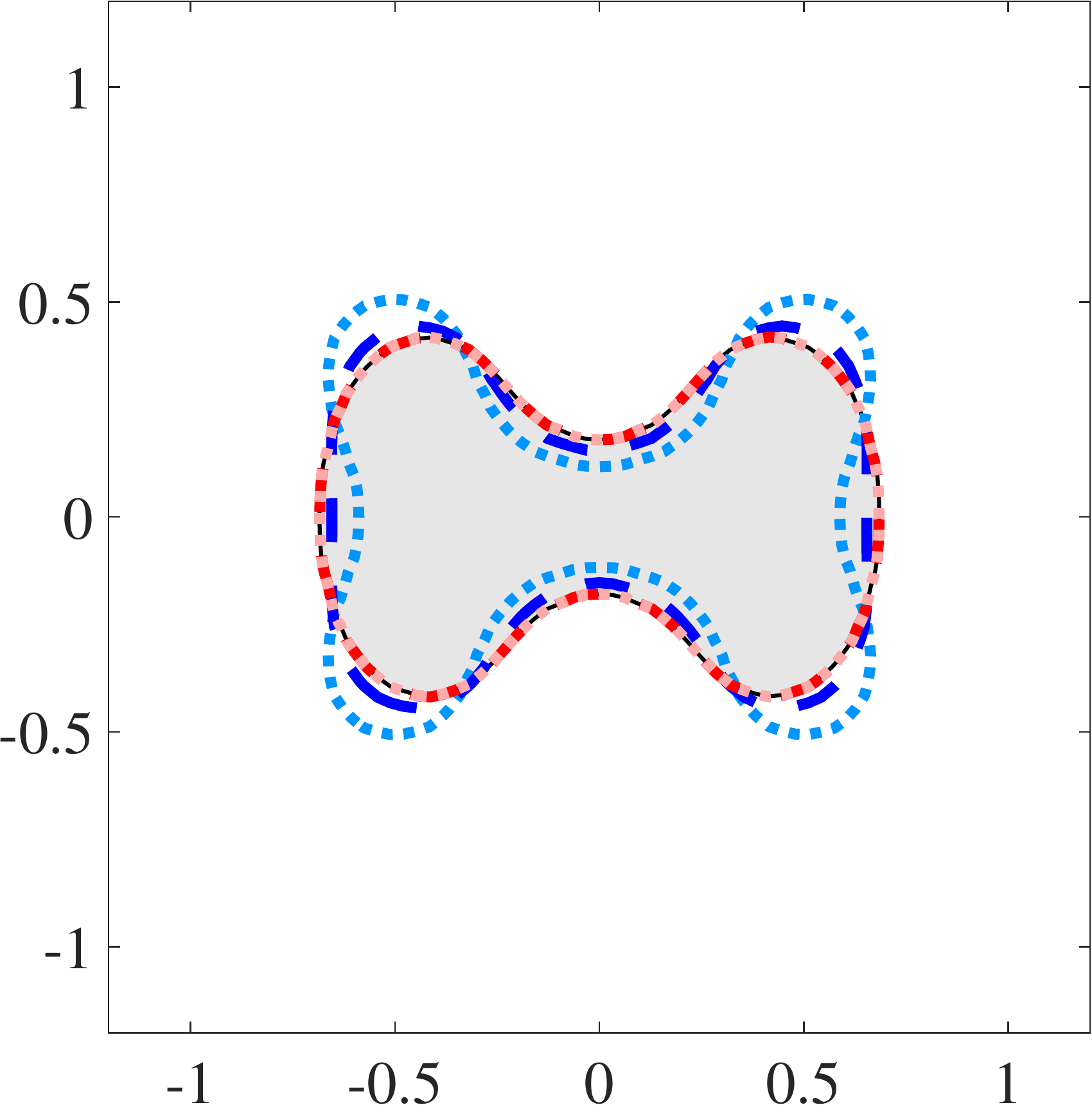}
\caption{$t = 1.0$} %
\end{subfigure}
\begin{subfigure}{0.245\columnwidth} \centering
\includegraphics[trim={0cm 0cm 0cm 0cm},clip,width=1\columnwidth]{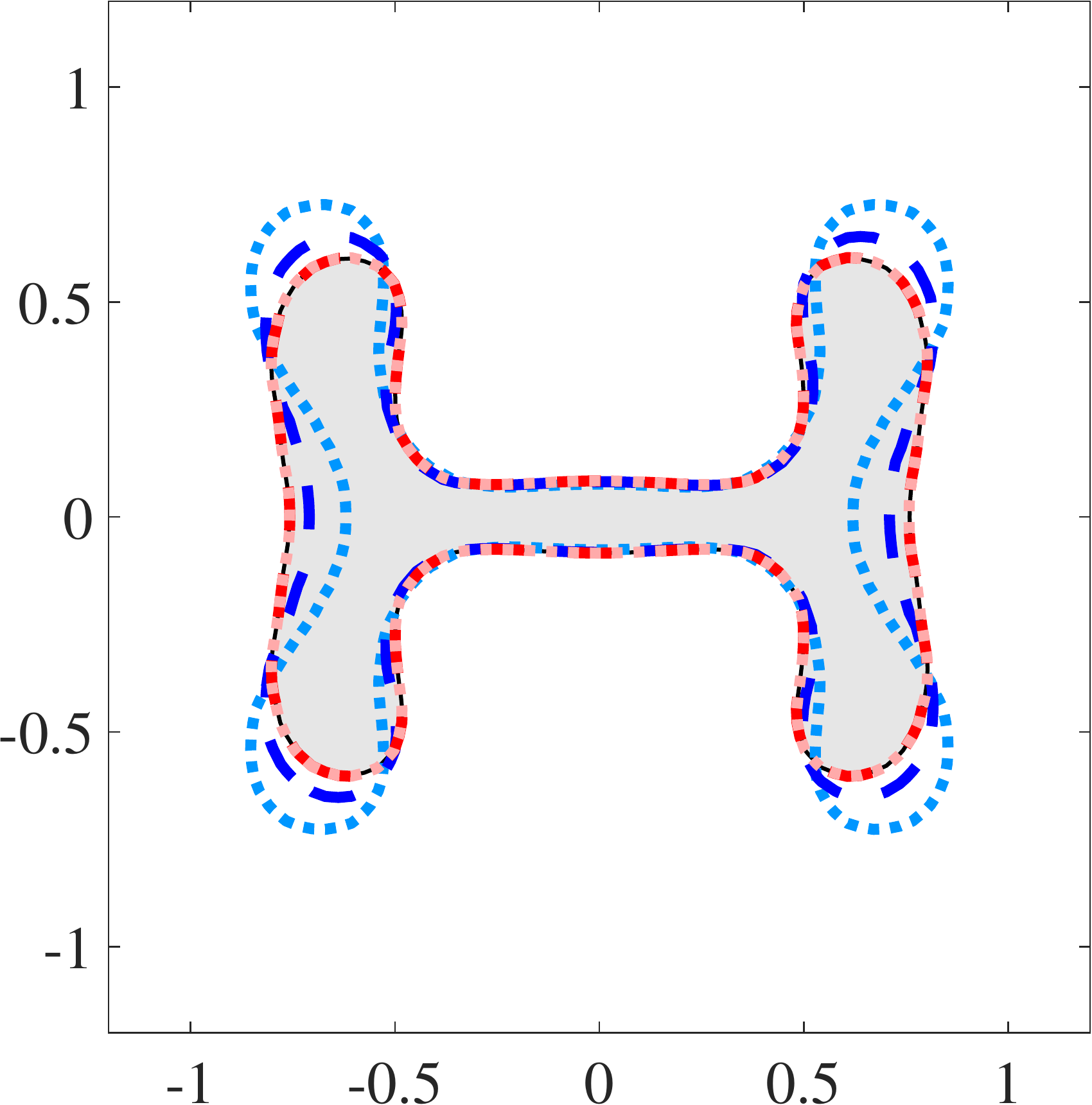}
\caption{$t = 1.5$} %
\end{subfigure}
\begin{subfigure}{0.245\columnwidth} \centering
\includegraphics[trim={0cm 0cm 0cm 0cm},clip,width=1\columnwidth]{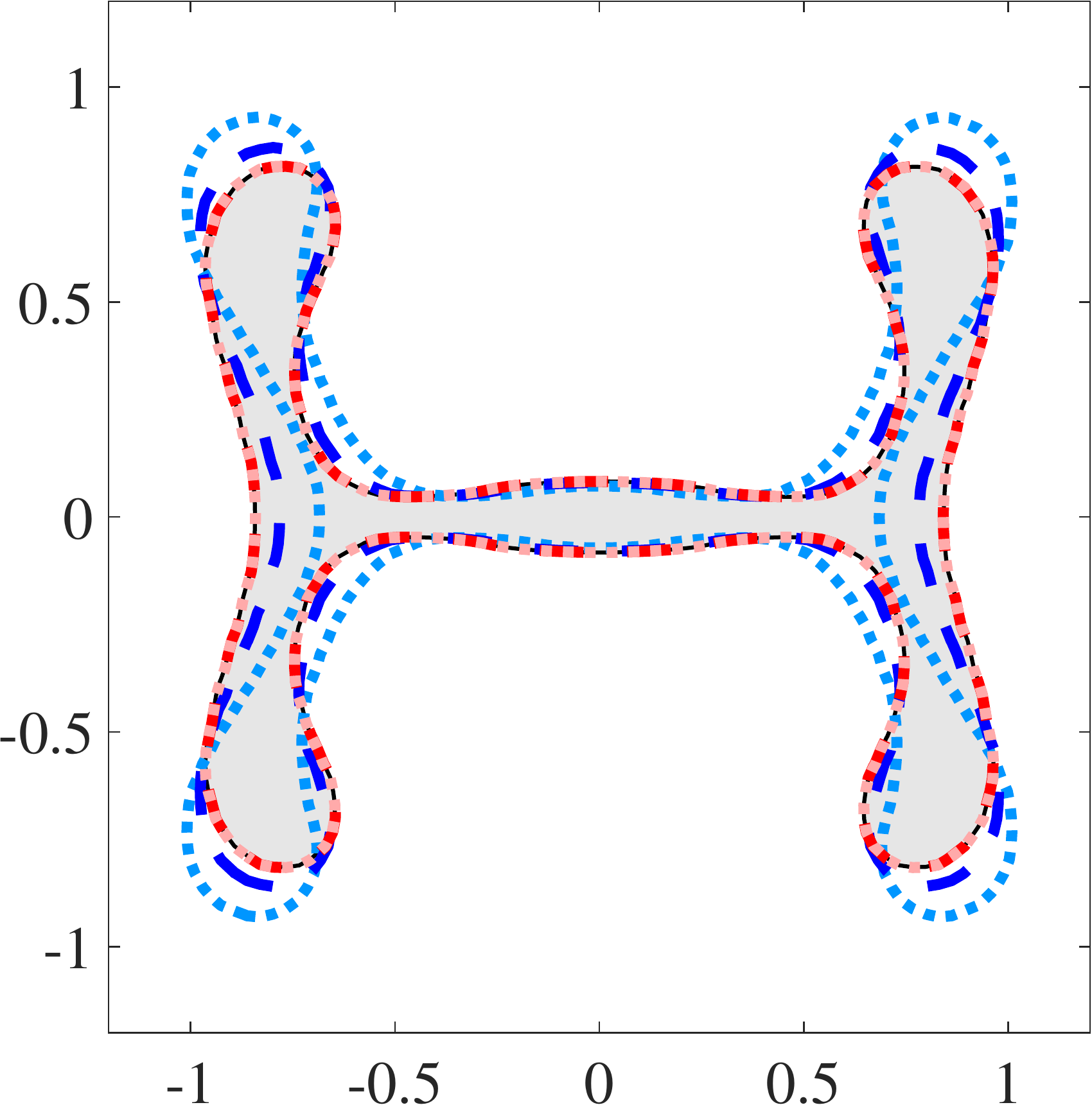}
\caption{$t = 2.0$} %
\end{subfigure}
\begin{subfigure}{0.245\columnwidth} \centering
\includegraphics[trim={0cm 0cm 0cm 0cm},clip,width=1\columnwidth]{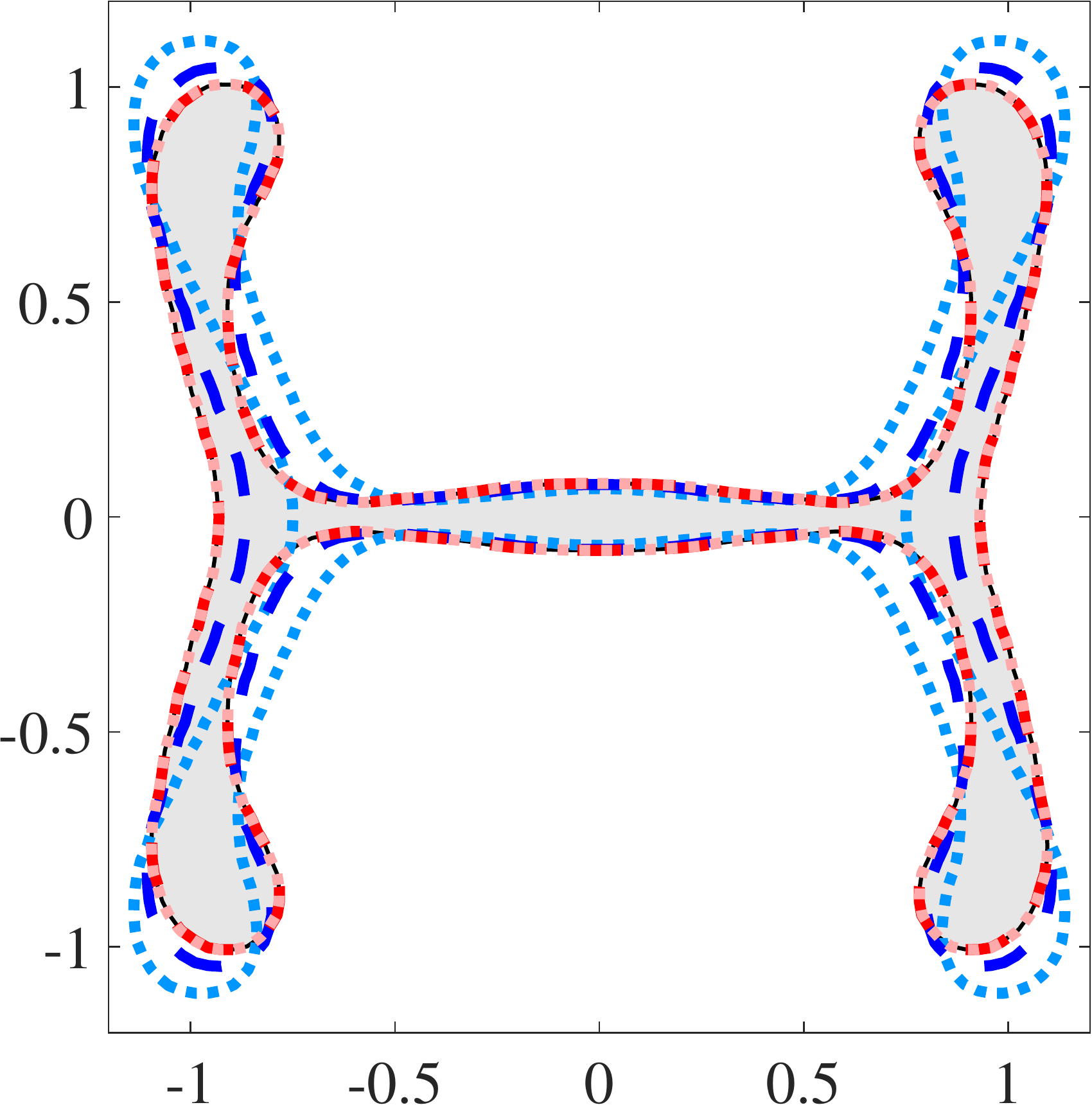}
\caption{$t = 2.5$} %
\end{subfigure}
\caption{Contour traced along the interface separating healthy and tumor tissue ($\phi \approx 0$)  for different uniform tensor product B-spline meshes with either $2^{10}$ or $2^{11}$ elements per side and polynomial degree either $p=2$ or $p=3$ at time 1, 1.5, 2.0, and 2.5.
The reference solution is shown as gray area representing the tumor shape, and it corresponds to the uniform tensor product B-spline mesh with $p = 4$ and $1024$ elements in each direction shown in Fig.~\ref{figure_1_final}.}
\label{figure_3_final}
\end{figure}

\begin{figure}[!ht]\centering
\begin{subfigure}{0.48\columnwidth} \centering
\includegraphics[trim={0cm 0cm 0cm 0cm},clip,width=0.49\columnwidth]{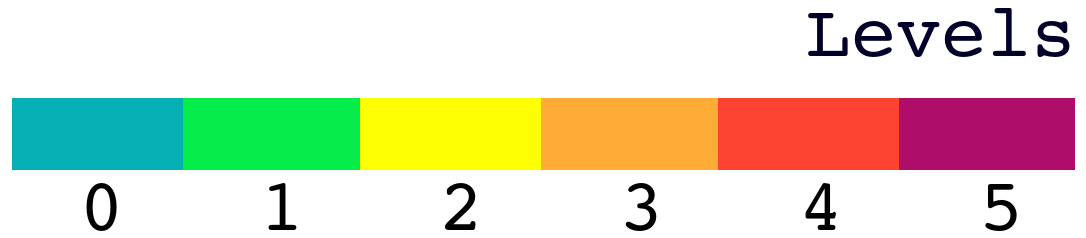}
\end{subfigure}
\begin{subfigure}{0.48\columnwidth} \centering
\includegraphics[trim={0cm 0cm 0cm 0cm},clip,width=0.55\columnwidth]{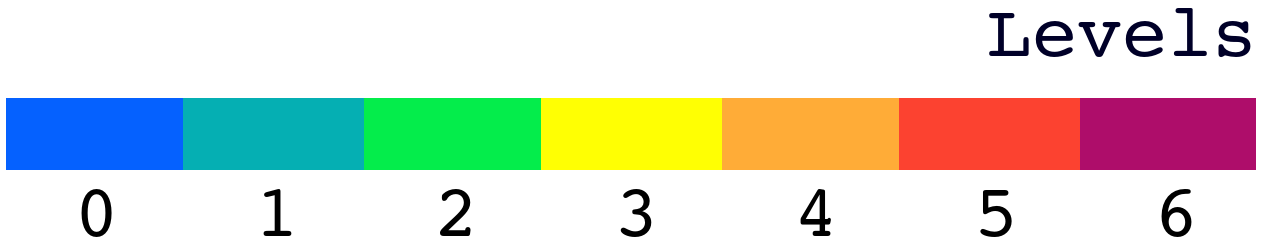}
\end{subfigure}
\begin{subfigure}{0.48\columnwidth} \centering
\includegraphics[trim={0cm 0cm 0cm 0cm},clip,width=0.49\columnwidth]{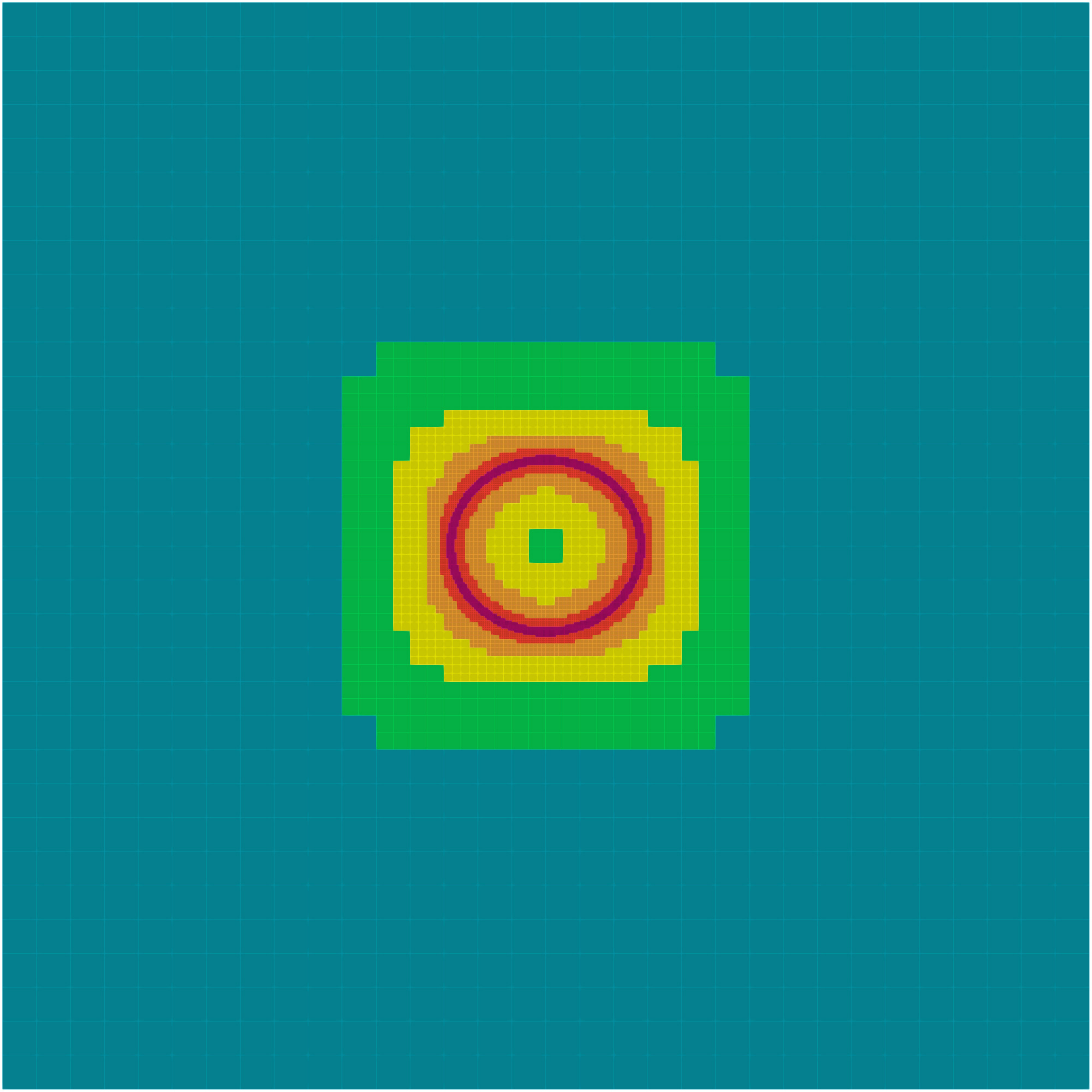}
\includegraphics[trim={0cm 0cm 0cm 0cm},clip,width=0.49\columnwidth]{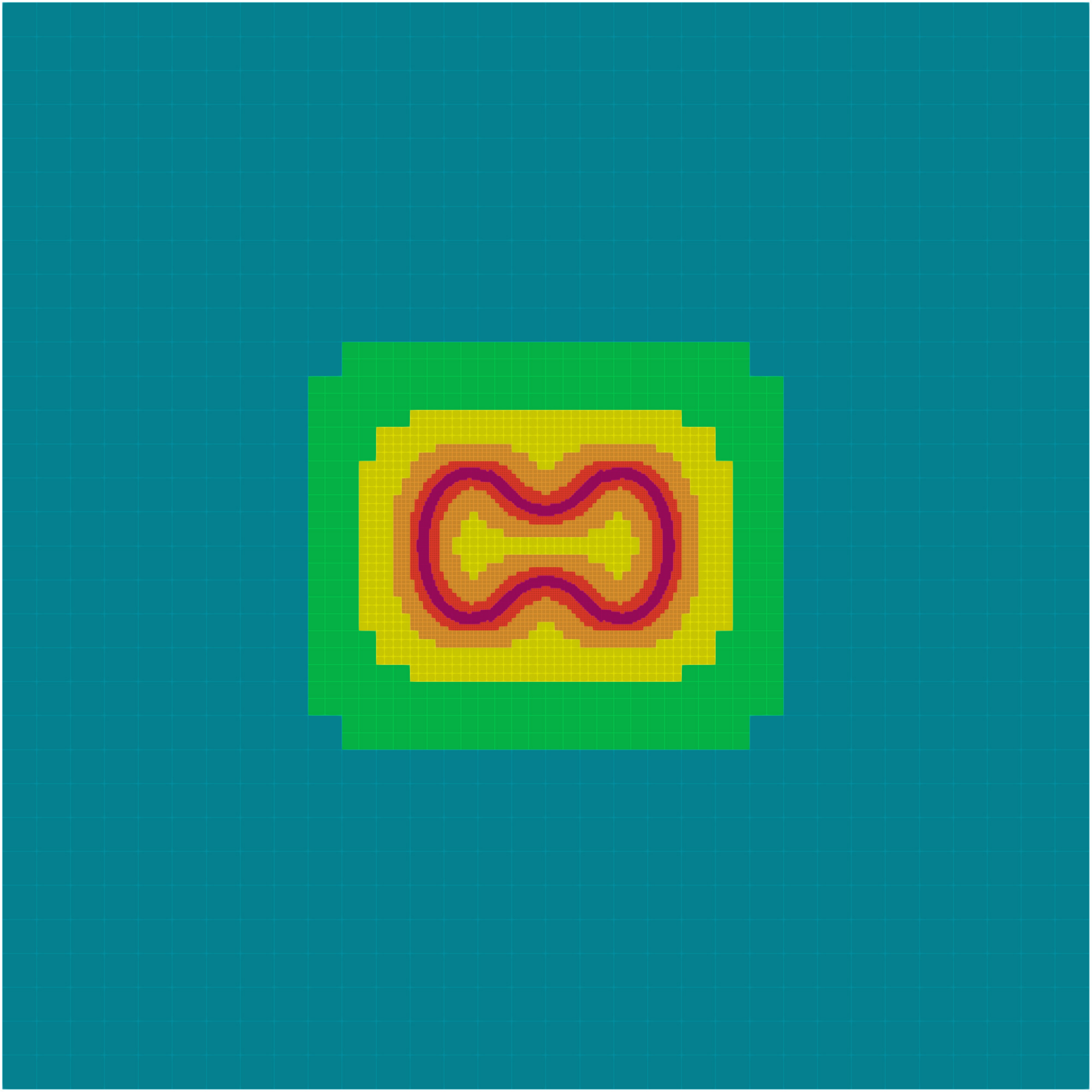}
\caption{$\ell = 5$, $m=2$, $\alpha = 0.1$, $\beta = 0.001$} %
\end{subfigure}
\begin{subfigure}{0.48\columnwidth} \centering
\includegraphics[trim={0cm 0cm 0cm 0cm},clip,width=0.49\columnwidth]{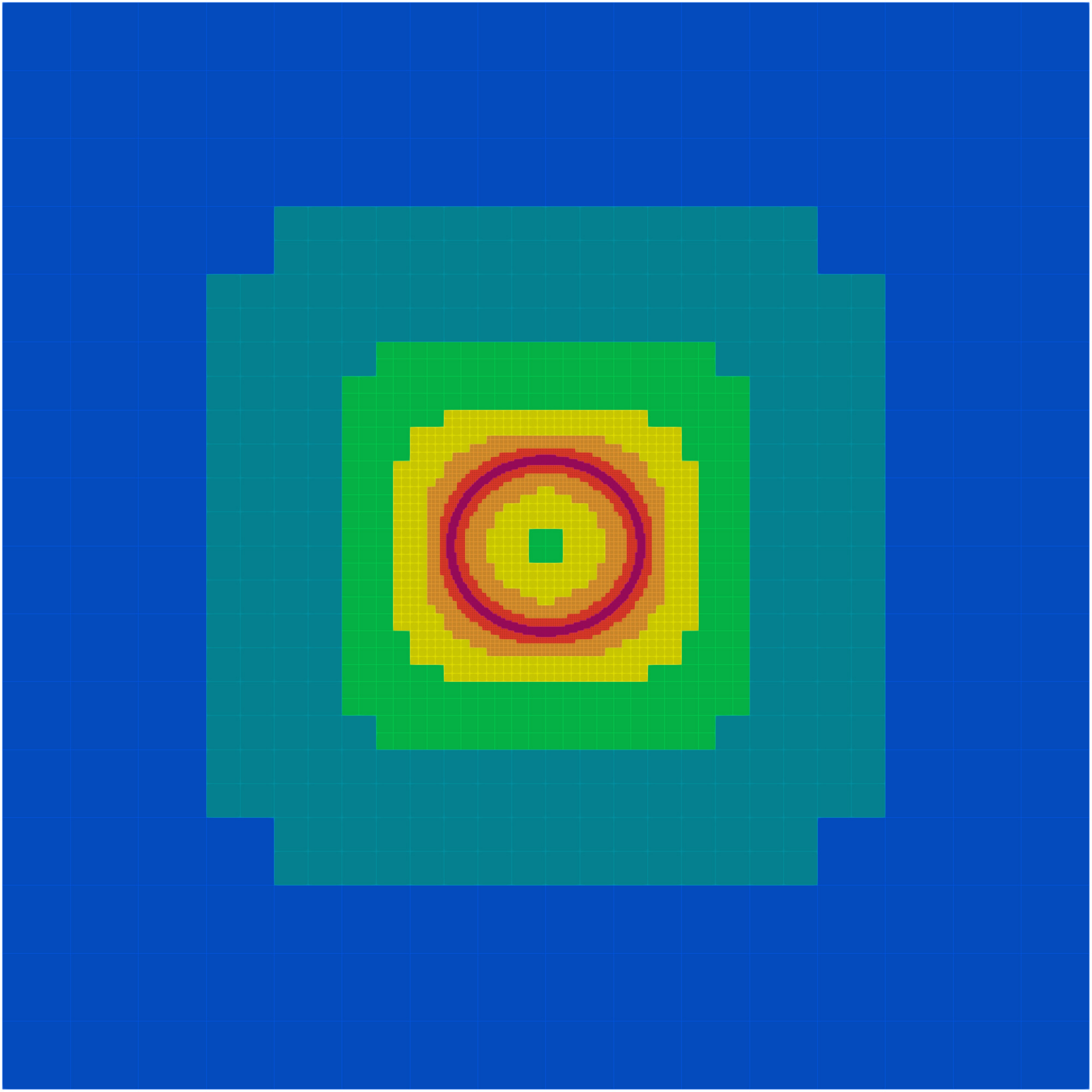}
\includegraphics[trim={0cm 0cm 0cm 0cm},clip,width=0.49\columnwidth]{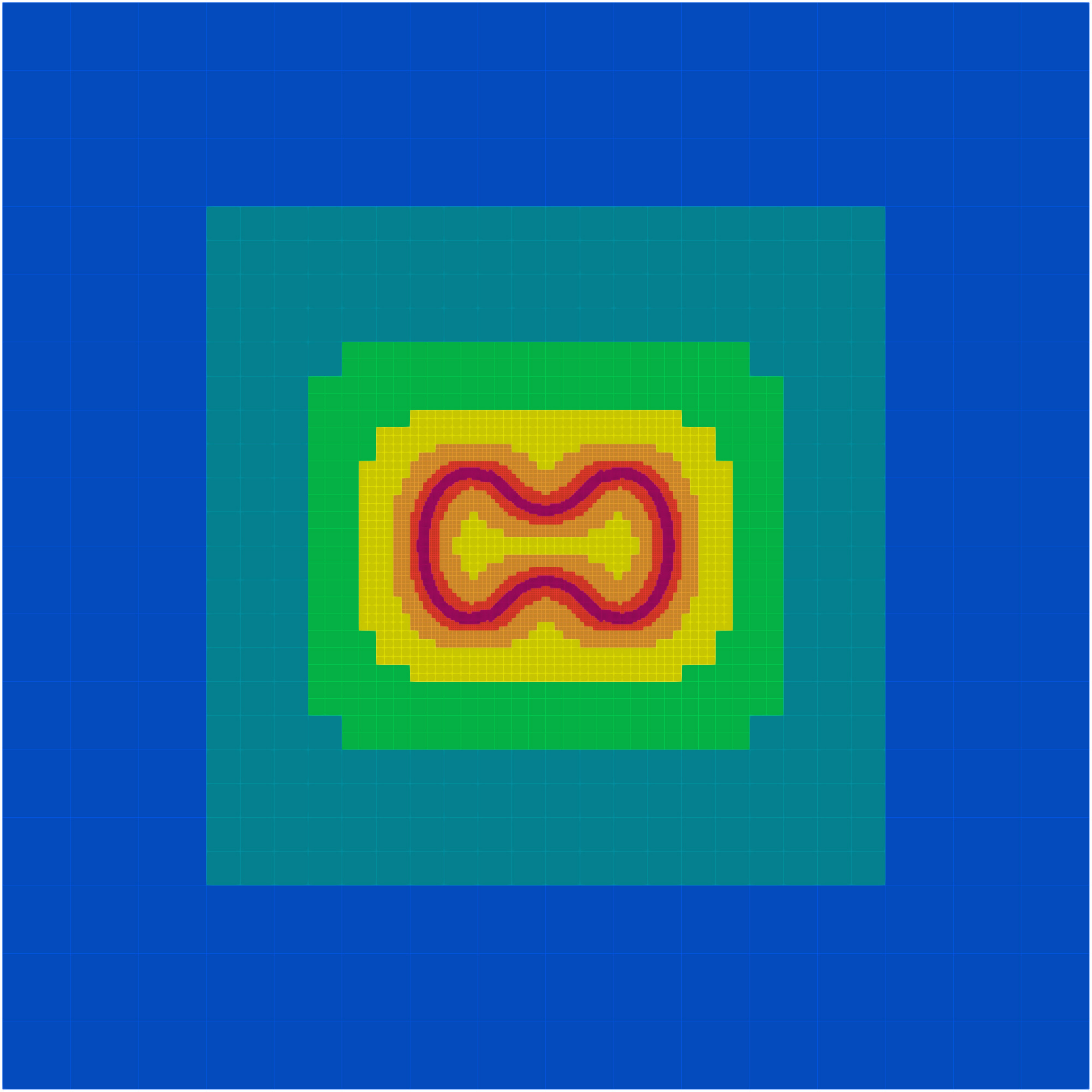}
\caption{$\ell = 6$, $m=2$, $\alpha = 0.1$, $\beta = 0.001$} %
\end{subfigure}
\vspace{0.15cm}
\\
\begin{subfigure}{0.48\columnwidth} \centering
\includegraphics[trim={0cm 0cm 0cm 0cm},clip,width=0.49\columnwidth]{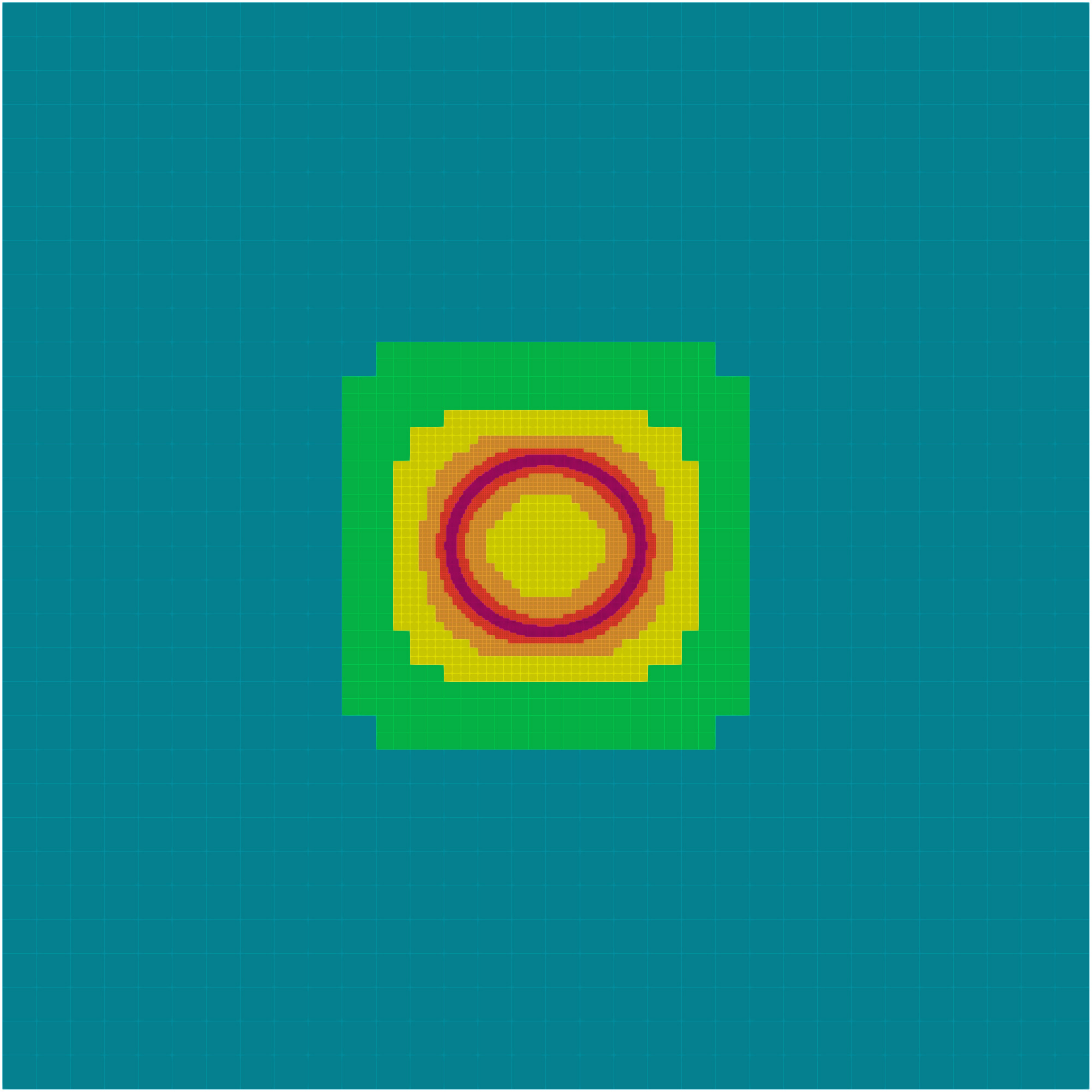}
\includegraphics[trim={0cm 0cm 0cm 0cm},clip,width=0.49\columnwidth]{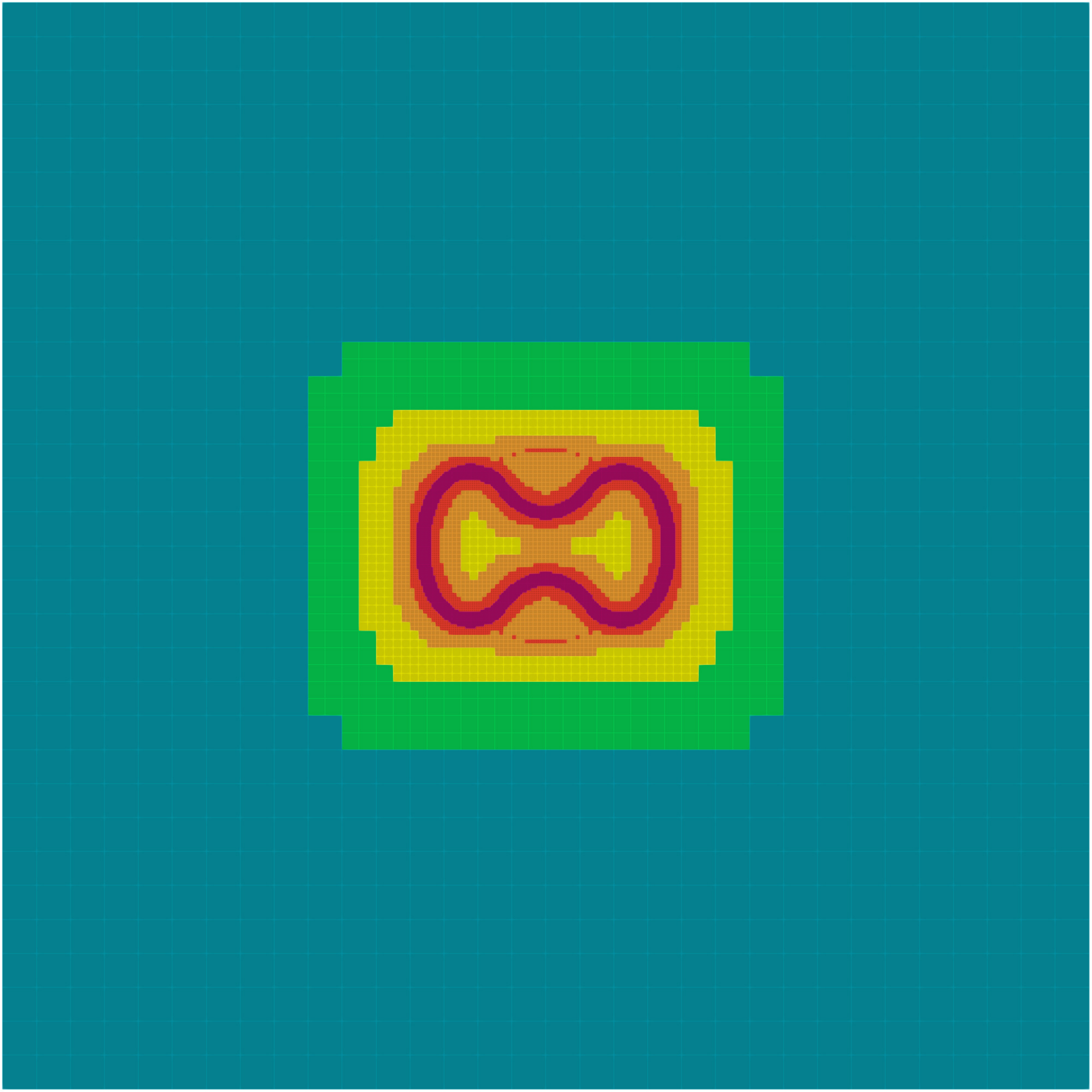}
\caption{$\ell = 5$, $m=2$, $\alpha = 0.01$, $\beta = 0.0001$} %
\end{subfigure}
\begin{subfigure}{0.48\columnwidth} \centering
\includegraphics[trim={0cm 0cm 0cm 0cm},clip,width=0.49\columnwidth]{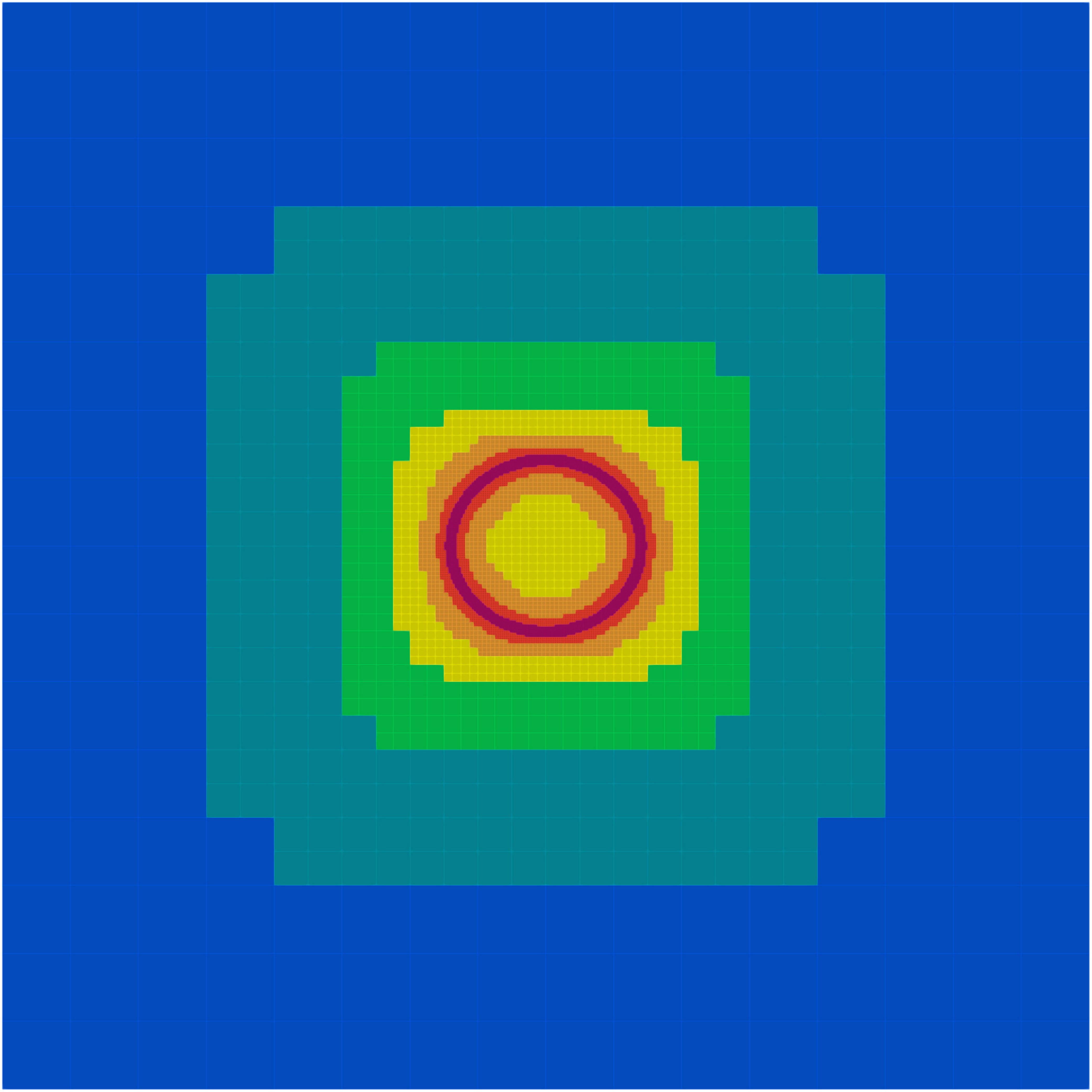}
\includegraphics[trim={0cm 0cm 0cm 0cm},clip,width=0.49\columnwidth]{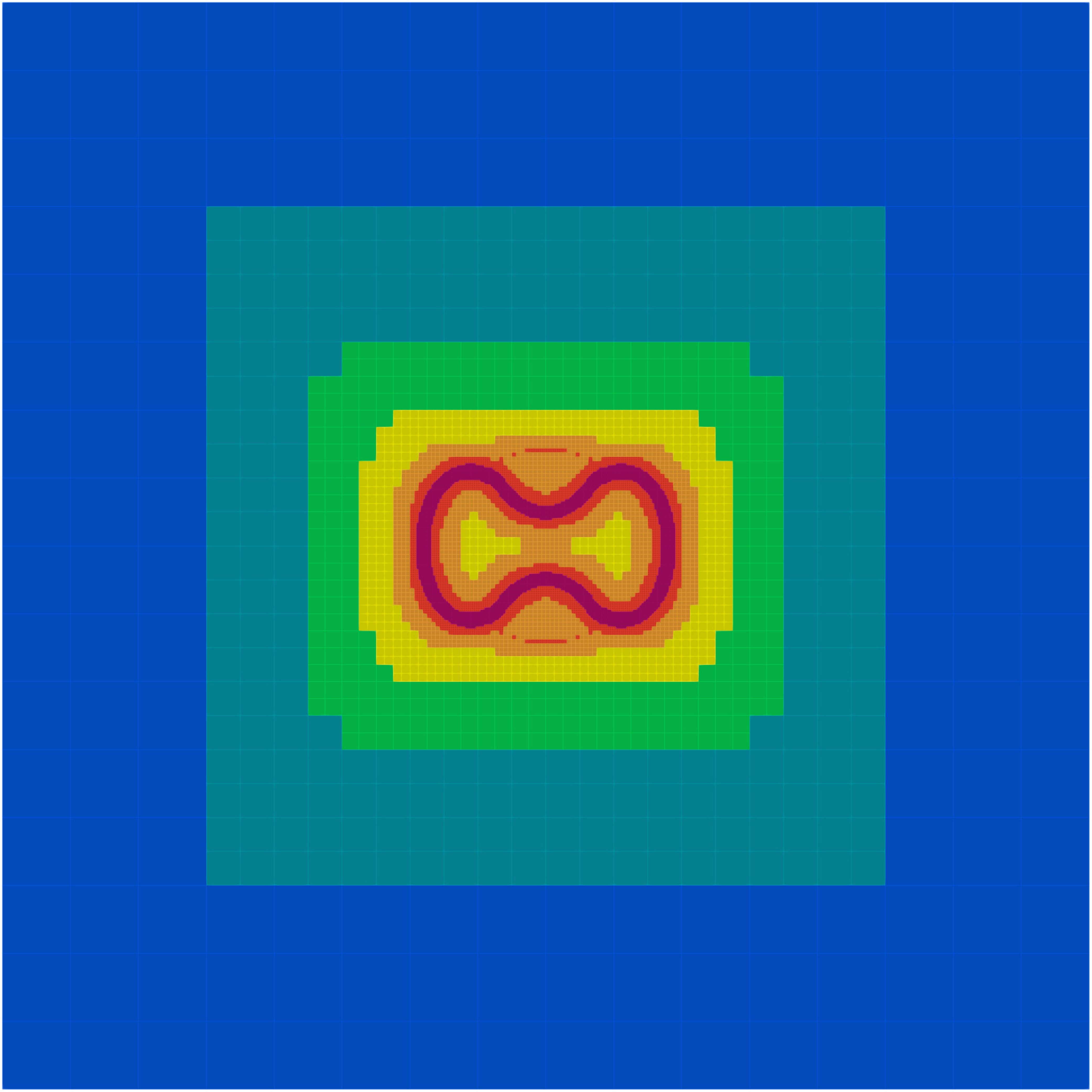}
\caption{$\ell = 6$, $m=2$, $\alpha = 0.01$, $\beta = 0.0001$} %
\end{subfigure}
\vspace{0.15cm}
\\
\begin{subfigure}{0.48\columnwidth} \centering
\includegraphics[trim={0cm 0cm 0cm 0cm},clip,width=0.49\columnwidth]{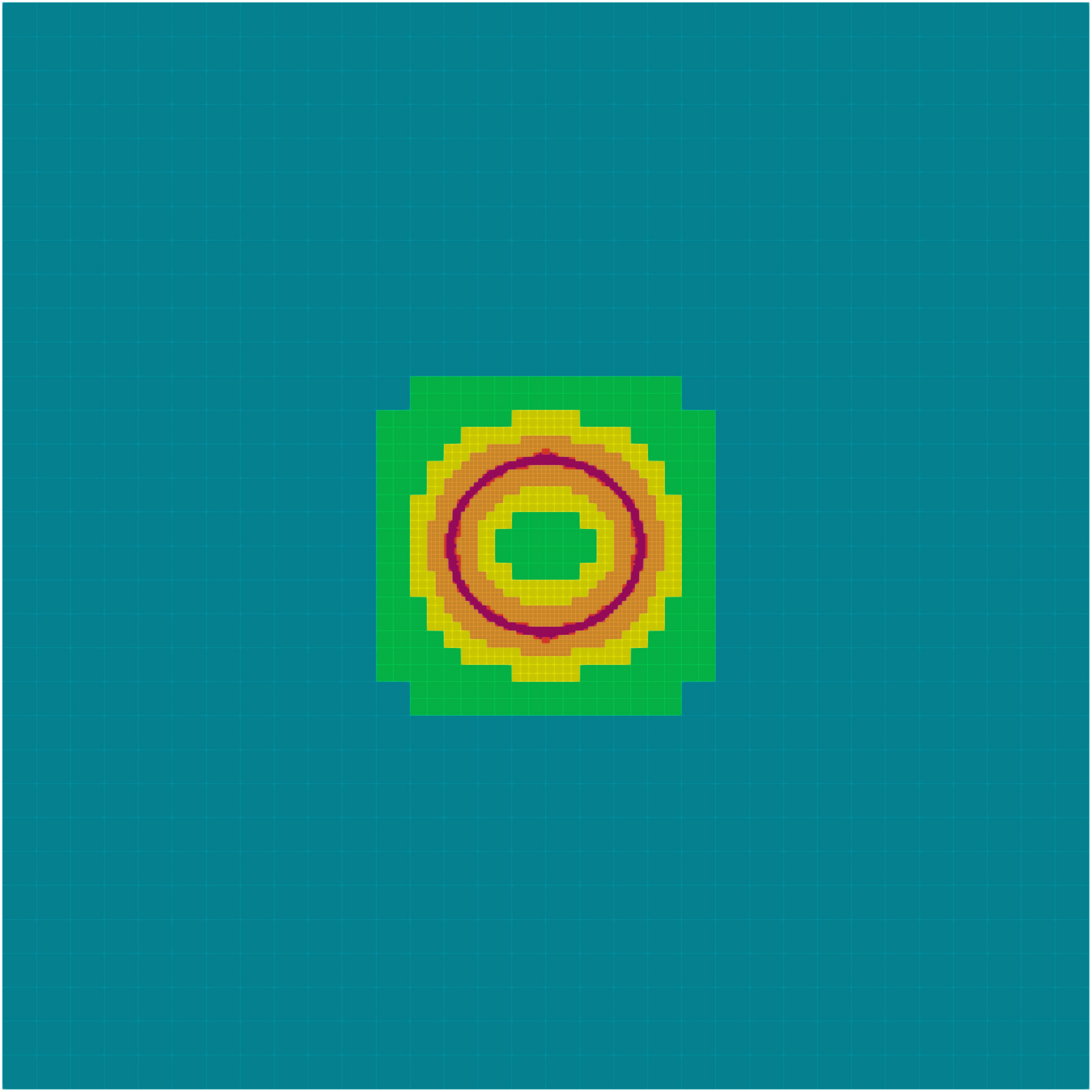}
\includegraphics[trim={0cm 0cm 0cm 0cm},clip,width=0.49\columnwidth]{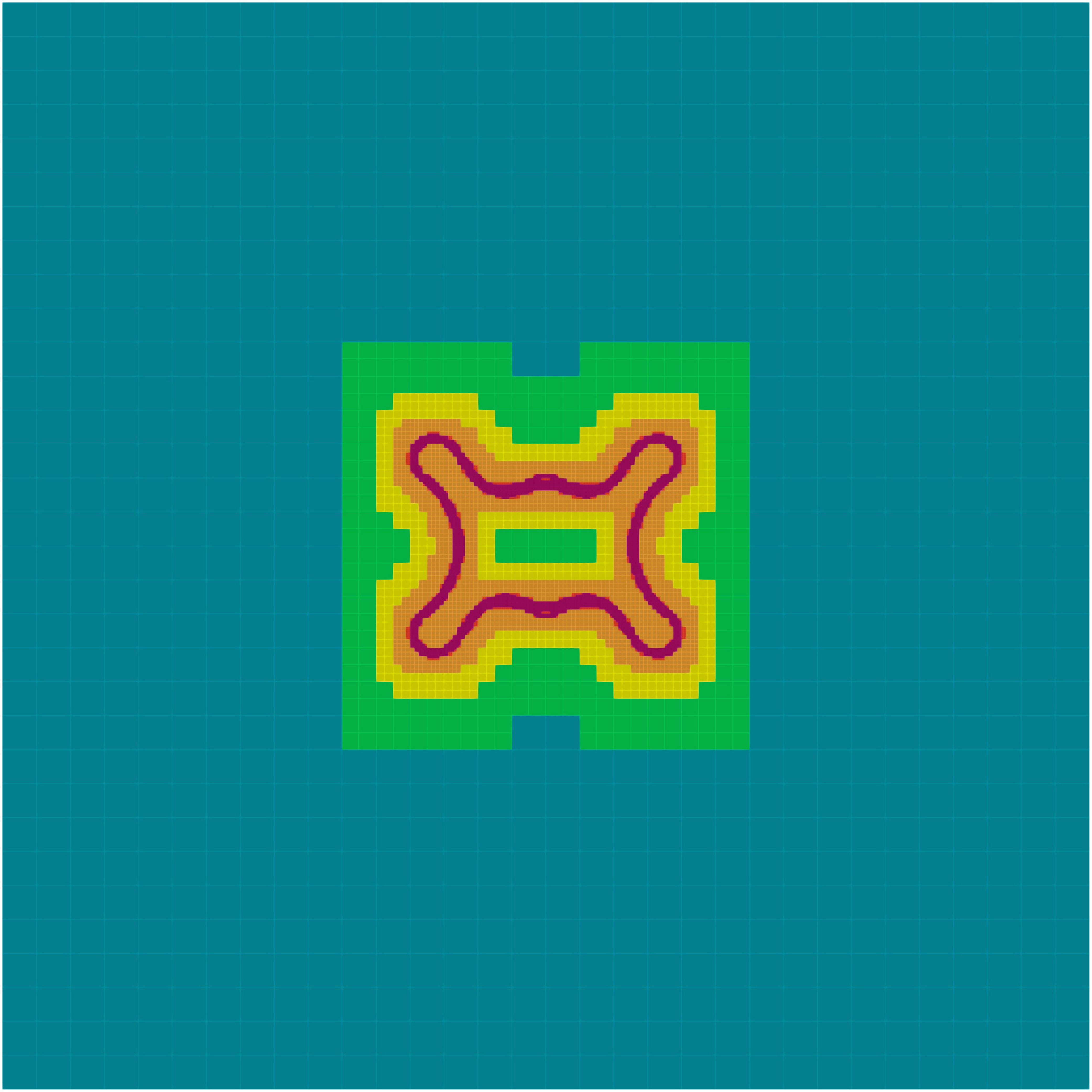}
\caption{$\ell = 5$, $m=3$, $\alpha = 0.1$, $\beta = 0.001$} %
\end{subfigure}
\begin{subfigure}{0.48\columnwidth} \centering
\includegraphics[trim={0cm 0cm 0cm 0cm},clip,width=0.49\columnwidth]{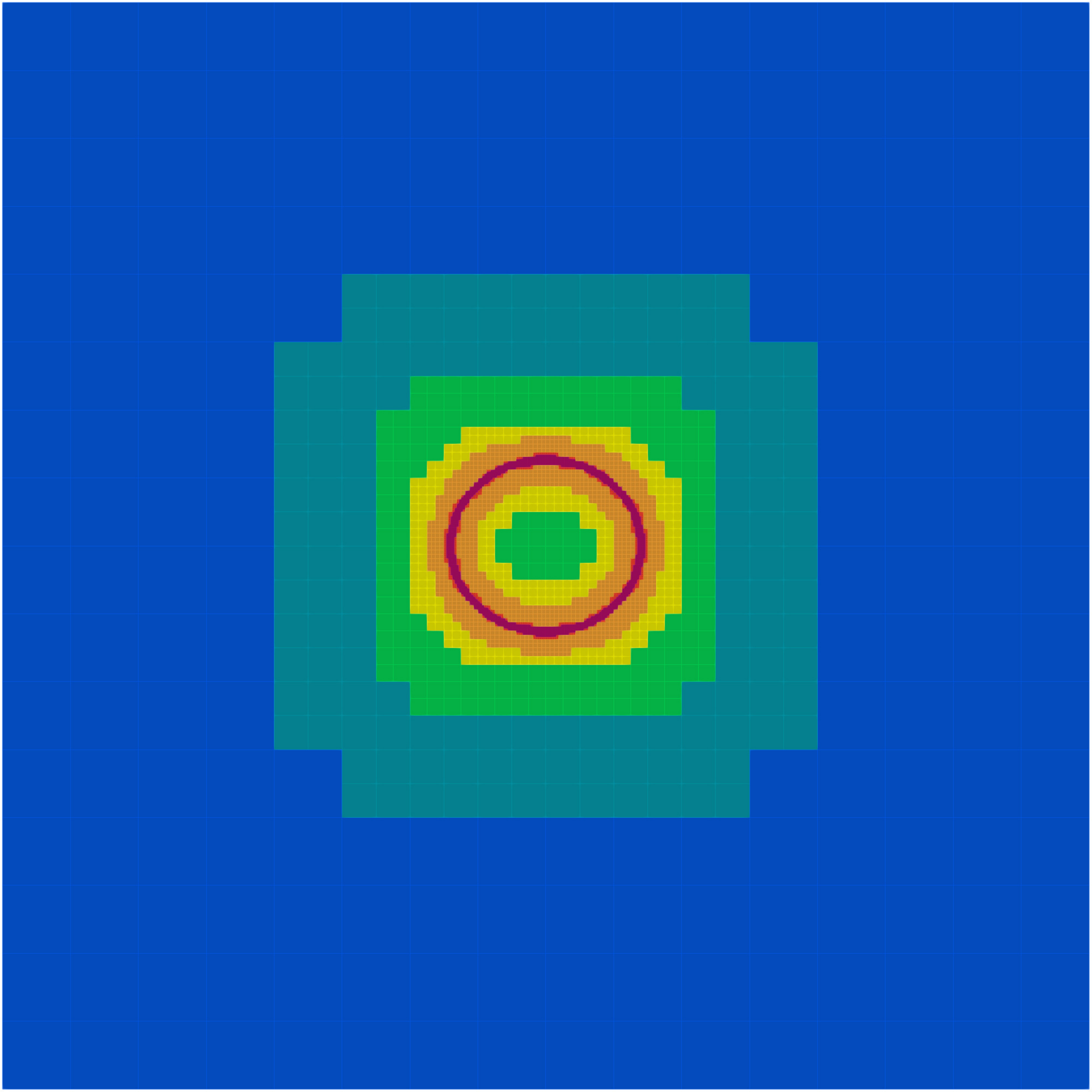}
\includegraphics[trim={0cm 0cm 0cm 0cm},clip,width=0.49\columnwidth]{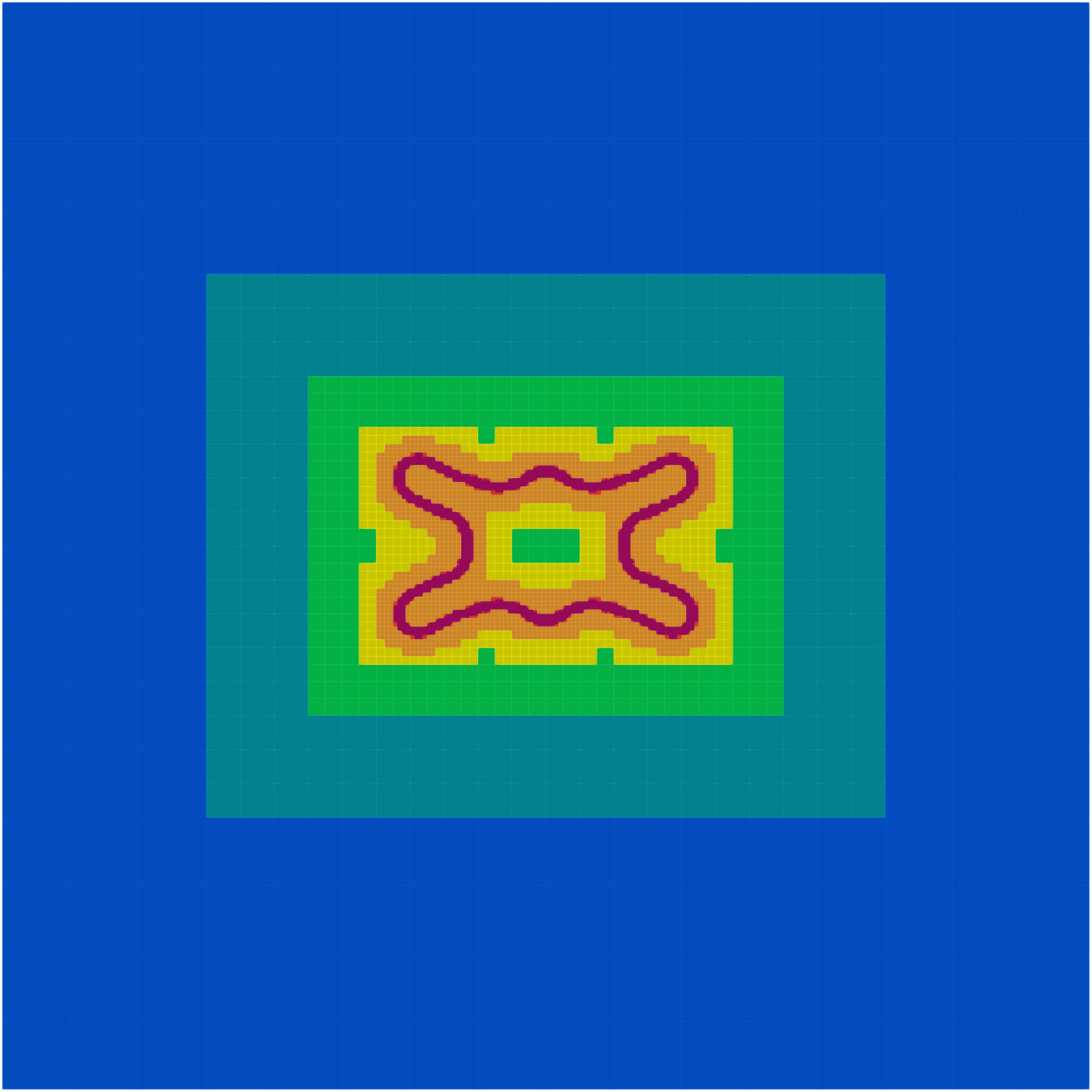}
\caption{$\ell = 6$, $m=3$, $\alpha = 0.1$, $\beta = 0.001$} %
\end{subfigure}
\vspace{0.15cm}
\\
\begin{subfigure}{0.48\columnwidth} \centering
\includegraphics[trim={0cm 0cm 0cm 0cm},clip,width=0.49\columnwidth]{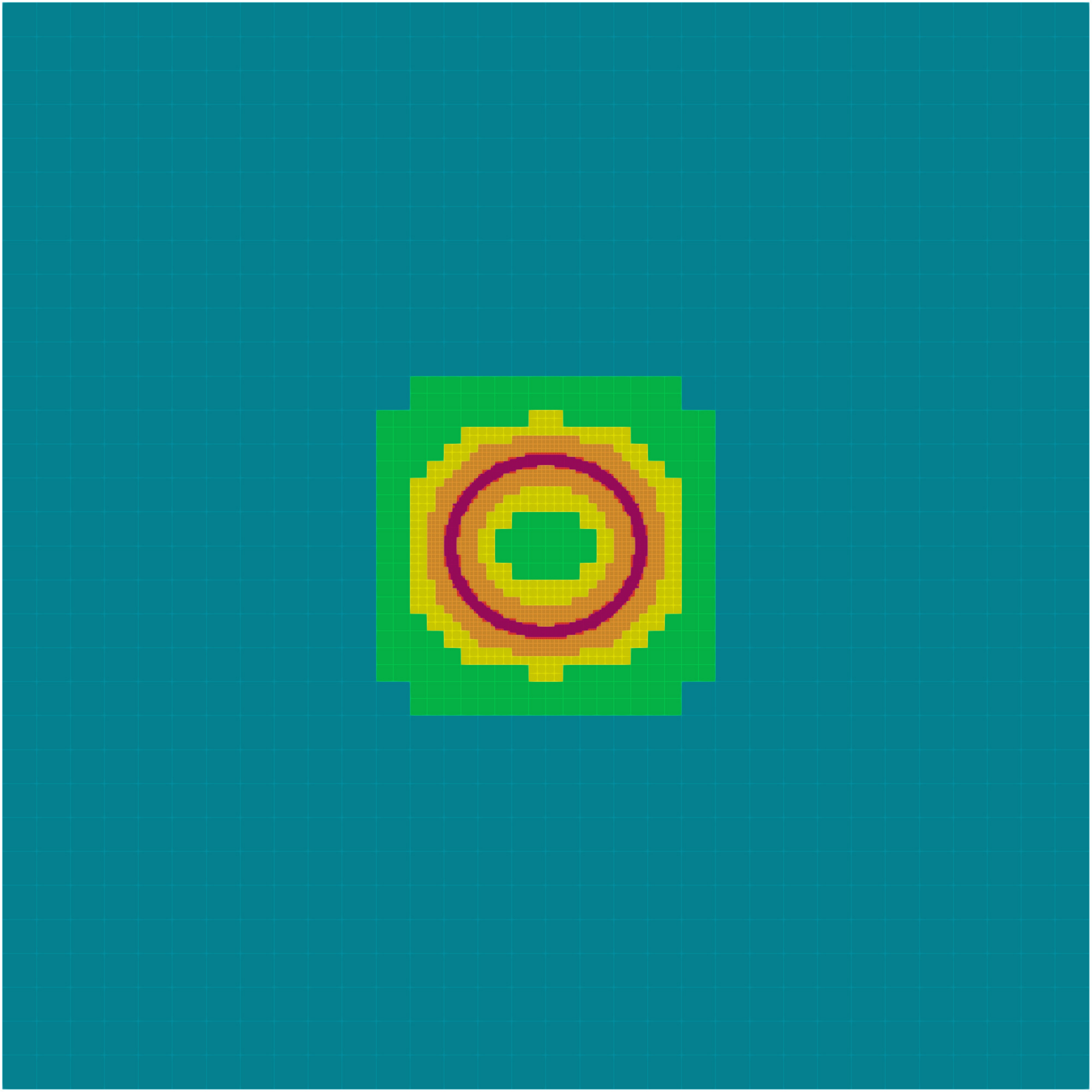}
\includegraphics[trim={0cm 0cm 0cm 0cm},clip,width=0.49\columnwidth]{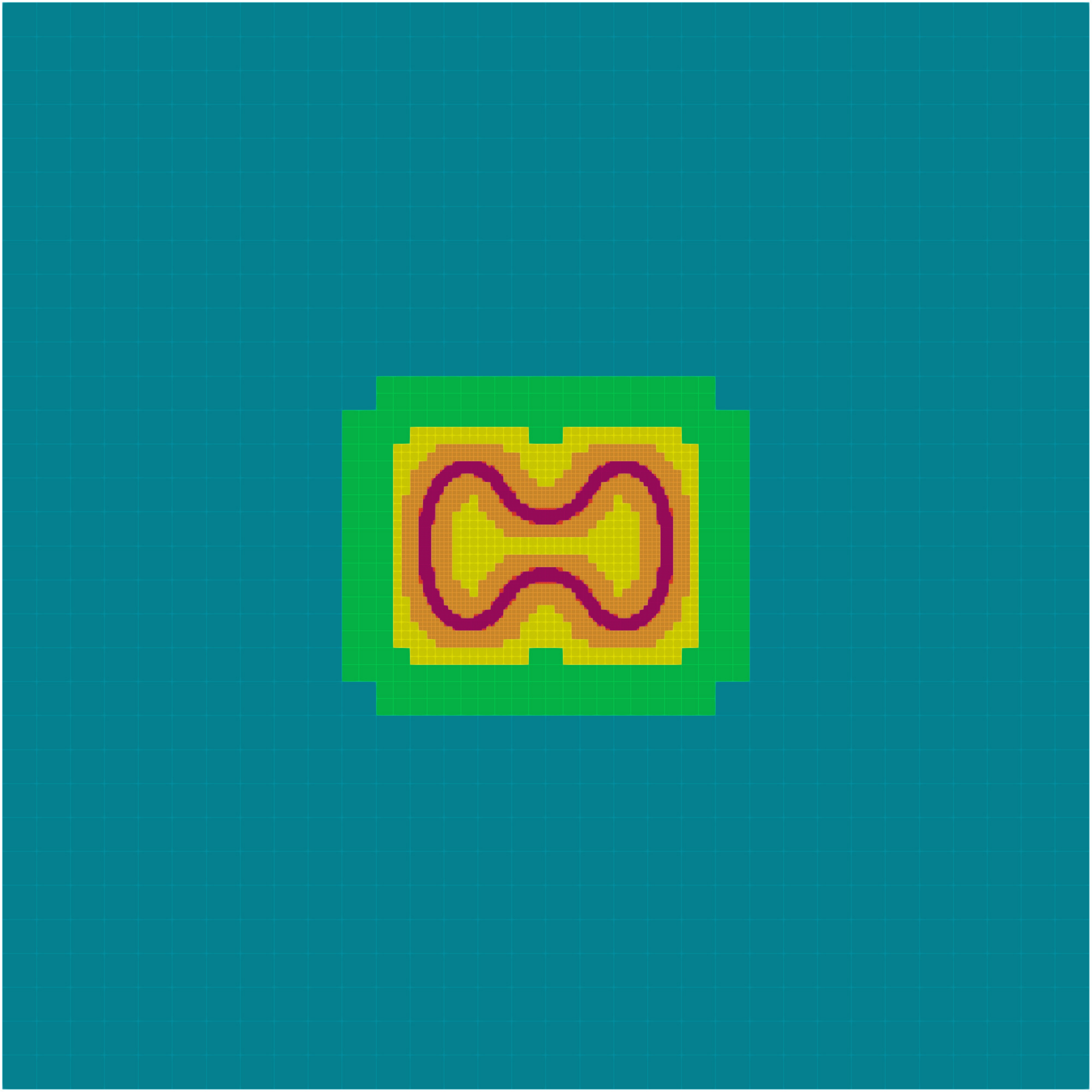}
\caption{$\ell = 5$, $m=3$, $\alpha = 0.01$, $\beta = 0.0001$} %
\end{subfigure}
\begin{subfigure}{0.48\columnwidth} \centering
\includegraphics[trim={0cm 0cm 0cm 0cm},clip,width=0.49\columnwidth]{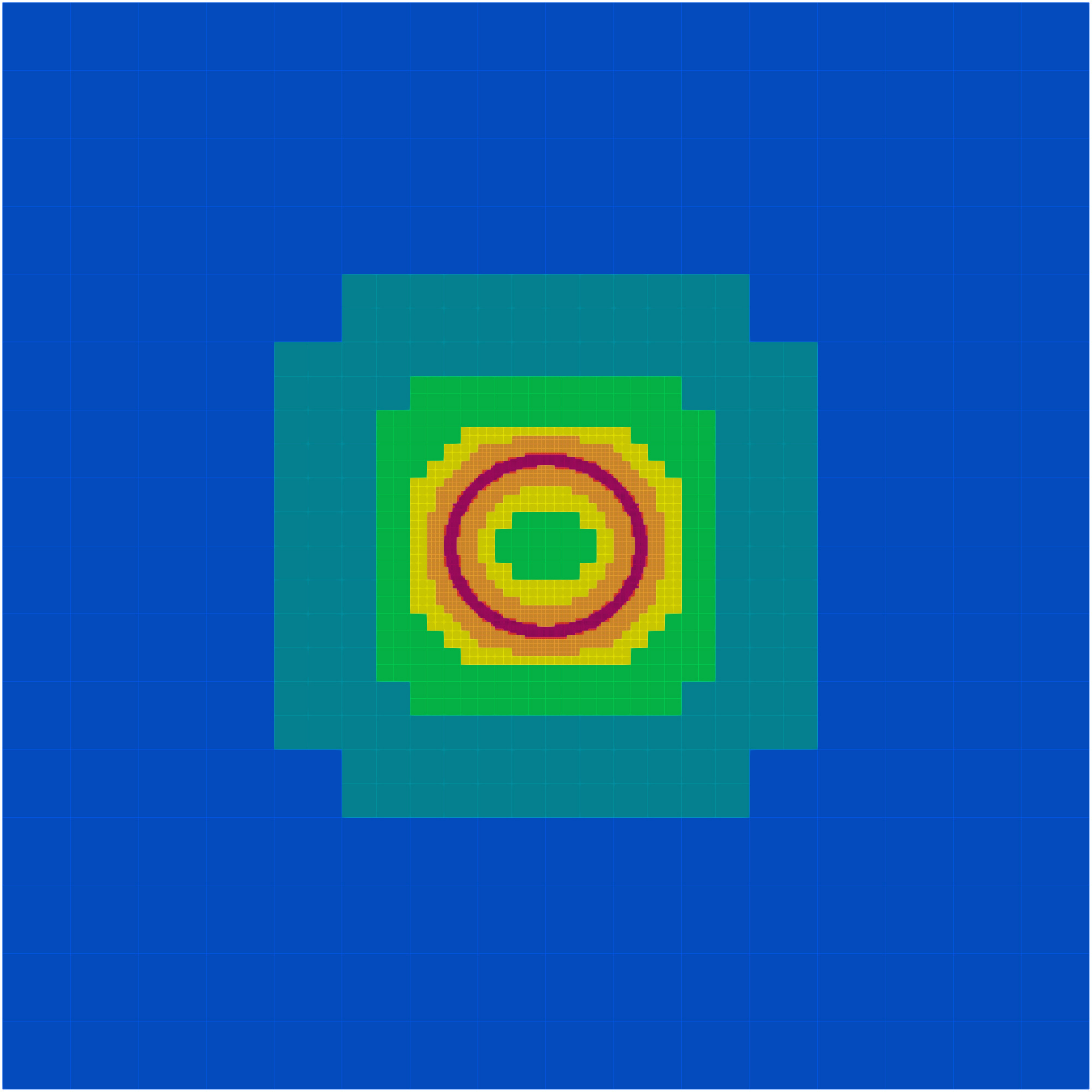}
\includegraphics[trim={0cm 0cm 0cm 0cm},clip,width=0.49\columnwidth]{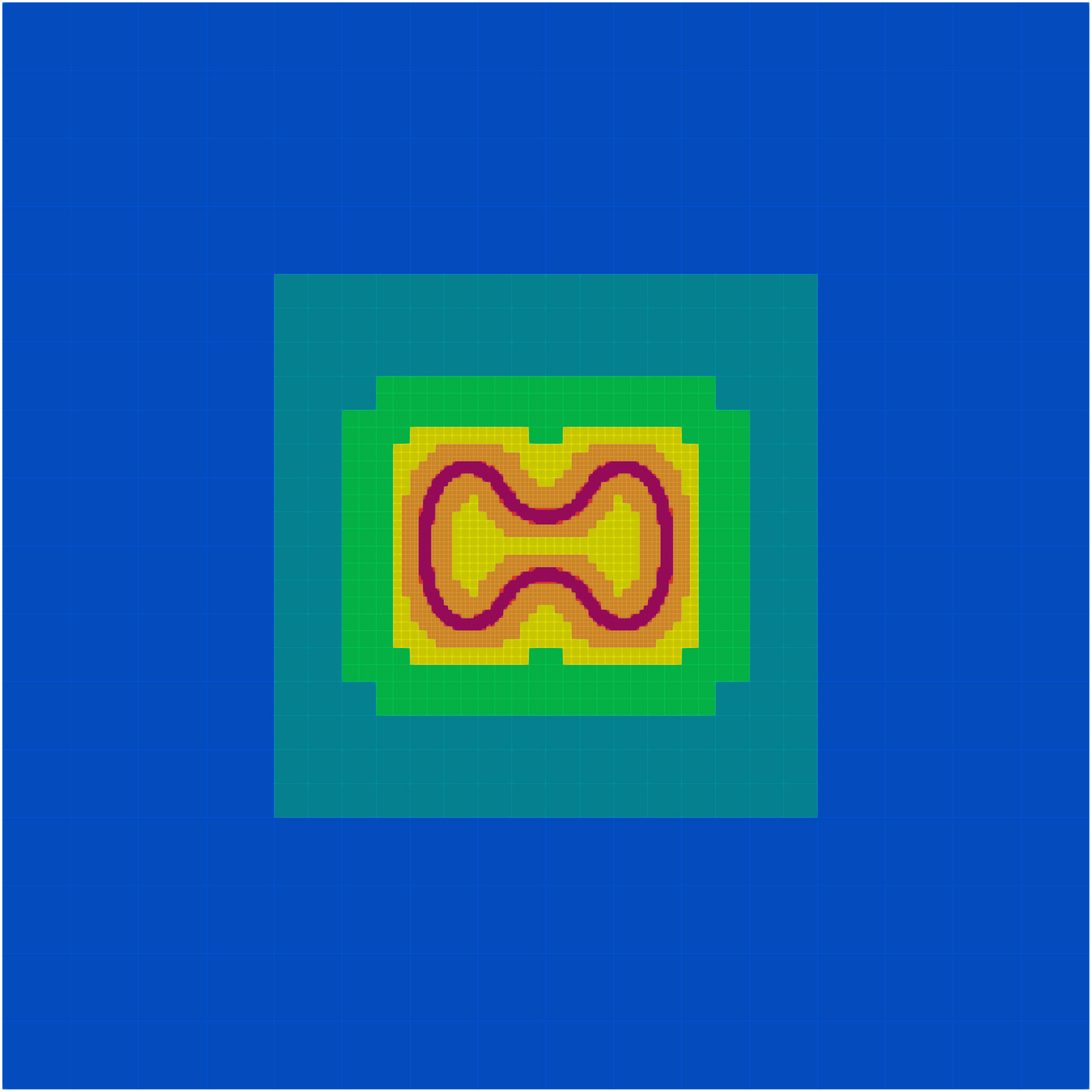}
\caption{$\ell = 6$, $m=3$, $\alpha = 0.01$, $\beta = 0.0001$} %
\end{subfigure}
\caption{Adaptive THB-spline mesh configurations of degree $3$ with varying hierarchical levels, admissibility class, and refinement and coarsening parameters at $t=0$ and $t=1$ (left and right figure in each panel, respectively). The finest level for all the meshes is fixed at $2^{10} \times 2^{10}$ (i.e., level 5 in panels (a), (c), (e), and (g) as well as level 6 in panels (b), (d), (f) and (h)). These configurations represent the candidate hierarchical meshes evaluated to select the most suited one for the simulation of tumor growth with the phase-field model.}
\label{figure_4_final}
\end{figure}

\section{Numerical studies}\label{Numerical_studies}
This section presents a series of numerical studies to assess the performance and capability of the CH-based tumor growth model. 
We begin with standard benchmark problems on a 2D square and a 3D cubic domain to study the established tumor growth morphologies and to explore the influence of key parameters on the resulting growth dynamics. 
The framework is then extended to an organ-scale, patient-specific breast geometry reconstructed from MRI data to examine model behavior in a realistic anatomical setting.
To quantify numerical accuracy, the relative $L^2(\Omega)$ error and the relative $H^1(\Omega)$ semi-norm of the error are evaluated as follows: 
\begin{equation}
    \| e \|_{L^2}^{\mathrm{rel}} :=
\frac{\left( \int_\Omega (\phi^{\mathrm{ref}} - \phi^{\mathrm{h}})^2 \, d\Omega \right)^{1/2}}
     {\left( \int_\Omega (\phi^{\mathrm{ref}})^2 \, d\Omega \right)^{1/2}} \\ \text{ and }
     | e |_{H^1}^{\mathrm{rel}} :=
\frac{\left( \int_\Omega |\nabla(\phi^{\mathrm{ref}} - \phi^{\mathrm{h}})|^2 \, d\Omega \right)^{1/2}}
     {\left( \int_\Omega |\nabla \phi^{\mathrm{ref}}|^2 \, d\Omega \right)^{1/2}},
\end{equation}
where $\phi^{\mathrm{ref}}$ denotes a reference solution and $\phi^{\mathrm{h}}$ is numerical solution obtained for a chosen mesh configuration. 
The ratio of the number of DOFs of a given mesh to that of the reference mesh is defined as $n_{\text{dof}}^{\text{rel}} := \frac{\#\text{DOFs}^\text{h}}{\#\text{DOFs}^{\text{ref}}}$.
The relative computational runtime is defined as $n_{\text{time}}^{\text{rel}} := \frac{T_{\mathrm{run}}^{h}} {T_{\mathrm{run}}^{\mathrm{ref}}} $, where $T_{\mathrm{run}}^{h}$ and $T_{\mathrm{run}}^{\mathrm{ref}}$ denote the runtimes of the given and reference meshes, respectively.
All examples presented here use the same degree spline basis in every direction, i.e., $\mathbf{p} = (p_1, \ldots, p_d)$ with $p_1 = \cdots = p_d = p$.

\subsection{Tumor growth in a square domain}\label{TG2D}

The present section shows the spatiotemporal tumor growth in a square domain defined by $\Omega = [-3, 3]^2$.
The model parameters used are given as follows \cite{Ebenbeck2021}: $
    \lambda = 0.0002,\: %
    E = 1, \:
    M = 1, \:
    \kappa = 2, \:
    \mathcal{P} = 0.1, \:
    \mathcal{A} = 0,
     D = 1, \:
     \mathcal{C} = 2, \:
     \mathcal{X} = 0.02, \:
     \mathcal{X}_\phi = 5,$ and $
     \sigma_B = 1.$
These parameter values are adopted throughout the subsequent numerical examples unless stated otherwise.
In the first numerical example, the phase-field is initialized through the following relation: 
\begin{equation} \label{initial_cindition}
\phi_0(\mathbf{x}) =
\begin{cases}
1 & r(\mathbf{x}) \le -\frac{1}{2}\pi\epsilon, \\
-\sin\!\left(\frac{r(\mathbf{x})}{\epsilon}\right) & |r(\mathbf{x})| < \frac{1}{2}\pi\epsilon, \\
-1 & r(\mathbf{x}) \ge \frac{1}{2}\pi\epsilon,
\end{cases}
\end{equation}
where $r(\mathbf{x}) = |\mathbf{x}| - \left[\frac{1}{2} + \frac{1}{40}\cos(2\theta) \right]$ and $\epsilon=0.02$.
The resulting initial condition $\phi_0(\bx)$ yields an elliptical profile centered at the origin of the domain, as shown in Fig.~\ref{fig1_initial_condition} (top row). 
The problem is solved over $t \in \left[ 0, 2.5 \right]$ with a time step $\Delta t = 0.001$.

\begin{figure}[!t]\centering
\begin{subfigure}{1\columnwidth} \centering
\includegraphics[trim={0cm -0.25cm 0cm 0cm},clip,width=0.96\columnwidth]{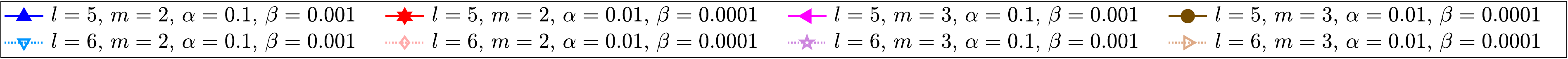}
\end{subfigure}
\\
\begin{subfigure}{0.49\columnwidth} \centering
\includegraphics[trim={0cm 0cm 0cm 0cm},clip,width=0.9\columnwidth]{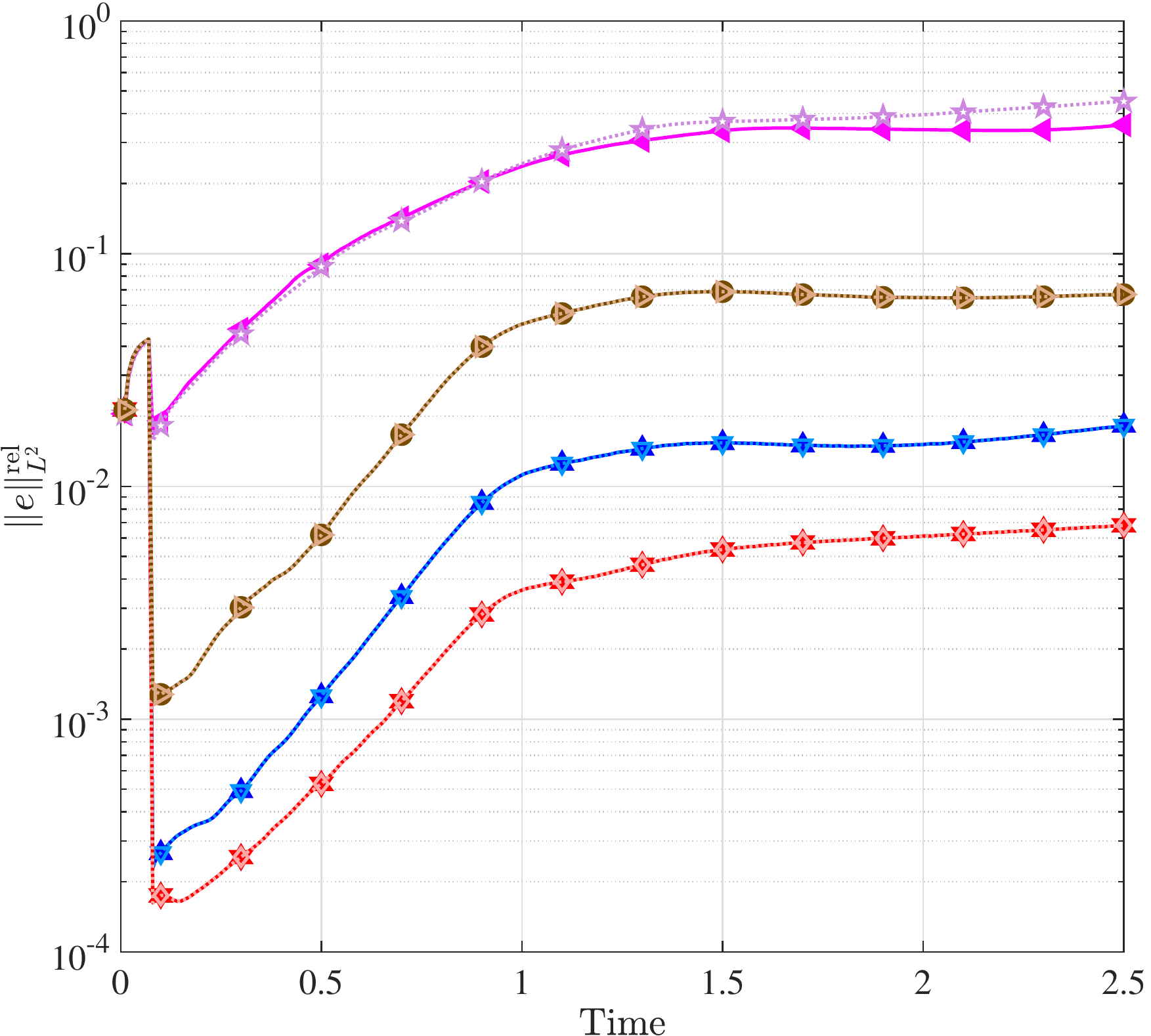}
\caption{$\| e \|_{L^2}^{\mathrm{rel}}$} \label{figure_5_l2err}
\end{subfigure}
\begin{subfigure}{0.49\columnwidth} \centering
\includegraphics[trim={0cm 0cm 0cm 0cm},clip,width=0.9\columnwidth]{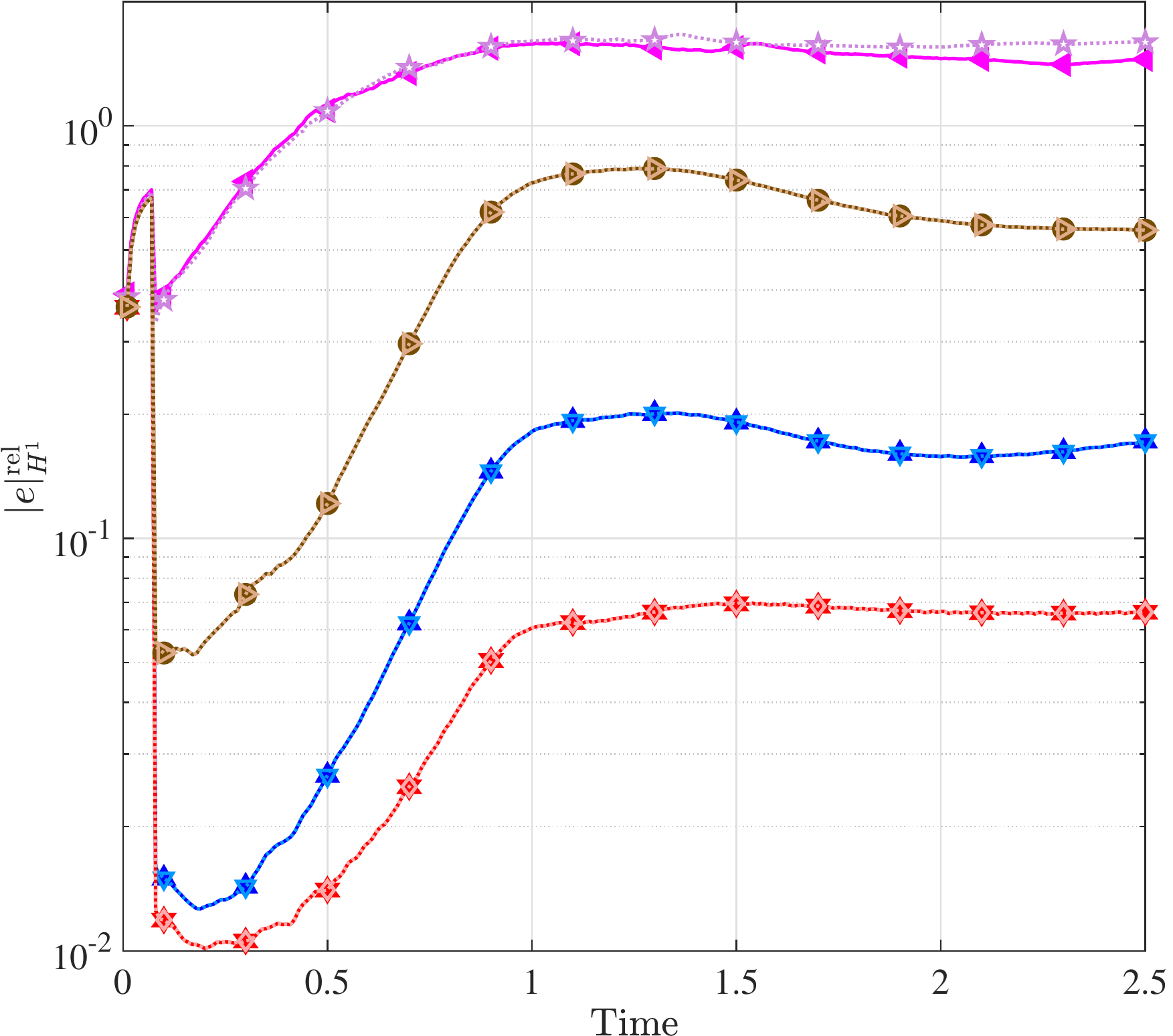}
\caption{$| e |_{H^1}^{\mathrm{rel}}$} \label{figure_5_h1err}
\end{subfigure}
\caption{Relative $\| e \|_{L^2}^{\mathrm{rel}}$ and $| e |_{H^1}^{\mathrm{rel}}$ errors over the time period $\left[0, 2.5\right]$ for different adaptive THB mesh configurations. The reference solution corresponds to the uniform tensor product B-spline mesh with $p = 4$ and $1024$ elements in each direction shown in Fig.~\ref{figure_1_final}}.
\label{figure_5_final}
\end{figure}

\begin{figure}[!t]\centering
\begin{subfigure}{1\columnwidth} \centering
\includegraphics[trim={0cm -0.25cm 0cm 0cm},clip,width=0.96\columnwidth]{Images/Problem1_2DSq_UniRef/Cub1024_THB_legends-eps-converted-to.pdf}
\end{subfigure}
\\
\begin{subfigure}{0.39\columnwidth} \centering
\includegraphics[trim={0cm 0cm 0cm 0cm},clip,width=0.98\columnwidth]{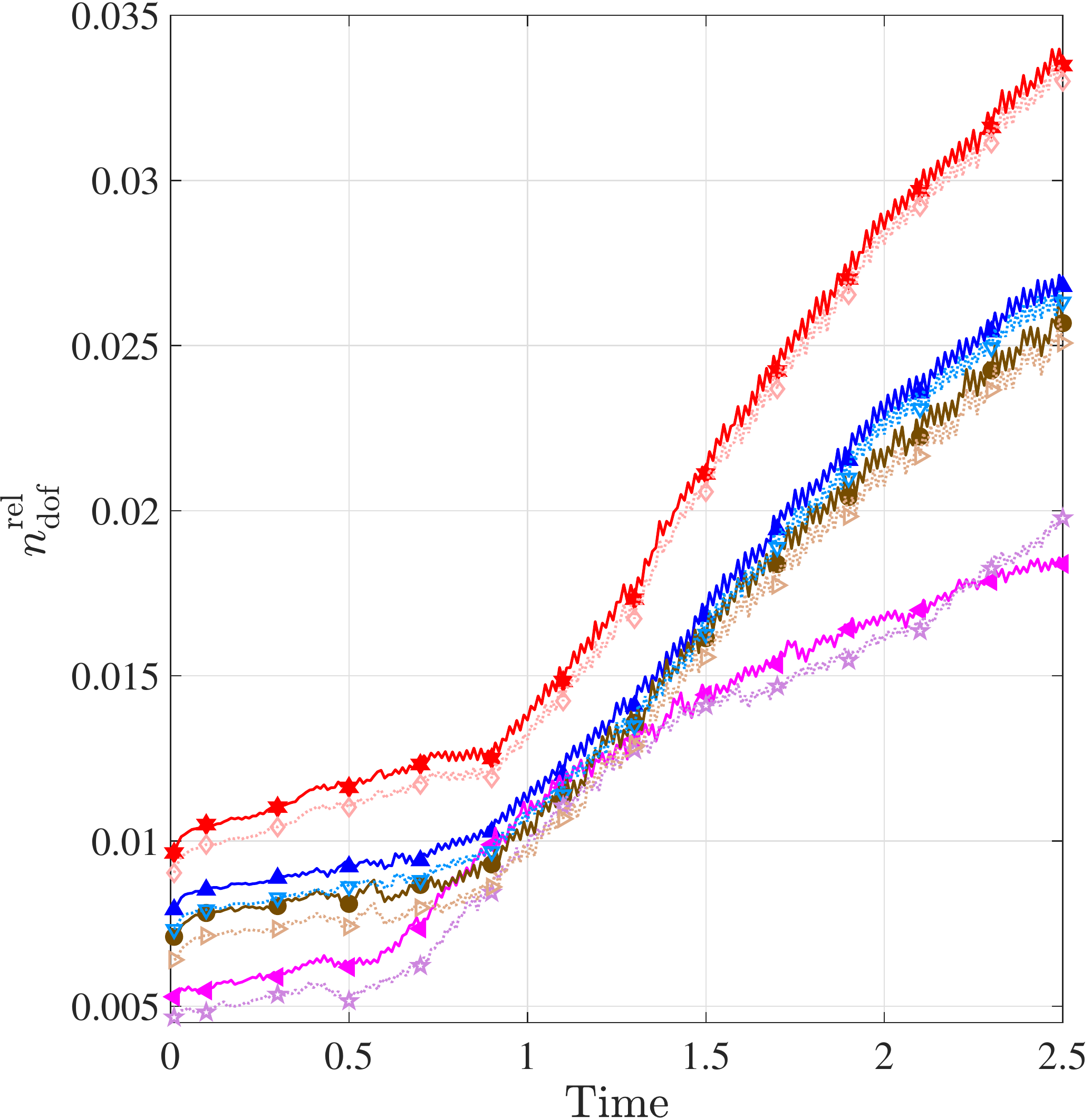}
\caption{$n_{\text{dof}}^{\text{rel}}$}  \label{figure_5_dof}
\end{subfigure}
\begin{subfigure}{0.59\columnwidth} \centering
\includegraphics[trim={0cm 0cm 0cm 0cm},clip,width=0.98\columnwidth]{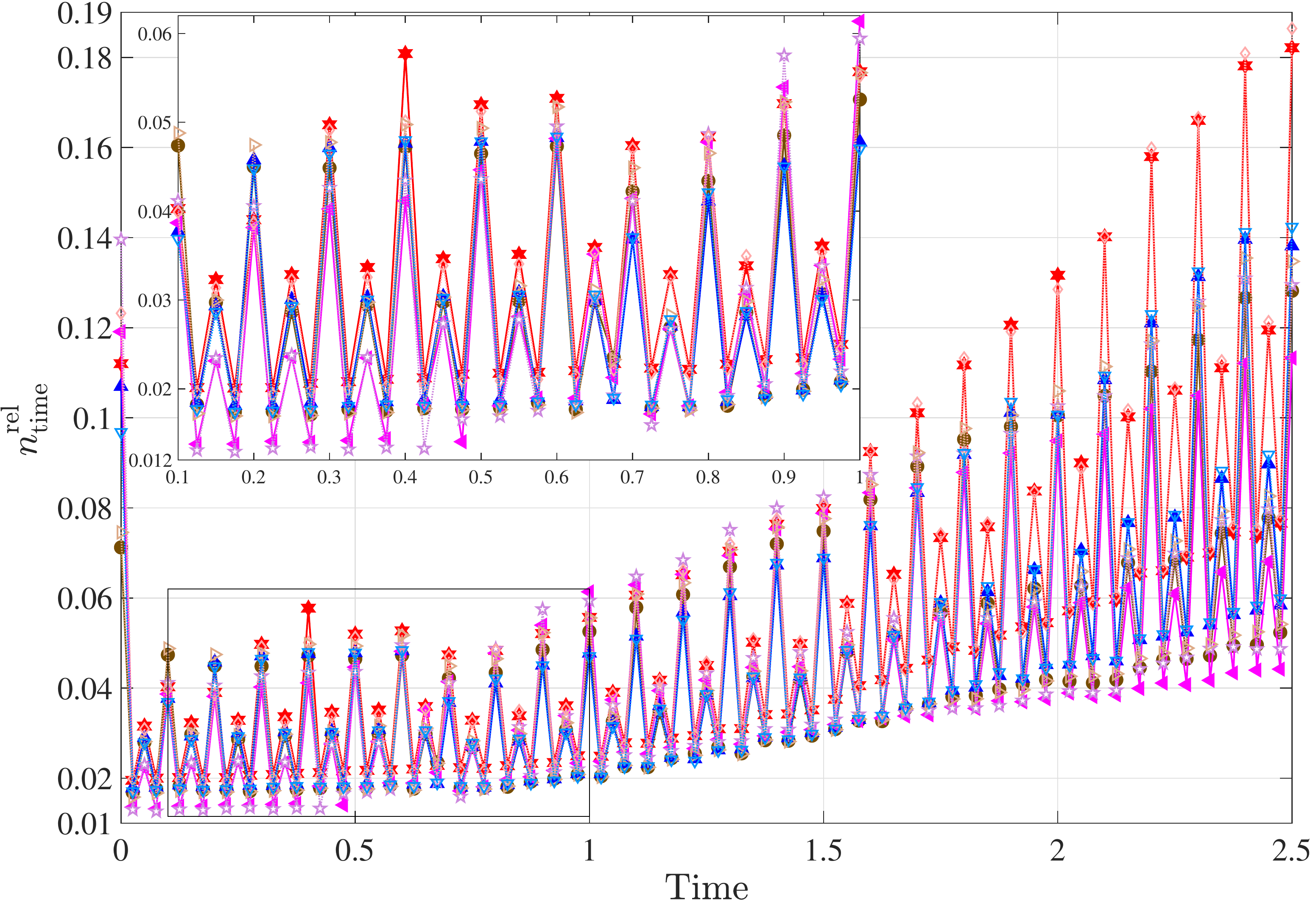}
\caption{$n_{\text{time}}^{\text{rel}}$} \label{figure_5_time}
\end{subfigure}
\caption{Relative $n_{\text{dof}}^{\text{rel}}$ and $n_{\text{time}}^{\text{rel}}$ (defined in Section~\ref{Numerical_studies}) over the time interval $\left[0, 2.5\right]$ for different adaptive THB mesh configurations. To ensure a fair comparison of the computational effort, the reference used in (a) and (b) correspond to a uniform tensor-product B-spline mesh with $p = 3$ and $1024$ elements in each direction, which provides a solution in close agreement with the reference solution (Fig.~\ref{figure_2_final}) shown in Fig.~\ref{figure_1_final}.}
\label{figure_5a_final}
\end{figure}

\begin{figure}[!ht]\centering
\begin{subfigure}{1\columnwidth} \centering
\includegraphics[trim={0cm -0.25cm 0cm 0cm},clip,width=0.98\columnwidth]{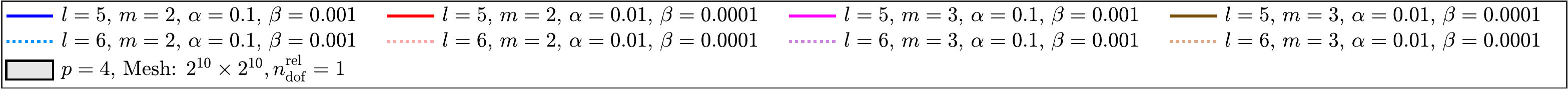}
\end{subfigure}
\begin{subfigure}{0.49\columnwidth} \centering
\includegraphics[trim={0cm 0cm 0cm 0cm},clip,width=0.58\columnwidth]{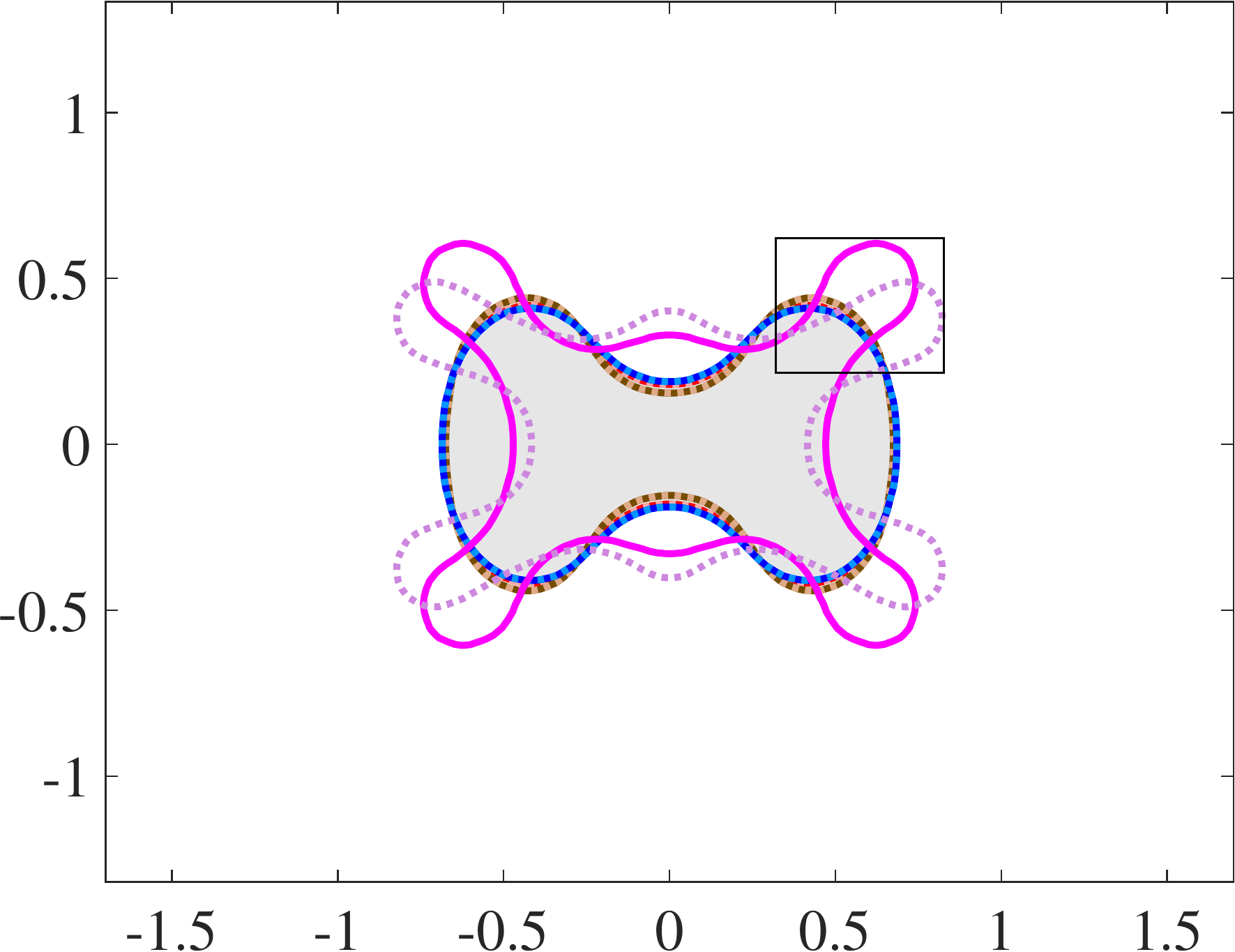}
\includegraphics[trim={0cm -3cm 0cm 0cm},clip,width=0.39\columnwidth]{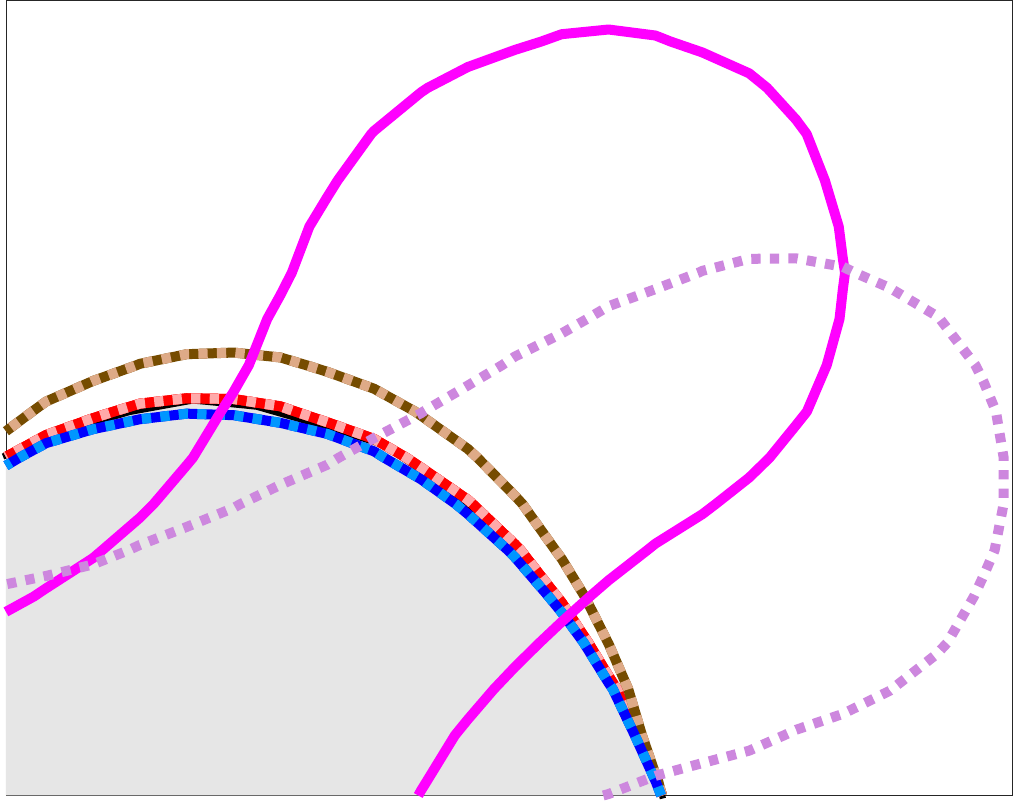}
\caption{$t = 1.0$} %
\end{subfigure}
\begin{subfigure}{0.49\columnwidth} \centering
\includegraphics[trim={0cm 0cm 0cm 0cm},clip,width=0.58\columnwidth]{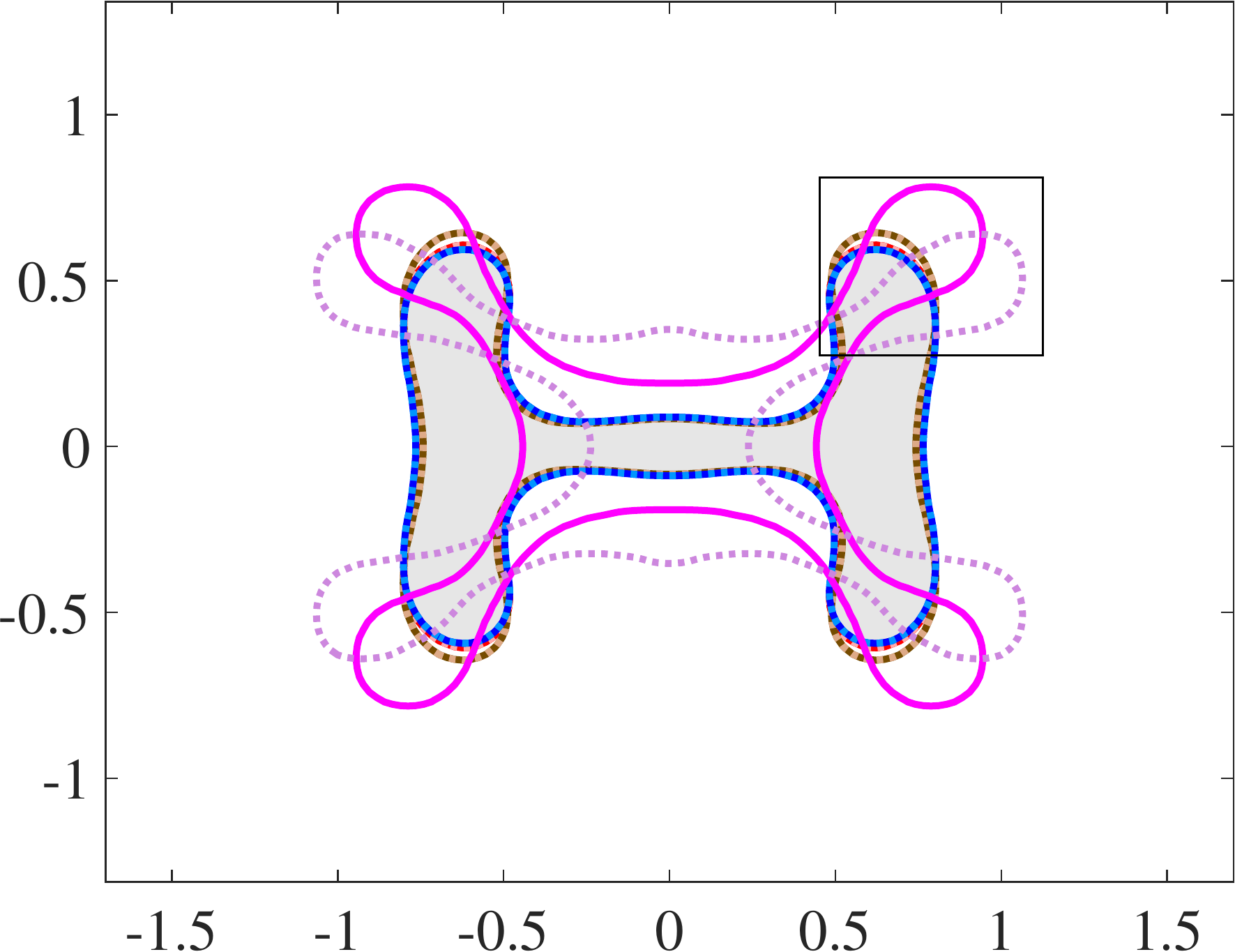}
\includegraphics[trim={0cm -3cm 0cm 0cm},clip,width=0.39\columnwidth]{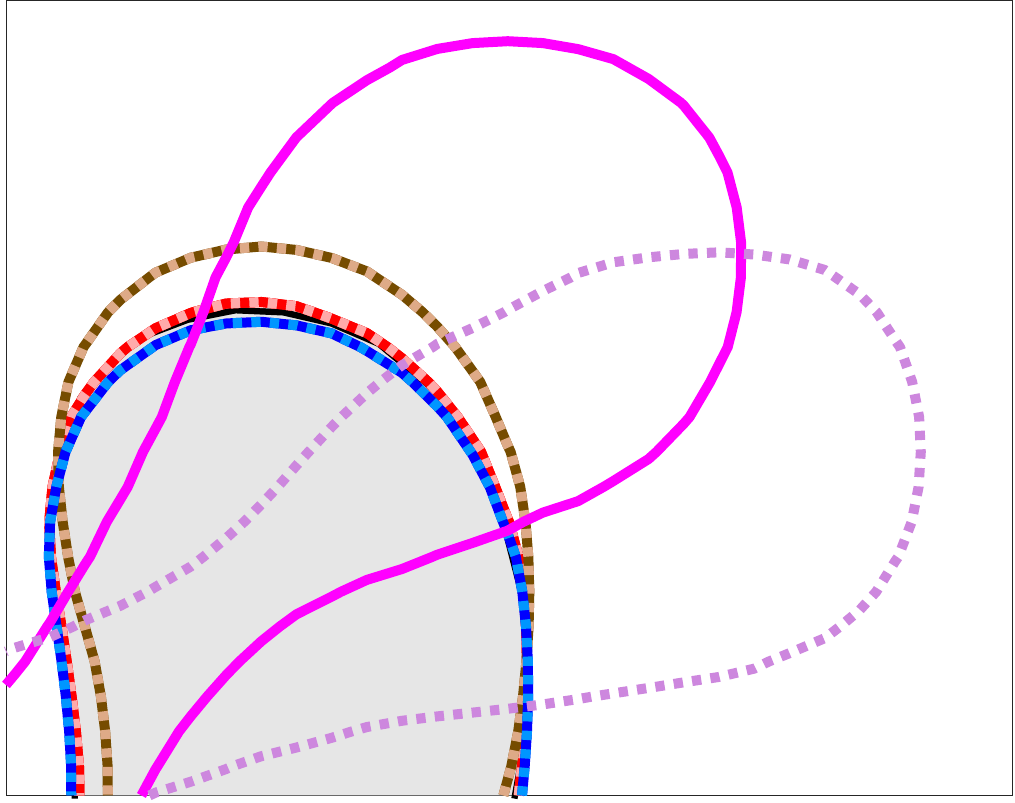}
\caption{$t = 1.5$} %
\end{subfigure}
\\
\begin{subfigure}{0.49\columnwidth} \centering
\includegraphics[trim={0cm 0cm 0cm 0cm},clip,width=0.58\columnwidth]{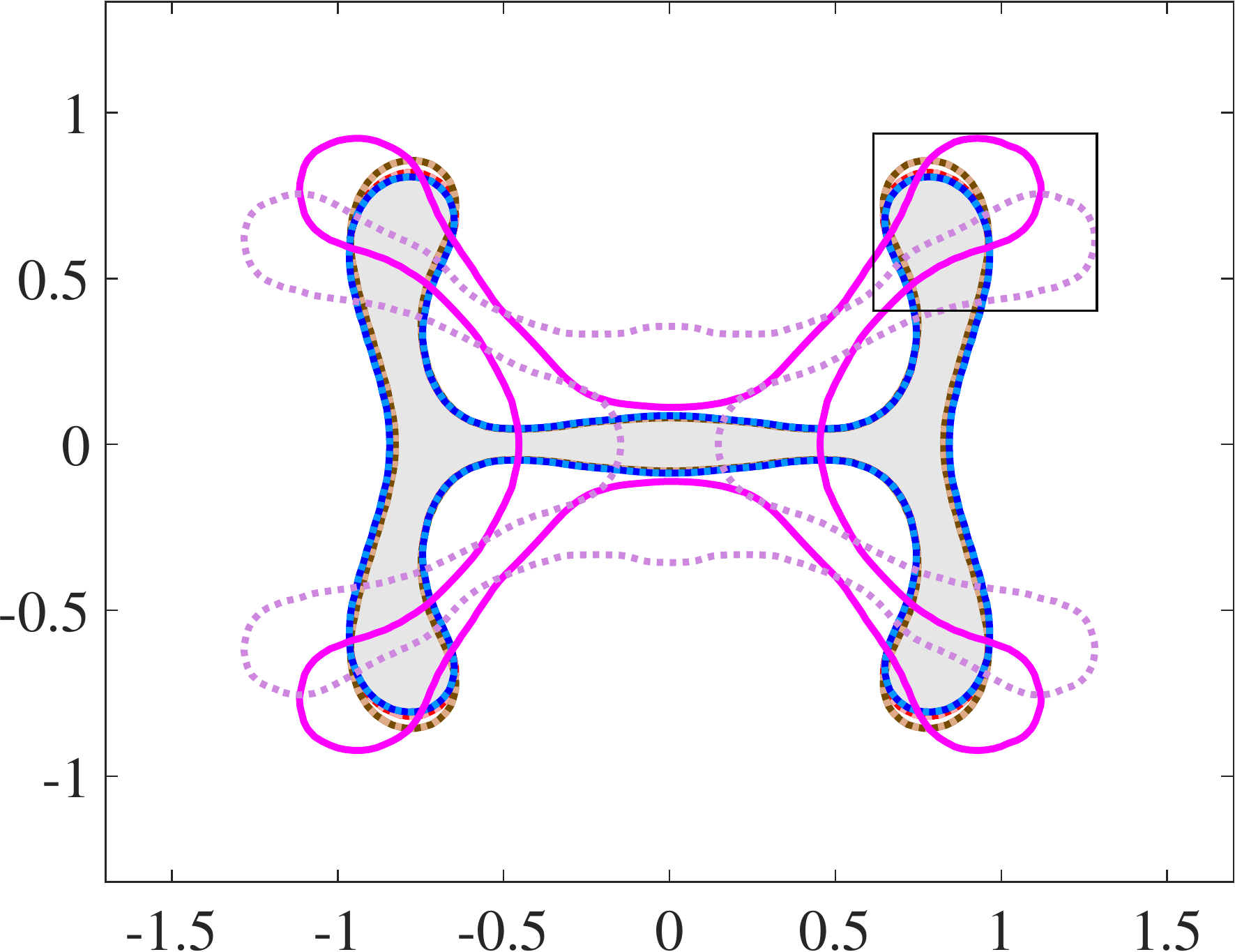}
\includegraphics[trim={0cm -3cm 0cm 0cm},clip,width=0.39\columnwidth]{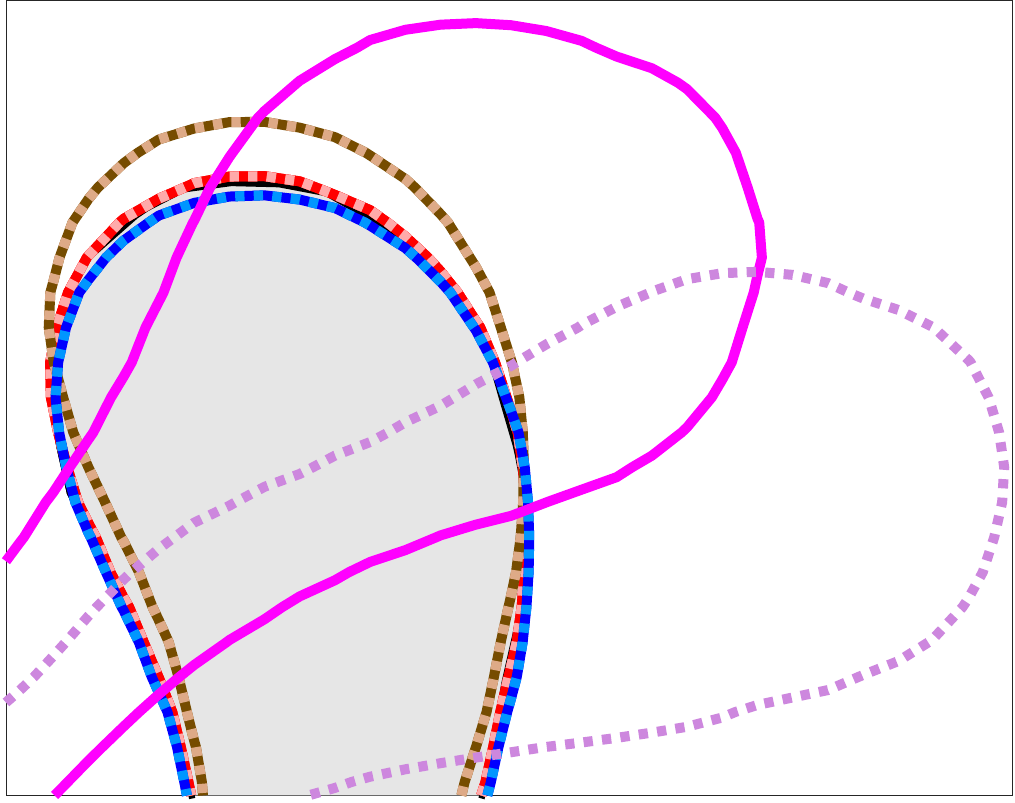}
\caption{$t = 2.0$} %
\end{subfigure}
\begin{subfigure}{0.49\columnwidth} \centering
\includegraphics[trim={0cm 0cm 0cm 0cm},clip,width=0.58\columnwidth]{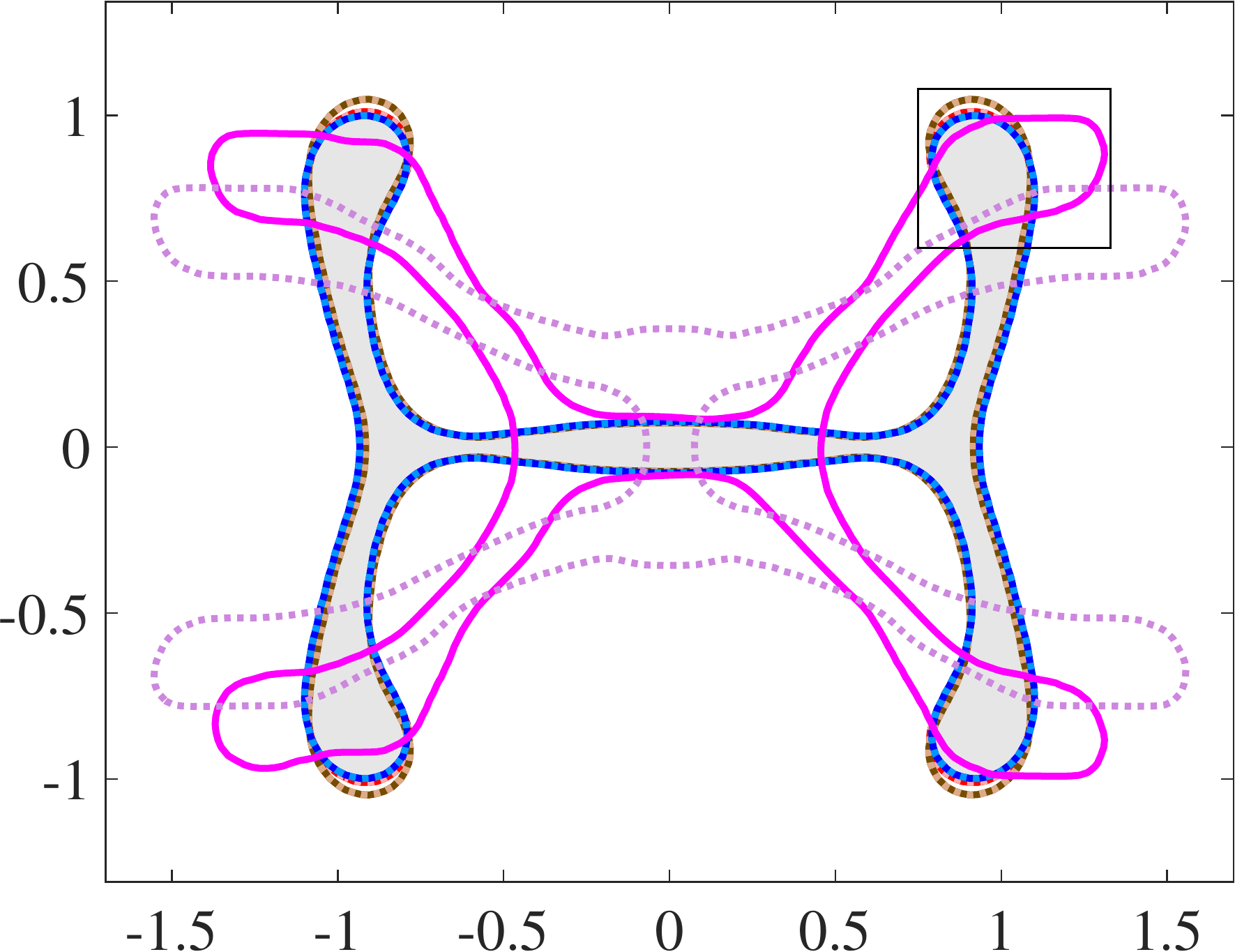}
\includegraphics[trim={0cm -3cm 0cm 0cm},clip,width=0.39\columnwidth]{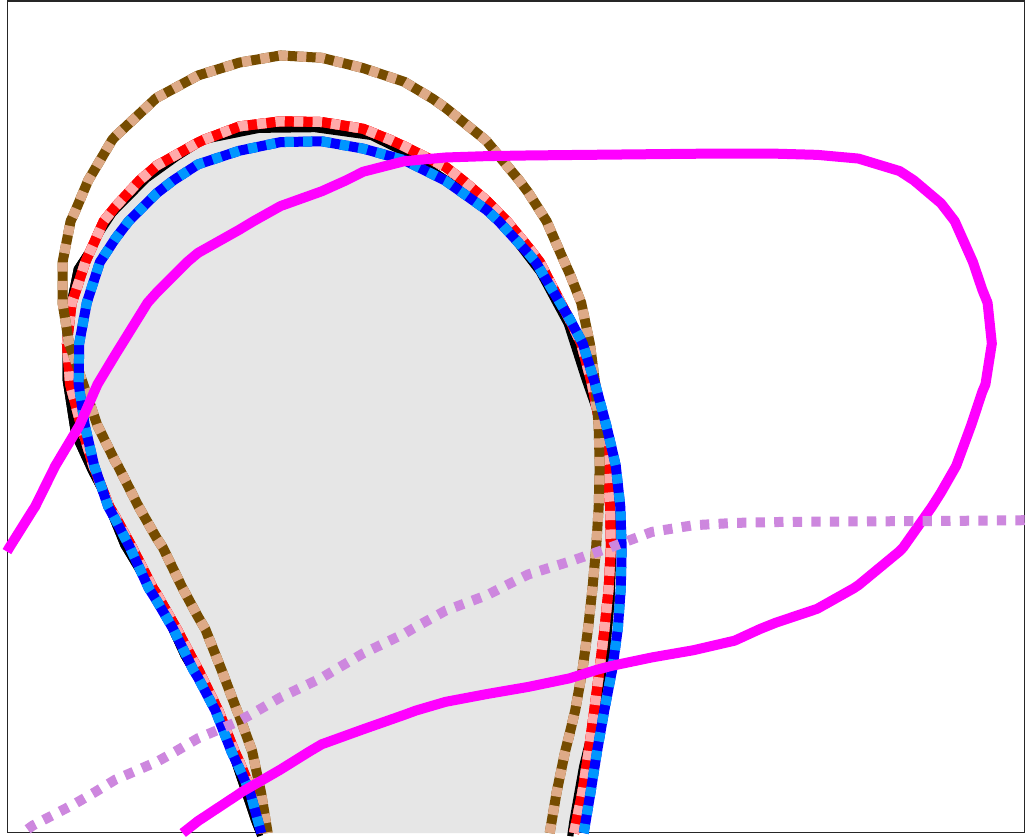}
\caption{$t = 2.5$} %
\end{subfigure}
\caption{Contour traced along the interface separating healthy and tumor tissue ($\phi \approx 0$)  for different adaptive THB-spline mesh configurations at time 1, 1.5, 2.0, and 2.5.
The reference solution is shown as gray area representing the tumor shape, and it corresponds to the uniform tensor product B-spline mesh with $p = 4$ and $1024$ elements in each direction shown in Fig.~\ref{figure_1_final}.}
\label{figure_6_final}
\end{figure}

\begin{figure}[!ht]\centering
\begin{subfigure}{0.185\columnwidth} \centering
\includegraphics[trim={0cm 0cm 0cm 0cm},clip,width=1\columnwidth]{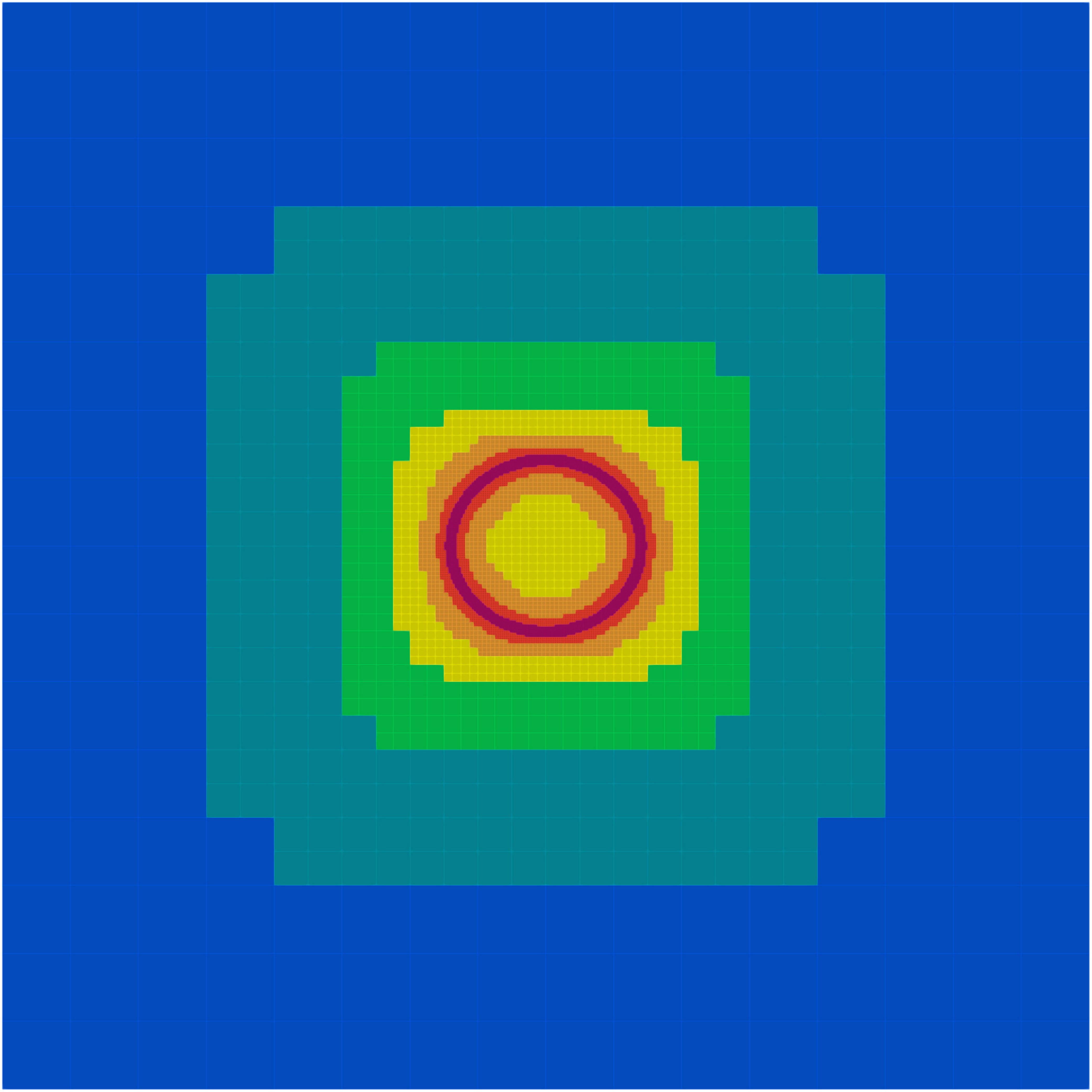}  \end{subfigure}
\begin{subfigure}{0.185\columnwidth} \centering
\includegraphics[trim={0cm 0cm 0cm 0cm},clip,width=1\columnwidth]{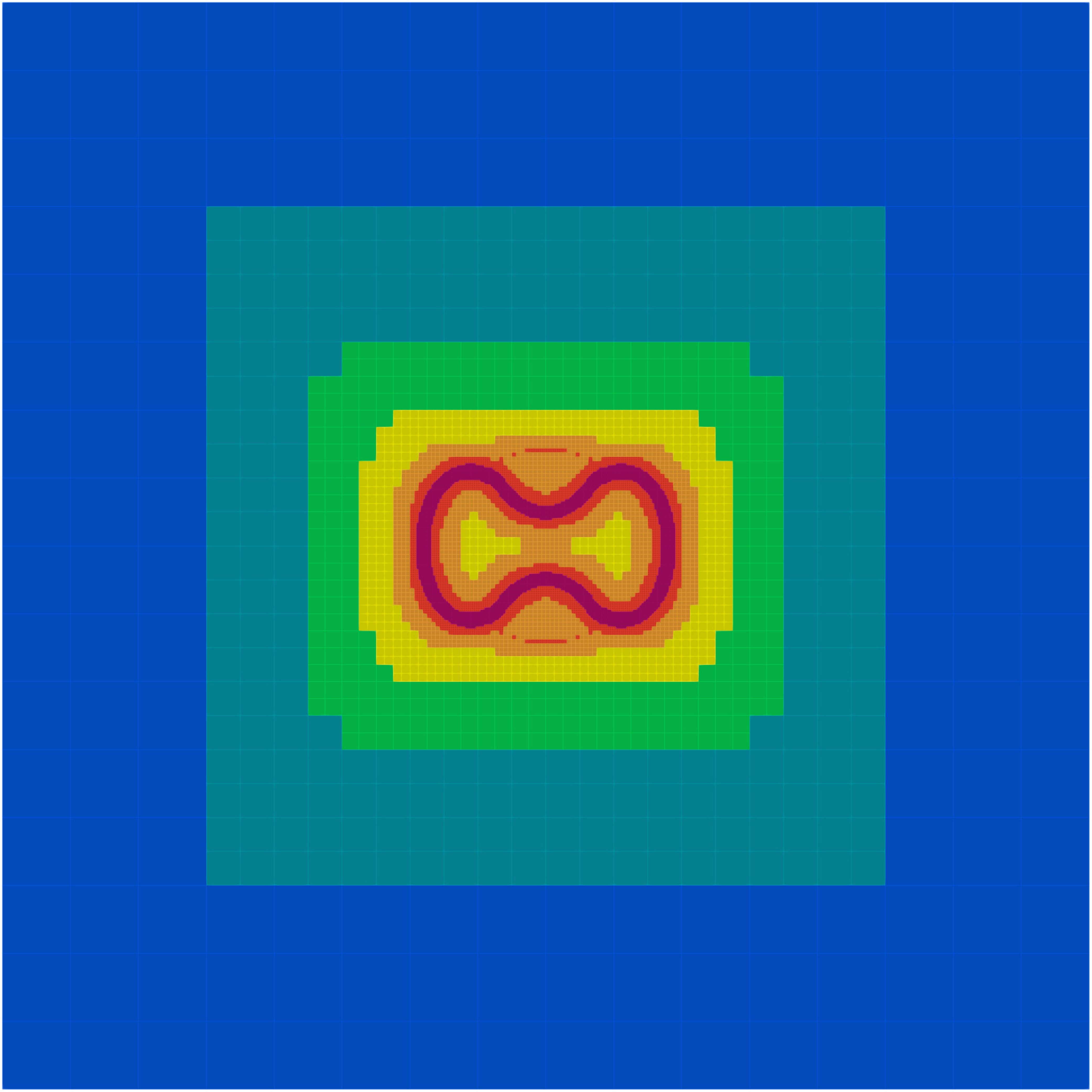}   \end{subfigure}
\begin{subfigure}{0.185\columnwidth} \centering
\includegraphics[trim={0cm 0cm 0cm 0cm},clip,width=1\columnwidth]{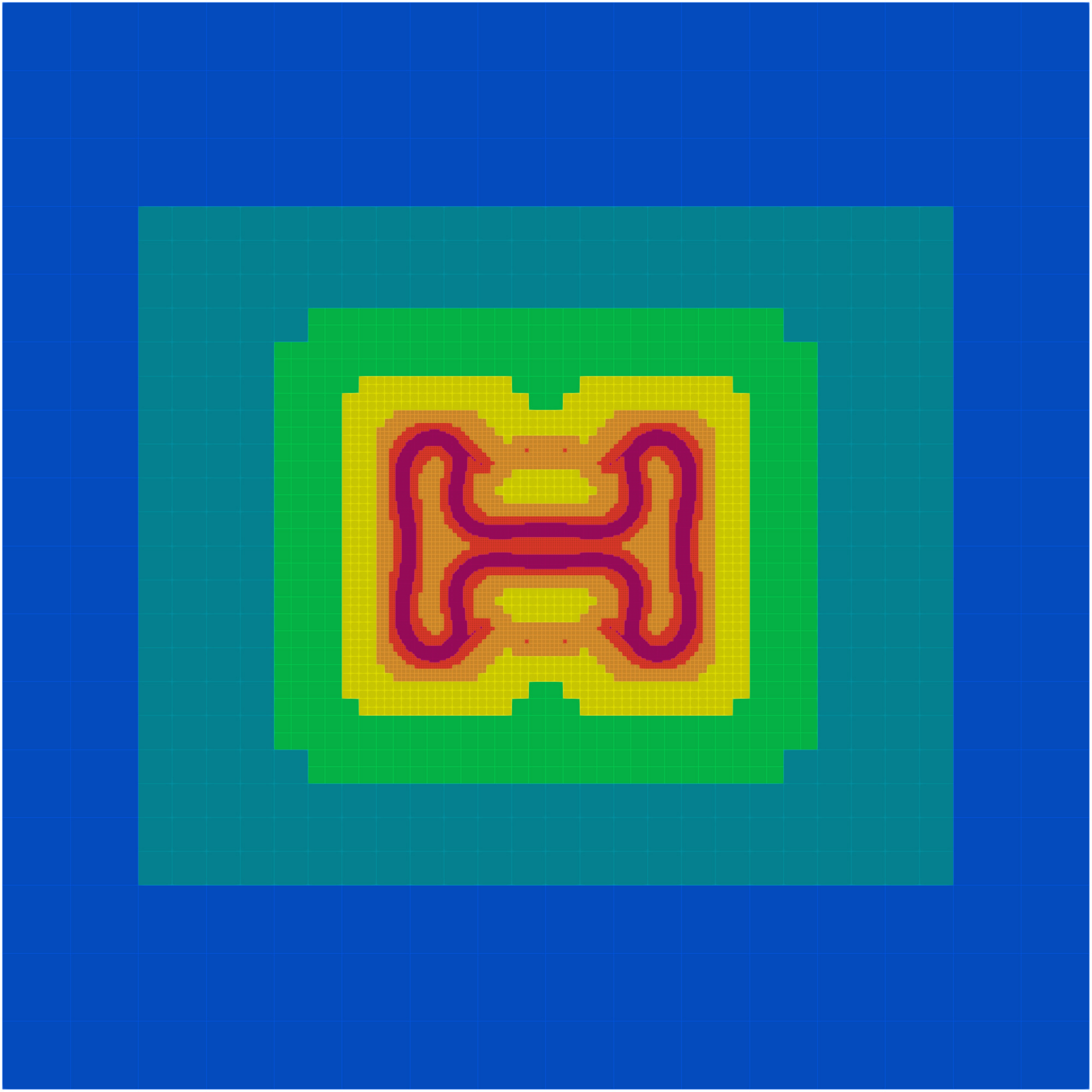} \end{subfigure}
\begin{subfigure}{0.185\columnwidth} \centering
\includegraphics[trim={0cm 0cm 0cm 0cm},clip,width=1\columnwidth]{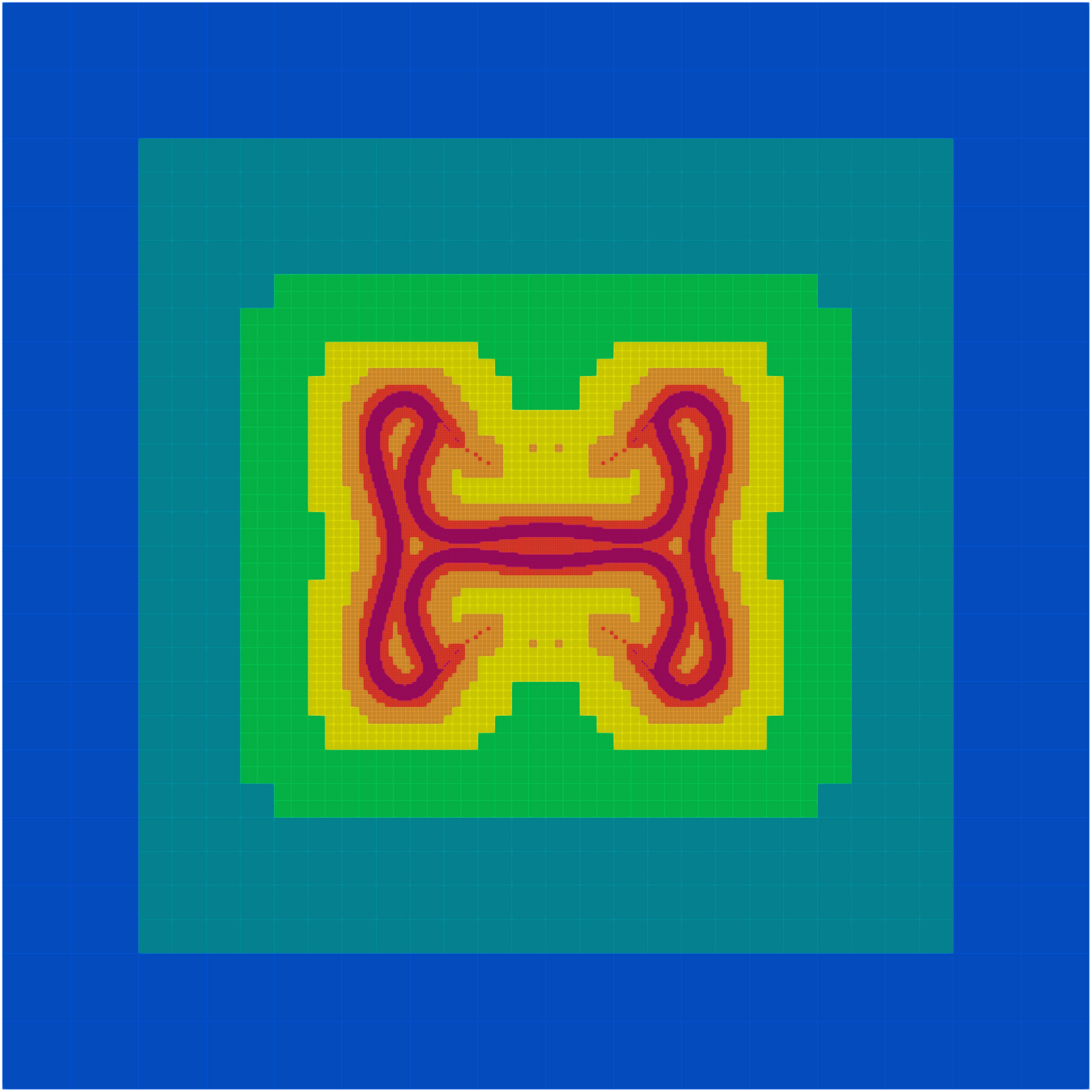} \end{subfigure}
\begin{subfigure}{0.185\columnwidth} \centering
\includegraphics[trim={0cm 0cm 0cm 0cm},clip,width=1\columnwidth]{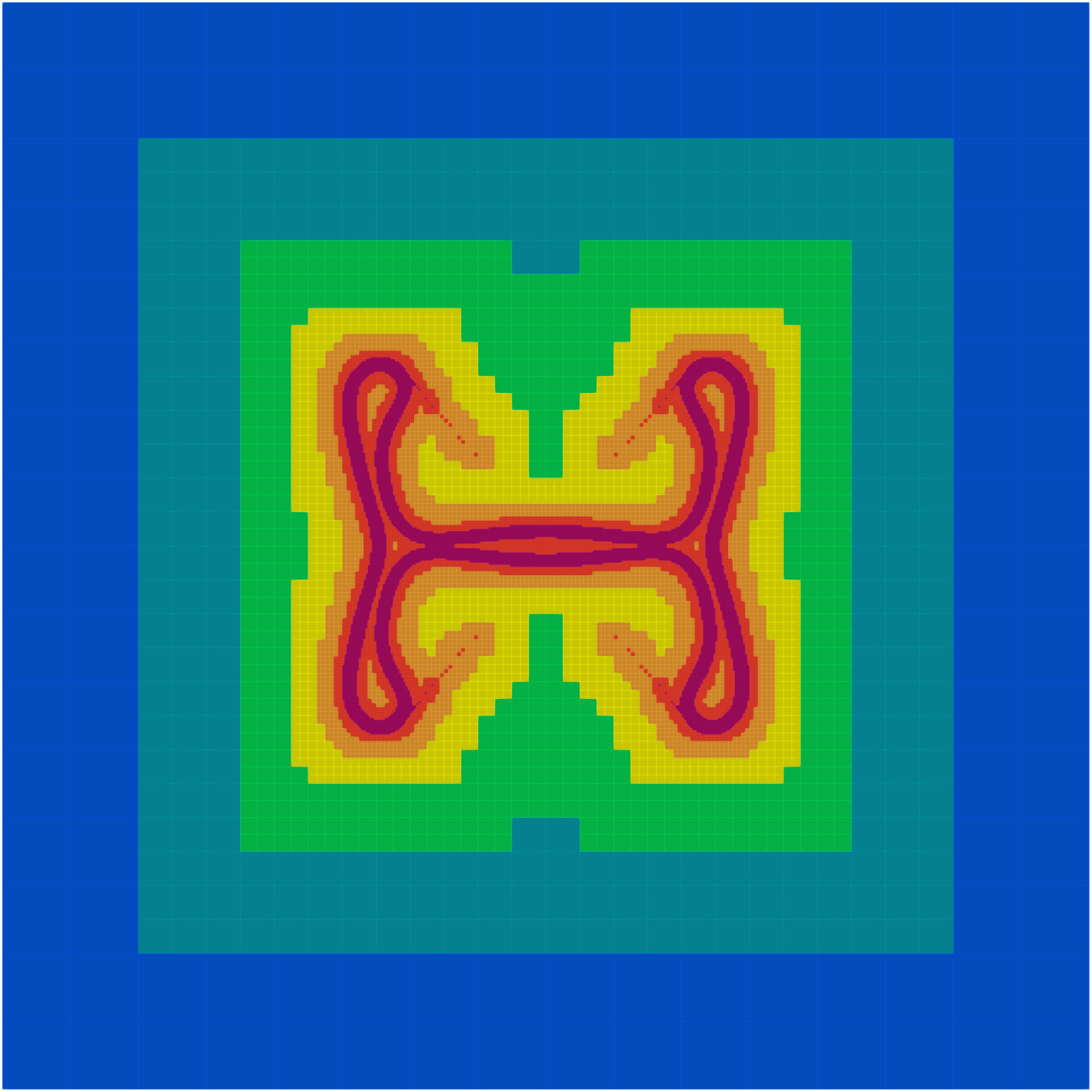}  \end{subfigure}
\begin{subfigure}{0.05\columnwidth} \centering
\includegraphics[trim={1cm 0cm 0cm 0cm},clip,width=0.84\columnwidth]{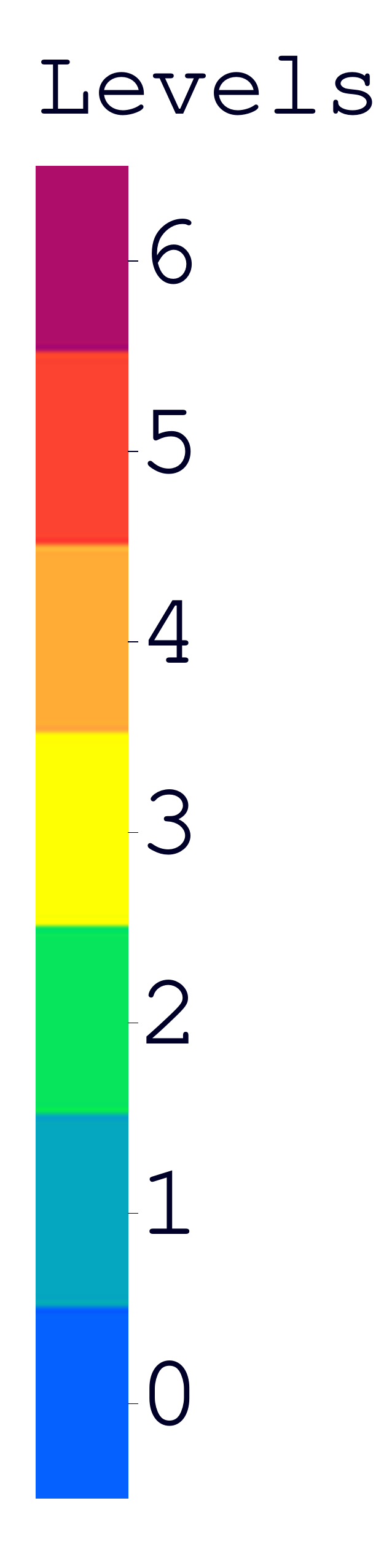} \end{subfigure}
\\ 
\begin{subfigure}{0.185\columnwidth} \centering
\includegraphics[trim={0cm 0cm 0cm 0cm},clip,width=1\columnwidth]{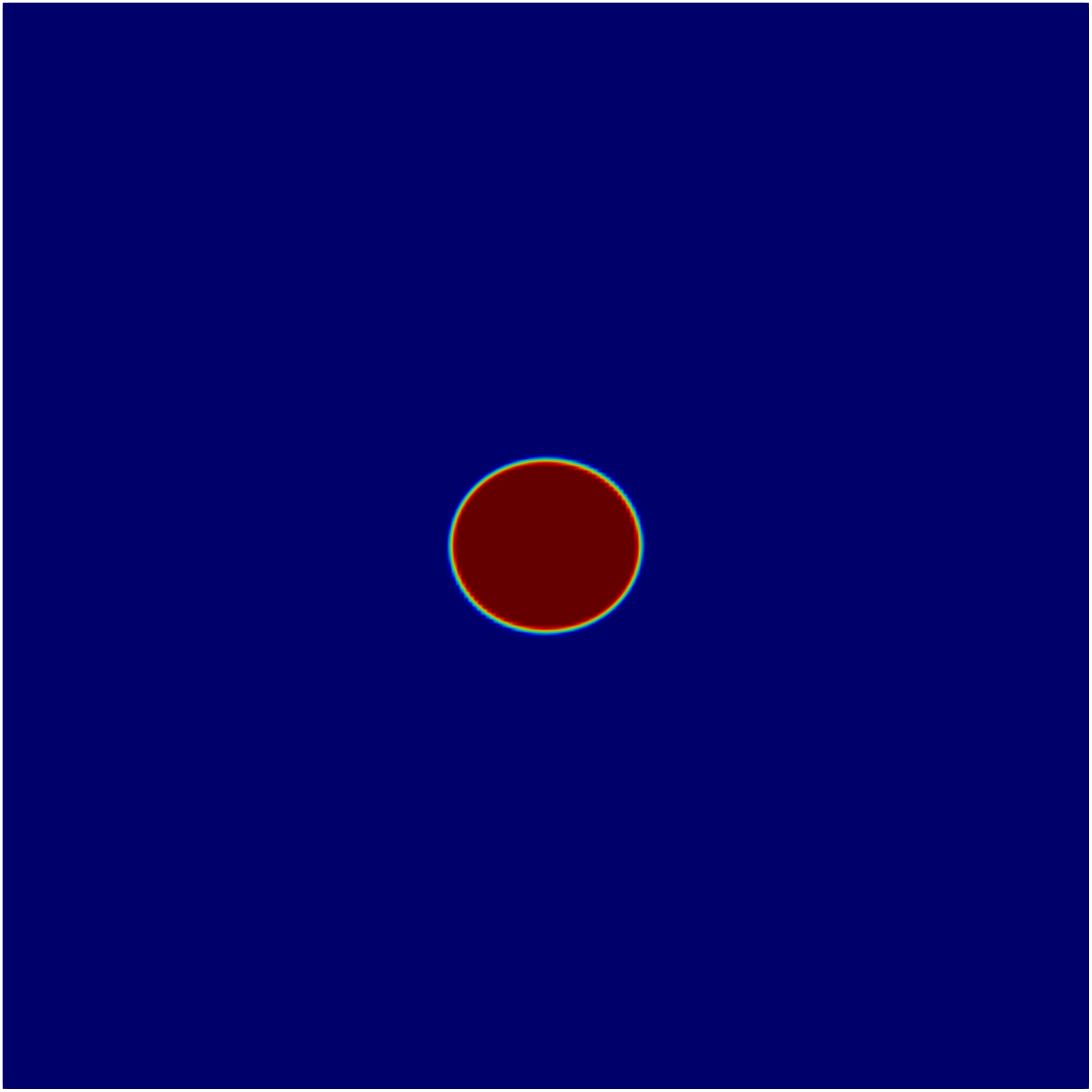}  \end{subfigure}
\begin{subfigure}{0.185\columnwidth} \centering
\includegraphics[trim={0cm 0cm 0cm 0cm},clip,width=1\columnwidth]{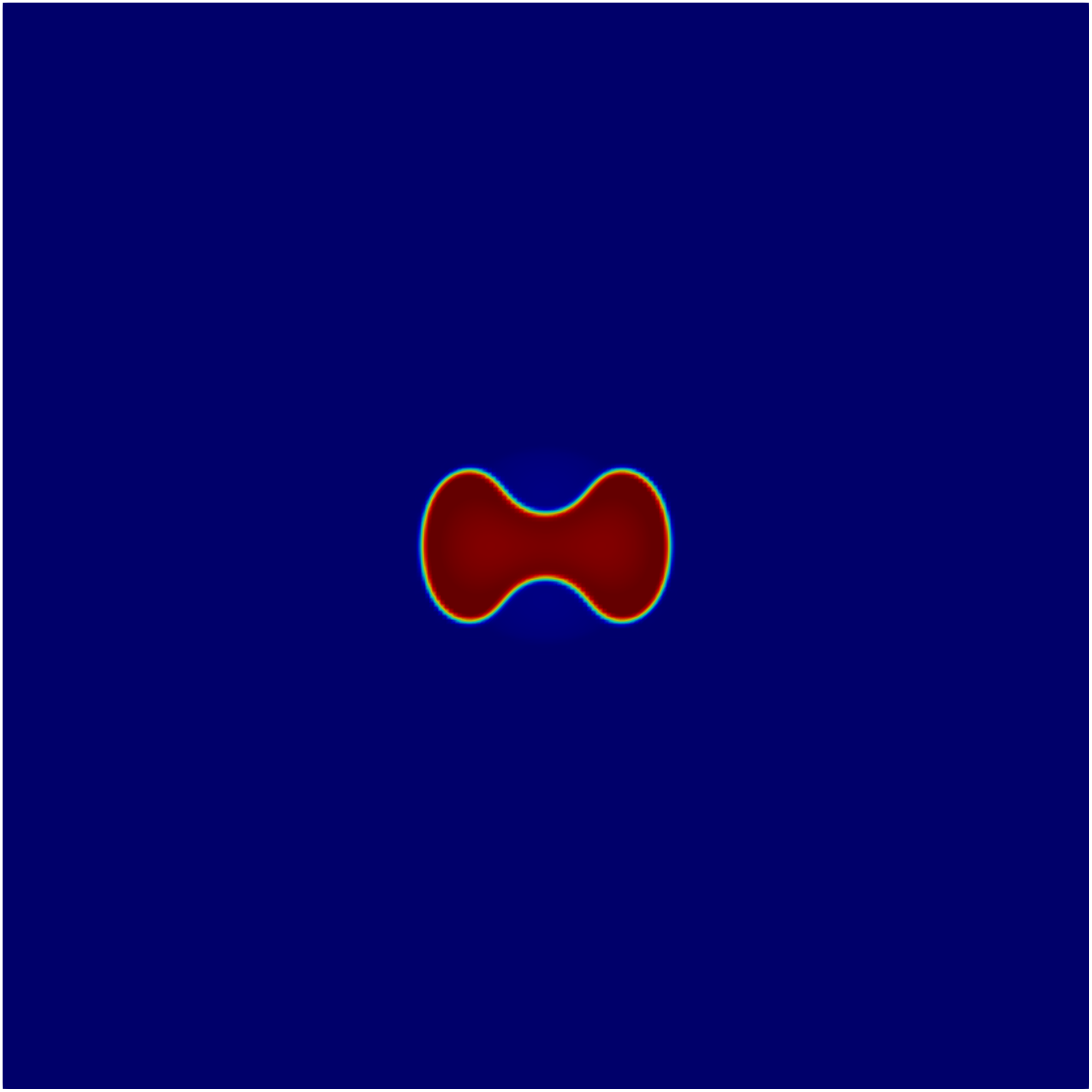}   \end{subfigure}
\begin{subfigure}{0.185\columnwidth} \centering
\includegraphics[trim={0cm 0cm 0cm 0cm},clip,width=1\columnwidth]{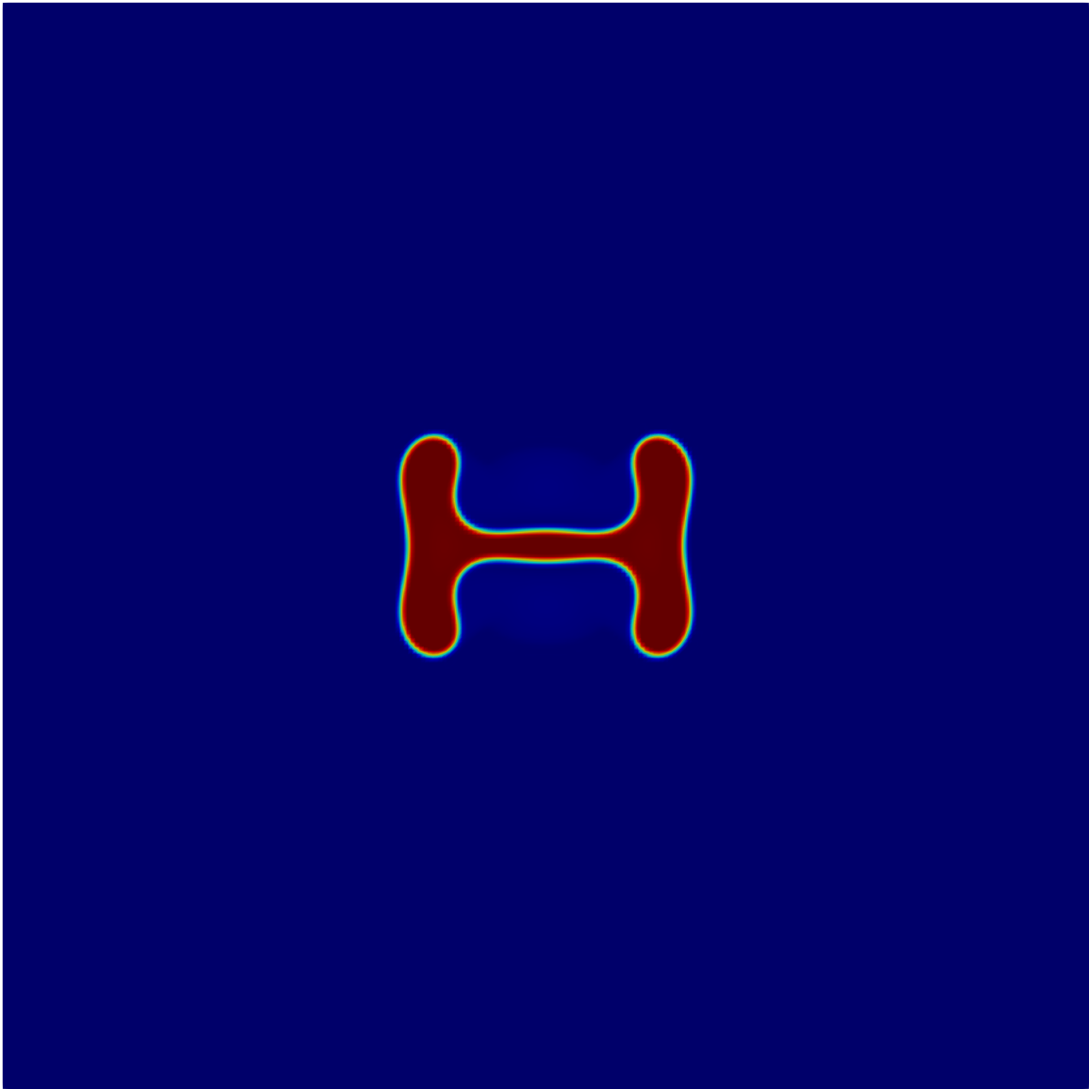} \end{subfigure}
\begin{subfigure}{0.185\columnwidth} \centering
\includegraphics[trim={0cm 0cm 0cm 0cm},clip,width=1\columnwidth]{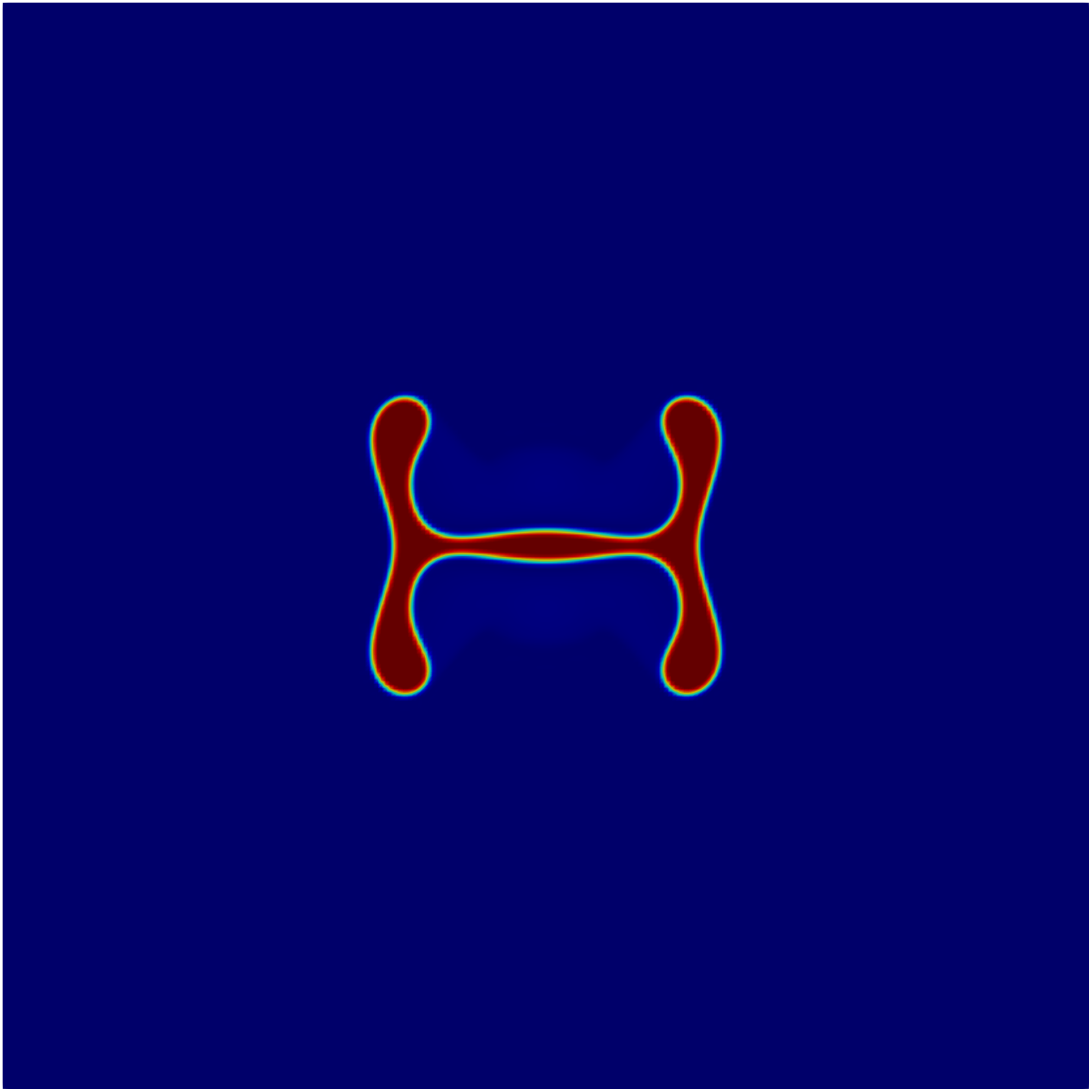} \end{subfigure}
\begin{subfigure}{0.185\columnwidth} \centering
\includegraphics[trim={0cm 0cm 0cm 0cm},clip,width=1\columnwidth]{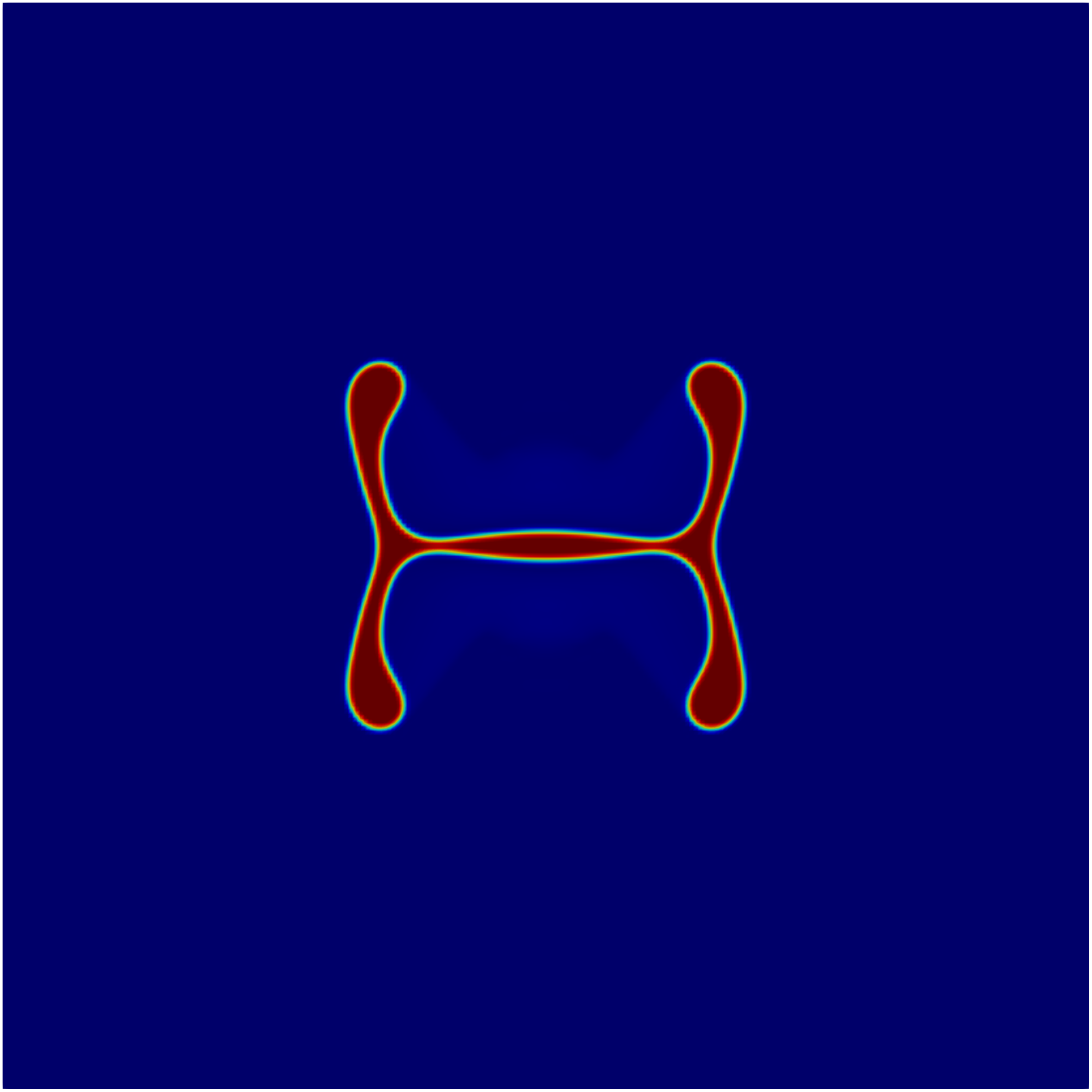}  \end{subfigure}
\begin{subfigure}{0.05\columnwidth} \centering
\includegraphics[trim={0cm 1cm 0cm 0cm},clip,width=1.2\columnwidth]{Images/Problem1_2DSq_UniRef/ColoarBar_Phi-eps-converted-to.pdf} \end{subfigure}
\\
\begin{subfigure}{0.185\columnwidth} \centering
\includegraphics[trim={0cm 0cm 0cm 0cm},clip,width=1\columnwidth]{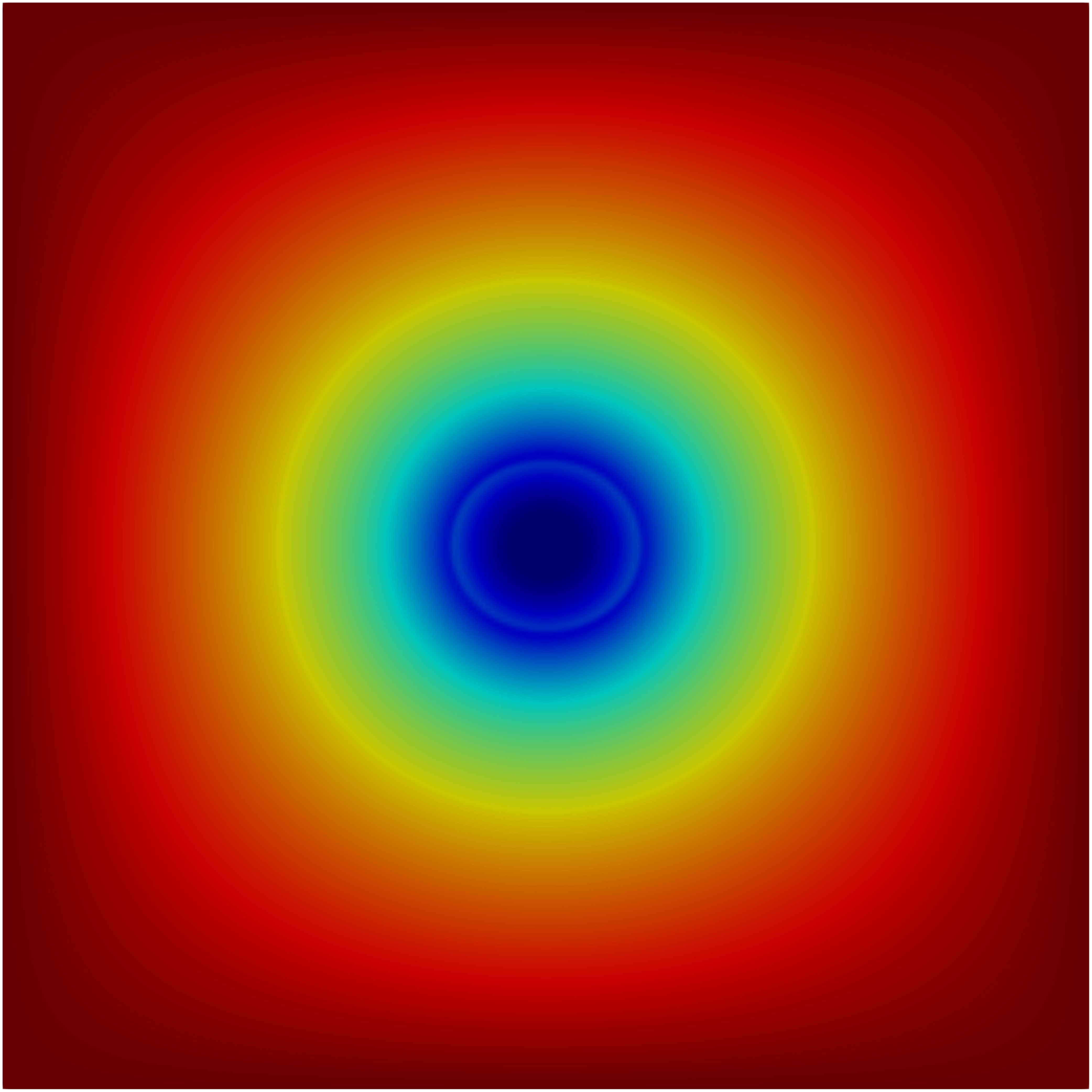} \caption{$t = 0$} \end{subfigure}
\begin{subfigure}{0.185\columnwidth} \centering
\includegraphics[trim={0cm 0cm 0cm 0cm},clip,width=1\columnwidth]{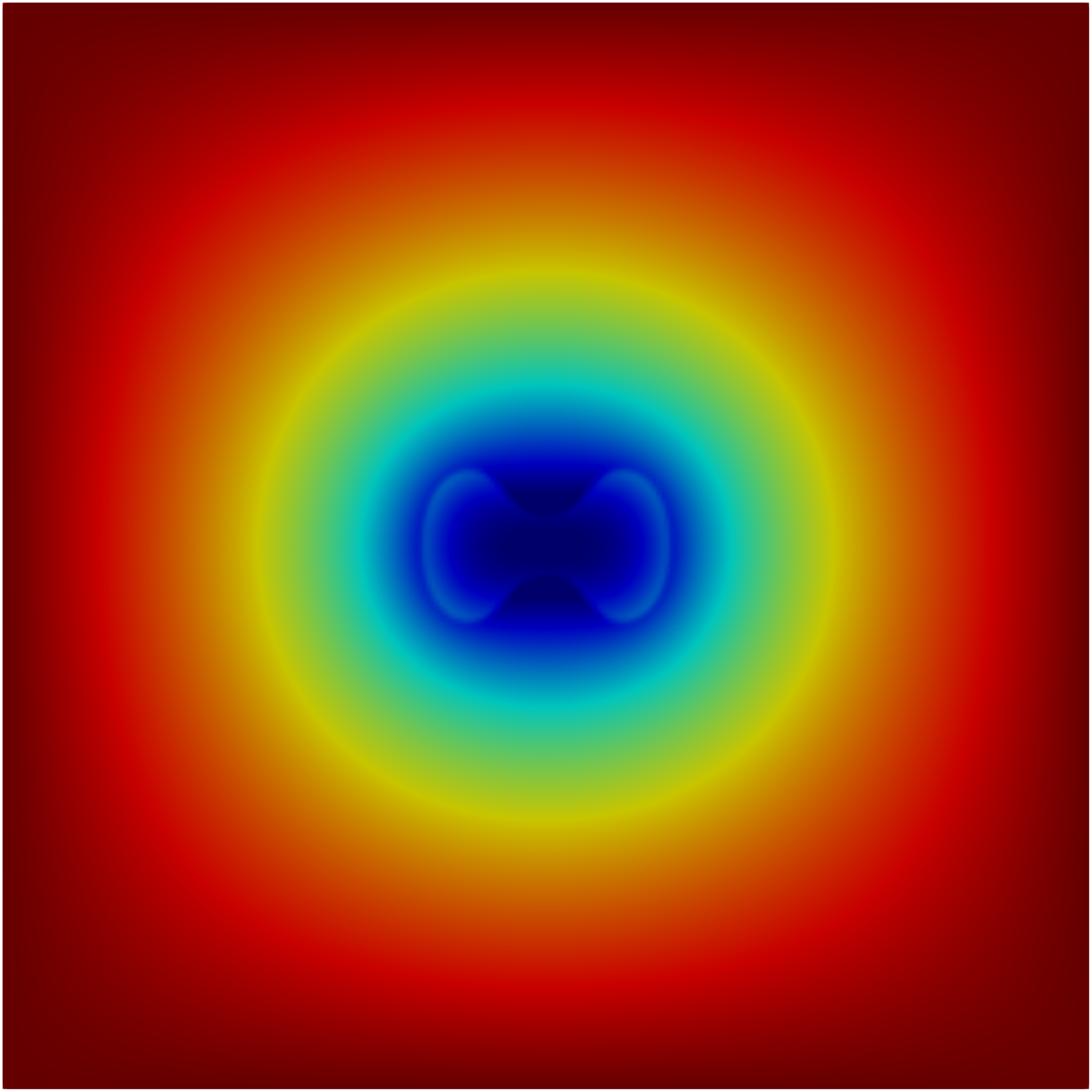} \caption{$t = 1$}  \end{subfigure}
\begin{subfigure}{0.185\columnwidth} \centering
\includegraphics[trim={0cm 0cm 0cm 0cm},clip,width=1\columnwidth]{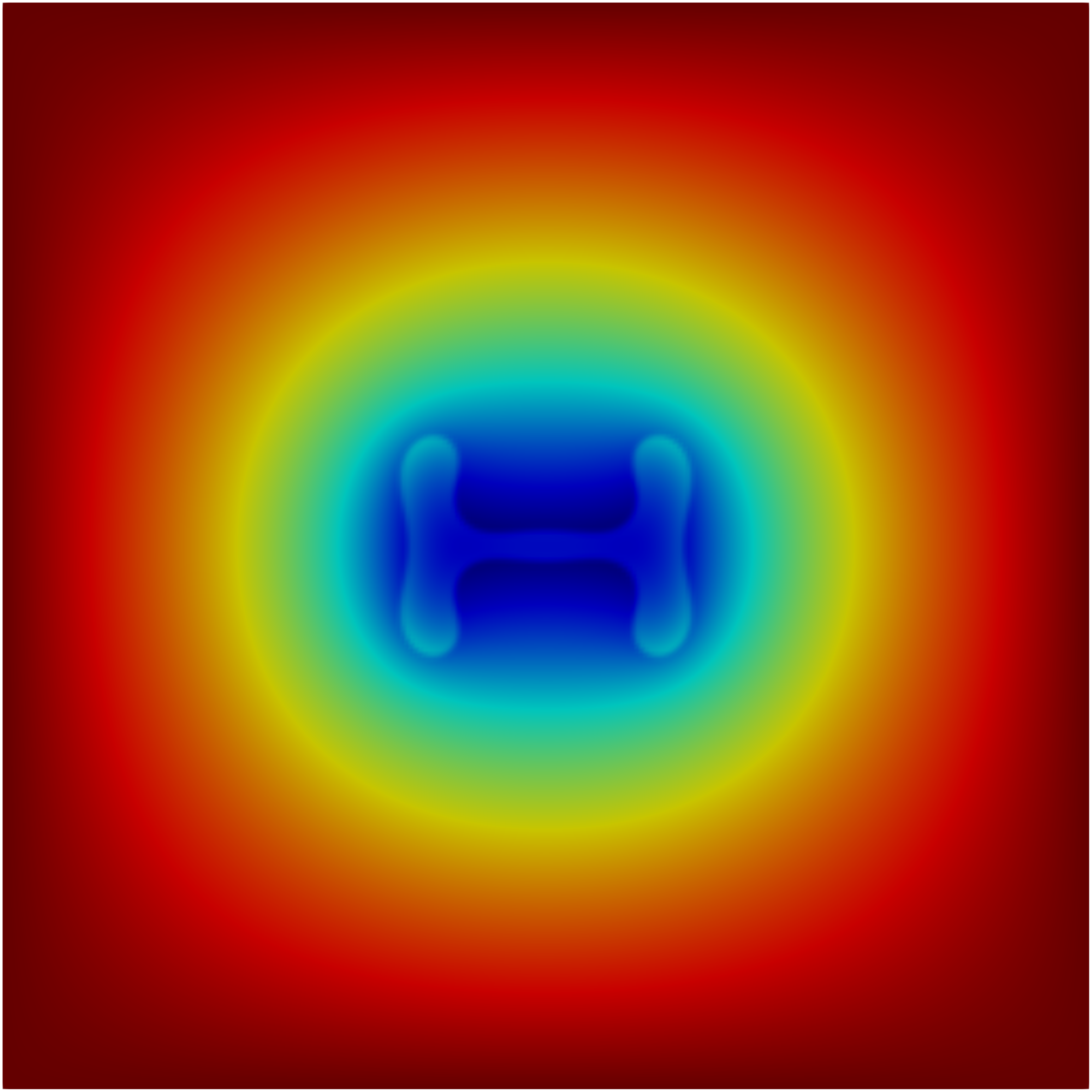} \caption{$t = 1.5$}\end{subfigure}
\begin{subfigure}{0.185\columnwidth} \centering
\includegraphics[trim={0cm 0cm 0cm 0cm},clip,width=1\columnwidth]{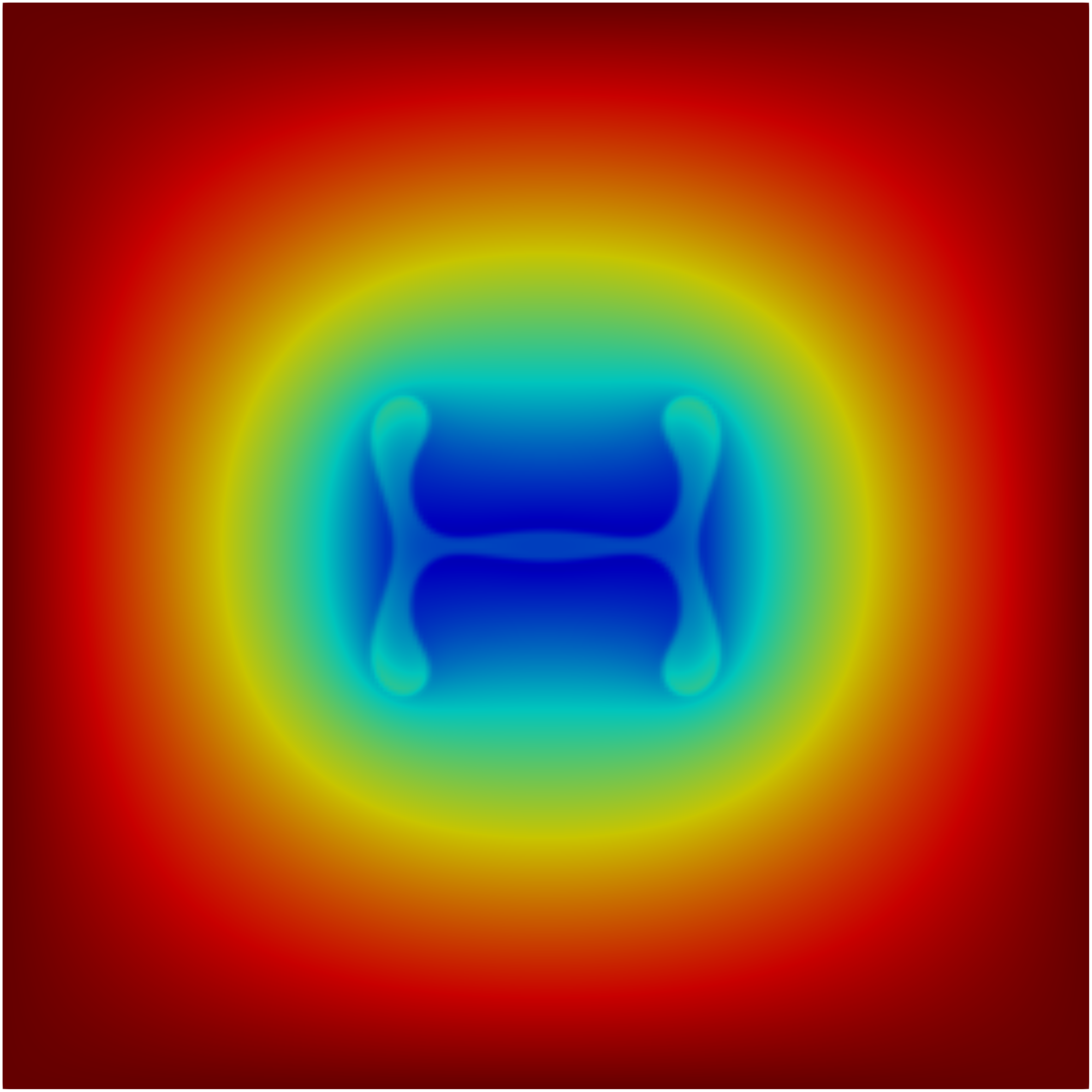} \caption{$t = 2$} \end{subfigure}
\begin{subfigure}{0.185\columnwidth} \centering
\includegraphics[trim={0cm 0cm 0cm 0cm},clip,width=1\columnwidth]{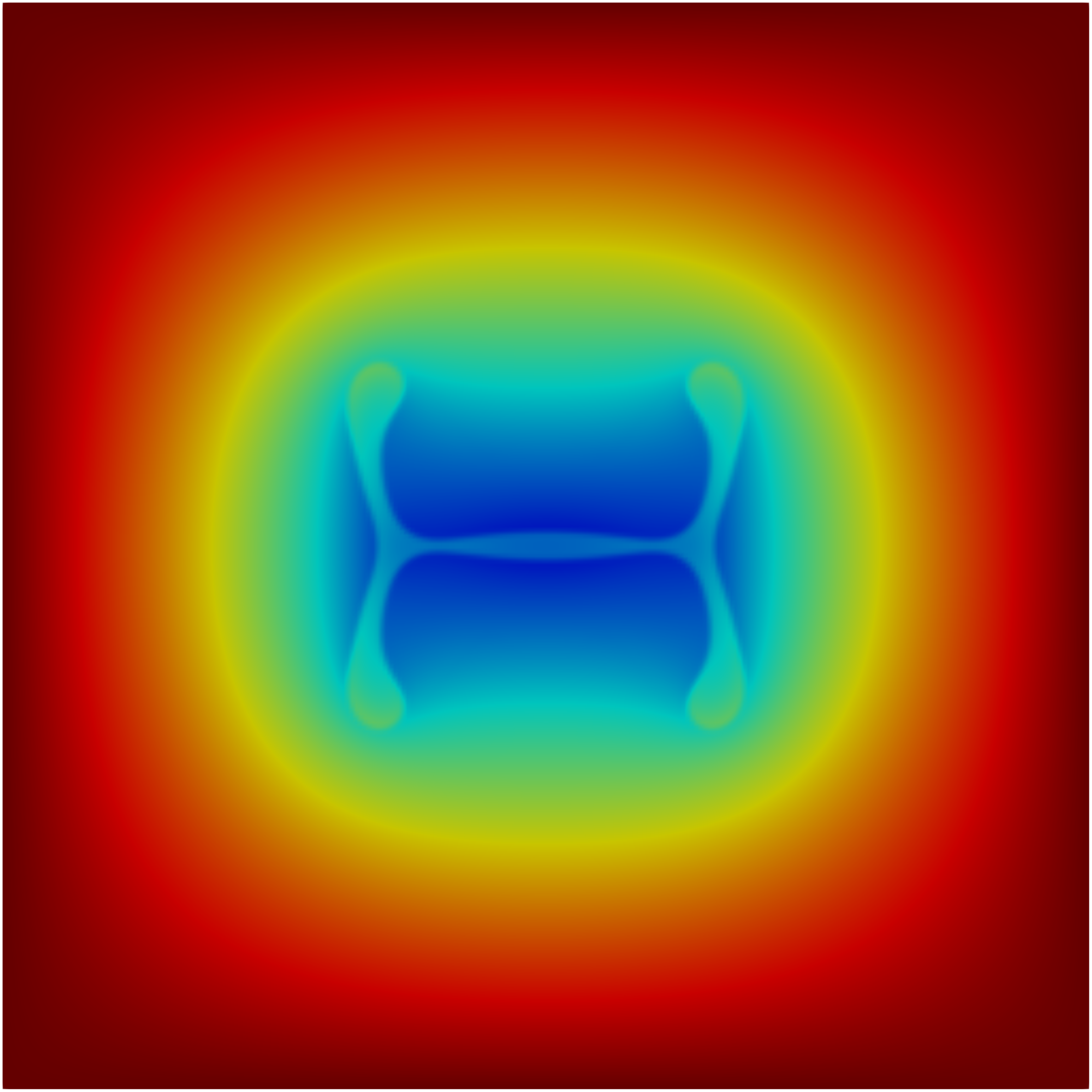} \caption{$t = 2.5$}\end{subfigure}
\begin{subfigure}{0.05\columnwidth} \centering
\includegraphics[trim={0cm -6cm 0cm 0cm},clip,width=1.2\columnwidth]{Images/Problem1_2DSq_UniRef/ColoarBar_Sigma-eps-converted-to.pdf} \end{subfigure}
\caption{Evolution of an elliptical tumor in square domain at time 0, 1, 1.5, 2, and 2.5. Top: adaptive THB-spline mesh configurations of degree $3$, $\ell = 6$, $m=2$, $\alpha = 0.01$, $\beta = 0.0001$ with a finest refinement level at $2^{10} \times 2^{10}$ ($h_e=6/1024$) mesh resolution. Middle: evolution of tumor geometry. Bottom: corresponding nutrient concentration.}
\label{figure_7_final}
\end{figure}

\begin{figure}[!t]\centering
\begin{subfigure}{0.185\columnwidth} \centering
\includegraphics[trim={0cm 0cm 0cm 0cm},clip,width=1\columnwidth]{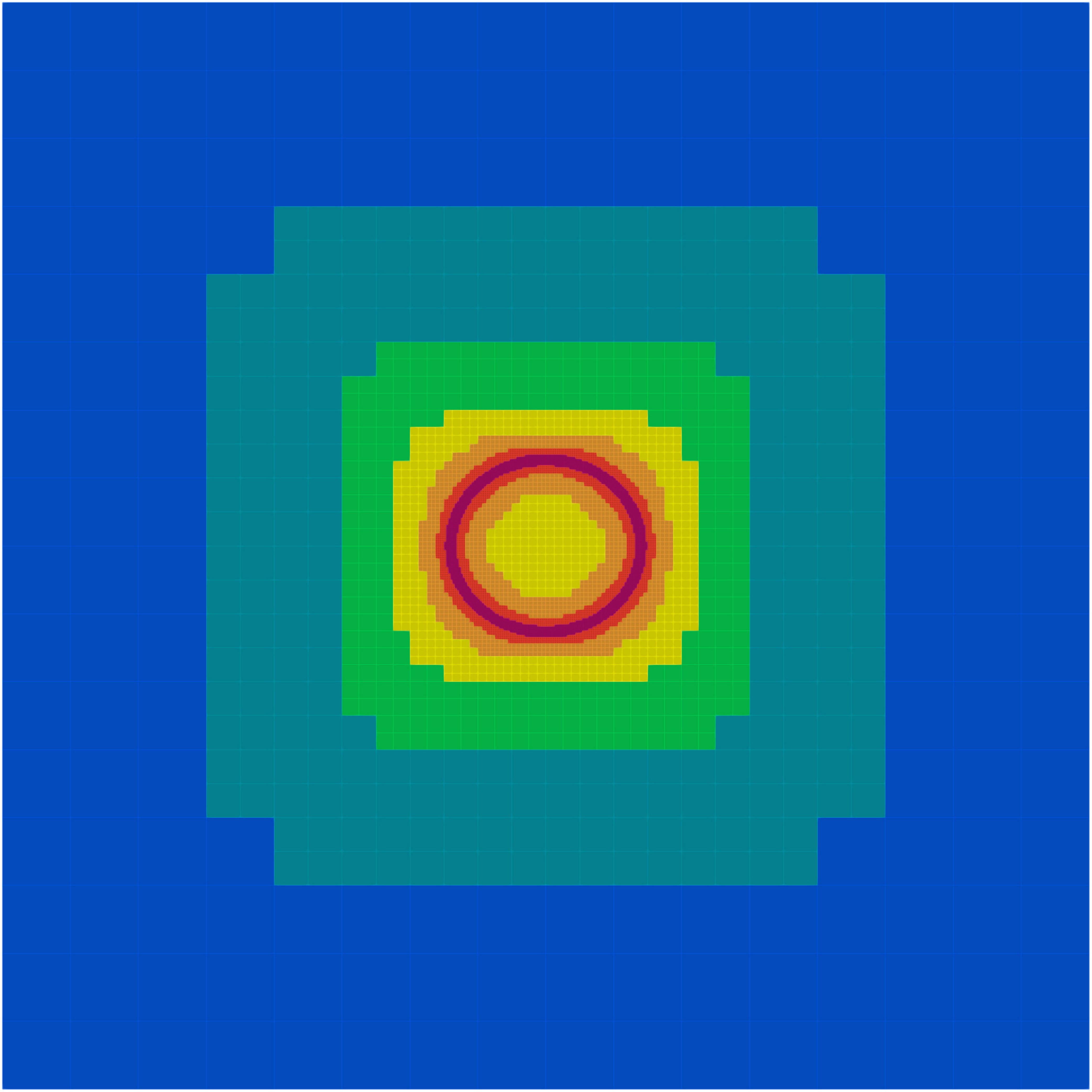}  \end{subfigure}
\begin{subfigure}{0.185\columnwidth} \centering
\includegraphics[trim={0cm 0cm 0cm 0cm},clip,width=1\columnwidth]{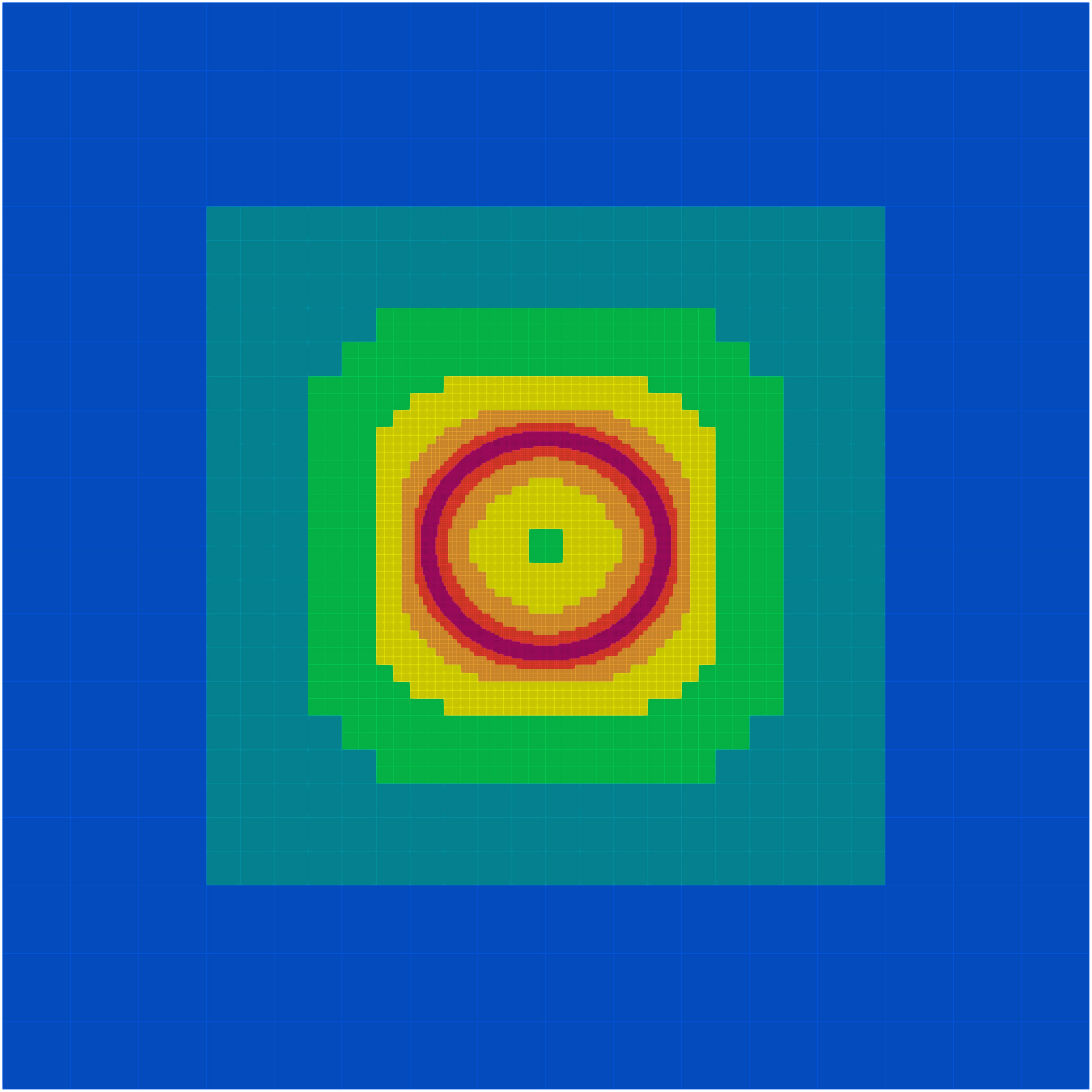}   \end{subfigure}
\begin{subfigure}{0.185\columnwidth} \centering
\includegraphics[trim={0cm 0cm 0cm 0cm},clip,width=1\columnwidth]{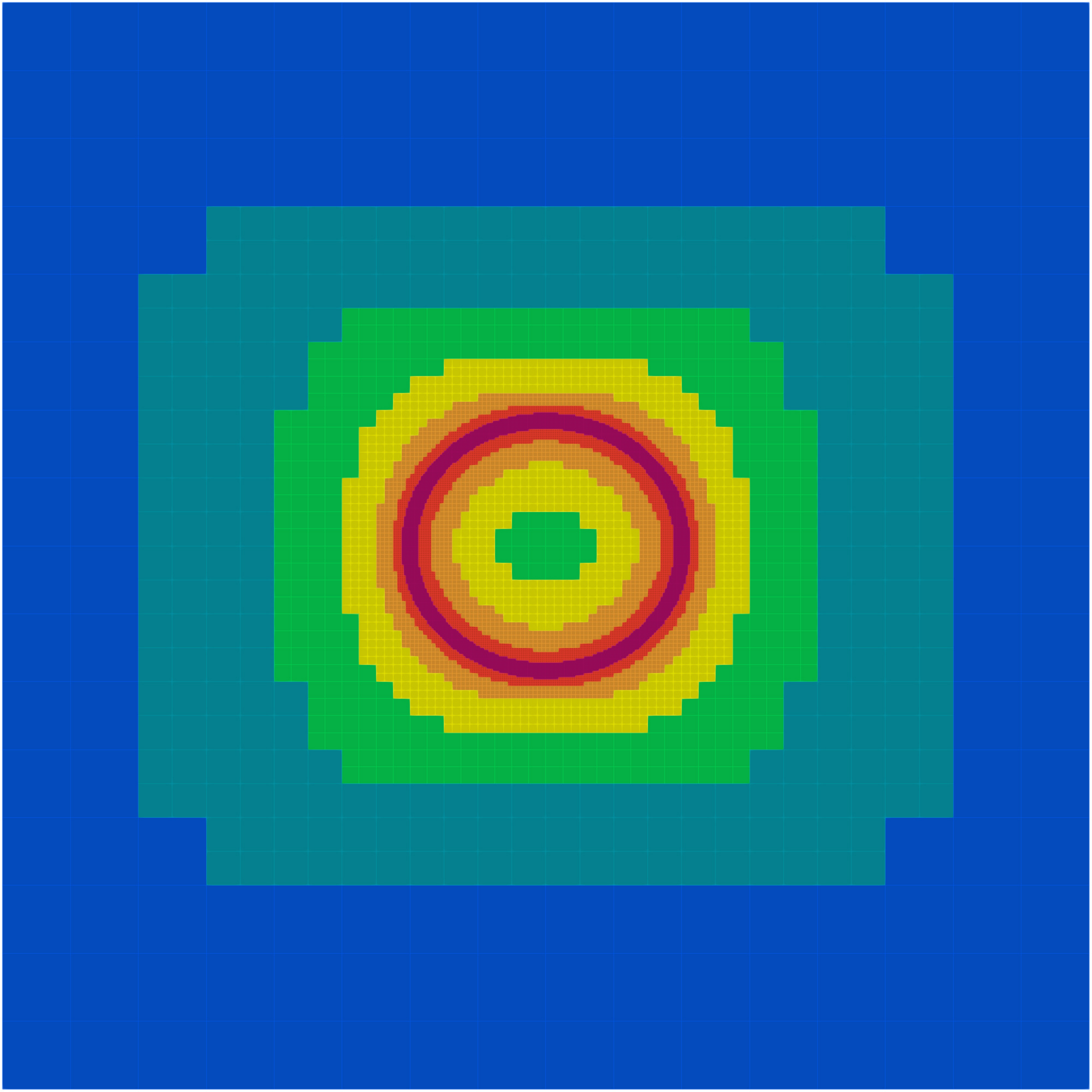} \end{subfigure}
\begin{subfigure}{0.185\columnwidth} \centering
\includegraphics[trim={0cm 0cm 0cm 0cm},clip,width=1\columnwidth]{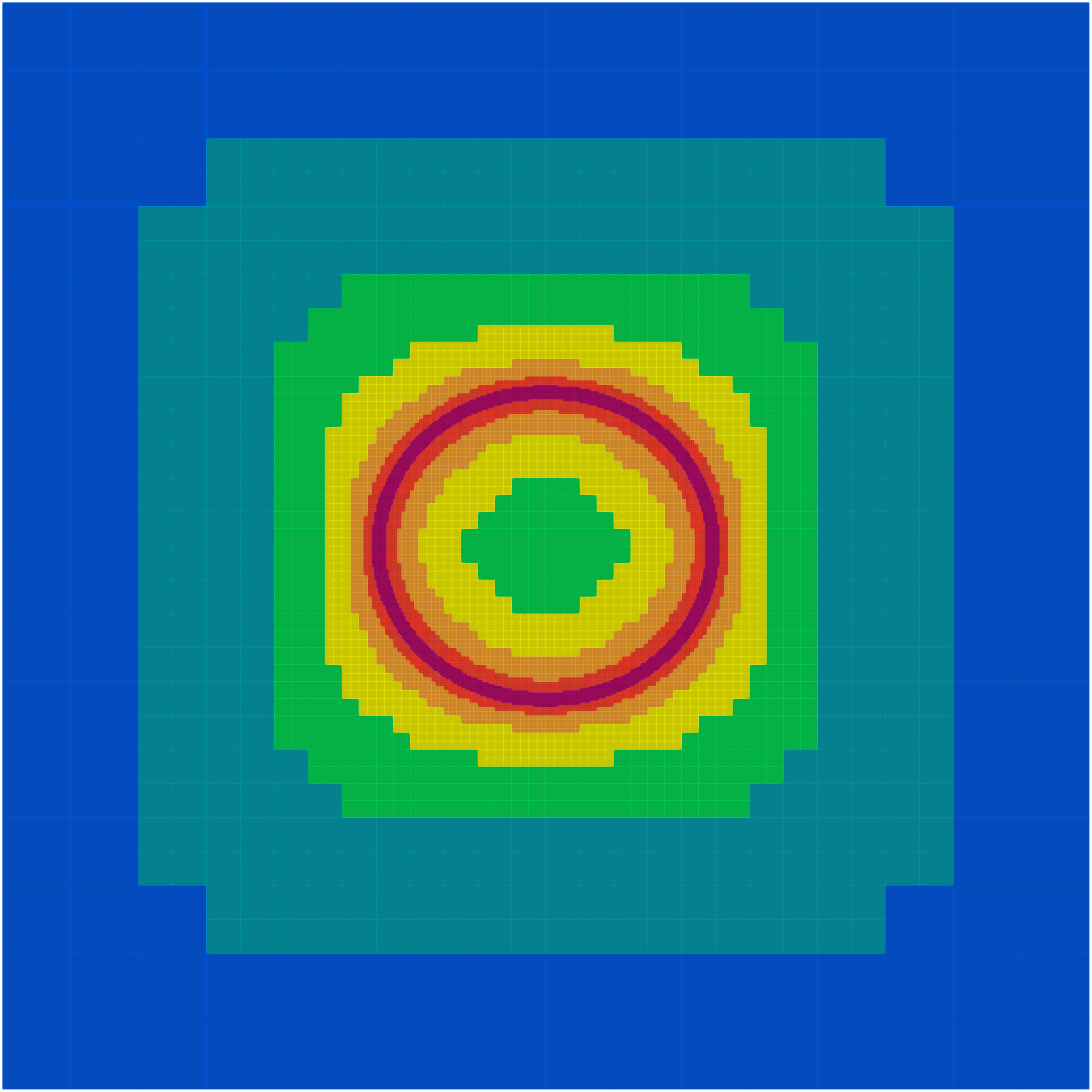} \end{subfigure}
\begin{subfigure}{0.185\columnwidth} \centering
\includegraphics[trim={0cm 0cm 0cm 0cm},clip,width=1\columnwidth]{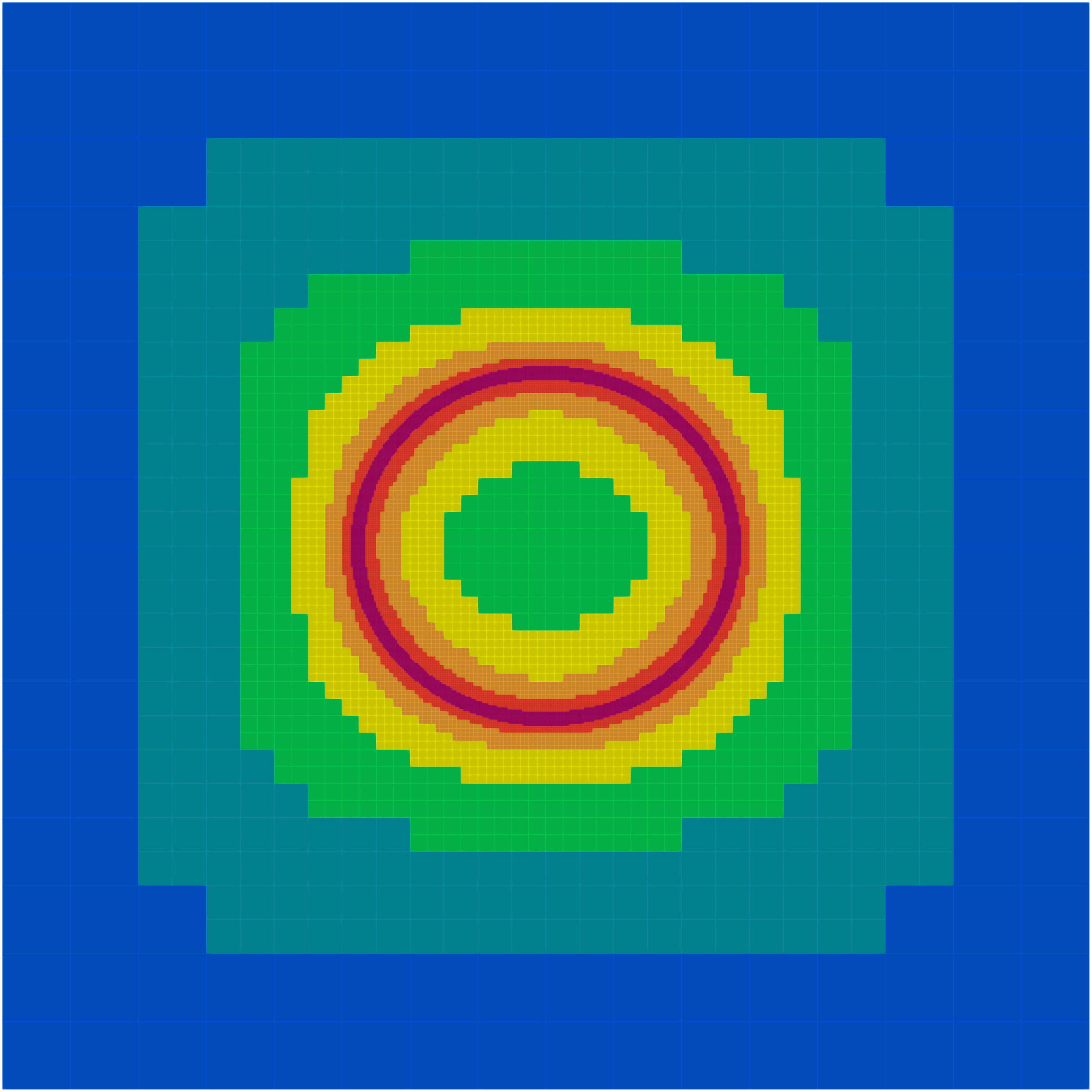}  \end{subfigure}
\begin{subfigure}{0.05\columnwidth} \centering
\includegraphics[trim={1cm 0cm 0cm 0cm},clip,width=0.84\columnwidth]{Images/Problem1_THBRes/ColorBar_Mesh_Lvl6-eps-converted-to.pdf} \end{subfigure}
\\ 
\begin{subfigure}{0.185\columnwidth} \centering
\includegraphics[trim={0cm 0cm 0cm 0cm},clip,width=1\columnwidth]{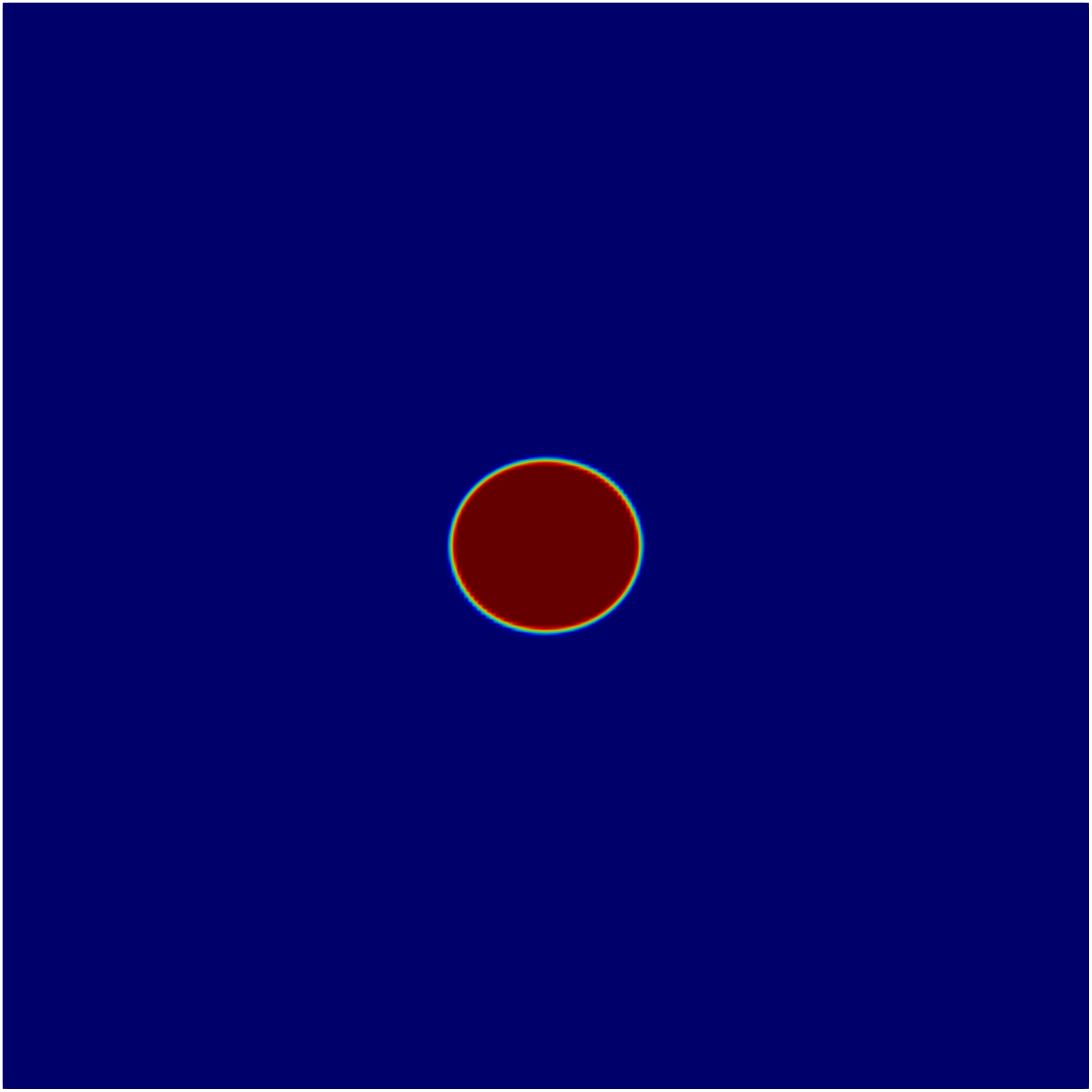}  \end{subfigure}
\begin{subfigure}{0.185\columnwidth} \centering
\includegraphics[trim={0cm 0cm 0cm 0cm},clip,width=1\columnwidth]{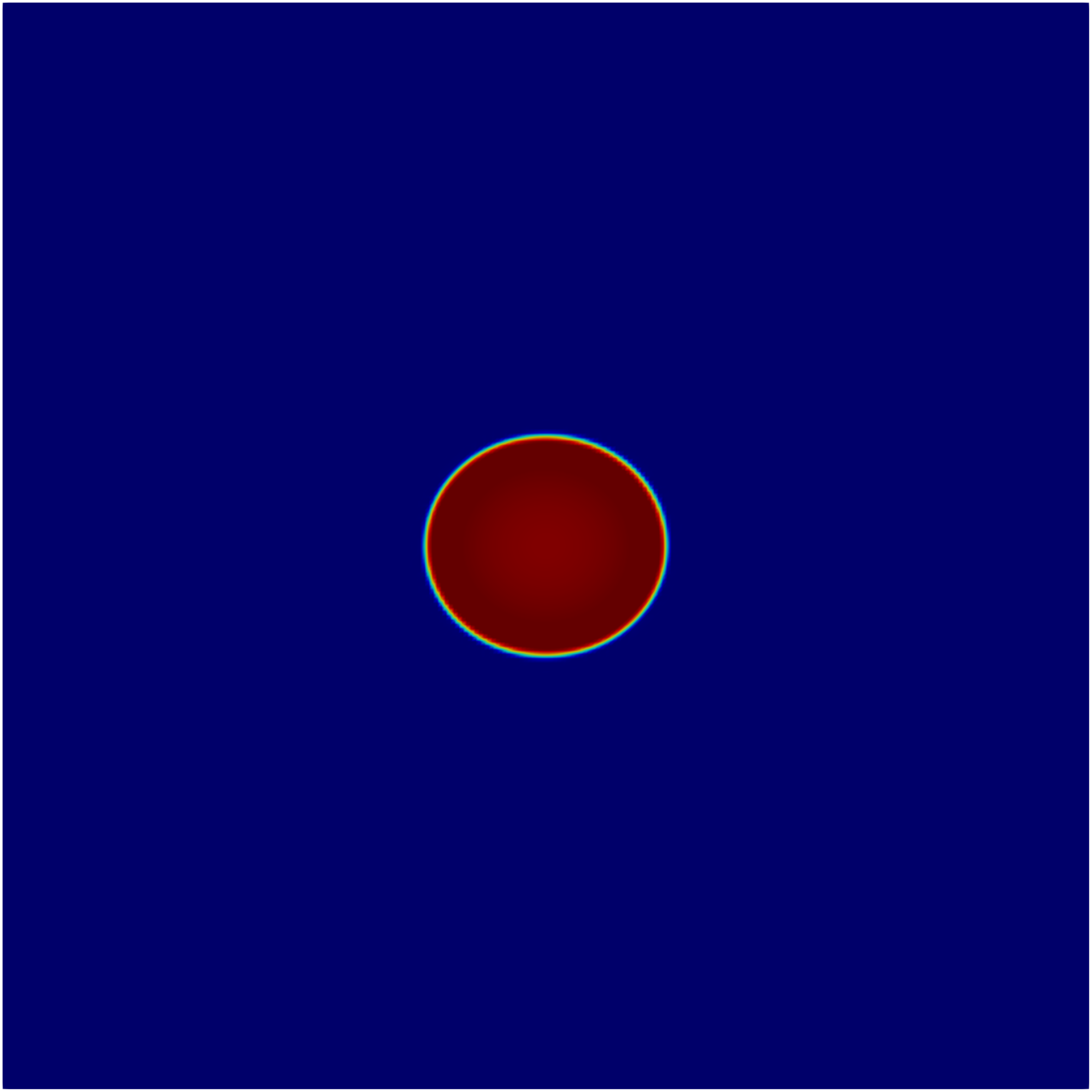}   \end{subfigure}
\begin{subfigure}{0.185\columnwidth} \centering
\includegraphics[trim={0cm 0cm 0cm 0cm},clip,width=1\columnwidth]{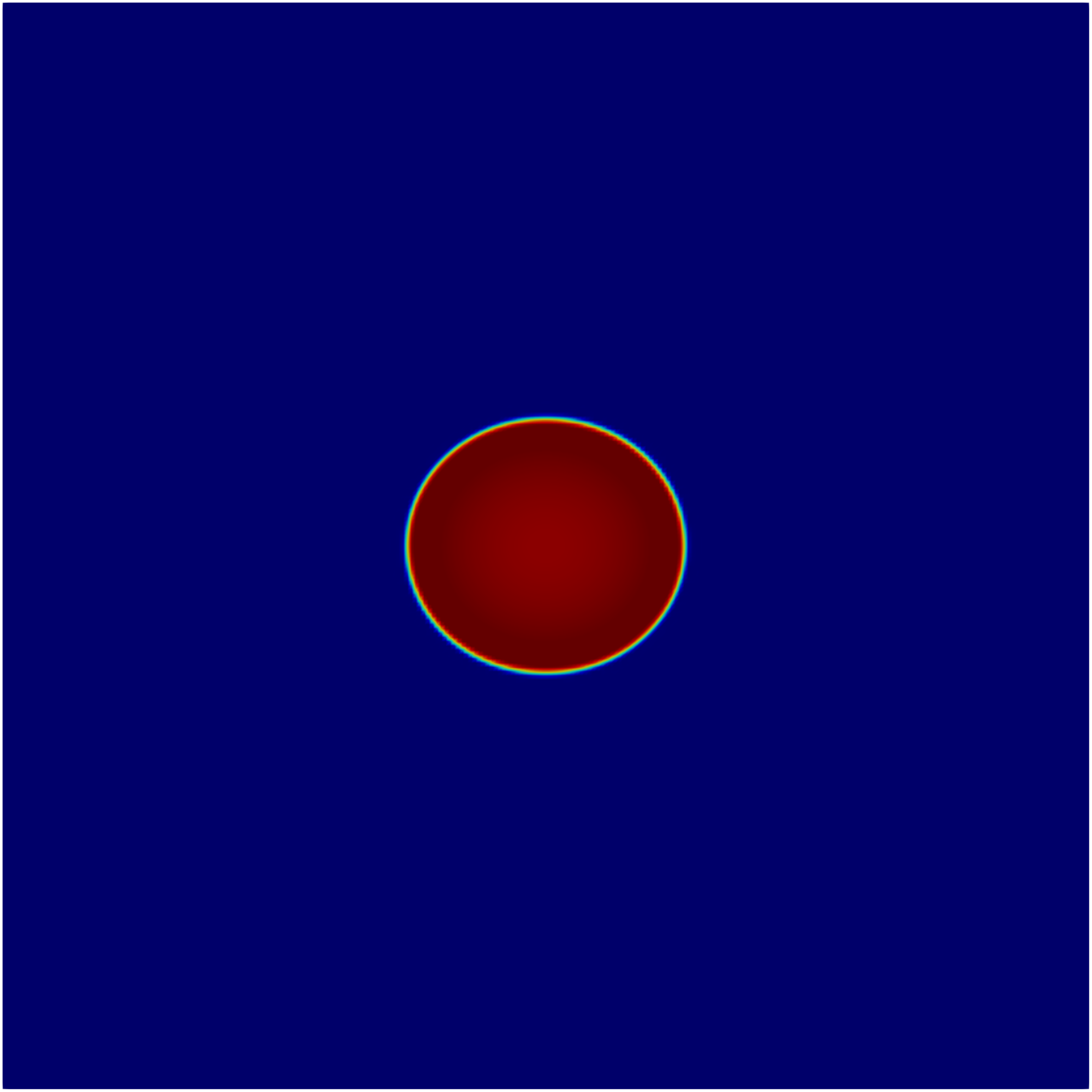} \end{subfigure}
\begin{subfigure}{0.185\columnwidth} \centering
\includegraphics[trim={0cm 0cm 0cm 0cm},clip,width=1\columnwidth]{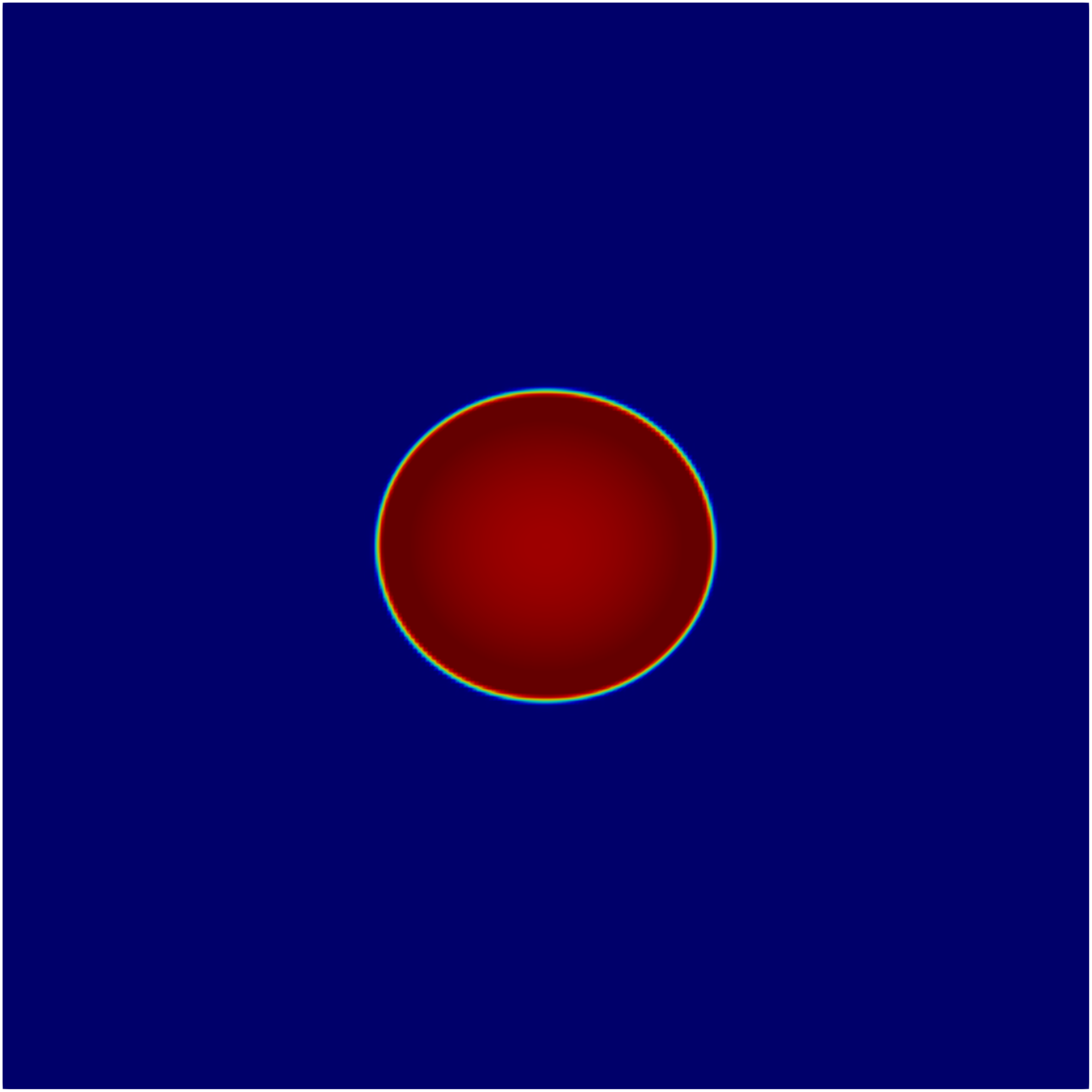} \end{subfigure}
\begin{subfigure}{0.185\columnwidth} \centering
\includegraphics[trim={0cm 0cm 0cm 0cm},clip,width=1\columnwidth]{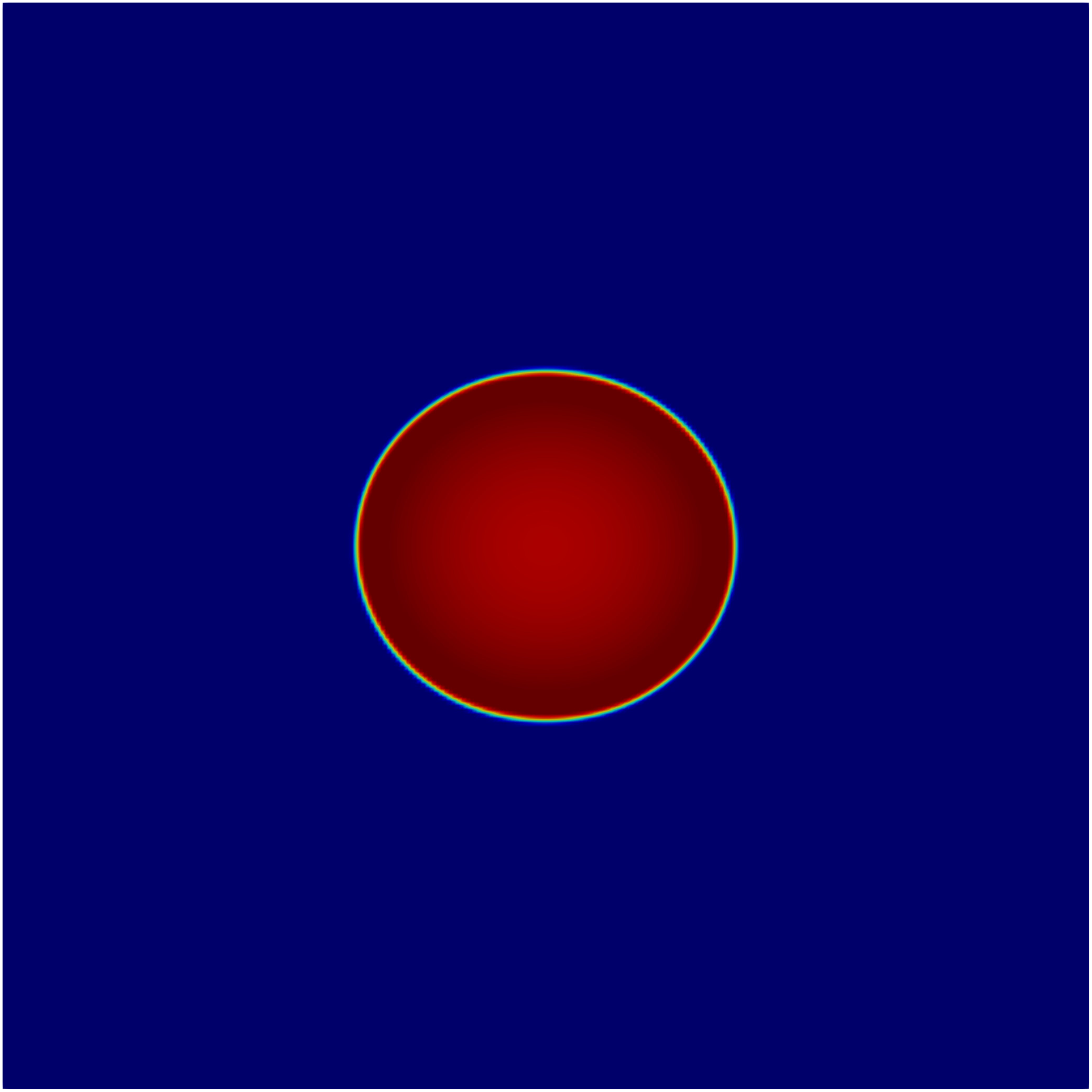}  \end{subfigure}
\begin{subfigure}{0.05\columnwidth} \centering
\includegraphics[trim={0cm 1cm 0cm 0cm},clip,width=1.18\columnwidth]{Images/Problem1_2DSq_UniRef/ColoarBar_Phi-eps-converted-to.pdf} \end{subfigure}
\\
\begin{subfigure}{0.185\columnwidth} \centering
\includegraphics[trim={0cm 0cm 0cm 0cm},clip,width=1\columnwidth]{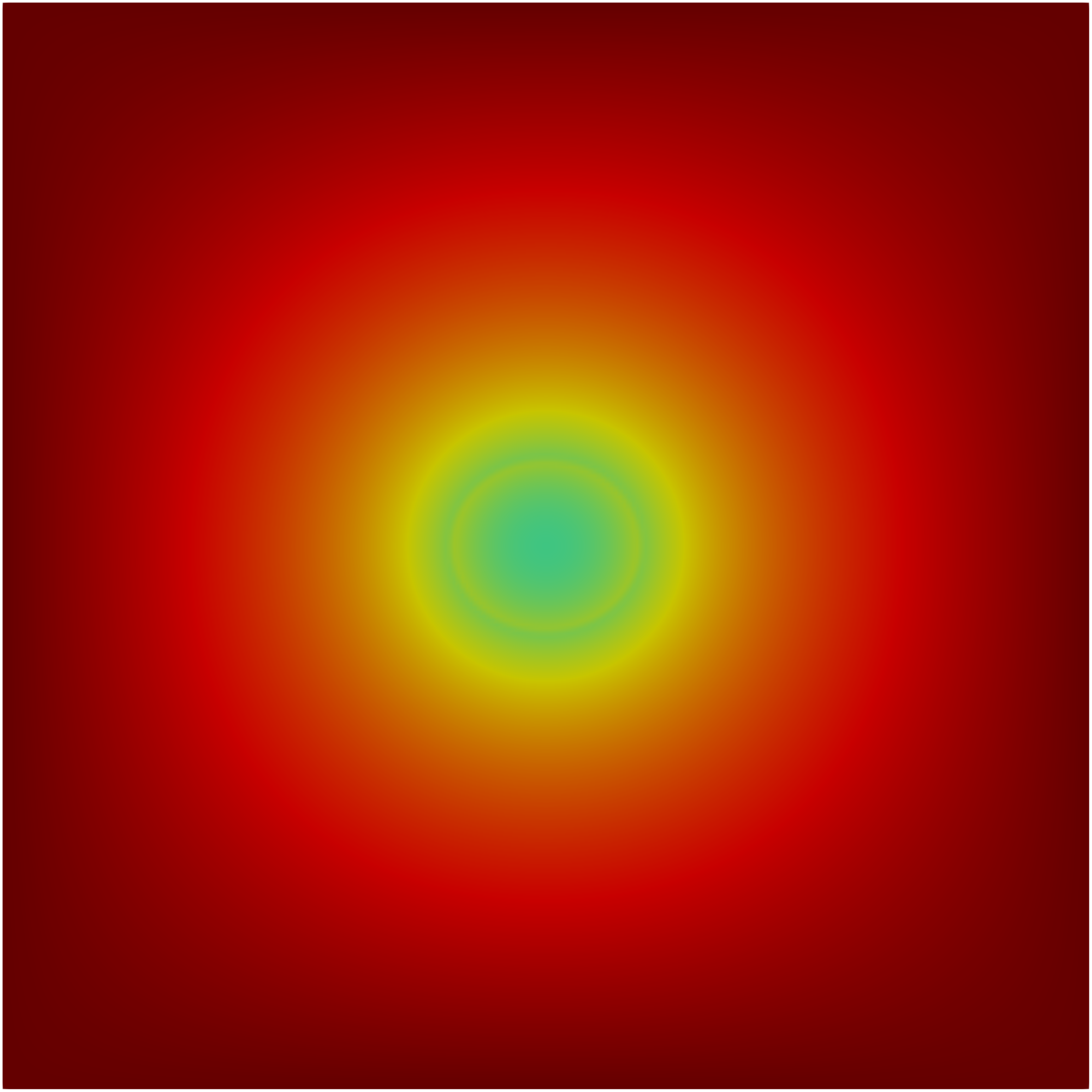} \caption{$t = 0$} \end{subfigure}
\begin{subfigure}{0.185\columnwidth} \centering
\includegraphics[trim={0cm 0cm 0cm 0cm},clip,width=1\columnwidth]{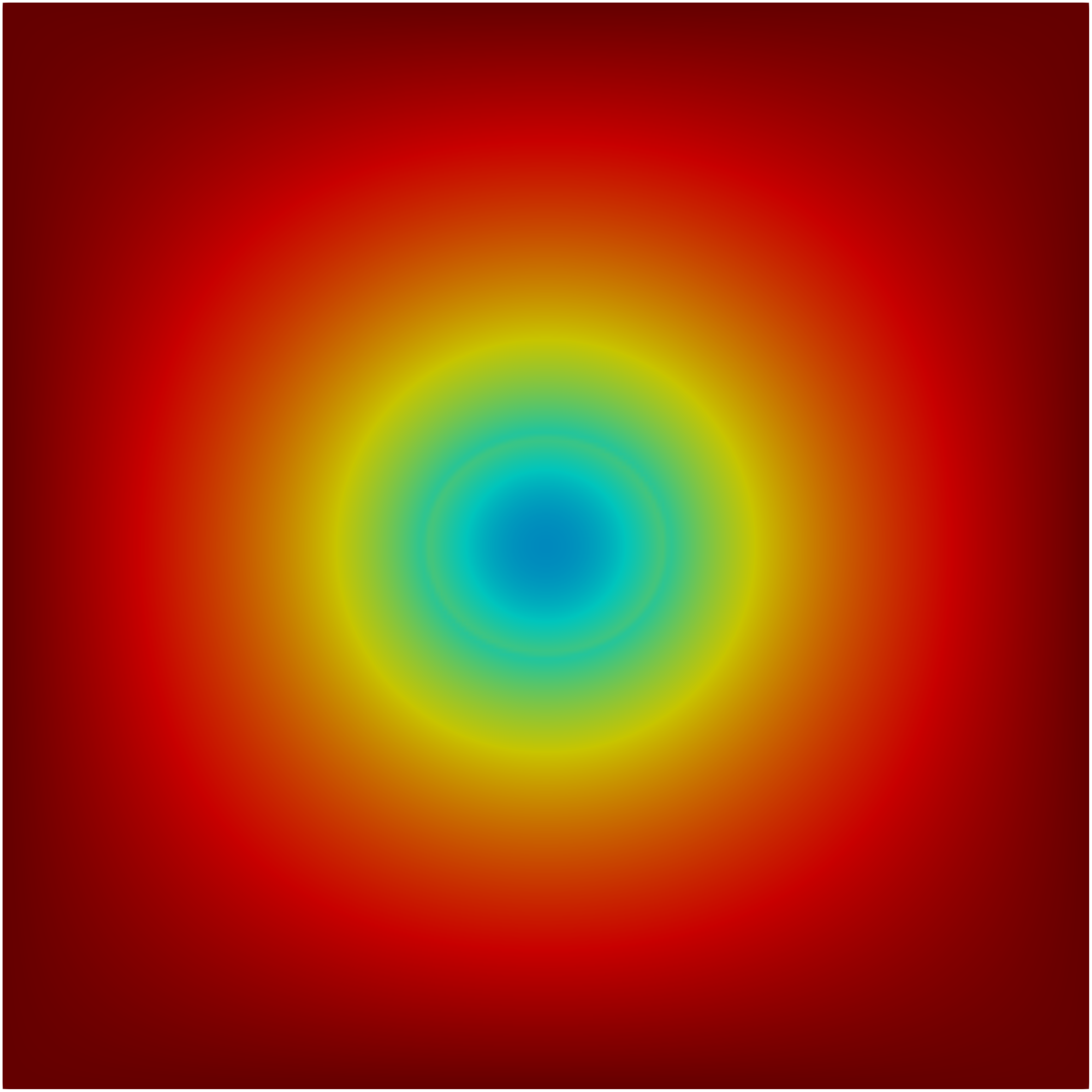} \caption{$t = 0.15$}  \end{subfigure}
\begin{subfigure}{0.185\columnwidth} \centering
\includegraphics[trim={0cm 0cm 0cm 0cm},clip,width=1\columnwidth]{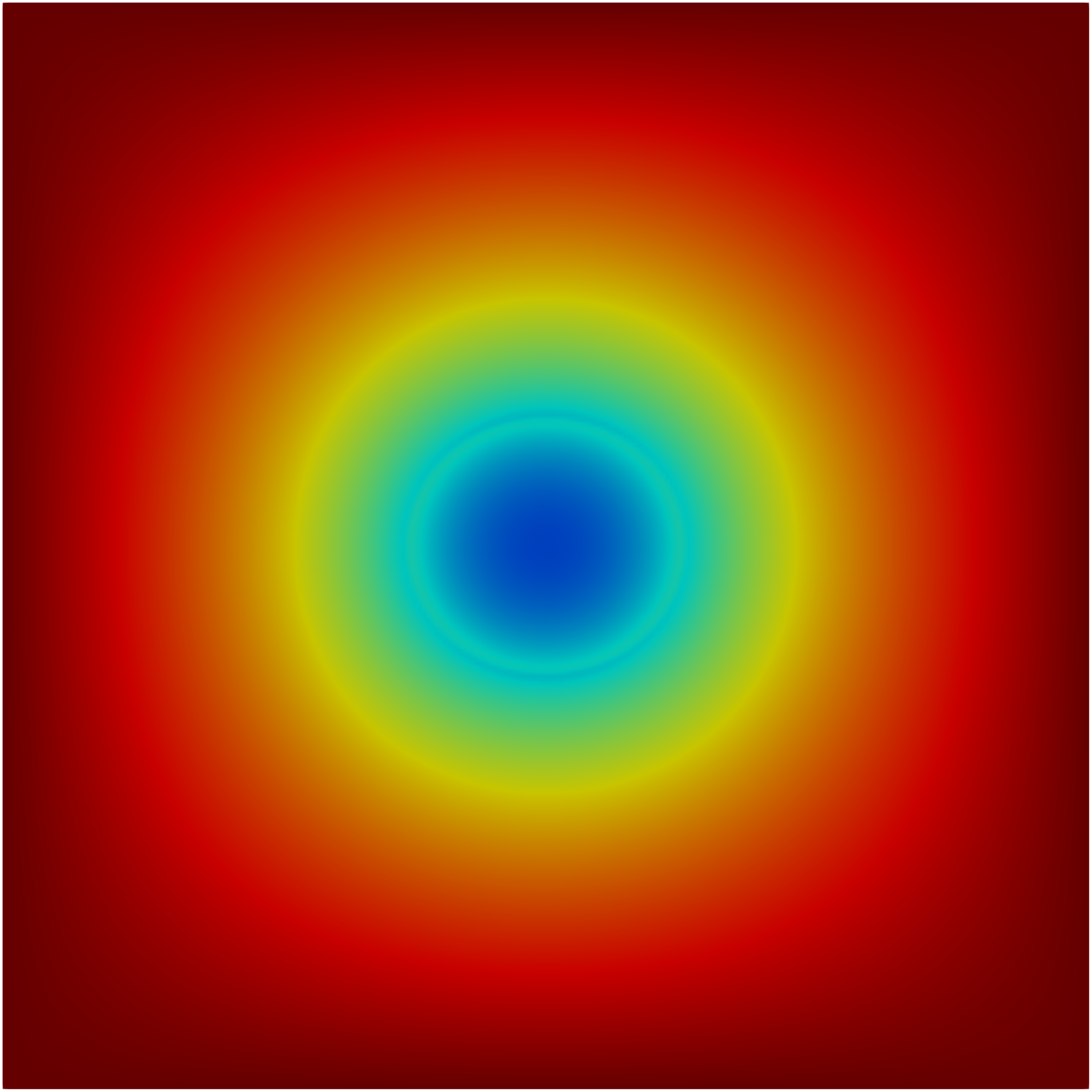} \caption{$t = 0.25$}\end{subfigure}
\begin{subfigure}{0.185\columnwidth} \centering
\includegraphics[trim={0cm 0cm 0cm 0cm},clip,width=1\columnwidth]{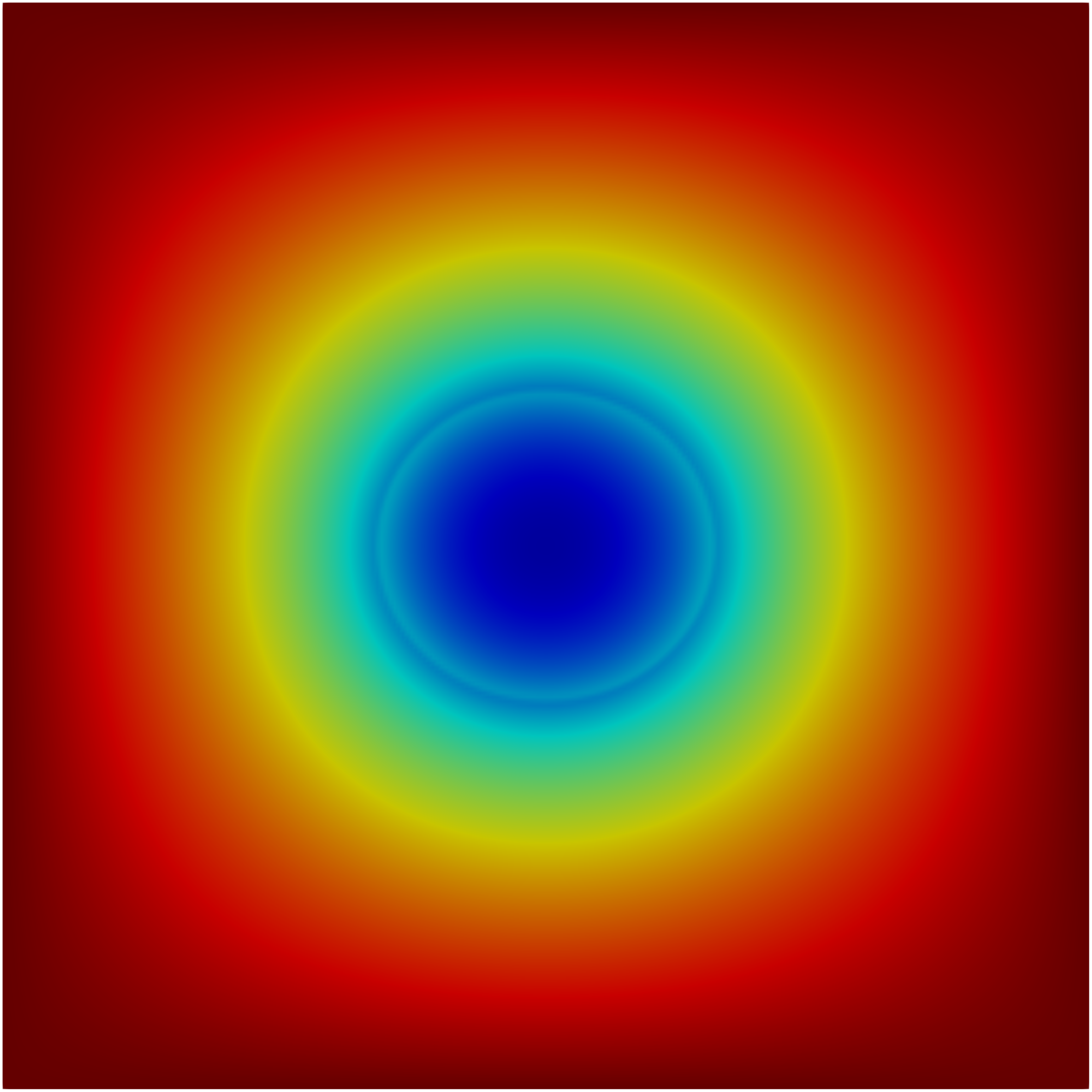} \caption{$t = 0.4$} \end{subfigure}
\begin{subfigure}{0.185\columnwidth} \centering
\includegraphics[trim={0cm 0cm 0cm 0cm},clip,width=1\columnwidth]{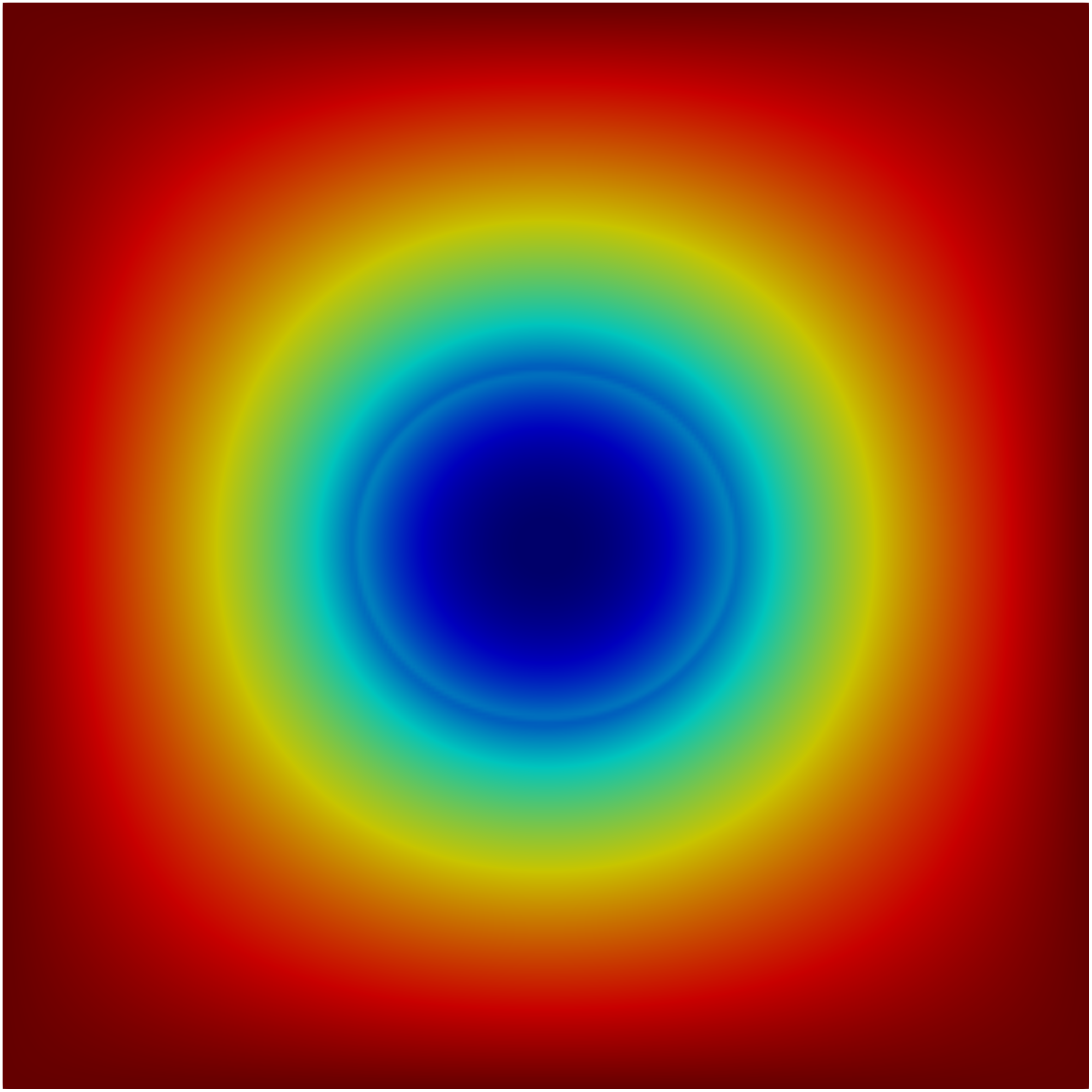} \caption{$t = 0.5$}\end{subfigure}
\begin{subfigure}{0.05\columnwidth} \centering
\includegraphics[trim={0cm -6cm 0cm 0cm},clip,width=1.08\columnwidth]{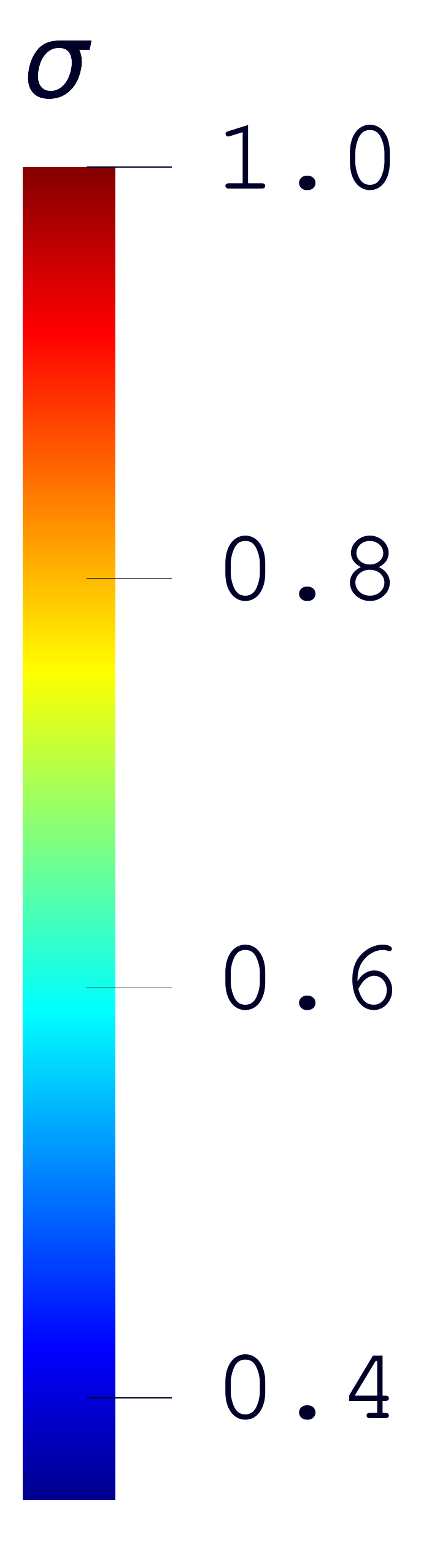} \end{subfigure}
\caption{Evolution of an elliptical tumor in square domain with $\mathcal{P} = 5$, while keeping all other parameters unchanged. Top: adaptive THB-spline mesh configurations of degree $3$, $\ell = 6$, $m=2$, $\alpha = 0.01$, $\beta = 0.0001$ with a finest refinement level at $2^{10} \times 2^{10}$ ($h_e = 6/1024$) mesh resolution. Middle: evolution of tumor geometry. Bottom: corresponding nutrient concentration.}
\label{figure_8_final}
\end{figure}

\begin{figure}[!t]\centering
\begin{subfigure}{0.185\columnwidth} \centering
\includegraphics[trim={0cm 0cm 0cm 0cm},clip,width=1\columnwidth]{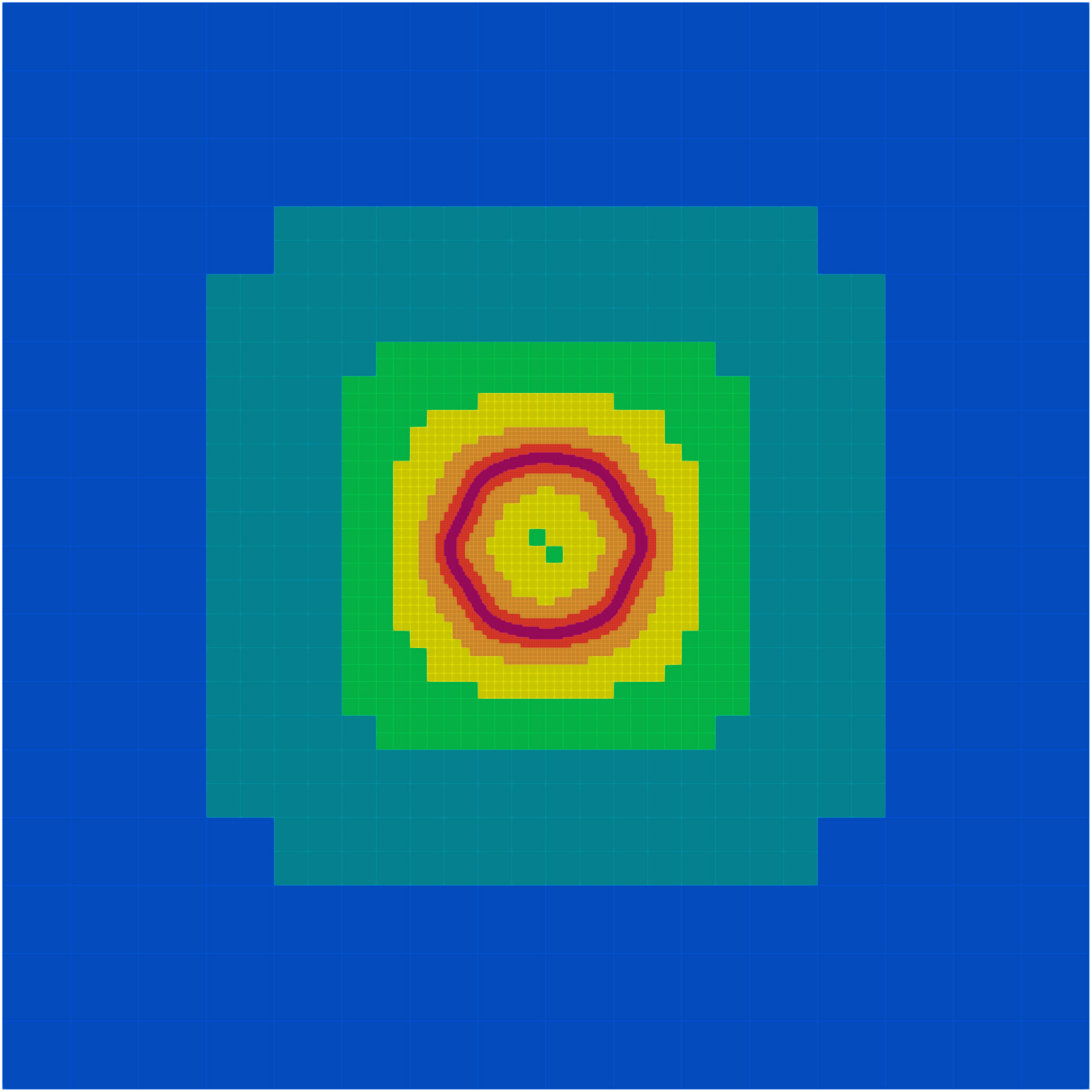}  \end{subfigure}
\begin{subfigure}{0.185\columnwidth} \centering
\includegraphics[trim={0cm 0cm 0cm 0cm},clip,width=1\columnwidth]{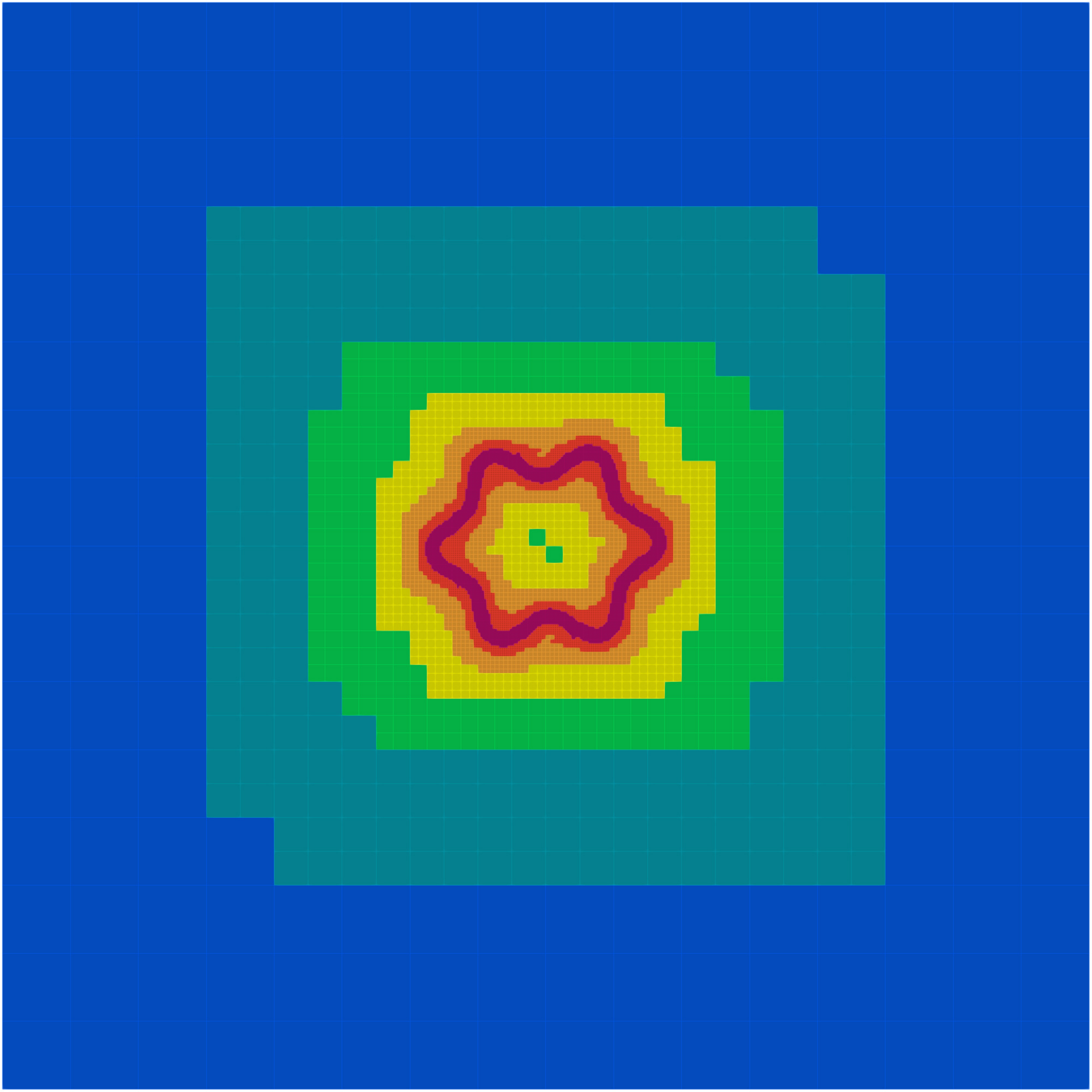}   \end{subfigure}
\begin{subfigure}{0.185\columnwidth} \centering
\includegraphics[trim={0cm 0cm 0cm 0cm},clip,width=1\columnwidth]{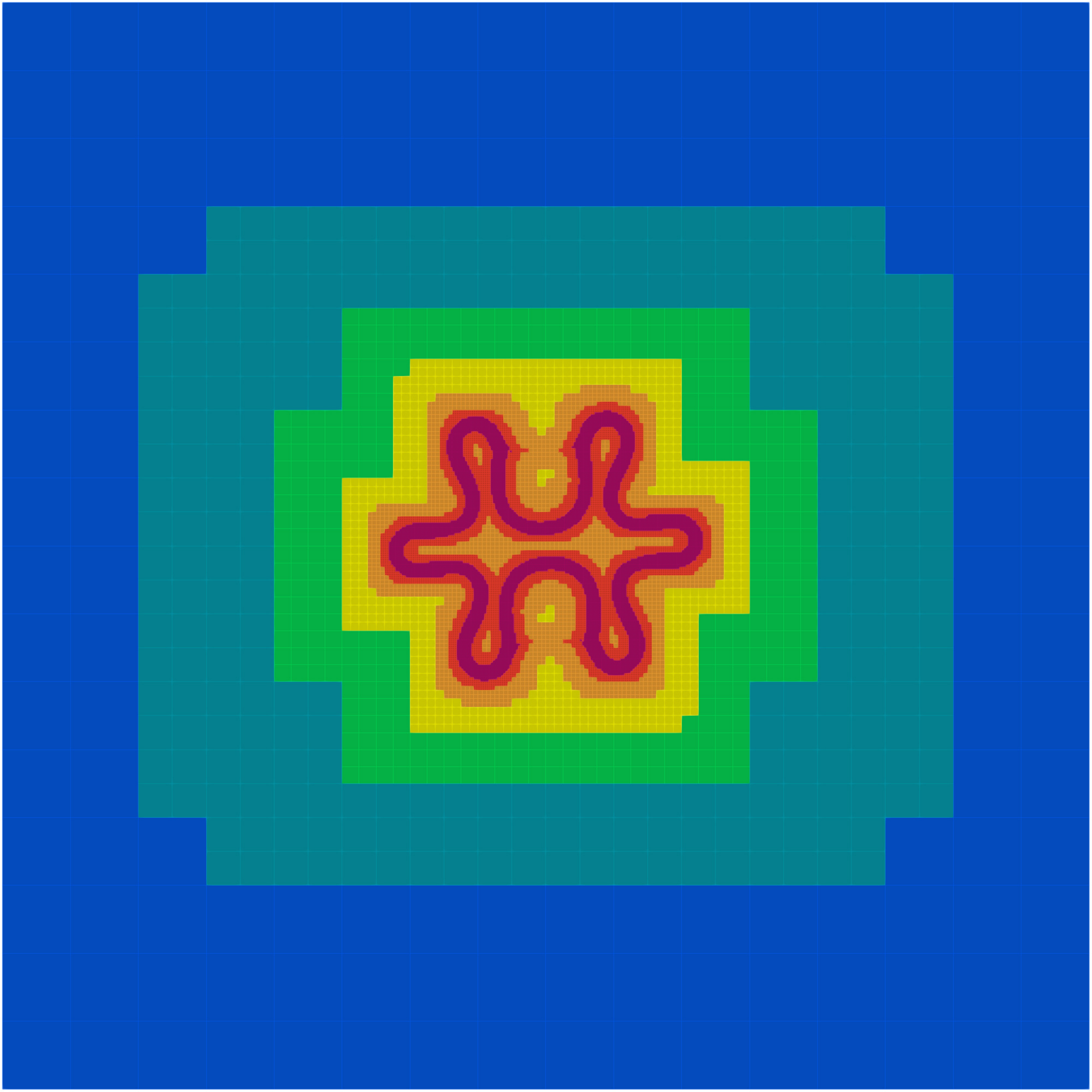} \end{subfigure}
\begin{subfigure}{0.185\columnwidth} \centering
\includegraphics[trim={0cm 0cm 0cm 0cm},clip,width=1\columnwidth]{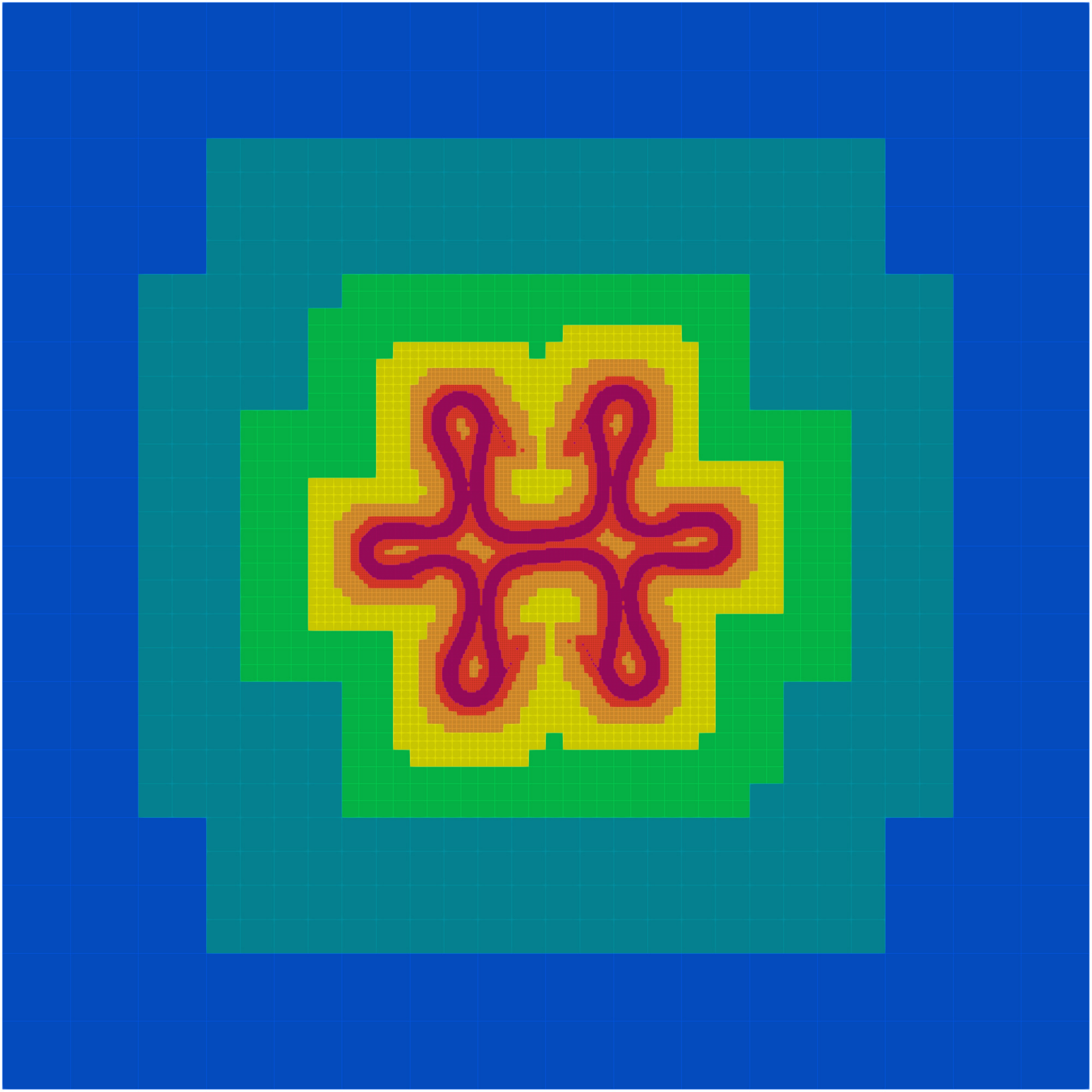} \end{subfigure}
\begin{subfigure}{0.185\columnwidth} \centering
\includegraphics[trim={0cm 0cm 0cm 0cm},clip,width=1\columnwidth]{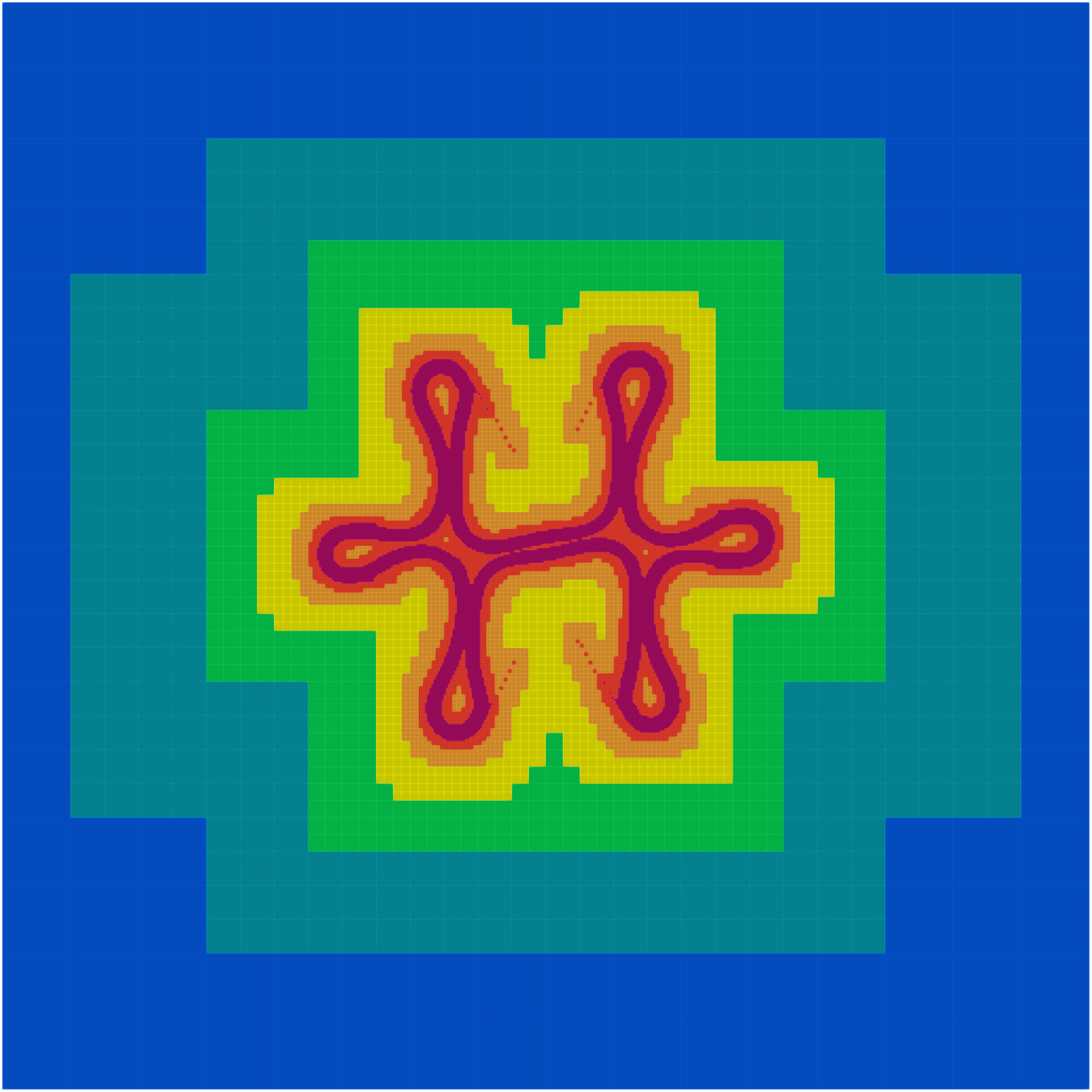}  \end{subfigure}
\begin{subfigure}{0.05\columnwidth} \centering
\includegraphics[trim={1cm 0cm 0cm 0cm},clip,width=0.84\columnwidth]{Images/Problem1_THBRes/ColorBar_Mesh_Lvl6-eps-converted-to.pdf} \end{subfigure}
\\ 
\begin{subfigure}{0.185\columnwidth} \centering
\includegraphics[trim={0cm 0cm 0cm 0cm},clip,width=1\columnwidth]{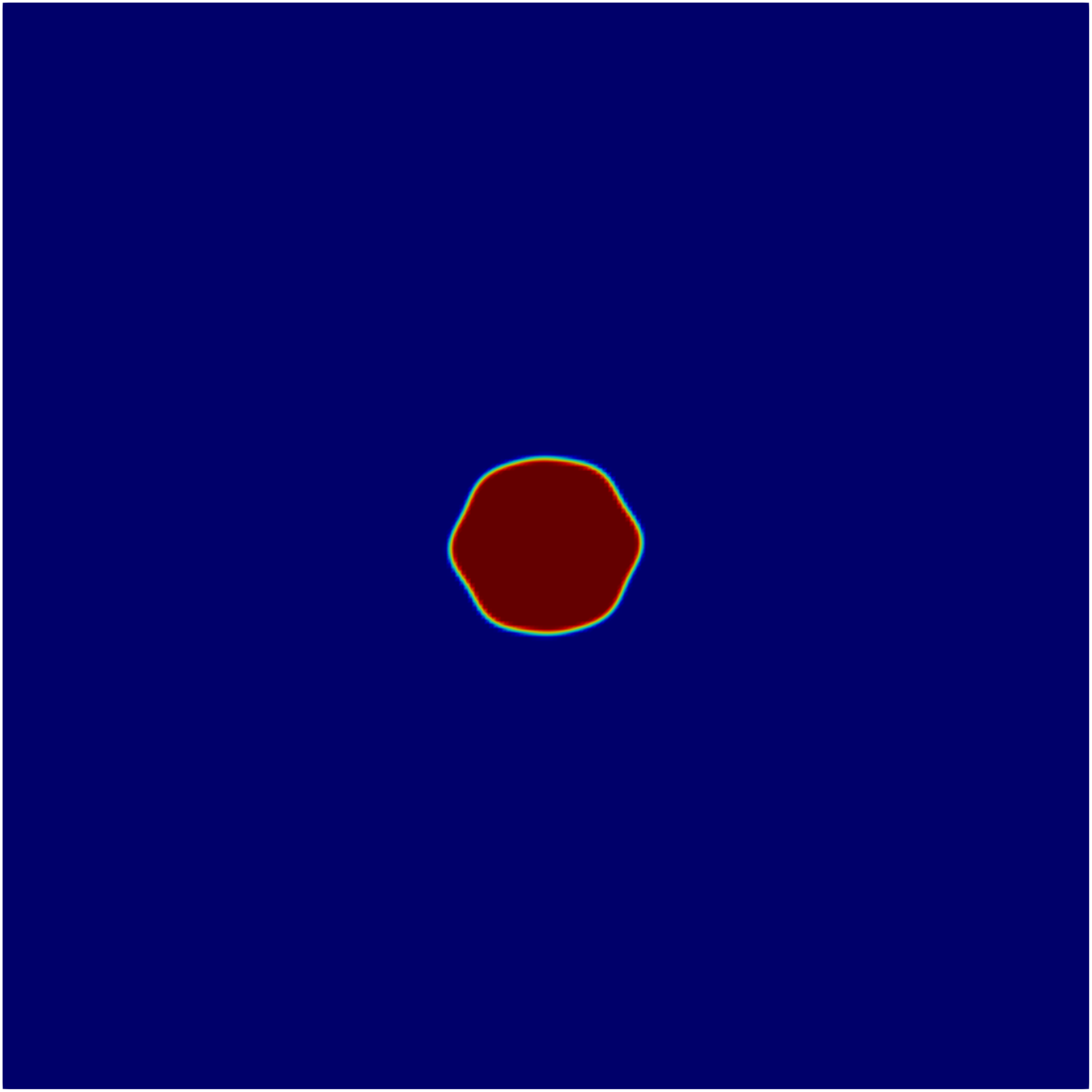}  \end{subfigure}
\begin{subfigure}{0.185\columnwidth} \centering
\includegraphics[trim={0cm 0cm 0cm 0cm},clip,width=1\columnwidth]{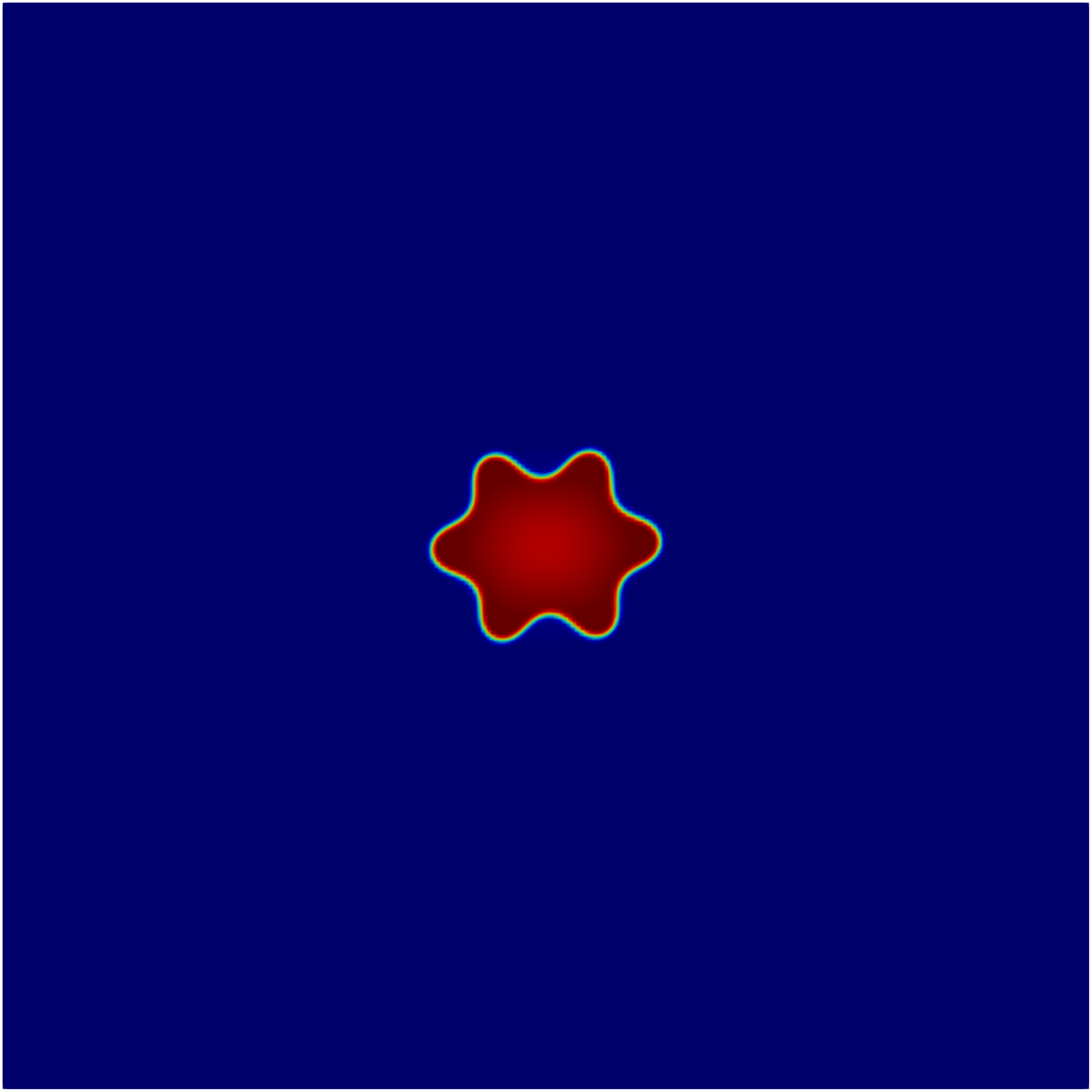}   \end{subfigure}
\begin{subfigure}{0.185\columnwidth} \centering
\includegraphics[trim={0cm 0cm 0cm 0cm},clip,width=1\columnwidth]{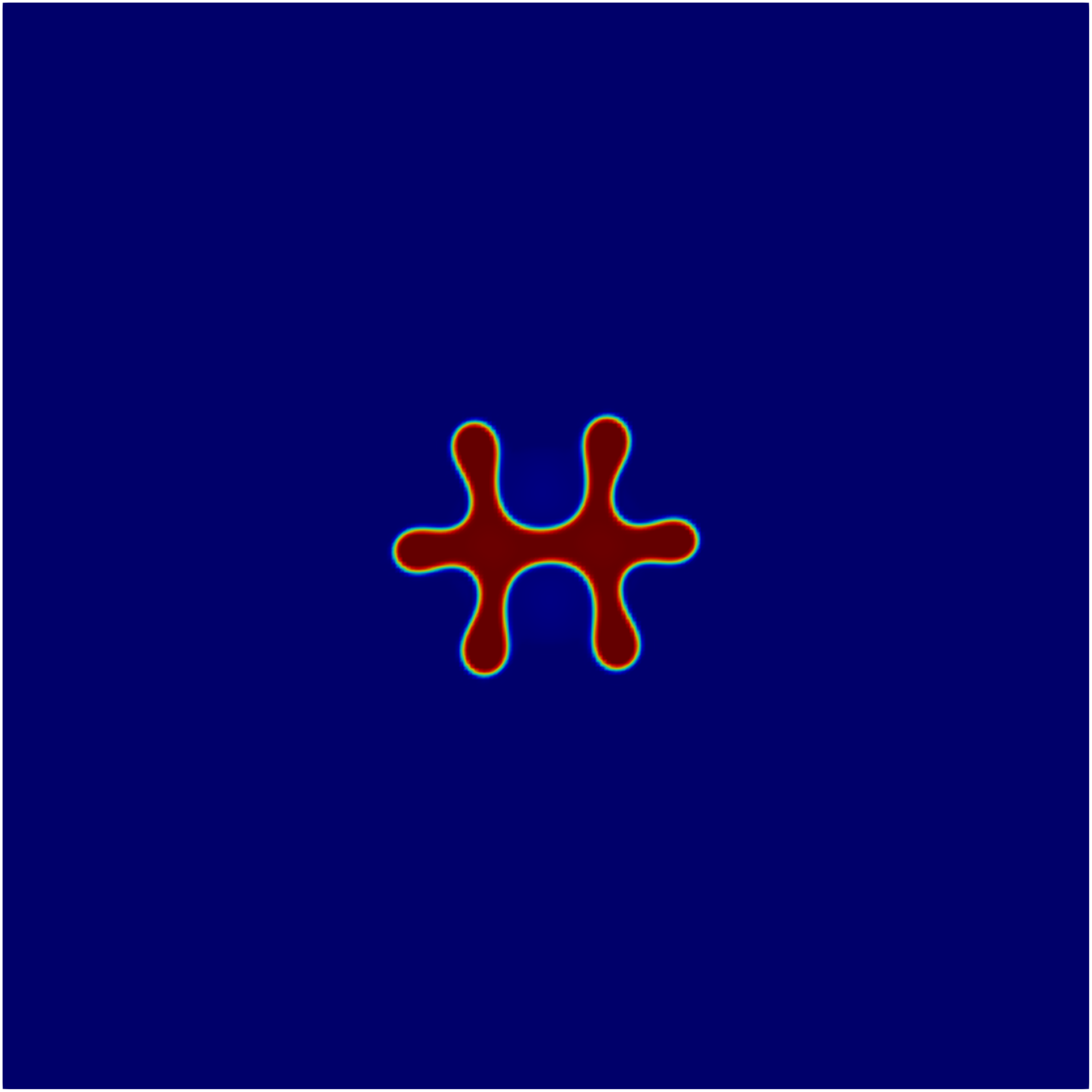} \end{subfigure}
\begin{subfigure}{0.185\columnwidth} \centering
\includegraphics[trim={0cm 0cm 0cm 0cm},clip,width=1\columnwidth]{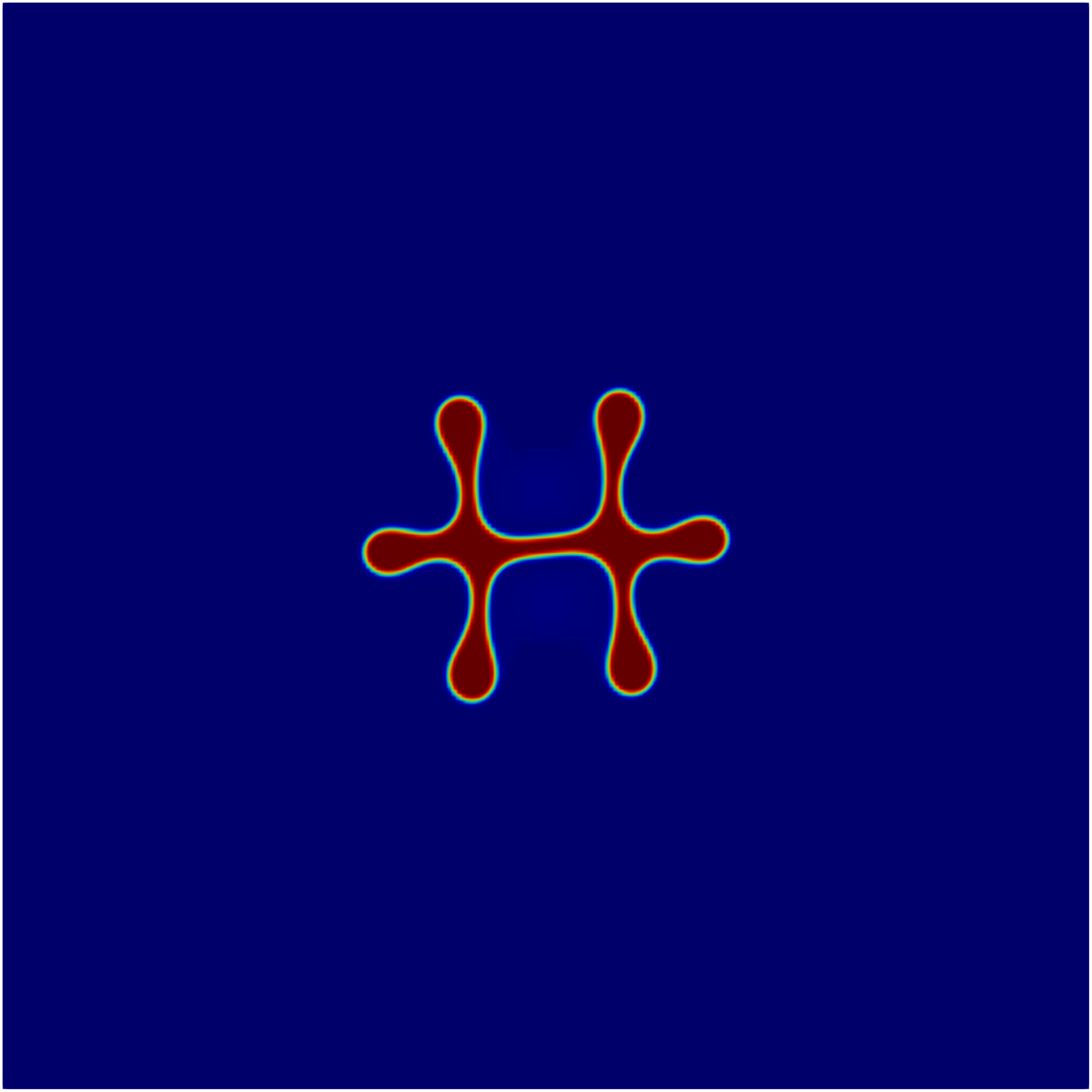} \end{subfigure}
\begin{subfigure}{0.185\columnwidth} \centering
\includegraphics[trim={0cm 0cm 0cm 0cm},clip,width=1\columnwidth]{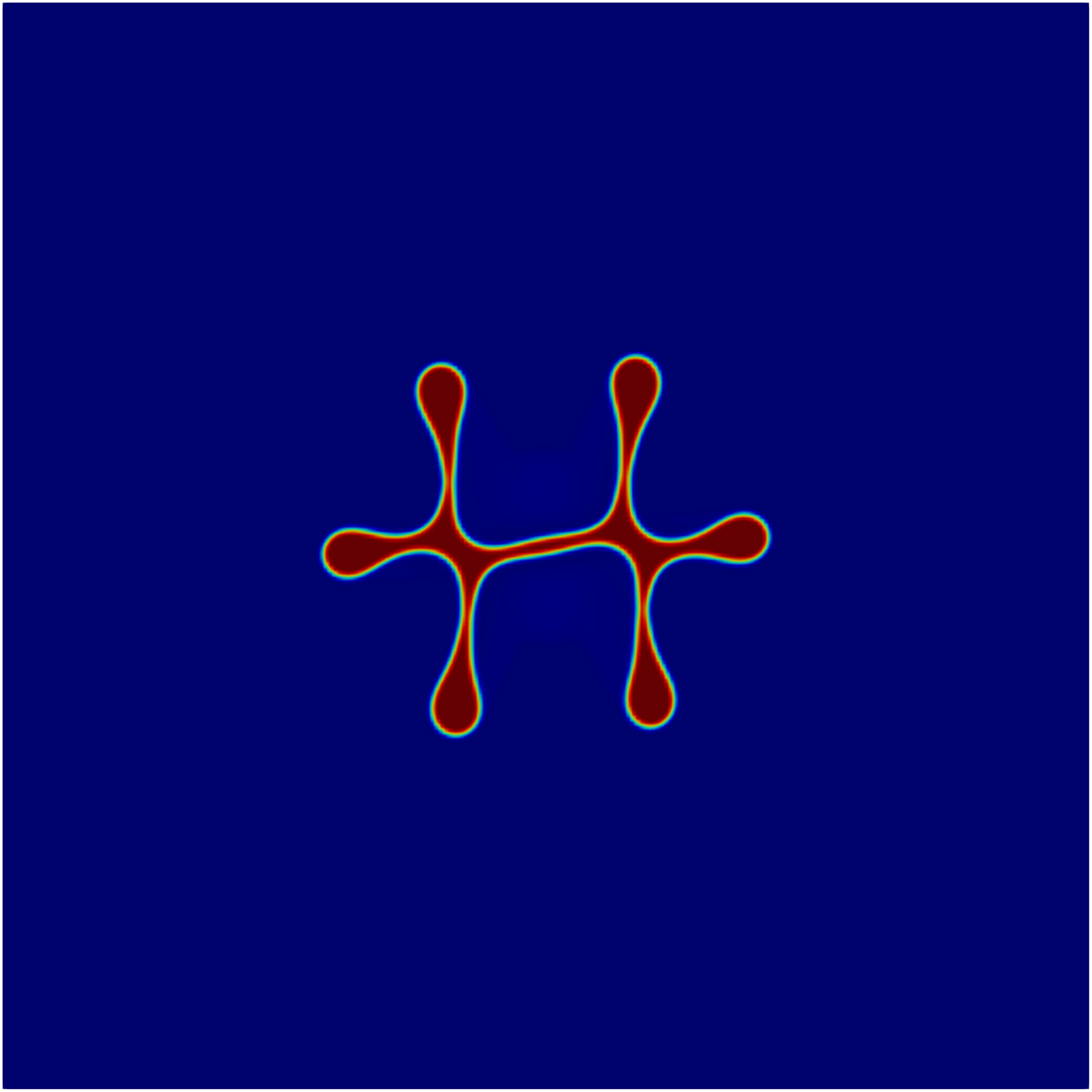}  \end{subfigure}
\begin{subfigure}{0.05\columnwidth} \centering
\includegraphics[trim={0cm 1cm 0cm 0cm},clip,width=1.2\columnwidth]{Images/Problem1_2DSq_UniRef/ColoarBar_Phi-eps-converted-to.pdf} \end{subfigure}
\\
\begin{subfigure}{0.185\columnwidth} \centering
\includegraphics[trim={0cm 0cm 0cm 0cm},clip,width=1\columnwidth]{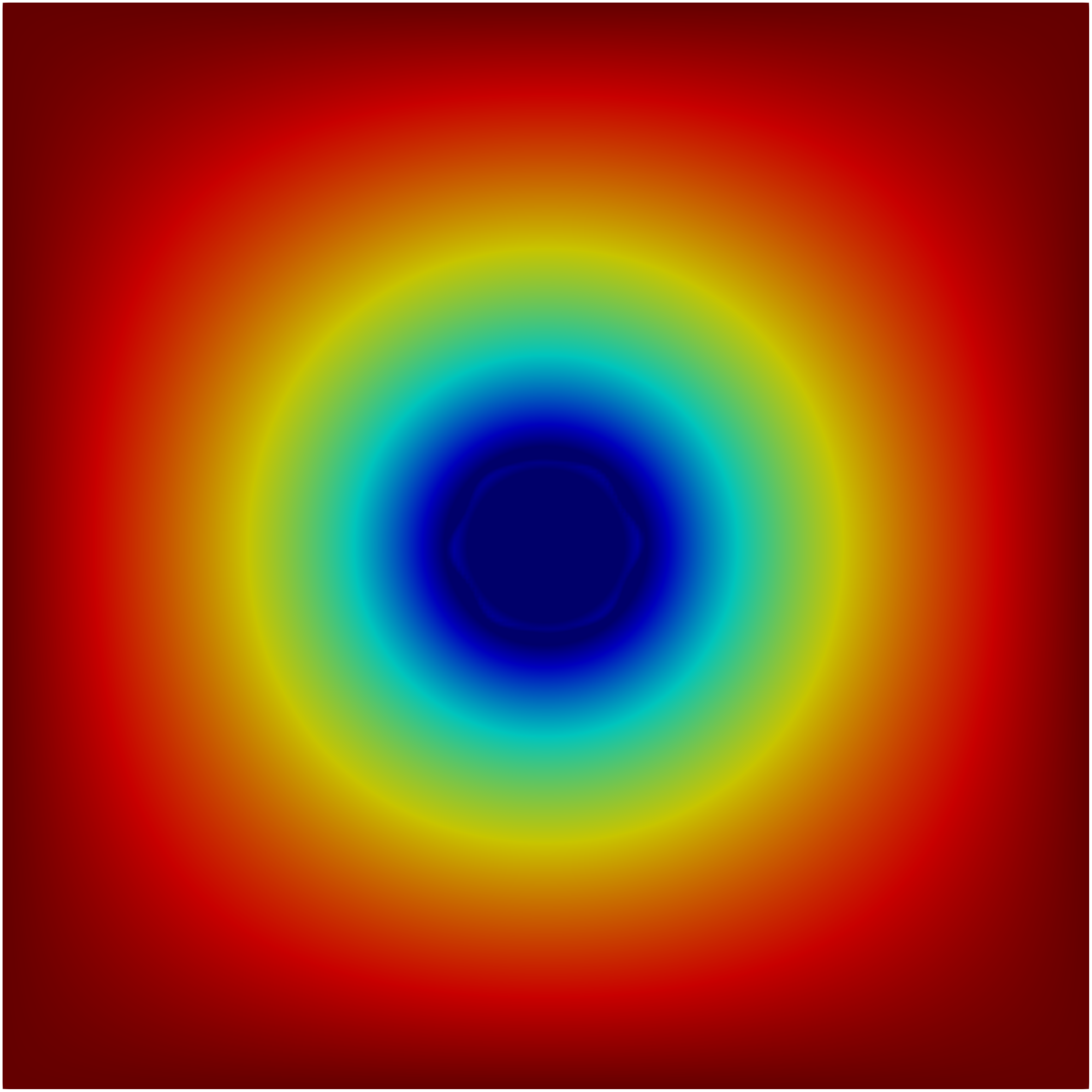} \caption{$t = 0$} \end{subfigure}
\begin{subfigure}{0.185\columnwidth} \centering
\includegraphics[trim={0cm 0cm 0cm 0cm},clip,width=1\columnwidth]{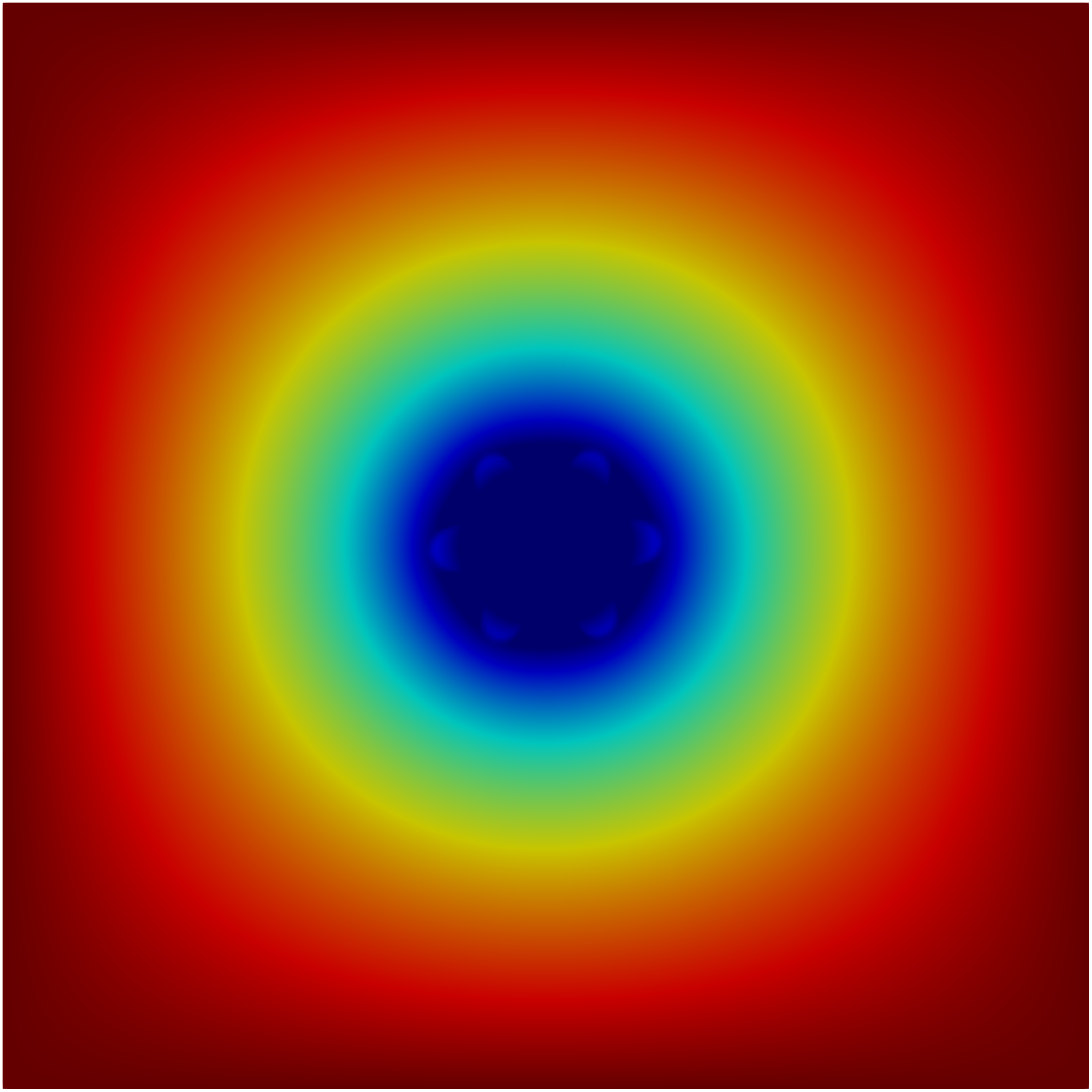} \caption{$t = 0.6$}  \end{subfigure}
\begin{subfigure}{0.185\columnwidth} \centering
\includegraphics[trim={0cm 0cm 0cm 0cm},clip,width=1\columnwidth]{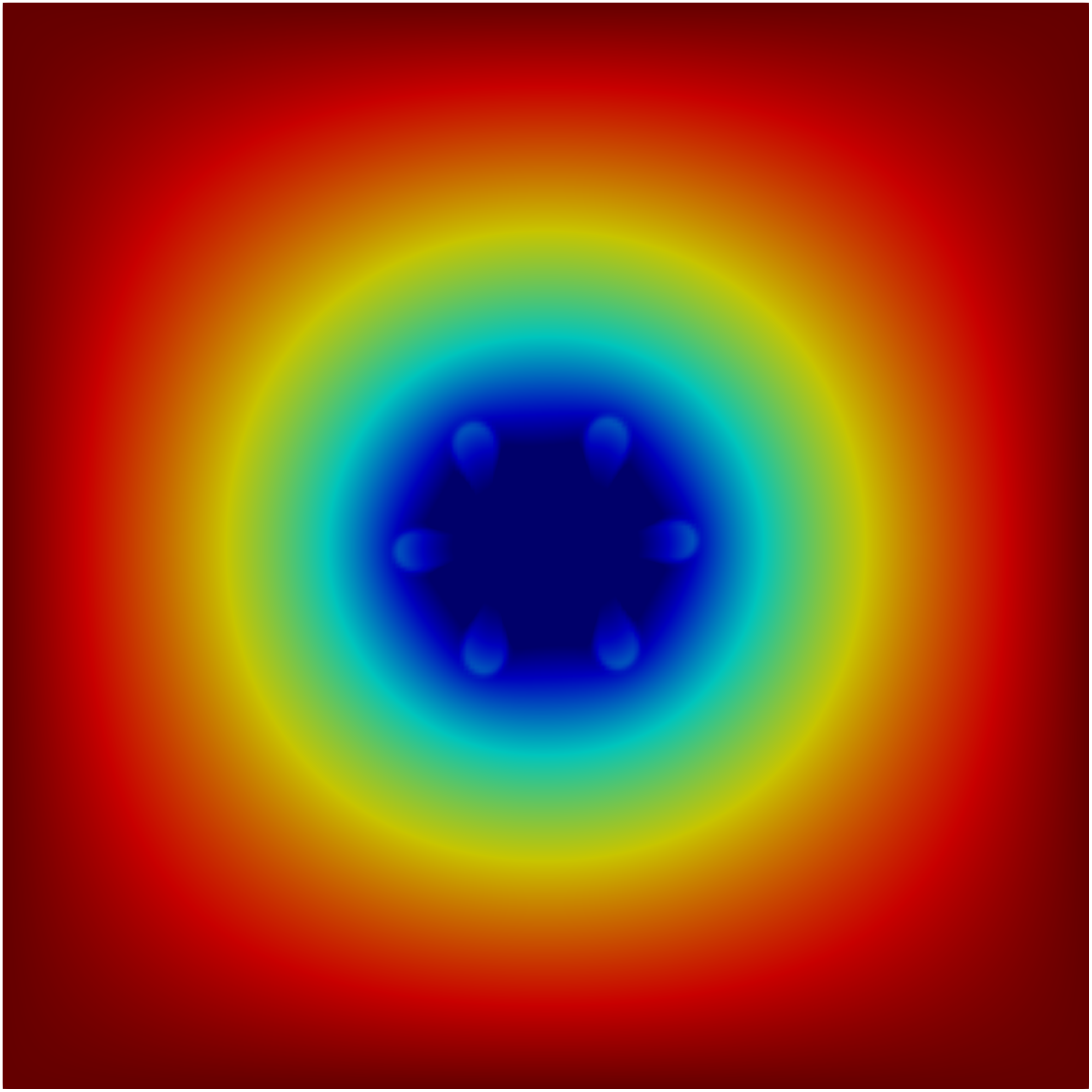} \caption{$t = 1$}\end{subfigure}
\begin{subfigure}{0.185\columnwidth} \centering
\includegraphics[trim={0cm 0cm 0cm 0cm},clip,width=1\columnwidth]{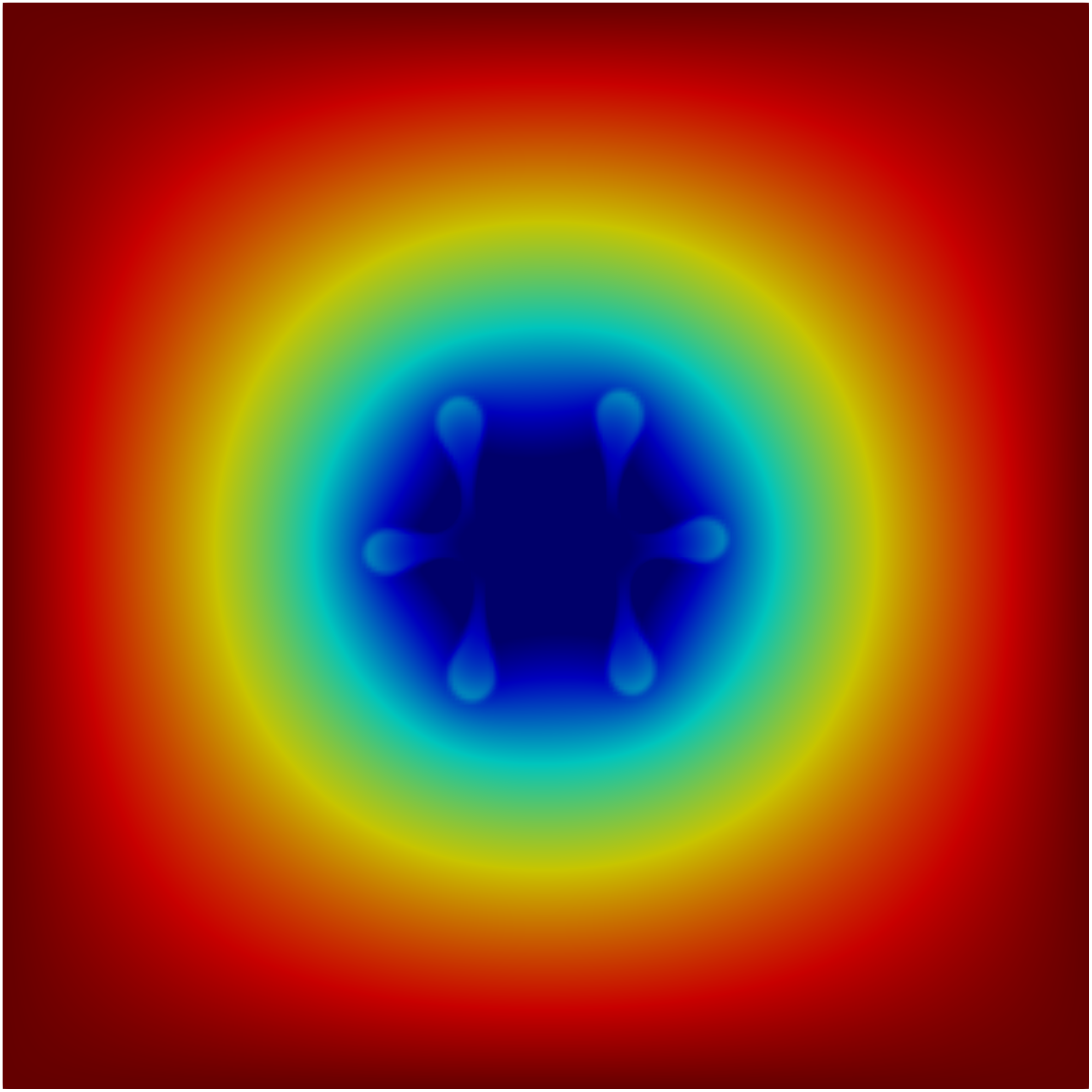} \caption{$t = 1.3$} \end{subfigure}
\begin{subfigure}{0.185\columnwidth} \centering
\includegraphics[trim={0cm 0cm 0cm 0cm},clip,width=1\columnwidth]{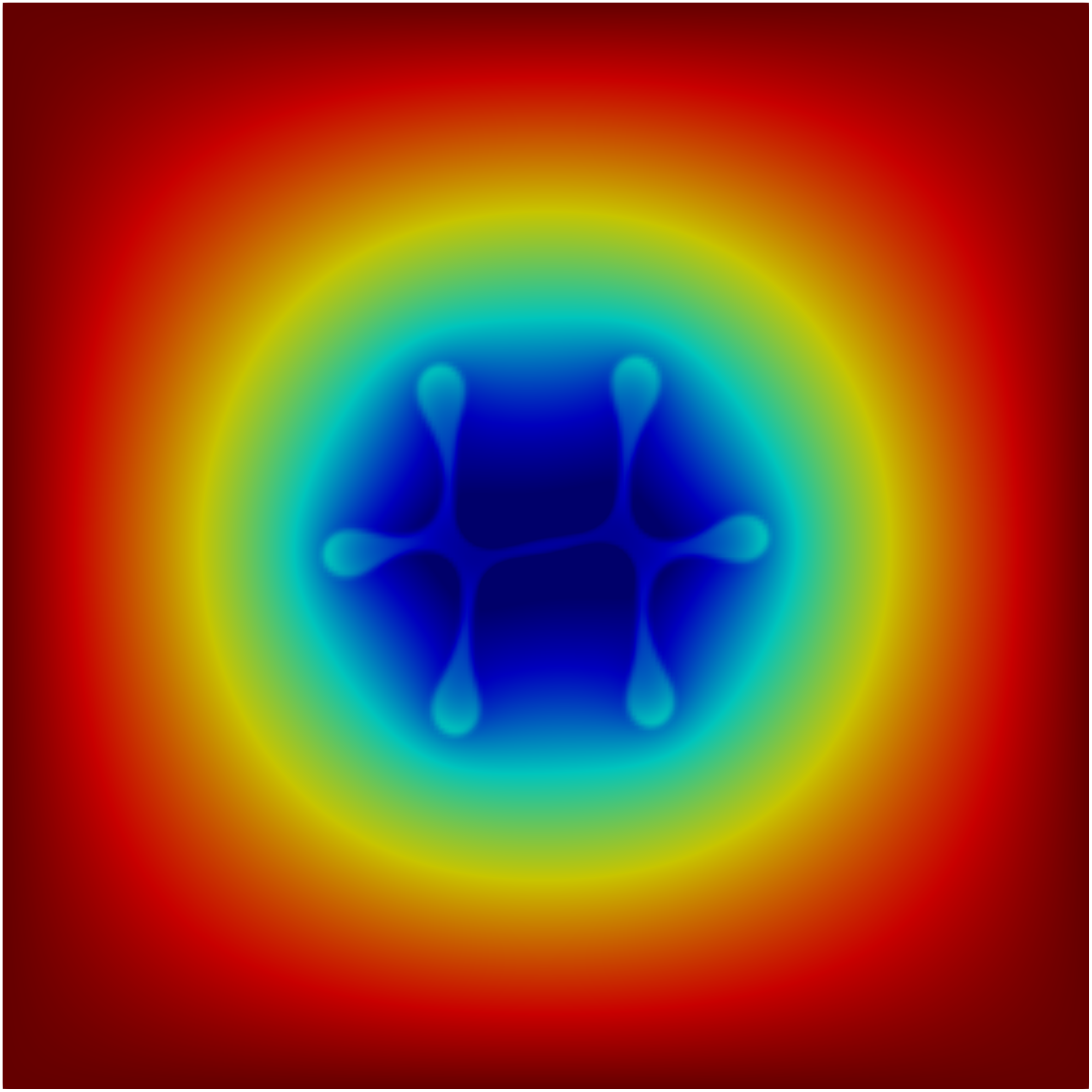} \caption{$t = 1.75$}\end{subfigure}
\begin{subfigure}{0.05\columnwidth} \centering
\includegraphics[trim={0cm -8cm 0cm 0cm},clip,width=1.06\columnwidth]{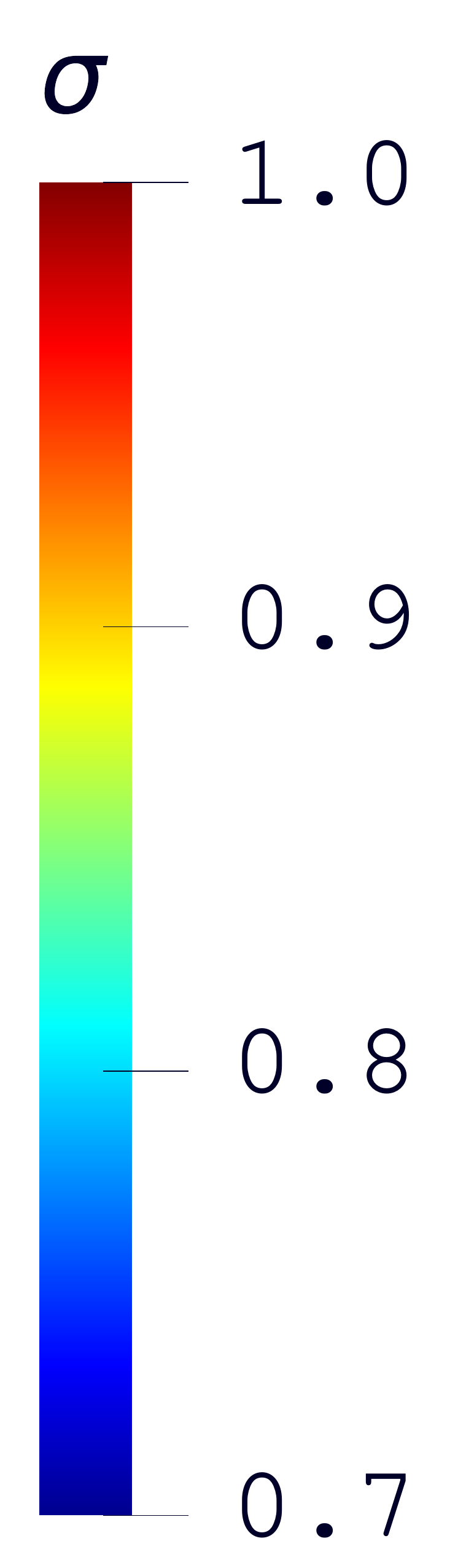} \end{subfigure}
\caption{Evolution of a tumor in a square domain with a perturbed elliptical initial profile defined by Eq.~\eqref{initial_cindition} with $r(\bx) = |\mathbf{x}| - \left( \frac{1}{2} + 10^{-3} \left[ \cos(2\theta) + \frac{5}{4}\cos\!\left(6\theta - \frac{\pi}{12}\right) + \frac{3}{4}\cos\!\left(8\theta - \frac{\pi}{7}\right) \right] \right) $. Top: adaptive THB-spline mesh configurations of degree $3$, $\ell = 6$, $m=2$, $\alpha = 0.01$, $\beta = 0.0001$ with a finest refinement level at $2^{10} \times 2^{10}$ mesh resolution. Middle: evolution of  tumor geometry. Bottom: corresponding nutrient concentration. }
\label{figure_9_final}
\end{figure}

To establish a reference solution, an overkill solution is computed using tensor product quartic B-splines ($p=4$)
with 1024 elements per direction (element size $h_e = 6/1024$).
The tumor morphology and corresponding nutrient concentration at  $t = 0, 1 , 1.5, 2, 2.5$ are shown in Fig.~\ref{figure_1_final}, which are in close agreement with those reported in \cite{Ebenbeck2021}.
Beyond this mesh, additional degree elevation or mesh refinement produce no perceptible variation in the solution. 

Subsequently, using the overkill solution as a reference, we proceed to investigate the minimum B-spline degree and the maximum admissible element size required to accurately resolve the steep gradients at the diffuse interface considering tensor-product uniform B-spline meshes, such that this resolution can be used in the deepest level of the THB spline framework.
To this end, we consider uniform meshes built with refinements using quadratic ($p=2$) and cubic ($p=3$) B-splines with 1024 and 2024 elements per direction.
The relative $L^2(\Omega)$ and $H^1(\Omega)$ errors for the considered mesh configurations are shown in Fig.~\ref{figure_2_final}, which provides a quantitative assessment of the solution accuracy.
Further, the spatial manifestation of errors in the evolving tumor morphology for each discretization against the reference solution across successive time instances is shown in Fig.~\ref{figure_3_final}.
The contours are traced along $\phi \approx 0$ to mark the boundary between the two phases, essentially outlining the tumor shape at each time instance.
Overlaying these contours with the reference solution  highlights the impact of mesh resolution and B-spline degree on the resulting tumor morphology.
The results in Figs.~\ref{figure_2_final} and \ref{figure_3_final} collectively 
establish that the mesh with $p=3$ B-splines with 1024 elements ($h_e = 6/1024$) per spatial direction is adequate for obtaining solutions in close agreement with the reference.
Thus, we choose this mesh configuration as the deepest refinement level along the tumor interface within the THB-spline framework.

After determining the required resolution at the deepest refinement level, we carry out a series of numerical experiments to determine the optimal THB-spline mesh configuration that balances accuracy and associated DOFs.
Eight configurations are considered, varying in hierarchical levels $(\ell = 5 \text{ and } 6)$, refinement and coarsening parameters $\left((\alpha, \beta) = (0.1, 0.001) \text{ and } (0.01, 0.0001)\right)$, and admissibility class $(m = 2 \text{ and } 3)$. 
The corresponding meshes at two representative time instances are shown in Fig.~\ref{figure_4_final}. 
Following the accuracy study on uniform meshes, the finest level for all THB-spline meshes is fixed at mesh resolution $2^{10} \times 2^{10}$ ($h_e=6/1024$) with $p=3$, which ensure that the narrow interfacial region is supported by sufficiently fine elements to accurately resolve the steep phase-field gradients. 
Furthermore, by varying the error indicators $\alpha$ and $\beta$, the extent of the finest refinement band along the interface is controlled. Lower values lead to a broader region of fine elements, while higher values result in a narrower refined zone. 
The admissibility class, governed by $m$, constrains the number of hierarchical levels over which basis functions remain active on any given element. 
We consider $m=2$ and $m=3$, which allow basis function support to span at most two and three hierarchical levels, respectively.
With the resolution at the finest level held fixed, varying the total number of levels provides control over the mesh resolution at the coarsest level.
Collectively, these variations cover a broad range of configurations, enabling a comprehensive assessment to identify the suitable setup for the analysis.

To quantify the error, Fig.~\ref{figure_5_final} shows the relative $L^2(\Omega)$ and $H^1(\Omega)$ error norms computed with respect to the reference solution shown in Fig.~\ref{figure_1_final}.
Additionally, the spatial manifestation of these errors in the evolving tumor morphology is shown in Fig.~\ref{figure_6_final} for each discretization and across successive time instances.
Among all configurations, the adaptive THB-spline meshes of degree 3 with $\ell=5 \text{ and } 6$, $m=2$, $\alpha = 0.01$, $\beta = 0.0001$ yield nearly identical $L^2(\Omega)$ and $H^1(\Omega)$ errors (see Fig.~\ref{figure_5_l2err} and \ref{figure_5_h1err}).
However, the $\ell = 6$ configuration requires fewer DOFs than the $\ell = 5$ configuration, as shown in Fig.~\ref{figure_5_dof}. 
Furthermore, Fig.~\ref{figure_5_time} compares the relative computational time. Each curve exhibits two recurring peaks per adaptation cycle. The first peak corresponds to the additional computational cost associated with mesh refinement, whereas the second reflects the cost of refinement followed by coarsening. The baseline between successive peaks represents the computational effort during time intervals in which the mesh topology remains unchanged.
The magnitude of the peaks varies depending on the number of elements marked for refinement and coarsening and refinement level.
Although the $\ell = 6$ configuration incurs a slightly higher computational cost at a few adaptation steps, %
the overall fitted trend shows that it requires nearly the same, or marginally lower, computational time than the $\ell = 5$ configuration over most of the simulation period.

Therefore, we select the mesh with $\ell = 6$ as the preferred hierarchical mesh. 
This configuration is subsequently used for all remaining numerical examples unless otherwise stated.
Finally, the detailed temporal evolution of an elliptical tumor within the square domain, simulated using the selected mesh configuration ($p=3$, $\ell = 6$, $m=2$, $\alpha = 0.01$, $\beta = 0.0001$) is presented in Fig.~\ref{figure_7_final}.
The figure shows the adaptive THB-spline mesh, the tumor morphology, and the corresponding nutrient concentration field at times $t = 0$, $1$, $1.5$, $2$, and $2.5$.
The tumor progressively develops finger-like protrusions, a morphological feature well documented in the literature and consistent with the results reported in \cite{Ebenbeck2021}, while maintaining close agreement with the overkill reference solution presented herein. 
These findings validate the proposed approach and establish the effectiveness of the locally adaptive THB-spline IGA framework for solving the CH phase-field tumor growth model in its primal form. 

To complete our initial study of the model and assess generalizability of our computational approach, we examine the influence of model parameters and initial conditions on tumor growth dynamics.
Beginning with the effect of model parameters, the proposed model is capable of reproducing not only the finger-like morphology observed in the preceding results but also spheroidal tumor growth, another well-documented tumor morphology in the literature.
This transition is achieved through an increase in the proliferation rate $\mathcal{P}$ from 0.1 to 5, while all remaining parameters are kept at their baseline values, thereby increasing the role of proliferation in driving the observed shift in tumor morphology. 
Hence, although we still use the same elliptical configuration for initial conditions, the tumor progressively evolves into a smooth spheroidal morphology rather than an irregular finger-like shape as shown in Fig.~\ref{figure_8_final}.
The morphological diversity is not limited to the two shapes discussed so far, as elongated tumor shapes can also be obtained through appropriate parameter adjustments. These results are presented in Fig.~\ref{figure_2_appendix} in \ref{Complementary_Numerical_Studies}. %
Collectively, these results highlight the sensitivity of parameter selection on the resulting tumor shapes.
Finally, the sensitivity of the tumor morphology to the initial condition is examined by introducing targeted perturbations to the initial phase-field profile $r(\bx)$ (see Eq.~\eqref{initial_cindition}) while keeping all model parameters unchanged. The results, presented in Fig.~\ref{figure_9_final} and complemented by an additional example in Fig.~\ref{figure_3_appendix} in \ref{Complementary_Numerical_Studies}, demonstrate that minimal changes to the initial configuration can give rise to entirely different growth morphologies, highlighting the dependence of the tumor evolution on its initial configuration.

\begin{figure}[!t]\centering
\begin{subfigure}{0.245\columnwidth} \centering
\includegraphics[trim={0cm 0cm 0cm 0cm},clip,width=0.82\columnwidth]{Images/Problem_1/leg6.png}
\end{subfigure}
\begin{subfigure}{0.245\columnwidth} \centering
\includegraphics[trim={0cm 14.5cm 0cm 15cm},clip,width=0.9\columnwidth]{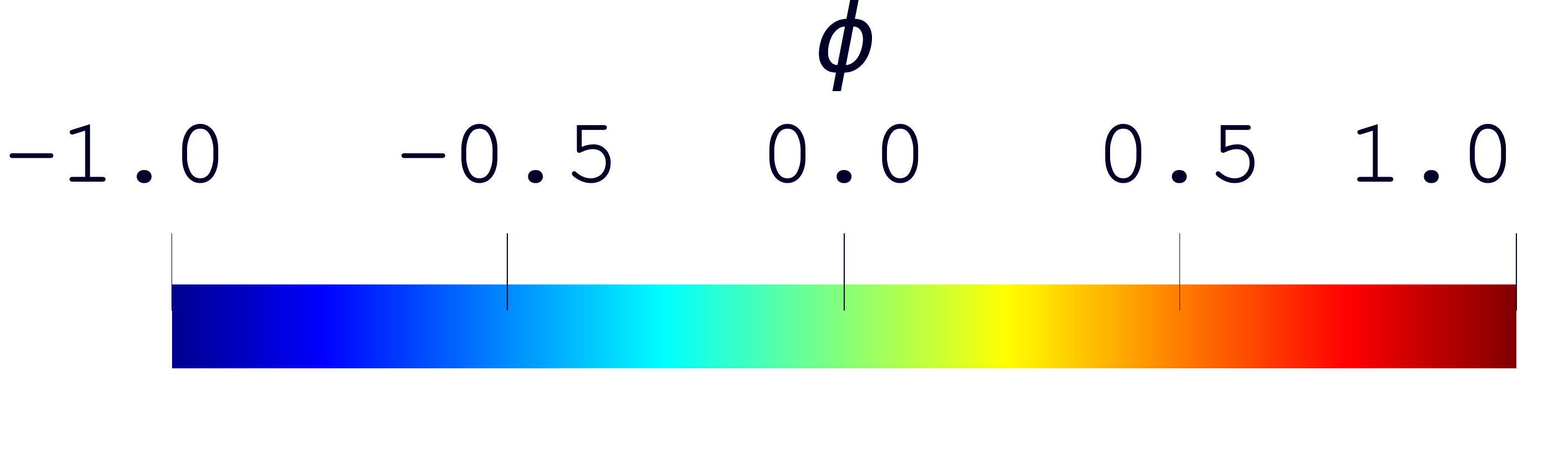}
\end{subfigure}
\begin{subfigure}{0.245\columnwidth} \centering
\includegraphics[trim={0cm 14.5cm 0cm 15cm},clip,width=0.9\columnwidth]{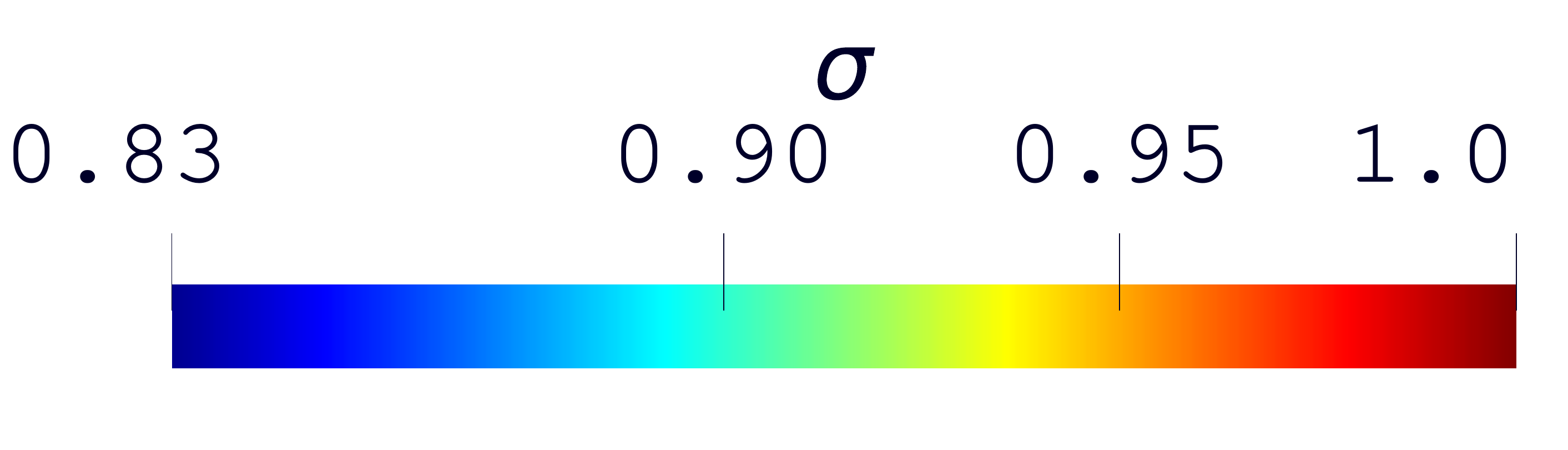}
\end{subfigure}
\begin{subfigure}{0.245\columnwidth} \centering
\mbox{}
\end{subfigure}
\\
\begin{subfigure}{1\columnwidth} \centering
\includegraphics[trim={0cm 0cm 0cm 0cm},clip,width=0.245\columnwidth]{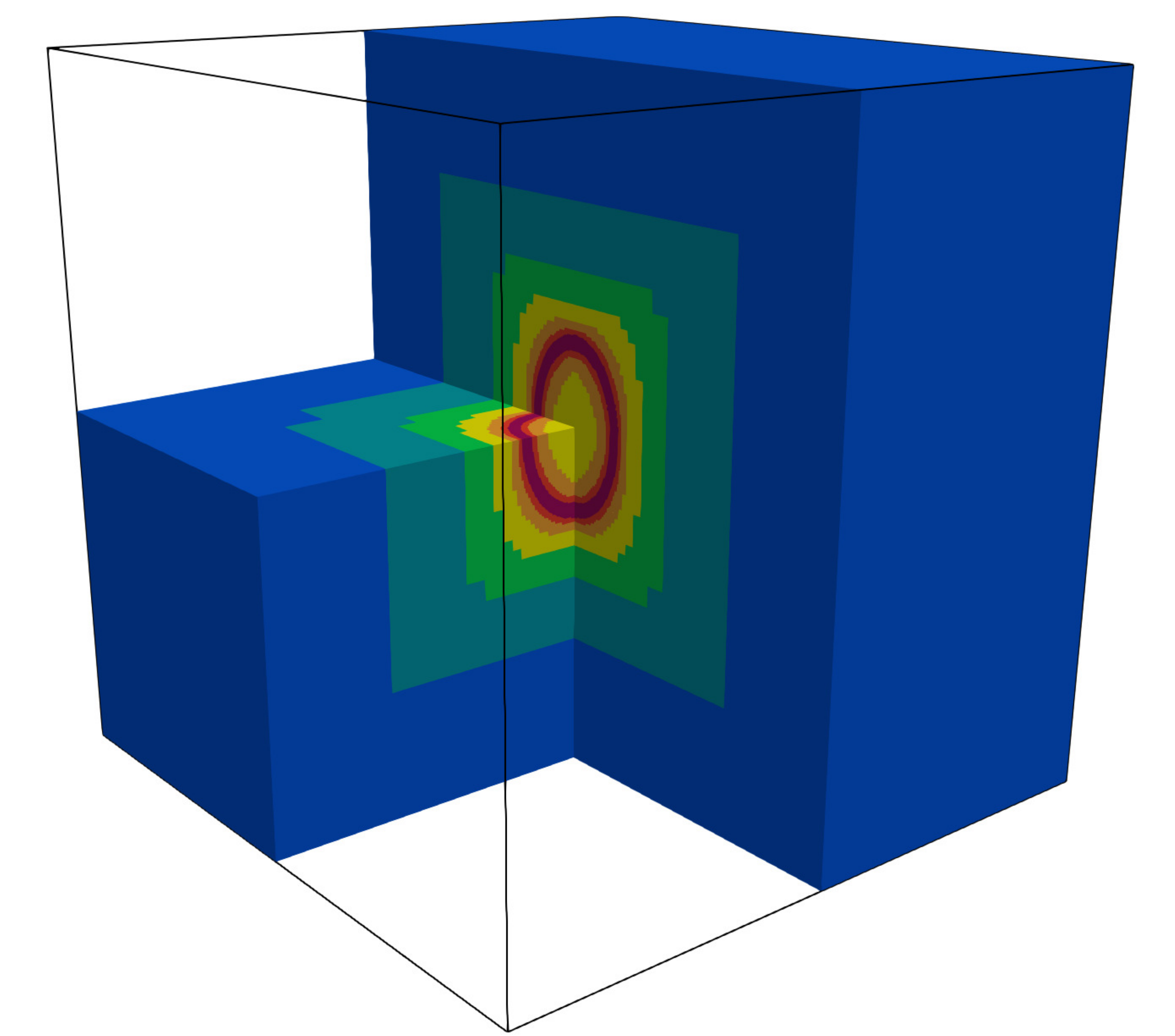}
\includegraphics[trim={0cm 0cm 0cm 0cm},clip,width=0.245\columnwidth]{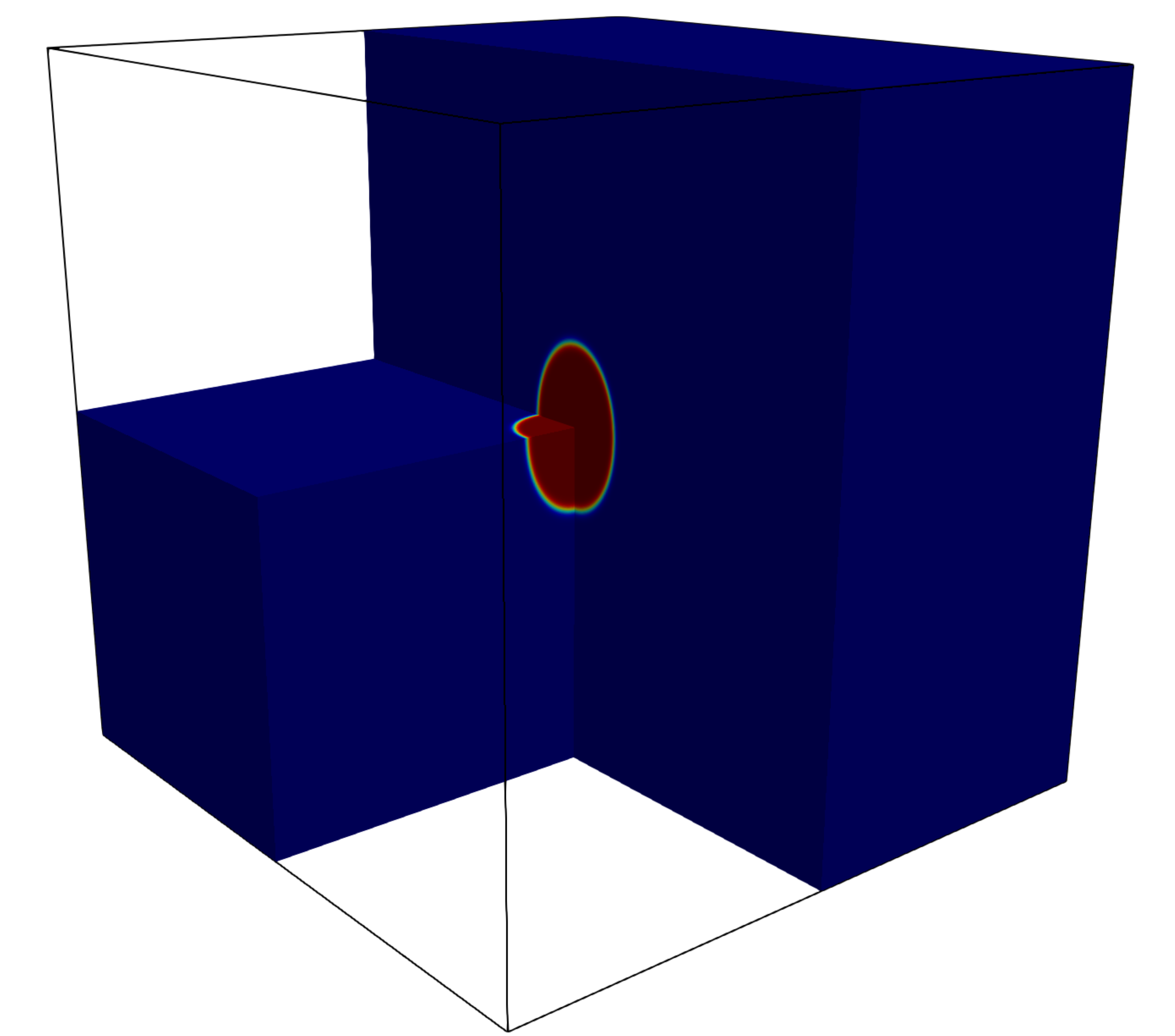}
\includegraphics[trim={0cm 0cm 0cm 0cm},clip,width=0.245\columnwidth]{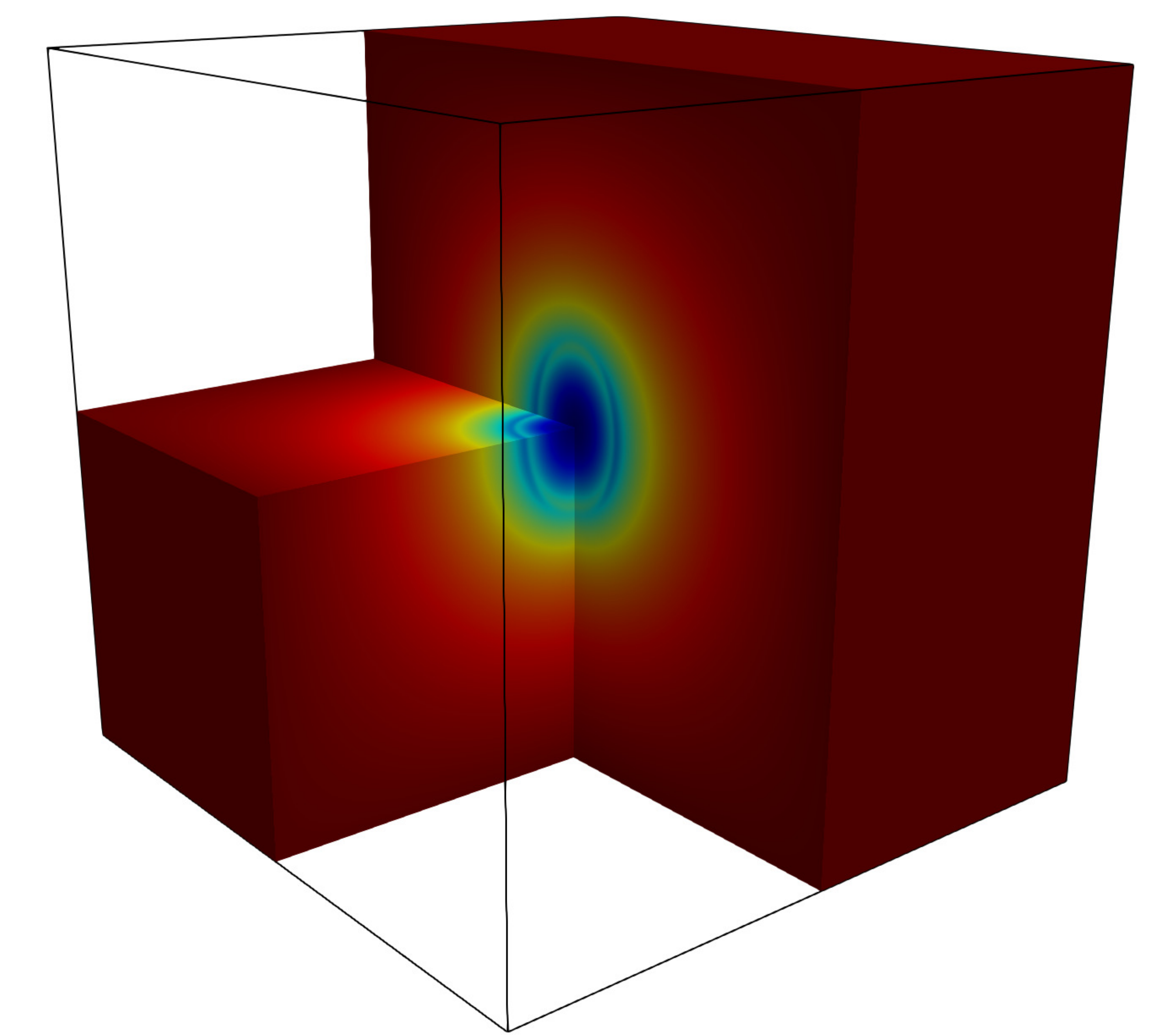}
\includegraphics[trim={0cm 0cm 0cm 0cm},clip,width=0.245\columnwidth]{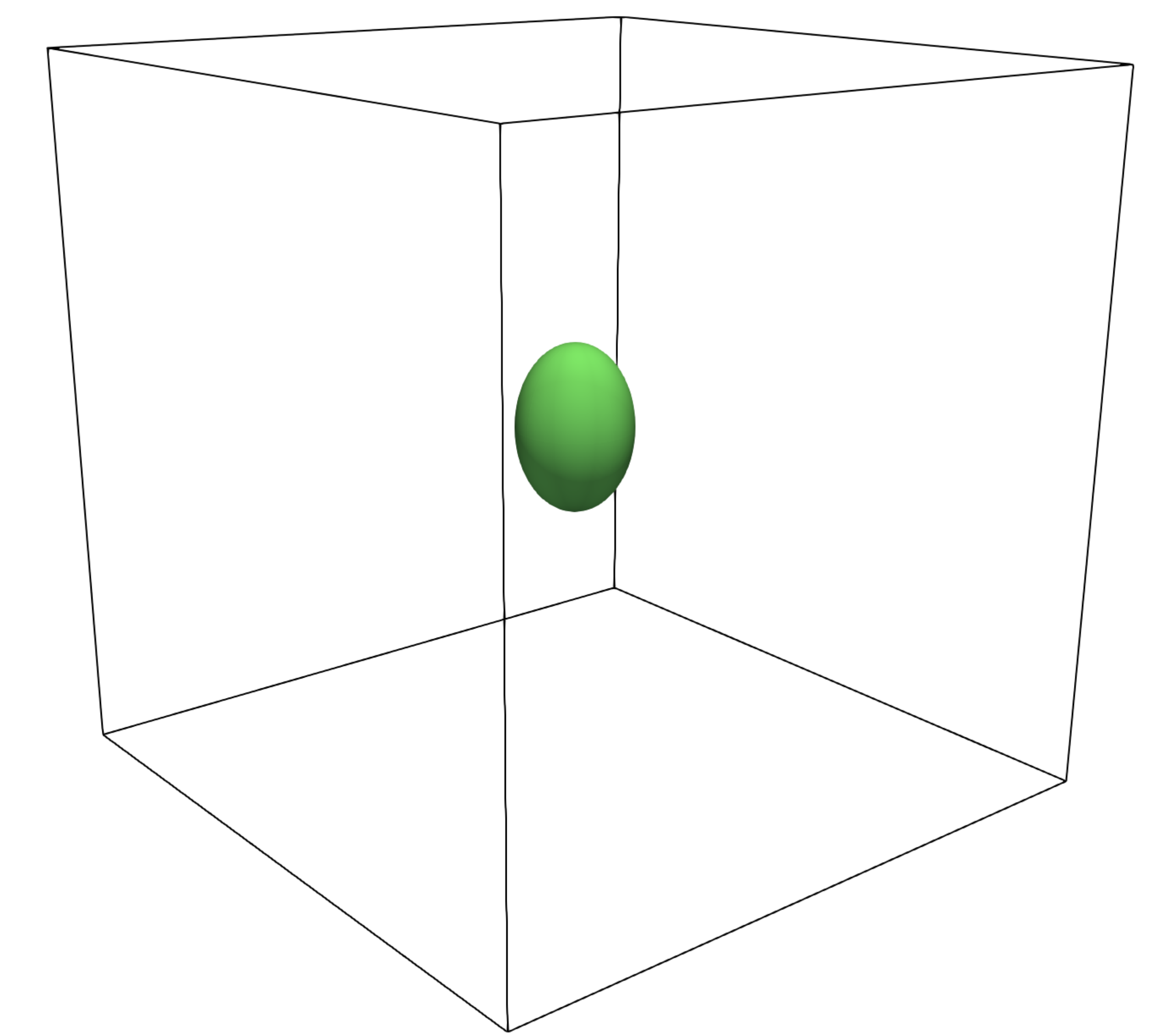}
\caption{$t=0$}
\end{subfigure}
\\
\begin{subfigure}{1\columnwidth} \centering
\includegraphics[trim={0cm 0cm 0cm 0cm},clip,width=0.245\columnwidth]{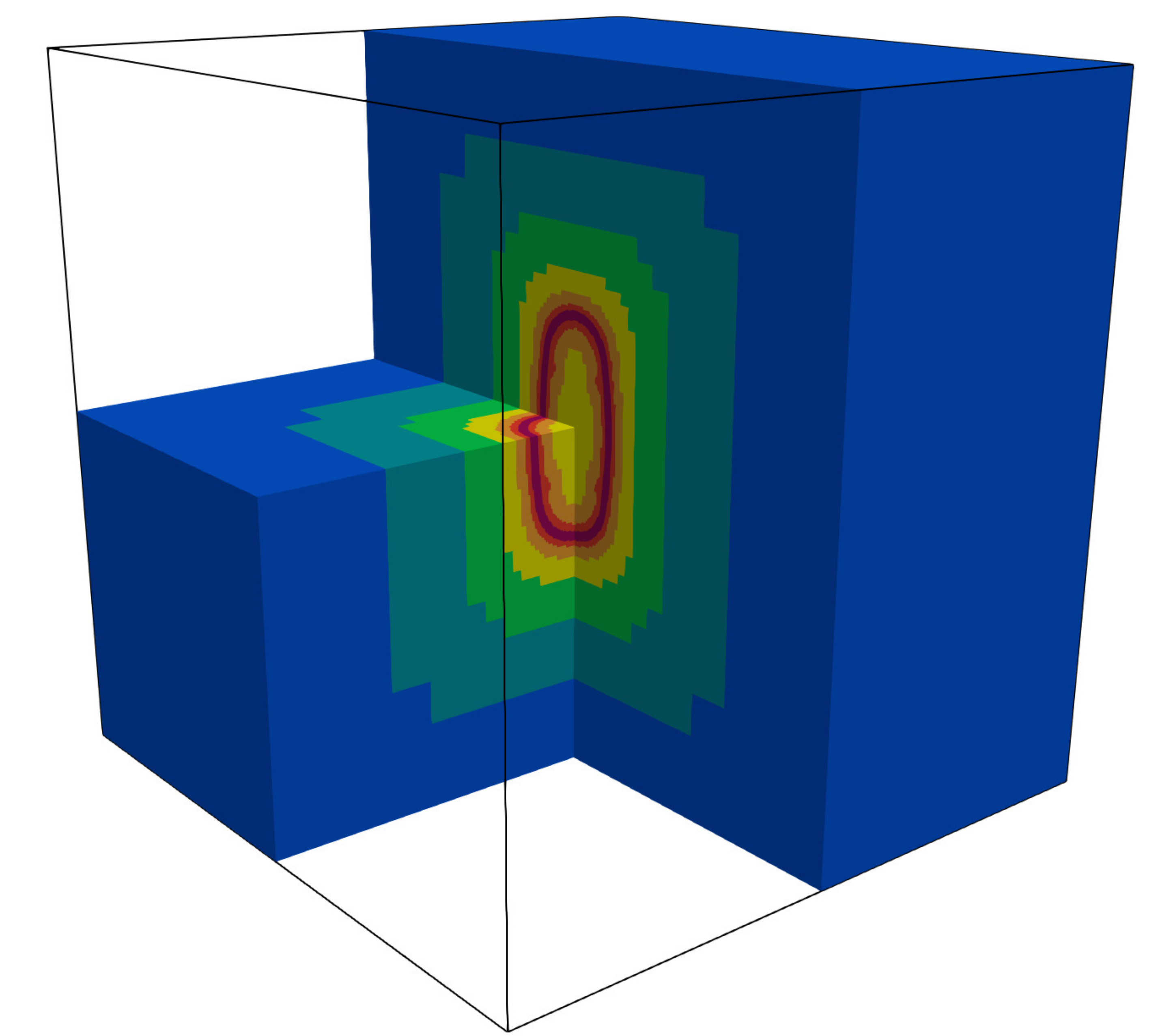}
\includegraphics[trim={0cm 0cm 0cm 0cm},clip,width=0.245\columnwidth]{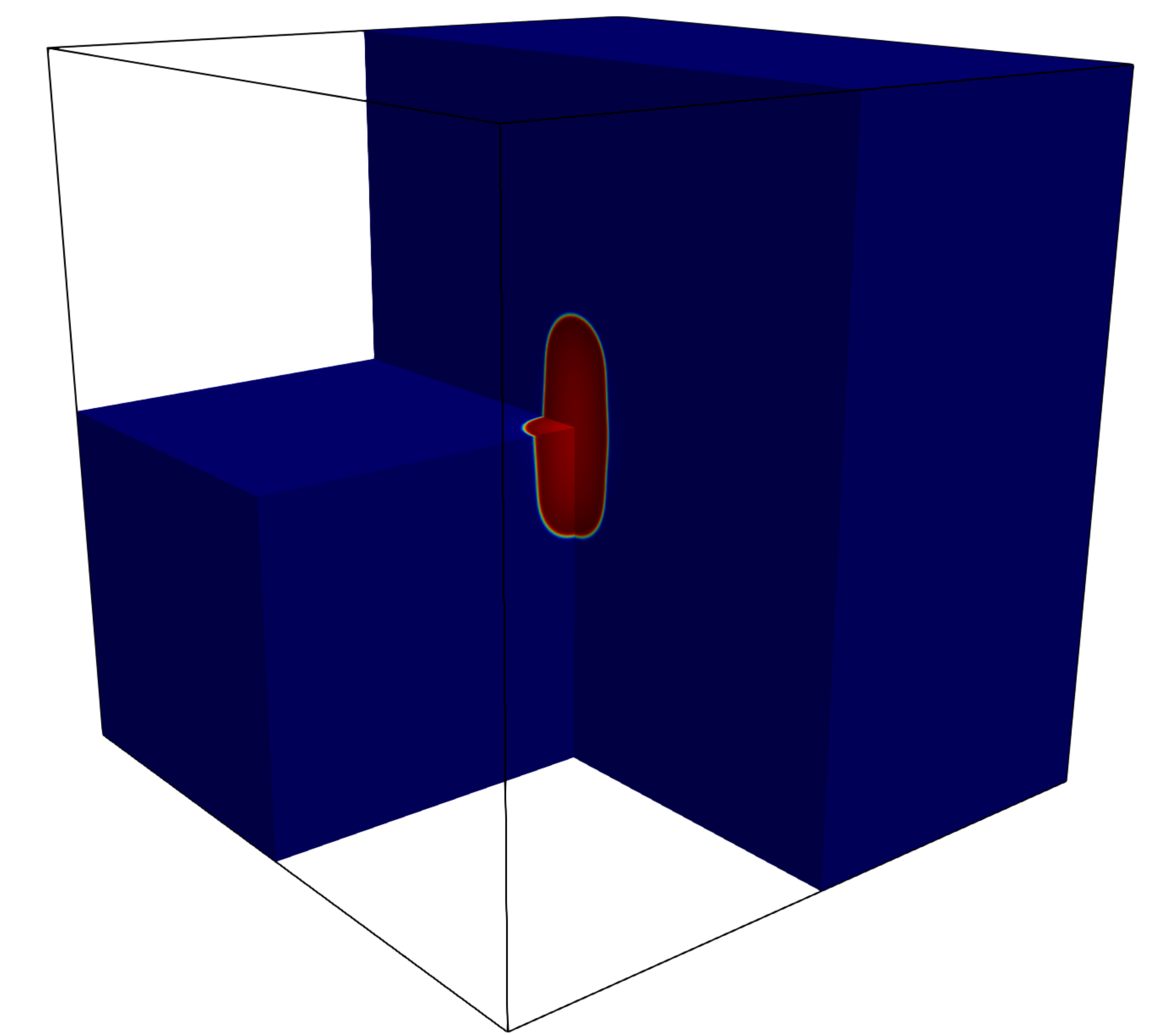}
\includegraphics[trim={0cm 0cm 0cm 0cm},clip,width=0.245\columnwidth]{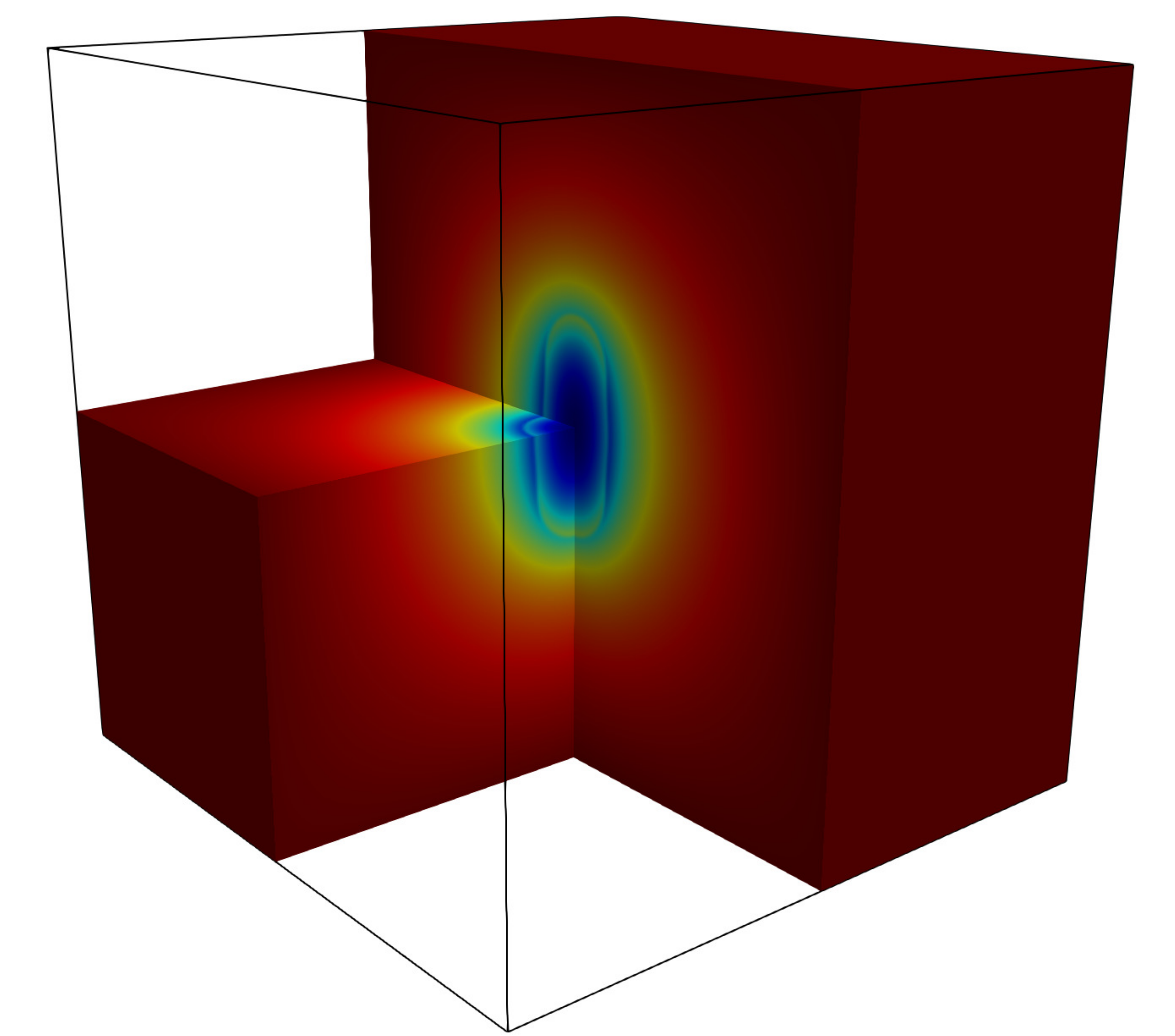}
\includegraphics[trim={0cm 0cm 0cm 0cm},clip,width=0.245\columnwidth]{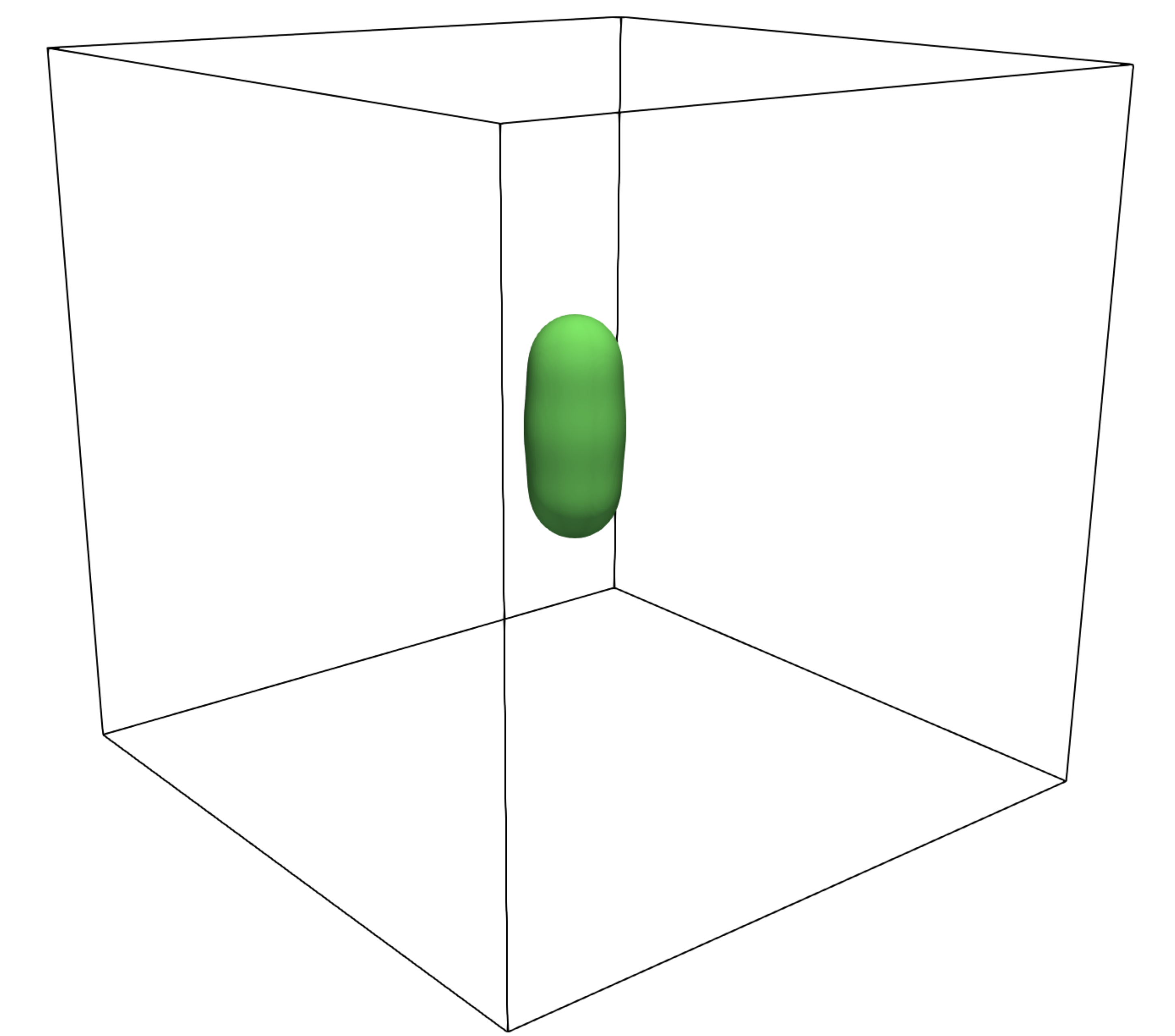}
\caption{$t=1$}
\end{subfigure}
\\
\begin{subfigure}{1\columnwidth} \centering
\includegraphics[trim={0cm 0cm 0cm 0cm},clip,width=0.245\columnwidth]{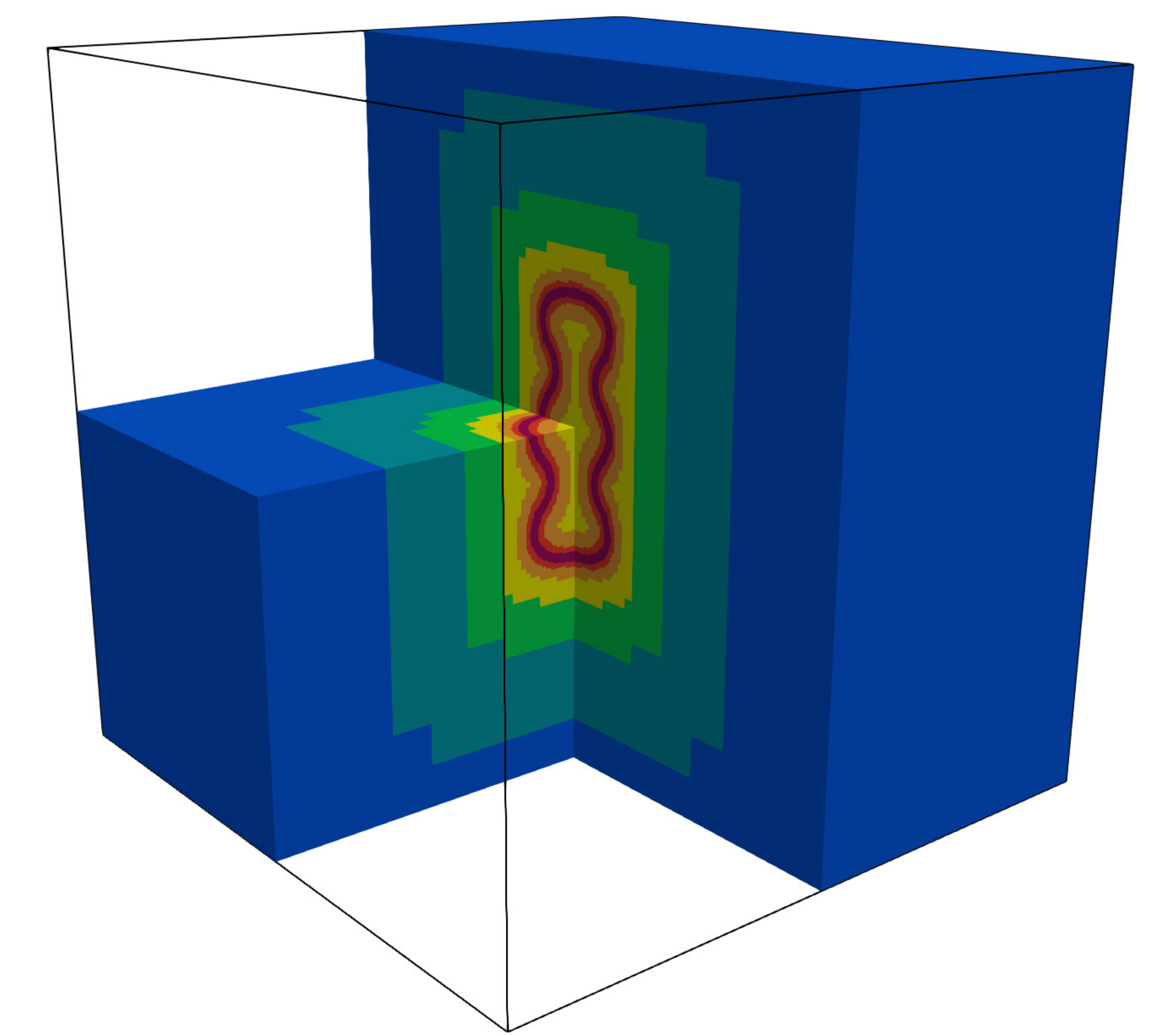}
\includegraphics[trim={0cm 0cm 0cm 0cm},clip,width=0.245\columnwidth]{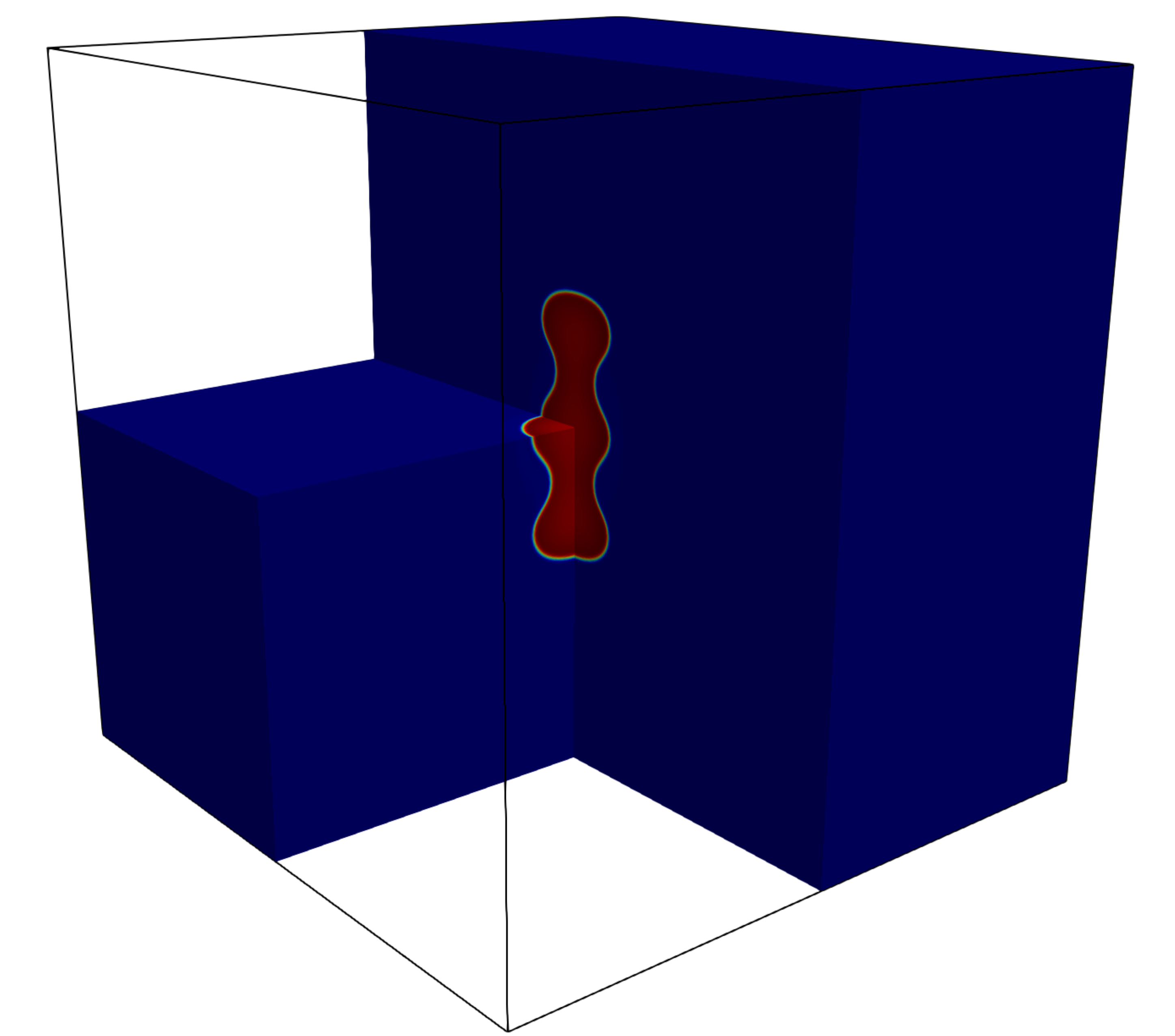}
\includegraphics[trim={0cm 0cm 0cm 0cm},clip,width=0.245\columnwidth]{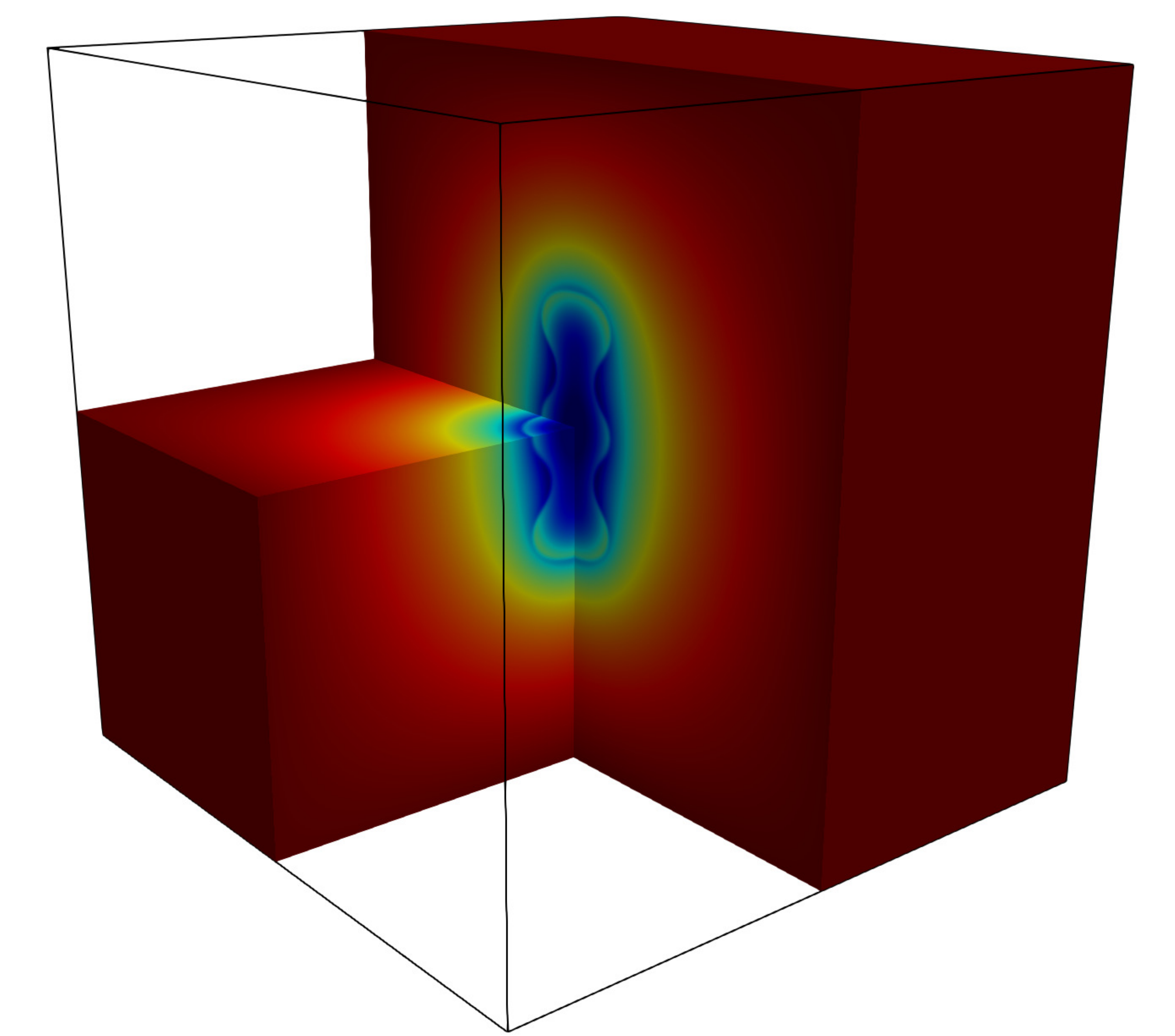}
\includegraphics[trim={0cm 0cm 0cm 0cm},clip,width=0.245\columnwidth]{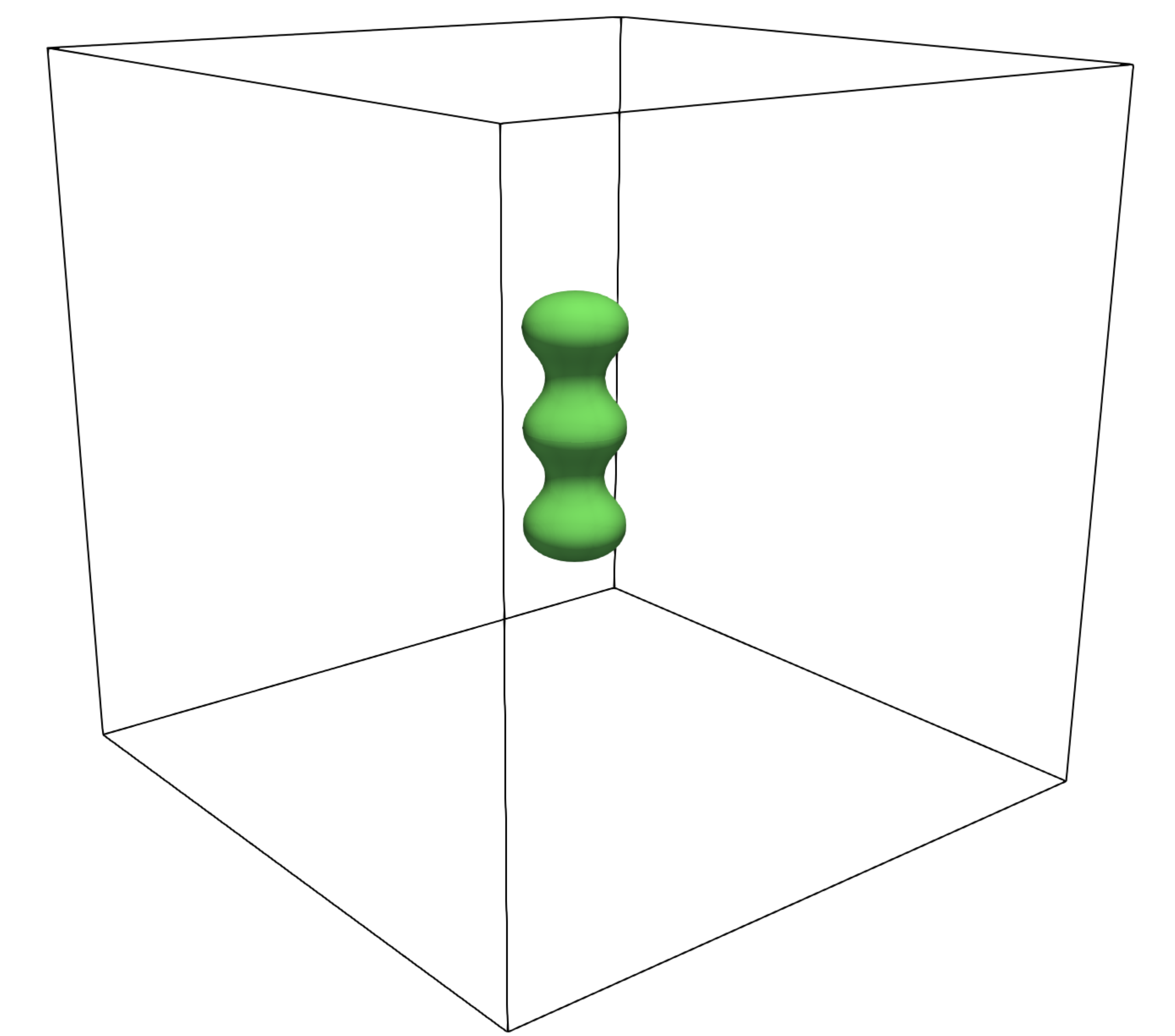}
\caption{$t=2$}
\end{subfigure}
\\
\begin{subfigure}{1\columnwidth} \centering
\includegraphics[trim={0cm 0cm 0cm 0cm},clip,width=0.245\columnwidth]{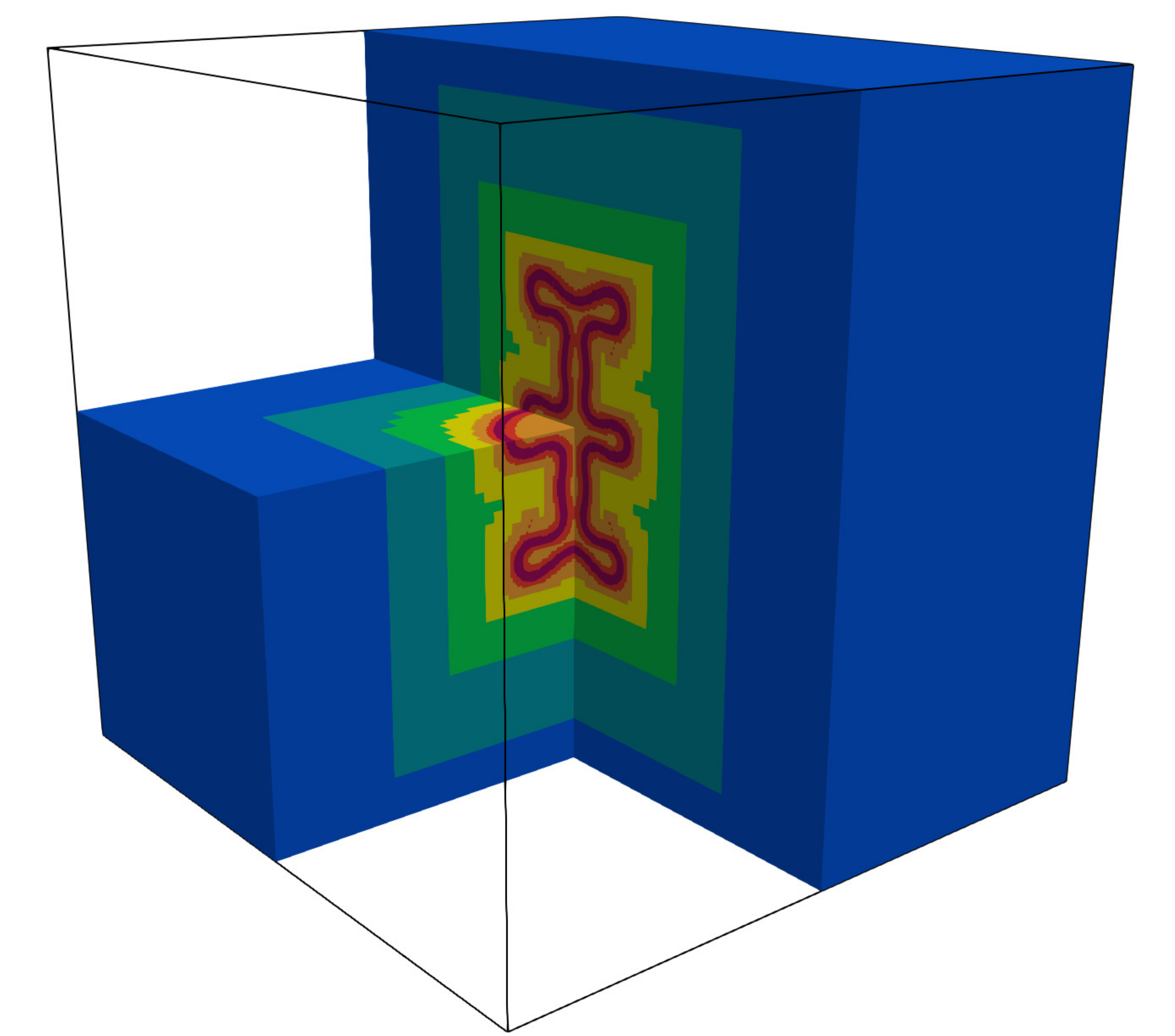}
\includegraphics[trim={0cm 0cm 0cm 0cm},clip,width=0.245\columnwidth]{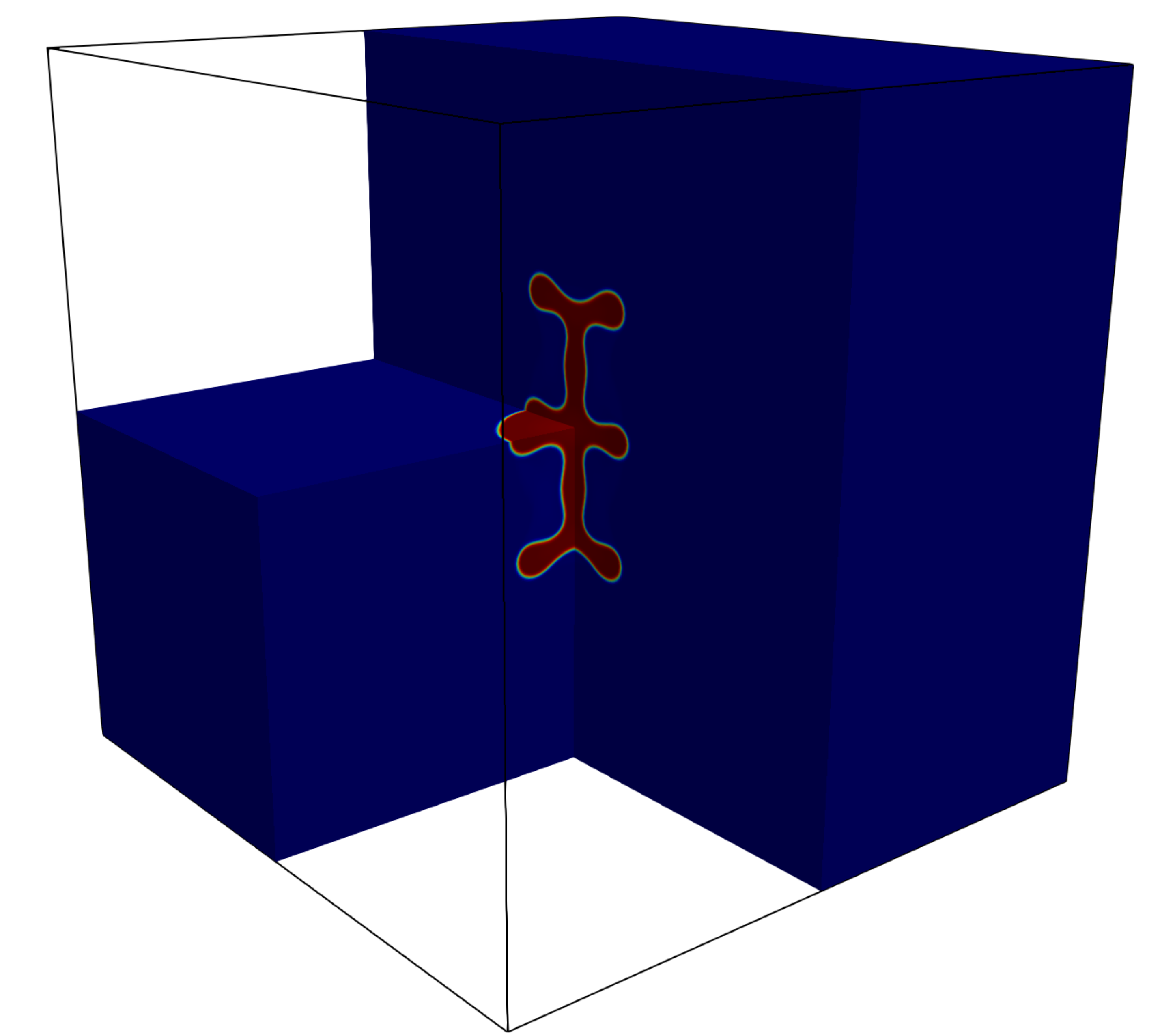}
\includegraphics[trim={0cm 0cm 0cm 0cm},clip,width=0.245\columnwidth]{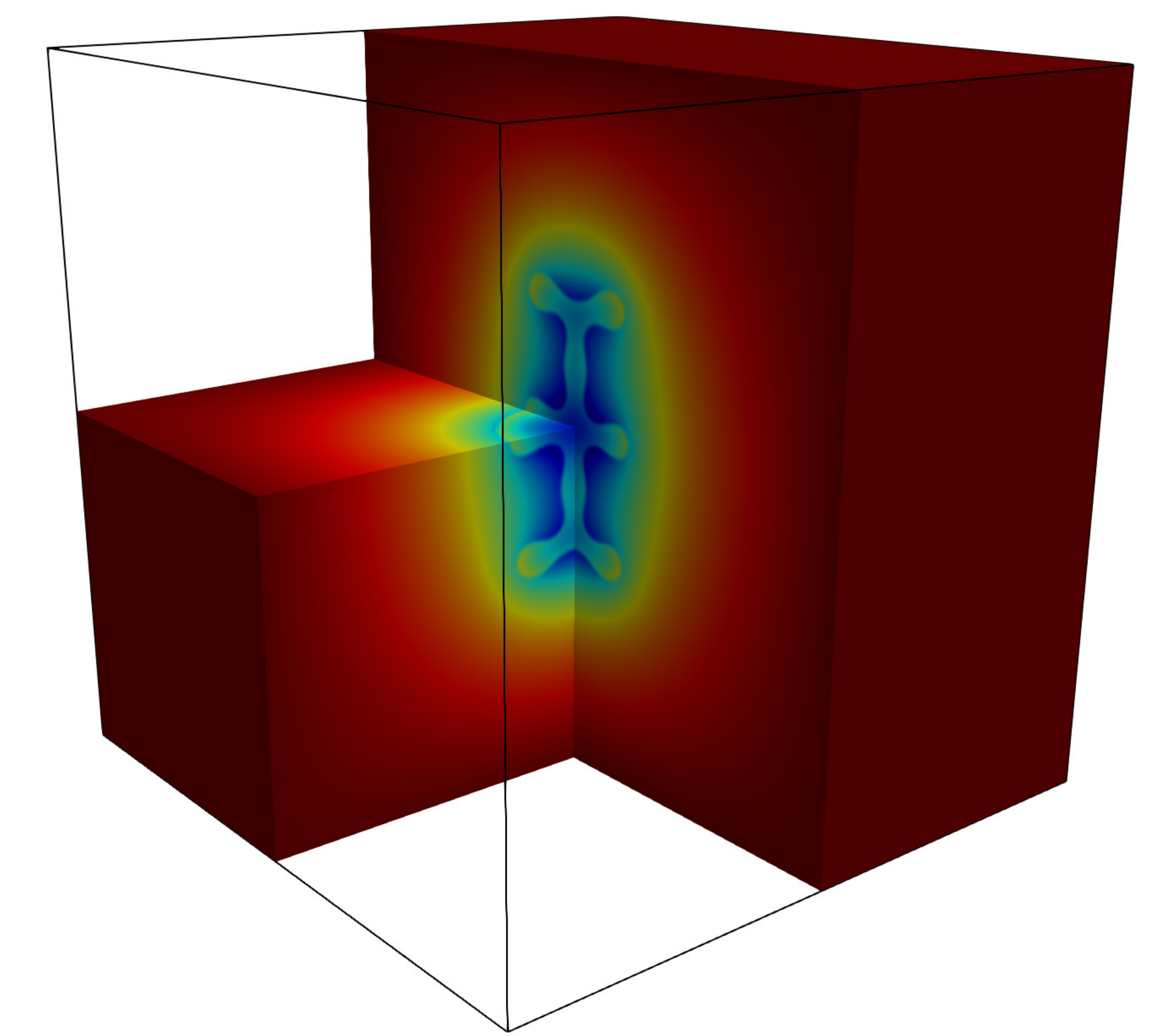}
\includegraphics[trim={0cm 0cm 0cm 0cm},clip,width=0.245\columnwidth]{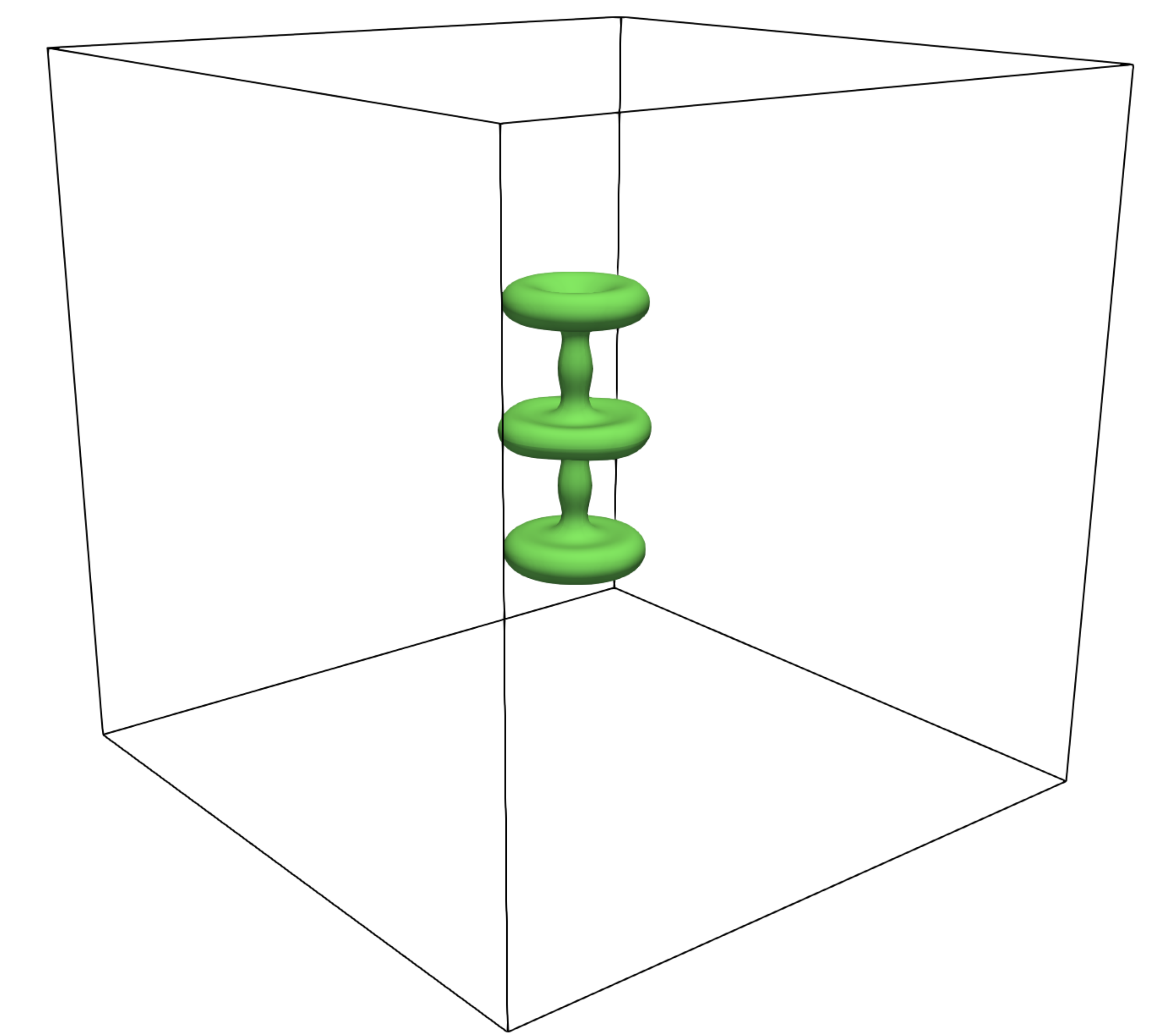}
\caption{$t=3$ }
\end{subfigure}
\caption{Spatiotemporal evolution of an ellipsoidal tumor within a cubic domain. From left to right: adaptive THB-spline mesh ($p =3$, $\ell = 6$, $m=2$, $\alpha = 0.01$, $\beta = 0.0001$), phase-field $\phi$, nutrient concentration field $\sigma$, and the resulting tumor morphology at times $t=$ 0, 1, 2, and 3.}
\label{figure_10_final}
\end{figure}

\begin{figure}[!t]\centering
\begin{subfigure}{0.245\columnwidth} \centering
\includegraphics[trim={0cm 0cm 0cm 0cm},clip,width=0.82\columnwidth]{Images/Problem_1/leg6.png}
\end{subfigure}
\begin{subfigure}{0.245\columnwidth} \centering
\includegraphics[trim={0cm 14.5cm 0cm 15cm},clip,width=0.9\columnwidth]{Images/Problem1_3D/ColorBar_phi_3D-eps-converted-to.pdf}
\end{subfigure}
\begin{subfigure}{0.245\columnwidth} \centering
\includegraphics[trim={0cm 14.5cm 0cm 15cm},clip,width=0.9\columnwidth]{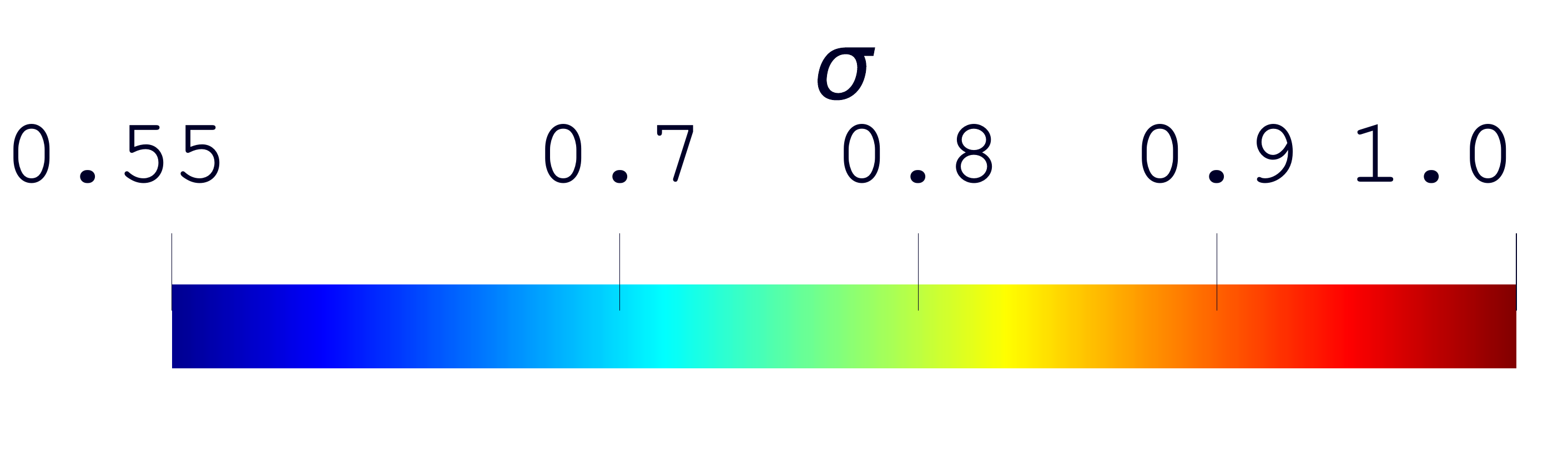}
\end{subfigure}
\begin{subfigure}{0.245\columnwidth} \centering
\mbox{}
\end{subfigure}
\\
\begin{subfigure}{1\columnwidth} \centering
\includegraphics[trim={0cm 0cm 0cm 0cm},clip,width=0.245\columnwidth]{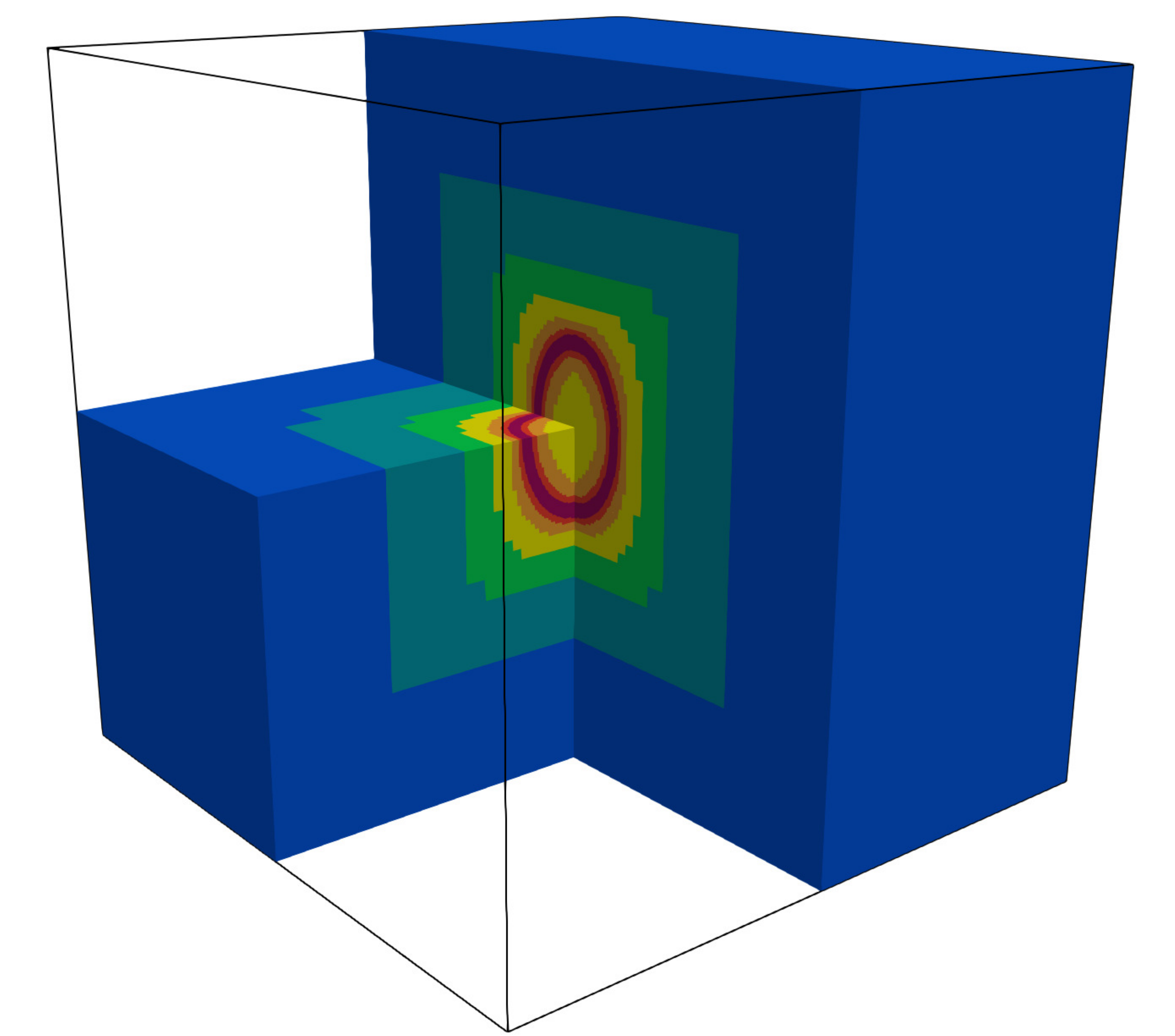}
\includegraphics[trim={0cm 0cm 0cm 0cm},clip,width=0.245\columnwidth]{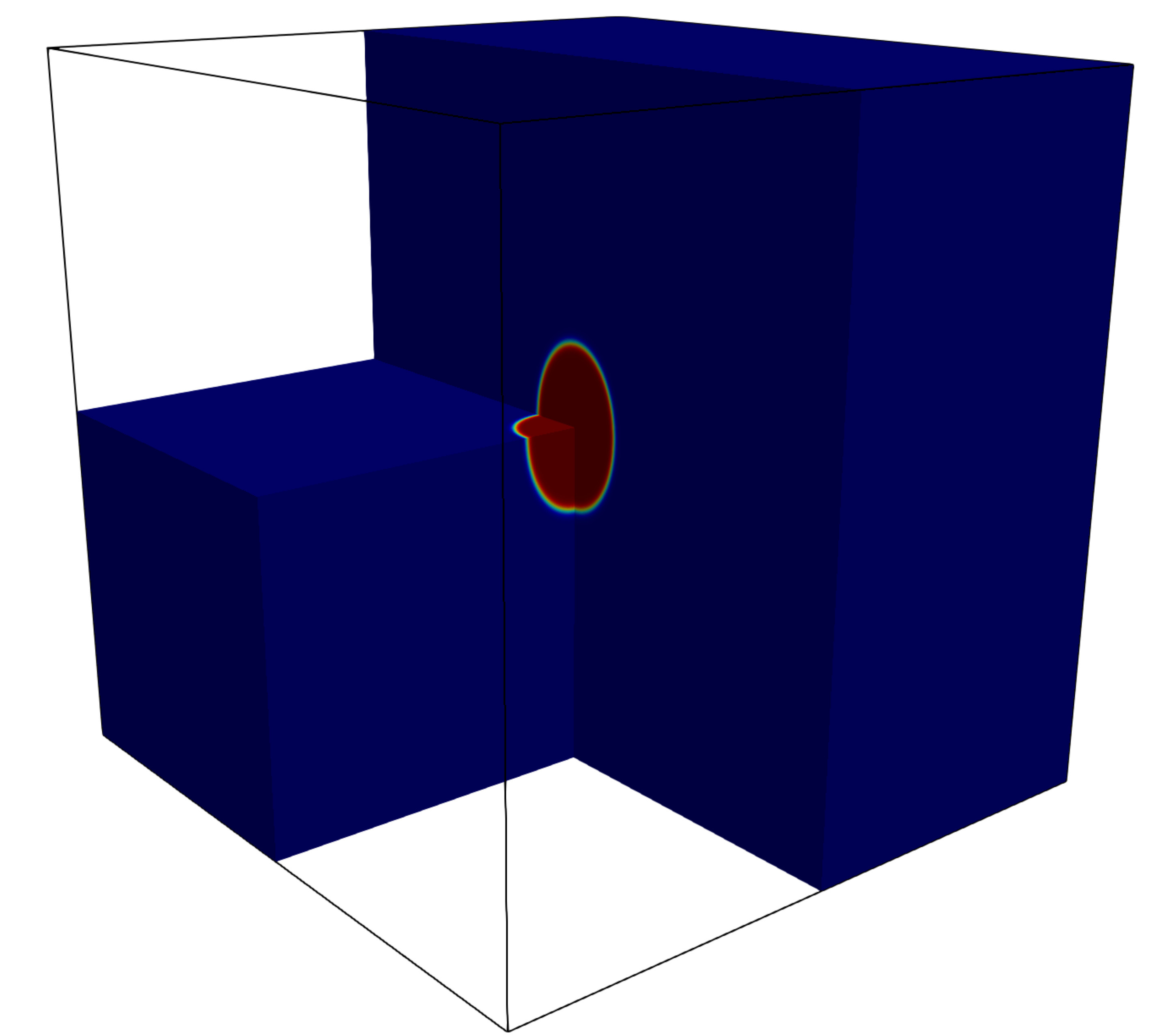}
\includegraphics[trim={0cm 0cm 0cm 0cm},clip,width=0.245\columnwidth]{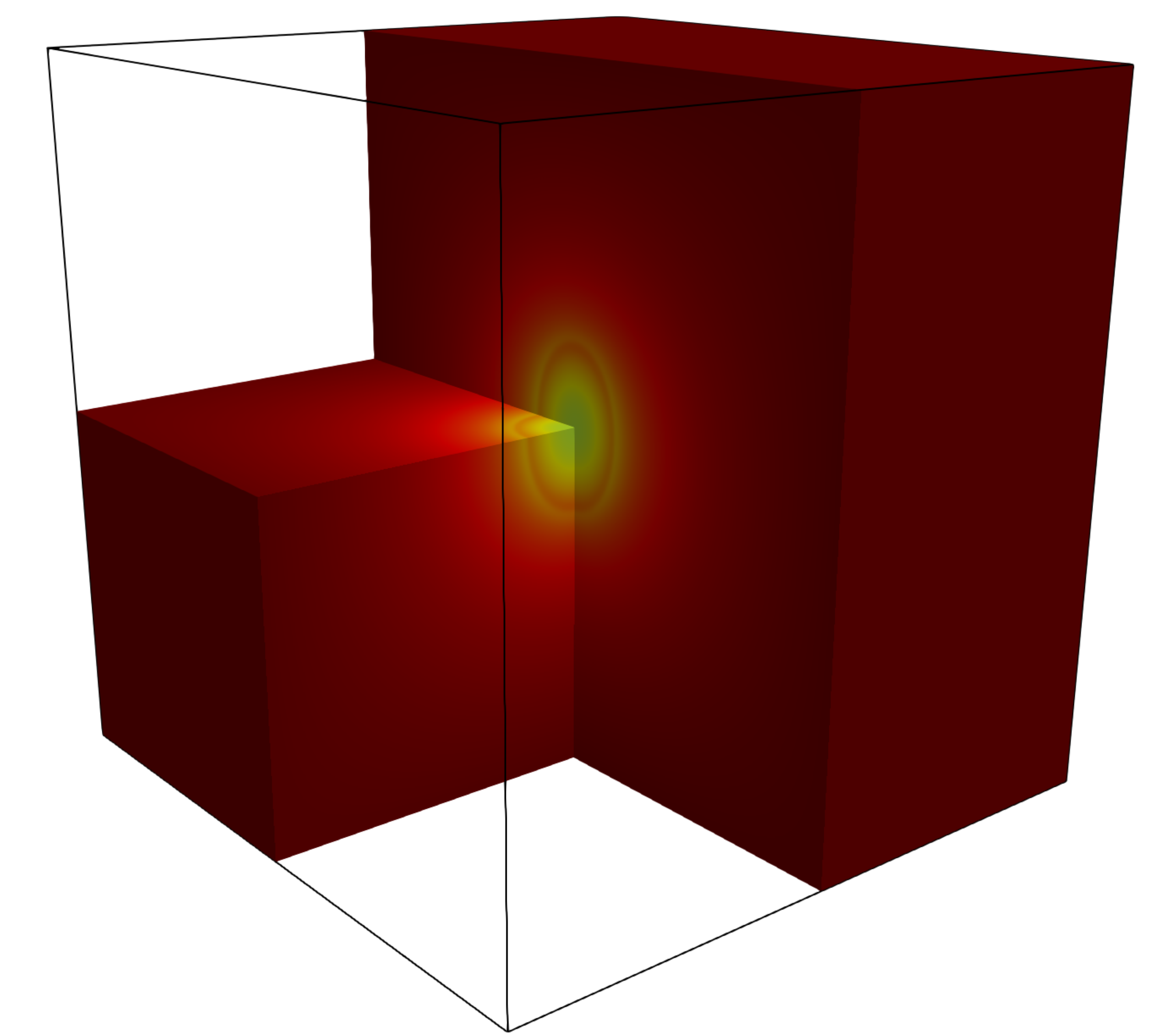}
\includegraphics[trim={0cm 0cm 0cm 0cm},clip,width=0.245\columnwidth]{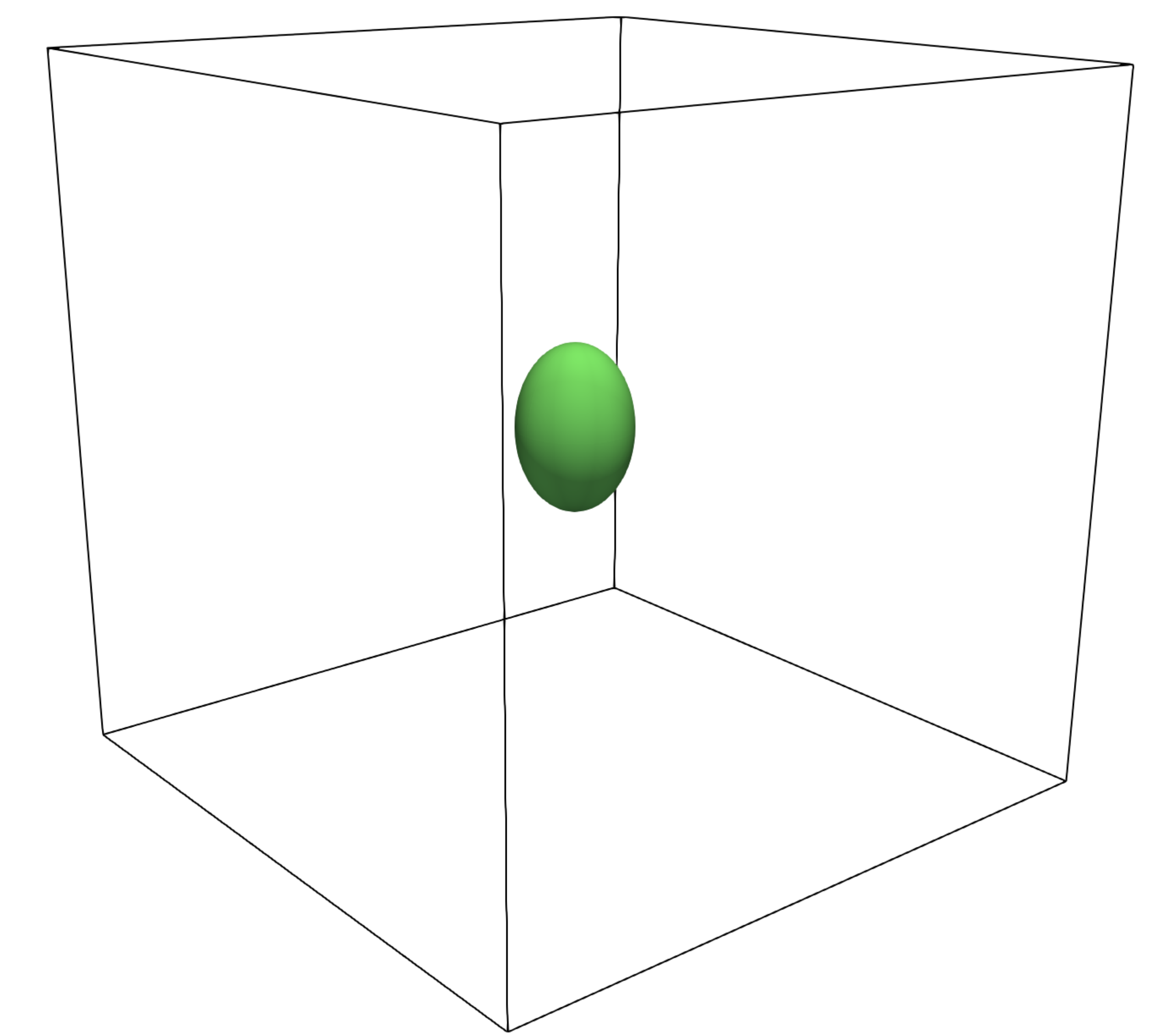}
\caption{$t=0$}
\end{subfigure}
\\
\begin{subfigure}{1\columnwidth} \centering
\includegraphics[trim={0cm 0cm 0cm 0cm},clip,width=0.245\columnwidth]{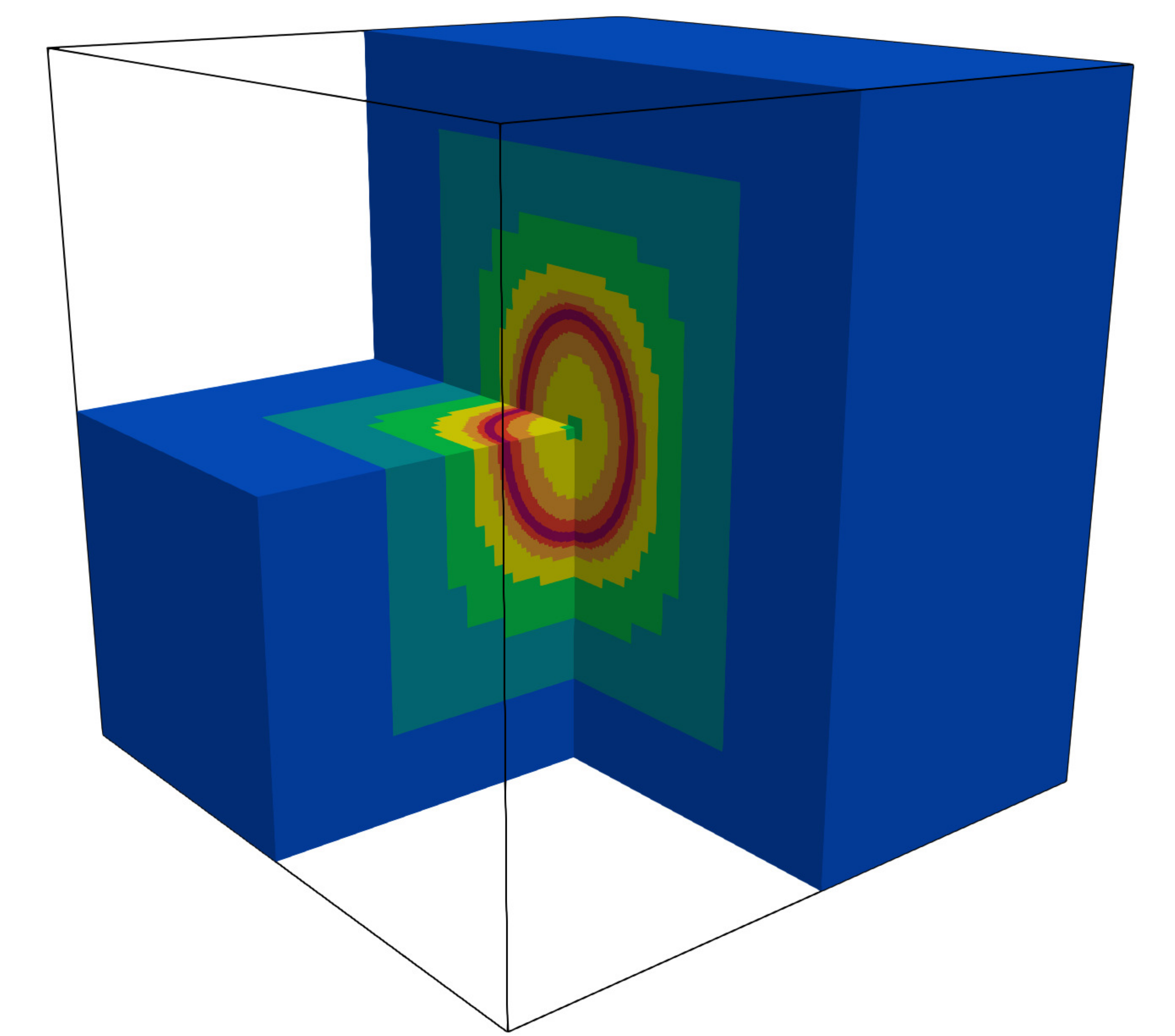}
\includegraphics[trim={0cm 0cm 0cm 0cm},clip,width=0.245\columnwidth]{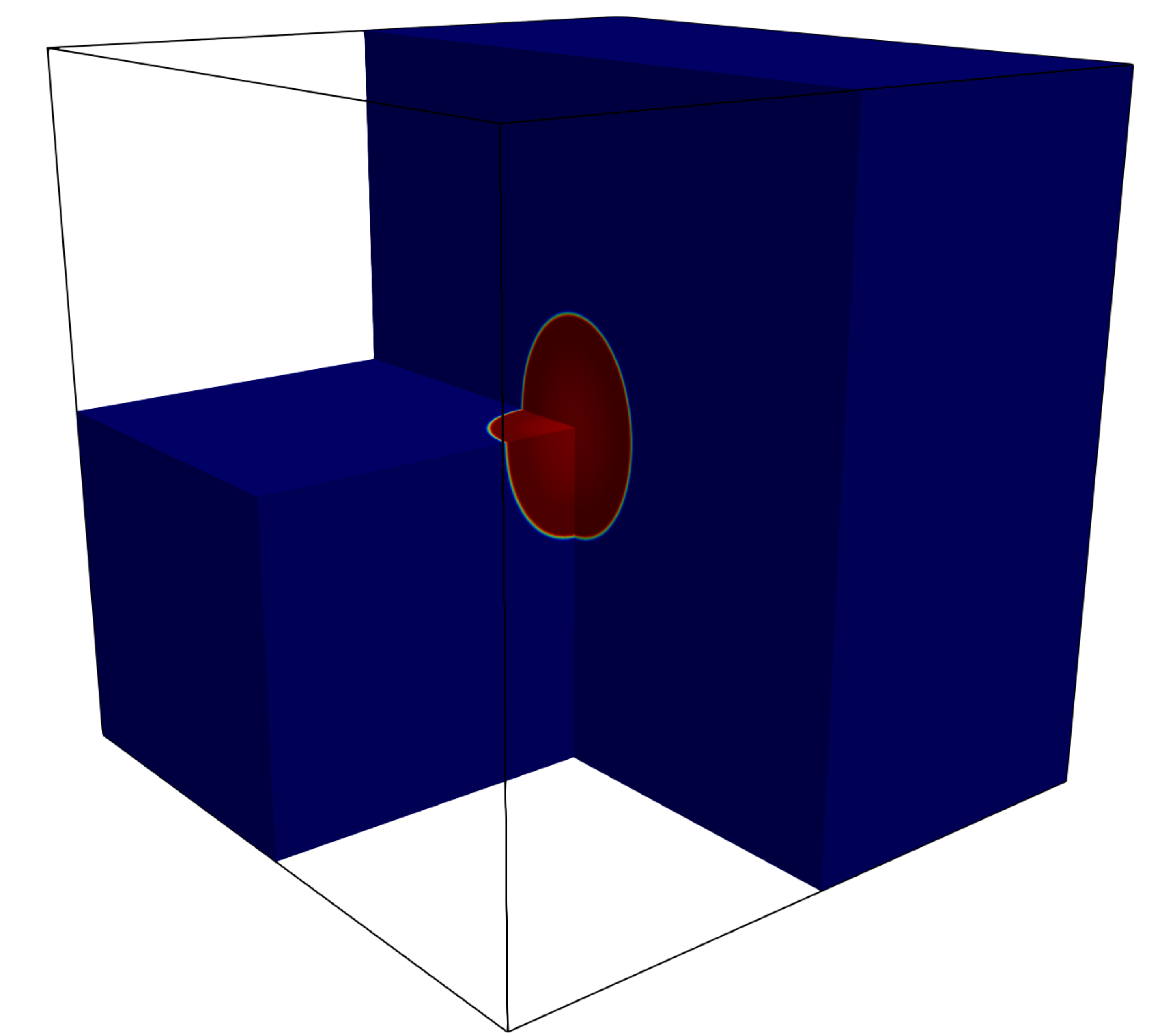}
\includegraphics[trim={0cm 0cm 0cm 0cm},clip,width=0.245\columnwidth]{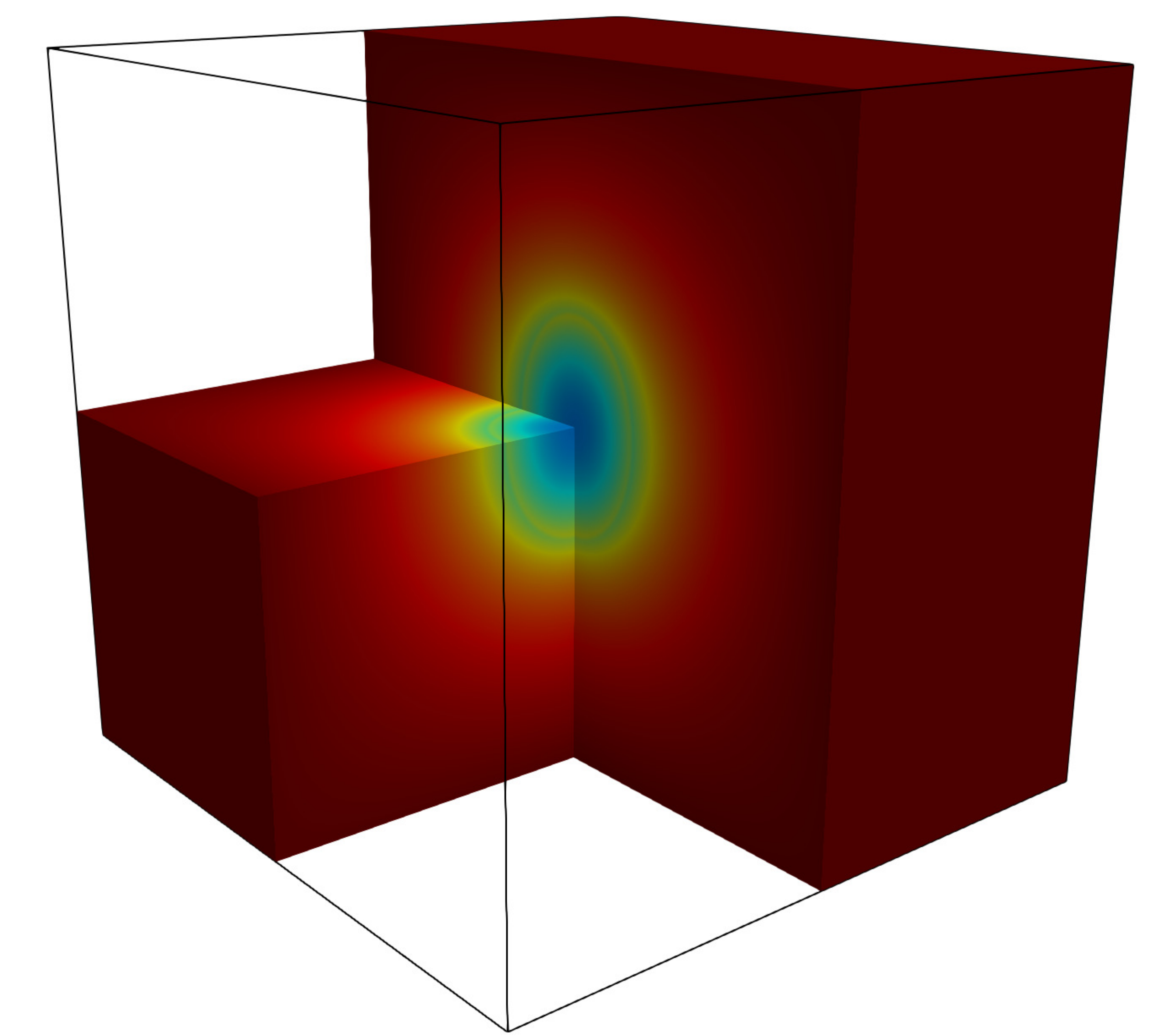}
\includegraphics[trim={0cm 0cm 0cm 0cm},clip,width=0.245\columnwidth]{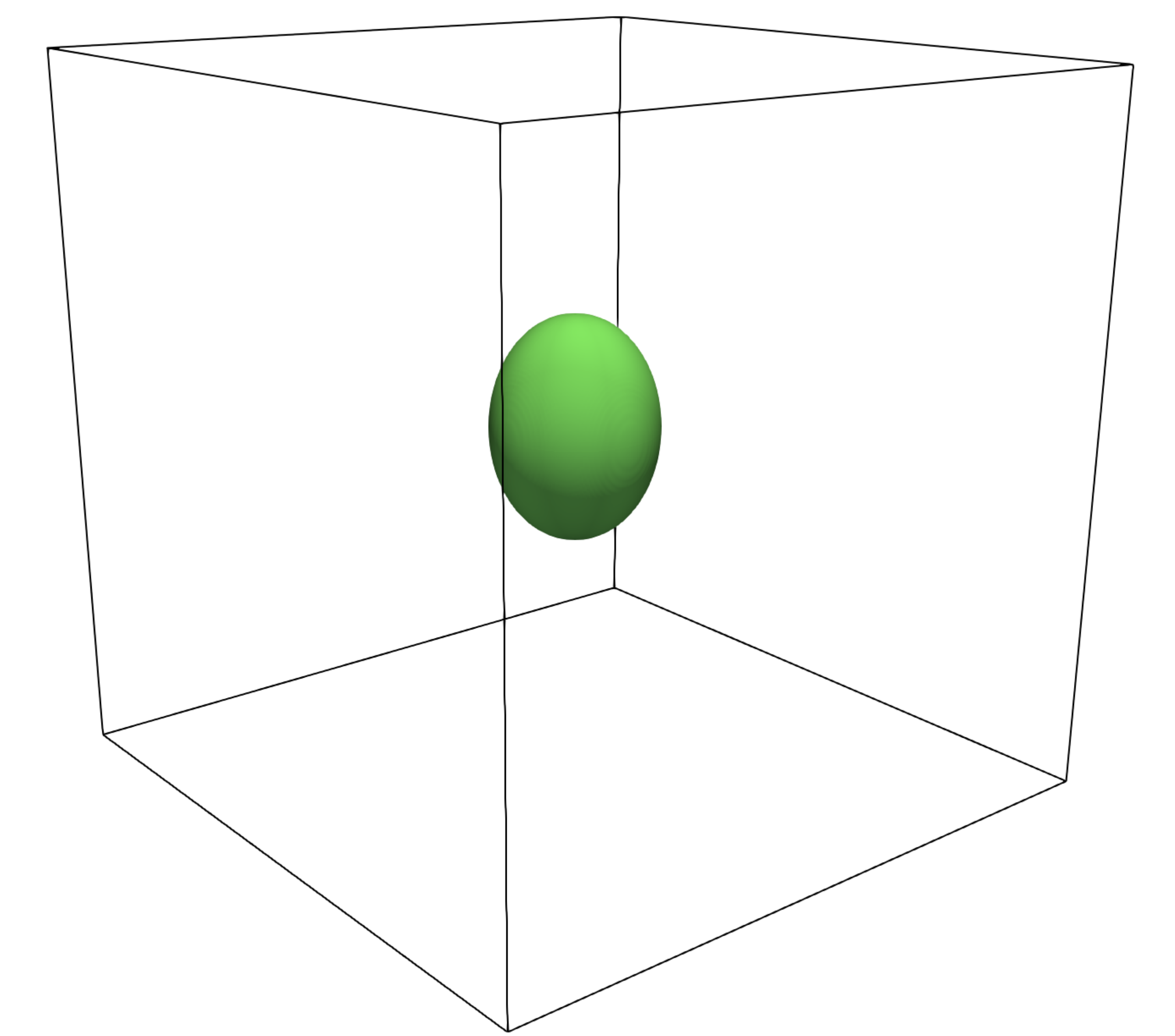}
\caption{$t=0.25$}
\end{subfigure}
\\
\begin{subfigure}{1\columnwidth} \centering
\includegraphics[trim={0cm 0cm 0cm 0cm},clip,width=0.245\columnwidth]{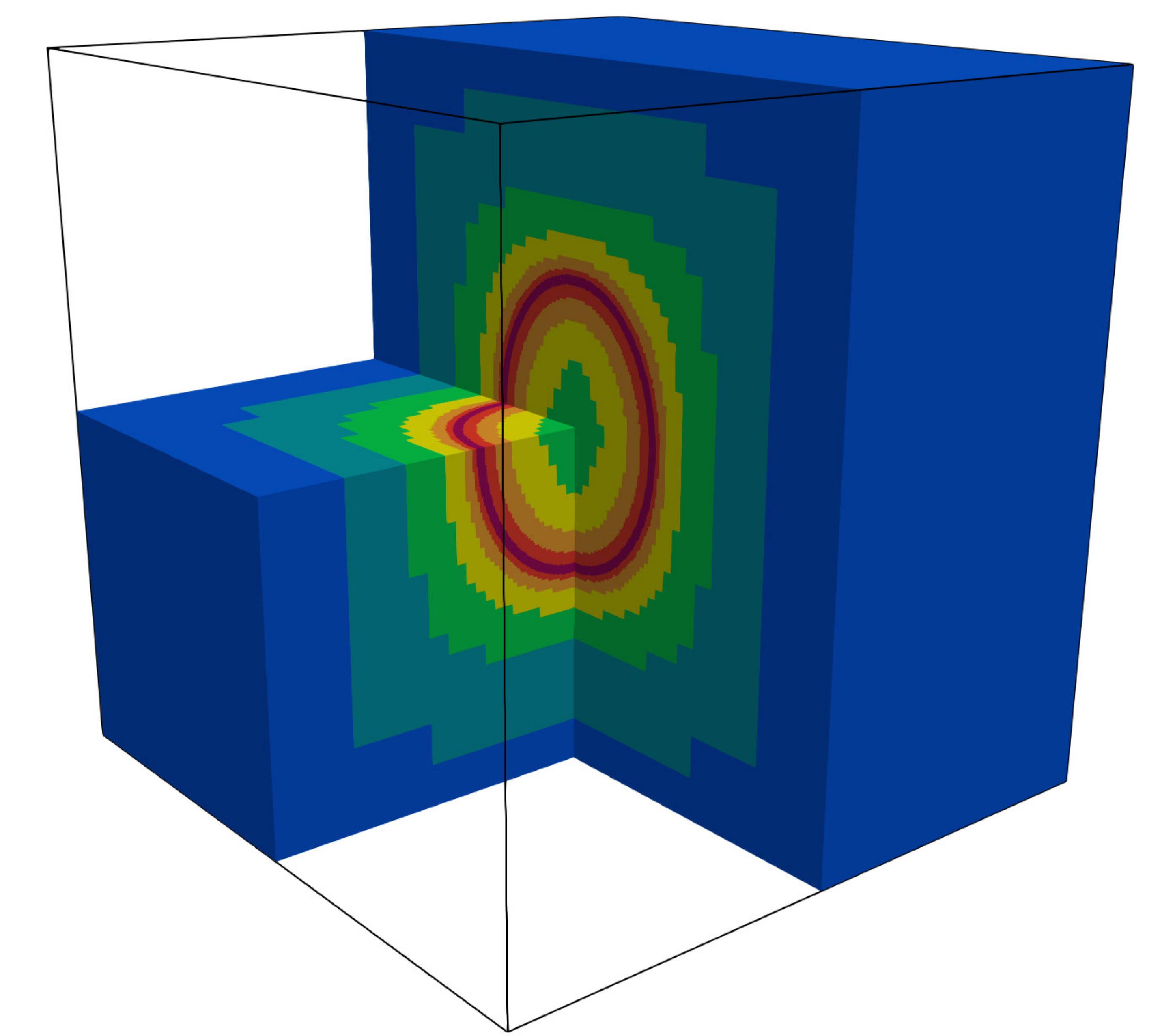}
\includegraphics[trim={0cm 0cm 0cm 0cm},clip,width=0.245\columnwidth]{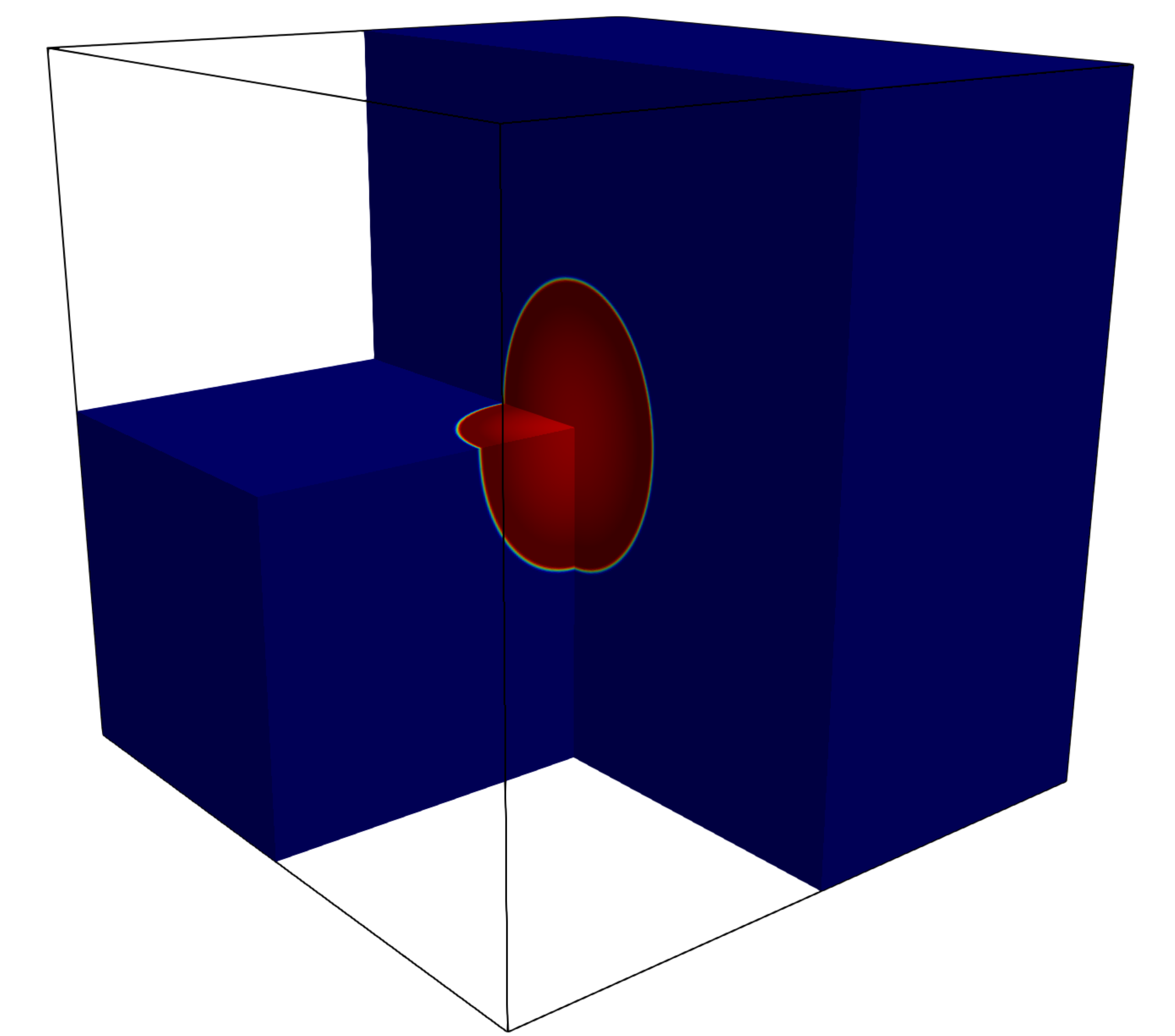}
\includegraphics[trim={0cm 0cm 0cm 0cm},clip,width=0.245\columnwidth]{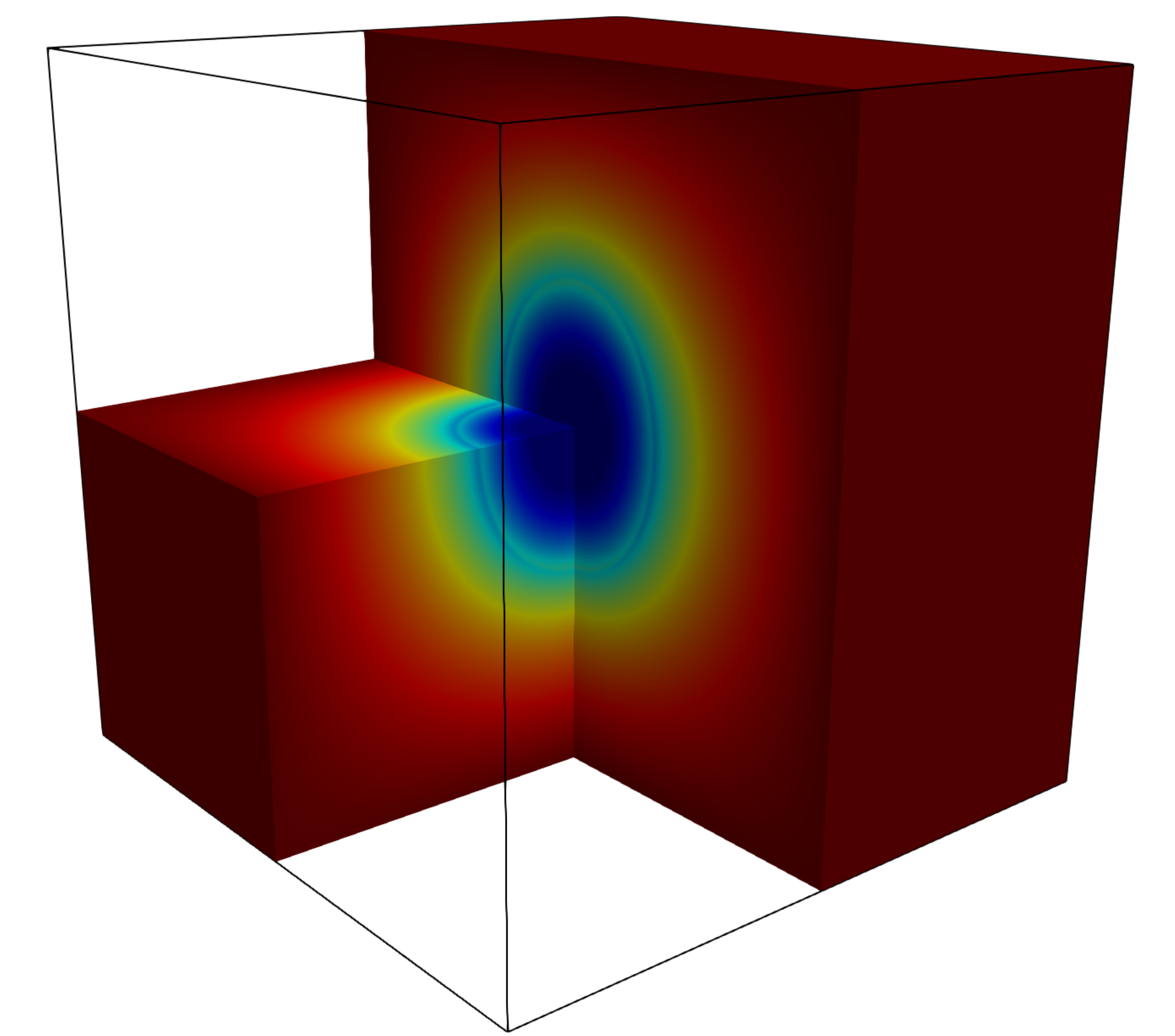}
\includegraphics[trim={0cm 0cm 0cm 0cm},clip,width=0.245\columnwidth]{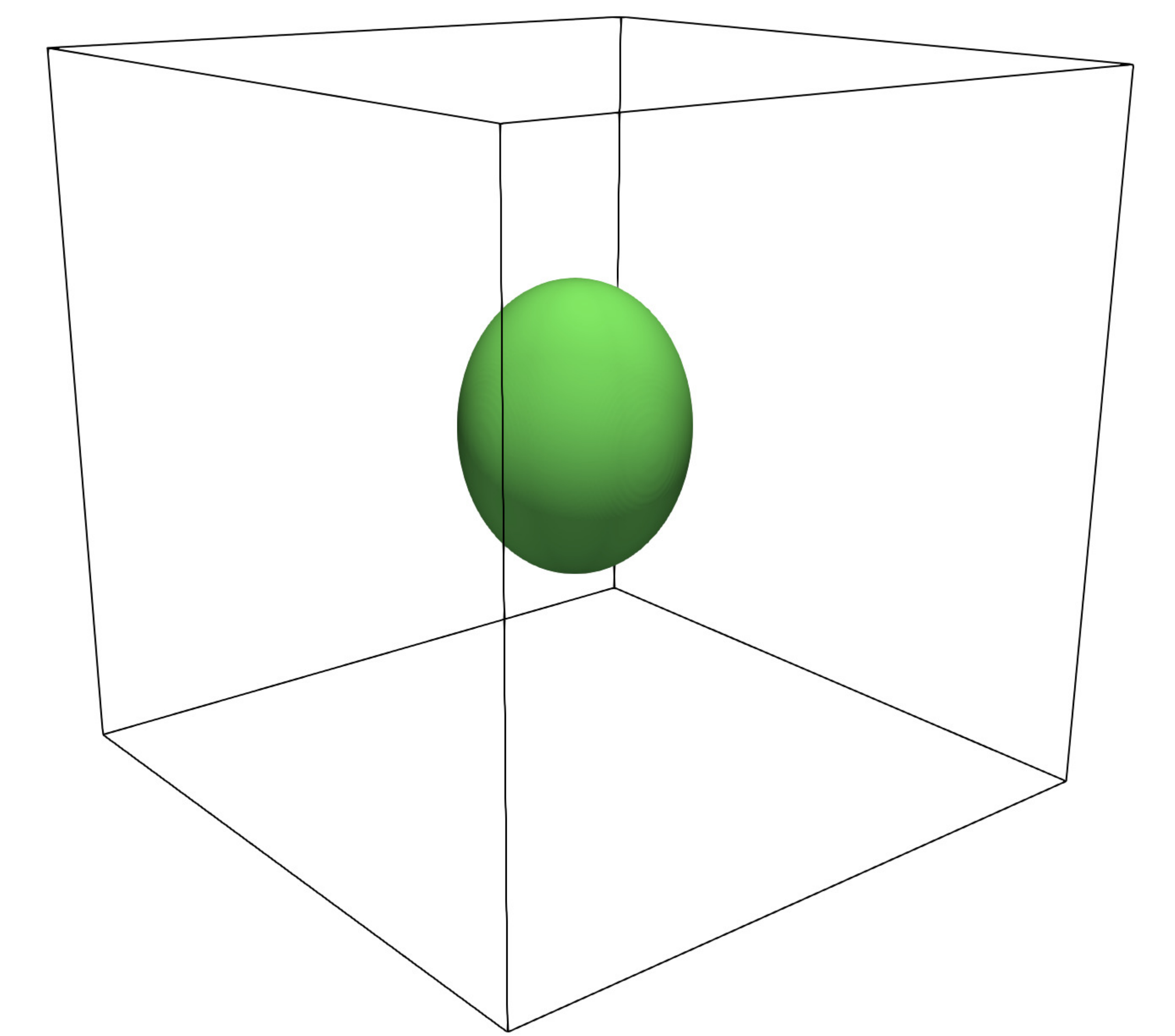}
\caption{$t=0.5$}
\end{subfigure}
\caption{Spatiotemporal evolution of an ellipsoidal tumor within a cubic domain with $\mathcal{P} = 5$, while keeping all other parameters unchanged. From left to right: adaptive THB-spline mesh ($p=3$, $\ell = 6$, $m=2$, $\alpha = 0.01$, $\beta = 0.0001$), phase-field $\phi$, nutrient concentration field $\sigma$, and the resulting tumor morphology at times 0, 0.25, and 0.5.}
\label{figure_11_final}
\end{figure}

\subsection{Tumor growth in cubic domain}\label{TG3D}

The present section extends the numerical investigation to 3D by studying the spatiotemporal evolution of tumor growth in a cubic domain defined by $\Omega = [-3, 3]^3$.
The tumor phase field is initialized as \cite{Garcke2024}:
\begin{equation*}
\phi_0(\mathbf{x})  = - \tanh\!\left( \frac{r(\mathbf{x})}{\sqrt{2}\,\epsilon} \right), \text{ where } r(\mathbf{x}) = \sqrt{\left(2x\right)^2 + \left(2y\right)^2 + \left(1.4z\right)^2} - 1.
\end{equation*}
The resulting initial configuration is an ellipsoid centered at the origin of the cube, with two semi-axes of length 0.5 and a slightly elongated third semi-axis of length 0.71429.
The mesh parameters established in the preceding study ($p =3$, $\ell = 6$, $m=2$, $\alpha = 0.01$, $\beta = 0.0001$) are retained for the problem.
Exploiting the symmetries of the problem, the deepest refinement level is reduced from $2^{10} \times 2^{10}$ to $2^{9} \times 2^{9}$, ensuring that the minimum element length ($h_e = 3/512$) remains consistent with the full domain discretization ($ h_e = 6/1024$) while reducing the overall computational efforts.

Firstly, we kept the model parameters consistent with those of the preceding 2D square domain.
The results are presented in Fig.~\ref{figure_10_final}, which shows the adaptive THB-spline mesh as it tracks the evolving interface between healthy and tumor tissues, alongside the phase field  $\phi$, the corresponding nutrient concentration $\sigma$, and the tumor morphology at $t = $ 0, 1, 2, and 3.
Starting from an ellipsoidal initial configuration, the tumor evolves into an umbrella-like morphology over time. 
While a direct one-to-one correspondence with the 2D results is not immediately apparent, a closer examination reveals that the 3D growth pattern can be understood as a superposition of two distinct growth tendencies: a fingered growth profile analogous to the 2D configuration observed in the vertical plane (Fig.~\ref{figure_7_final}), and a circular growth profile in the horizontal plane (Fig.~\ref{figure_8_final}), the interplay of which gives rise to the composite morphology as shown in Fig.~\ref{figure_10_final}.
Given the absence of reference results for this specific model configuration, the predicted tumor morphology can be qualitatively validated against the findings reported in \cite{Garcke2024}.

The influence of model parameters is further investigated by increasing the proliferation rate from 0.1 to 5, while retaining the same initial ellipsoidal condition and all other parameters unchanged. 
As in the corresponding 2D scenario, this modification significantly alters the growth dynamics. The elevated proliferation rate drives a transition toward a spheroidal tumor morphology dominated by volumetric growth rather than interfacial instability. 
The nutrient concentration field adapts accordingly, with a steeper depletion profile in the tumor interior reflecting the increased consumption associated with higher proliferation.
The corresponding results, showing the adaptive mesh, $\phi$, $\sigma$, and tumor morphology at $t = $ 0, 0.25, and 0.5 are presented in Fig.~\ref{figure_11_final}.

\begin{figure}[!t]\centering
\begin{subfigure}{0.245\columnwidth} \centering
\includegraphics[trim={0cm 0cm 0cm 0cm},clip,width=1\columnwidth]{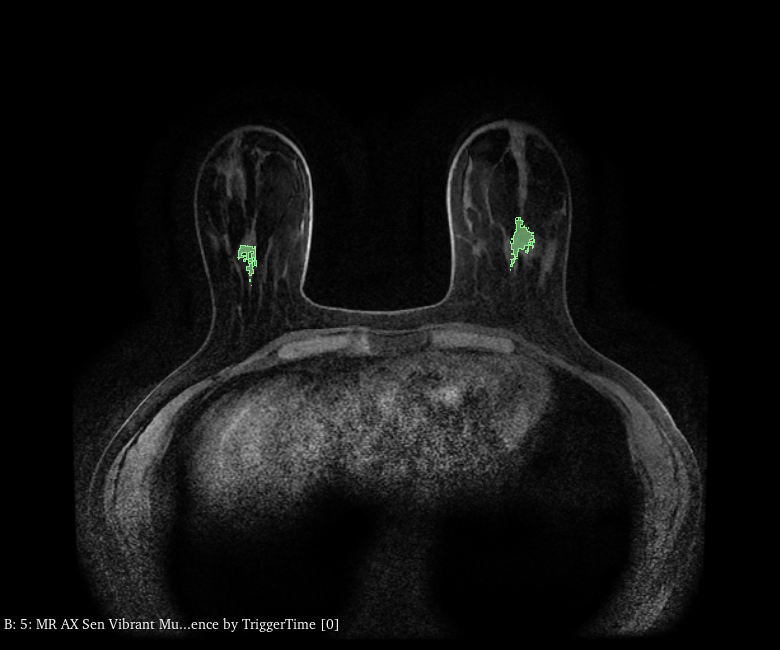}  \caption{Axial plane} \end{subfigure}
\begin{subfigure}{0.245\columnwidth} \centering
\includegraphics[trim={0cm 0cm 0cm 0cm},clip,width=1\columnwidth]{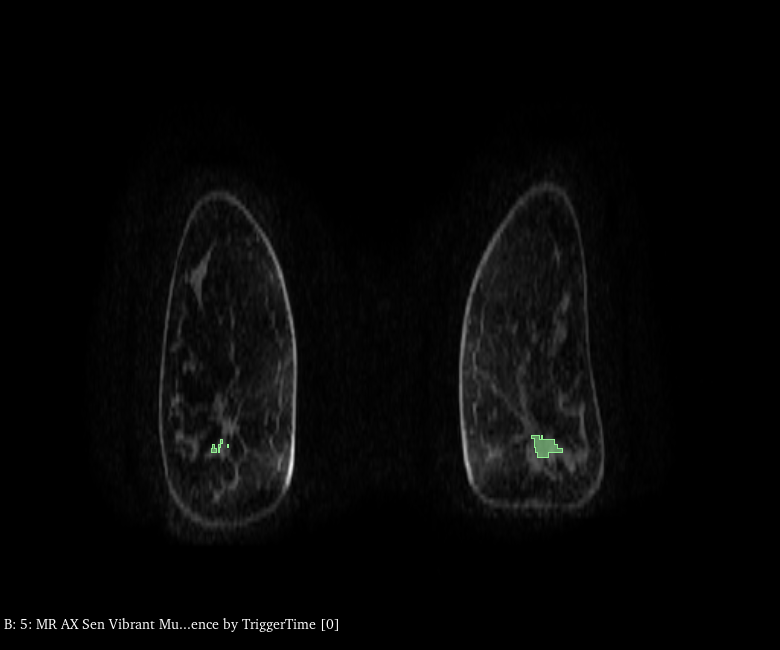}  \caption{Coronal plane} \end{subfigure}
\begin{subfigure}{0.245\columnwidth} \centering
\includegraphics[trim={0cm 0cm 0cm 0cm},clip,width=1\columnwidth]{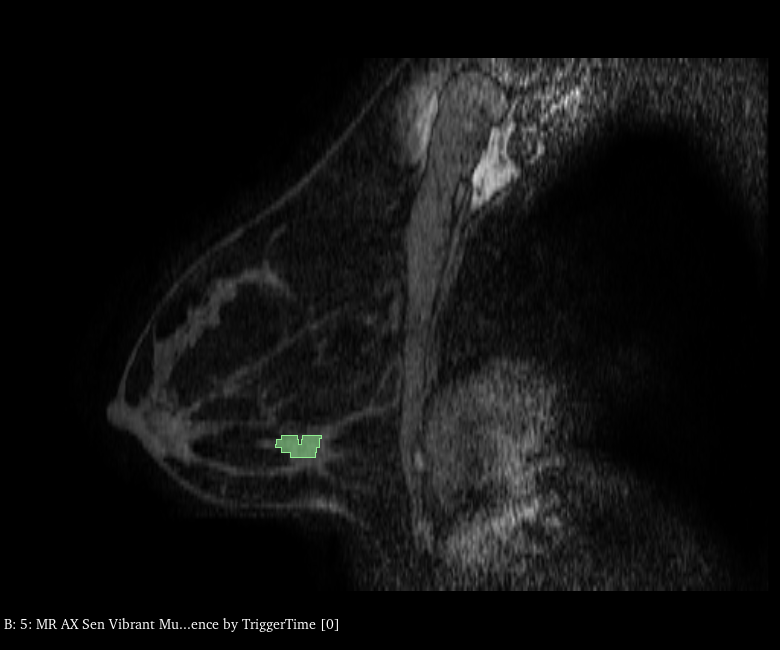}  \caption{Sagittal plane} \end{subfigure}
\begin{subfigure}{0.245\columnwidth} \centering
\includegraphics[trim={0cm 0cm 0cm 0cm},clip,width=1\columnwidth]{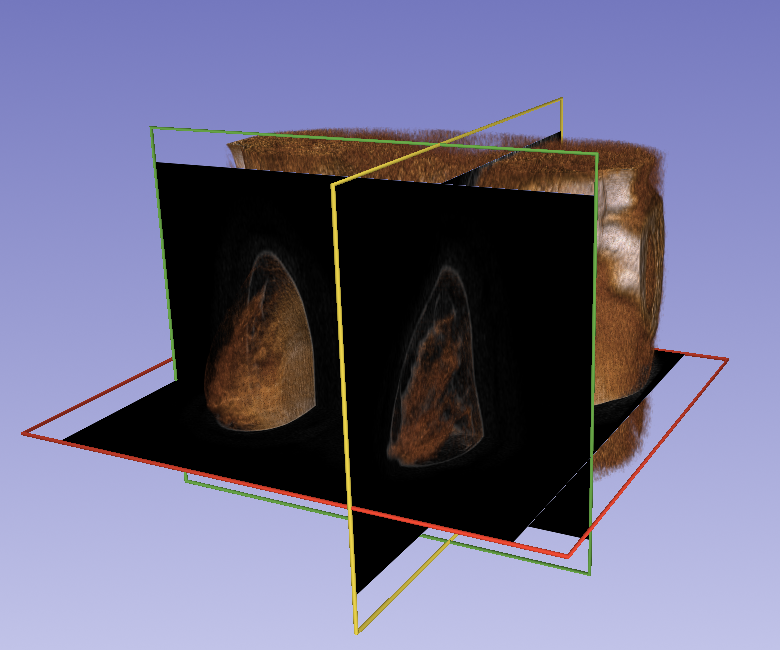}  \caption{Rendered volume} \end{subfigure}
\caption{Representative MRI slices of the breast of patient AMBL-001 visualized in three orthogonal planes: (a) axial, (b) coronal, and (c) sagittal, with the tumor mass highlighted in green. Panel (d) shows the corresponding 3D volume rendering with the three slice planes superimposed at the tumor location, where the red, green, and yellow planes correspond to the axial, coronal, and sagittal orientations, respectively. }
\label{fig:MRI_slices}
\end{figure}

\begin{figure}[!t]\centering
\begin{subfigure}{0.145\columnwidth} \centering
\includegraphics[trim={0cm 0cm 0cm 0cm},clip,width=0.99\columnwidth]{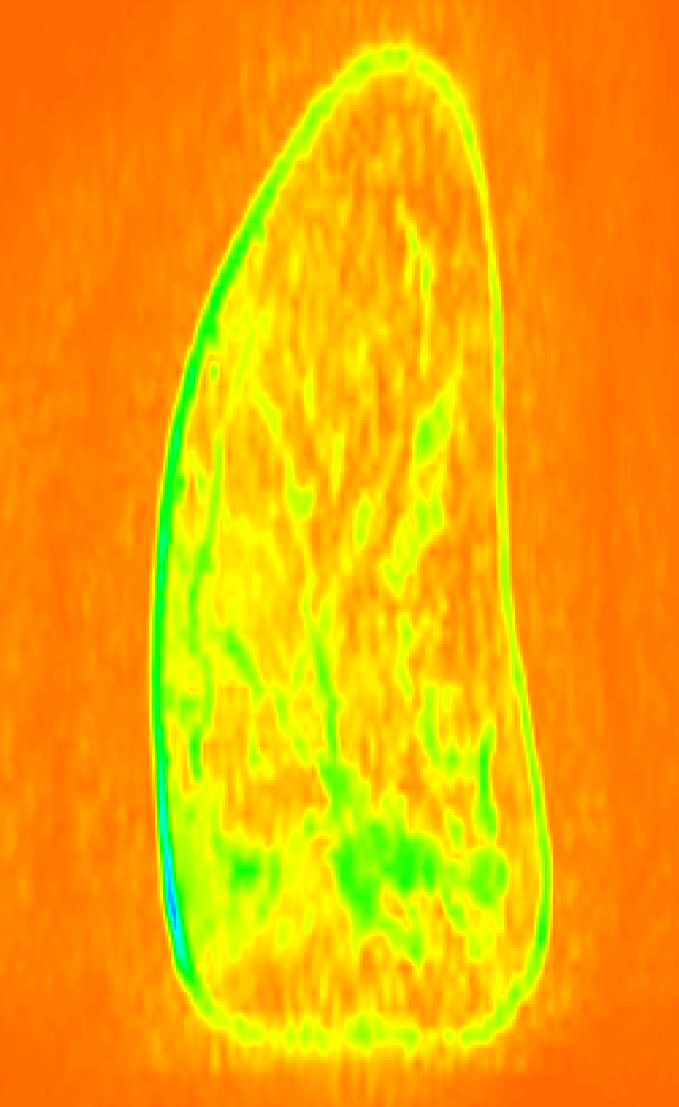}  \caption{MRI image} \end{subfigure}
\begin{subfigure}{0.145\columnwidth} \centering
\includegraphics[trim={0cm 0cm 0cm 0cm},clip,width=0.99\columnwidth]{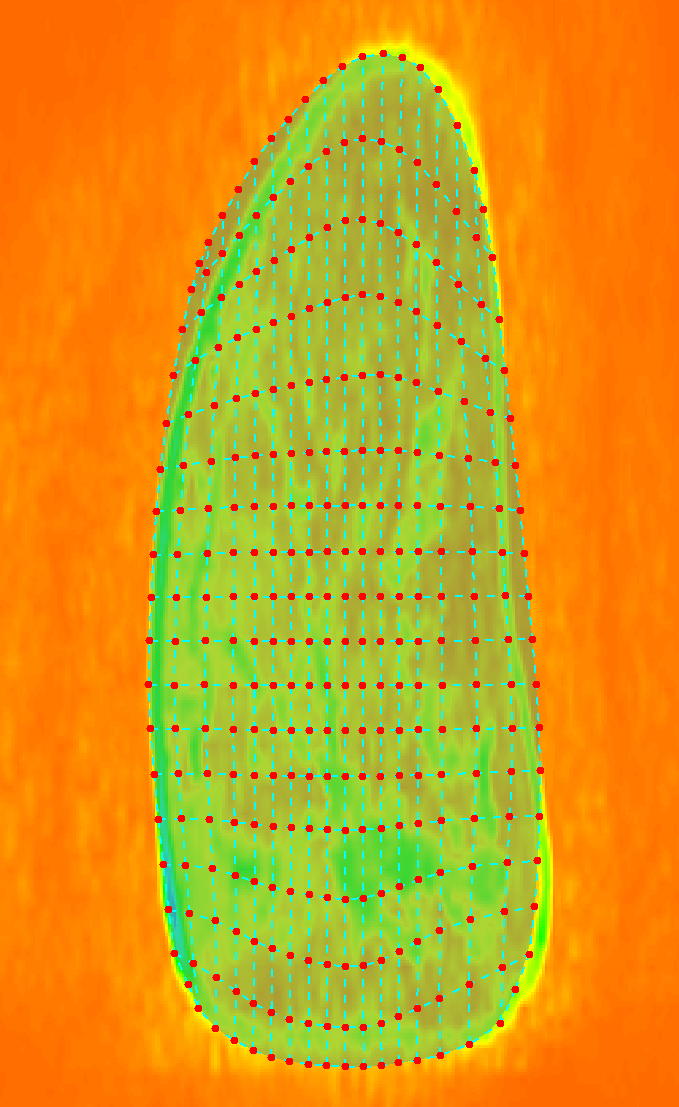}  \caption{Fitted boundary} \end{subfigure}
\begin{subfigure}{0.23\columnwidth} \centering
\includegraphics[trim={0cm 0cm 0cm 0cm},clip,width=0.85\columnwidth]{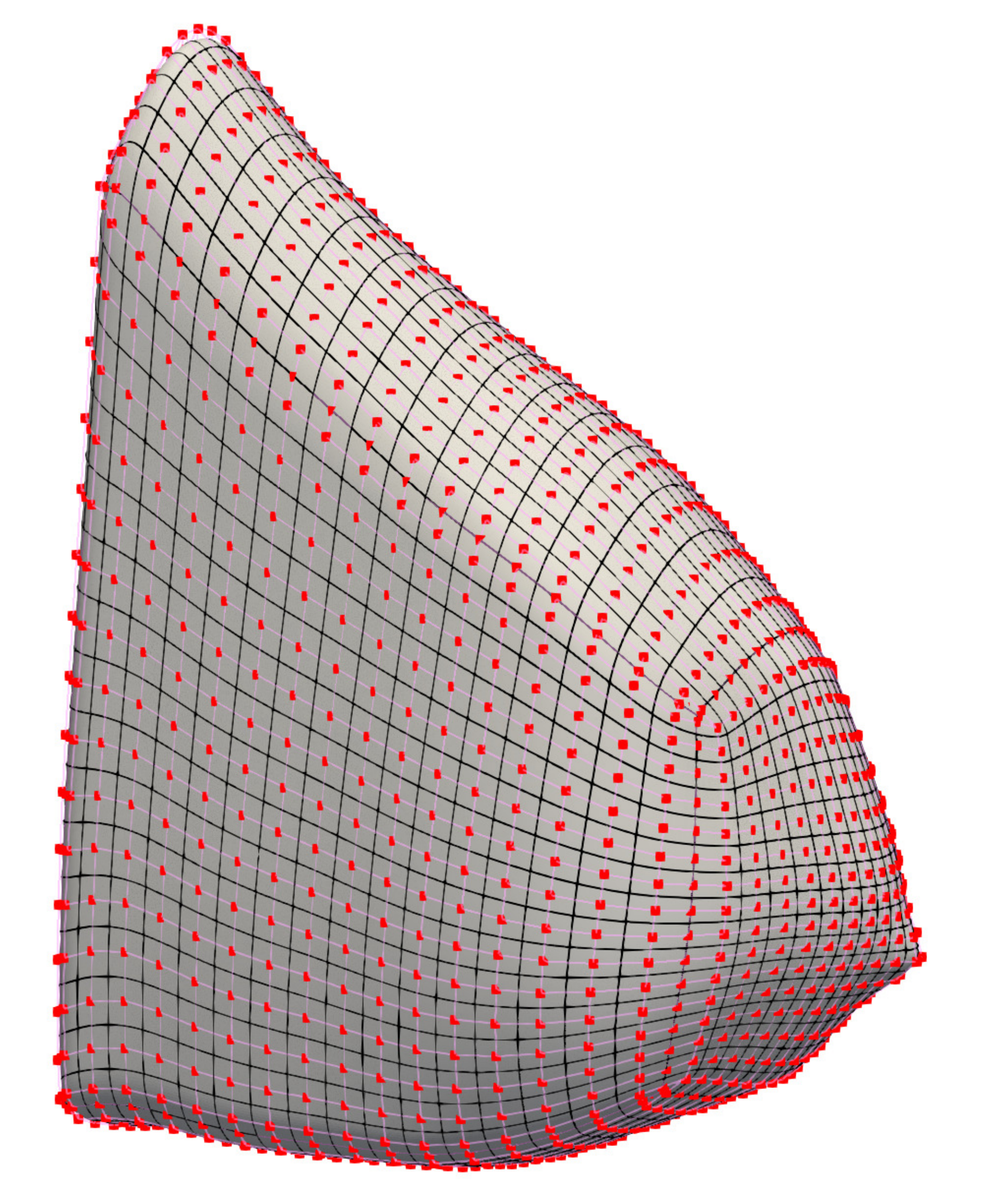}  \caption{Coarse mesh and control point net} \label{IGA_Compatible_model1} \end{subfigure}
\begin{subfigure}{0.23\columnwidth} \centering
\includegraphics[trim={0cm 0cm 0cm 0cm},clip,width=0.85\columnwidth]{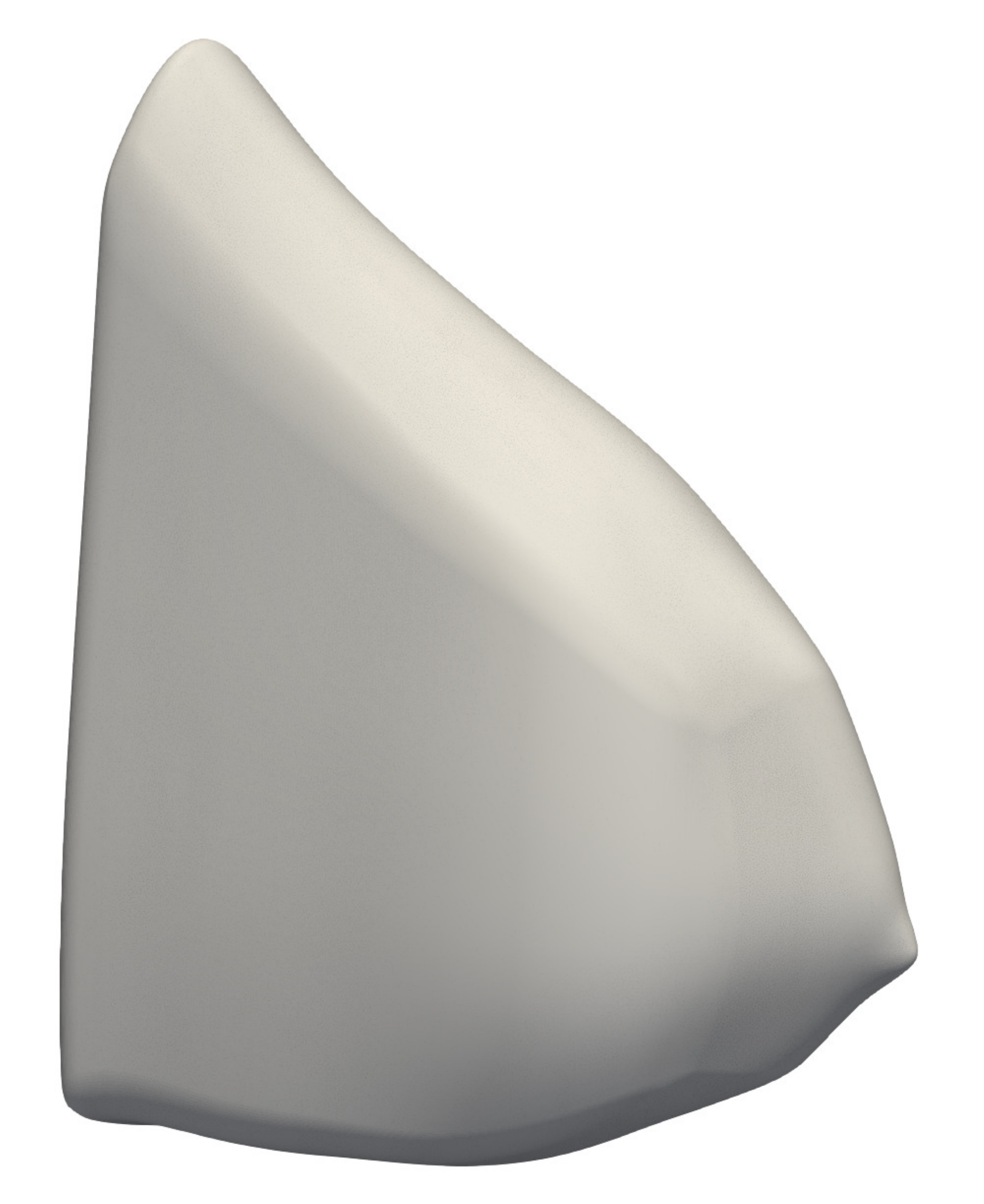}  \caption{IGA-compatible geometric model} \label{IGA_compatible_model2} \end{subfigure}
\begin{subfigure}{0.23\columnwidth} \centering
\includegraphics[trim={0cm 0cm 0cm 0cm},clip,width=0.85\columnwidth]{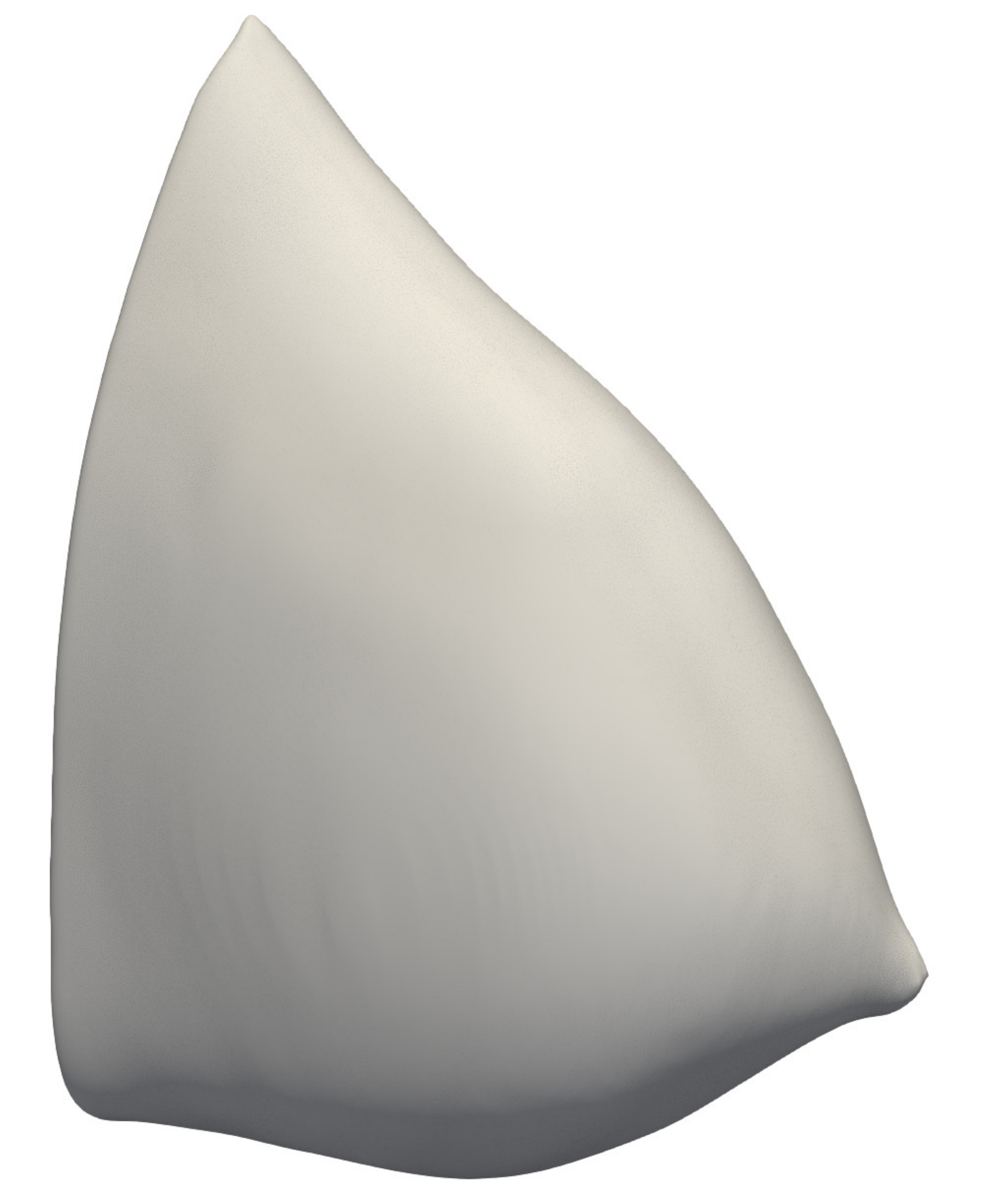}  \caption{Slicer reference model} \label{Slicer_model} \end{subfigure}
\caption{Workflow for constructing the patient-specific geometric model from MRI data. (a) Selected anatomic MRI data in the coronal plane used as the reference image. (b) Structured mesh (in light red) of control points (in red) overlaid on the MRI image, defined using quadratic B-splines ($p=2$) with 16 elements per direction (in black), interactively repositioned to trace the outer boundary of the breast tissue.  (c) Resulting 3D net of control points and corresponding quadratic B-spline elements obtained by repeating the process across 18 MRI slices and assembling the cross-sectional profiles into a full organ-scale geometry. (d) IGA-compatible volumetric model (e) Corresponding 3D surface reconstruction generated using 3D Slicer for comparative reference.}
\label{fig:3D_model}
\end{figure}

\begin{figure}[!t]\centering
\begin{tikzpicture}
\node (topimg) at (0,0)
{
\begin{subfigure}{0.49\columnwidth} \centering
\includegraphics[trim={0cm 5cm 0cm 10cm},clip,width=0.12\columnwidth]
{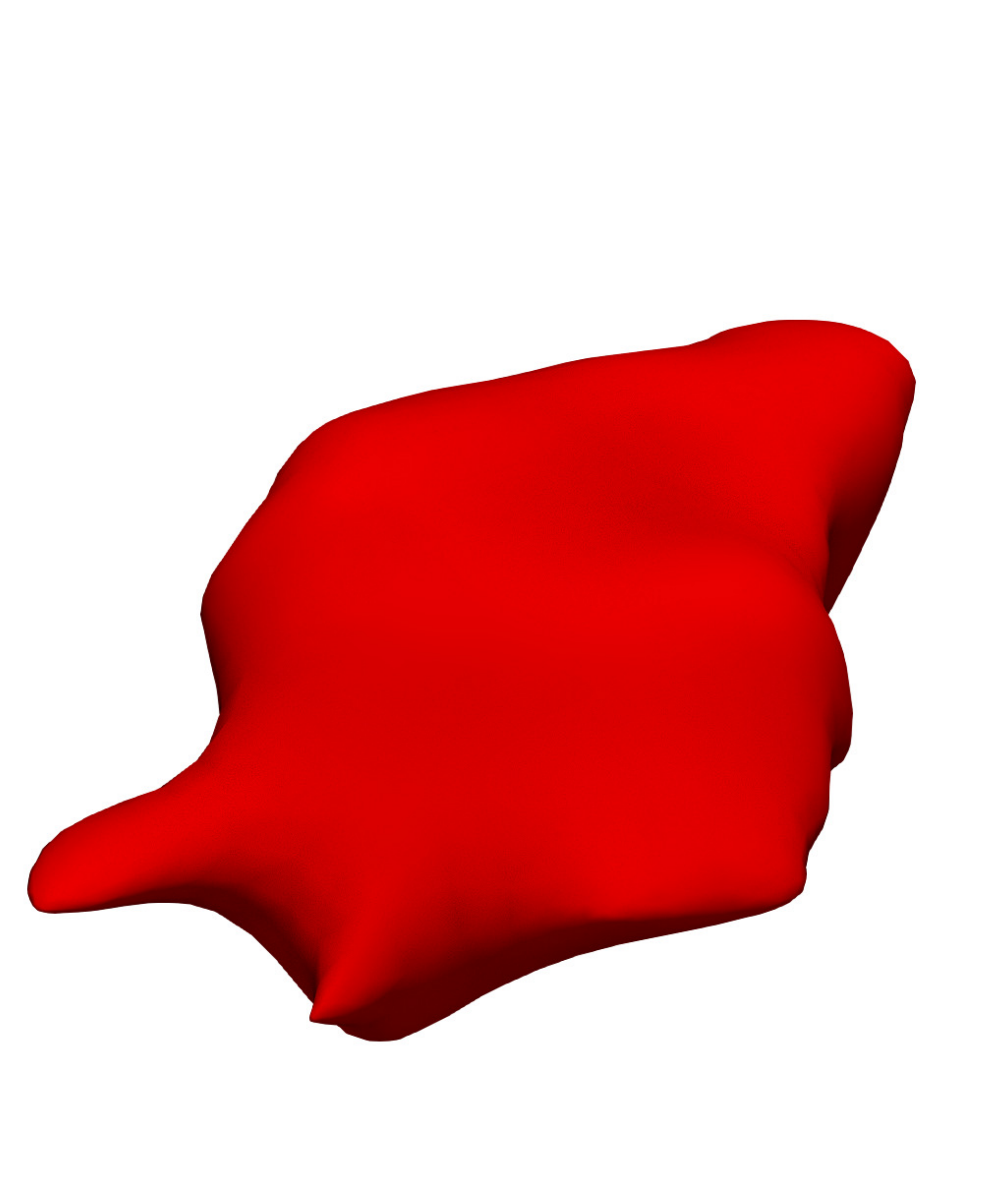}
\end{subfigure}
};
\node (bottomimg) at (0,-3.5)   %
{
\begin{subfigure}{0.49\columnwidth} \centering
\includegraphics[width=0.45\columnwidth]
{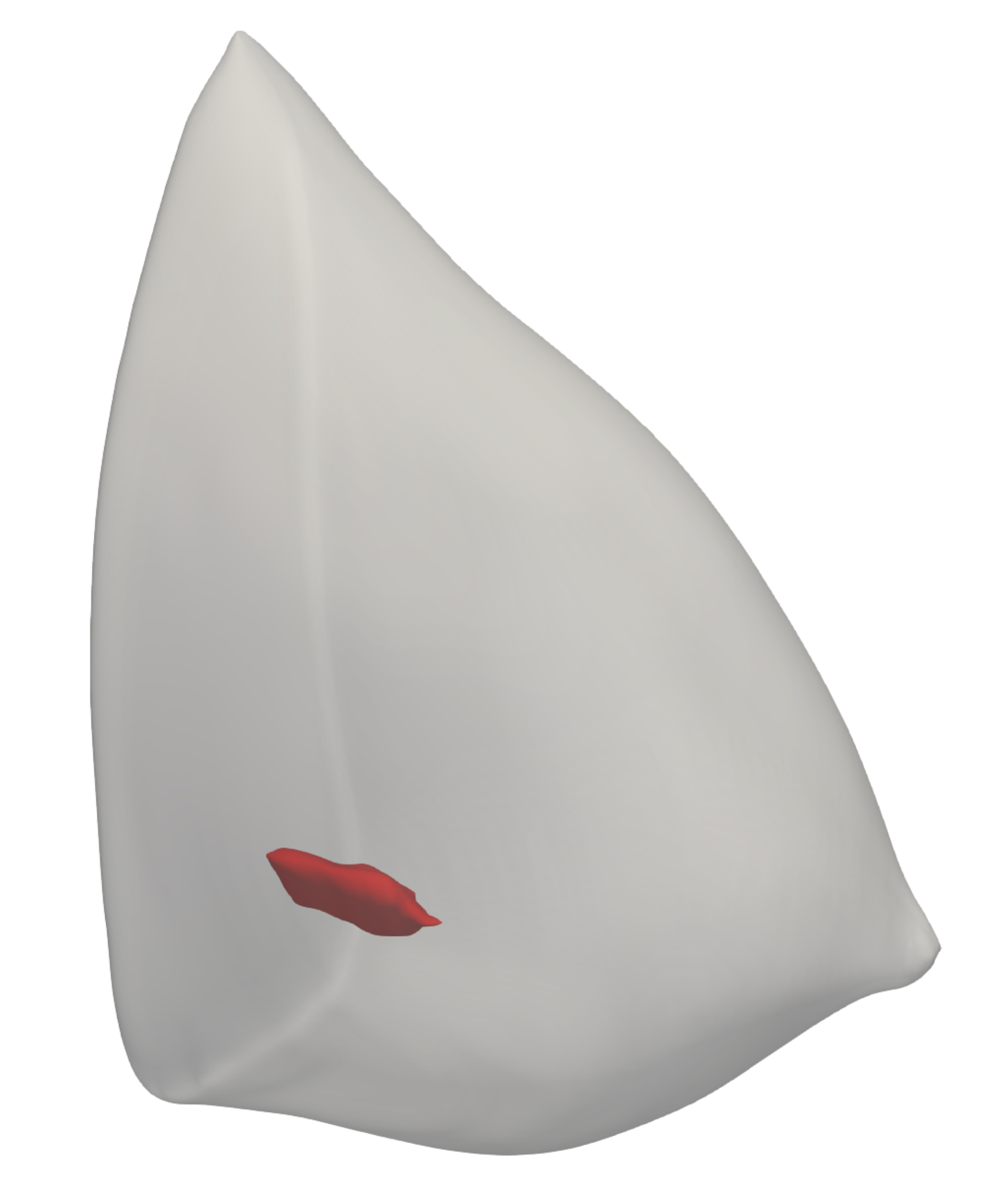}
\includegraphics[width=0.45\columnwidth]
{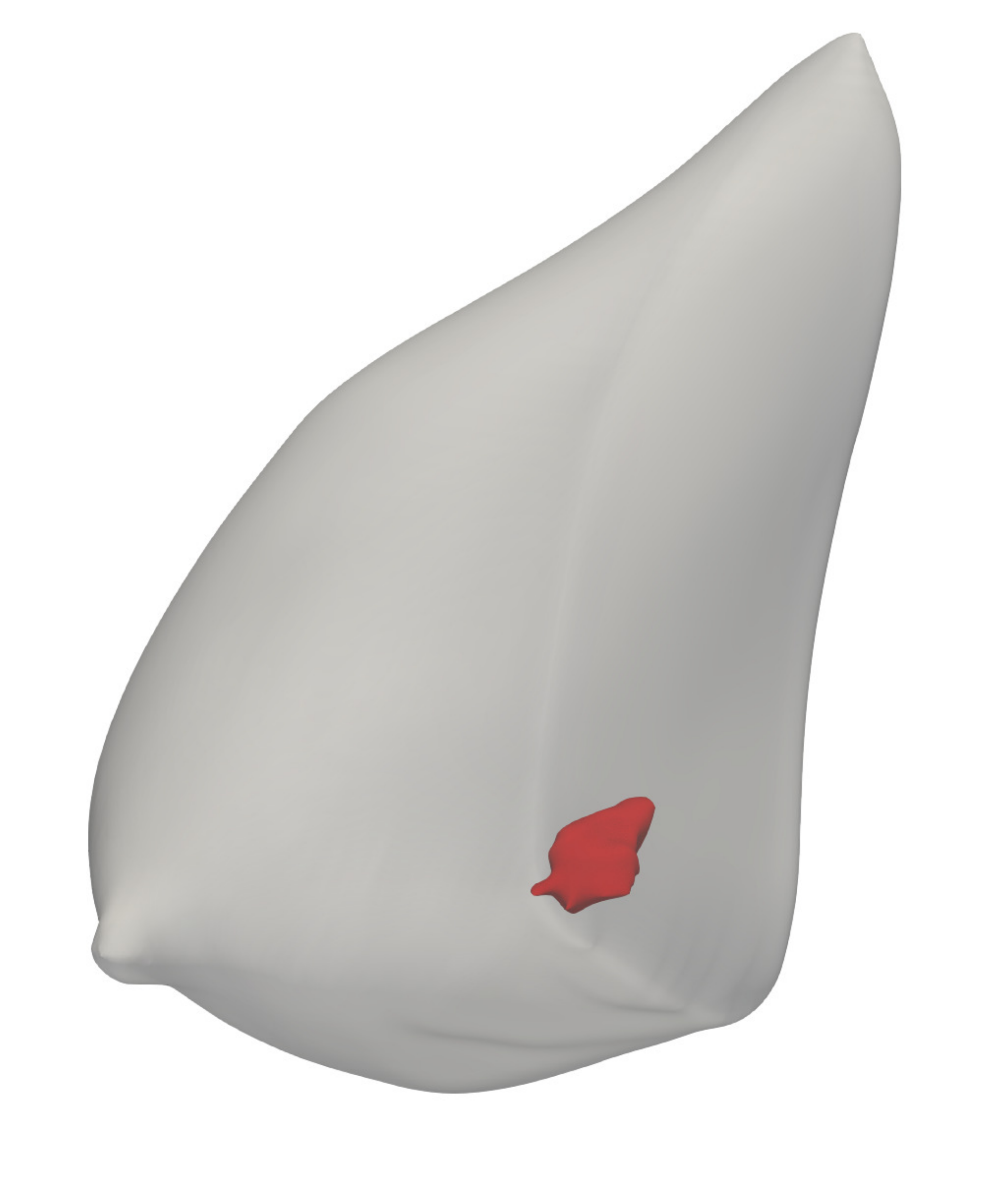}
\caption{Tumor mass and location (from MRI data)}
\end{subfigure}
};
\draw[->, thick, black]
(-2.4,-4.2) -- (-0.5,-0.5);
\draw[->, thick, black]
(2.1,-4.0) -- (0.5,-0.5);
\node[fill=white, draw=white]
at ($(topimg.south)+(1.5,0.8)$)
{Tumor mass};
\node[]
at ($(-2.6,-2.2)$)
{View 1};
\node[]
at ($(2.6,-2.2)$)
{View 2};
\end{tikzpicture}
\begin{tikzpicture}

\node (igaimg)
{
\begin{subfigure}{0.47\columnwidth}
\centering
\includegraphics[
trim={0cm 4cm 0cm 0cm},
clip,
width=0.85\columnwidth]
{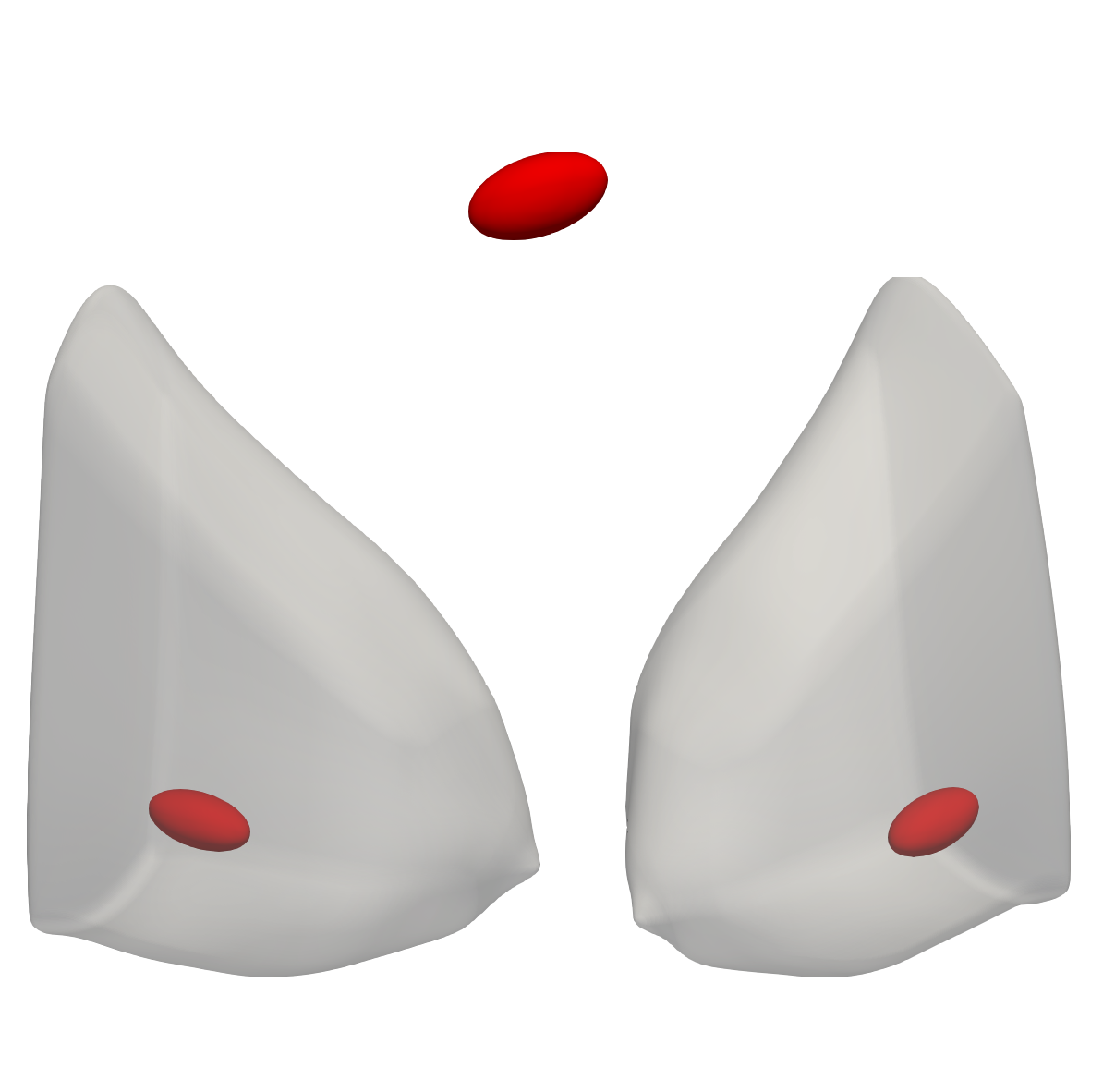}
\caption{IGA-ready breast model and approximated tumor mass}
\end{subfigure}
};

\draw[->, thick, black]
(-2.0,-1.5) -- (-0.5,1.6);

\draw[->, thick, black]
(2.1,-1.5) -- (0.3,1.6);

\node
at (-2.3,0.4)
{View 1};

\node
at (2.1,0.4)
{View 2};

\node[fill=white]
at (1.2,2.5)
{Tumor mass};

\end{tikzpicture}

\caption{
Breast model reconstructed from MRI data in two views, with the tumor mass shown in red.
Note that the image smoothing during reconstruction results in the removal of fine tumor branches.
}
\label{fig:tumorMass}
\end{figure}

\subsection{Tumor growth in patient-specific organ-scale breast model}\label{sec:BrCa_patient_spec}

To complete our simulation study, we test our THB spline isogeometric approach on an organ-scale scenario for untreated breast cancer (BCa). This disease remains the most prevalent cancer among women worldwide, with 2,308,897 new cases (11.6\% incidence) and 665,684 deaths (6.9\% mortality) reported globally in 2022 \cite{Bray2024}. 
While advances in screening and treatment have improved survival rates \cite{Giaquinto2024, Peintinger2019}, interest is increasingly shifting towards personalized disease management which demands predictive tools capable of capturing tumor growth dynamics at clinically relevant scales. 
This motivates the application of the CH-based tumor growth model developed in this work to BCa, which, to the best of our knowledge, represents the first attempt to model BCa growth using a higher-order CH-based formulation in a realistic, organ-scale, three-dimensional setting.

A patient-specific geometry is developed using clinical MRI data obtained from The Cancer Imaging Archive \cite{Daniels2024AdvancedMRI}.
The dataset is publicly available and anonymized, with all personal identifiers removed in compliance with data privacy standards. 
The selected patient (AMBL-001) was 47.4 years of age at the time of MRI acquisition, with a tumor detected in the left breast and another in the right breast.
Figure~\ref{fig:MRI_slices} shows an MRI slice in three orthogonal planes (axial, coronal, and sagittal) with the tumor region highlighted in green. 
Since the purpose of our study is to demonstrate the viability of our THB-spline isogeometric method in a realistic organ-scale scenario, the subsequent geometric modeling and simulation are restricted to the left breast domain.

We construct a patient-specific geometric model directly from the MRI scans using a dedicated in-house graphical user interface (GUI) developed for this purpose. 
First, a coronal anatomic $T_1$-weighted MRI slice is loaded into the GUI, a structured rectangular mesh of control points is overlaid, and the control points are interactively repositioned to trace the outer boundary of the breast tissue in real time. 
Once a satisfactory approximation of the boundary is achieved, the geometry is smoothed to yield a smooth B-spline curve representation of the cross-sectional profile. 
This process is repeated independently for 18 MRI slices distributed along the anterior-posterior longitudinal axis of the breast, providing a sufficiently dense set of cross-sectional profiles to capture the 3D geometry closely. 
The individual slice geometries are subsequently assembled and patched together to construct a volumetric 3D IGA-compatible geometric model, defined in terms of control points and knot vectors, as shown in Fig.~\ref{fig:3D_model}.
This figure also shows a comparison between the surface model reconstructed using 3D Slicer (Fig.~\ref{Slicer_model}) and the IGA-compatible geometric model (Figs.~\ref{IGA_Compatible_model1} and \ref{IGA_compatible_model2}) developed for the simulation, demonstrating close agreement between the two representations.

\begin{figure}[!t]\centering
\begin{subfigure}{\columnwidth}
\begin{subfigure}{0.328\columnwidth} \centering
\includegraphics[trim={0cm 8cm 0cm 0cm},clip,width=0.45\columnwidth]{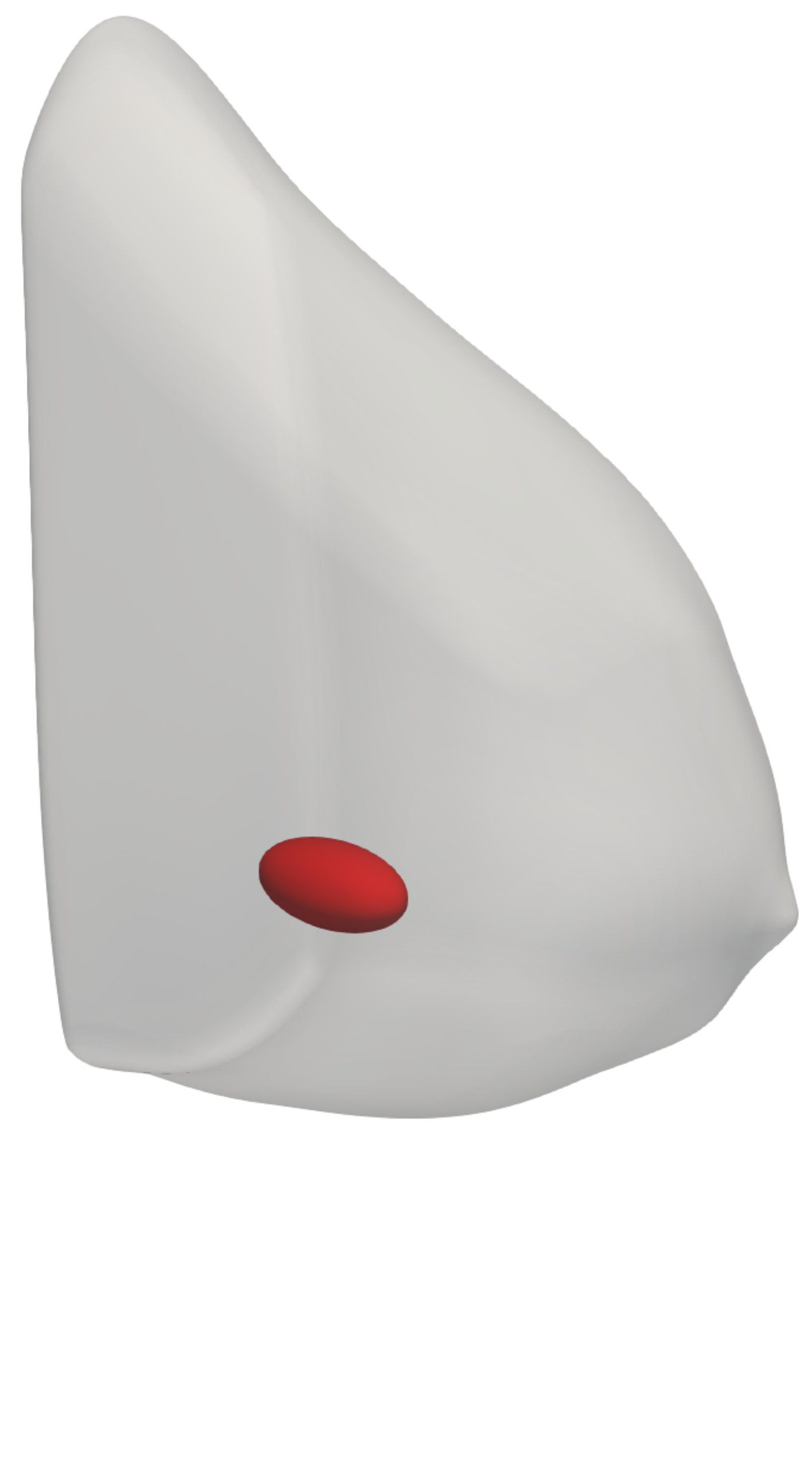}  
\includegraphics[trim={0cm 8cm 0cm 0cm},clip,width=0.45\columnwidth]{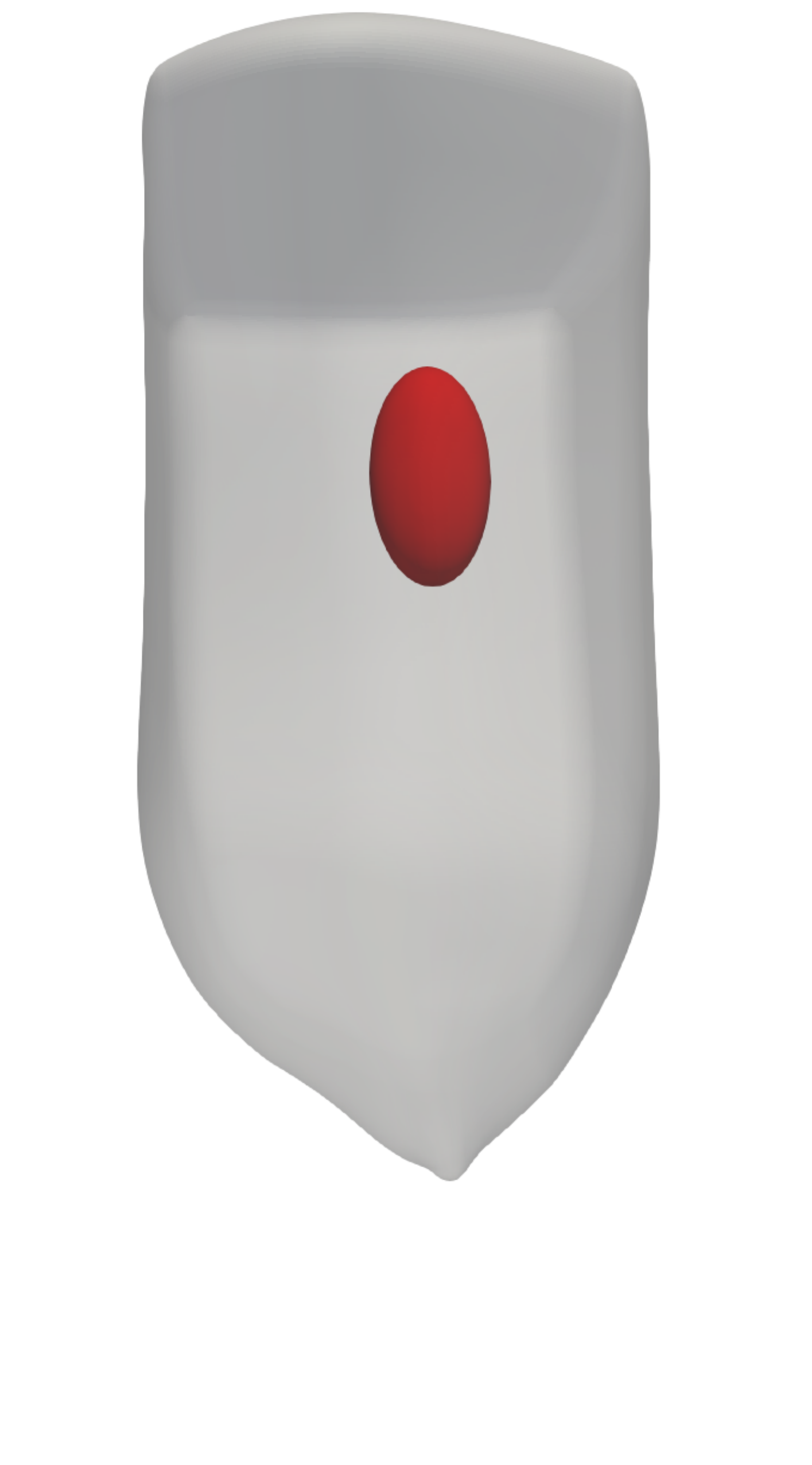}  
\end{subfigure}
\begin{subfigure}{0.328\columnwidth} \centering
\includegraphics[trim={0cm 8cm 0cm 0cm},clip,width=0.45\columnwidth]{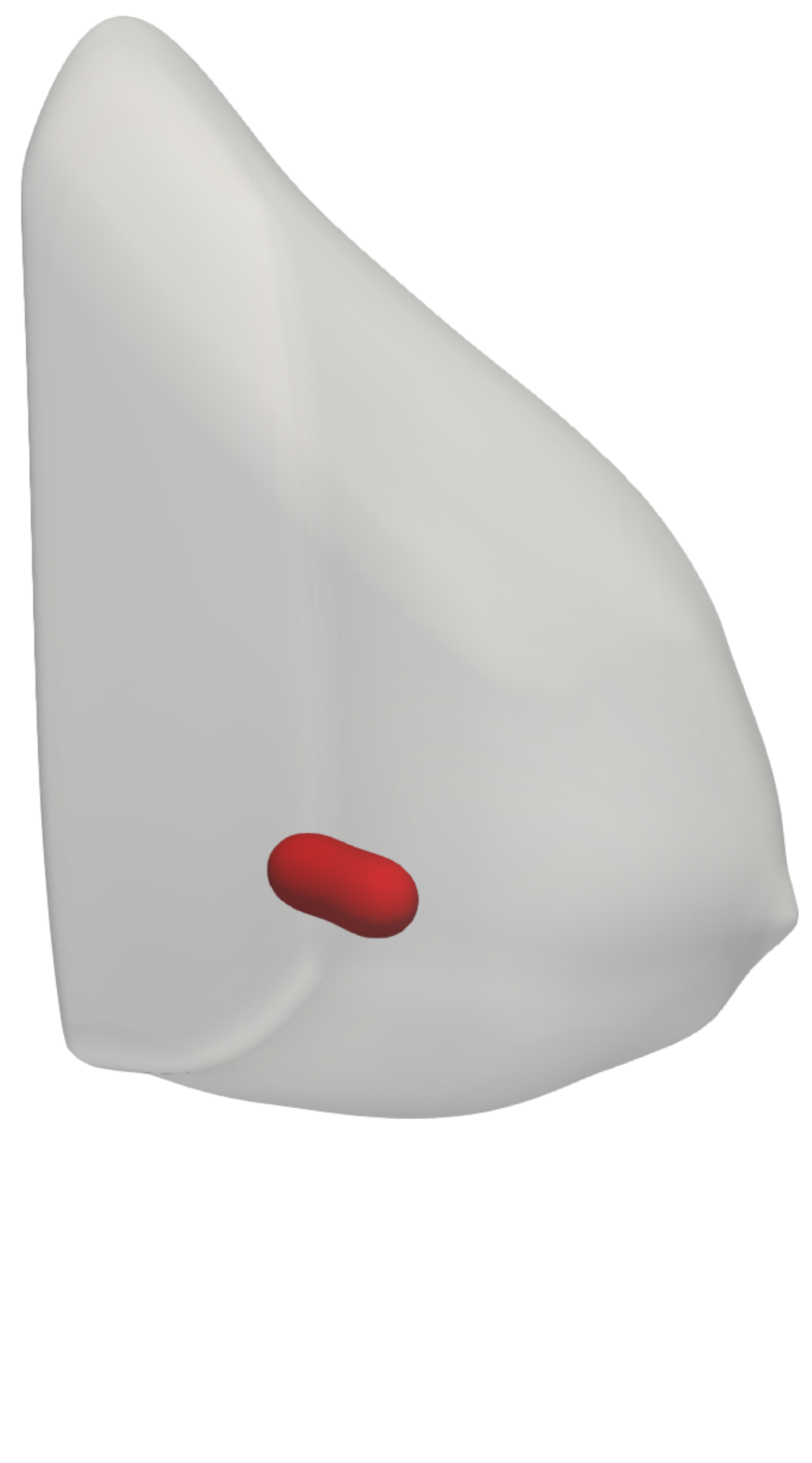}  
\includegraphics[trim={0cm 8cm 0cm 0cm},clip,width=0.45\columnwidth]{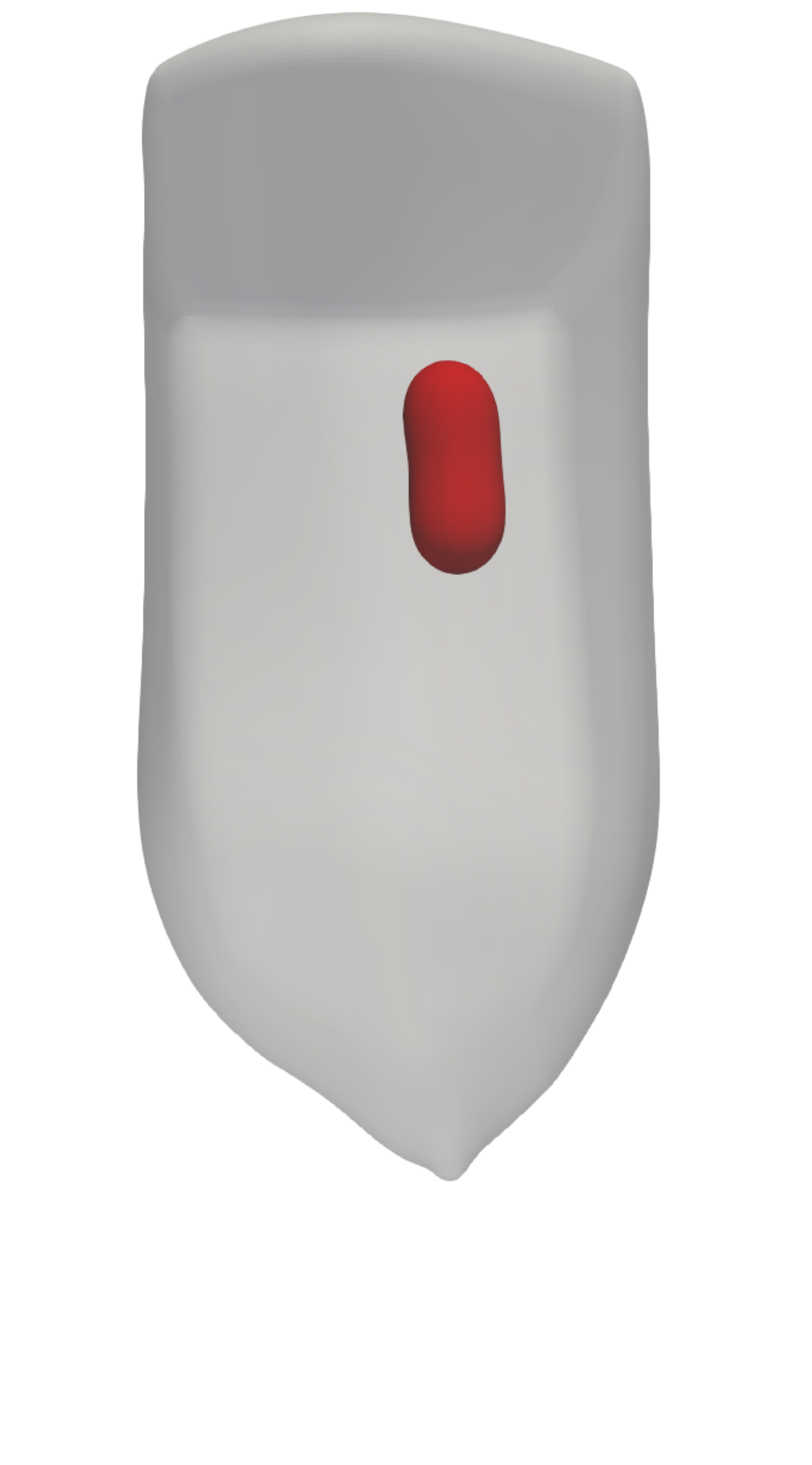}  
\end{subfigure}
\begin{subfigure}{0.328\columnwidth} \centering
\includegraphics[trim={0cm 8cm 0cm 0cm},clip,width=0.45\columnwidth]{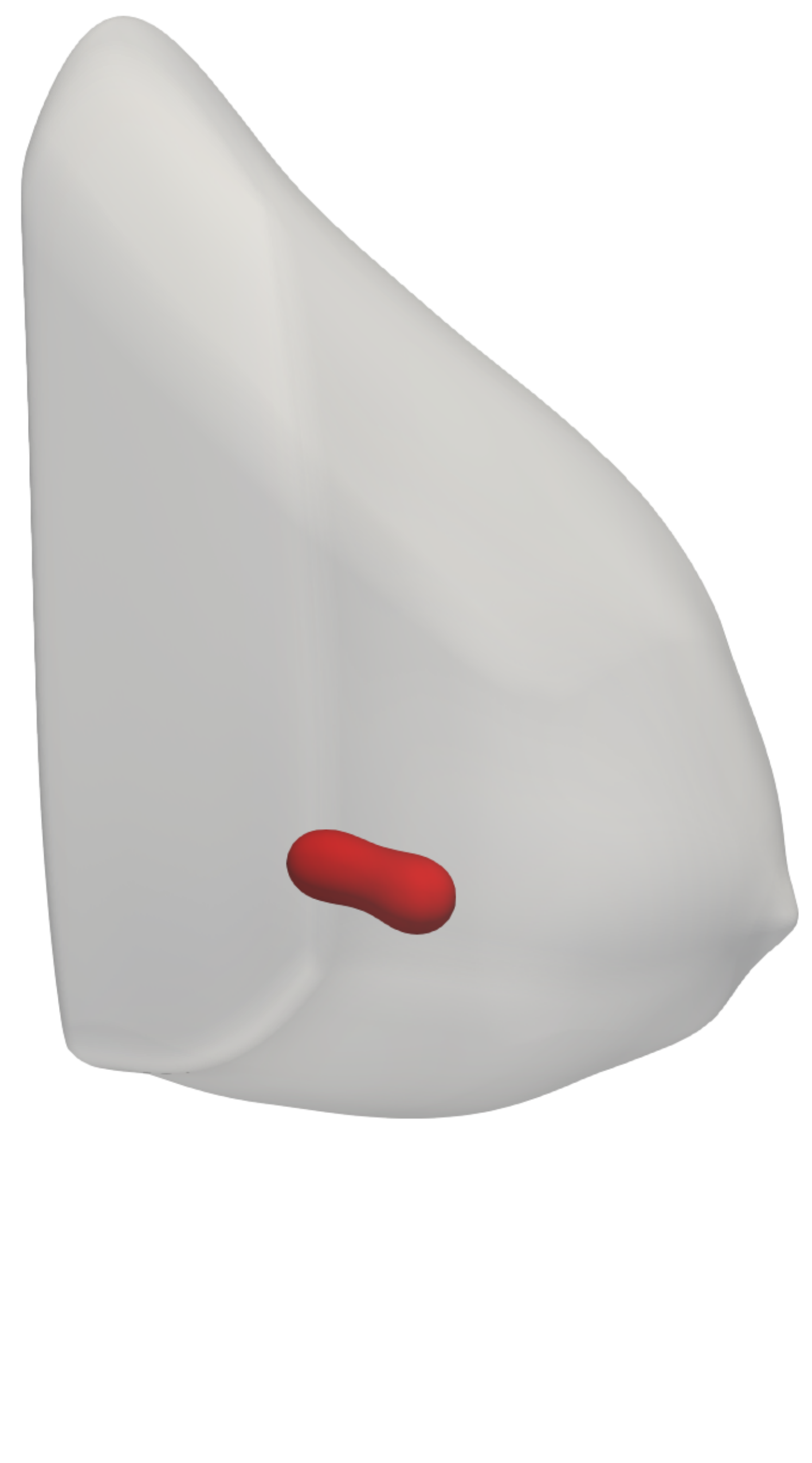}  
\includegraphics[trim={0cm 8cm 0cm 0cm},clip,width=0.45\columnwidth]{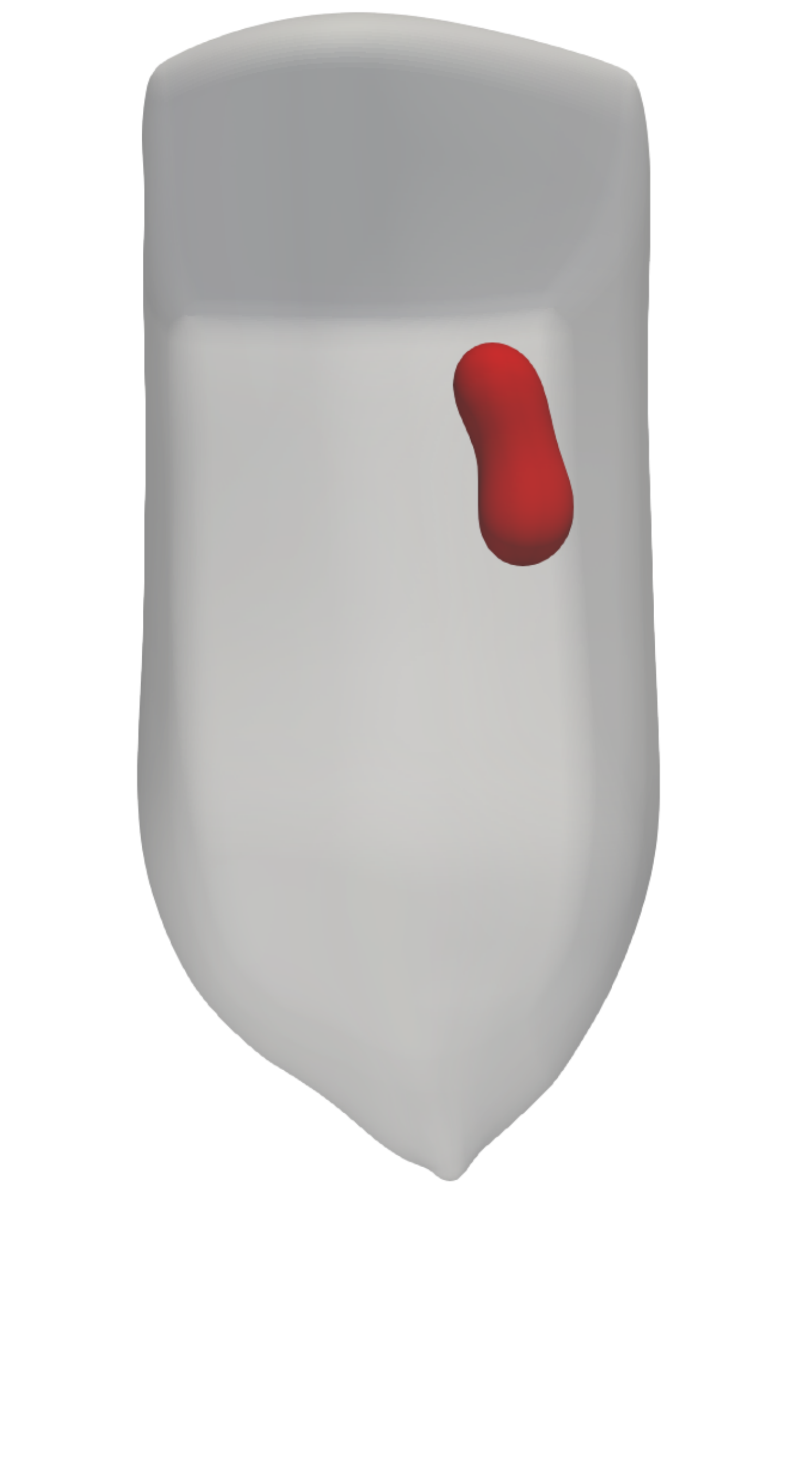}  
 \end{subfigure}
\caption{Tumor morphology} \end{subfigure}
\\
\begin{subfigure}{\columnwidth}
\begin{subfigure}{0.328\columnwidth} \centering
\includegraphics[trim={0cm 8cm 0cm 0cm},clip,width=0.55\columnwidth]{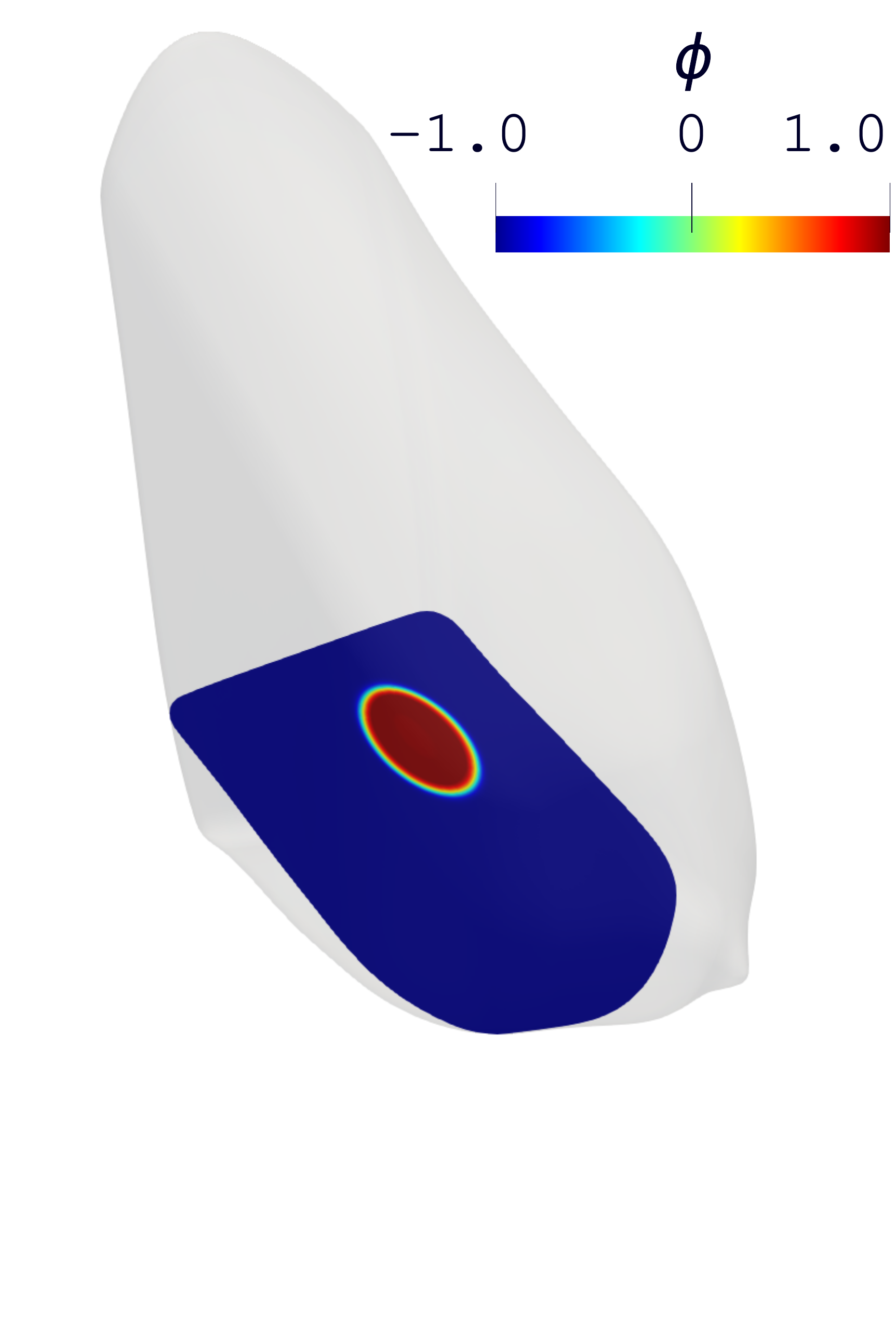}  
\includegraphics[trim={0cm 8cm 0cm 0cm},clip,width=0.35\columnwidth]{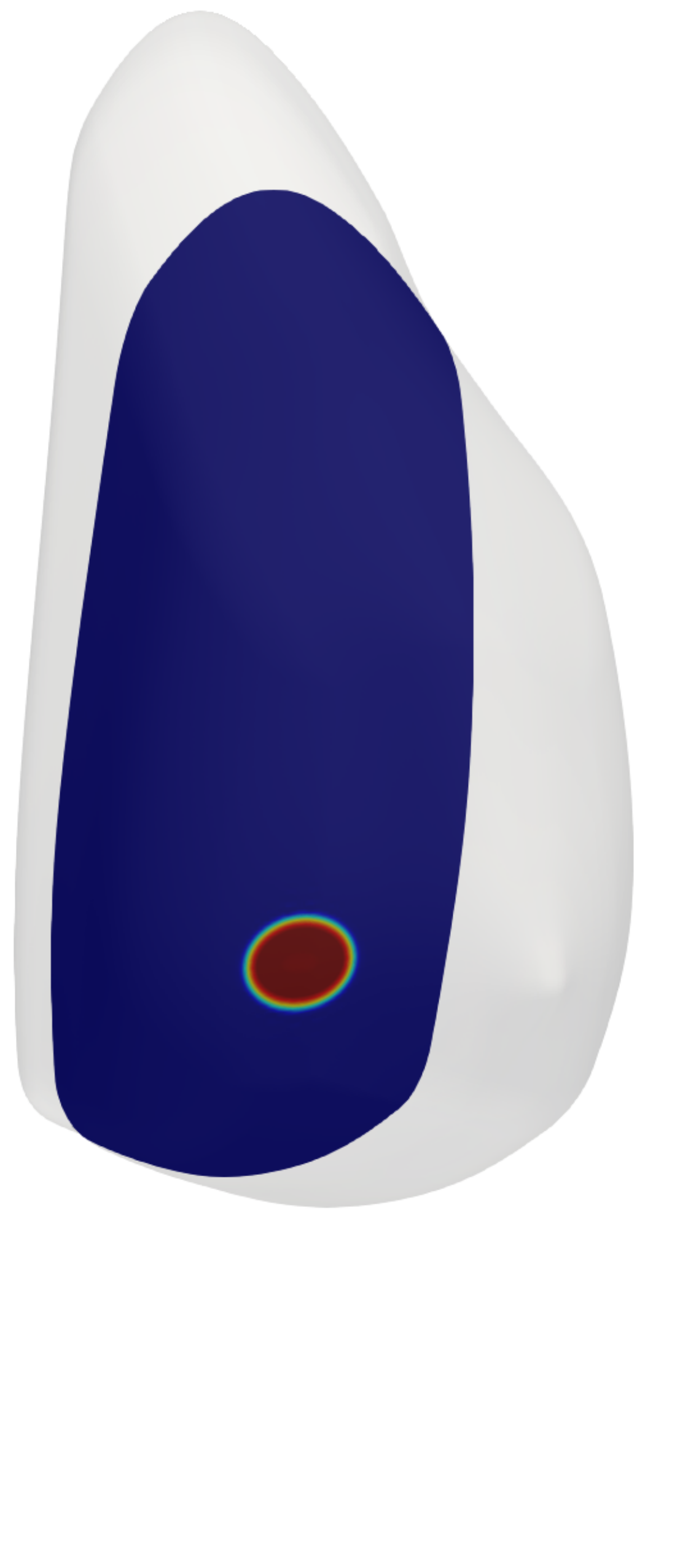}  
 \end{subfigure}
\begin{subfigure}{0.328\columnwidth} \centering
\includegraphics[trim={0cm 8cm 0cm 0cm},clip,width=0.55\columnwidth]{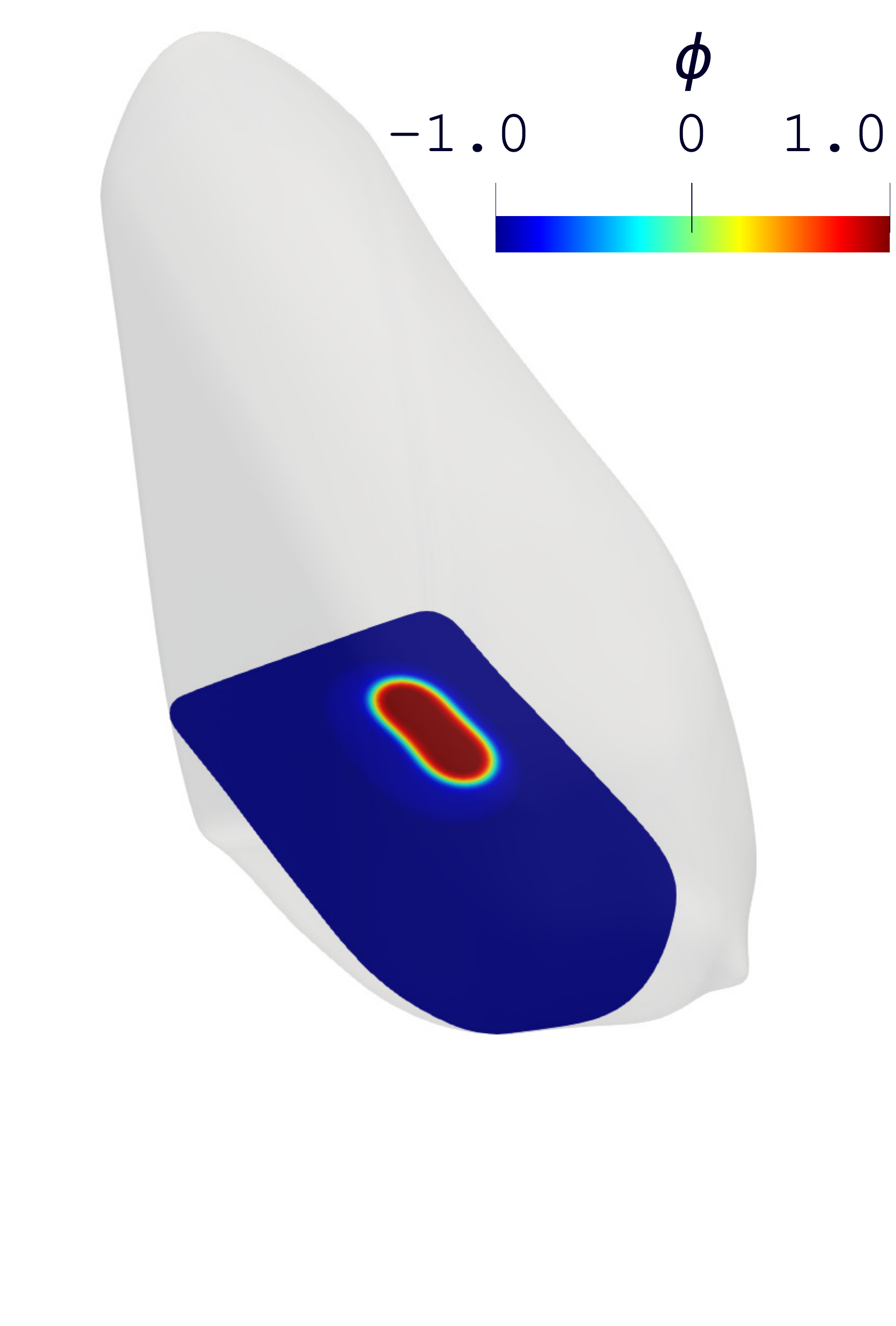}  
\includegraphics[trim={0cm 8cm 0cm 0cm},clip,width=0.35\columnwidth]{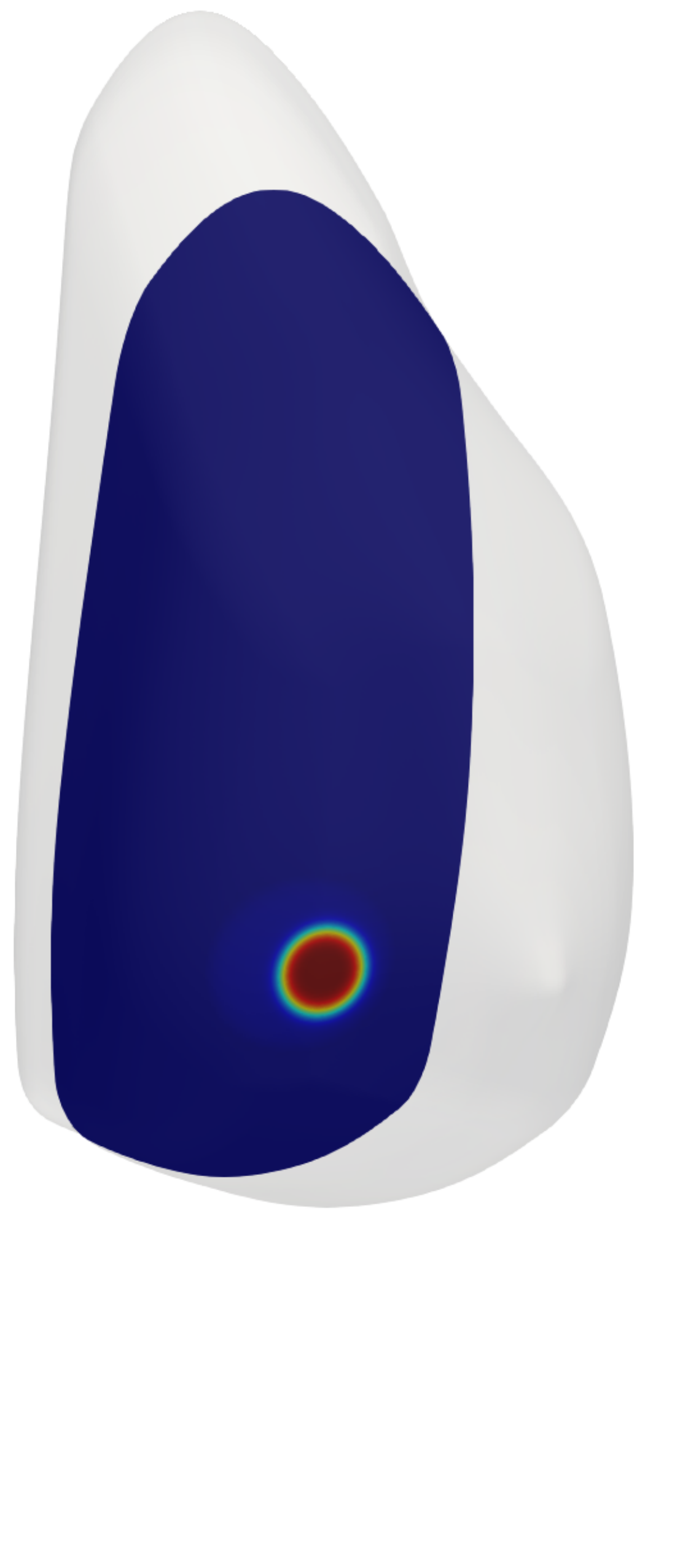}  
\end{subfigure}
\begin{subfigure}{0.328\columnwidth} \centering
\includegraphics[trim={0cm 8cm 0cm 0cm},clip,width=0.55\columnwidth]{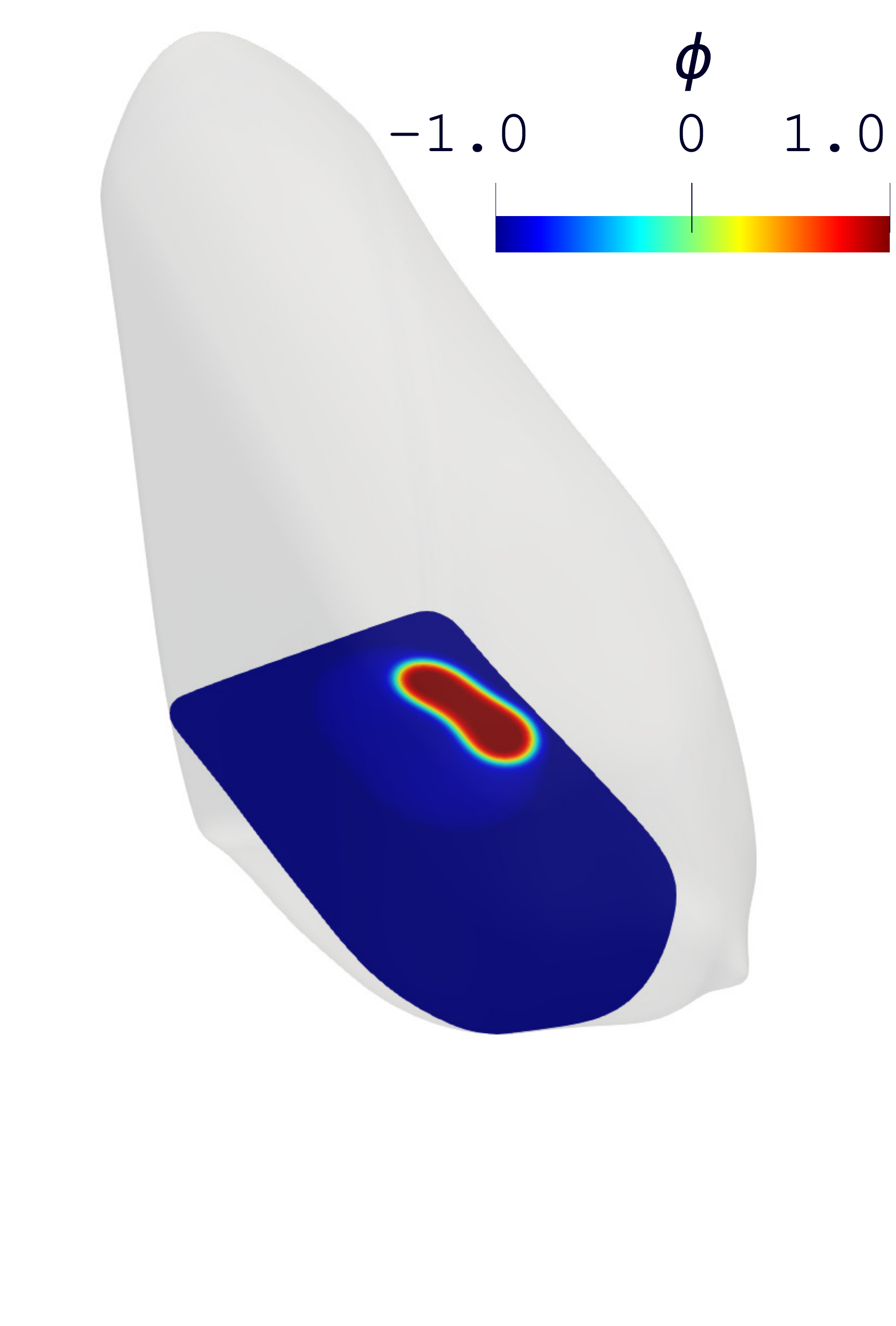}  
\includegraphics[trim={0cm 8cm 0cm 0cm},clip,width=0.35\columnwidth]{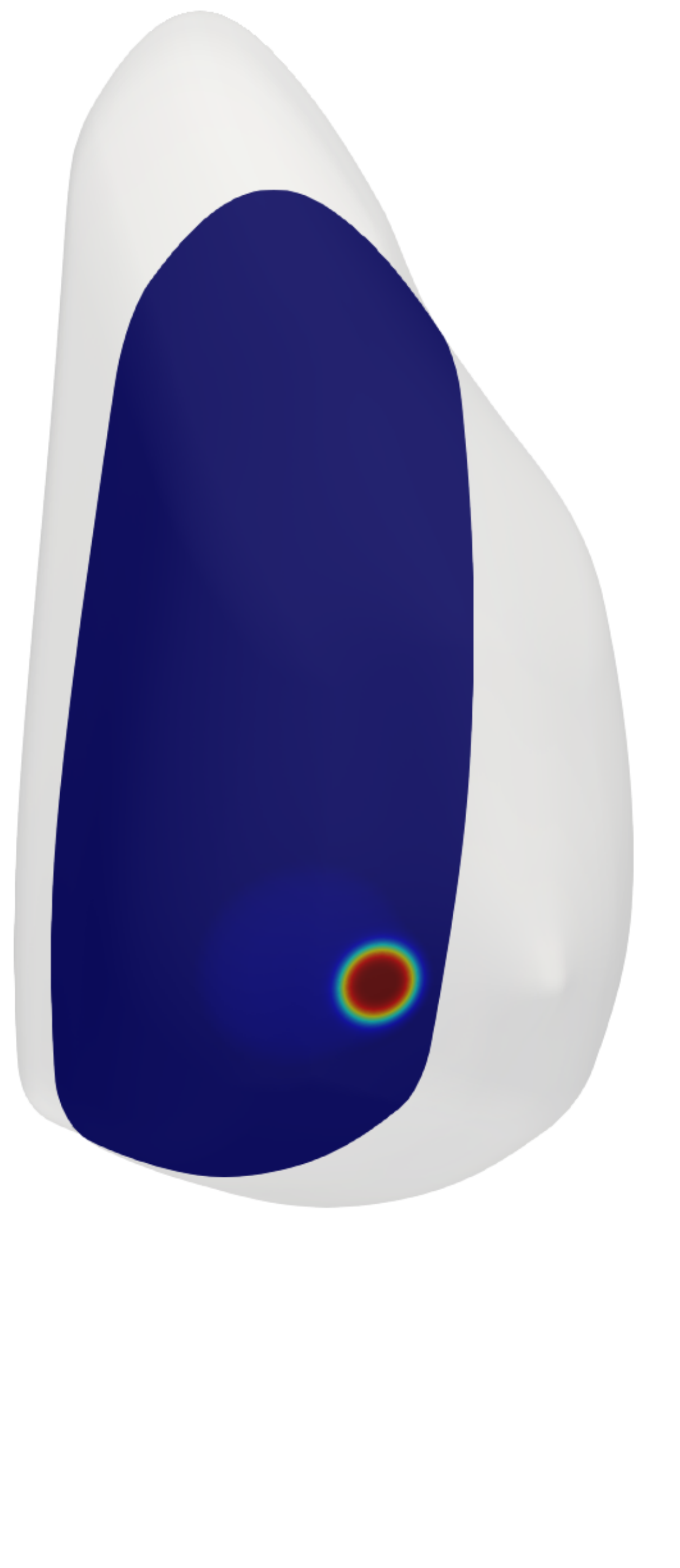}   \end{subfigure}
\caption{Phase-field}
\end{subfigure}
\\
\begin{subfigure}{\columnwidth}
\begin{subfigure}{0.328\columnwidth} \centering
\includegraphics[trim={0cm 8cm 0cm 0cm},clip,width=0.55\columnwidth]{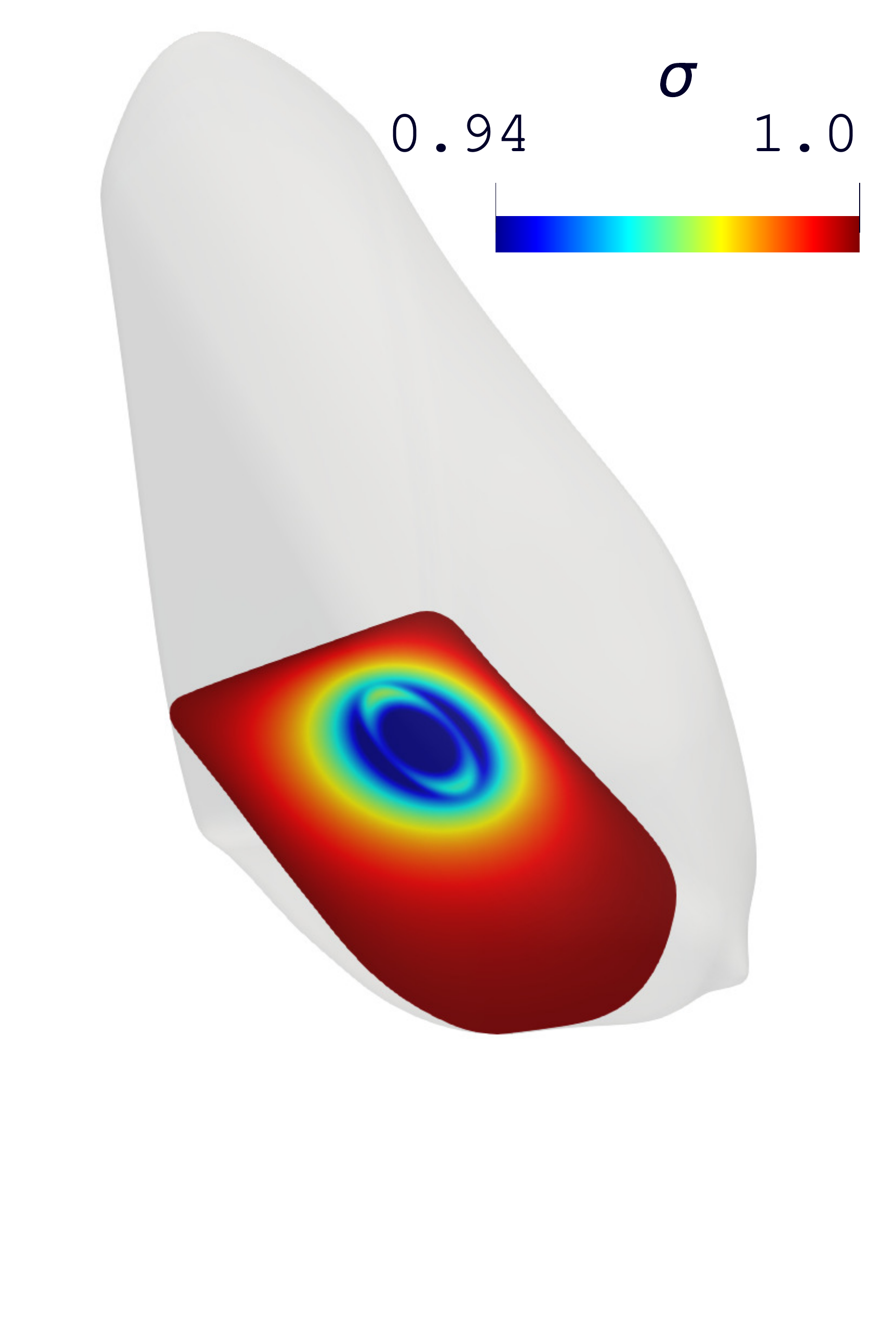}  
\includegraphics[trim={0cm 8cm 0cm 0cm},clip,width=0.35\columnwidth]{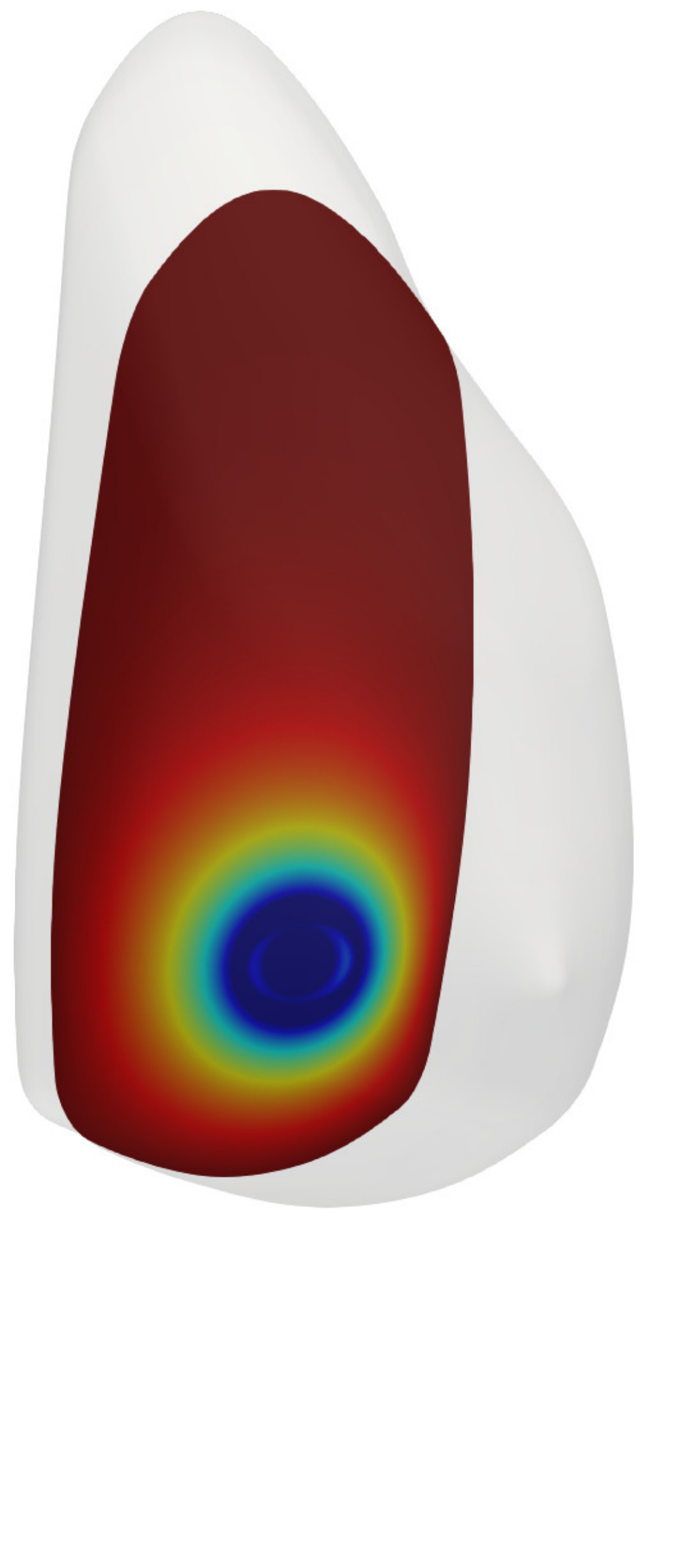}  
 \end{subfigure}
\begin{subfigure}{0.328\columnwidth} \centering
\includegraphics[trim={0cm 8cm 0cm 0cm},clip,width=0.55\columnwidth]{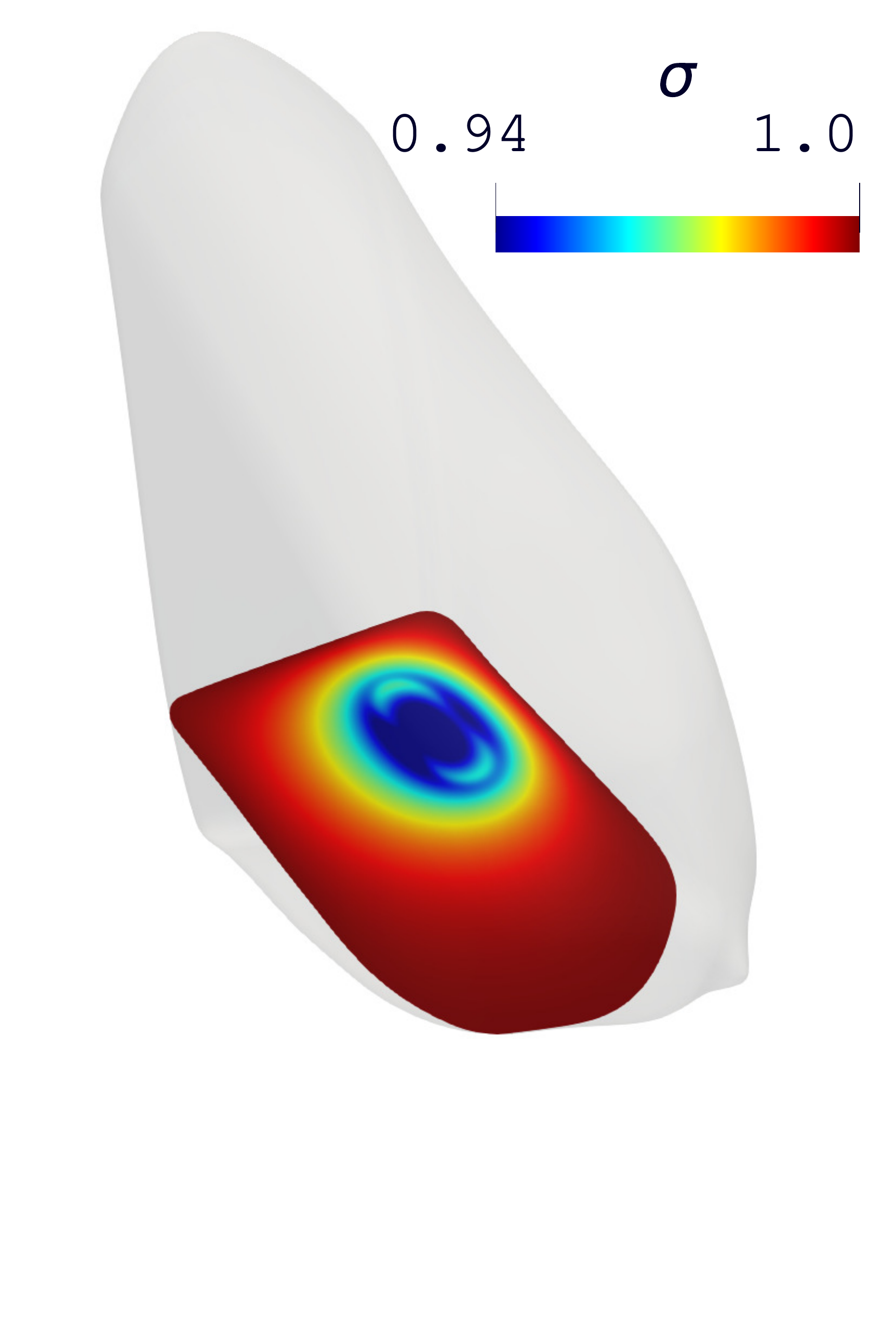}  
\includegraphics[trim={0cm 8cm 0cm 0cm},clip,width=0.35\columnwidth]{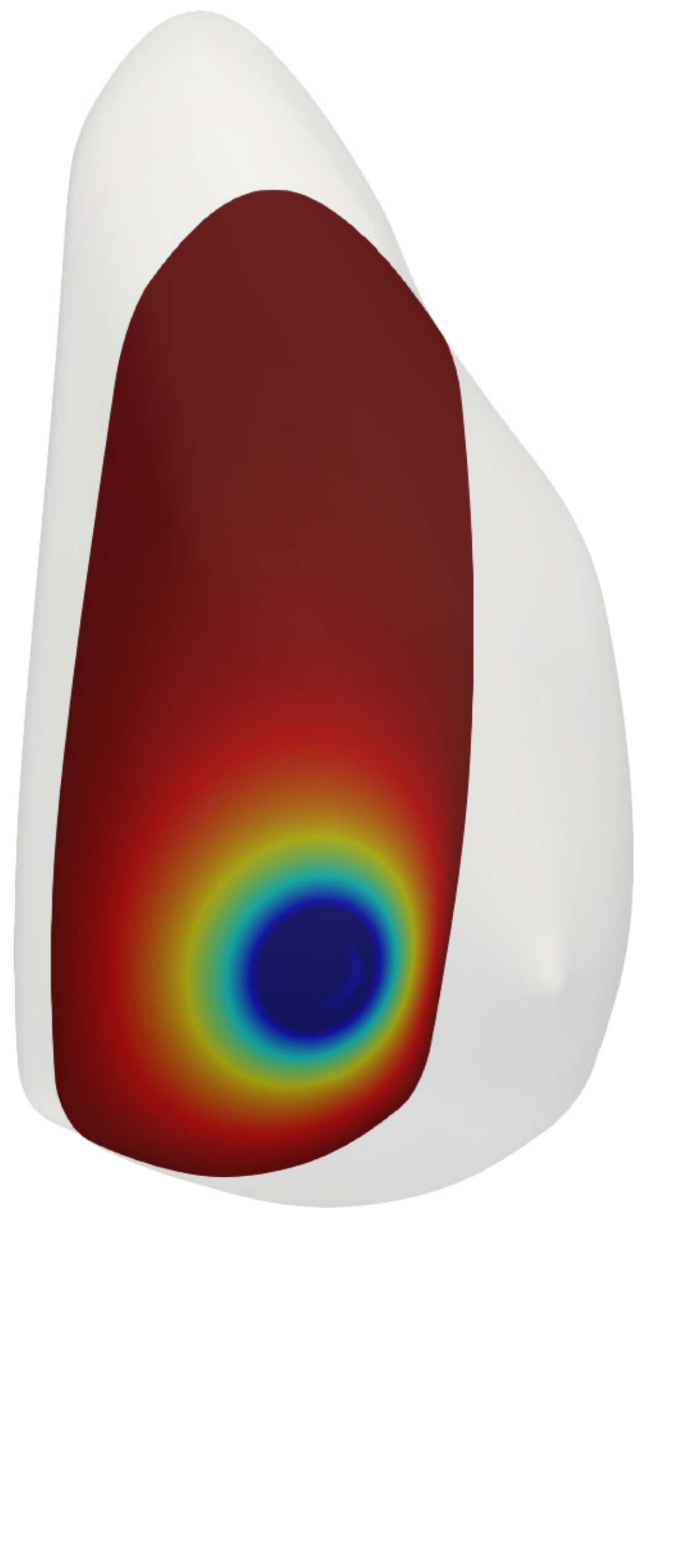}  
\end{subfigure}
\begin{subfigure}{0.328\columnwidth} \centering
\includegraphics[trim={0cm 8cm 0cm 0cm},clip,width=0.55\columnwidth]{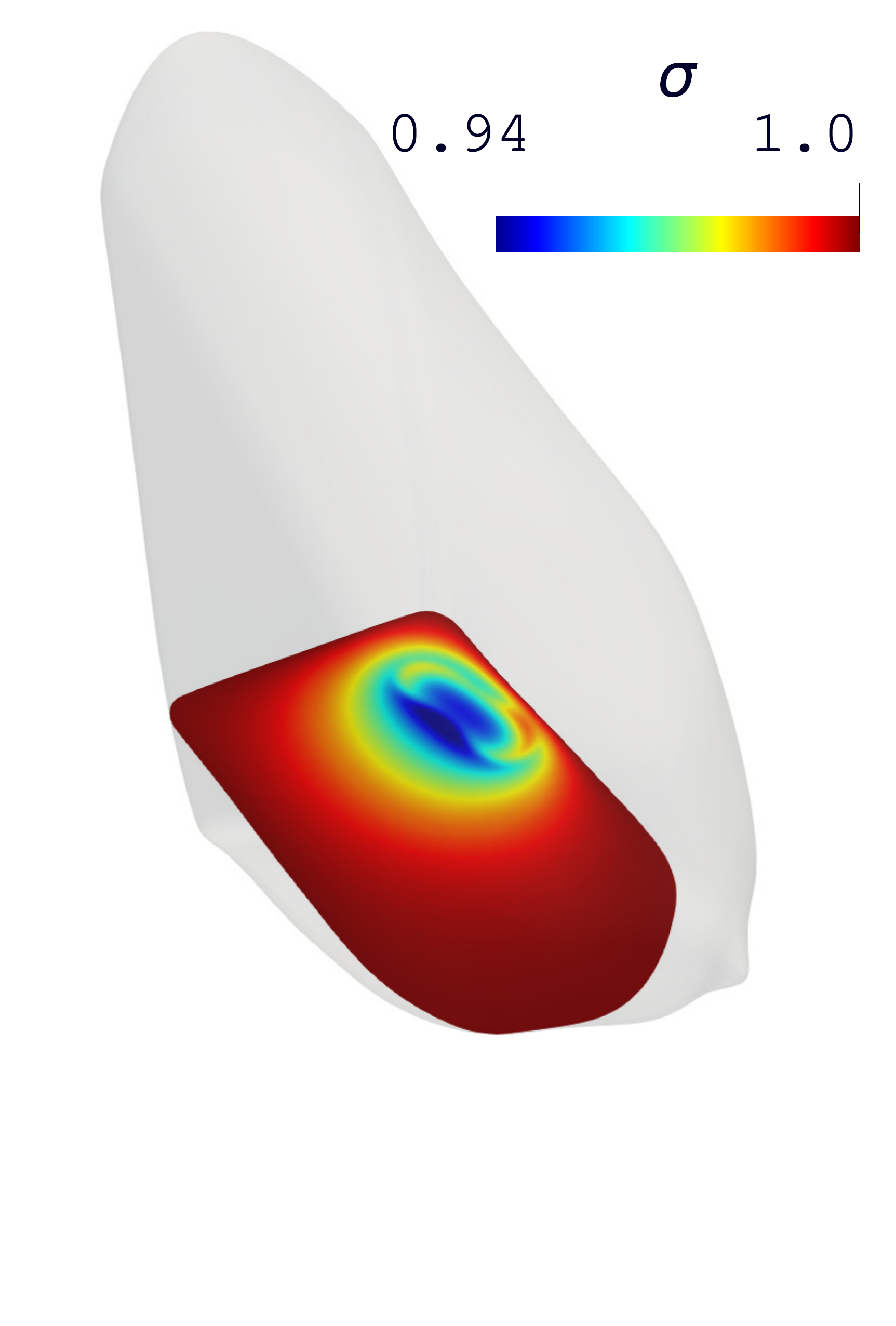}  
\includegraphics[trim={0cm 8cm 0cm 0cm},clip,width=0.35\columnwidth]{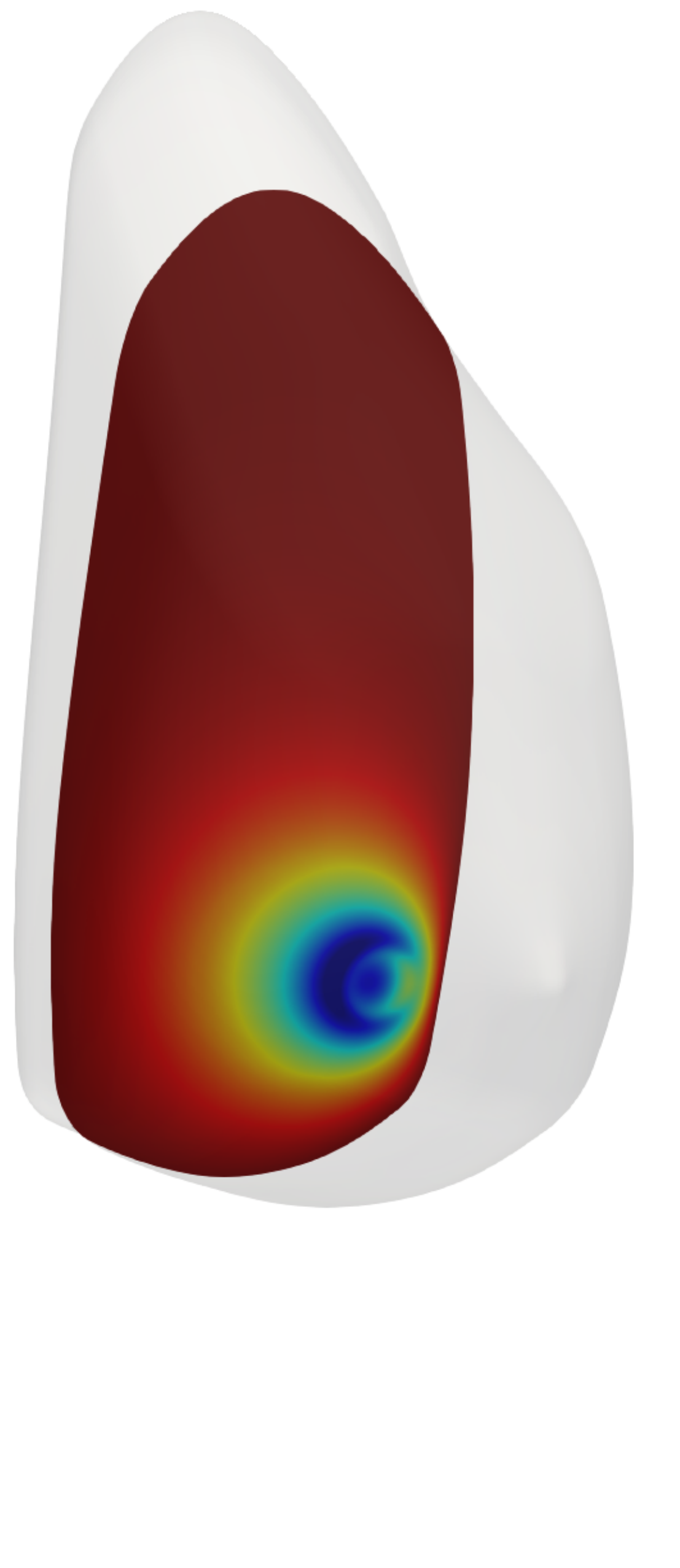}   \end{subfigure}
\caption{Nutrient concentration}
\end{subfigure}
\\
\begin{subfigure}{\columnwidth}
\begin{subfigure}{0.328\columnwidth} \centering
\includegraphics[trim={0cm 0cm 0cm 0cm},clip,width=0.42\columnwidth]{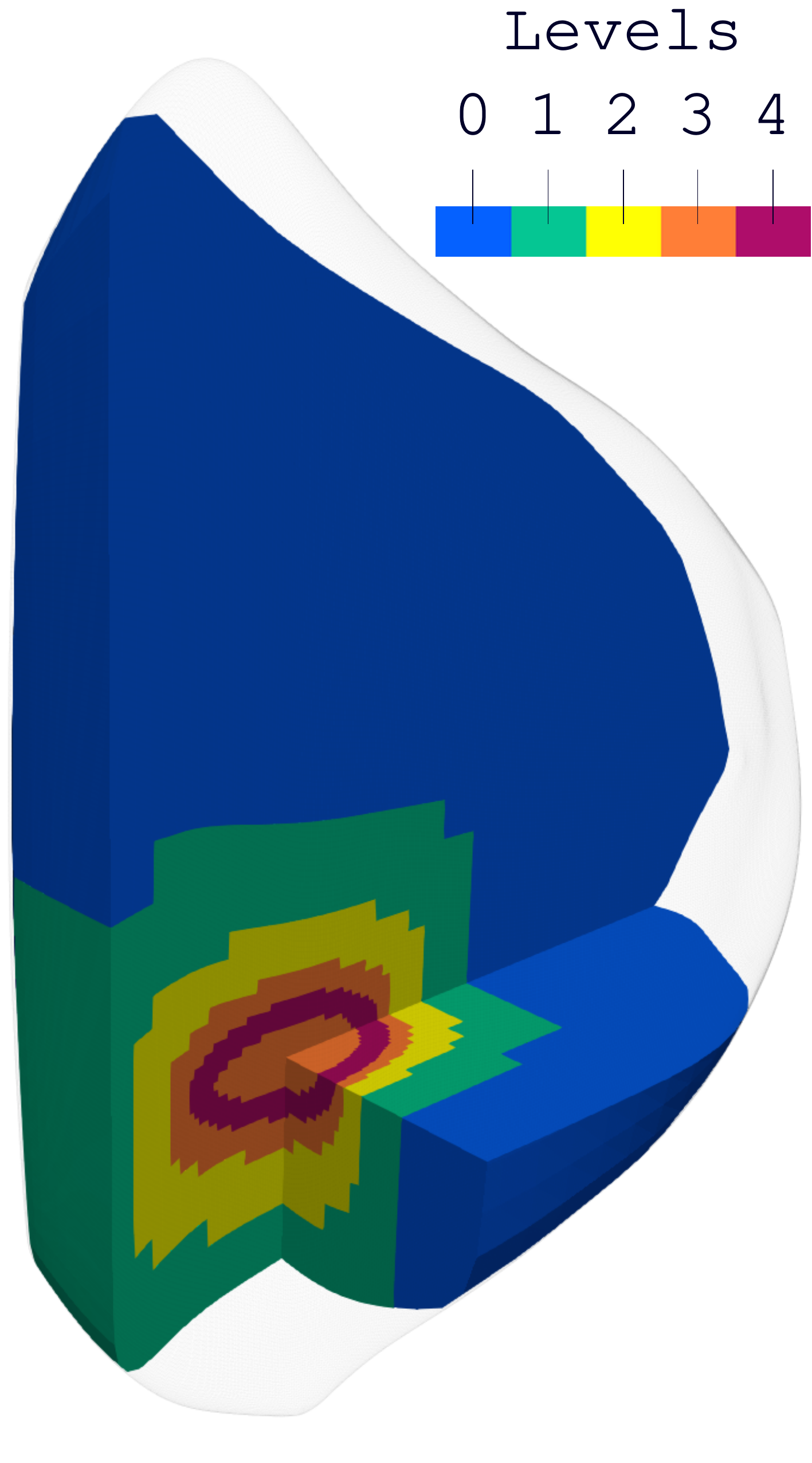}  
\includegraphics[trim={0cm 0cm 0cm 0cm},clip,width=0.42\columnwidth]{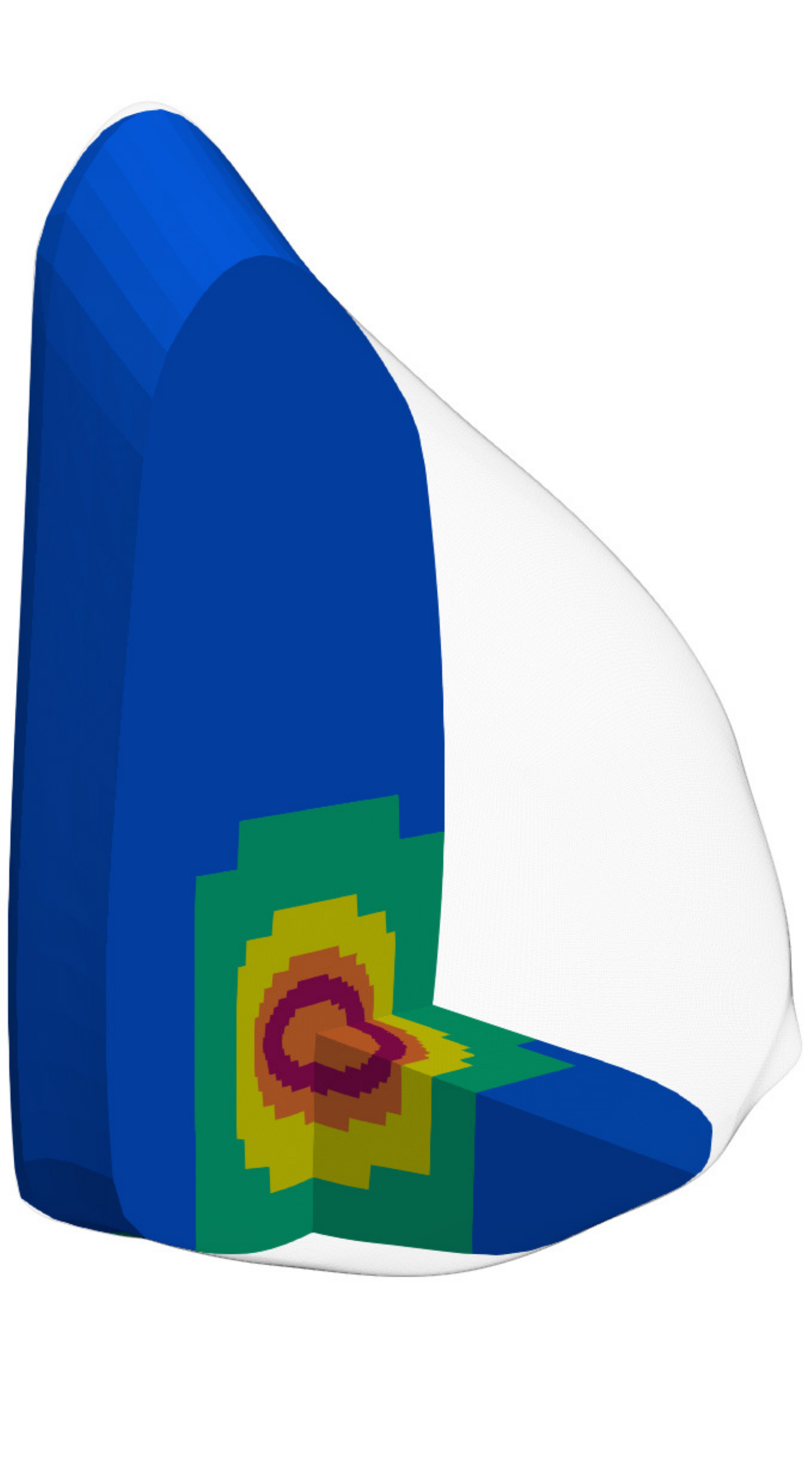}  
 \end{subfigure}
\begin{subfigure}{0.328\columnwidth} \centering
\includegraphics[trim={0cm 0cm 0cm 0cm},clip,width=0.42\columnwidth]{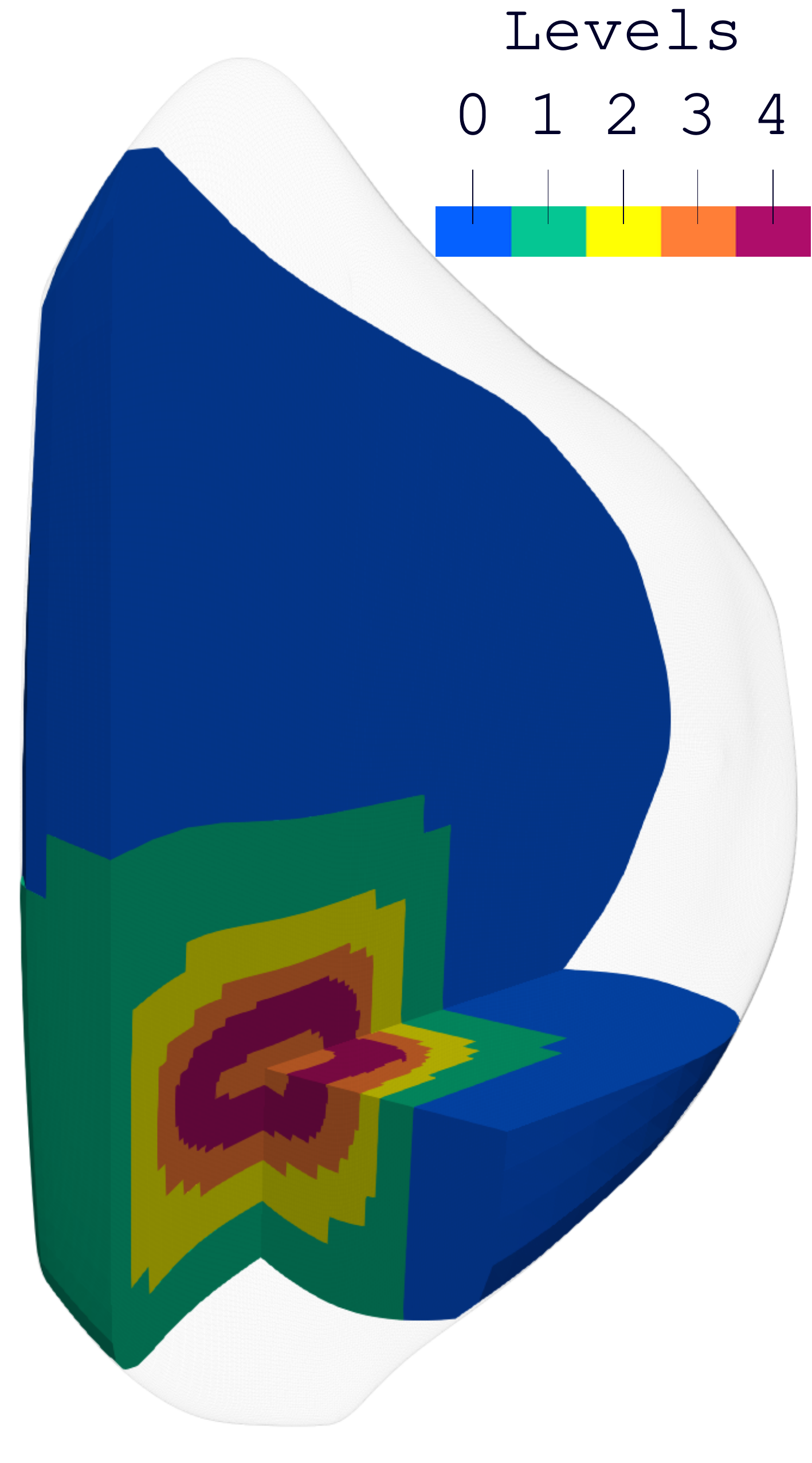}  
\includegraphics[trim={0cm 0cm 0cm 0cm},clip,width=0.42\columnwidth]{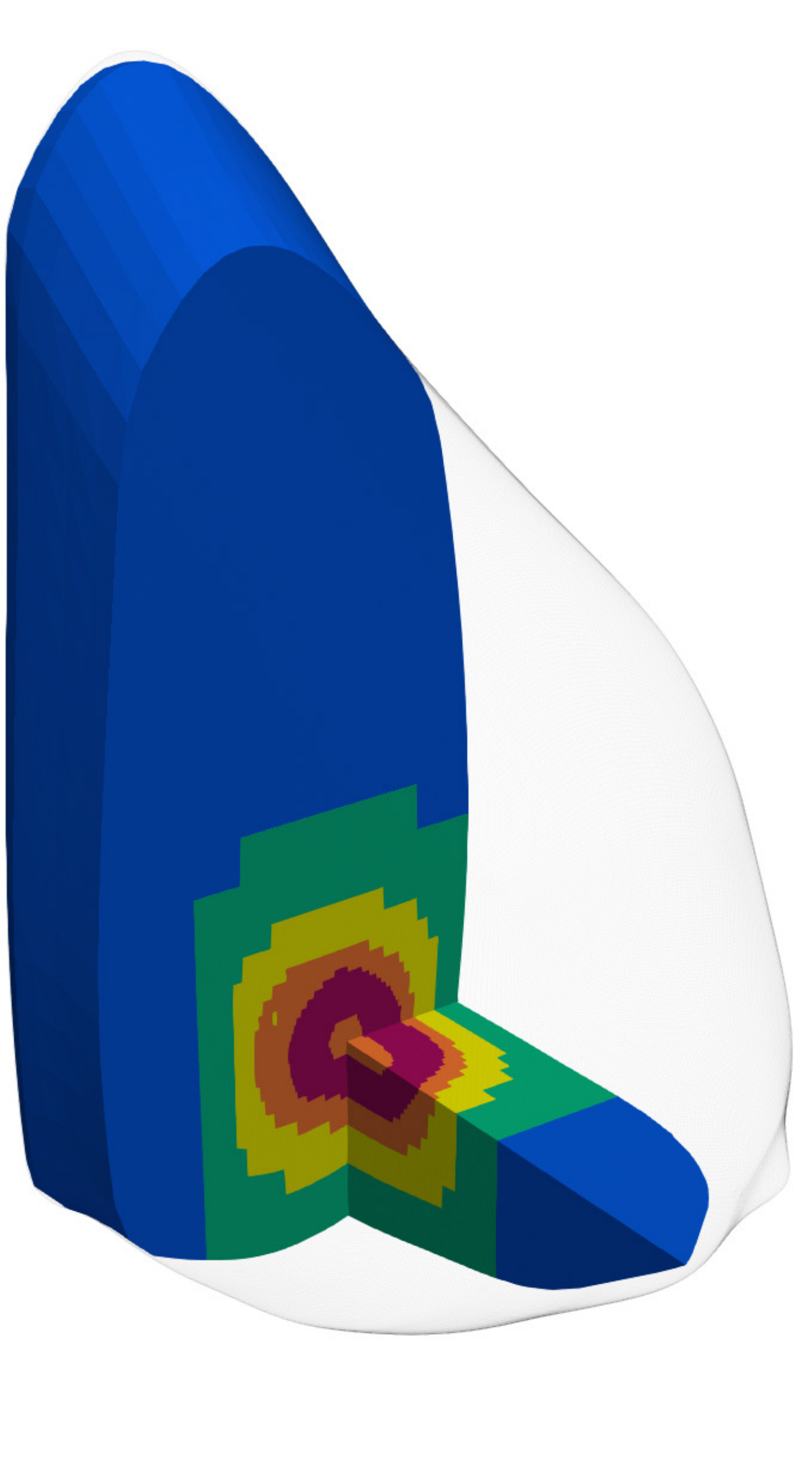}  
\end{subfigure}
\begin{subfigure}{0.328\columnwidth} \centering
\includegraphics[trim={0cm 0cm 0cm 0cm},clip,width=0.42\columnwidth]{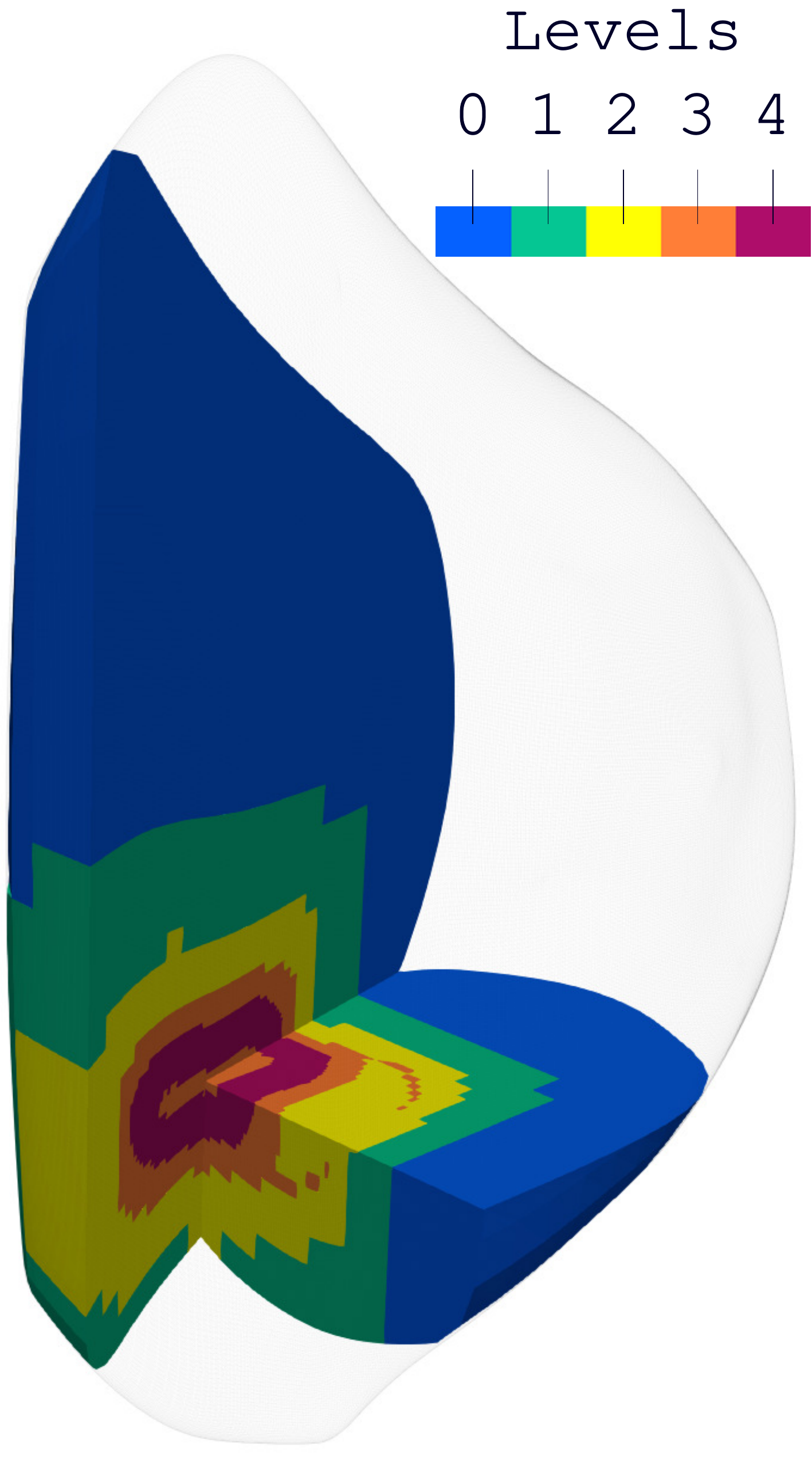}  
\includegraphics[trim={0cm 0cm 0cm 0cm},clip,width=0.42\columnwidth]{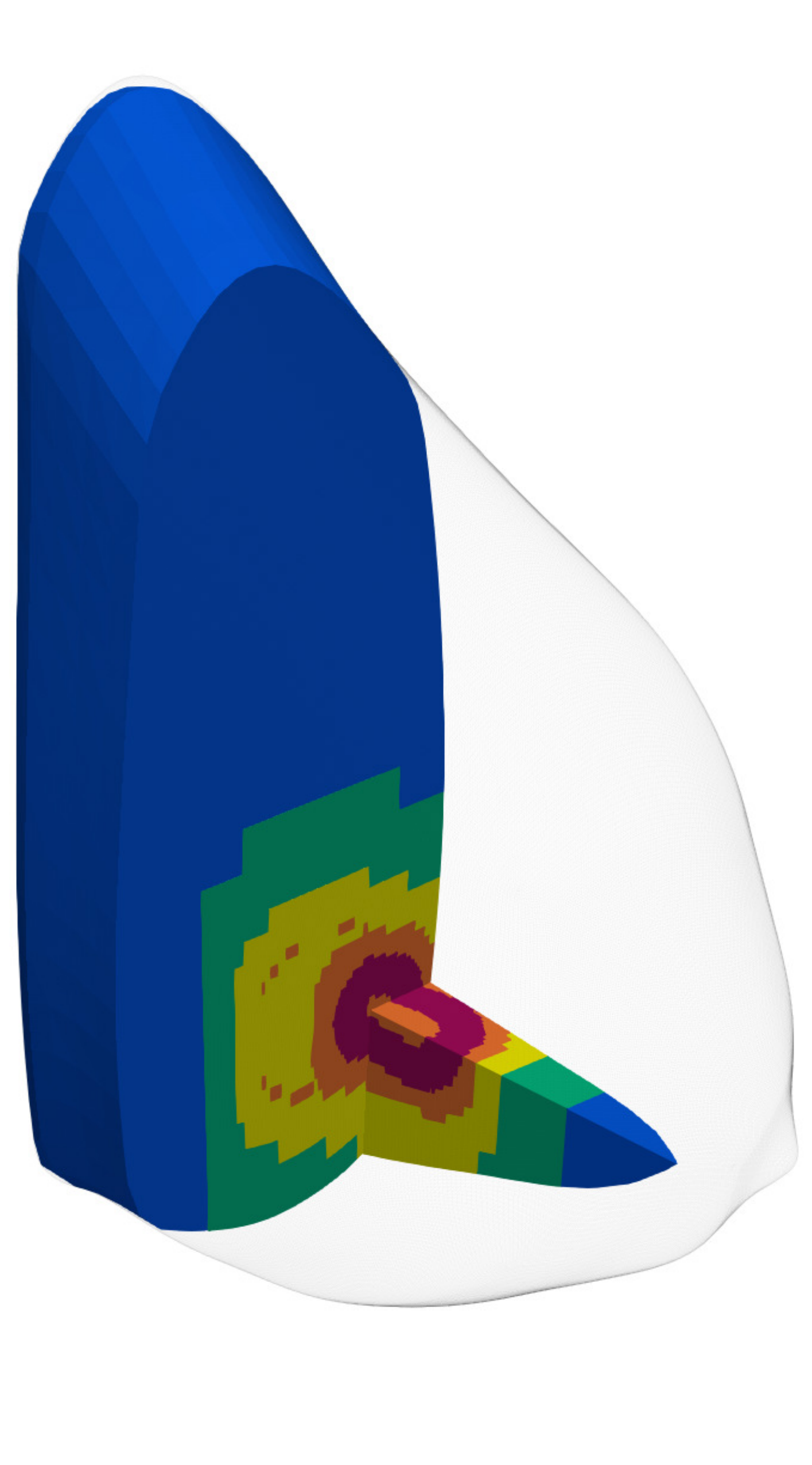}  
\end{subfigure}
\caption{Adaptive THB-spline mesh}
\end{subfigure}
\caption{Spatiotemporal evolution of an ellipsoidal tumor in a patient-specific breast geometry at times $t=$ 0, 0.135, and 0.27 (from left to right). THB-spline mesh configuration: $p=3$, $\ell = 4$, $m=2$, $\alpha = 0.01$ and $\beta = 0.0001$, with deepest refinement level at $2^8\times2^8\times2^8$.}
\label{fig:breast1_model}
\end{figure}

A representative anatomical length scale of $L=100$ mm is used to non-dimensionalize the breast model. 
With this scaling, the computational domain ($\Omega$) has approximate maximum and minimum base diameters of 1.67 and 0.63, respectively, and a maximum height of 1.18. 
Owing to the reduced size of the resulting computational domain $\Omega$, the deepest refinement level is selected to maintain a mesh resolution consistent with that adopted in the previous benchmark studies in Sections~\ref{TG2D} and \ref{TG3D}. 
Accordingly, the deepest level of the hierarchy consists of a mesh with $2^8 \times 2^8 \times 2^8$ elements. 
This provides a conservative characteristic element size of approximately $1.67/256 = 6.52 \cdot 10^{-3}$, resulting in a similar level of discretization across the phase-field transition region as in the 2D and 3D cases in Sections~\ref{TG2D} and \ref{TG3D}.
The coarsest refinement level is selected to consist of 16 elements in each direction. 
Therefore, the results presented here employ a mesh configuration with $p=3$, $\ell = 4$, $m=2$, $\alpha = 0.01$, and $\beta = 0.0001$.
To accommodate the organ-scale scenario within our modeling framework, we select $\lambda = 0.0001$, $D=0.05$, $M=2$, and $\mathcal{P} = 0.5$, while all other parameters remain unchanged from Section~\ref{TG2D}.
The initial tumor field is initialized using a hyperbolic tangent profile of an ellipsoidal shape, with semi-axes 0.08, 0.065, and 0.145, and centered at the patient-specific tumor location.

The spatiotemporal evolution of the tumor within the patient-specific breast geometry is presented in Fig.~\ref{fig:breast1_model}, which shows the 3D tumor morphology together with axial cuts of the breast providing the $\phi$ and $\sigma$ fields over the time interval $[0, 0.27]$. 
It also shows the adaptive THB-spline mesh as it dynamically tracks the evolving interface between the healthy and tumor tissues.
The tumor has an elongated shape close to the boundary, which exposes a large portion of the tumor interface to the nutrient gradient in the direction of higher nutrient concentrations.
This leads to migration-dominated evolution towards the boundaries with limited morphological changes.
While more realistic tumor dynamics could be achieved by calibrating the model parameters to longitudinal MRI data \cite{Wu2022}, the qualitative results shown in Fig.~\ref{fig:breast1_model} establish that the proposed computational pipeline can be successfully applied to patient-specific breast geometries to govern the coupled tumor–nutrient dynamics in clinically relevant scenarios. 
For completeness, a second representative example is presented in Fig.~\ref{fig:breast2_model} in the Appendix, where the tumor is initialized near the center of the domain. Using the same mesh configuration and model parameters, this example highlights the resulting morphological evolution under reduced boundary influence.

\section{Conclusion} \label{Conclusion}

In this work, we investigated the use of THB spline isogeometric methods to solve a CH, nutrient-driven phase-field model of tumor growth in its primal form.
The hierarchical mesh is dynamically refined along the evolving diffuse interface between healthy and tumor tissues, where solution gradients are steepest, and coarsened in regions of low activity.
A systematic series of numerical experiments was first conducted to establish the appropriate mesh parameters (i.e., polynomial degree, the maximum admissible element size, and the admissibility class) required to accurately resolve the interfacial dynamics of the CH system. 
Benchmark problems in both 2D and 3D were solved to capture complex tumor morphologies and investigate their sensitivity to variations in model parameters (e.g., proliferation rate, mobility, and interfacial thickness parameter) as well as in the choice of the initial condition of the tumor phase field.
The results confirm that the proposed framework is capable of reproducing a broad range of clinically relevant tumor morphologies highlighting the versatility and predictive scope of the model.
Finally, the model was applied to a clinically motivated 3D untreated BCa scenario, employing a patient-specific geometry reconstructed from anatomic MRI data. 
The results demonstrate that the adaptive THB-spline framework is applicable to organ-scale tumor growth simulations, which could facilitate the predictive modeling of BCa progression in patient-specific settings using CH-based tumor growth models to capture intricate cancer geometries. 

A compelling direction for future work is the calibration of model parameters against longitudinal MRI data, which would enable the transition from qualitative morphological prediction to quantitative, patient-specific tumor growth modeling with direct clinical relevance. 
Further extensions may include the incorporation of mechanical feedback from the surrounding breast tissue \cite{lorenzo2024global, Cristini2009, Garcke2021} and the modeling of treatment response such as the effect of chemotherapy or radiotherapy on tumor proliferation and nutrient consumption \cite{lima2017, colli2020, lorenzo2024global}. 
Together, these developments aimed at advancing the predictive fidelity of the framework to clinical applicability.

\section*{Declaration of competing interest} 

The authors declare that they have no known competing financial interests or personal relationships that could have appeared to influence the work reported in this paper.

\section*{Acknowledgments}

The authors gratefully acknowledge the support from the European Union’s Horizon Europe (HORIZON) under the Marie Skłodowska-Curie grant agreement No 101208599 Call: HORIZON-MSCA-2024-PF-01. 
GL also acknowledges grant PID2023-146347OA-I00 funded by MICIU/AEI/10.13039/501100011033 and ERDF/EU, as well as grant RYC2022-036010-I funded by MICIU/AEI/10.13039/501100011033 and ESF+. CG and AR are members of INdAM-GNCS, whose support is gratefully acknowledged.

\clearpage

\appendix 

\section{Implementation workflow and pseudocode} \label{Appendix-flowchart}
The key implementation details of the CH-based tumor growth model considered herein using adaptive THB-spline framework are systematically outlined in the pseudocode presented in Algorithm~\ref{seudocode_timemarchNR}. 
It traces the computational workflow, from nonlinear solution via the NR method and time integration via the generalised-$\alpha$ scheme to the adaptive mesh refinement and coarsening cycles.

\setcounter{algocf}{0}
\renewcommand{\thealgocf}{A\arabic{algocf}}

\begin{algorithm}[!ht]
\DontPrintSemicolon
\caption{Pseudocode to solve phase-field-nutrient system.}\label{seudocode_timemarchNR}
\textbf{Initialize:} $\tilde{\boldsymbol{\phi}}_{0}$, $\dot{\tilde{\boldsymbol{\phi}}}_0 = \boldsymbol{0}$, $\tilde{\boldsymbol{\sigma}}_0 $, $\Delta t$, $\left[ 0, T \right]$ \tcp*[r]{$\tilde{\boldsymbol{\sigma}}_0 \gets \bR_\sigma(\tilde{\boldsymbol{\sigma}}_{0}, \tilde{\boldsymbol{\phi}}_{0}) = \bf{0} $ (Eq.~\eqref{R_sigma})}
\While{$t < T$}{
    $t_{n+1} \gets t_n + \Delta t$\;
    $\tilde{\boldsymbol{\phi}}_{n+1}^0 = \tilde{\boldsymbol{\phi}}_{n}$, \quad $\dot{\tilde{\boldsymbol{\phi}}}_{n+1}^0 = \dfrac{\gamma - 1}{\gamma}\dot{\tilde{\boldsymbol{\phi}}}_{n}$  \tcp*[r]{Predictor (Eqs.~\eqref{predictor1} and \eqref{predictor2})}
    
\For(\tcp*[f]{$i$: iteration counter, $i_{\max}$: maximum iterations}){$i$ = 0, 1, $\hdots$, $i_{max}$ }{ 
$\tilde{\boldsymbol{\phi}}_{n+\alpha_f}^{(i+1)}  = \tilde{\boldsymbol{\phi}}_{n} + \alpha_f \left( \tilde{\boldsymbol{\phi}}_{n+1}^{(i)} - \tilde{\boldsymbol{\phi}}_{n} \right)$ \tcp*[r]{ Eqs.~\eqref{phi_nalpha} and \eqref{phi_nplus1}}
$\dot{\tilde{\boldsymbol{\phi}}}_{n+\alpha_m}^{(i+1)}  = \dot{\tilde{\boldsymbol{\phi}}}_n + \alpha_m \left( \dot{\tilde{\boldsymbol{\phi}}}_{n+1}^{(i)} - \dot{\tilde{\boldsymbol{\phi}}}_n \right)$ \tcp*[r]{ Eq.~\eqref{phidot_nalpha}}

\If(\tcp*[f]{$\varepsilon_{\text{tol}}$ = NR tolerance for convergence}){$\|\bR_\phi(\tilde{\boldsymbol{\phi}}_{n+\alpha_f}^{(i+1)}, 
\dot{\tilde{\boldsymbol{\phi}}}_{n+\alpha_m}^{(i+1)}, 
\tilde{\boldsymbol{\sigma}}_n)\| < \varepsilon_{\text{tol}}$}{
    \textbf{break}\;}
\Else{
    Solve: $\bK_{\text{NR}}^{(i+1)} \Delta\dot{\tilde{\boldsymbol{\phi}}}_{n+1}^{(i+1)} 
    = - \bR_\phi(\tilde{\boldsymbol{\phi}}_{n+\alpha_f}^{(i+1)}, 
    \dot{\tilde{\boldsymbol{\phi}}}_{n+\alpha_m}^{(i+1)}, 
    \tilde{\boldsymbol{\sigma}}_n)$ \tcp*[r]{ Eqs.~\eqref{R_phi}, \eqref{knablaphi_R} and \eqref{tangentMatK}} 
    Update: 
    $\dot{\tilde{\boldsymbol{\phi}}}_{n+1}^{(i+1)} = \dot{\tilde{\boldsymbol{\phi}}}_{n+1}^{(i)} + \Delta\dot{\tilde{\boldsymbol{\phi}}}_{n+1}^{(i+1)},$ \quad
    $\tilde{\boldsymbol{\phi}}_{n+1}^{(i+1)} = \tilde{\boldsymbol{\phi}}_{n+1}^{(i)} + \gamma \Delta t \Delta\dot{\tilde{\boldsymbol{\phi}}}_{n+1}^{(i+1)}$ \tcp*[r]{ Eqs.~\eqref{phidot_nplus1_i} and \eqref{phi_nplus1_i}} 
}
}
Solve: $\tilde{\boldsymbol{\sigma}}_{n+1} \gets \bR_\sigma(\tilde{\boldsymbol{\sigma}}_{n+1}, \tilde{\boldsymbol{\phi}}_{n+1}) = \bf{0} $ \tcp*[r]{Eq.~\eqref{R_sigma}}
\If(\tcp*[f]{$n_{\mathrm{t}}$ = time-step counter}){$n_{\mathrm{t}} \bmod n_{\mathrm{r}} = 0$ }{
    Evaluate: $\varepsilon_{\phi}^e \gets  \frac{ 1 }{\text{Vol} (\Omega_e)} \int_{\Omega_e} \| \nabla\phi_{n+1} \|^2 \diff \Omega$, $\: e = 1 \text{ to } n_{\mathrm{el}}$  \tcp*[r]{Eq.~\eqref{error_estimator}}
    
    \If{$\varepsilon_{\phi}^e > \alpha \cdot \varepsilon_{\max}$}{
    Mark element $e$ for refinement: $Q_r$ \\
    Ensure admissibility class $m$: $Q_r \gets Q_r \cup {\cal N}_r({\cal Q},{Q}_r,m)$ \\
    }
    Refine all marked elements: $\mathcal{Q}^\ell \to \mathcal{Q}^{\ell+1}$, $\: \mathcal{T}(\mathcal{Q}^\ell) \to \mathcal{T}(\mathcal{Q}^{\ell+1})$\; 
    Transfer solution: Find updated $\tilde{\boldsymbol{\phi}}_{n+1}^{*}$ such that: $\phi_{n+1}^{*} = \phi_{n+1}$, $\phi^{*}_{n+1}\in \mathcal{T}(\mathcal{Q}^{\ell+1})$\;
    Update: $\tilde{\boldsymbol{\phi}}_{n+1} \gets \tilde{\boldsymbol{\phi}}^{*}_{n+1}$ \\
}

\If{$n_{\mathrm{t}} \bmod n_{\mathrm{c}} = 0$}{
    Evaluate: $\varepsilon_{\phi}^e \gets  \frac{ 1 }{\text{Vol} (\Omega_e)} \int_{\Omega_e} \| \nabla\phi_{n+1} \|^2 \diff \Omega$, $\: e = 1 \text{ to } n_{\mathrm{el}}$  \tcp*[r]{Eq.~\eqref{error_estimator}}
    
    \If{$\varepsilon_{\phi}^e < \beta \cdot \varepsilon_{\max}$}{
    Mark element $e$ for coarsening: $Q_c$ \\
    Ensure admissibility class $m$: $Q_c \gets Q_c \cup {\cal N}_c({\cal Q},{Q}_c,m)$ \\
    }
    Coarsen all marked elements: $\mathcal{Q}^\ell \to \mathcal{Q}^{\ell-1}$, $\: \mathcal{T}(\mathcal{Q}^\ell) \to \mathcal{T}(\mathcal{Q}^{\ell-1})$ \\
    $L^2$ projection: Find updated $\tilde{\boldsymbol{\phi}}^{*}_{n+1}$ such that:  $\int_{\Omega} {\phi}^{*}_{n+1} \, \delta{\phi}^*_{n+1} \, \mathrm{d}\Omega = \int_{\Omega} {\phi_{n+1}} \, \delta{\phi}^*_{n+1} \, \mathrm{d}\Omega \quad \forall \, \phi^*_{n+1}, \delta{\phi}^*_{n+1} \in \mathcal{T}(\mathcal{Q}^{\ell-1})$\;
    Update: $\tilde{\boldsymbol{\phi}}_{n+1} \gets \tilde{\boldsymbol{\phi}}^{*}_{n+1}$ \\
}

}
\textbf{Result:} $\phi(\bx,t)$ and $\sigma(\bx,t)$ $\in \mathbb{R} $ $: (\bx,t) \in \bar{\Omega} \times \left[0, T \right]$
\end{algorithm}

\section{Extended numerical studies} \label{Complementary_Numerical_Studies}
\setcounter{figure}{0}
This section presents the morphological diversity achievable through selective variation of model parameters and initial conditions for the problem described in Section~\ref{TG2D}. 
The computational domain is the same square $\Omega = [-3,3]^2$. 
For the first variation, the model parameters are set to $E = 5$, $\lambda = 0.002$, and $M = 5$, while all remaining parameters are kept unchanged as stated in Section~\ref{TG2D}. 
Despite starting from the same initial condition, the tumor evolves into an elliptical morphology under this parameter choice (see Fig.~\ref{figure_2_appendix}), which shows the THB-spline mesh, the tumor evolution, and the corresponding nutrient concentration. 
The same mesh parameters are used, i.e. $p=3$, $\ell = 6$, $m=2$, $\alpha = 0.01$, $\beta = 0.0001$, with the finest refinement level corresponding to a $2^{10} \times 2^{10}$ ($h_e = 6/1024$) mesh resolution. 
We note that since $\lambda$ is increased from 0.0002 to 0.002, the diffuse interface is relatively broader, and the adopted mesh resolution is marginally conservative for this configuration.

\begin{figure}[!h]\centering
\begin{subfigure}{0.185\columnwidth} \centering
\includegraphics[trim={0cm 0cm 0cm 0cm},clip,width=1\columnwidth]{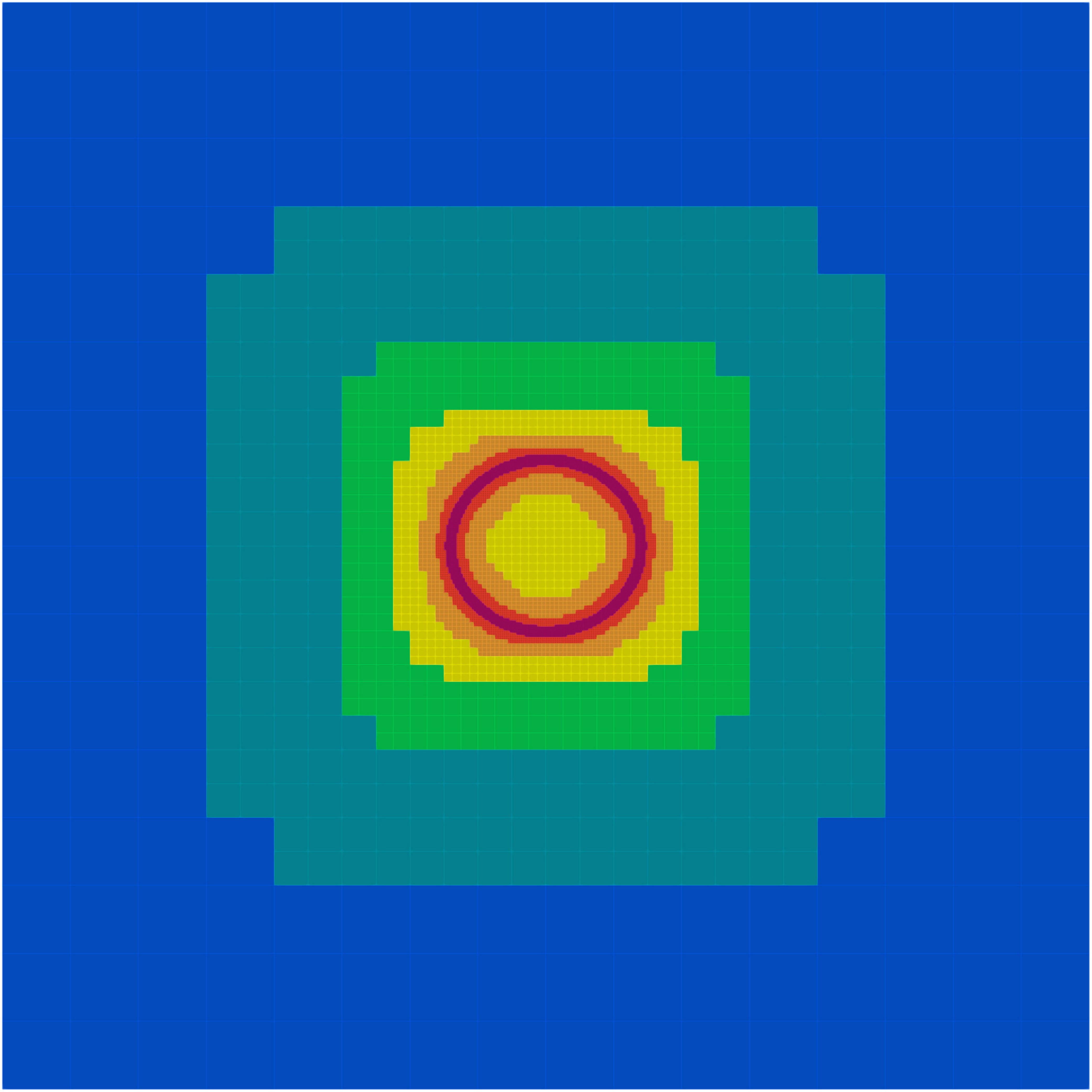}  \end{subfigure}
\begin{subfigure}{0.185\columnwidth} \centering
\includegraphics[trim={0cm 0cm 0cm 0cm},clip,width=1\columnwidth]{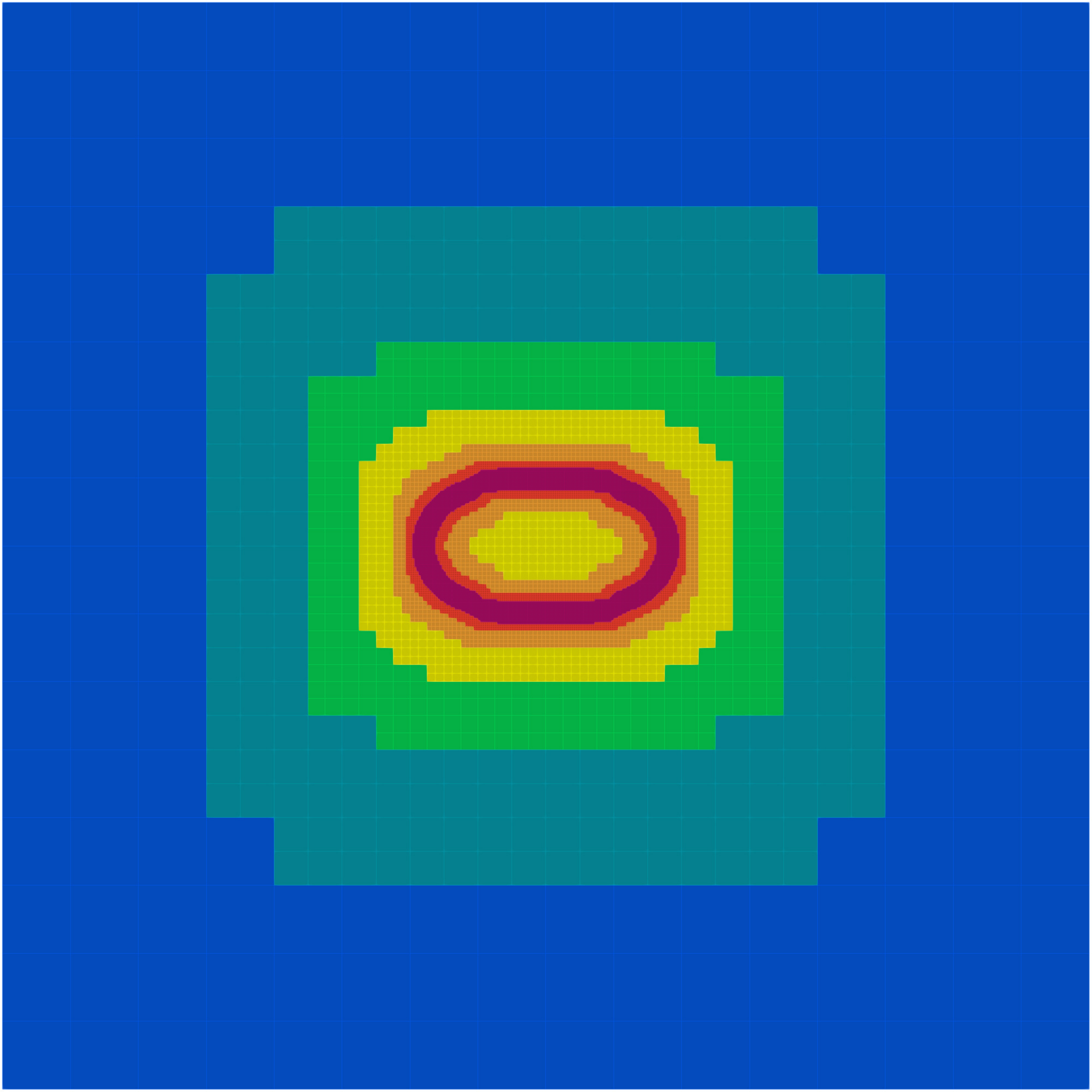}   \end{subfigure}
\begin{subfigure}{0.185\columnwidth} \centering
\includegraphics[trim={0cm 0cm 0cm 0cm},clip,width=1\columnwidth]{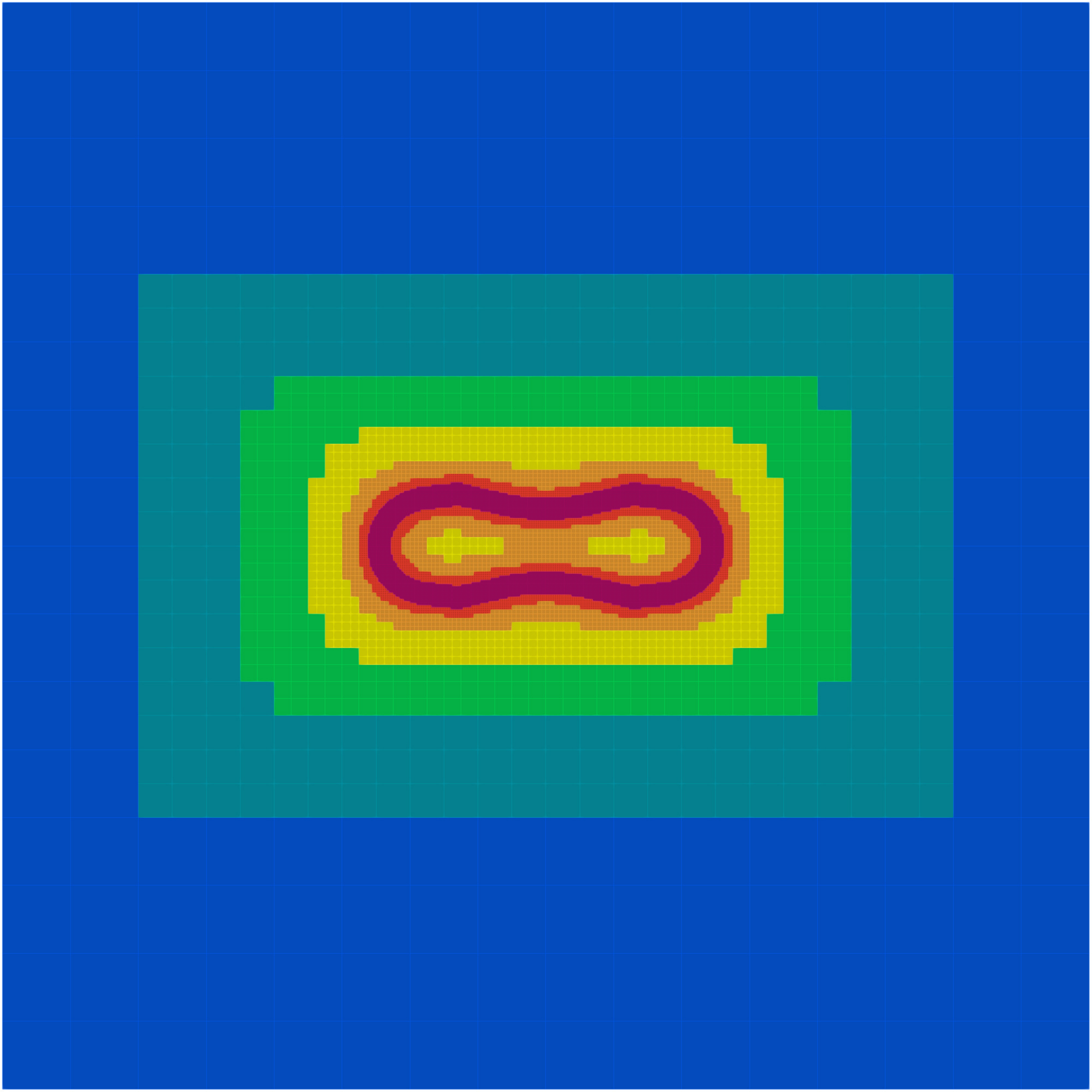} \end{subfigure}
\begin{subfigure}{0.185\columnwidth} \centering
\includegraphics[trim={0cm 0cm 0cm 0cm},clip,width=1\columnwidth]{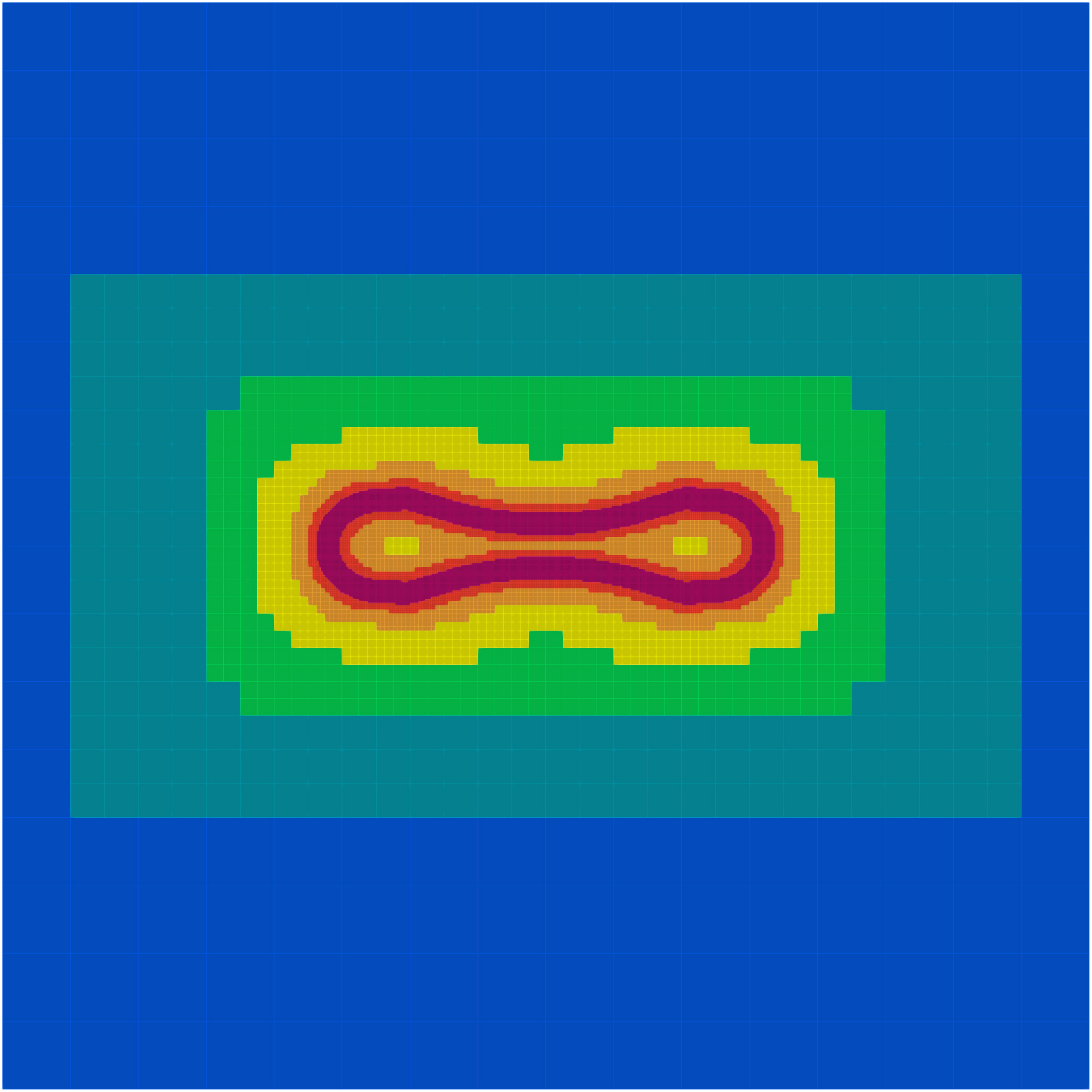} \end{subfigure}
\begin{subfigure}{0.185\columnwidth} \centering
\includegraphics[trim={0cm 0cm 0cm 0cm},clip,width=1\columnwidth]{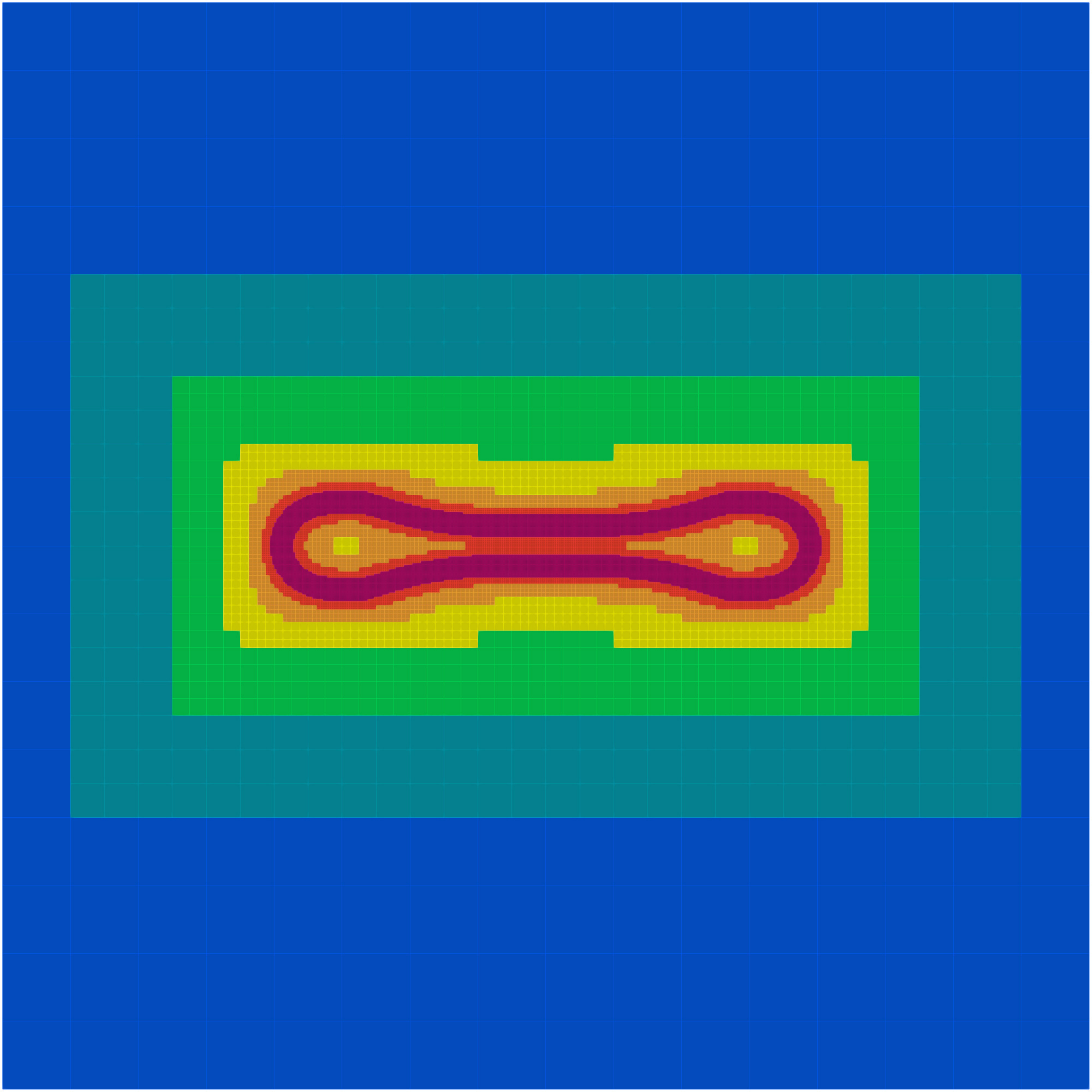}  \end{subfigure}
\begin{subfigure}{0.05\columnwidth} \centering
\includegraphics[trim={1cm 0cm 0cm 0cm},clip,width=0.84\columnwidth]{Images/Problem1_THBRes/ColorBar_Mesh_Lvl6-eps-converted-to.pdf} \end{subfigure}
\\ 
\begin{subfigure}{0.185\columnwidth} \centering
\includegraphics[trim={0cm 0cm 0cm 0cm},clip,width=1\columnwidth]{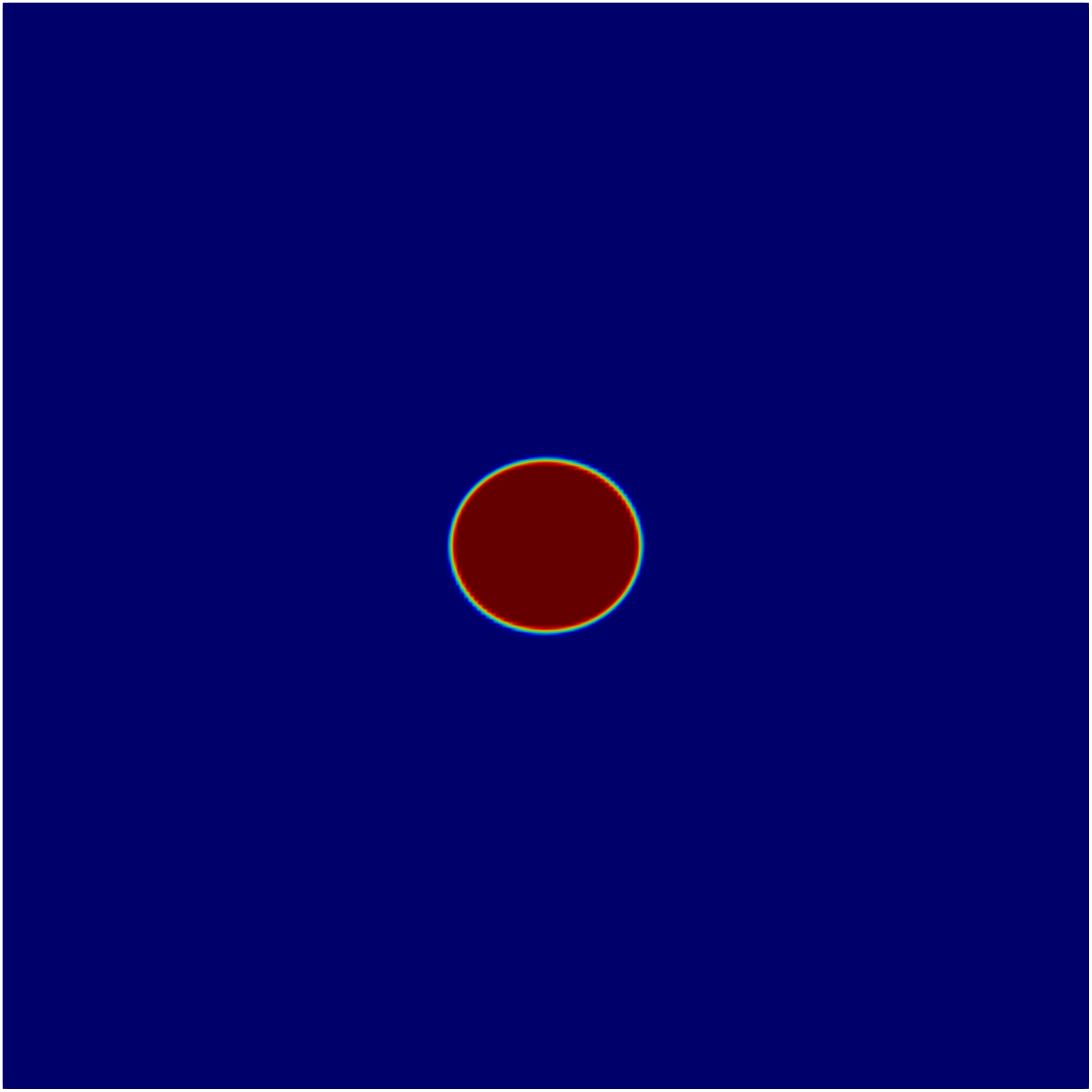}  \end{subfigure}
\begin{subfigure}{0.185\columnwidth} \centering
\includegraphics[trim={0cm 0cm 0cm 0cm},clip,width=1\columnwidth]{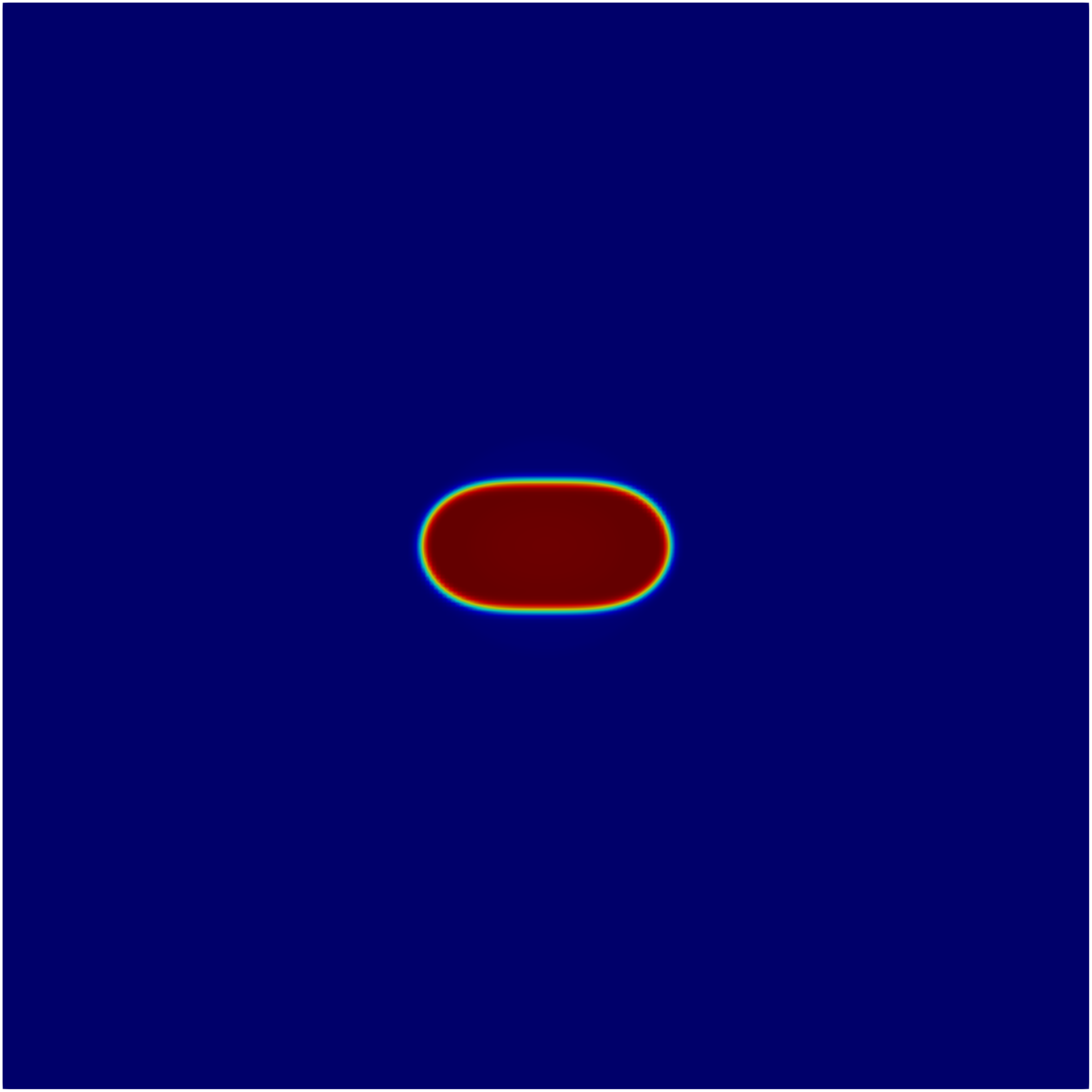}   \end{subfigure}
\begin{subfigure}{0.185\columnwidth} \centering
\includegraphics[trim={0cm 0cm 0cm 0cm},clip,width=1\columnwidth]{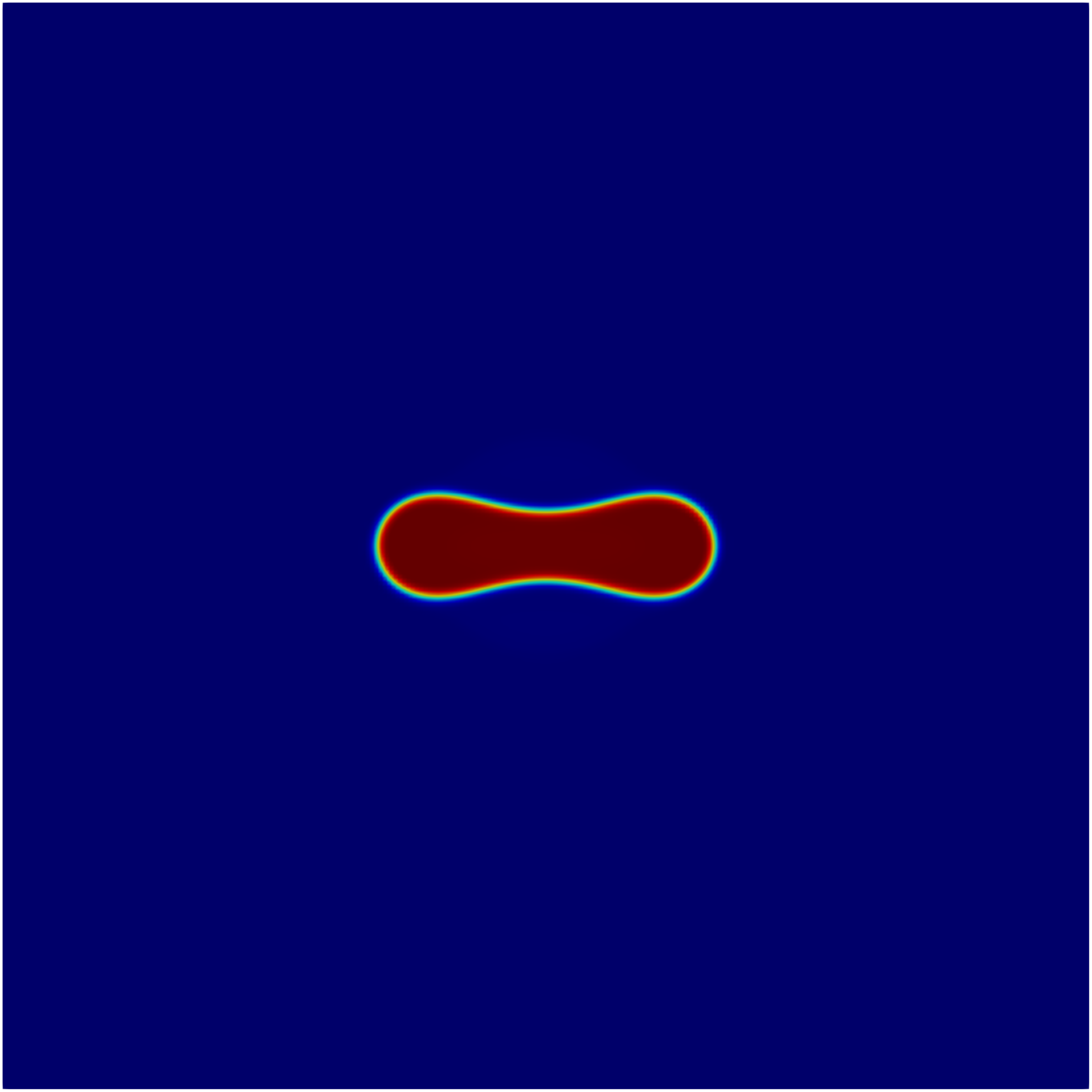} \end{subfigure}
\begin{subfigure}{0.185\columnwidth} \centering
\includegraphics[trim={0cm 0cm 0cm 0cm},clip,width=1\columnwidth]{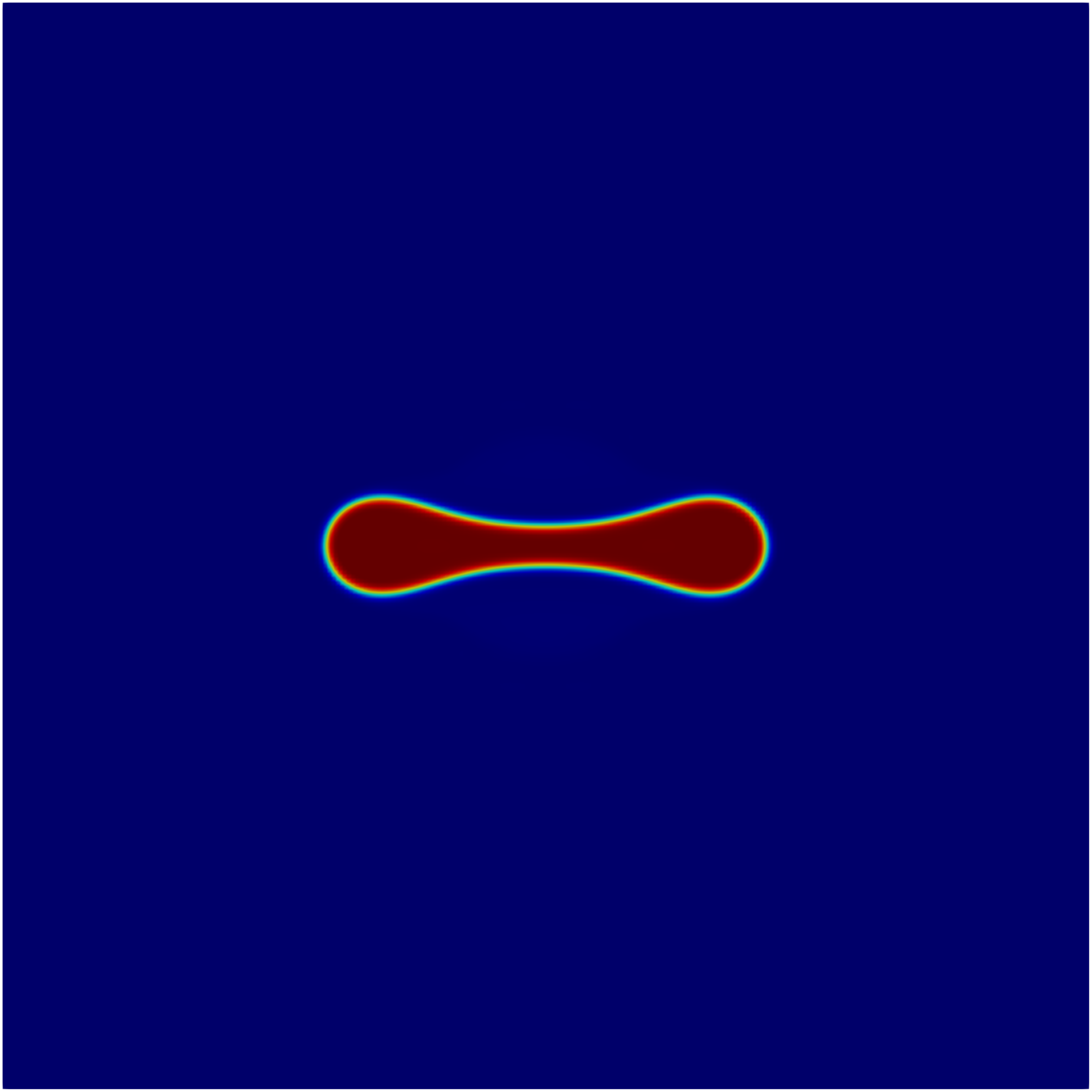} \end{subfigure}
\begin{subfigure}{0.185\columnwidth} \centering
\includegraphics[trim={0cm 0cm 0cm 0cm},clip,width=1\columnwidth]{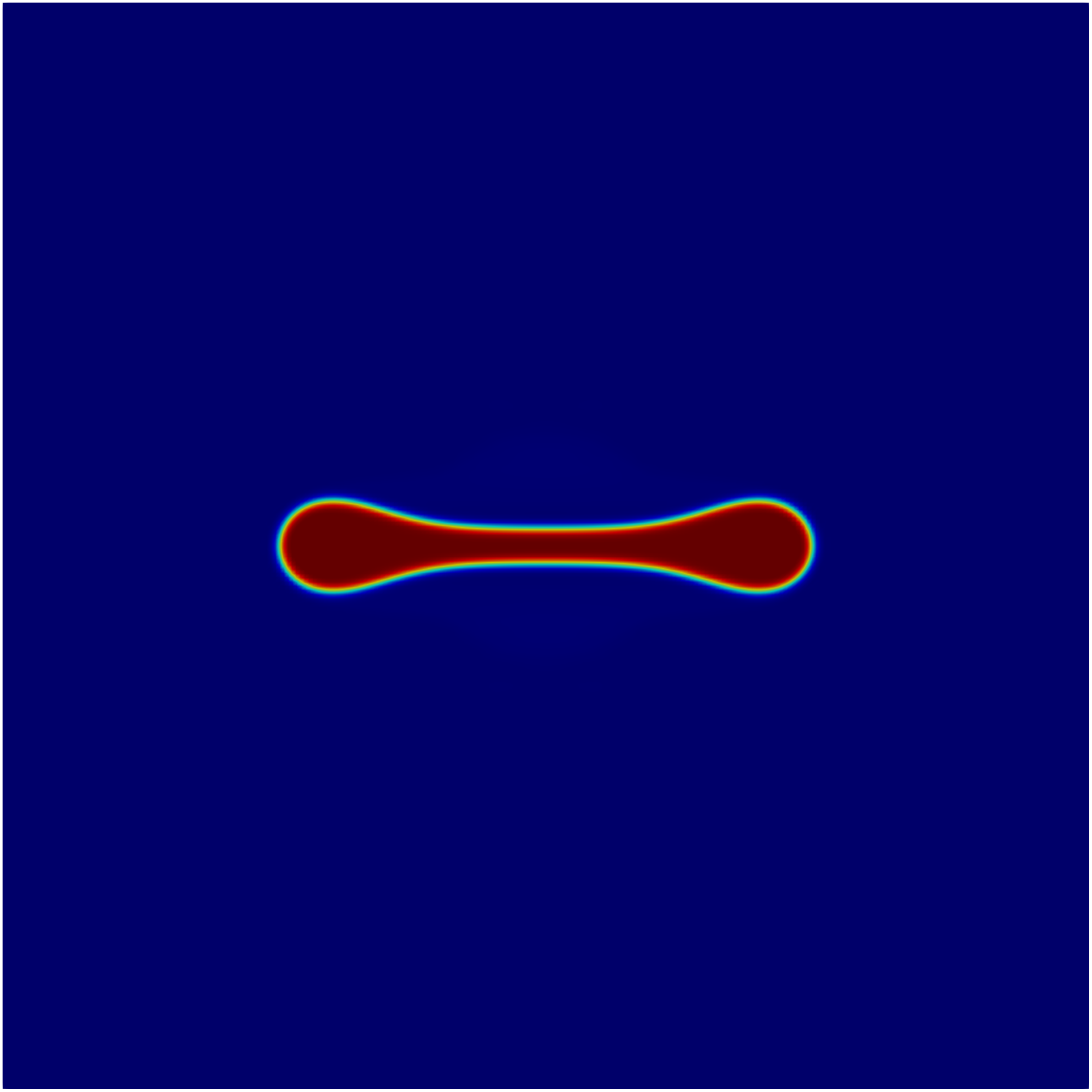}  \end{subfigure}
\begin{subfigure}{0.05\columnwidth} \centering
\includegraphics[trim={0cm 1cm 0cm 0cm},clip,width=1.2\columnwidth]{Images/Problem1_2DSq_UniRef/ColoarBar_Phi-eps-converted-to.pdf} \end{subfigure}
\\
\begin{subfigure}{0.185\columnwidth} \centering
\includegraphics[trim={0cm 0cm 0cm 0cm},clip,width=1\columnwidth]{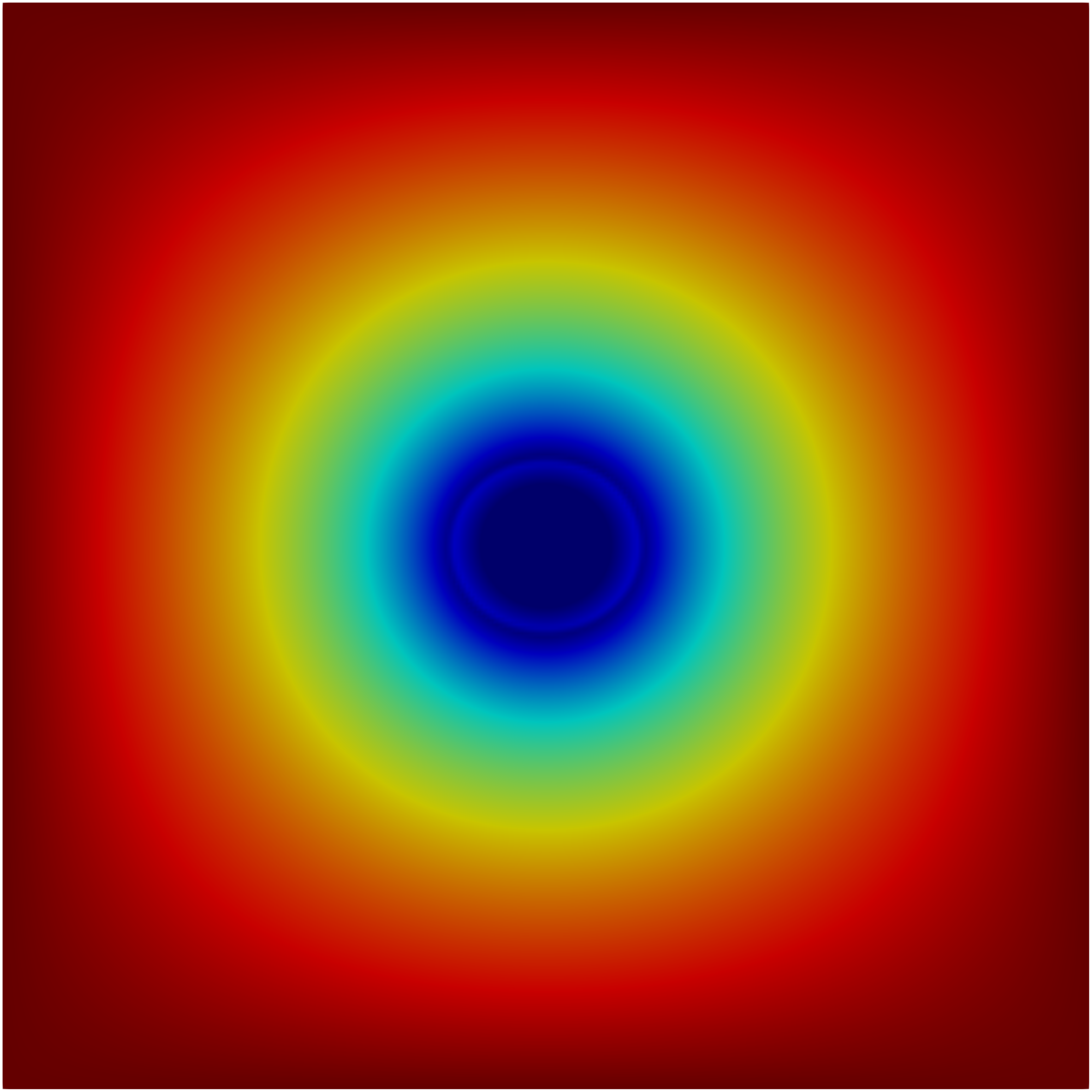} \caption{$t = 0$} \end{subfigure}
\begin{subfigure}{0.185\columnwidth} \centering
\includegraphics[trim={0cm 0cm 0cm 0cm},clip,width=1\columnwidth]{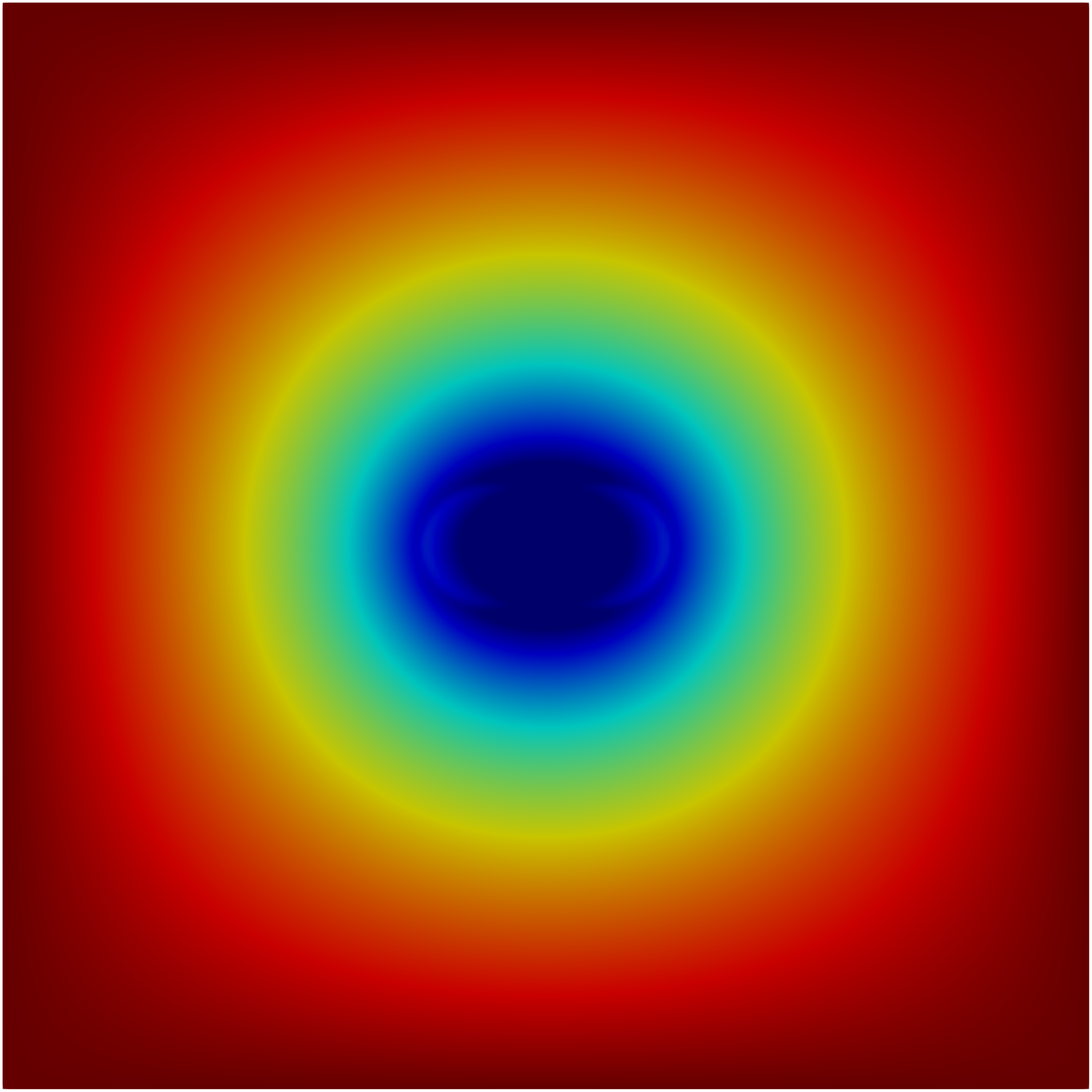} \caption{$t = 1$}  \end{subfigure}
\begin{subfigure}{0.185\columnwidth} \centering
\includegraphics[trim={0cm 0cm 0cm 0cm},clip,width=1\columnwidth]{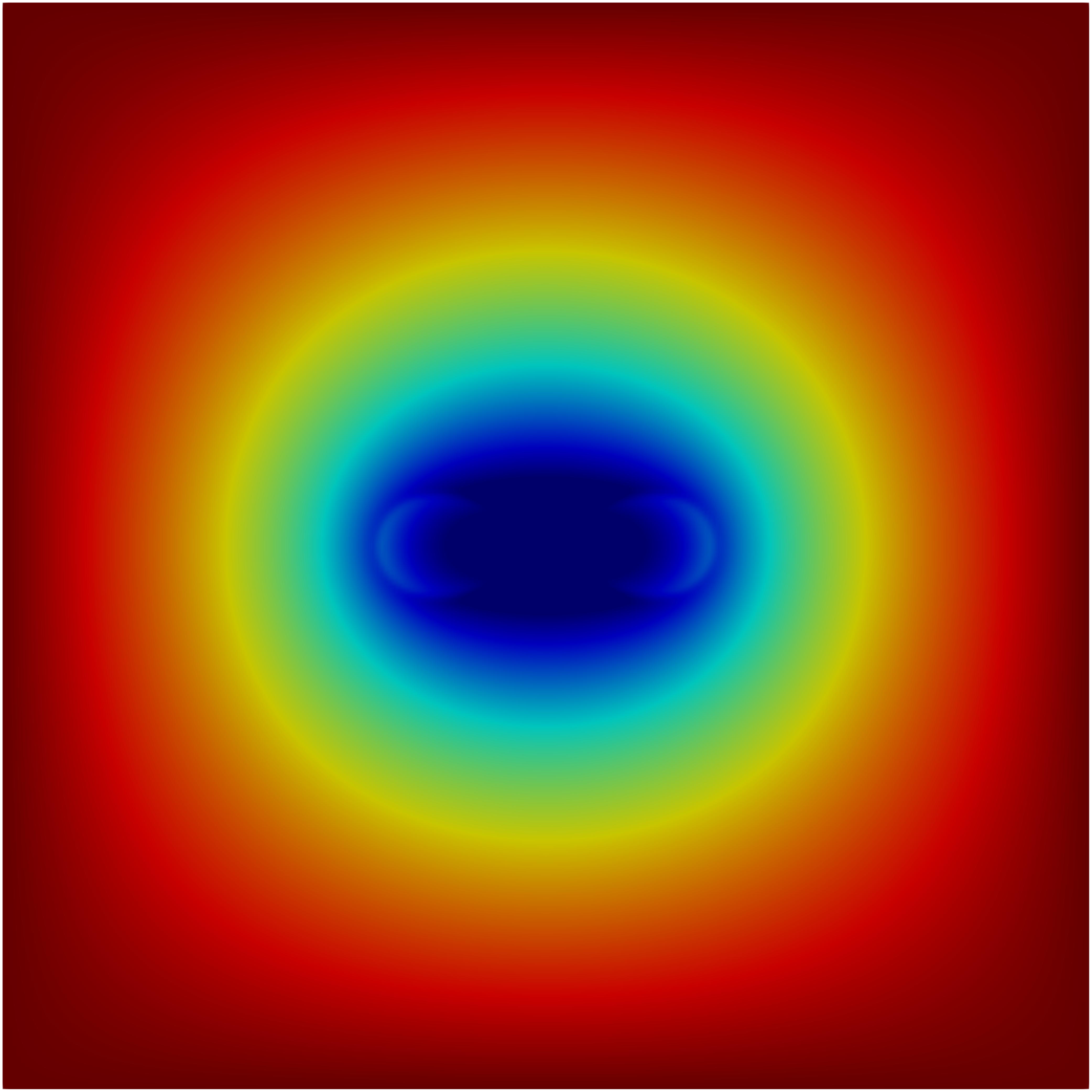} \caption{$t = 1.5$}\end{subfigure}
\begin{subfigure}{0.185\columnwidth} \centering
\includegraphics[trim={0cm 0cm 0cm 0cm},clip,width=1\columnwidth]{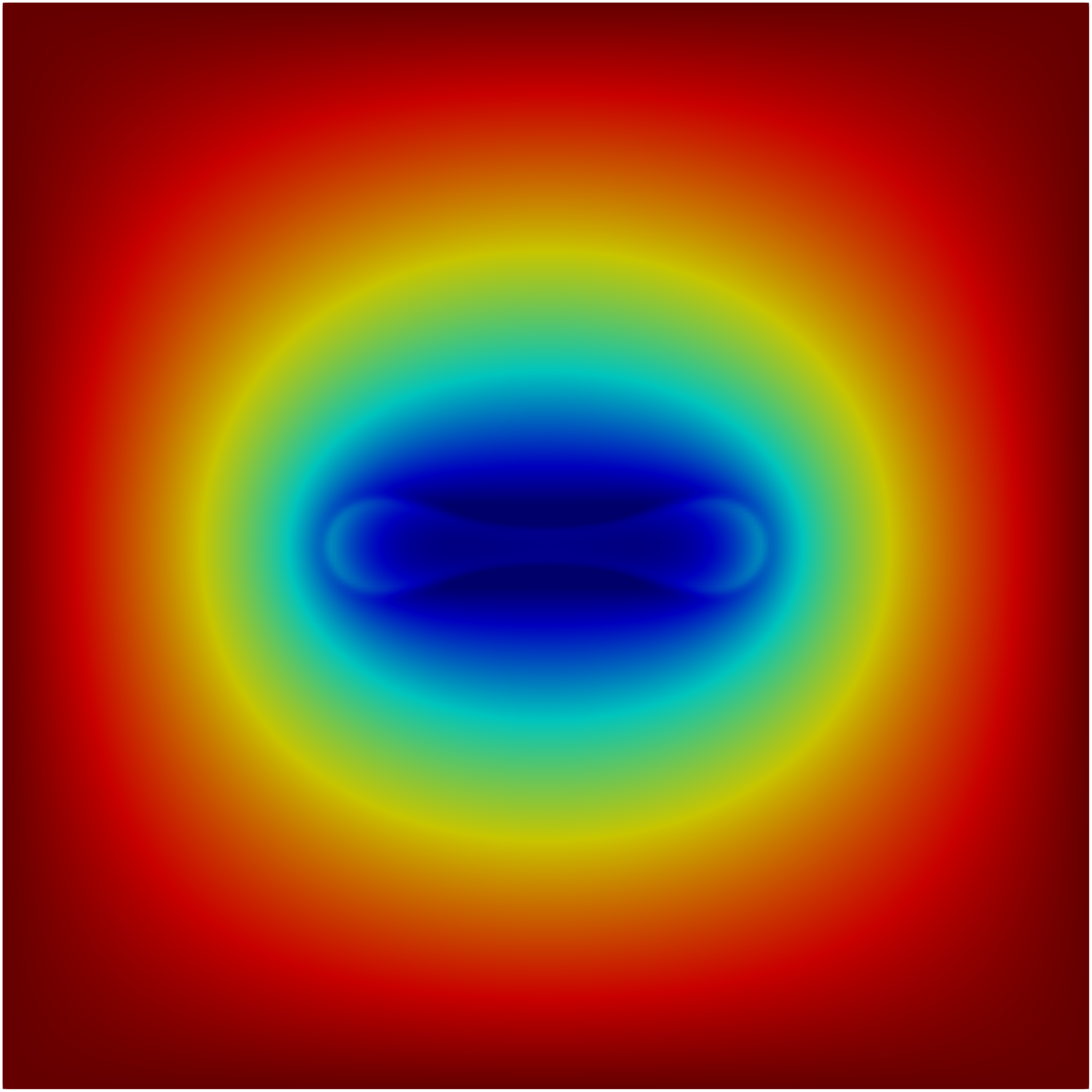} \caption{$t = 2$} \end{subfigure}
\begin{subfigure}{0.185\columnwidth} \centering
\includegraphics[trim={0cm 0cm 0cm 0cm},clip,width=1\columnwidth]{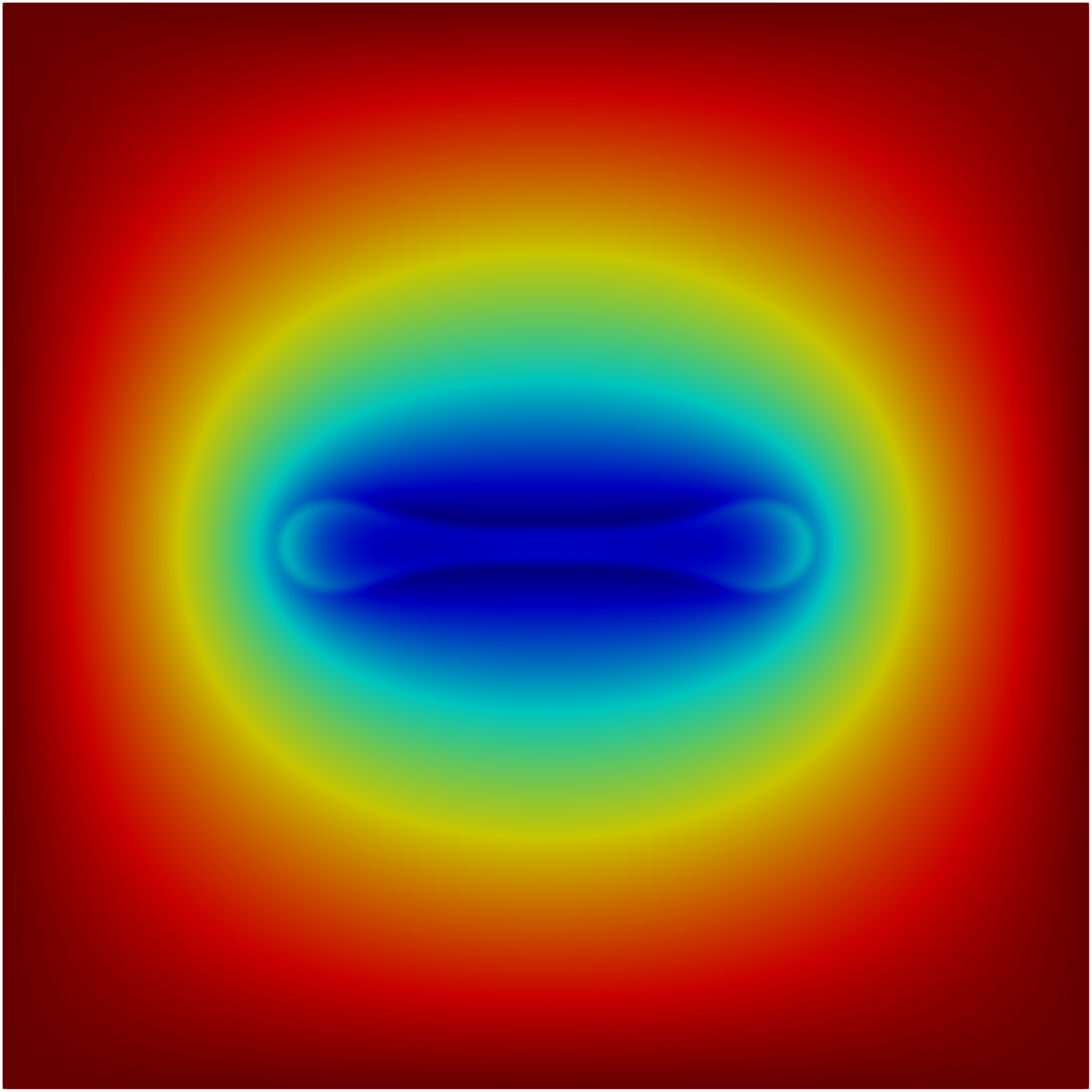} \caption{$t = 2.5$}\end{subfigure}
\begin{subfigure}{0.05\columnwidth} \centering
\includegraphics[trim={0cm -6cm 0cm 0cm},clip,width=1.07\columnwidth]{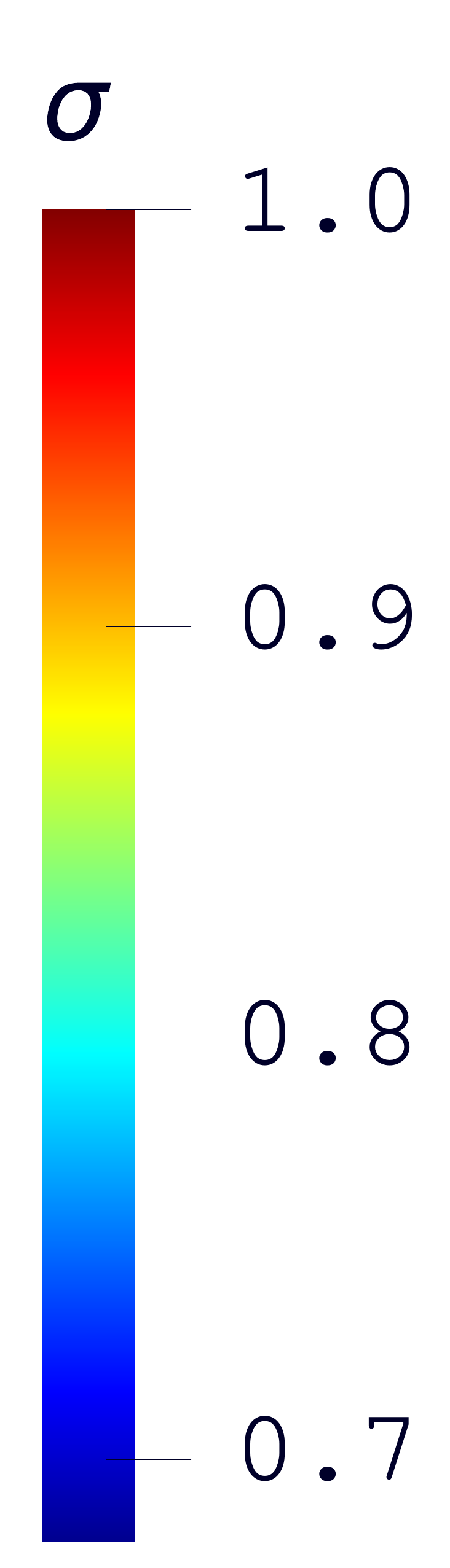} \end{subfigure}
\caption{Evolution of an elliptical tumor in square domain with $E = 5$, $\lambda = 0.002$, and $M = 5$, while keeping all other parameters unchanged (see Section~\ref{TG2D}). Top: adaptive THB-spline mesh configurations of degree $3$, $\ell = 6$, $m=2$, $\alpha = 0.01$, $\beta = 0.0001$ with a finest refinement level at $2^{10} \times 2^{10}$ mesh resolution. Middle: evolution of tumor geometry. Bottom: corresponding nutrient concentration. }
\label{figure_2_appendix}
\end{figure}

\begin{figure}[!h]\centering
\begin{subfigure}{0.185\columnwidth} \centering
\includegraphics[trim={0cm 0cm 0cm 0cm},clip,width=1\columnwidth]{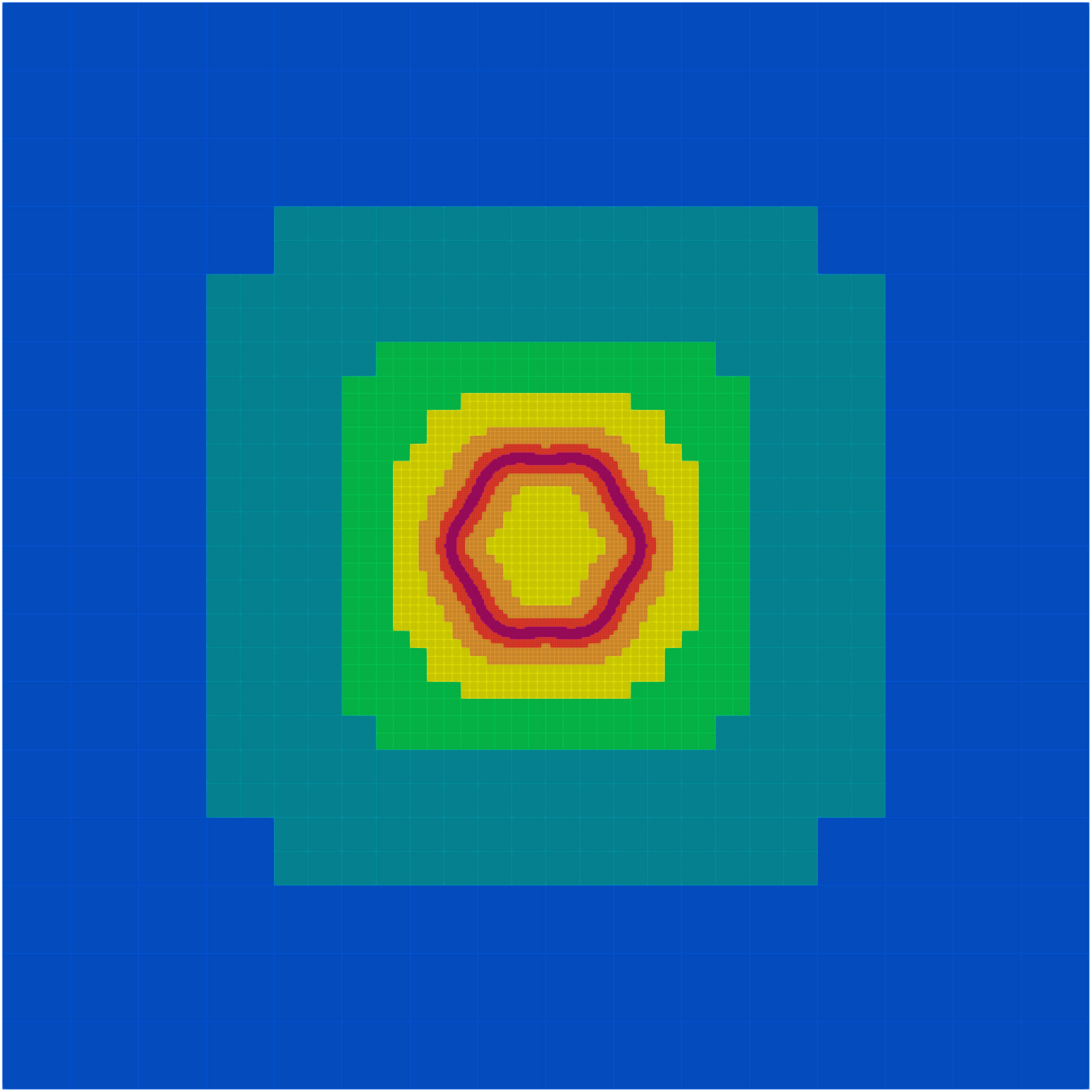}  \end{subfigure}
\begin{subfigure}{0.185\columnwidth} \centering
\includegraphics[trim={0cm 0cm 0cm 0cm},clip,width=1\columnwidth]{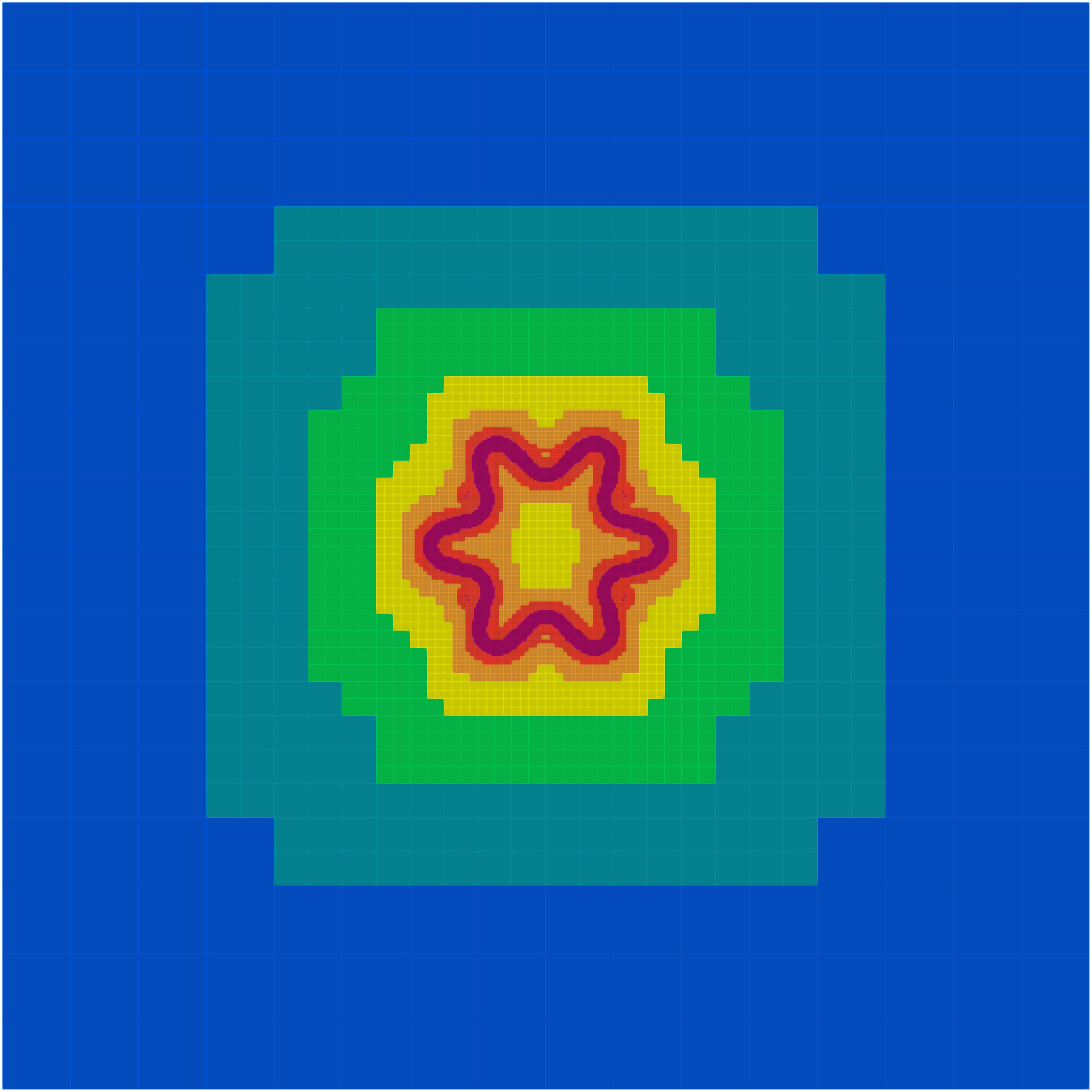}   \end{subfigure}
\begin{subfigure}{0.185\columnwidth} \centering
\includegraphics[trim={0cm 0cm 0cm 0cm},clip,width=1\columnwidth]{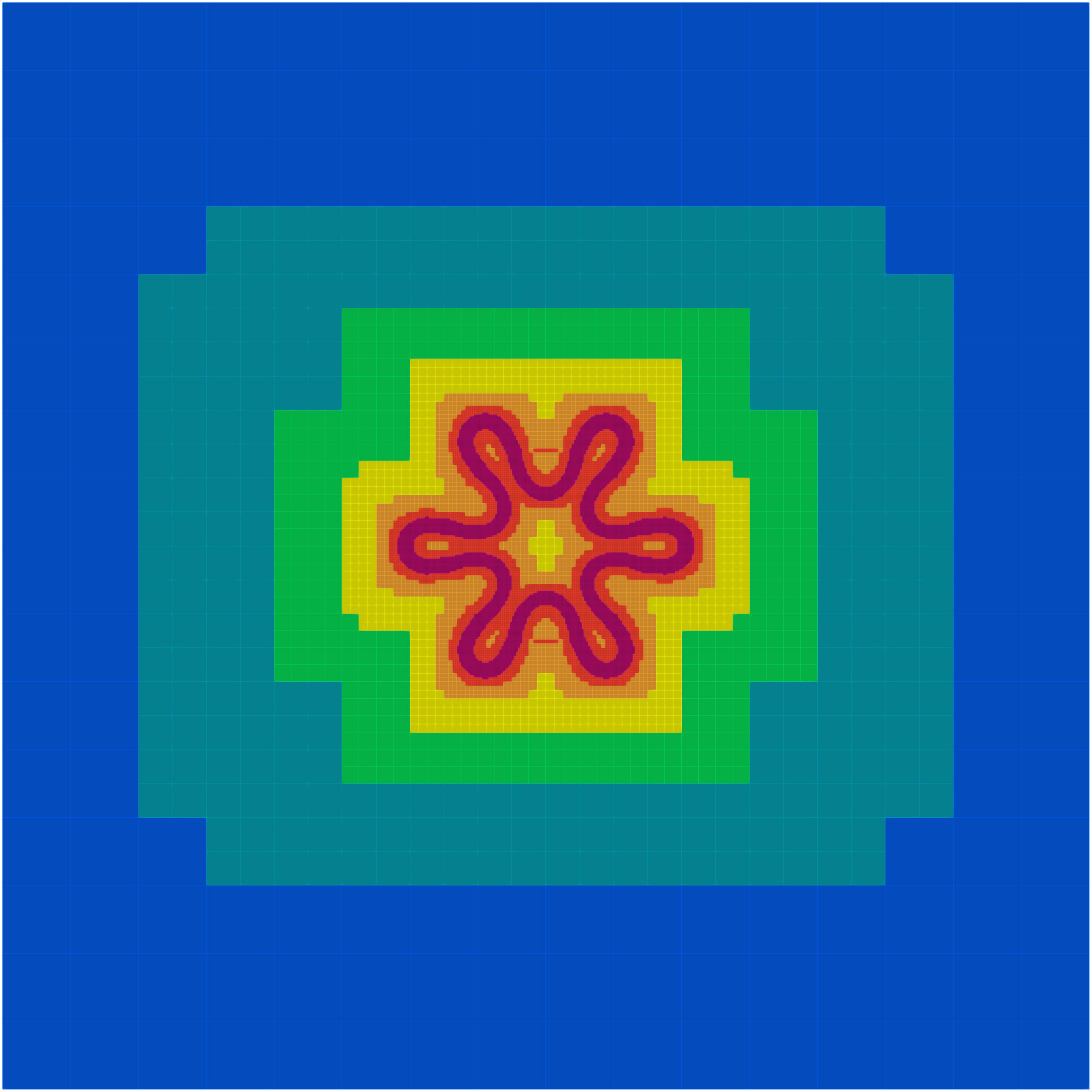} \end{subfigure}
\begin{subfigure}{0.185\columnwidth} \centering
\includegraphics[trim={0cm 0cm 0cm 0cm},clip,width=1\columnwidth]{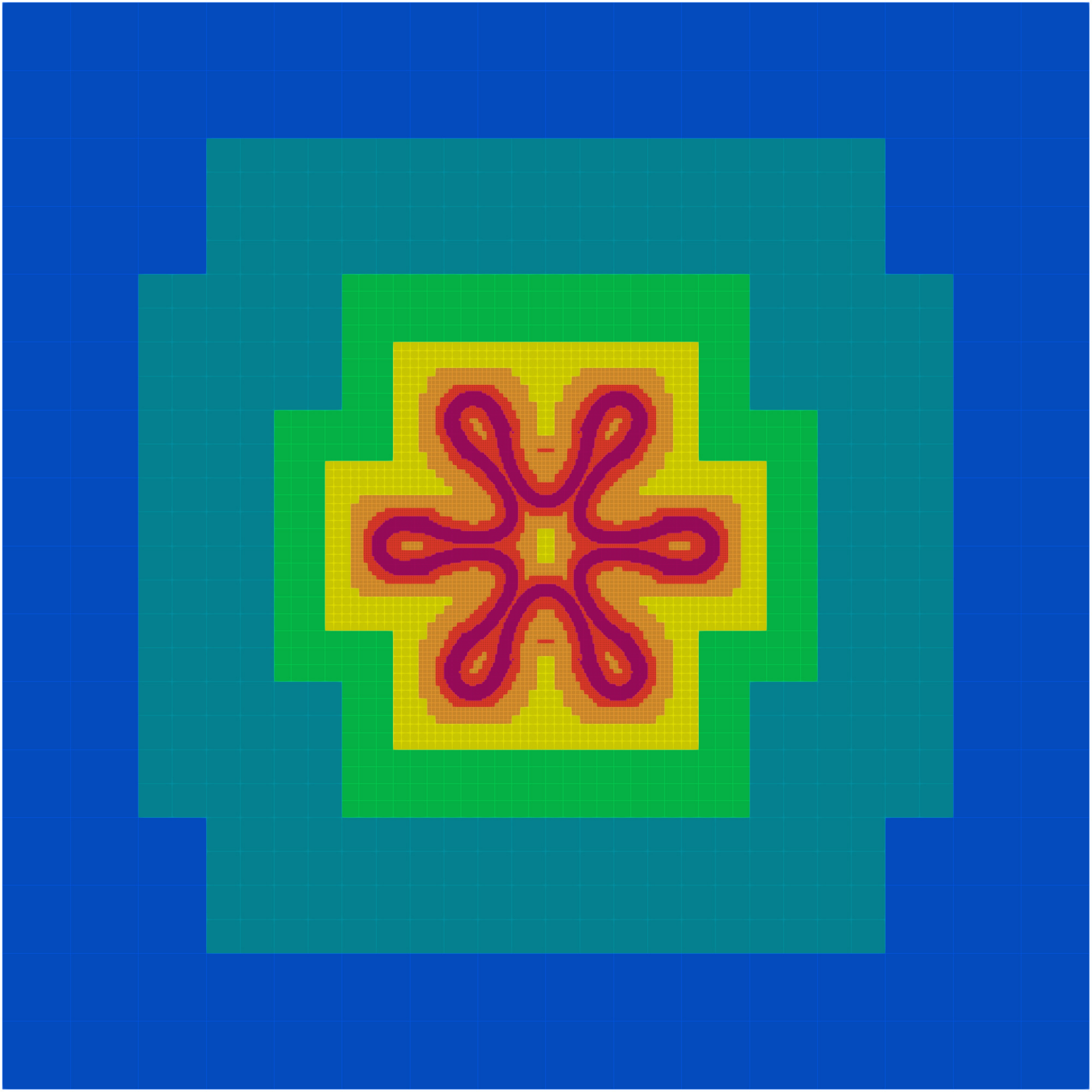} \end{subfigure}
\begin{subfigure}{0.185\columnwidth} \centering
\includegraphics[trim={0cm 0cm 0cm 0cm},clip,width=1\columnwidth]{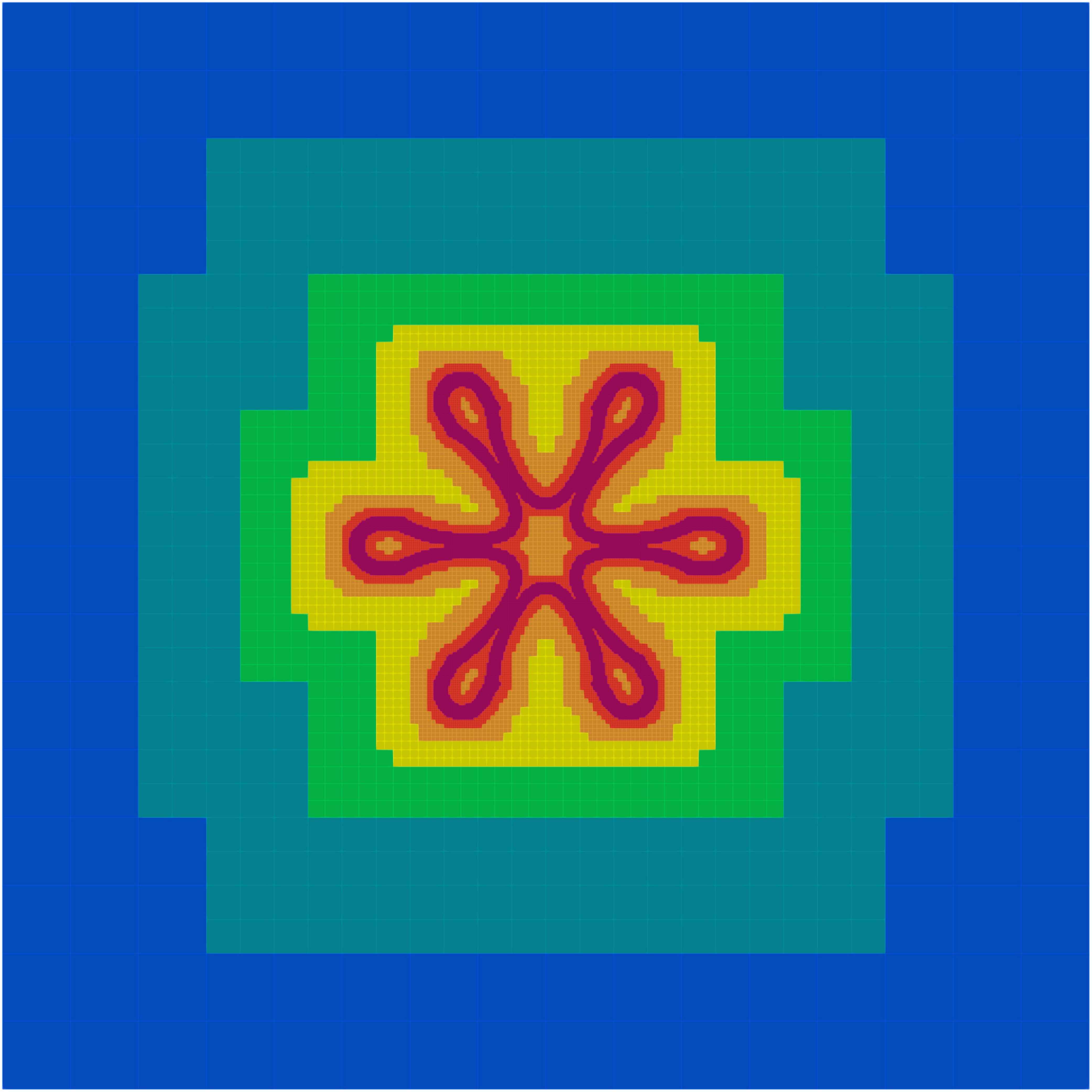}  \end{subfigure}
\begin{subfigure}{0.05\columnwidth} \centering
\includegraphics[trim={1cm 0cm 0cm 0cm},clip,width=0.84\columnwidth]{Images/Problem1_THBRes/ColorBar_Mesh_Lvl6-eps-converted-to.pdf} \end{subfigure}
\\ 
\begin{subfigure}{0.185\columnwidth} \centering
\includegraphics[trim={0cm 0cm 0cm 0cm},clip,width=1\columnwidth]{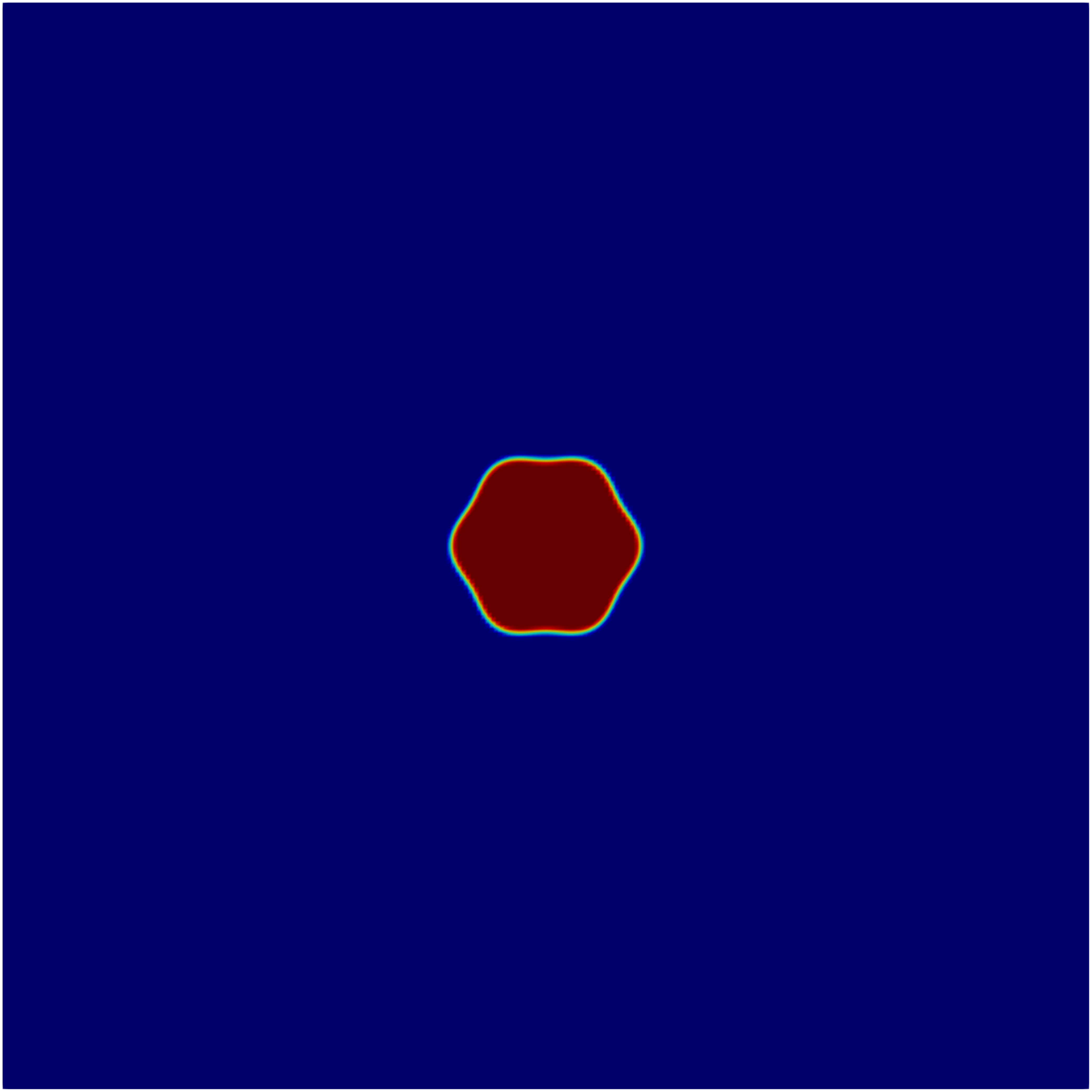}  \end{subfigure}
\begin{subfigure}{0.185\columnwidth} \centering
\includegraphics[trim={0cm 0cm 0cm 0cm},clip,width=1\columnwidth]{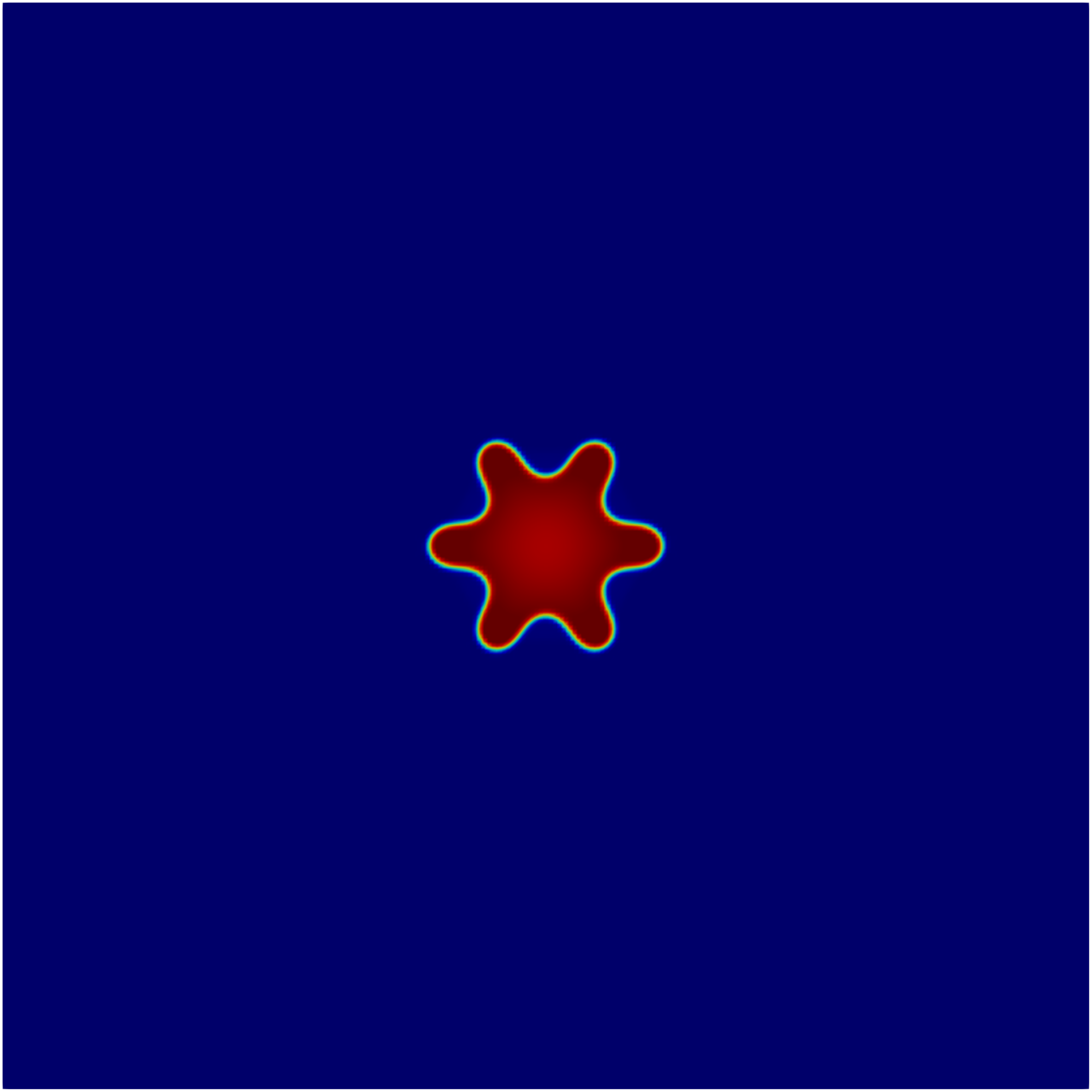}   \end{subfigure}
\begin{subfigure}{0.185\columnwidth} \centering
\includegraphics[trim={0cm 0cm 0cm 0cm},clip,width=1\columnwidth]{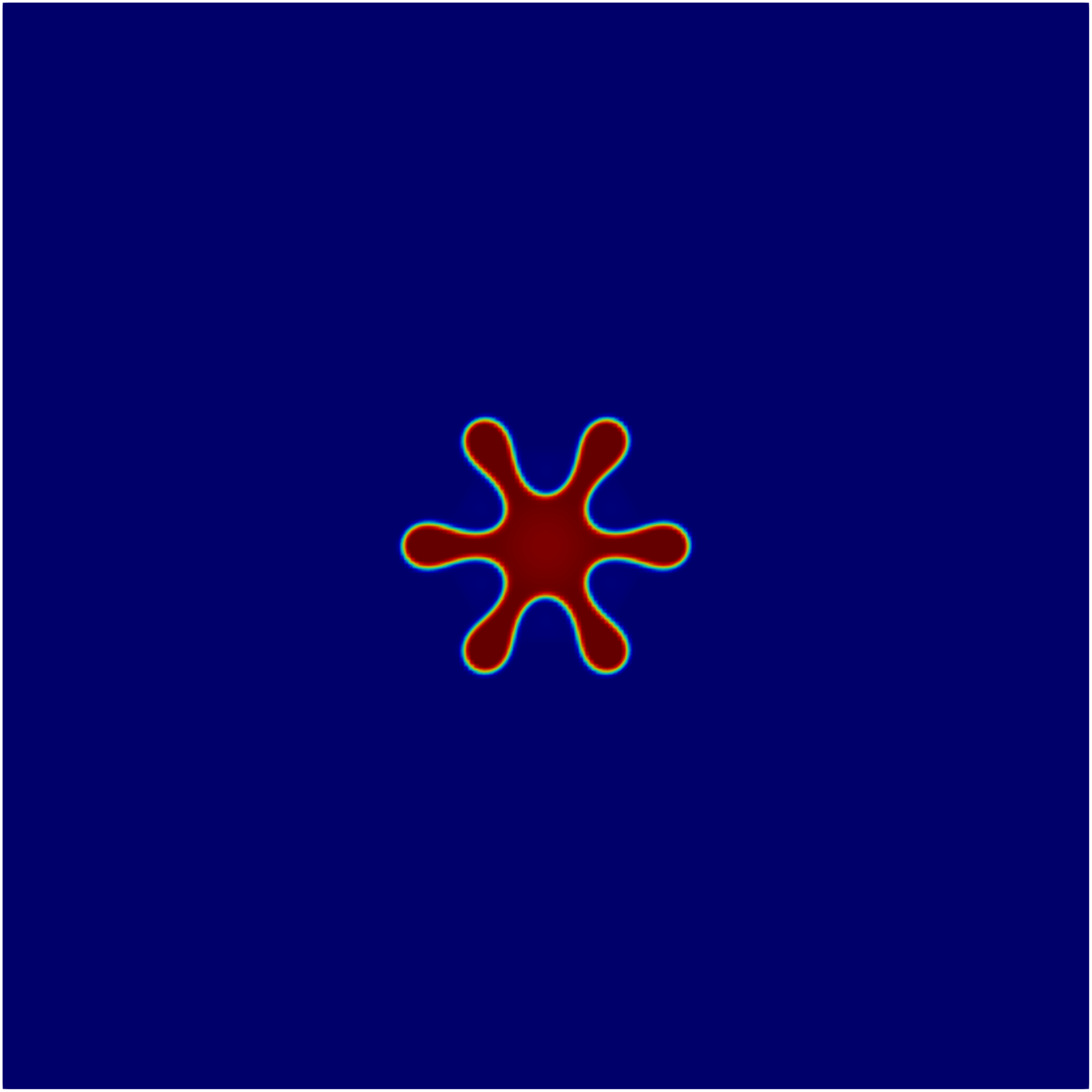} \end{subfigure}
\begin{subfigure}{0.185\columnwidth} \centering
\includegraphics[trim={0cm 0cm 0cm 0cm},clip,width=1\columnwidth]{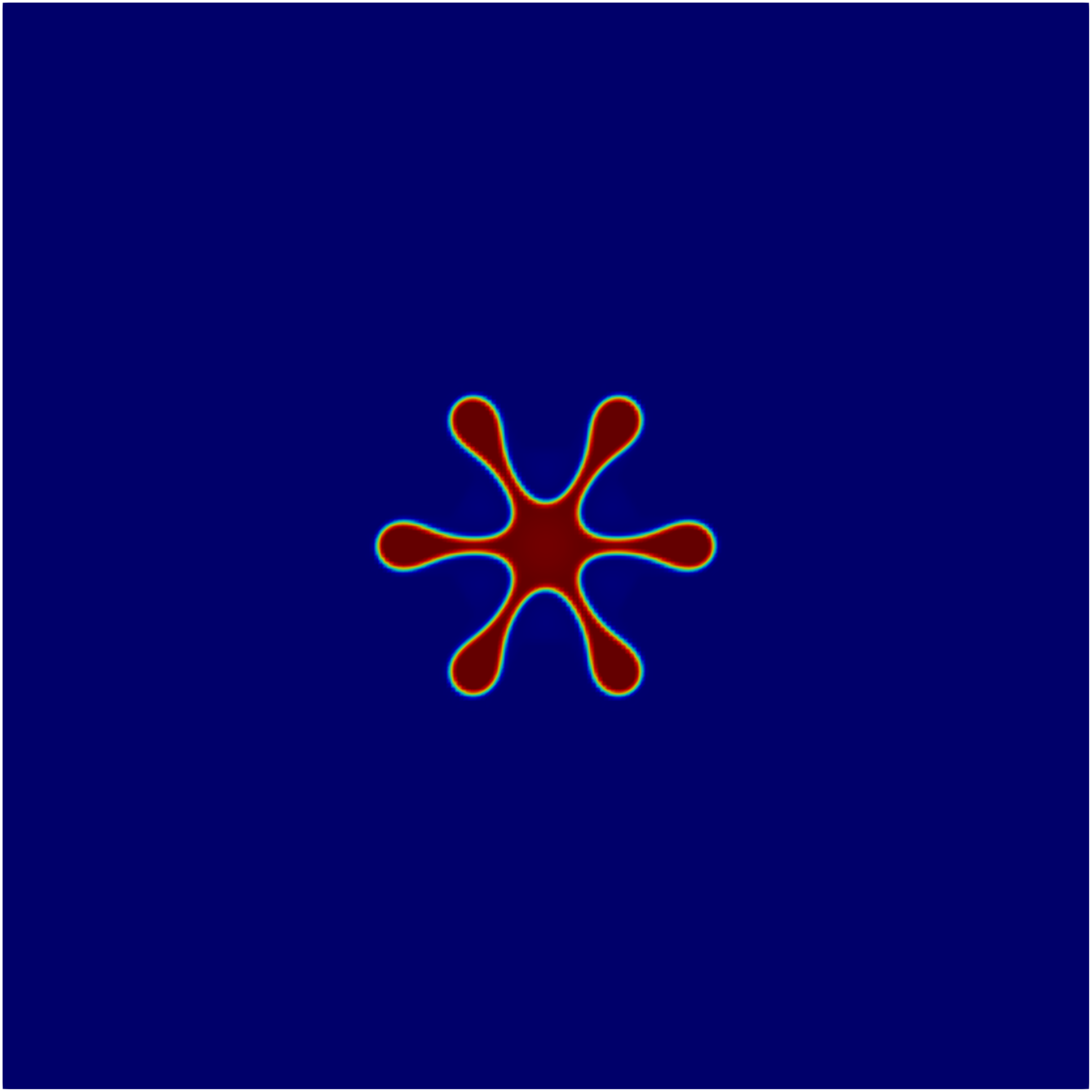} \end{subfigure}
\begin{subfigure}{0.185\columnwidth} \centering
\includegraphics[trim={0cm 0cm 0cm 0cm},clip,width=1\columnwidth]{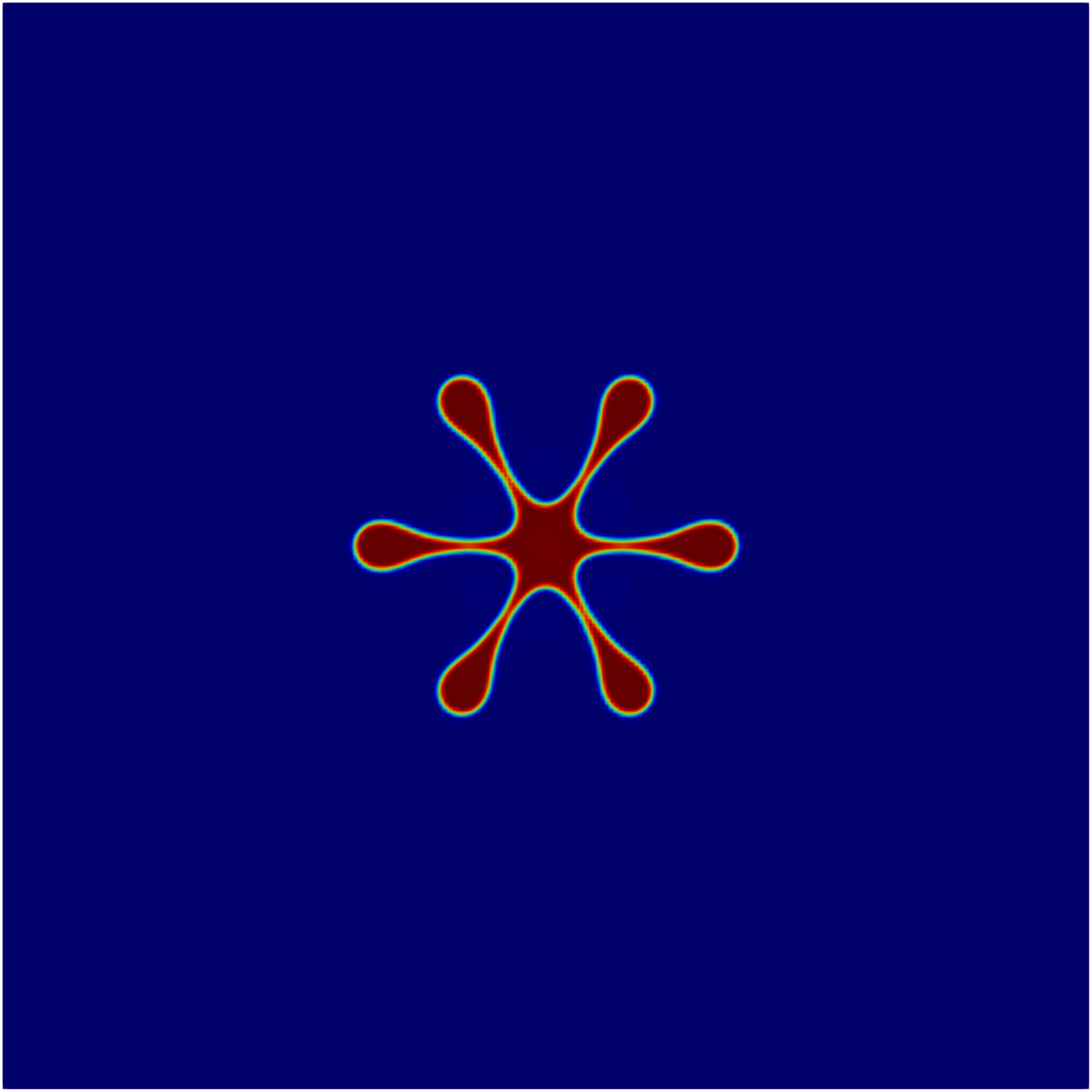}  \end{subfigure}
\begin{subfigure}{0.05\columnwidth} \centering
\includegraphics[trim={0cm 1cm 0cm 0cm},clip,width=1.2\columnwidth]{Images/Problem1_2DSq_UniRef/ColoarBar_Phi-eps-converted-to.pdf} \end{subfigure}
\\
\begin{subfigure}{0.185\columnwidth} \centering
\includegraphics[trim={0cm 0cm 0cm 0cm},clip,width=1\columnwidth]{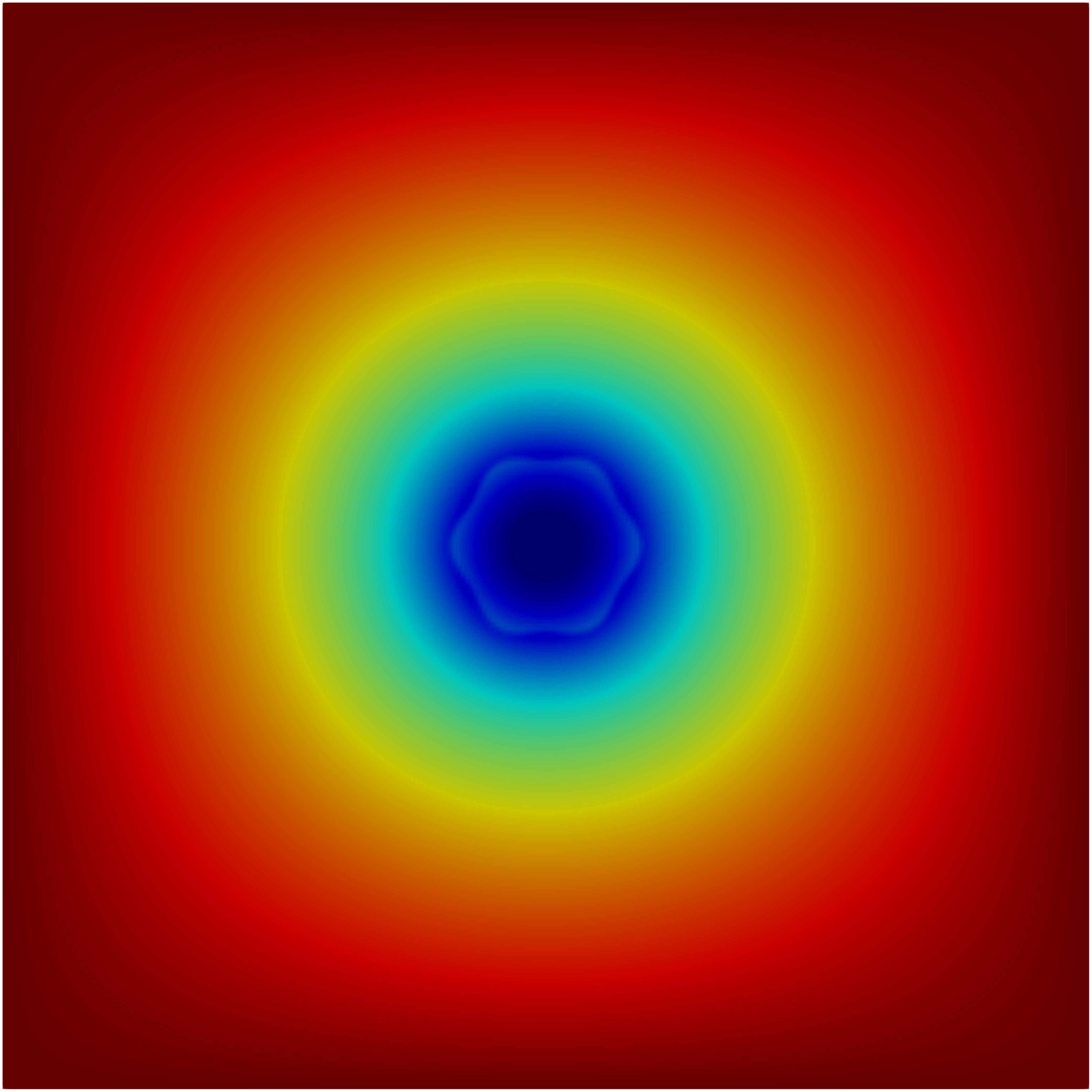} \caption{$t = 0$} \end{subfigure}
\begin{subfigure}{0.185\columnwidth} \centering
\includegraphics[trim={0cm 0cm 0cm 0cm},clip,width=1\columnwidth]{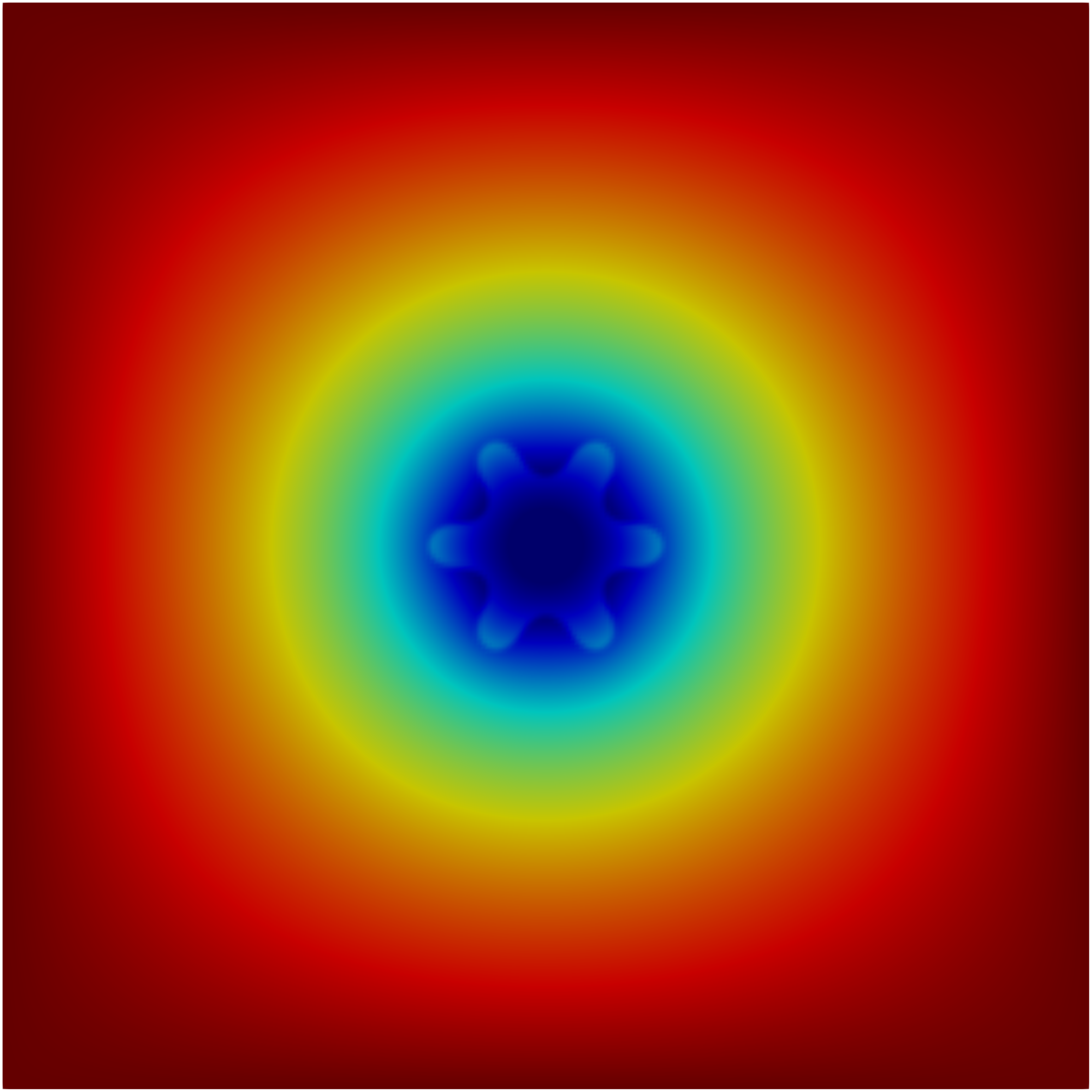} \caption{$t = 0.5$}  \end{subfigure}
\begin{subfigure}{0.185\columnwidth} \centering
\includegraphics[trim={0cm 0cm 0cm 0cm},clip,width=1\columnwidth]{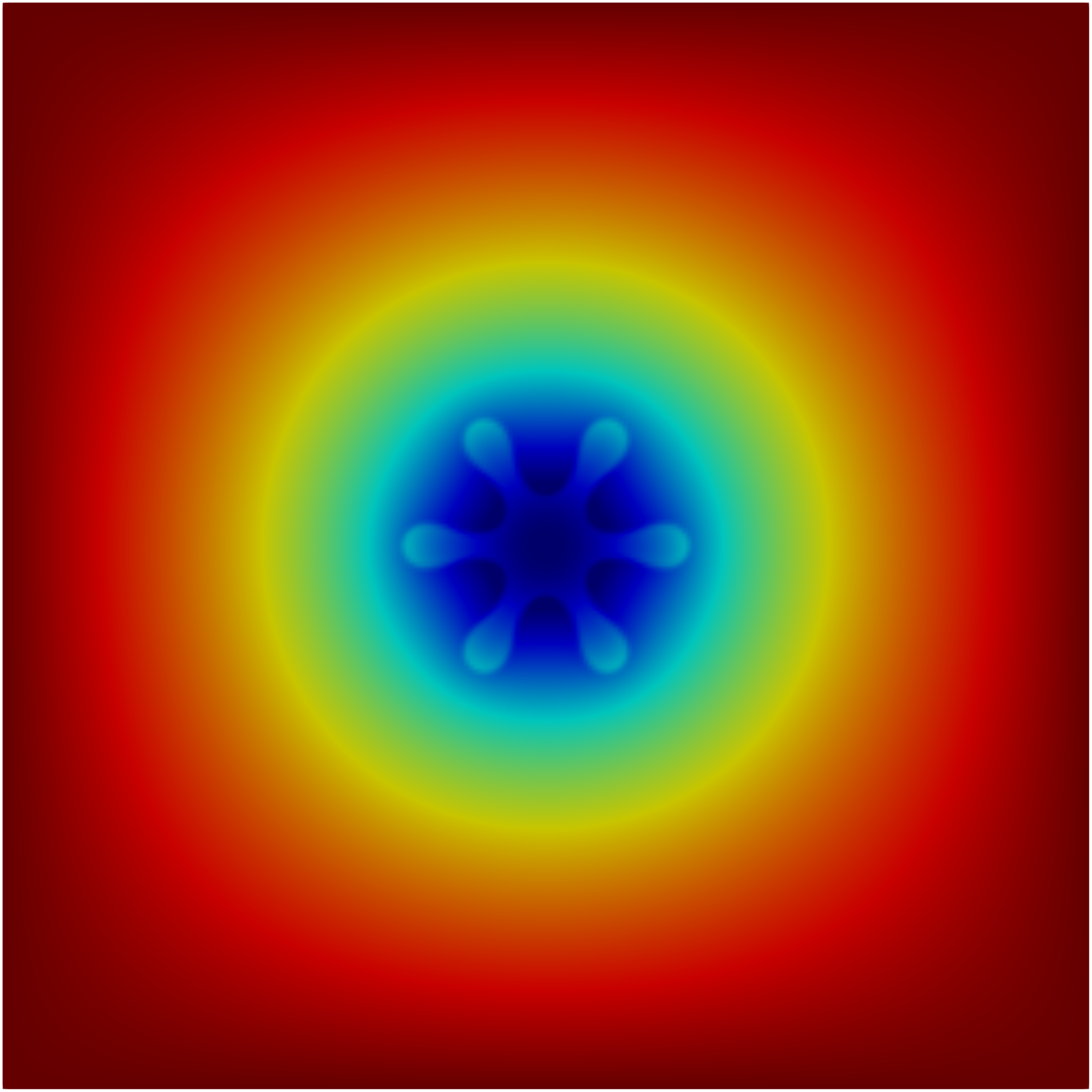} \caption{$t = 0.75$}\end{subfigure}
\begin{subfigure}{0.185\columnwidth} \centering
\includegraphics[trim={0cm 0cm 0cm 0cm},clip,width=1\columnwidth]{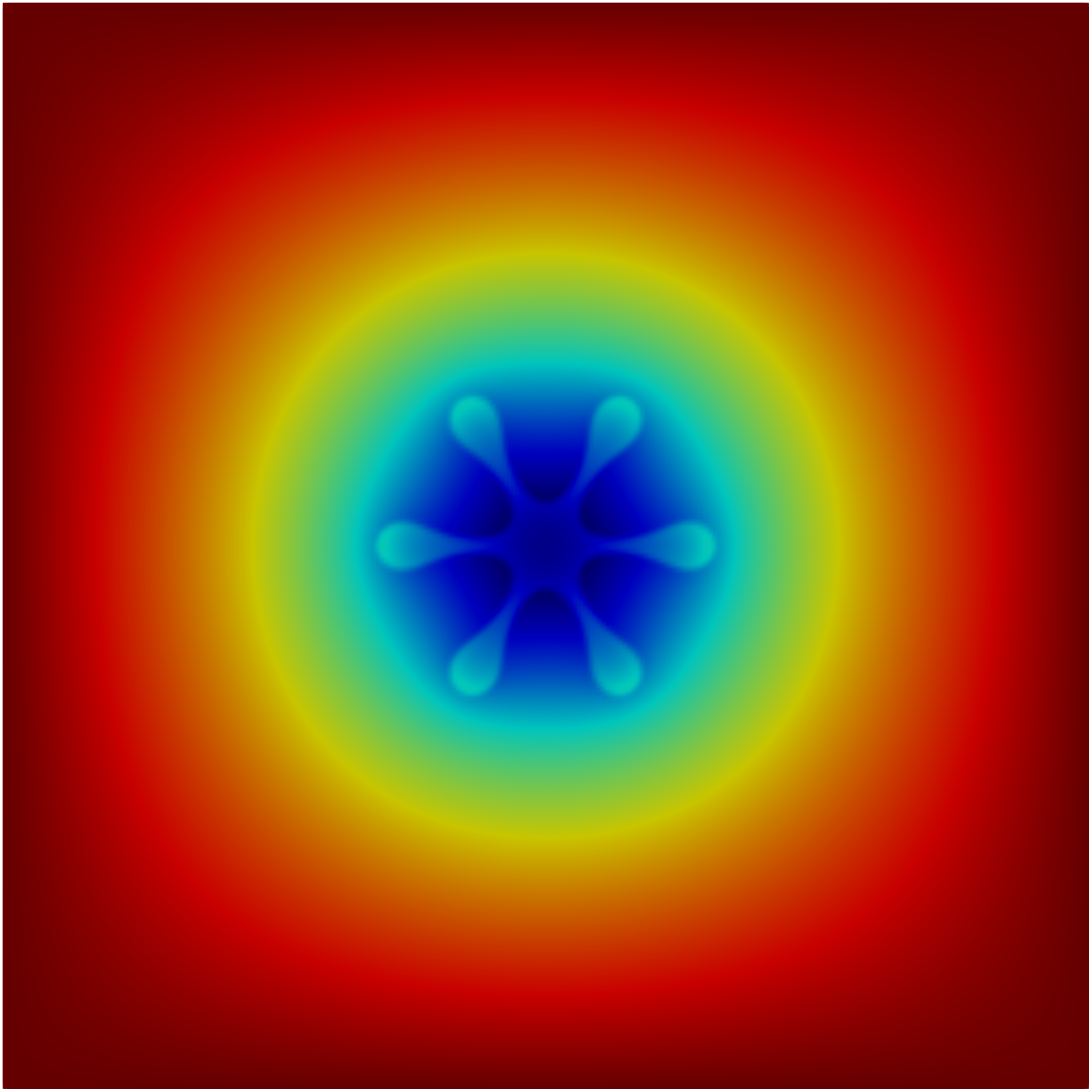} \caption{$t = 1$} \end{subfigure}
\begin{subfigure}{0.185\columnwidth} \centering
\includegraphics[trim={0cm 0cm 0cm 0cm},clip,width=1\columnwidth]{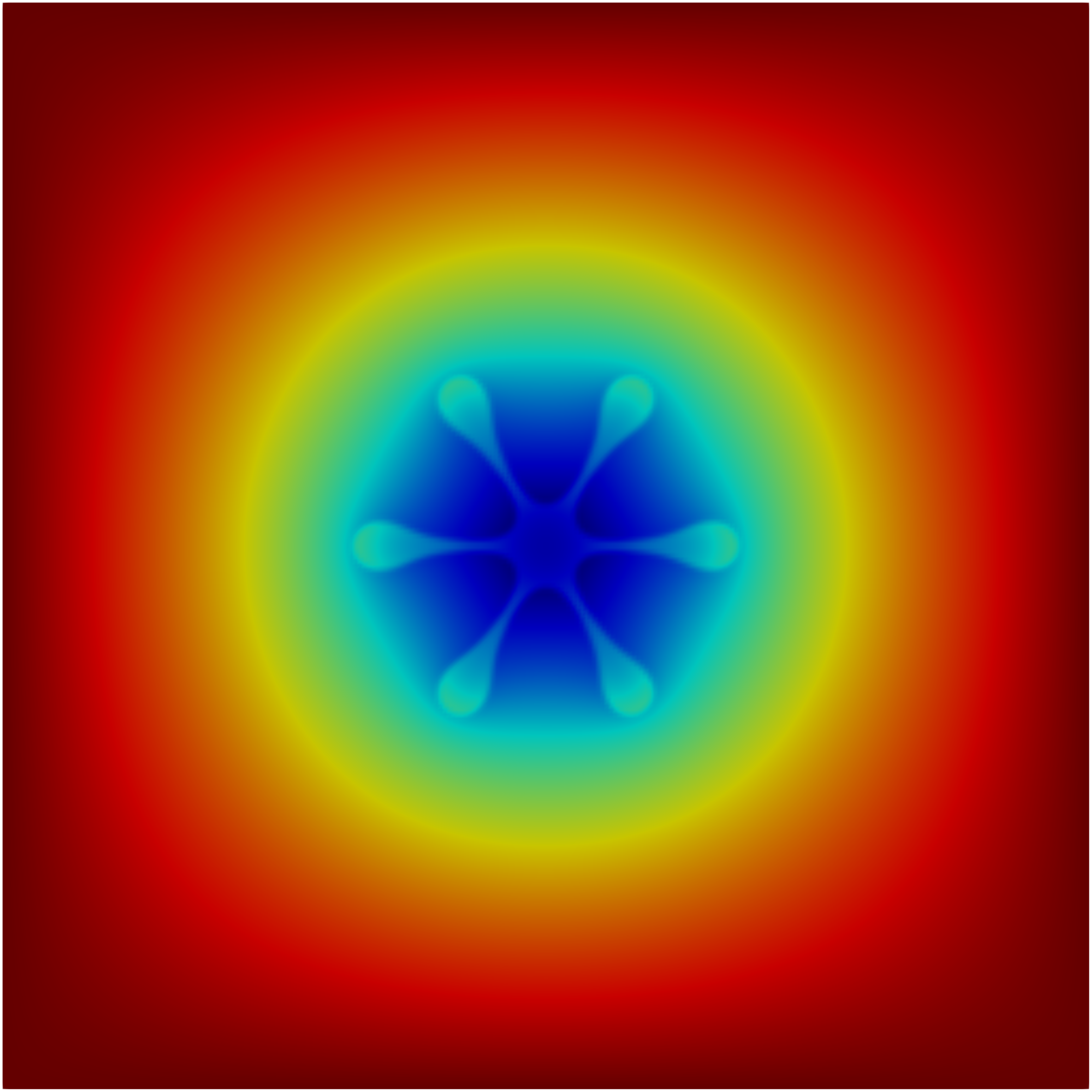} \caption{$t = 1.25$}\end{subfigure}
\begin{subfigure}{0.05\columnwidth} \centering
\includegraphics[trim={0cm -8cm 0cm 0cm},clip,width=1.15\columnwidth]{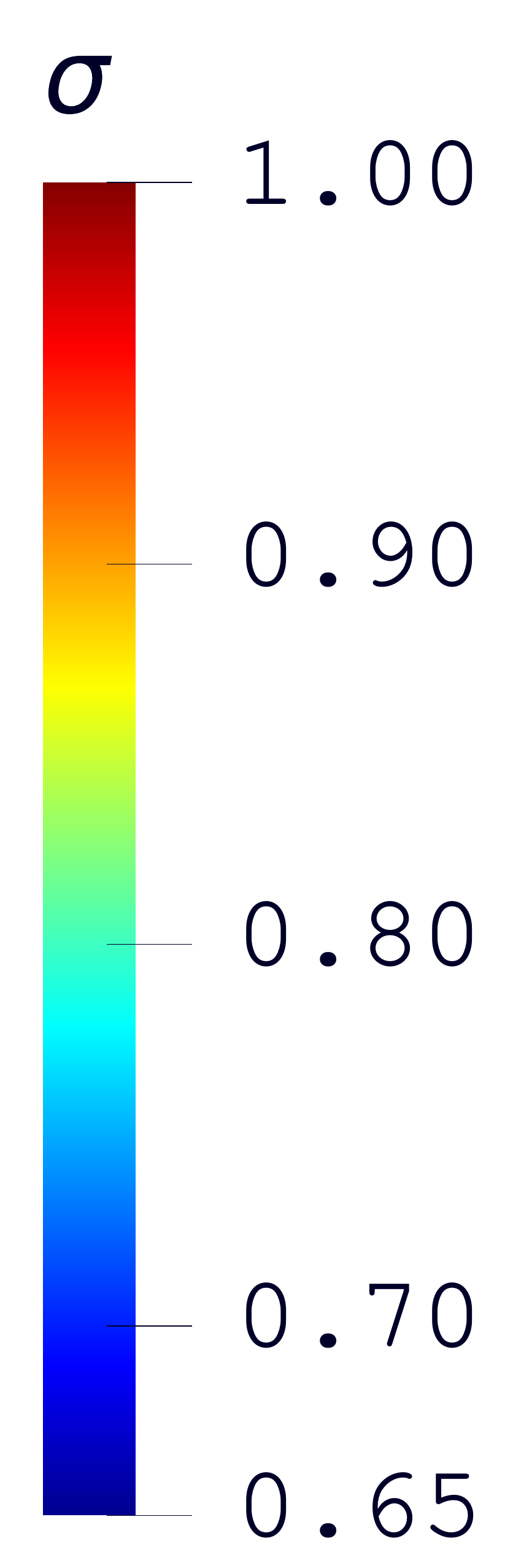} \end{subfigure}
\caption{Evolution of a tumor in a square domain with a perturbed elliptical initial profile defined by Eq.~\eqref{initial_cindition}, where $r(\mathbf{x}) = |\mathbf{x}| - \left[\frac{1}{2} + \frac{1}{40}\cos(6\theta) \right]$. Top: adaptive THB-spline mesh configurations of degree $3$, $\ell = 6$, $m=2$, $\alpha = 0.01$, $\beta = 0.0001$ with a finest refinement level at $2^{10} \times 2^{10}$ mesh resolution. Middle: evolution of tumor geometry. Bottom: corresponding nutrient concentration. }
\label{figure_3_appendix}
\end{figure}

In the second variation, the model parameters are kept identical to those stated in Section~\ref{TG2D}, and instead the initial condition is perturbed. 
Specifically, the initial condition follows the general form given in Eq.~\eqref{initial_cindition}, with the signed distance function modified to
\[
r(\mathbf{x}) = |\mathbf{x}| - \left[\frac{1}{2} + \frac{1}{40}\cos(6\theta) \right]
\]
The cosine perturbation introduces a six-fold angular modulation of amplitude $1/40$ onto an circular interface, breaking the radial symmetry of the initial tumor geometry. 
This structured perturbation triggers a qualitatively distinct growth morphology.
The tumor evolution and the corresponding THB-spline mesh adaptation and nutrient concentration are presented in Fig.~\ref{figure_3_appendix}, where the interface instabilities are amplified by the underlying dynamics, leading to a fingering pattern that differs from the morphologies observed in Figs.~\ref{figure_7_final} and \ref{figure_9_final}.

\begin{figure}[!h]\centering
\begin{subfigure}{1\columnwidth} \centering
\includegraphics[trim={0cm -0.25cm 0cm 0cm},clip,width=0.94\columnwidth]{Images/Problem1_2DSq_UniRef/Cub1024_THB_legends-eps-converted-to.pdf}
\end{subfigure}
\\
\begin{subfigure}{0.33\columnwidth} \centering
\includegraphics[trim={0cm 0cm 0cm 0cm},clip,width=0.99\columnwidth]{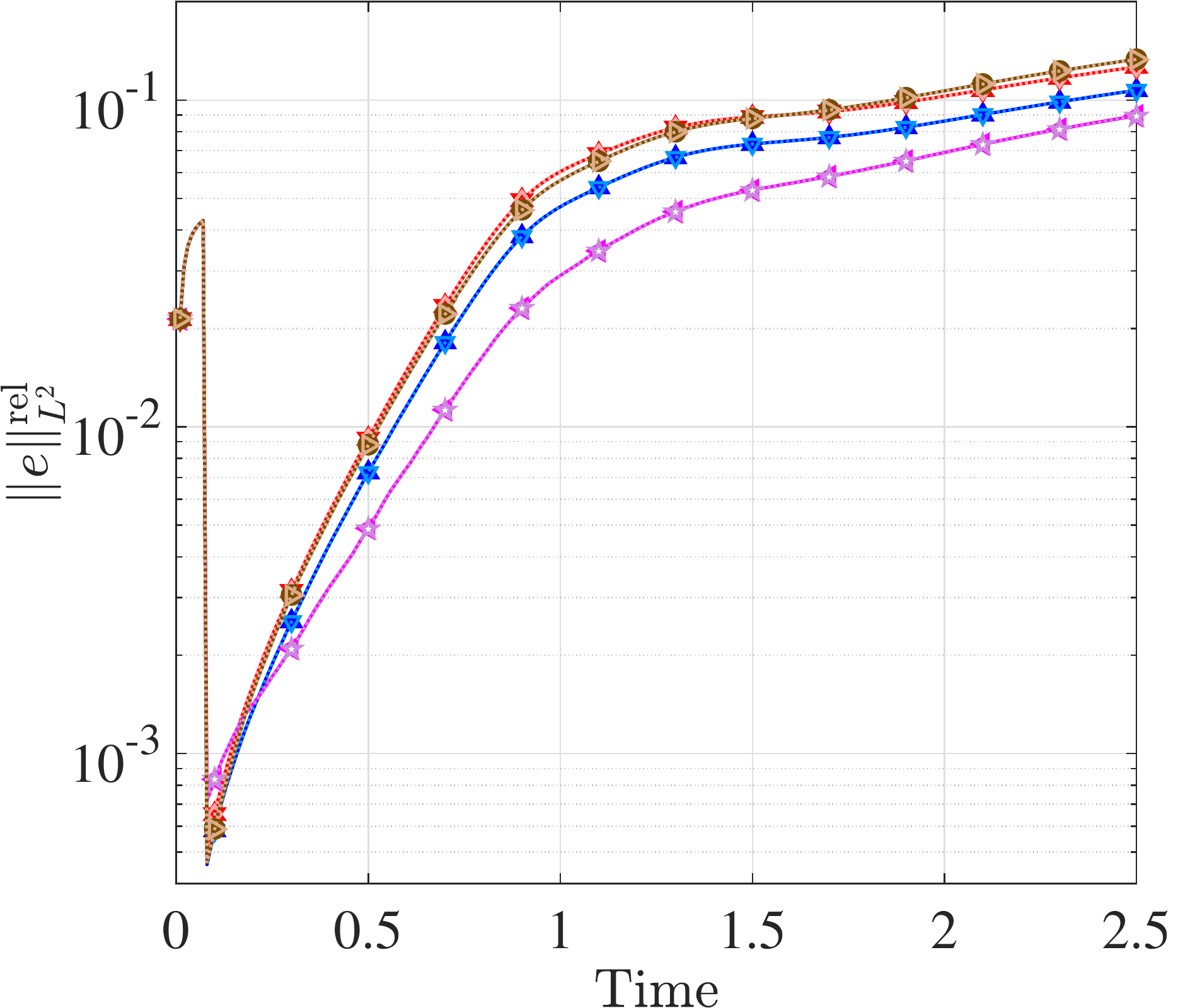}
\caption{$\| e \|_{L^2}^{\mathrm{rel}}$}
\end{subfigure}
\begin{subfigure}{0.33\columnwidth} \centering
\includegraphics[trim={0cm 0cm 0cm 0cm},clip,width=0.99\columnwidth]{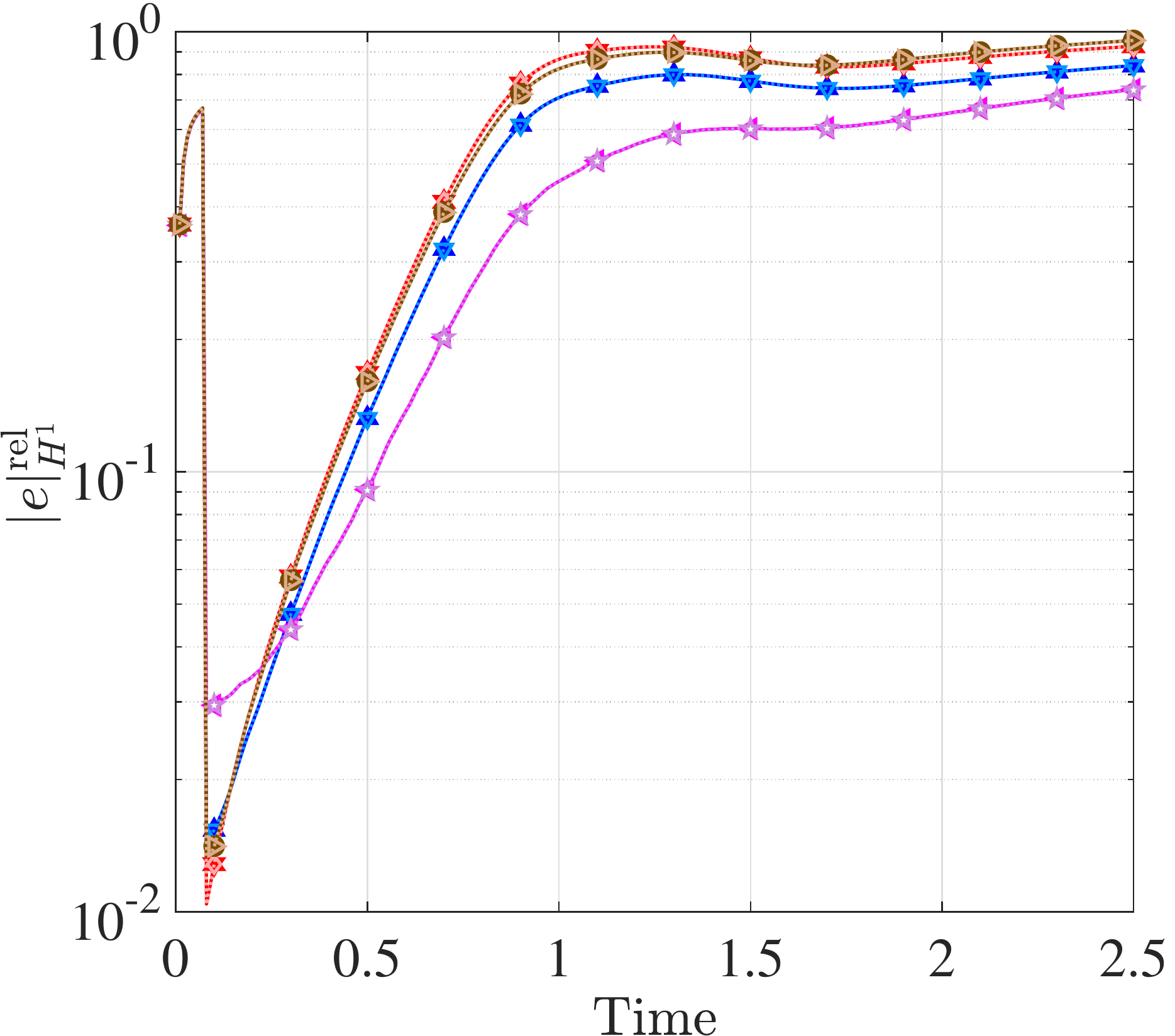}
\caption{$| e |_{H^1}^{\mathrm{rel}}$}
\end{subfigure}
\begin{subfigure}{0.33\columnwidth} \centering
\includegraphics[trim={0cm 0cm 0cm 0cm},clip,width=0.96\columnwidth]{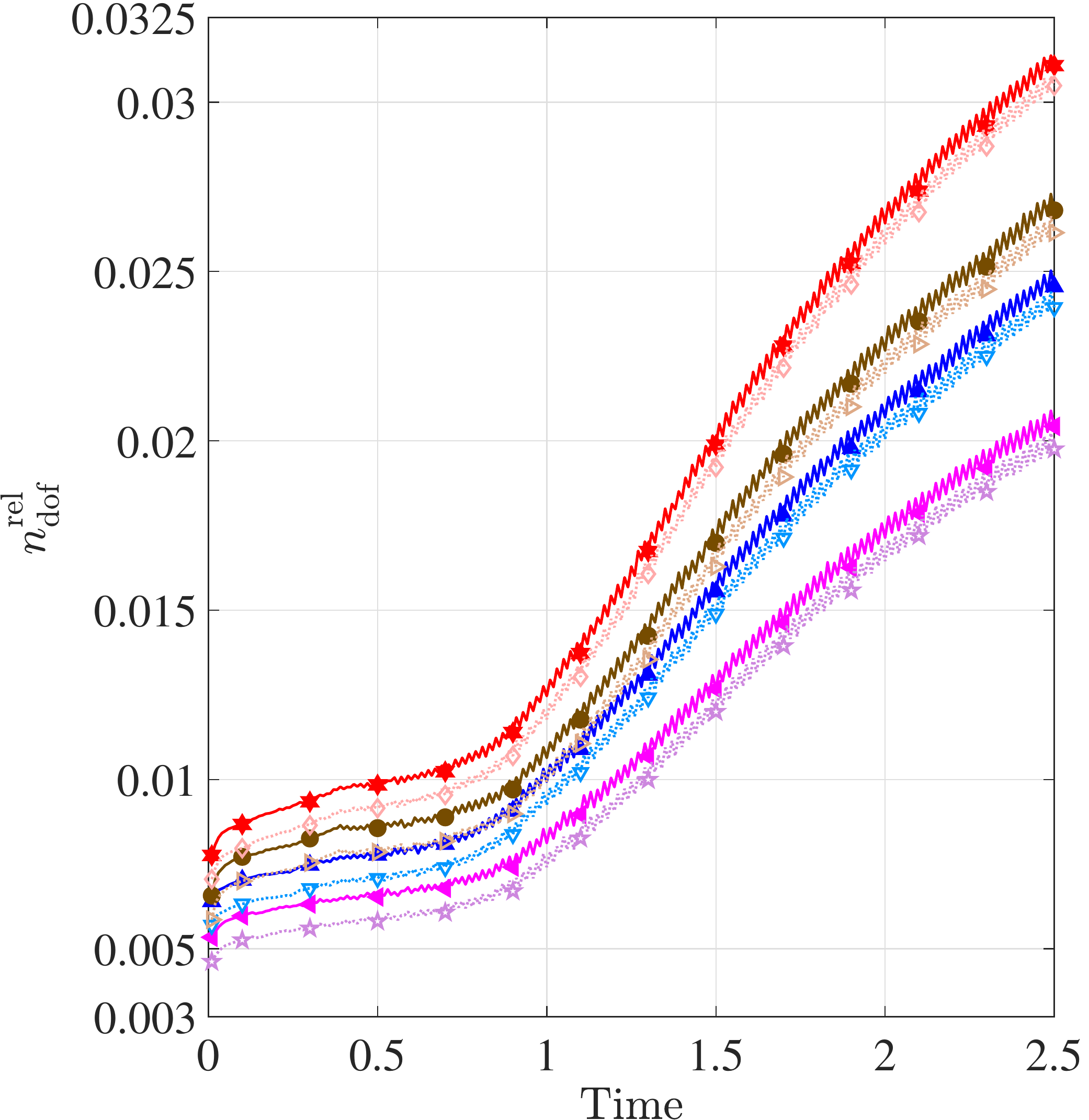}
\caption{$n_{\text{dof}}^{\text{rel}}$} %
\end{subfigure}
\caption{Relative $\| e \|_{L^2}^{\mathrm{rel}}$, $| e |_{H^1}^{\mathrm{rel}}$ , and $n_{\text{dof}}^{\mathrm{rel}}$ over the time period $\left[0, 2.5\right]$ for different adaptive THB mesh configurations of degree $2$ with varying hierarchical levels, admissibility class, and refinement and coarsening parameters. The finest level for all meshes is fixed at $2^{11}\times2^{11}$. The reference solution used to assess the accuracy in (a) and (b) corresponds to the uniform tensor-product B-spline mesh with $p=4$ and $1024$ elements in each direction as shown in Fig.~\ref{figure_1_final}, whereas the reference in (c) corresponds to a uniform tensor-product B-spline mesh with $p=2$ and $2048$ elements in each direction, indicating the relative computational effort, in terms of DOFs, required when the discretization is limited to quadratic B-splines on a uniform mesh.}
\label{figure_4_appendix}
\end{figure}

As a complementary remark, the performance of THB-meshes with quadratic B-splines ($p = 2$) is also examined in Fig.~\ref{figure_4_appendix} for readers who may prefer lower-order discretizations. 
Eight configurations are considered, varying in hierarchical levels $(\ell = 5 \text{ and } 6)$, refinement and coarsening parameters $\left((\alpha, \beta) = (0.1, 0.001) \text{ and } (0.01, 0.0001)\right)$, and admissibility class $(m = 2 \text{ and } 3)$. 
Following the tensor product mesh analyses in Section~\ref{TG2D}, the finest level for all THB-spline meshes is fixed at mesh resolution $2^{11} \times 2^{11}$ ($h_e=6/2048$), which is twice that of used for the cubic case. 
The results indicate that, irrespective of the mesh refinement level, quadratic B-splines are unable to match the accuracy achieved by cubic or quartic counterparts, and the errors reported in Fig.~\ref{figure_4_appendix} represent the best attainable solution quality within the quadratic approximation space for this configuration, see Fig.~\ref{figure_4_appendix}. 
This highlights the limitation of low-order bases in resolving high-gradient interfacial features, and may serve as a practical guideline for those willing to accept a moderate reduction in accuracy while working exclusively with quadratic elements. 
More broadly, this observation also highlights one of the advantages of IGA, where through $k$-refinement, higher orders of convergence and superior approximation per DOFs can be attained. 

\begin{figure}[!t]\centering
\begin{subfigure}{0.245\columnwidth}
\begin{subfigure}{\columnwidth} \centering
\includegraphics[trim={0cm 8cm 0cm 0cm},clip,width=0.59\columnwidth]{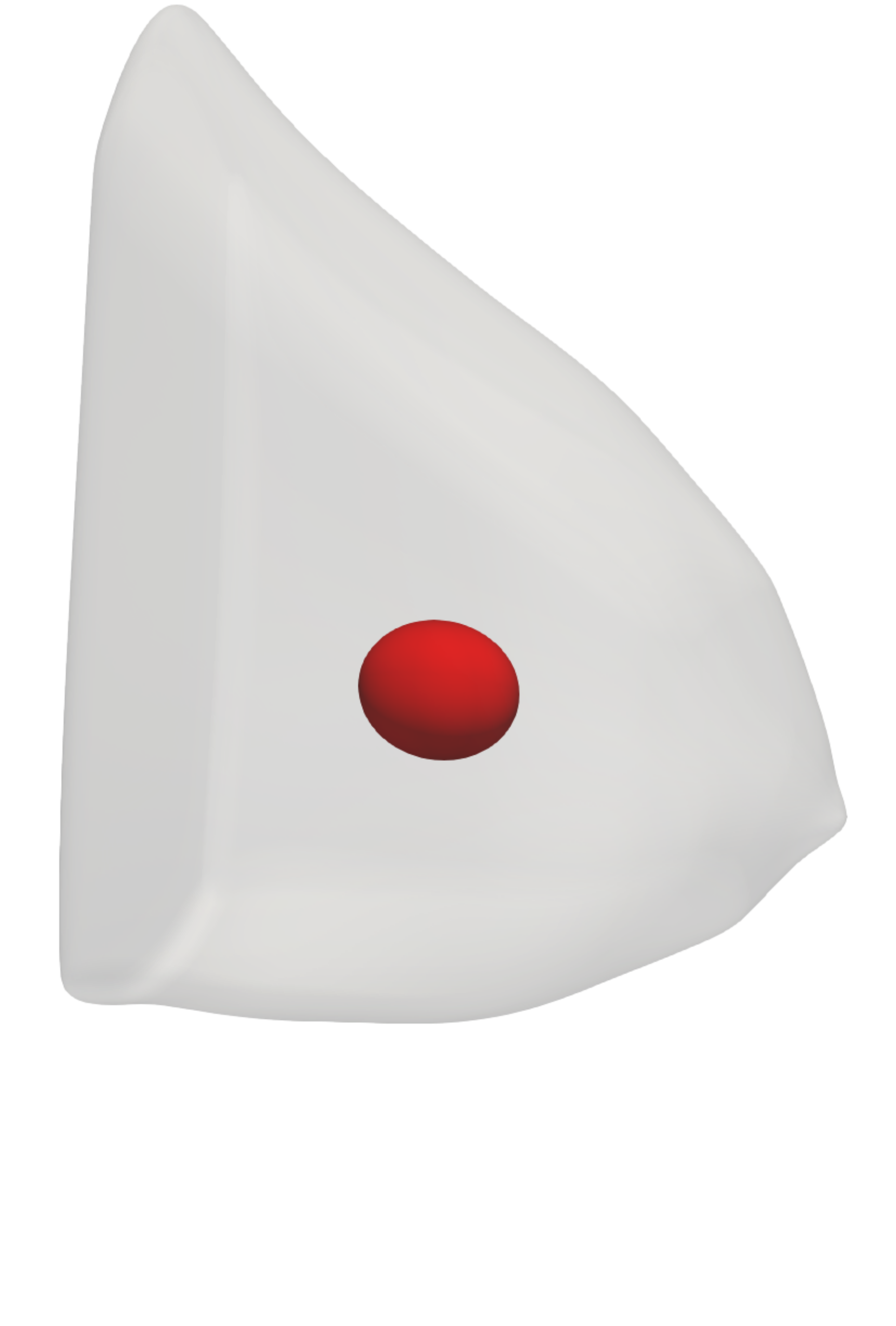}  
\includegraphics[trim={0cm 8cm 0cm 0cm},clip,width=0.39\columnwidth]{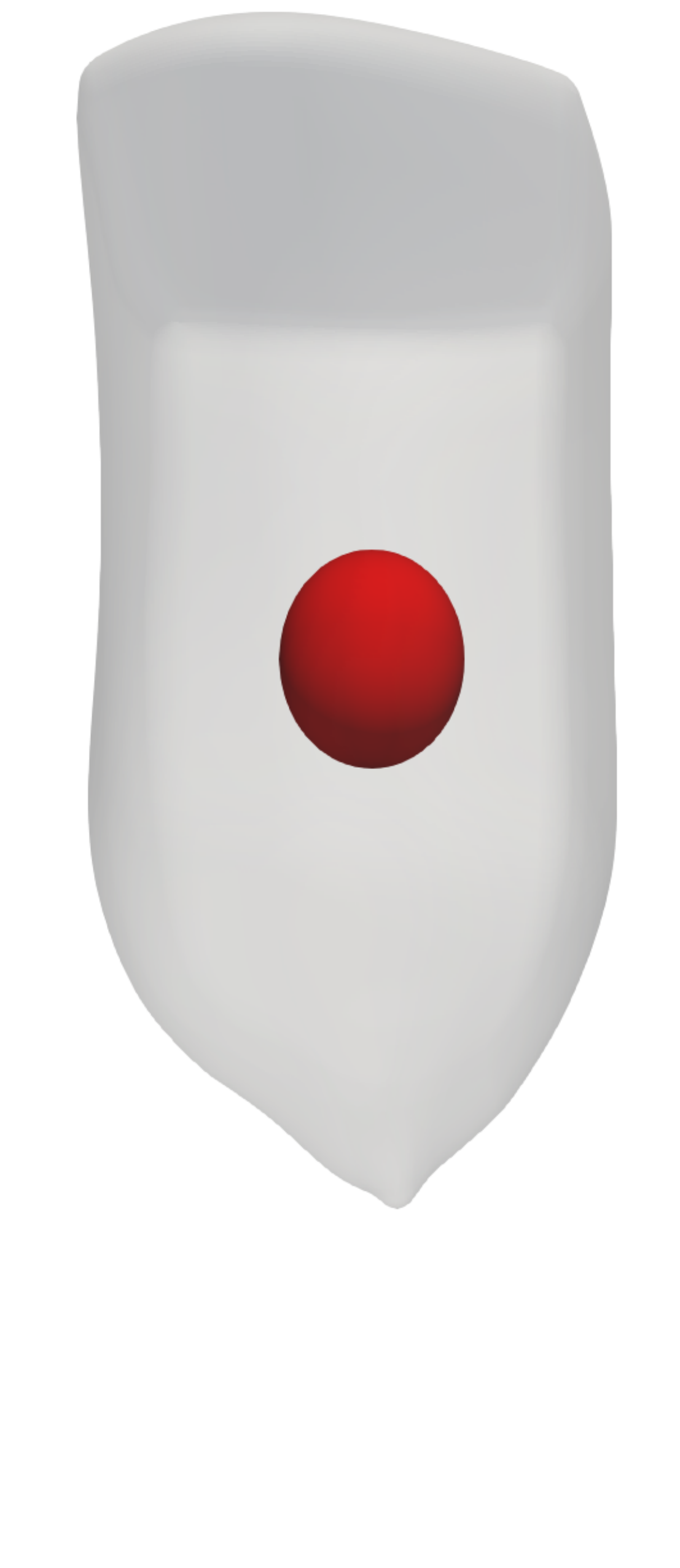}  
\end{subfigure}
\\
\begin{subfigure}{\columnwidth} \centering
\includegraphics[trim={0cm 8cm 0cm 0cm},clip,width=0.59\columnwidth]{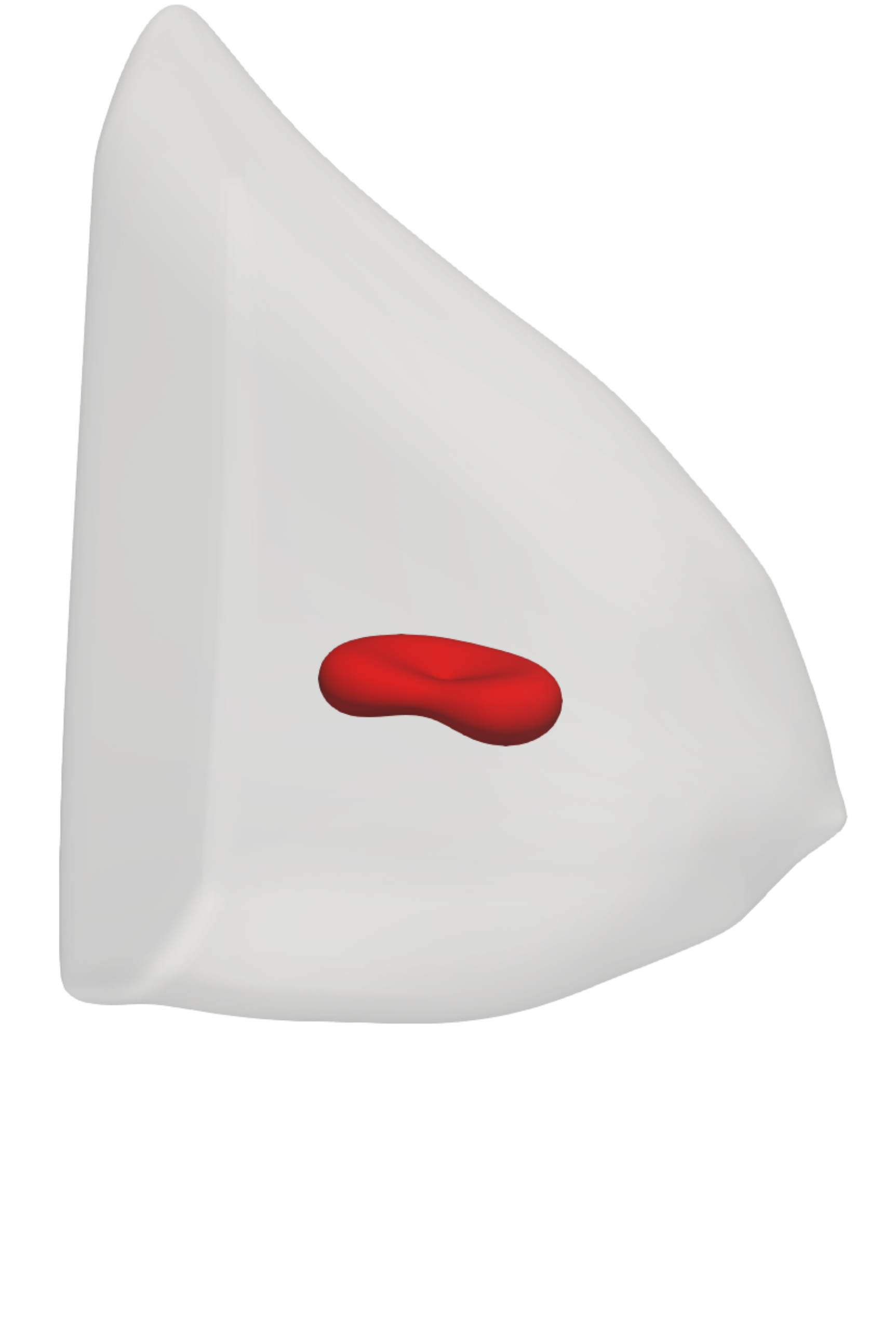}  
\includegraphics[trim={0cm 8cm 0cm 0cm},clip,width=0.39\columnwidth]{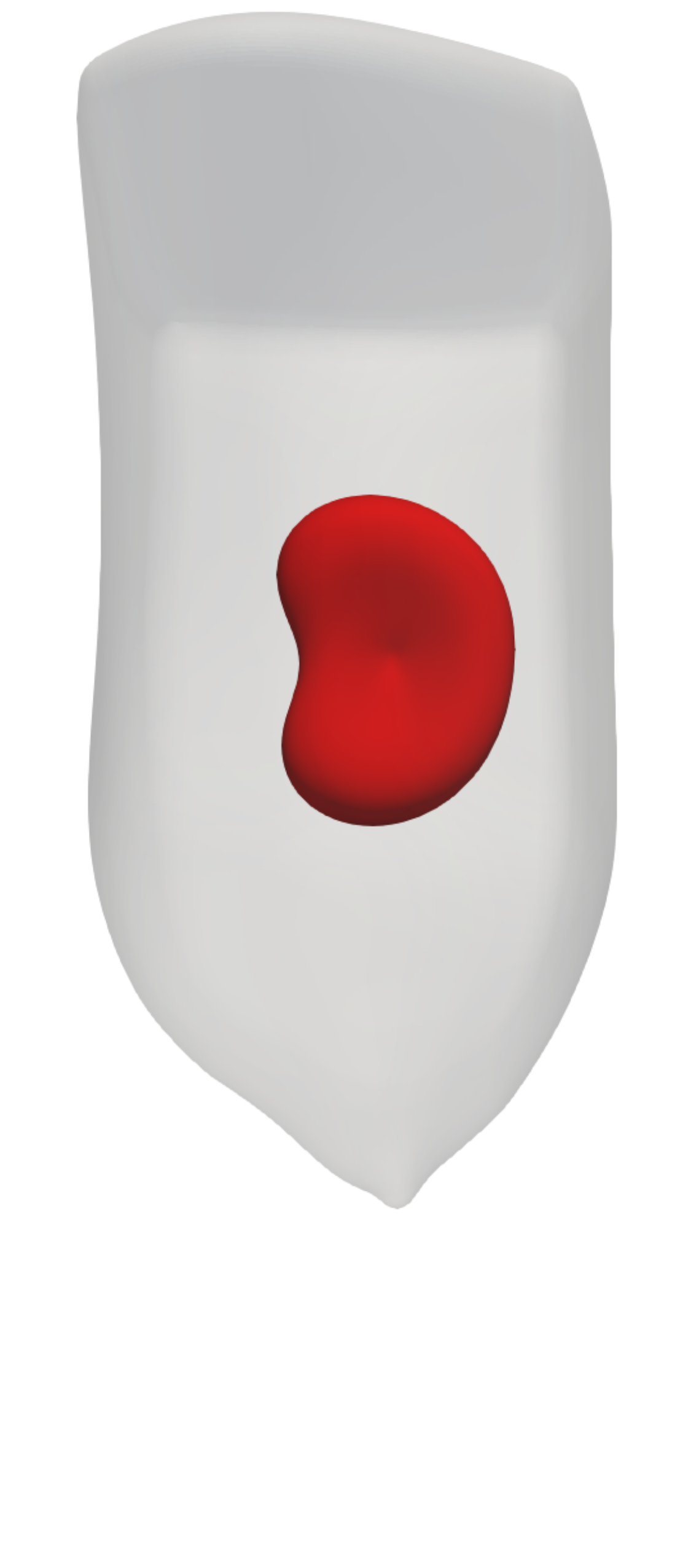}  
\end{subfigure}
\caption{Tumor morphology} \label{fig:breast2_modela} \end{subfigure}
\begin{subfigure}{0.245\columnwidth}
\begin{subfigure}{\columnwidth} \centering
\includegraphics[trim={0cm 8cm 0cm 0cm},clip,width=0.59\columnwidth]{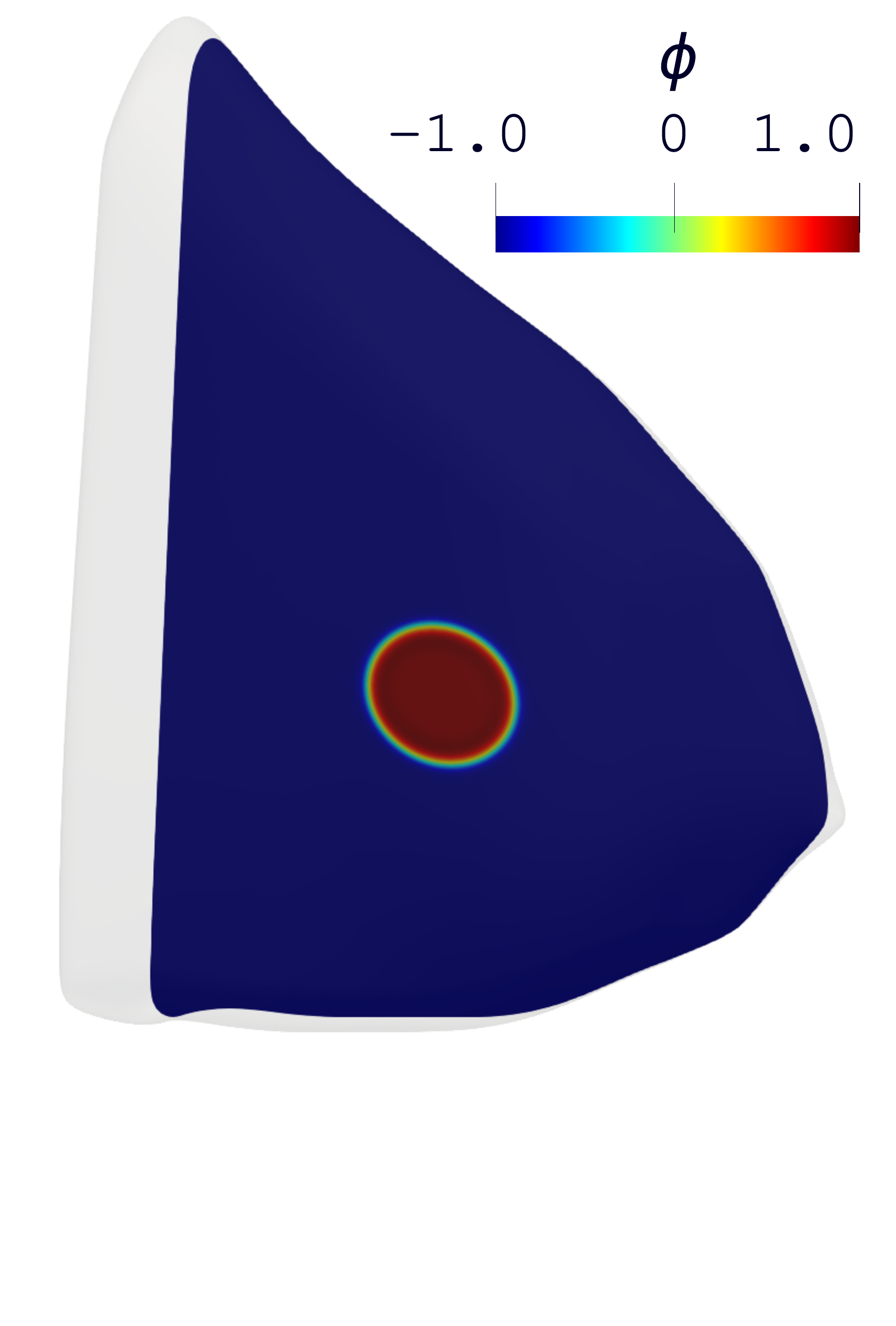}  
\includegraphics[trim={0cm 8cm 0cm 0cm},clip,width=0.39\columnwidth]{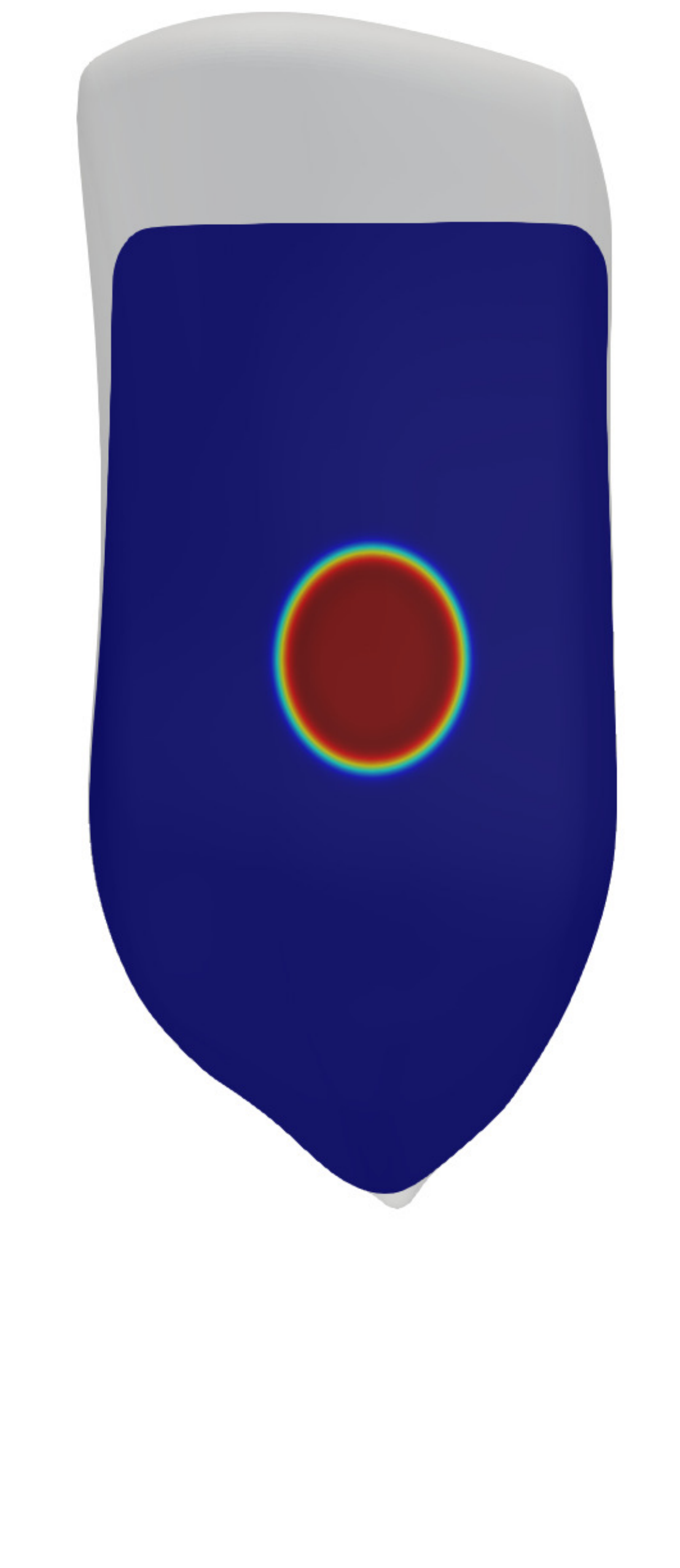}  
\end{subfigure}
\\
\begin{subfigure}{\columnwidth} \centering
\includegraphics[trim={0cm 8cm 0cm 0cm},clip,width=0.59\columnwidth]{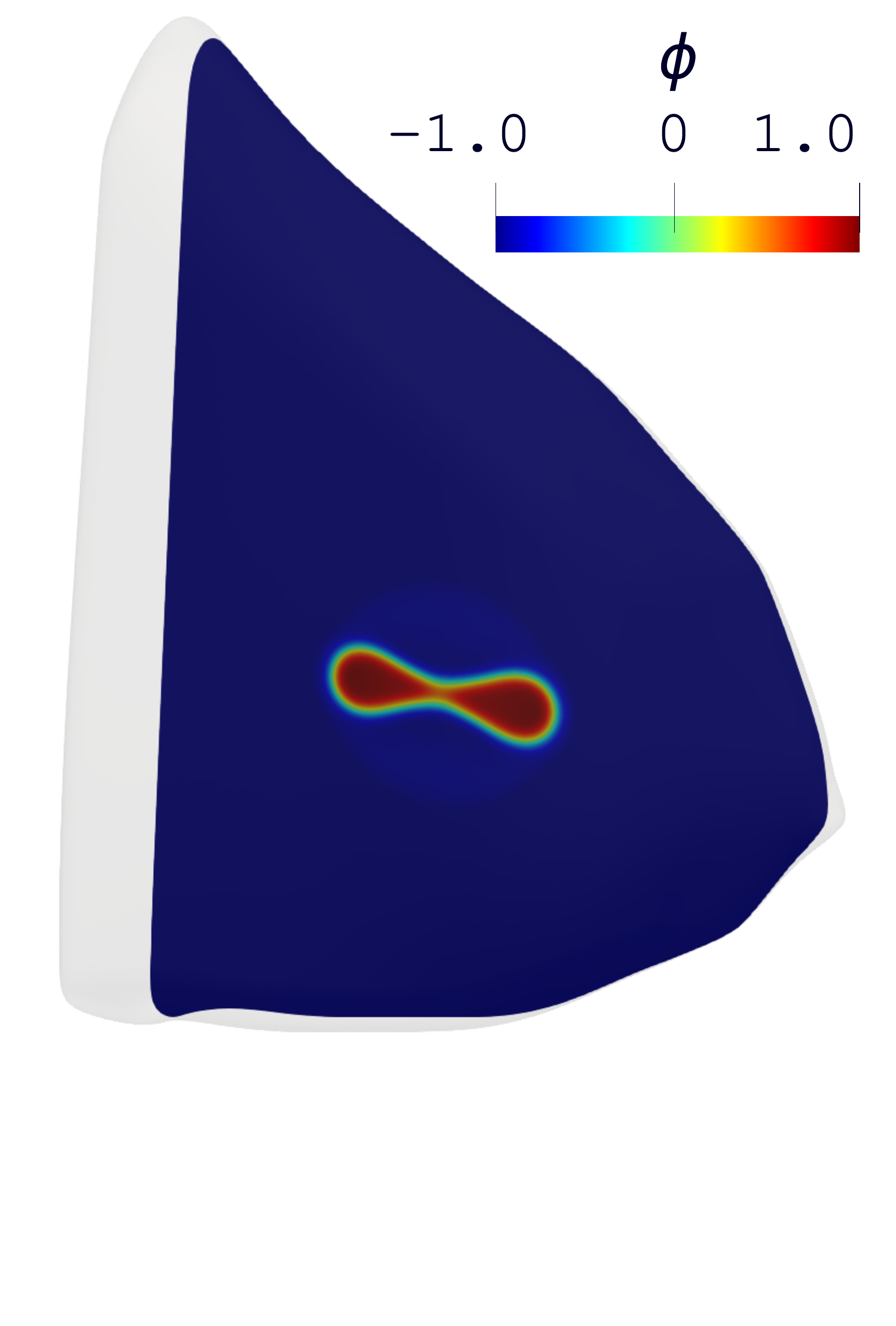}  
\includegraphics[trim={0cm 8cm 0cm 0cm},clip,width=0.39\columnwidth]{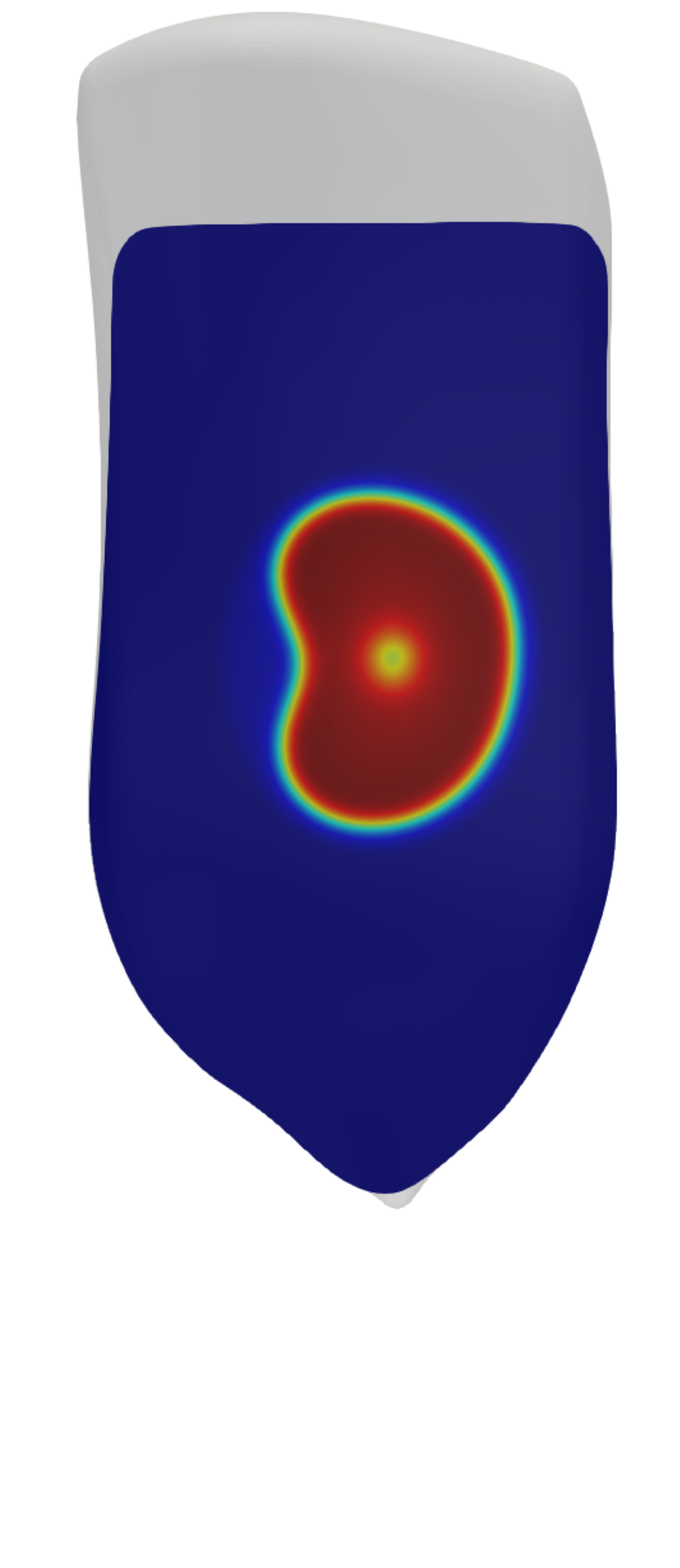}  
\end{subfigure}
\caption{Phase-field} \label{fig:breast2_modelb} \end{subfigure}
\begin{subfigure}{0.245\columnwidth}
\begin{subfigure}{\columnwidth} \centering
\includegraphics[trim={0cm 8cm 0cm 0cm},clip,width=0.59\columnwidth]{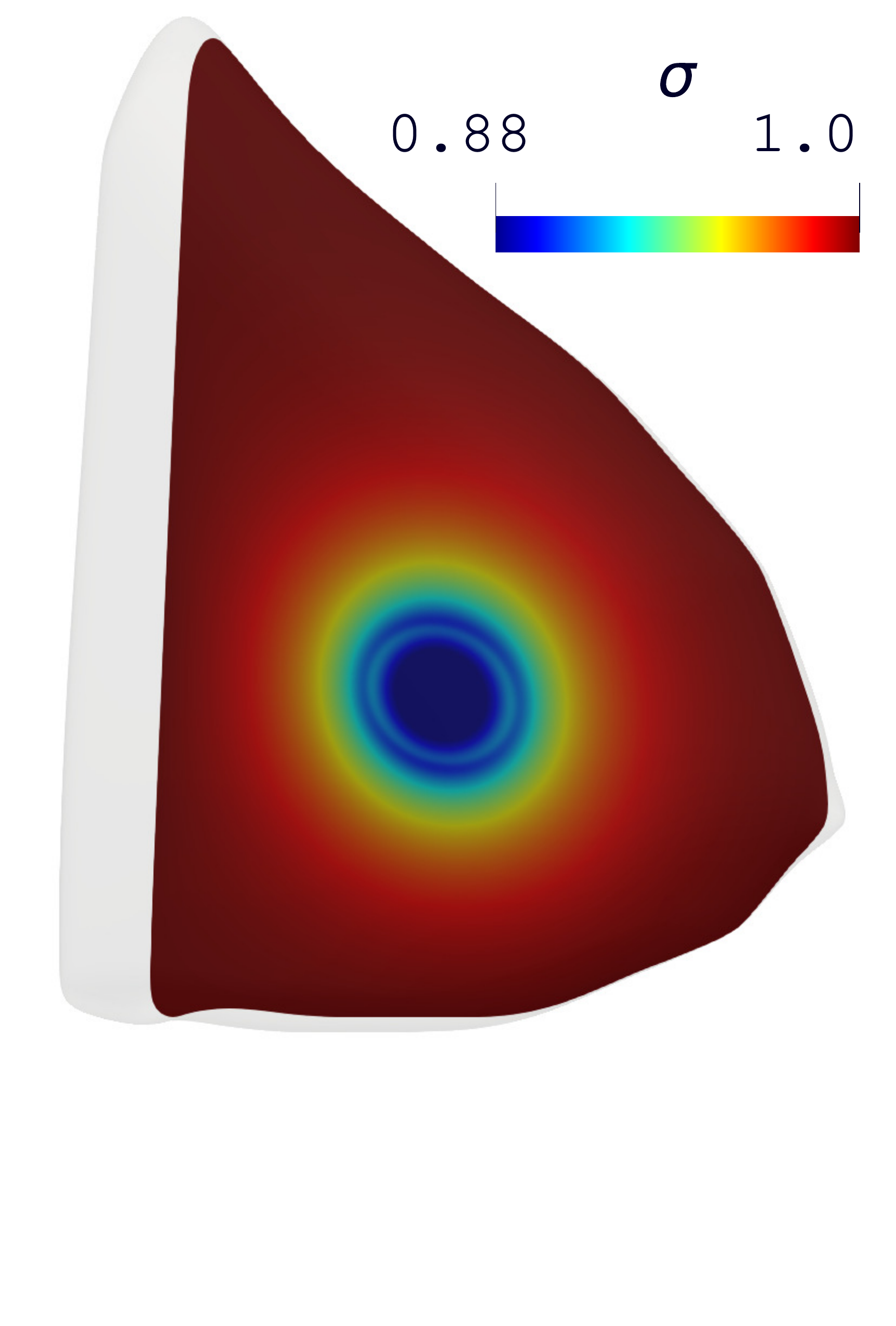}  
\includegraphics[trim={0cm 8cm 0cm 0cm},clip,width=0.39\columnwidth]{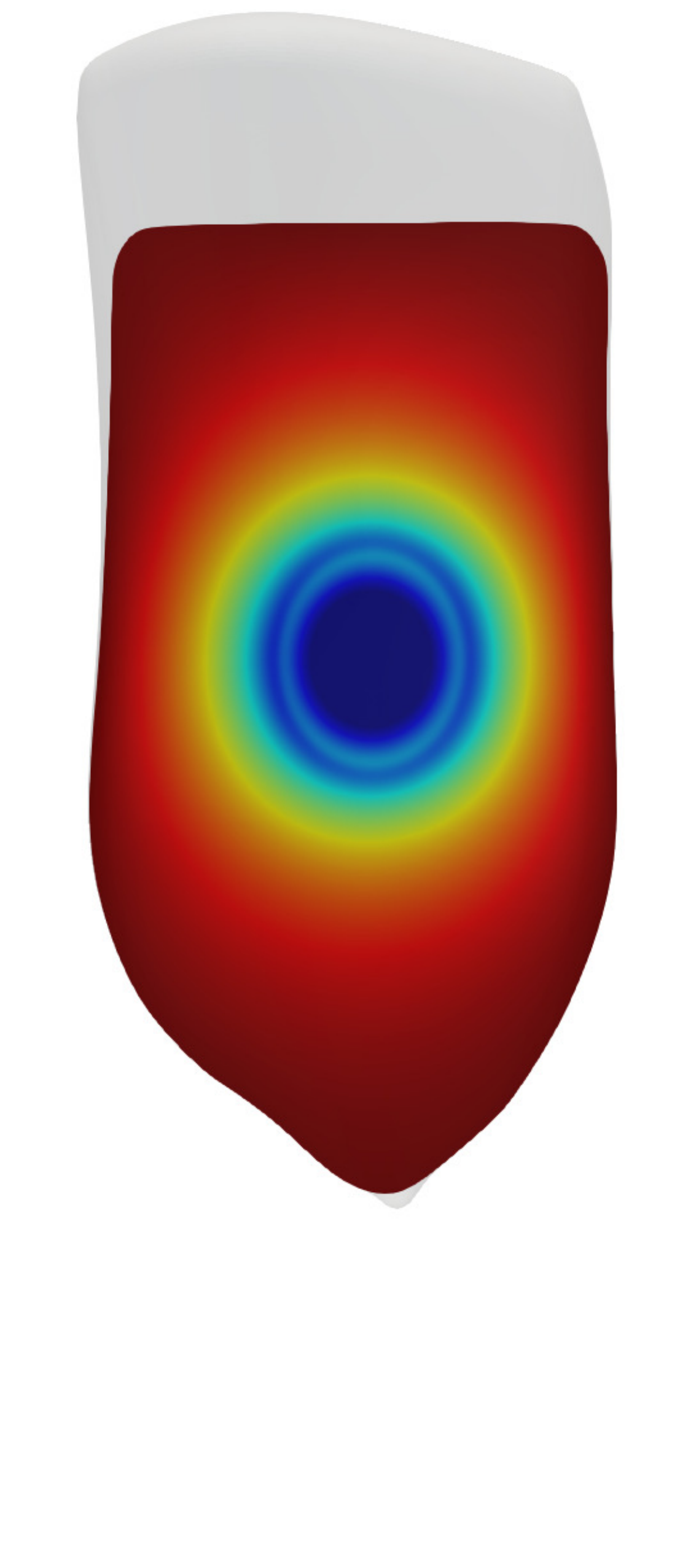}  
\end{subfigure}
\\
\begin{subfigure}{\columnwidth} \centering
\includegraphics[trim={0cm 8cm 0cm 0cm},clip,width=0.59\columnwidth]{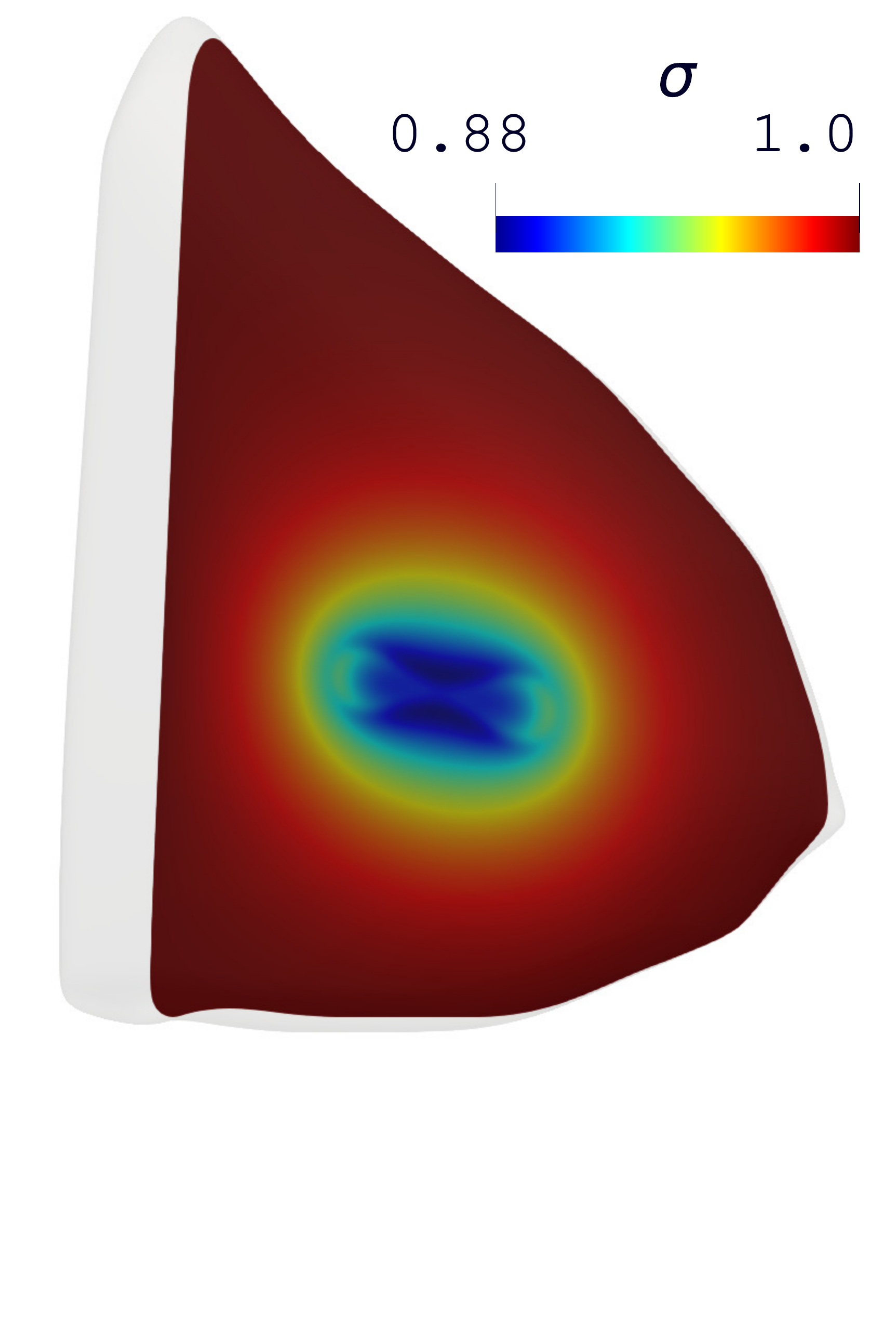}  
\includegraphics[trim={0cm 8cm 0cm 0cm},clip,width=0.39\columnwidth]{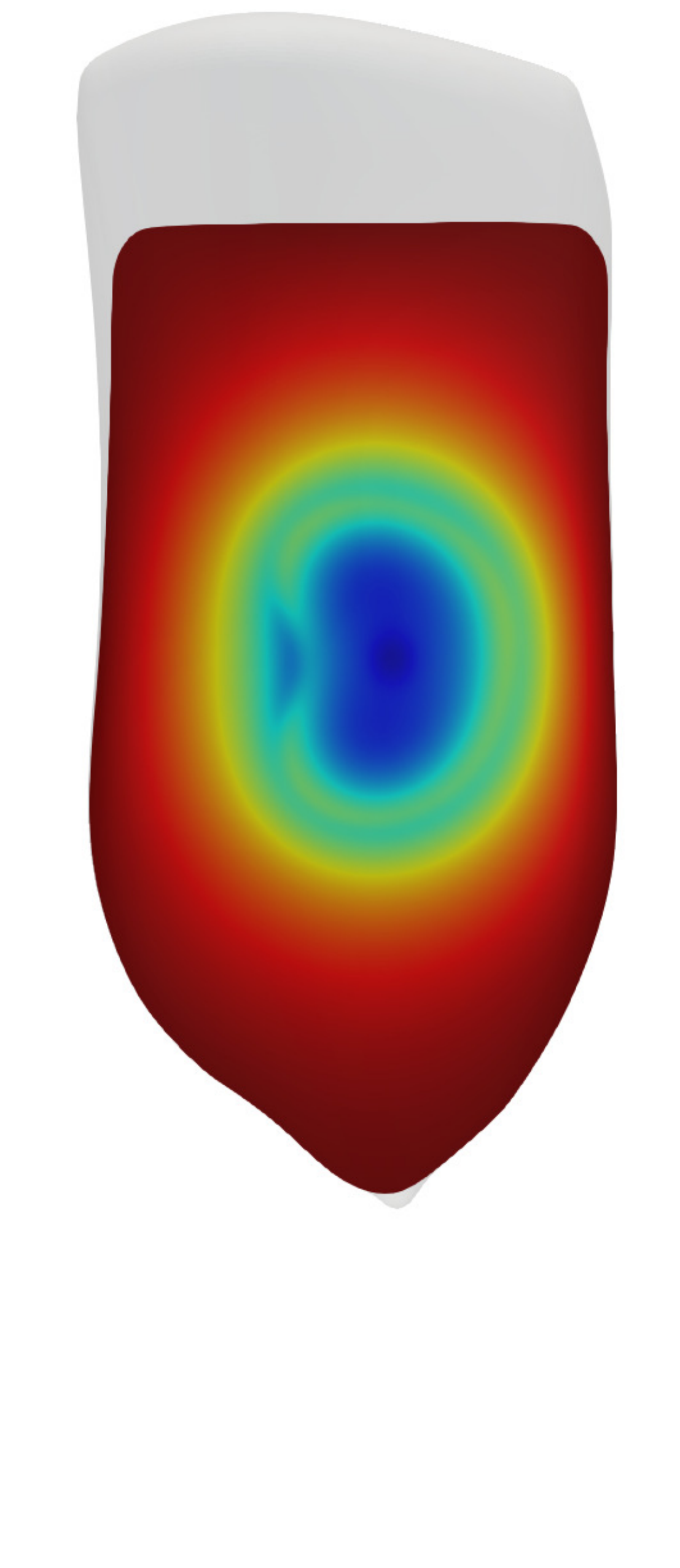}  
\end{subfigure}
\caption{Nutrient concentration} \label{fig:breast2_modelc} \end{subfigure}
\begin{subfigure}{0.245\columnwidth}
\begin{subfigure}{\columnwidth} \centering
\includegraphics[trim={0cm 0cm 0cm 0cm},clip,width=0.425\columnwidth]{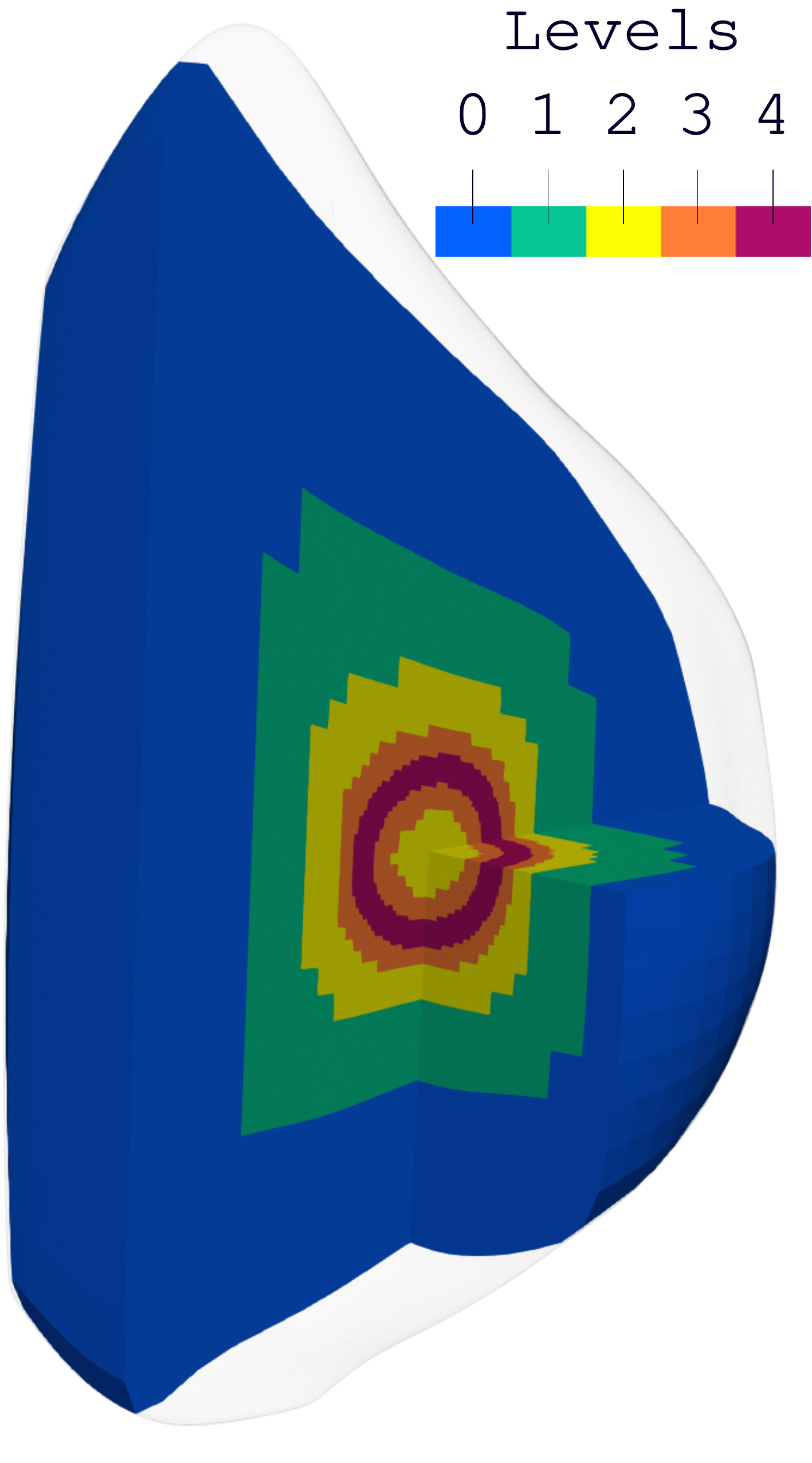}  
\includegraphics[trim={0cm 0cm 0cm 0cm},clip,width=0.425\columnwidth]{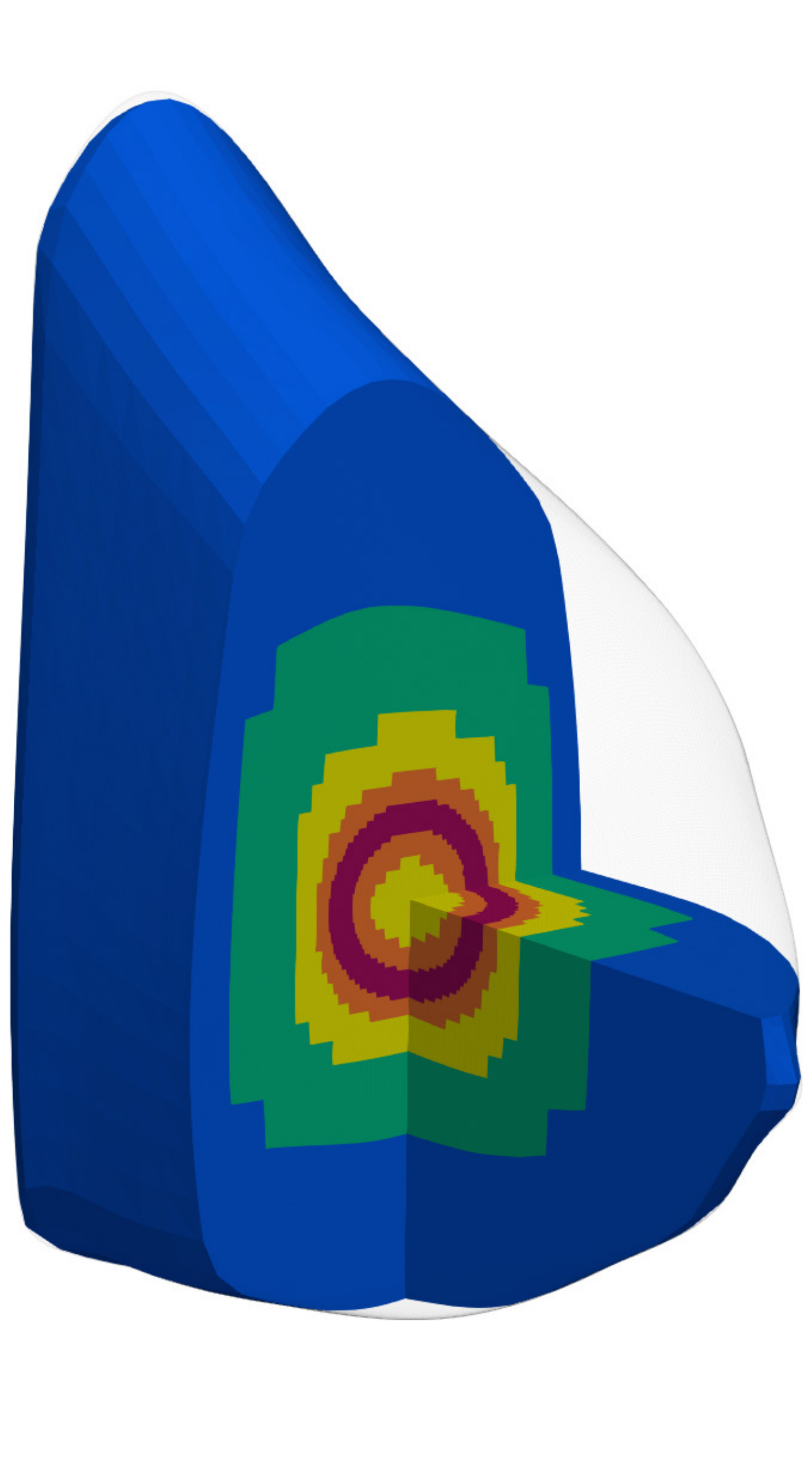}  
\end{subfigure}
\\
\begin{subfigure}{\columnwidth} \centering
\includegraphics[trim={0cm 0cm 0cm 0cm},clip,width=0.425\columnwidth]{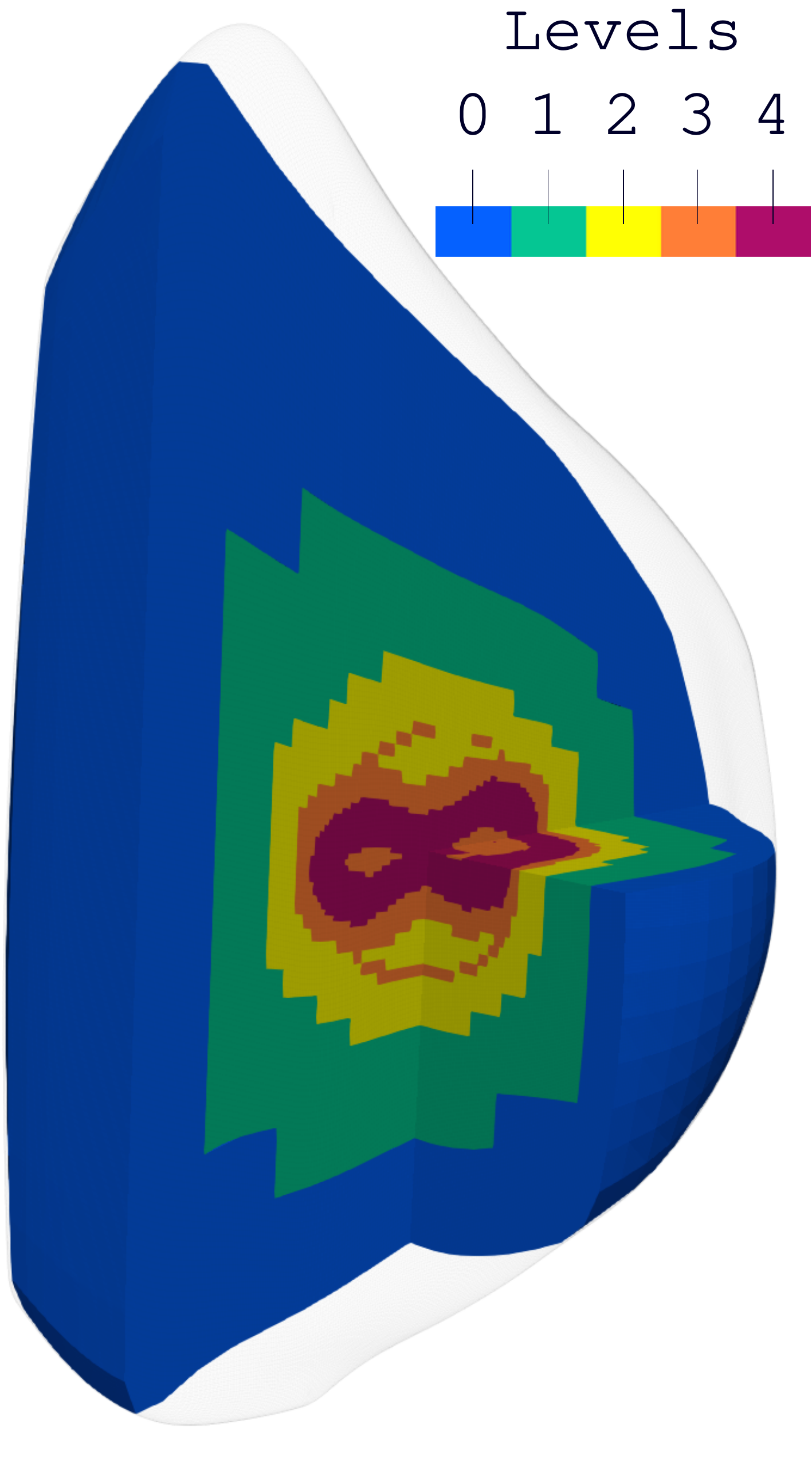}  
\includegraphics[trim={0cm 0cm 0cm 0cm},clip,width=0.425\columnwidth]{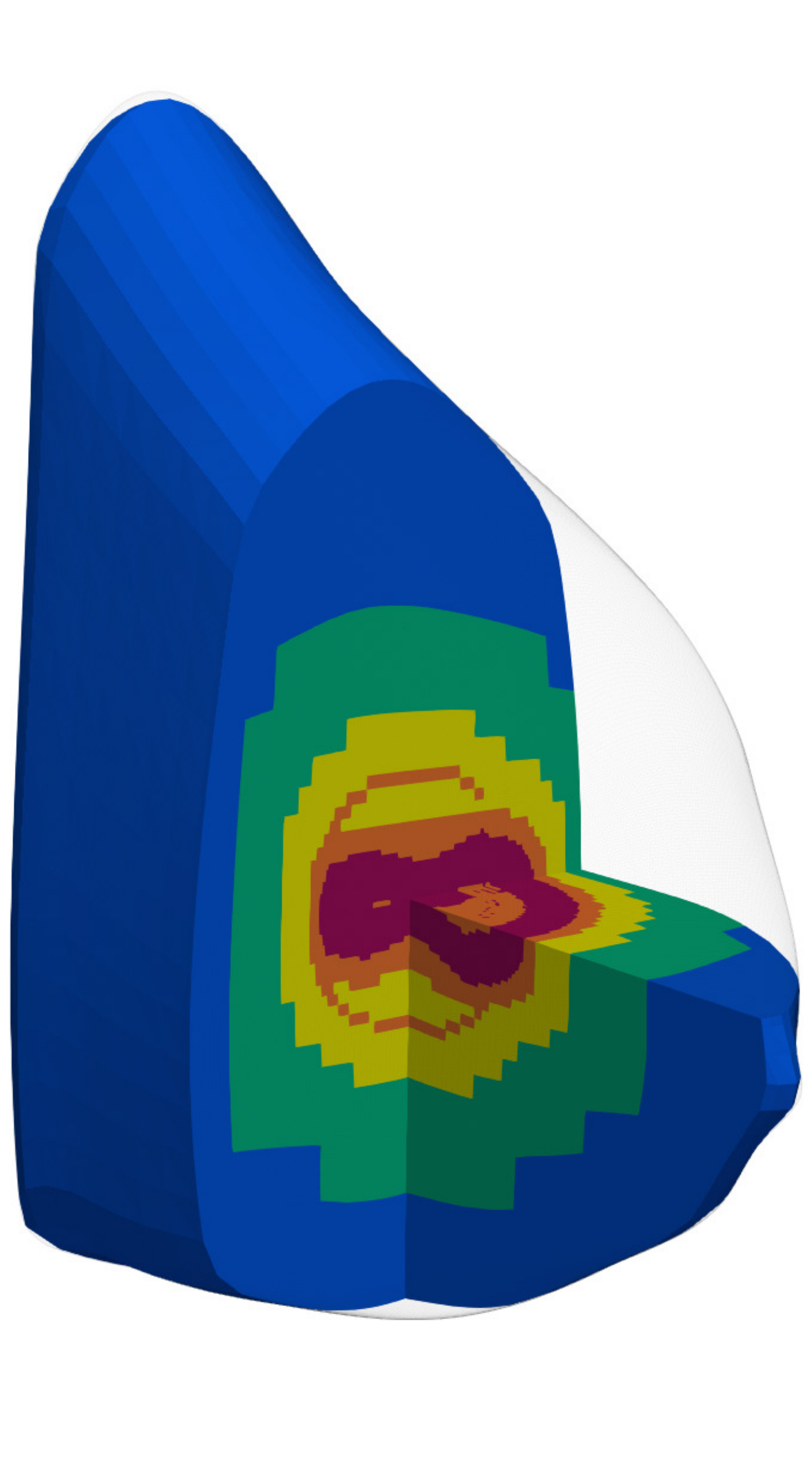}  
\end{subfigure}
\caption{Adaptive THB-spline mesh} \label{fig:breast2_modeld} \end{subfigure}
\caption{Spatiotemporal evolution of an ellipsoidal tumor in a patient-specific breast geometry at times $t=$ 0 (top row) and 0.125 (bottom row). THB-spline mesh configuration: $p=3$, $\ell = 4$, $m=2$, $\alpha = 0.01$ and $\beta = 0.0001$, with deepest refinement level at $2^8\times2^8\times2^8$. Model parameters are same as Section~\ref{sec:BrCa_patient_spec}. The initial tumor field is initialized using a hyperbolic tangent profile with an ellipsoidal shape, having two equal semi-axes of 0.11 and a third semi-axis of 0.13 and centered within the interior region of the breast geometry.}
\label{fig:breast2_model}
\end{figure}

Finally, an additional case is considered to investigate tumor evolution within the patient-specific breast geometry reconstructed from MRI data, as described in Section~\ref{sec:BrCa_patient_spec}.
Conversely to the case shown in Fig.~\ref{fig:breast1_model}, the tumor is initialized away from the domain boundaries using a hyperbolic tangent profile with an ellipsoidal shape, having two equal semi-axes of 0.11 and a third semi-axis of 0.13, and centered at the interior region of the breast geometry.
The mesh configuration and model parameters are the same as in Section~\ref{sec:BrCa_patient_spec}.
The resulting morphological changes under reduced influence from the nutrient-rich boundaries are shown in Fig.~\ref{fig:breast2_modela}.
Furthermore, Figs.~\ref{fig:breast2_modelb} and \ref{fig:breast2_modelc} show the phase field and the corresponding nutrient concentration, respectively, on two orthogonal planes passing through the center of the tumor.
Unlike the case shown in Fig.~\ref{fig:breast1_model}, the reduced influence of the nutrient-rich boundaries allows the growth-driven morphological changes to become more apparent, with migration no longer being the dominant mechanism.

\def\bibsection{\section*{References}}

\bibliographystyle{elsarticle-num} 
\bibliography{Library}

\end{document}